\pdfoutput=1

\input{tex/preamble/switches}

\documentclass[
  leqno,
  a4paper,
]{scrbook}

\usepackage[ngerman,american]{babel}

\usepackage[utf8]{inputenc}

\usepackage{csquotes}

\usepackage[markcase=upper]{scrlayer-scrpage}

\usepackage{mathtools}

\usepackage{amsthm}

\usepackage{thmtools}

\usepackage{geometry}

\iftoggle{checkMode}{
  \usepackage[mark]{lua-check-hyphen}
}{}

\iftoggle{debugMode}{
  \usepackage{lua-visual-debug}
}{}

\usepackage[
  type1,
  scale=0.88,
]{GoMono}

\usepackage{tgheros}

\usepackage[T1]{fontenc}

\makeatletter
\edef\bfseries@rm{b}
\edef\mdseries@rm{m}
\makeatother

\usepackage[bitstream-charter]{mathdesign}

\usepackage[
  final,
  tracking=smallcaps,
  stretch=10,
  shrink=10,
]{microtype}

\usepackage{lastpage}

\usepackage{refcount}

\usepackage{graphicx}

\usepackage[table]{xcolor}

\usepackage{contour}

\usepackage{cellspace}

\usepackage{subcaption}

\usepackage{sidecap}

\usepackage{array}

\usepackage{tabularx}

\usepackage{booktabs}

\usepackage{multirow}

\usepackage{enumitem}

\usepackage{expl3}
\ExplSyntaxOn
\msg_redirect_name:nnn{siunitx}{moved-cellspace-column}{none}
\ExplSyntaxOff

\usepackage[binary-units=true]{siunitx}

\usepackage{xspace}

\usepackage{everypage}

\usepackage[framemethod=tikz]{mdframed}

\usepackage{tikz}

\usepackage{tocbasic}

\usepackage{etoc}

\usepackage{wrapfig}

\usepackage{lettrine}

\usepackage{pdfpages}

\usepackage{nameref}

\usepackage[onehalfspacing]{setspace}

\usepackage{silence}

\usepackage{regexpatch}

\usepackage[hyphens]{url}

\usepackage{hyperref}

\usepackage{algorithmicx}

\usepackage[noend]{algpseudocode}

\usepackage{bookmark}

\iftoggle{partialCompileMode}{
  \usepackage{substr}
  \usepackage{longtable}
}{
  \usepackage[
    toc,
    nopostdot,
    indexonlyfirst,
    ucmark,
    nohypertypes=main,
  ]{glossaries}
}

\usepackage[
  giveninits=true,
  date=year,
  style=alphabetic,
  maxalphanames=1,
]{biblatex}

\usepackage[
  capitalise,
  nameinlink,
]{cleveref}

\input{tex/preamble/settings}

\newcommand*{\ifempty}[2]{\if\relax\detokenize{#1}\relax#2\fi}
\newcommand*{\ifnotempty}[2]{\if\relax\detokenize{#1}\relax\else#2\fi}
\newcommand*{\appendwithcomma}[1]{\ifnotempty{#1}{,#1}}

\newnotationcommand{\nat}{\mathbb{N}}{N30}{$\nat$}{
  Natural numbers without zero ($1, 2, 3, \dotsc$)
}

\newnotationcommand{\natz}{\mathbb{N}_0}{N40}{$\natz$}{
  Natural numbers with zero ($\nat \cup \{0\}$)
}

\newnotationcommand{\integer}{\mathbb{Z}}{Z}{$\integer$}{
  Integer numbers ($\dotsc, -2, -1, 0, 1, 2, \dotsc$)
}

\newnotationcommand{\real}{\mathbb{R}}{R}{$\real$}{
  Real numbers
}

\newnotationcommand{\posreal}{\mathbb{R}_{{} > 0}}{R>0}{$\posreal$}{
  Positive real numbers
}

\newnotationcommand{\nonnegreal}{\mathbb{R}_{{} \ge 0}}{R>=0}{$\nonnegreal$}{
  Non-negative real numbers
}

\newnotation{l}{$\*l$}{
  Hierarchical level $\*l \in \natz^d$
}

\newnotation{i}{$\*i$}{
  Hierarchical index $\*i = \*0, \dotsc, \*2^\*l$
}

\newnotation{kT00}{$\*k$}{
  Continuously enumerated index $\*k \in \natz^d$
  of the level-index pairs $(\*l, \*i)$
}

\newnotation{d}{$d$}{
  Dimensionality $d \in \nat$
}

\newnotation{n10}{$n$}{
  Level $\in \natz$ of full or sparse grid
}

\newnotationcommand{\ngp}{N}{N20}{$\ngp$}{
  Number of grid points
  for a finite set $\sgset \subset \clint{\*0, \*1}$ of grid points
}

\newnotationcommandoptarg{1}{\vlinin}{%
  \*u\ifnotempty{#1}{^{#1}}}{u}{$\vlinin$%
}{
  Input of the linear operator $\linop$
}

\newcommand*{\linin}[2][]{u_{#2}\ifnotempty{#1}{^{#1}}}

\newnotationcommandoptarg{1}{\vlinout}{%
  \*y\ifnotempty{#1}{^{#1}}}{y}{$\vlinout$%
}{
  Output of the linear operator $\linop$
}

\newcommand*{\linout}[2][]{y_{#2}\ifnotempty{#1}{^{#1}}}

\newnotation{@0}{$\*0$}{
  $(0, \dotsc, 0) \in \natz^d$
}

\newnotation{@1}{$\*1$}{
  $(1, \dotsc, 1) \in \nat^d$
}

\newnotation{kT10}{$\*k_{-t}$}{
  Vector of all entries $k_{t'}$ except the $t$-th
}

\newnotation{kT20}{$\*k_T$}{
  Vector of all dimensions that are contained in $T \in \{1, \dotsc, d\}^j$
}

\newnotation{kT30}{$\*k_{-T}$}{
  Vector of all dimensions that are not contained in $T \in \{1, \dotsc, d\}^j$
}

\newnotationcommand[2]{\range}{{{#1}:{#2}}}{kT40}{$\*k_{\range{a}{b}}$}{
  Vector of all dimensions that are contained in $a, a + 1, \dotsc, b$
}

\newnotationcommand[1]{\stdbasis}{\*e_{#1}}{et}{$\stdbasis{t}$}{
  $t$-th standard basis vector
  $\stdbasis{t} \ceq (0, \dotsc, 0, 1, 0, \dotsc, 0) \in \real^d$
}

\hidenextnotation
\newnotationcommandoptarg[\*k]{2}{\chain}{#1^{(#2)}}{kt}{$\chain{t_{j'}}$}{
  Point on chain $(\chain{t_1}, \dotsc, \chain{t_j})$
  between $\*k'$ and $\*k''$
}

\newcommand*{\chainuv}[1]{\chain[k]{#1}}

\newnotationcommand[1]{\gp}{%
  \containsvector{#1}{\*x}{x}_{#1}%
}{xli}{$\gp{\*l,\*i}$}{
  Grid point $\gp{\*l,\*i} \ceq \*i \cdot \ms{\*l}$
}

\newcommand*{\vgp}[1]{\*x_{#1}}

\newnotationcommand[1]{\ccgp}{%
  \containsvector{#1}{\*x}{x}_{#1}^\cc%
}{xlicc}{$\ccgp{\*l,\*i}$}{
  Clenshaw--Curtis grid point
}

\newnotationcommand[1]{\ms}{%
  \containsvector{#1}{\*h}{h}_{#1}%
}{hl}{$\ms{\*l}$}{
  Mesh size $\ms{\*l} \ceq \*2^{-\*l}$
}

\newcommand*{\cc}{\mathrm{cc}}
\newcommand*{\modified}{\mathrm{mod}}
\newcommand*{\ntrl}{\mathrm{nat}}
\newcommand*{\nak}{\mathrm{nak}}
\newcommand*{\sparse}{\mathrm{s}}
\newcommand*{\ct}{\mathrm{ct}}
\newcommand*{\hft}{\mathrm{hft}}
\newcommand*{\tift}{\mathrm{tift}}
\newcommand*{\ext}{\mathrm{ext}}
\newcommand*{\fs}{\mathrm{fs}}
\newcommand*{\wfs}{\mathrm{wfs}}
\newcommand*{\chol}{\mathrm{chol}}
\newcommand*{\reference}{\mathrm{ref}}
\newcommand*{\opt}{\mathrm{opt}}

\newnotation{zzzz1}{$(\cdot)^1$}{
  Superscript for ``Piecewise linear''
}

\newnotation{zzzzp}{$(\cdot)^p$}{
  Superscript for ``B-splines of degree $p$''
}

\newnotation{zzzzs}{$(\cdot)^\sparse$}{
  Superscript for ``Sparse grid''
}

\newnotation{zzzzsb}{$(\cdot)^{\sparse(b)}$}{
  Superscript for ``Sparse grid with boundary parameter $b$''
}

\newnotation{zzzzmod}{$(\cdot)^\modified$}{
  Superscript for ``Modified''
}

\newnotation{zzzzcc}{$(\cdot)^\cc$}{
  Superscript for ``Clenshaw--Curtis''
}

\newnotation{zzzznak}{$(\cdot)^\nak$}{
  Superscript for ``Not-a-knot''
}

\newnotation{zzzzfs}{$(\cdot)^\fs$}{
  Superscript for ``Fundamental spline''
}

\newnotation{zzzzwfs}{$(\cdot)^\wfs$}{
  Superscript for ``Weakly fundamental spline''
}

\newnotation{zzzzopt}{$(\cdot)^\opt$}{
  Superscript for ``Optimal''
}

\newnotationcommand[1]{\clint}{[#1]}{zzzzzzab10}{$\clint{\*a, \*b}$}{
  Closed hyper-rectangle $\{\*x \in \real^d \mid \*a \le \*x \le \*b\}$
}

\newnotationcommand[1]{\opint}{\mathopen]#1\mathclose[}{zzzzzzab20}{%
  $\opint{\*a, \*b}$%
}{
  Open hyper-rectangle $\{\*x \in \real^d \mid \*a < \*x < \*b\}$
}

\newnotationcommand[1]{\hopint}{\mathopen[#1\mathclose[}{zzzzzzab30}{%
  $\hopint{\*a, \*b}$%
}{
  Half-open hyper-rectangle $\{\*x \in \real^d \mid \*a \le \*x < \*b\}$
}

\hidenextnotation
\newnotationcommand[1]{\pohint}{\mathopen]#1\mathclose]}{zzzzzzab40}{%
  $\pohint{\*a, \*b}$%
}{
  Half-open hyper-rectangle $\{\*x \in \real^d \mid \*a < \*x \le \*b\}$
}

\newcommand*{\opintscaled}[1]{\left]#1\right[}

\newnotationcommand[1]{\normone}{\norm[1]{#1}}{zzzzzznorm1}{$\normone{\cdot}$}{
  $\ell_1$ norm $\normone{\*x} \ceq \sum_{t=1}^d \abs{x_t}$
}

\newcommand*{\Ltwo}{L^2}
\newcommand*{\Linfty}{L^\infty}

\newnotationcommand[1]{\normLtwo}{\norm[\Ltwo]{#1}}{zzzzzznormL2}{%
  $\normLtwo{\cdot}$%
}{
  $\Ltwo$ norm
  $\normLtwo{\objfun} \ceq \sqrt{\int_\Omega \objfun(\*x)^2 \diff{}\*x}$
  of a function $\objfun\colon \Omega \to \real$
}

\newcommand*{\normLtwoscaled}[1]{\normscaled[\Ltwo]{#1}}

\newnotationcommand[1]{\normLinfty}{\norm[\Linfty]{#1}}{zzzzzznormLinfty}{%
  $\normLinfty{\cdot}$%
}{
  $\Linfty$ norm
  $\normLinfty{\objfun} \ceq \max_{\*x \in \Omega} \abs{\objfun(\*x)}$
  of a continuous function $\objfun\colon \Omega \to \real$
}

\newcommand*{\normLinftyscaled}[1]{\normscaled[\Linfty]{#1}}

\let\equivorig\equiv
\let\equiv\undefined
\newnotationcommand{\equiv}{\equivorig}{zzzzzzequiv}{$\equivorig$}{
  Equality of functions everywhere on their domain
  (i.e., $\fa{\*x}{\objfun_1(\*x) = \objfun_2(\*x)}$)
}

\newcommand*{\gridset}[3][]{%
  \if\relax\detokenize{#1}\relax\Omega\else#1{\Omega}\fi_{#2}^{#3}%
}

\newcommand*{\interiorgrid}[1]{\mathring{#1}}

\newcommand*{\interiorsgset}{\gridset[\interiorgrid]{}{\sparse}}

\makecommandnotation{\interiorgrid}{zzZs2}{$\interiorsgset$}{
  Set $\ceq \sgset \cap \opint{\*0, \*1}$ of
  interior grid points
  for a finite set $\sgset \subset \clint{\*0, \*1}$ of grid points
}

\newnotationcommandoptarg{2}{\fgset}{%
  \Omega_{#2}\ifnotempty{#1}{^{#1}}%
}{zzZl}{$\fgset{\*l}$}{
  Set of full grid points of level $\*l$
}

\newnotationcommandoptarg{1}{\sgset}{%
  \gridset{}{\sparse\appendwithcomma{#1}}%
}{zzZs}{$\sgset$}{
  Arbitrary sparse grid (possibly spatially adaptive)
}

\newnotationcommandoptarg{3}{\regsgset}{%
  \gridset{#2,#3}{\sparse\appendwithcomma{#1}}%
}{zzZsnd}{$\regsgset{n}{d}$}{
  Set of regular sparse grid points of level $n$, dimensionality $d$
}

\newnotationcommandoptarg{4}{\coarseregsgset}{%
  \gridset{#2,#3}{\sparse(#4)\appendwithcomma{#1}}%
}{zzZsndb}{$\coarseregsgset{n}{d}{b}$}{
  Set of regular sparse grid points of level $n$, dim. $d$,
  boundary parameter $b$
}

\newcommand*{\interiorregsgset}[3][]{%
  \gridset[\interiorgrid]{#2,#3}{\sparse\appendwithcomma{#1}}%
}

\newcommand*{\gridspace}[3][]{%
  \if\relax\detokenize{#1}\relax V\else#1{V}\fi_{#2}^{#3}%
}

\newnotationcommandoptarg{2}{\ns}{%
  V_{#2}\ifnotempty{#1}{^{#1}}%
}{Vl}{$\ns{\*l}$}{
  Nodal space of level $\*l$
}

\newcommand*{\nsbspl}[3][]{\ns[#3\appendwithcomma{#1}]{#2}}

\newnotation{Vnd}{$\ns{n,d}$}{
  Multivariate nodal space
  $\ceq \ns{n \cdot \*1}$ of level $n$ with dimensionality $d$
}

\newnotationcommandoptarg{2}{\hs}{%
  W_{#2}\ifnotempty{#1}{^{#1}}%
}{Wl}{$\hs{\*l}$}{
  Hierarchical subspace of level $\*l$
}

\newcommand*{\hsbspl}[3][]{\hs[#3\appendwithcomma{#1}]{#2}}

\newnotationcommand{\sgspace}{\gridspace{}{\sparse}}{Vs}{$\sgspace$}{
  Arbitrary sparse grid space (possibly spatially adaptive)
}

\newnotationcommandoptarg{3}{\regsgspace}{%
  \gridspace{#2,#3}{\sparse\appendwithcomma{#1}}%
}{Vnds}{$\regsgspace{n}{d}$}{
  Regular sparse grid space of level $n$ with dimensionality $d$
}

\hidenextnotation
\newnotationcommand[2]{\nonunifsplspace}{S_{#1}^{#2}}{Slop}{%
  $\nonunifsplspace{\knotseq}{p}$%
}{
  Spline space of degree $p$ with knots $\knotseq$
}

\hidenextnotation
\newnotationcommand[2]{\restrictedsplspace}{S_{#1}^{#2}}{Slp}{%
  $\restrictedsplspace{l}{p}$%
}{
  Spline space $\ceq \nonunifsplspace{\nodalknotseq{l}{p}}{p}$
  (spanned by the uniform nodal B-splines of level $l$, degree $p$)
}

\newnotationcommand[2]{\wholesplspace}{S_{#1}^{#2,\clint{0,1}}}{Slp01}{%
  $\wholesplspace{l}{p}$%
}{
  Spline space of degree $p$ on the grid
  $\{\gp{l,i} \mid i = 0, \dotsc, 2^l\}$ of level $l$
}

\hidenextnotation
\newnotationcommand[2]{\naksplspace}{S_{#1}^{#2,\nak}}{Slpnak}{%
  $\naksplspace{l}{p}$%
}{
  Spline space $\ceq \nonunifsplspace{\nodalknotseq[\nak]{l}{p}}{p}$
  (spanned by the nodal not-a-knot B-splines of level $l$, degree $p$)
}

\newnotationcommand[1]{\polyspace}{P^{#1}}{Pp}{$\polyspace{\*p}$}{
  Space of all $d$-variate polynomials of
  coordinate degree $\le \*p$ on $\clint{\*0, \*1}$
}

\newnotationcommand[2]{\restrictspace}{#1|_{#2}}{V|D}{%
  $\restrictspace{V}{D}$%
}{
  Restriction $\ceq \{\restrictfcn{\objfun}{D} \mid \objfun \in V\}$
  onto a sub-domain $D$ for a function space $V$
}

\newnotationcommand[1]{\landauO}{\mathcal{O}(#1)}{O}{$\landauO{\objfun(x)}$}{
  Big-$\mathcal{O}$ Landau notation
}

\newnotationcommand[1]{\landauOmega}{\Omega(#1)}{zzZ}{%
  $\landauOmega{\objfun(x)}$%
}{
  Big-$\Omega$ Landau notation
}

\newnotationcommand[1]{\landauTheta}{\Theta(#1)}{zzH}{%
  $\landauTheta{\objfun(x)}$%
}{
  Big-$\Theta$ Landau notation
}

\makecommandnotation{\dim}{dim}{$\dim$}{
  Vector space dimension
}

\hidenextnotation
\newnotationcommand[1]{\abs}{|#1|}{zzzzzzabs}{$\abs{\cdot}$}{
  Absolute value of a scalar
}

\hidenextnotation
\newnotationcommand[1]{\vabs}{%
  \boldsymbol{|}#1\boldsymbol{|}%
}{zzzzzzabs2}{$\vabs{\cdot}$}{
  Absolute value of a vector (entry-wise)
}

\hidenextnotation
\newnotationcommand[1]{\setsize}{|#1|}{zzzzzzcount}{$\setsize{\cdot}$}{
  Number of elements of a set
}

\newcommand*{\bndrydomain}[1]{\mathop{}\!\partial#1}
\let\bndrydomainorig\bndrydomain
\makecommandnotation{\bndrydomain}{zzZd}{$\bndrydomainorig{\Omega}$}{
  Topological boundary of a set $\Omega \subset \real^d$
}

\newnotationcommand{\linop}{\mathfrak{L}}{L2}{$\linop$}{
  Linear operator $\linop\colon \real^{\ngp} \to \real^{\ngp}$
  on grid point data
}

\newnotationcommand{\intpmat}{\mat{A}}{A}{$\intpmat$}{
  Interpolation matrix with entries $\basis{\*l',\*i'}(\gp{\*l,\*i})$
}

\newcommand*{\intpmatuv}[1]{\mat{A}^{(#1)}}

\hidenextnotation
\newnotationcommand[1]{\upop}{\linop^{(#1)}}{Lt1}{%
  $\upop{t_1, \dotsc, t_j}$%
}{
  Operator $\upop{t_1, \dotsc, t_j}\colon \real^{\ngp} \to \real^{\ngp}$
  corresponding to the application of the unidirectional principle in the
  dimensions $t_1, \dotsc, t_j$ on the $d$-dimensional sparse grid
}

\hidenextnotation
\newnotationcommand[2]{\upopuv}{\linop^{(#1),#2}}{LtKpole}{%
  $\upopuv{t}{\lisetpole}$%
}{
  One-dimensional operator $\upopuv{t}{\lisetpole}\colon
  \real^{\setsize{\lisetpole}} \to \real^{\setsize{\lisetpole}}$
  corresponding to the application of the unidirectional principle in the
  $t$-th dimension on a pole $\lisetpole$ of the sparse grid
}

\newcommand*{\linopinv}{\linop^{-1}}
\newcommand*{\intpmatinv}{\intpmat^{-1}}
\newcommand*{\intpmatuvinv}[1]{(\intpmatuv{#1})^{-1}}
\newcommand*{\upopinv}[1]{(\upop{#1})^{-1}}

\newnotationcommand{\idop}{\mathrm{id}}{id}{$\idop$}{
  Identity operator $\idop\colon \real^{\ngp} \to \real^{\ngp}$
}

\hidenextnotation
\newnotationcommand[2]{\intpmatentry}{A_{#1,#2}}{Akk'}{$\intpmatentry{k}{k'}$}{
  $(k, k')$-th entry of $\intpmat$
}

\DeclareMathOperator{\supp}{supp}
\let\supporig\supp
\makecommandnotation{\supp}{supp}{$\supporig \objfun$}{
  Support of a function (i.e., the closure of $\interiorsupporig \objfun$)
}

\DeclareMathOperator{\interiorsupp}{\mathring{\supporig}}
\let\interiorsupporig\interiorsupp
\makecommandnotation{\interiorsupp}{suppo}{$\interiorsupporig \objfun$}{
  Interior of the support of a continuous function
  (i.e., $\{\*x \mid \objfun(\*x) \not= 0\}$)
}

\DeclareMathOperator{\spn}{span}
\let\spnorig\spn
\makecommandnotation{\spn}{span}{$\spnorig$}{
  Linear span (set of all linear combinations)
}

\hidenextnotation
\makecommandnotation{\deg}{deg}{$\deg \objfun$}{
  Degree of the polynomial $\objfun$
}

\DeclareMathOperator*{\vecmin}{\mathbf{min}}
\DeclareMathOperator*{\vecmax}{\mathbf{max}}

\DeclareMathOperator*{\argmin}{arg\,min}
\DeclareMathOperator*{\vecargmin}{\mathbf{arg\,min}}
\let\vecargminorig\vecargmin
\hidenextnotation
\makecommandnotation{\vecargmin}{argmin}{$\vecargminorig$}{
  Point where the minimum is attained
}

\DeclareMathOperator*{\vecargmax}{\mathbf{arg\,max}}
\let\vecargmaxorig\vecargmax
\hidenextnotation
\makecommandnotation{\vecargmax}{argmax}{$\vecargmaxorig$}{
  Point where the maximum is attained
}

\newcommand*{\floor}[1]{\lfloor#1\rfloor}

\newnotationcommand[1]{\ceil}{\lceil#1\rceil}{zzzzzzceil}{%
  $\floor{\cdot}$, $\ceil{\cdot}$%
}{
  Floor/ceiling function
  (greatest/smallest integer $\le$/$\ge$ than $\cdot$)
}

\hidenextnotation
\newnotationcommandoptarg{3}{\deriv}{%
  \frac{\diff\ifnotempty{#1}{^{#1}}}{\diff{}#2\ifnotempty{#1}{^{#1}}} #3%
}{zzzzzzderivative}{$\deriv[q]{x}{\objfun}$}{
  $q$-th derivative of $\objfun$ with respect to $x$
}

\hidenextnotation
\newnotationcommandoptarg{3}{\tderiv}{%
  \tfrac{\diff\ifnotempty{#1}{^{#1}}}{\diff{}#2\ifnotempty{#1}{^{#1}}} #3%
}{zzzzzzderivativesmall}{$\tderiv[q]{x}{\objfun}$}{
  $q$-th derivative of $\objfun$ with respect to $x$
}

\hidenextnotation
\newnotationcommandoptarg{3}{\partialderiv}{%
  \frac{\partialdiff\ifnotempty{#1}{^{#1}}}{#2} #3%
}{zzzzzzpartialderivative}{$\partialderiv{\partialdiff{}x}{\objfun}$}{
  Partial derivative of $\objfun$ with respect to $x$
}

\hidenextnotation
\newnotationcommandoptarg{3}{\tpartialderiv}{%
  \tfrac{\partialdiff\ifnotempty{#1}{^{#1}}}{#2} #3%
}{zzzzzzpartialderivativesmall}{$\tpartialderiv{\partialdiff{}x}{\objfun}$}{
  Partial derivative of $\objfun$ with respect to $x$
}

\newnotationcommand[2]{\gradient}{%
  \mathop{\*\nabla_{\hspace*{-0.1em}#1}} #2%
}{zzzzzzgradient}{$\gradient{\*x}{\objfun}$}{
  Gradient of a function $\objfun$ with respect to $\*x$
  (transposed Jacobian)
}

\newnotationcommand[2]{\hessian}{%
  \mathop{\*\nabla_{\hspace*{-0.1em}#1}^2} #2%
}{zzzzzzhessian}{$\hessian{\*x}{\objfun}$}{
  Hessian of a function $\objfun$ with respect to $\*x$
}

\newnotationcommand[1]{\tr}{#1^\mathrm{T}}{zzzzT}{$\tr{(\cdot)}$}{
  Transpose of a vector or matrix
}

\newnotationcommand{\eye}{\mat{I}}{I}{$\eye$}{
  Identity matrix
}

\newnotationcommand[1]{\nonnegpart}{(#1)_+}{zzzzzz.+}{$\nonnegpart{\cdot}$}{
  Non-negative part $\max({\cdot}, 0)$
}

\newcommand*{\expectationsign}[1][]{\mathbb{E}\ifnotempty{#1}{_{#1}}}

\newnotationcommandoptarg{2}{\expectation}{%
  \mathbb{E}\ifnotempty{#1}{_{#1}}\bracket*{#2}%
}{EX}{$\expectation{X}$}{
  Expectation of the random variable $X$
}

\let\oplusorig\oplus
\let\oplus\undefined
\newnotationcommand{\oplus}{\oplusorig}{zzzzzzoplus}{$\oplus$}{%
  Internal direct sum of vector spaces
}

\hidenextnotation
\newnotationcommand{\convolution}{\ast}{zzzzzzconvolution}{$\convolution$}{
  Convolution $f \ast g$ of two functions $f$ and $g$
}

\newnotationcommand[2]{\kronecker}{\delta_{#1,#2}}{zzdAB}{$\kronecker{A}{B}$}{
  Kronecker delta $\ceq 1$ if $A = B$ and $\ceq 0$ otherwise
  ($A, B$ arbitrary objects)
}

\DeclareMathOperator{\xor}{xor}
\let\xororig\xor
\hidenextnotation
\makecommandnotation{\xor}{xor}{$\xororig$}{
  Bitwise ``exclusive or''
}

\newnotationcommand[2]{\restrictfcn}{#1|_{#2}}{fD}{$\restrictfcn{\objfun}{D}$}{
  Restriction $\restrictfcn{\objfun}{D}\colon D \to \real$,
  $\restrictfcn{\objfun}{D}(x) \ceq \objfun(x)$,
  onto a sub-domain $D$
}

\newnotationcommand{\dotcup}{\mathbin{\dot{\cup}}}{zzzzzzunion}{$\dotcup$}{
  Disjoint union of sets
}

\DeclareMathOperator*{\bigdotcup}{\dot{\bigcup}}

\hidenextnotation
\newnotationcommand{\eq}{\sim}{zzzzzzeq}{$\eq$}{
  Equivalence relation where
  $a \eq b$ denotes that two elements $a, b$ are equivalent
}

\hidenextnotation
\newnotationcommand[1]{\samepole}{\eq_{#1}}{zzzzzzsimt}{$\samepole{t}$}{
  Pole equivalence relation with respect to dimension $t = 1, \dotsc, d$
  ($\*k \samepole{t} \*k'$ if $\*k$ and $\*k'$ are on the same pole
  $\lisetpole \in \eqclasses{\liset}{\samepole{t}}$ in dimension $t$)
}

\newnotationcommand[2]{\eqclasses}{{#1}/{#2}}{zzzzzzeqclasses}{%
  $\eqclasses{\cdot}{\eq}$%
}{
  Set of equivalence classes for an equivalence relation $\eq$ on
  a set $\cdot$
}

\newnotationcommand[2]{\eqclass}{[#1]_{#2}}{zzzzzzeqclass}{%
  $\eqclass{\cdot}{\eq}$%
}{
  Equivalence class of $\cdot$ (set of all elements equivalent to $\cdot$)
  with respect to $\eq$
}

\newnotationcommand[1]{\hiset}{I_{#1}}{Il}{$\hiset{\*l}$}{
  Set of odd hierarchical indices of level $\*l$
}

\newnotationcommand{\liset}{K}{K}{$\liset$}{
  Finite set of hierarchical level-index pairs $(\*l, \*i)$
  or continuous indices $\*k \in \natz^d$
}

\newnotationcommandoptarg{1}{\lisetpole}{%
  \liset_\mathrm{pole}\ifnotempty{#1}{^{#1}}%
}{Kpole}{$\lisetpole$}{
  Pole $\lisetpole \in \eqclasses{\liset}{\samepole{t}}$
  in some dimension $t$ ($\lisetpole$ is a subset of $\liset$)
}

\newnotationcommand{\levelset}{L}{L1}{$\levelset$}{
  Finite subset $\levelset \subset \natz^d$ of levels
}

\hidenextnotation
\newnotationcommandoptarg{4}{\coarselevelset}{%
  \levelset_{#2,#3}^{\sparse(#4)\appendwithcomma{#1}}%
}{Lndbs}{$\coarselevelset{n}{d}{b}$}{
  Finite subset $\coarselevelset{n}{d}{b} \subset \natz^d$ of levels
  for the grid $\coarseregsgset{n}{d}{b}$
}

\hidenextnotation
\newnotationcommandoptarg{2}{\relindexset}{%
  J_{#2}\ifnotempty{#1}{^{#1}}%
}{Jl}{$\relindexset{l}$}{
  Set of relevant indices for the translation-invariant
  fundamental transformation
}

\hidenextnotation
\newnotationcommandoptarg{2}{\fcnval}{%
  f\ifnotempty{#1}{^{#1}}(\gp{#2})%
}{fli}{$\fcnval{\*l,\*i}$}{
  Value of the objective function $\objfun$ at the grid point $\gp{\*l,\*i}$
}

\newnotationcommandoptarg{2}{\surplus}{%
  \alpha_{#2}\ifnotempty{#1}{^{#1}}%
}{zzali}{$\surplus{\*l,\*i}$}{
  Hierarchical surpluses
}

\newcommand*{\vsurplus}{\*\alpha}
\newcommand*{\surplustilde}[2][]{\tilde{\alpha}_{#2}\ifnotempty{#1}{^{#1}}}

\newnotationcommandoptarg{2}{\interpcoeff}{%
  c_{#2}\ifnotempty{#1}{^{#1}}%
}{cli}{$\interpcoeff{\*l,\*i}$}{
  Full grid interpolation coefficients
}

\newcommand*{\interpcoefftilde}[2][]{\tilde{c}_{#2}\ifnotempty{#1}{^{#1}}}

\newcommand*{\fundsplcoeff}[3][]{c_{#2,#3}\ifnotempty{#1}{^{#1}}}
\newcommand*{\fundsplcutoff}[1]{n_{#1}}

\newcommand*{\wfundsplcoeff}[3][]{c_{#2,#3}\ifnotempty{#1}{^{#1}}}

\newnotationcommand{\econst}{\mathrm{e}}{e}{$\econst$}{
  Euler constant $\exp(1)$
}

\hidenextnotation
\newnotationcommand{\ngpMax}{\ngp_\mathrm{max}}{Nmax}{$\ngpMax$}{
  Maximum number of grid points when generating spatially adaptive sparse grids
}

\hidenextnotation
\newnotationcommand{\refinetol}{\kappa}{zzk}{$\refinetol$}{
  Refinement tolerance
}

\hidenextnotation
\newnotationcommand{\error}{\varepsilon}{zze}{$\error$}{
  Error measure
}

\newnotationcommand{\xopt}{\*x^\opt}{xopt}{$\xopt$}{
  Solution of an optimization problem of the form
  $\xopt = \vecargmin \objfun(\*x)$
}

\newnotationcommand{\xoptappr}{\*x^{\opt,\ast}}{xopt*}{$\xoptappr$}{
  Approximation for $\xopt = \vecargmin \objfun(\*x)$
}

\hidenextnotation
\newnotationcommandoptarg{1}{\xoptscaled}{%
  \bar{\*x}^\opt%
}{xoptbar}{$\xoptscaled$}{
  Scaled solution of an optimization problem
}

\hidenextnotation
\newnotationcommandoptarg{1}{\xscaled}{%
  \bar{\*x}%
}{xbar}{$\xscaled$}{
  Scaled $\xscaled \in \clint{\*a, \*b}$
  for the definition of test functions
}

\hidenextnotation
\newnotationcommand[1]{\xoptscaledentry}{%
  \bar{x}^\opt_{#1}%
}{xoptbart}{$\xoptscaledentry{t}$}{
  Scaled solution entry of an optimization problem
}

\hidenextnotation
\newnotationcommand[1]{\xscaledentry}{%
  \bar{x}_{#1}%
}{xbart}{$\xscaledentry{t}$}{
  Scaled $\xscaledentry{t} \in \clint{a, b}$
  for the definition of test functions
}

\newnotationcommandoptarg{2}{\basis}{%
  \varphi_{#2}\ifnotempty{#1}{^{#1}}%
}{zzvli}{$\basis{\*l,\*i}$}{
  Hierarchical basis function of level $\*l$, index $\*i$
}

\hidenextnotation
\newnotationcommand[1]{\fundbasis}{%
  \varphi_{#1}^\mathrm{f}%
}{zzvlif}{$\fundbasis{\*l,\*i}$}{
  Fundamental hierarchical basis function of level $\*l$, index $\*i$
  (satisfying the fundamental property)
}

\hidenextnotation
\newnotationcommand[1]{\wfundbasis}{%
  \varphi_{#1}^\mathrm{wf}%
}{zzvliwf}{$\wfundbasis{\*l,\*i}$}{
  Weakly fundamental hierarchical basis function of level $\*l$, index $\*i$
  (satisfying the weakly fundamental property)
}

\newnotationcommandoptarg{1}{\parentfcn}{%
  \varphi\ifnotempty{#1}{^{#1}}%
}{zzv}{$\parentfcn$}{
  Parent function $\parentfcn\colon \clint{0, 1} \to \real$
  ($\basis{l,i}$ is scaled translate of $\parentfcn$)
}

\newnotationcommandoptarg{3}{\bspl}{%
  \varphi_{#2}^{#3\appendwithcomma{#1}}%
}{zzvlip}{$\bspl{\*l,\*i}{\*p}$}{
  Hierarchical standard B-spline basis function of
  level $\*l$, index $\*i$, degree $\*p$
}

\newnotationcommand[1]{\cardbspl}{b^{#1}}{bp}{$\cardbspl{p}$}{
  Cardinal B-spline of degree $p$
}

\hidenextnotation
\newnotationcommandoptarg{2}{\parentbspl}{%
  \parentfcn[#2\appendwithcomma{#1}]%
}{zzvp}{$\parentbspl{p}$}{
  Centralized cardinal B-spline
  $\parentbspl{p} \ceq \cardbspl{p}({\cdot} + \tfrac{p+1}{2})$
  of degree $p$
}

\hidenextnotation
\newnotationcommandoptarg{2}{\parentfundspl}{%
  \parentfcn[#2,\fs\appendwithcomma{#1}]%
}{zzvpfs}{$\parentfundspl{p}$}{
  Fundamental spline of degree $p$ (parent function)
}

\hidenextnotation
\newnotationcommandoptarg{2}{\parentwfundspl}{%
  \parentfcn[#2,\wfs\appendwithcomma{#1}]%
}{zzvpwfs}{$\parentwfundspl{p}$}{
  Weakly fundamental spline of degree $p$ (parent function)
}

\hidenextnotation
\newnotationcommand[2]{\nonunifbspl}{b_{#1}^{#2}}{bkop}{%
  $\nonunifbspl{k,\knotseq}{p}$%
}{
  Non-uniform B-spline of degree $p$ with knots $\knotseq$ and index $k$
}

\hidenextnotation
\newnotationcommand{\spl}{s}{s}{$s$}{
  Spline (piecewise polynomial)
}

\hidenextnotation
\newnotationcommand[1]{\charfun}{\chi_{#1}}{zzxA}{$\charfun{A}$}{
  Characteristic function of the set $A \subset \real$
}

\newnotationcommand{\objfun}{f}{f}{$\objfun$}{
  Objective function $\objfun\colon \clint{\*0, \*1} \to \real$
}

\newnotationcommandoptarg{1}{\ineqconfun}{%
  \ifempty{#1}{\*g}\ifnotempty{#1}{g_{#1}}%
}{g}{$\ineqconfun$}{
  Inequality constraint function
  $\ineqconfun\colon \clint{\*0, \*1} \to \real^{m_{\ineqconfun}}$
  (constraint $\ineqconfun(\*x) \le \*0$)
}

\hidenextnotation
\newnotationcommandoptarg{1}{\eqconfun}{
  \ifempty{#1}{\*h}\ifnotempty{#1}{h_{#1}}%
}{h}{$\eqconfun$}{
  Equality constraint function
  $\eqconfun\colon \clint{\*0, \*1} \to \real^{m_{\eqconfun}}$
  (optimal solutions $\*x$ must satisfy $\eqconfun(\*x) = \*0$)
}

\hidenextnotation
\newnotationcommand{\objfunscaled}{\bar{\objfun}}{fbar}{$\objfunscaled$}{
  Scaled objective function $\objfunscaled\colon \clint{\*a, \*b} \to \real$
}

\hidenextnotation
\newnotationcommandoptarg{1}{\ineqconfunscaled}{%
  \ifempty{#1}{\bar{\*g}}\ifnotempty{#1}{\bar{g}_{#1}}%
}{gbar}{$\ineqconfunscaled$}{
  Scaled inequality constraint function
  $\ineqconfunscaled\colon \clint{\*a, \*b} \to \real^{m_{\ineqconfun}}$
}

\hidenextnotation
\newnotationcommand[1]{\testobjfunscaled}{%
  \bar{\objfun}_{\mathrm{#1}}%
}{fbarTest}{$\testobjfunscaled{Test}$}{
  Scaled test objective function
  $\testobjfunscaled{Test}\colon \clint{\*a, \*b} \to \real$
  with name ``Test''
}

\hidenextnotation
\newnotationcommandoptarg{2}{\testineqconfunscaled}{%
  \ifempty{#1}{\bar{\*g}_{\mathrm{#2}}}%
  \ifnotempty{#1}{\bar{g}_{\mathrm{#2},#1}}%
}{gbarTest}{$\testineqconfunscaled{Test}$}{
  Scaled test inequality constraint function
  $\testineqconfunscaled{Test}\colon
  \clint{\*a, \*b} \to \real^{m_{\ineqconfun}}$
  with name ``Test''
}

\newnotationcommand[1]{\fgintp}{f_{#1}}{fl}{$\fgintp{\*l}$}{
  Full grid interpolant of $\objfun$ in $\ns{\*l}$
}

\newnotationcommandoptarg{1}{\sgintp}{%
  f^{\sparse\appendwithcomma{#1}}%
}{fs}{$\sgintp$}{
  Sparse grid interpolant of $\objfun$ in $\sgspace$
  (on some grid $\sgset$)
}

\newcommand*{\vsgintp}[1][]{\*f^{\sparse\appendwithcomma{#1}}}

\newnotationcommandoptarg{3}{\regsgintp}{%
  f_{#2,#3}^{\sparse\appendwithcomma{#1}}%
}{fsnd}{$\regsgintp{n}{d}$}{
  Regular sparse grid interpolant of $\objfun$ in $\regsgspace{n}{d}$
}

\newnotationcommandoptarg{2}{\lagrangepoly}{%
  L_{#2}\ifnotempty{#1}{^{#1}}%
}{Lli}{$\lagrangepoly{l,i}$}{
  Lagrange polynomial of level $l$, index $i$
}

\hidenextnotation
\newnotationcommand[1]{\knot}{\xi_{#1}}{zzo10}{$\knot{k}$}{
  Knot $\in \real$ of a spline space
  (point where the splines are not smooth)
}

\hidenextnotation
\newnotationcommand{\knotseq}{\*\xi}{zzo20}{$\knotseq$}{
  Knot sequence $\ceq (\knot{0}, \dotsc, \knot{m+p})$
}

\hidenextnotation
\newnotationcommandoptarg{3}{\nodalknot}{%
  \xi_{#2}^{#3\appendwithcomma{#1}}%
}{zzolkp}{$\nodalknot{l,k}{p}$}{
  Knot for uniform nodal B-splines of level $l$, degree $p$
  ($k$-th entry of $\nodalknotseq{l}{p}$)
}

\hidenextnotation
\newnotationcommandoptarg{3}{\nodalknotseq}{%
  \*\xi_{#2}^{#3\appendwithcomma{#1}}%
}{zzolp}{$\nodalknotseq{l}{p}$}{
  Knot sequence for uniform nodal B-splines of level $l$, degree $p$
}

\hidenextnotation
\newnotationcommand[2]{\spldomain}{D_{#1}^{#2}}{D10}{%
  $\spldomain{\knotseq}{p}$%
}{
  Spline interpolation domain of $\nonunifsplspace{\knotseq}{p}$
}

\hidenextnotation
\newnotationcommand[2]{\rspldomain}{D_{#1}^{#2}}{D20}{%
  $\rspldomain{l}{p}$%
}{
  Spline interpolation domain of $\restrictedsplspace{l}{p}$
  (uniform nodal B-splines of level $l$, degree $p$)
}

\hidenextnotation
\newnotationcommandoptarg{2}{\fuzzy}{%
  \tilde{#2}\ifnotempty{#1}{^{#1}}%
}{xtilde}{$\fuzzy{x}$}{
  Fuzzy set
}

\hidenextnotation
\newnotationcommandoptarg{2}{\memfun}{%
  \mu_{\fuzzy[#1]{#2}}%
}{zzmtilde}{$\memfun{x}$}{
  Membership function $\memfun{x}\colon X \to \clint{0, 1}$ of the
  fuzzy set $\fuzzy{x}$
}

\hidenextnotation
\newnotationcommandoptarg{3}{\acut}{%
  (\tilde{#2}\ifnotempty{#1}{_{#1}})_{#3}%
}{xtildea}{$\acut{x}{\alpha}$}{
  $\alpha$-cut of $\fuzzy{x}$ for $\alpha \in \clint{0, 1}$
}

\hidenextnotation
\newnotationcommand{\dimobjdomain}{\tilde{d}}{dtilde}{$\dimobjdomain$}{
  \#dimensions of $\objdomain$
}

\hidenextnotation
\newnotationcommand{\objdomain}{\tilde{\Omega}}{zzZtilde}{$\objdomain$}{
  Object domain
}

\hidenextnotation
\newnotationcommand{\densglobal}{\tilde{\varrho}}{zzr2}{$\densglobal$}{
  Global density fcn.
}

\hidenextnotation
\newnotationcommandoptarg{1}{\denscell}{%
  \varrho\ifnotempty{#1}{^{#1}}%
}{zzr}{$\denscell$}{
  Micro-cell density fcn.
}

\hidenextnotation
\newnotationcommandoptarg{1}{\compliance}{%
  J\ifnotempty{#1}{^{#1}}%
}{Jr}{$\compliance(\densglobal)$}{
  Compliance fcn.
}

\hidenextnotation
\newnotationcommand[1]{\vol}{\mathrm{vol}(#1)}{volZ}{$\vol{\objdomain}$}{
  Total volume
}

\hidenextnotation
\newnotationcommand[2]{\voldens}{\mathrm{vol}_{#1}(#2)}{volrZ}{
  $\voldens{\densglobal}{\objdomain}$
}{%
  Volume w.r.t.\ $\densglobal$
}

\hidenextnotation
\newnotationcommand{\densub}{\varrho^\ast}{r*}{$\densub$}{
  Density bound
}

\hidenextnotation
\newnotationcommand{\densmean}{\bar{\varrho}}{rbar}{$\densmean$}{
  Mean density
}

\hidenextnotation
\newnotationcommand{\force}{\*F}{F}{$\force$}{
  External force
}

\hidenextnotation
\newnotationcommand{\displacement}{\*u}{u2}{$\displacement$}{
  Displacement fcn.
}

\hidenextnotation
\newnotationcommandoptarg{1}{\etensor}{%
  \mat{E}\ifnotempty{#1}{^{#1}}%
}{E}{$\etensor$}{
  Elasticity tensor
}

\hidenextnotation
\newnotationcommandoptarg{2}{\etensorentry}{%
  E_{#2}\ifnotempty{#1}{^{#1}}%
}{Ej}{$\etensorentry{j}$}{
  $j$-th entry of the elasticity tensor
}

\hidenextnotation
\newnotationcommandoptarg{1}{\cholfactor}{%
  \mat{R}\ifnotempty{#1}{^{#1}}%
}{R2}{$\cholfactor$}{
  Cholesky factor
}

\hidenextnotation
\newnotationcommandoptarg{1}{\etensorchol}{%
  \mat{E}^{\chol\appendwithcomma{#1}}%
}{Echol}{$\etensor$}{
  Elasticity tensor retrieved by Cholesky factorization
}

\newcommand*{\complianceintp}[1][]{%
  \compliance[\sparse\appendwithcomma{#1}]%
}
\newcommand*{\denscellintp}[1][]{%
  \denscell[\sparse\appendwithcomma{#1}]%
}
\newcommand*{\etensorintp}[1][]{%
  \etensor[\sparse\appendwithcomma{#1}]%
}
\newcommand*{\etensorentryintp}[2][]{%
  \etensorentry[\sparse\appendwithcomma{#1}]{#2}%
}
\newcommand*{\cholfactorintp}[1][]{%
  \cholfactor[\sparse\appendwithcomma{#1}]%
}
\newcommand*{\etensorcholintp}[1][]{%
  \etensorchol[\sparse\appendwithcomma{#1}]%
}

\hidenextnotation
\newnotationcommandoptarg{2}{\mcp}{%
  \*x^{\ifnotempty{#1}{#1,}(#2)}%
}{xq}{$\mcp{q}$}{
  Micro-cell param.
}

\newcommand*{\mcpopt}[1]{\mcp[\opt]{#1}}
\newcommand*{\mcpoptappr}[1]{\mcp[\opt,\ast]{#1}}

\hidenextnotation
\newnotationcommandoptarg{3}{\mcpentry}{%
  x^{\ifnotempty{#1}{#1,}(#2)}_{#3}%
}{xqt}{$\mcpentry{q}{t}$}{
  Micro-cell param. entry
}

\hidenextnotation
\newnotationcommand{\elbang}{\theta}{zzh}{$\elbang$}{
  Elbow angle
}

\hidenextnotation
\newnotationcommand[1]{\equielbang}{%
  \elbang_{#1}%
}{zzhFL}{$\equielbang{\forceL}$}{
  Equil. elbow angle
}

\hidenextnotation
\newnotationcommand{\tarelbang}{\elbang^\ast}{zzh*}{$\tarelbang$}{
  Target elbow angle
}

\hidenextnotation
\newnotationcommand{\forceX}{F_X}{FX}{$\forceX$}{
  Force
}

\hidenextnotation
\newnotationcommand{\forceT}{F_\mathrm{T}}{FT}{$\forceT$}{
  Triceps force
}

\hidenextnotation
\newnotationcommand{\forceB}{F_\mathrm{B}}{FB}{$\forceB$}{
  Biceps force
}

\hidenextnotation
\newnotationcommand{\forceL}{F_\mathrm{L}}{FL}{$\forceL$}{
  Load force
}

\hidenextnotation
\newnotationcommand{\arm}{r}{r}{$\arm$}{
  Lever arm
}

\hidenextnotation
\newnotationcommand{\armX}{\arm_X}{rX}{$\armX$}{
  Lever arm
}

\hidenextnotation
\newnotationcommand{\armT}{\arm_\mathrm{T}}{rT}{$\armT$}{
  Triceps lever arm
}

\hidenextnotation
\newnotationcommand{\armB}{\arm_\mathrm{B}}{rB}{$\armB$}{
  Biceps lever arm
}

\hidenextnotation
\newnotationcommand{\armL}{\arm_\mathrm{L}}{rL}{$\armL$}{
  Load lever arm
}

\hidenextnotation
\newnotationcommand{\act}{\beta}{zzb}{$\act$}{
  Muscle activation
}

\hidenextnotation
\newnotationcommand{\actX}{\act_X}{zzbX}{$\actX$}{
  Muscle activation
}

\hidenextnotation
\newnotationcommand{\actT}{\act_\mathrm{T}}{zzbT}{$\actT$}{
  Triceps activation
}

\hidenextnotation
\newnotationcommand{\actB}{\act_\mathrm{B}}{zzbB}{$\actB$}{
  Biceps activation
}

\hidenextnotation
\newnotationcommand{\moment}{M}{M}{$\moment$}{
  Moment
}

\hidenextnotation
\newnotationcommand[1]{\actdomain}{%
  D_{\equielbang{#1}}%
}{$DhL$}{$\actdomain{\forceL}$}{
  Domain of $\equielbang{\forceL}$
}

\hidenextnotation
\newnotationcommandoptarg{1}{\forceXintp}{%
  F_X^{\sparse\ifnotempty{#1}{,#1}}%
}{FXs}{$\forceXintp$}{
  Muscle force interpolated with sparse grids
}

\hidenextnotation
\newnotationcommandoptarg{1}{\forceTintp}{%
  \forceT^{\sparse\ifnotempty{#1}{,#1}}%
}{FTs}{$\forceTintp$}{
  Triceps force interpolated with sparse grids
}

\hidenextnotation
\newnotationcommandoptarg{1}{\forceBintp}{%
  \forceB^{\sparse\ifnotempty{#1}{,#1}}%
}{FBs}{$\forceTBntp$}{
  Biceps force interpolated with sparse grids
}

\hidenextnotation
\newnotationcommandoptarg{1}{\momentintp}{%
  \moment^{\sparse\ifnotempty{#1}{,#1}}%
}{Ms}{$\momentintp$}{
  Moment obtained through sparse grid interpolation
}

\hidenextnotation
\newnotationcommandoptarg{2}{\equielbangintp}{%
  \elbang_{#2}^{\sparse\ifnotempty{#1}{,#1}}%
}{zzhFLs}{$\equielbangintp{\forceL}$}{
  Equil. elbow angle obtained through sparse grid interpolation
}

\hidenextnotation
\newnotationcommandoptarg{2}{\actdomainintp}{%
  D_{\equielbangintp{#2}}\ifnotempty{#1}{^{#1}}%
}{$DhLs$}{$\actdomainintp{\forceL}$}{
  Domain of $\equielbangintp{\forceL}$
}

\hidenextnotation
\newnotationcommand{\forceXref}{\forceX^{\reference}}{FXref}{$\forceXref$}{
  Reference force
}

\hidenextnotation
\newnotationcommand{\forceTref}{\forceT^{\reference}}{FTref}{$\forceTref$}{
  Reference triceps force
}

\hidenextnotation
\newnotationcommand{\forceBref}{\forceB^{\reference}}{FBref}{$\forceTref$}{
  Reference biceps force
}

\hidenextnotation
\newnotationcommand{\momentref}{M^{\reference}}{Mref}{$\momentref$}{
  Reference moment
}

\hidenextnotation
\newnotationcommand[1]{\equielbangref}{%
  \elbang_{#1}^{\mathrm{ref}}%
}{zzhFLref}{$\equielbangref{\forceL}$}{
  Reference equil. elbow angle
}

\newnotation{t}{$t$}{
  Time
}

\hidenextnotation
\newnotationcommand{\state}{\*x}{xt}{$\state_t$}{
  State variables
}

\hidenextnotation
\newnotationcommand{\stateentry}{x}{xtj}{$\stateentry_{t,j}$}{
  $j$-th state variable ($j = 1, \dotsc, d$)
}

\hidenextnotation
\newnotationcommand{\policy}{\*y}{yt}{$\policy_t$}{
  Policy variables
}

\hidenextnotation
\newnotationcommand{\policyentry}{y}{ytj}{$\policyentry_{t,j}$}{
  $j$-th policy variable ($j = 1, \dotsc, m_{\policy}$)
}

\hidenextnotation
\newnotationcommand{\stochastic}{\*\omega}{zzzt}{$\stochastic_t$}{
  Stochastic variables
}

\hidenextnotation
\newnotationcommand{\stochdomain}{\Omega}{zzZ2}{$\stochdomain$}{
  Domain of the stochastic variables $\stochastic_t$
}

\hidenextnotation
\newnotationcommand{\discrstate}{\*\theta}{zzht}{$\discrstate_t$}{
  Discrete state variables
}

\hidenextnotation
\newnotationcommand{\discrstdomain}{\Theta}{zzH2}{$\discrstdomain$}{
  Domain of the discrete state variables $\discrstate_t$
}

\hidenextnotation
\newnotationcommand{\wealth}{w}{wt}{$\wealth_t$}{
  Wealth
}

\hidenextnotation
\newnotationcommand{\consume}{c}{ct}{$\consume_t$}{
  Consumption
}

\hidenextnotation
\newnotationcommand{\bond}{b}{bt}{$\bond_t$}{
  Bond investment
}

\hidenextnotation
\newnotationcommand{\stock}{s}{st}{$\stock_{t,j}$}{
  Stock holding
}

\hidenextnotation
\newnotationcommandoptarg{1}{\buy}{%
  \delta^{+\appendwithcomma{#1}}%
}{zzdtj+}{$\buy_{t,j}$}{
  Stock buy
}

\hidenextnotation
\newnotationcommandoptarg{1}{\sell}{%
  \delta^{-\appendwithcomma{#1}}%
}{zzdtj-}{$\sell_{t,j}$}{
  Stock sell
}

\hidenextnotation
\newnotationcommandoptarg{1}{\buysell}{%
  \delta^{\pm\appendwithcomma{#1}}%
}{zzdtj+-}{$\buysell_{t,j}$}{
  Stock buy/sell
}

\hidenextnotation
\newnotationcommand{\normstate}{\hat{\*x}}{xt2}{$\normstate_t$}{
  Normalized state variables
}

\hidenextnotation
\newnotationcommand{\normpolicy}{\hat{\*y}}{yt2}{$\normpolicy_t$}{
  Normalized policy variables
}

\hidenextnotation
\newnotationcommand{\normconsume}{\hat{c}}{c^t}{$\consume_t$}{
  Normalized consumption
}

\hidenextnotation
\newnotationcommand{\normbond}{\hat{b}}{b^t}{$\bond_t$}{
  Normalized bond investment
}

\hidenextnotation
\newnotationcommandoptarg{1}{\normbuy}{%
  \hat{\delta}^{+\appendwithcomma{#1}}%
}{zzd^tj+}{$\buy_{t,j}$}{
  Normalized stock buy
}

\hidenextnotation
\newnotationcommandoptarg{1}{\normsell}{%
  \hat{\delta}^{-\appendwithcomma{#1}}%
}{zzd^tj-}{$\sell_{t,j}$}{
  Normalized stock sell
}

\newcommand*{\vstock}{\*s}
\newcommand*{\vbuy}[1][]{\*\delta^{+\appendwithcomma{#1}}}
\newcommand*{\vsell}[1][]{\*\delta^{-\appendwithcomma{#1}}}
\newcommand*{\vbuysell}[1][]{\*\delta^{\pm\appendwithcomma{#1}}}
\newcommand*{\vnormbuy}[1][]{\hat{\*\delta}^{+\appendwithcomma{#1}}}
\newcommand*{\vnormsell}[1][]{\hat{\*\delta}^{-\appendwithcomma{#1}}}
\newcommand*{\vnormbuysell}[1][]{\hat{\*\delta}^{\pm\appendwithcomma{#1}}}

\hidenextnotation
\newnotationcommand{\wealthratio}{\eta}{zzft}{$\wealthratio_t$}{
  Wealth ratio
}

\hidenextnotation
\newnotationcommand{\tac}{\tau}{zzt}{$\tac$}{
  Transaction cost rate
}

\hidenextnotation
\newnotationcommand{\bondreturn}{r}{rt}{$\bondreturn_t$}{
  Bond return rate
}

\hidenextnotation
\newnotationcommand{\stockreturn}{\lambda}{zzltj}{$\stockreturn_{t,j}$}{
  Stock return rate
}

\newcommand*{\vstockreturn}{\*\lambda}

\hidenextnotation
\newnotationcommand{\valuefcn}{J}{Jt}{$\valuefcn_t$}{
  Value fcn.
}

\hidenextnotation
\newnotationcommand{\optpolicyfcn}{\policy^\opt}{poptt}{$\optpolicyfcn_t$}{
  Optimal policy fcn.
}

\hidenextnotation
\newnotationcommand{\optnormpolicyfcn}{%
  \hat{\policy}^\opt%
}{p^optt}{$\optnormpolicyfcn_t$}{
  Normalized optimal policy fcn.
}

\hidenextnotation
\newnotationcommandoptarg{1}{\valueintp}{%
  J^{\sparse\ifnotempty{#1}{,#1}}%
}{Jts}{$\valueintp_t$}{
  Interpolated value fcn.
}

\hidenextnotation
\newnotationcommandoptarg{1}{\optpolicyintp}{%
  \policy^{\opt,\sparse\ifnotempty{#1}{,#1}}%
}{poptts}{$\optpolicyintp_t$}{
  Interpolated optimal policy fcn.
}

\hidenextnotation
\newnotationcommandoptarg{1}{\optnormpolicyintp}{%
  \hat{\policy}^{\opt,\sparse\ifnotempty{#1}{,#1}}%
}{p^optts}{$\optnormpolicyintp_t$}{
  Interpolated normalized optimal policy fcn.
}

\hidenextnotation
\newnotationcommandoptarg{1}{\cetvalueintp}{%
  \tilde{J}^{\sparse\ifnotempty{#1}{,#1}}%
}{Jts2}{$\cetvalueintp_t$}{
  Interpolated certainty-equivalent-transf.\ value fcn.
}

\hidenextnotation
\newnotationcommandoptarg{1}{\normcetvaluefcn}{%
  \hat{\tilde{J}}^{\ifnotempty{#1}{,#1}}%
}{Jt3}{$\normcetvaluefcn_t$}{
  Normalized certainty-equivalent-transf.\ value fcn.
}

\hidenextnotation
\newnotationcommandoptarg{1}{\normcetvalueintp}{%
  \hat{\tilde{J}}^{\sparse\ifnotempty{#1}{,#1}}%
}{Jts3}{$\normcetvalueintp_t$}{
  Interpolated normalized certainty-equivalent-transf.\ value fcn.
}

\hidenextnotation
\newnotationcommand{\optpolicymean}{\bar{\policy}^\opt}{popttbar}{%
  $\optpolicymean_t$%
}{
  Optimal policy mean
}

\hidenextnotation
\newnotationcommand{\statefcn}{\*\psi}{zzyt}{$\statefcn_t$}{
  State transition fcn.
}

\hidenextnotation
\newnotationcommandoptarg{1}{\normstatefcn}{%
  \hat{\*\psi}\ifnotempty{#1}{^{\raisebox{-0.5em}{\scriptsize$#1$}}}%
}{zzyt2}{$\normstatefcn_t$}{
  Normalized state transition fcn.
}

\hidenextnotation
\newnotationcommand{\utilityfcn}{u}{u3}{$\utilityfcn$}{
  Utility fcn.
}

\hidenextnotation
\newnotationcommand{\riskav}{\gamma}{zzc}{$\riskav$}{
  Risk aversion
}

\hidenextnotation
\newnotationcommand{\patience}{\varrho}{zzr3}{$\patience$}{
  Patience factor
}

\hidenextnotation
\newnotationcommand{\norefine}{q}{qt}{$\norefine_t$}{
  Number of grid refines
}

\hidenextnotation
\newnotationcommand{\quadweight}{\zeta}{zzezztj}{$\quadweight_t^{(j)}$}{
  $j$-th quadrature weight
}

\hidenextnotation
\newnotationcommand{\sumfcn}{\mathop{\Sigma}}{zzS}{$\sumfcn$}{
  Sum function $\sumfcn(\*x) \ceq \tr{\*1} \*x$
}

\hidenextnotation
\newnotationcommand{\compmult}{\odot}{zzzzzz.}{$\*a \compmult \*b$}{
  Component-wise product of vectors $\*a$ and $\*b$
}

\hidenextnotation
\newnotationcommand{\normcropfactor}{\hat{\beta}}{zzb2}{$\normcropfactor$}{
  State cropping factor
}

\hidenextnotation
\newnotationcommand{\lagrangian}{\mathcal{L}}{Lt}{$\lagrangian$}{
  Lagrangian
}

\hidenextnotation
\newnotationcommand{\multiplier}{\mu}{zzm}{$\multiplier$}{
  Lagrangian multiplier
}

\hidenextnotation
\newnotationcommand{\eulererror}{%
  \error^{\mathrm{Eu}}%
}{zzetwEu}{$\eulererror_t$}{
  Euler equation error
}

\hidenextnotation
\newnotationcommand{\weightedeulererror}{%
  \error^{\mathrm{w},\mathrm{Eu}}%
}{zzetEu}{$\weightedeulererror_t$}{
  Weighted Euler equation error
}

\hidenextnotation
\newnotationcommand{\weightedeulererrorLtwo}{%
  \error^{\mathrm{w},\mathrm{Eu},\Ltwo}%
}{zzetEuL2}{$\weightedeulererrorLtwo_t$}{
  Normalized $\Ltwo$ norm of weighted Euler equation error
}

\newcommand*{\term}[1]{\emph{#1}}

\newcommand*{\containsvector}[3]{%
  \protect\IfSubStr{\detokenize{#1}}{\detokenize{\*}}{#2}{#3}%
}

\renewcommand*{\subset}{\subseteq}
\renewcommand*{\supset}{\supseteq}

\newcommand*{\ceq}{\coloneqq}

\newcommand*{\diff}{\mathop{}\!\mathrm{d}}
\newcommand*{\dx}{\diff{}x}
\newcommand*{\partialdiff}{\mathop{}\!\partial}

\renewcommand*{\vec}[1]{{\boldsymbol{#1}}}
\def\*#1{\vec{#1}}
\newcommand*{\veclog}{\mathop{\vec{\log}}}

\newcommand*{\mat}[1]{{\boldsymbol{#1}}}

\newcommand*{\fa}[2]{\forall_{#1}\;#2}
\newcommand*{\ex}[2]{\exists_{#1}\;#2}
\newcommand*{\fafa}[3]{\forall_{#1} \forall_{#2}\;#3}
\newcommand*{\faex}[3]{\forall_{#1} \exists_{#2}\;#3}
\newcommand*{\exfa}[3]{\exists_{#1} \forall_{#2}\;#3}

\newcommand*{\falarge}[2]{\forall_{#1}\;\;#2}
\newcommand*{\exlarge}[2]{\exists_{#1}\;\;#2}
\newcommand*{\fafalarge}[3]{\forall_{#1} \forall_{#2}\;\;#3}

\newcommand*{\exfalarge}[3]{\exists_{#1} \forall_{#2}\;\;#3}

\DeclarePairedDelimiter{\braced}{\{}{\}}
\DeclarePairedDelimiter{\bracket}{[}{]}
\DeclarePairedDelimiter{\paren}{(}{)}

\newcommand*{\innerprod}[3][]{%
  \langle{#2,#3}\rangle\ifnotempty{#1}{_{#1}}
}

\newcommand*{\norm}[2][]{%
  \lVert{#2}\rVert\ifnotempty{#1}{_{#1}}%
}
\newcommand*{\normscaled}[2][]{%
  \left\lVert{#2}\right\rVert\ifnotempty{#1}{_{#1}}%
}

\newcommand*{\largesum}[2][]{%
  \;\;\sum_{\mathclap{#2}}\ifnotempty{#1}{^{\mathclap{#1}}}\;\;%
}

\newcommand*{\notni}{\not\ni}

\newcommand*{\pop}{\texttt{pop}}
\newcommand*{\push}{\texttt{push}}

\newcommand*{\yes}{%
  \tikz[x=0.7em,y=0.7em]{
    \draw[C4,line width=1.4pt] (0,0.3) -- (0.3,0) -- (1,1);
  }%
}
\newcommand*{\no}{%
  \tikz[x=0.7em,y=0.7em]{
    \draw[C1,line width=1.4pt] (0,0) -- (1,1) (0,1) -- (1,0);
  }%
}

\newcommand*{\sgpp}{SG\textsuperscript{++}\xspace}

\newgacronym{bfs}{breadth-first search}
\newgacronym{crra}{constant relative risk aversion}
\newgacronym[DAG]{dagr}{directed acyclic graph}
\newgacronym{fem}{finite element method}
\newgacronym[HFT]{hftr}{hierarchical fundamental transformation}
\newgacronym{iga}{isogeometric analysis}
\newgacronym{lhs}{left-hand side}
\newgacronym{nurbs}{non-uniform rational B-splines}
\newgacronym{pde}{partial differential equation}
\newgacronym{rhs}{right-hand side}
\newgacronym{spd}{symmetric positive definite}
\newgacronym[TIFT]{tiftr}{translation-invariant fundamental transformation}
\newgacronym{up}{unidirectional principle}

\hidegacronym{iga}
\hidegacronym{nurbs}

\begin{document}
\pagenumbering{roman}
\setcounter{page}{1}

\includepdf{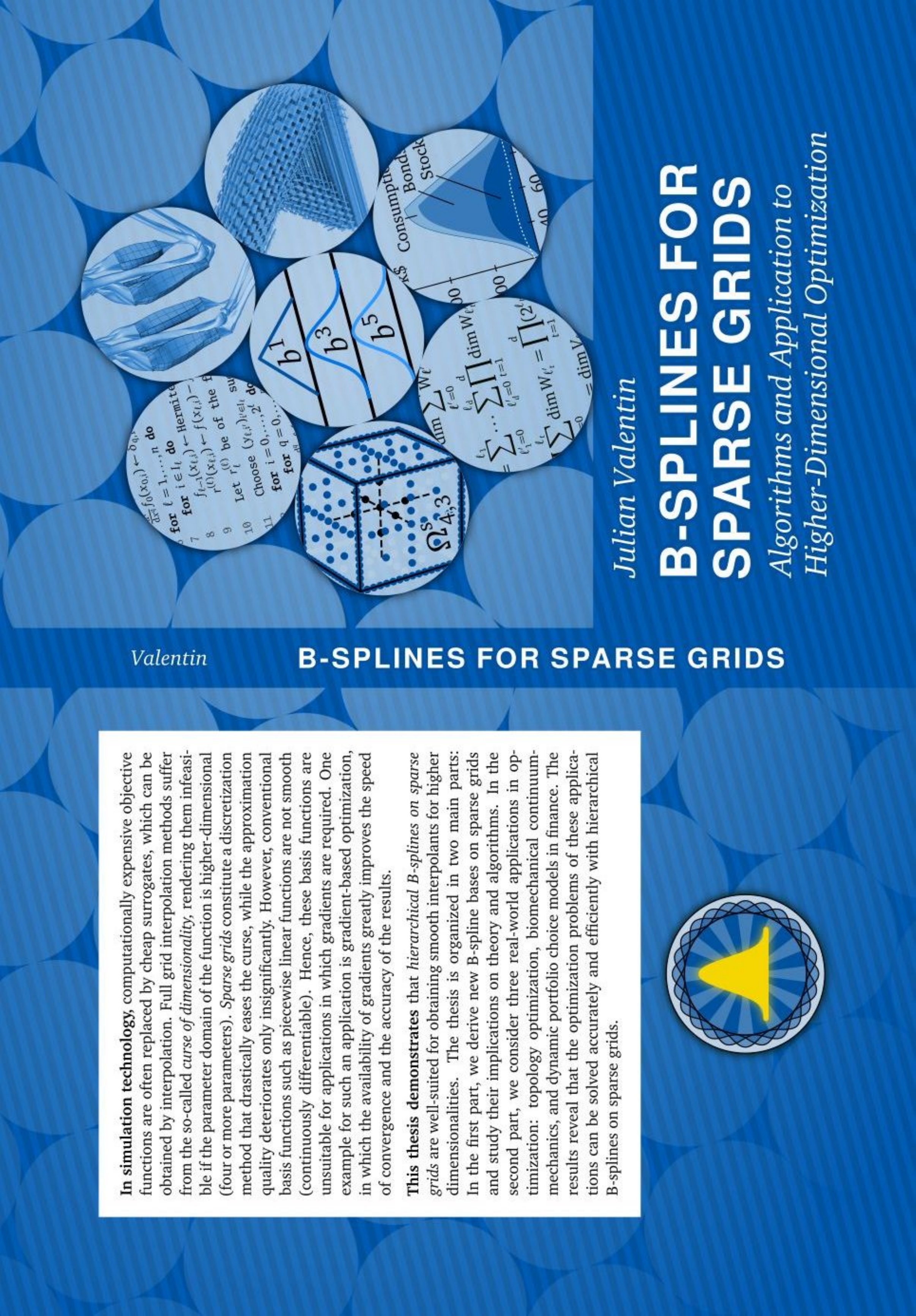}

\cleardoublepage

\pagenumbering{arabic}
\setcounter{page}{1}

  \iftoggle{partialCompileMode}{}{
    \begin{titlepage}
  \begin{spacing}{1}
    \begin{center}
      \begin{otherlanguage}{ngerman}
        \setlength{\parindent}{0pt}
        
        {\bfseries\huge\thetitle\par}
        
        \vfill
        
        \theapproval
        
        \vfill
        
        Vorgelegt von
        
        {\bfseries\Large\theauthor\par}
        
        aus \thebirthplace
        
        \vfill
        
        \begin{tabular}{ll}
          Hauptberichter:&
          \theadvisor\\[0.5em]
          Mitberichter:&
          \theexamineri\\
          &\theexaminerii\\[1em]
          \multicolumn{2}{l}{%
            Tag der mündlichen Prüfung:\quad%
            \thedefensedate%
          }
        \end{tabular}
        
        \vfill
        
        \includegraphics[scale=1.1]{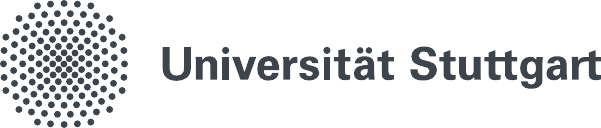}
        
        \vspace{2em}
        
        \theinstitute{} der \theuniversity
        
        \vspace{1em}
        
        \theyear
      \end{otherlanguage}
    \end{center}
  \end{spacing}
\end{titlepage}

\thispagestyle{empty}

{%
  \setlength{\parindent}{0pt}%
  \small
  
  \begin{center}
    \includegraphics[scale=1.1]{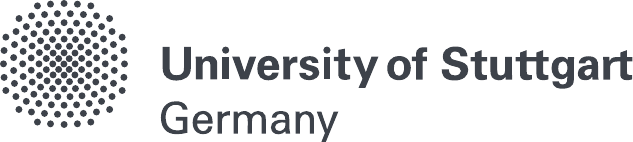}%
    
    \vspace{1em}
    
    Submitted to the University of Stuttgart
  \end{center}
  
  \begin{tabular}{@{}p{0.57\textwidth}@{}p{0.43\textwidth}@{}}
    \emph{Involved institutions and departments:}%
    \vspace{0.6mm}\newline%
    Cluster of Excellence in Simulation Technology%
    \vspace{0.6mm}\newline%
    Institute for Parallel and Distributed Systems%
    \vspace{0.6mm}\newline%
    Chair of Simulation Software Engineering%
    \vspace{0.6mm}\newline%
    Chair of Simulation of Large Systems&
    \raisebox{-0.4\height}{%
      \includegraphics[scale=1.0]{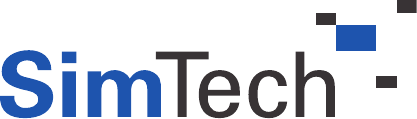}%
    }%
    \hspace{5mm}%
    \raisebox{-0.5\height}{%
      \includegraphics[scale=1.0]{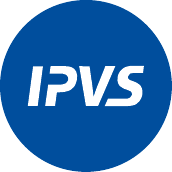}%
    }%
    \vspace{1.3mm}\newline%
    \raisebox{-0.3888\height}{%
      \includegraphics[scale=1.0]{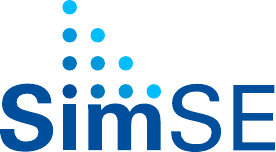}%
    }%
    \hspace{6mm}%
    \raisebox{-0.5\height}{%
      \includegraphics[height=19mm]{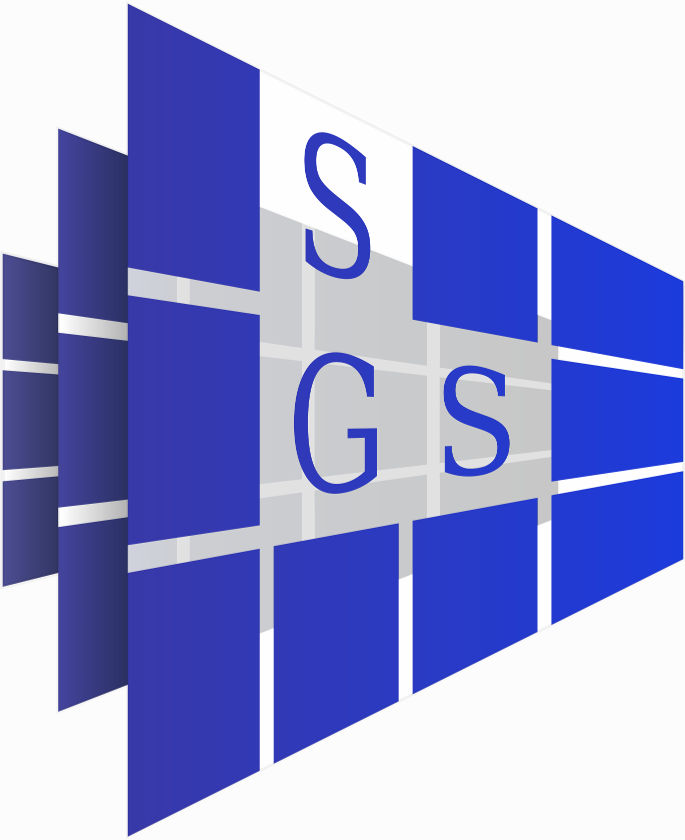}%
    }
  \end{tabular}
  
  \vfill
  
  Julian Valentin\\
  Simulation Software Engineering\\
  Institute for Parallel and Distributed Systems\\
  University of Stuttgart\\
  Universitätsstr. 38\\
  70569 Stuttgart\\
  Germany
  
  \vfill
  
  D\,93 (dissertation)
  
  \vspace{1em}
  
  Compiled as version \compileCounterText{} on \currentTimeLong{}.\\
  Committed as \gitCommitText{} on \gitCommitTimeLong{}.
  
  \vspace{1em}
  
  Typeset using \LaTeX{} and cover design by the author.
  
  Copyright \copyright{} \theyear{} \theauthor{}.
  
  \vspace{1em}
  
  \hangindent=29mm%
  \hangafter=0%
  This work is licensed under the\\
  \smash{%
    \rlap{%
      \hspace*{-29mm}\includegraphics[width=25mm]{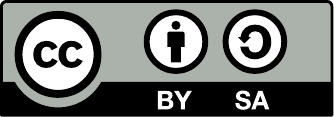}%
    }%
  }%
  Creative Commons Attribution-ShareAlike 4.0 International License.
  
  To view a copy of this license, visit
  \href{%
    https://creativecommons.org/licenses/by-sa/4.0/%
  }{%
    \texttt{https://creativecommons.org/licenses/by-sa/4.0/}%
  }
  or send a letter to
  Creative Commons, PO Box 1866, Mountain View, CA 94042, USA.
  
  \vspace{1em}
  
  Although this thesis was written with utmost care,
  it cannot be ruled out that it contains errors.\\
  Please send any corrections and mistakes to
  \href{mailto:thesis@bsplines.org}{\texttt{thesis@bsplines.org}}.
}

\cleardoublepage

    \vspace*{\fill}

\begin{center}
  \includegraphics[width=0.6\textwidth]{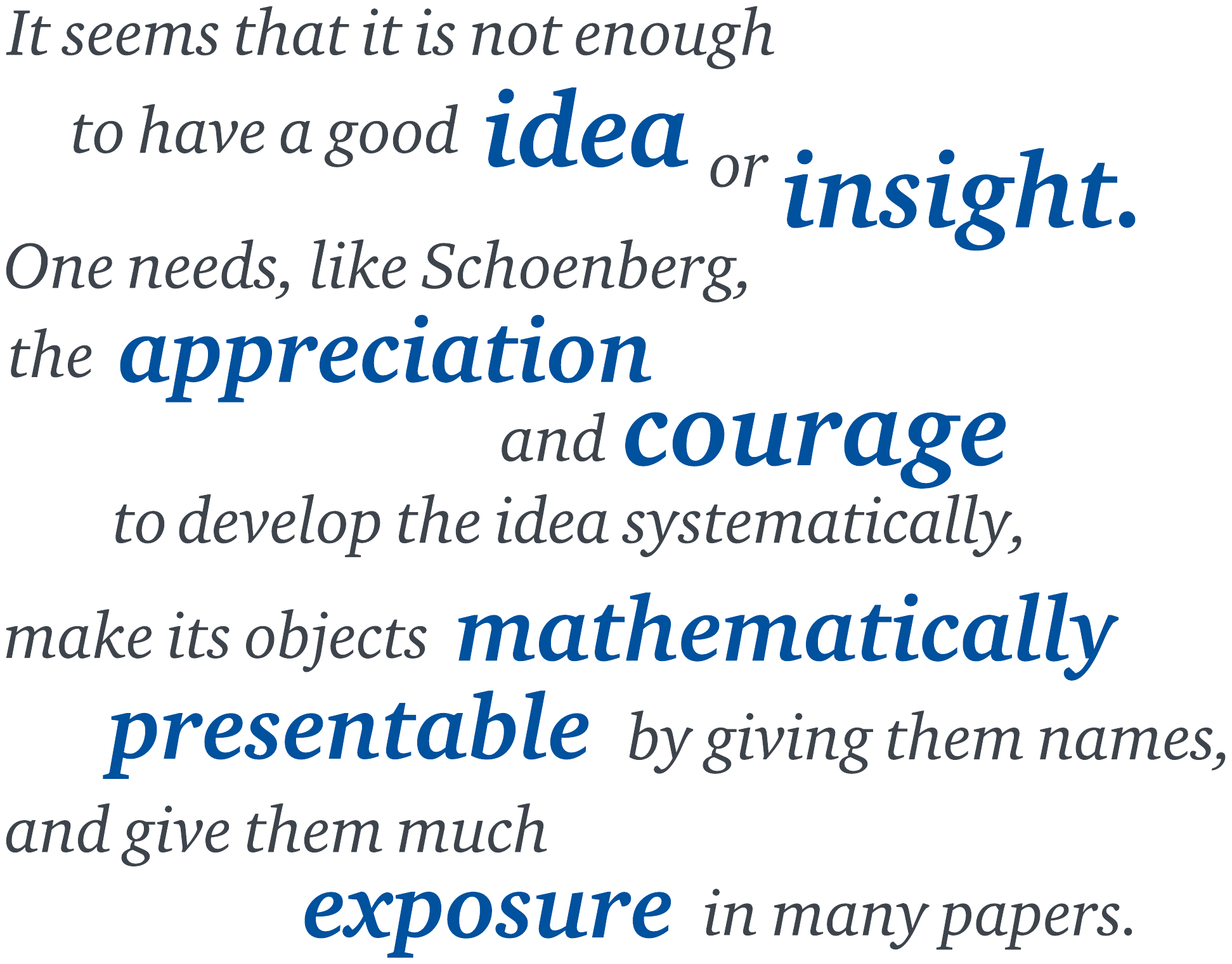}
  
  \begin{minipage}{0.6\textwidth}%
    \begin{flushright}
      \small--- Carl de~Boor \cite{Boor16Comment}
    \end{flushright}
  \end{minipage}
\end{center}

\vspace*{\fill}

\begin{center}
  \includegraphics{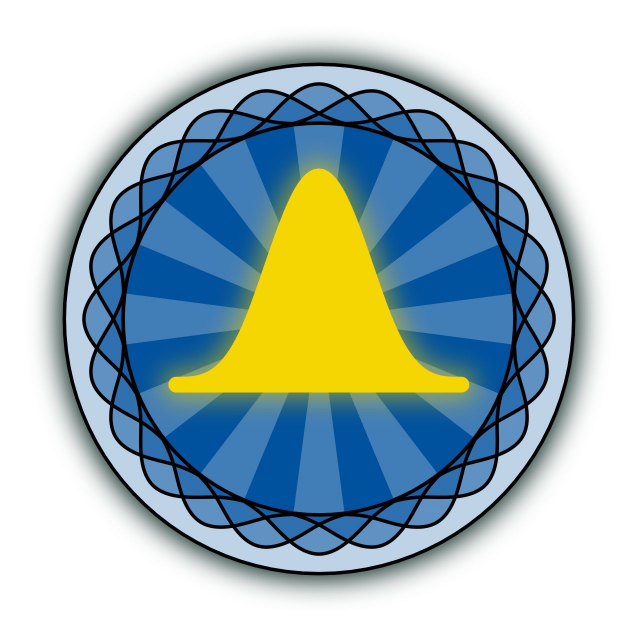}%
\end{center}

\vspace*{\fill}

\clearpage

\vspace*{\fill}

{%
  \setlength{\parindent}{0pt}
  
  Disclaimer:
  This version of the thesis slightly differs from the
  original published version (\texttt{published-v1})
  that was submitted to the University of Stuttgart.
  Only typos have been fixed, and no additional changes have been made.
  The original version can be obtained as a book or as a PDF file from
  the library of the University of Stuttgart.%
}

\cleardoublepage

  }
  
  \iftoggle{draftMode}{
    \input{tex/document/progress}
  }{}
  
  \iftoggle{partialCompileMode}{}{
\pdfbookmark[section]{\contentsname}{toc}
\tableofcontents

\cleardoublepage

    \addlongchap{%
  Lists of Figures, Tables, Algorithms, and Theorems%
}{%
  Lists of Figures, Tables,\texorpdfstring{\\}{ }Algorithms, and Theorems%
}{%
  Lists of Figures, Tables, Algorithms, and Theorems%
}

\disableornamentsfornextheadingtrue
\listoffigures
\disableornamentsfornextheadingtrue
\listoftables
\disableornamentsfornextheadingtrue
\listofalgorithms
\disableornamentsfornextheadingtrue
\listofmytheorems

\cleardoublepage

    \printglossary[title={List of Symbols and Acronyms}]

\cleardoublepage

    \addchap{Abstract/\foreignlanguage{ngerman}{Kurzzusammenfassung}}

\disableornamentsfornextheadingtrue
\section*{Abstract}

In simulation technology, computationally expensive objective functions
are often replaced by cheap surrogates,
which can be obtained by interpolation.
Full grid interpolation methods suffer from the
so-called curse of dimensionality,
rendering them infeasible if the parameter domain of the function
is higher-dimensional (four or more parameters).
Sparse grids constitute a discretization method that drastically eases the
curse, while the approximation quality deteriorates only insignificantly.
However, conventional basis functions such as piecewise linear functions
are not smooth (continuously differentiable).
Hence, these basis functions are unsuitable for applications
in which gradients are required.
One example for such an application is gradient-based optimization,
in which the availability of gradients greatly improves the speed of
convergence and the accuracy of the results.

This thesis demonstrates that hierarchical B-splines on sparse grids are
well-suited for obtaining smooth interpolants for higher dimensionalities.
The thesis is organized in two main parts:
In the first part, we derive new B-spline bases on sparse grids and study
their implications on theory and algorithms.
In the second part, we consider three real-world applications in optimization:
topology optimization, biomechanical continuum-mechanics, and
dynamic portfolio choice models in finance.
The results reveal that the optimization problems of these applications
can be solved accurately and efficiently with hierarchical B-splines on
sparse grids.

\newpage

\begin{otherlanguage}{ngerman}
  \disableornamentsfornextheadingtrue
  \section*{Kurzzusammenfassung}
  
  In der Simulationstechnik werden zeitaufwendige Zielfunktionen
  oft durch einfache Surrogate ersetzt, die durch Interpolation
  gewonnen werden können.
  Vollgitter-Interpola\-tions\-methoden leiden unter dem
  sogenannten Fluch der Dimensionalität,
  der sie unbrauchbar macht, falls der Parameterbereich der Funktion
  höherdimensional ist (vier oder mehr Parameter).
  Dünne Gitter sind eine Diskretisierungsmethode, die den Fluch drastisch
  lindert und die Approximationsqualität nur leicht verschlechtert.
  Leider sind konventionelle Basisfunktionen wie die stückweise
  linearen Funktionen nicht glatt (stetig differenzierbar).
  Daher sind sie für Anwendungen ungeeignet, in denen Gradienten
  benötigt werden.
  Ein Beispiel für eine solche Anwendung ist gradientenbasierte Optimierung,
  in der die Verfügbarkeit von Gradienten die Konvergenzgeschwindigkeit und
  die Ergebnisgenauigkeit deutlich verbessert.
  
  Diese Dissertation demonstriert, dass hierarchische B-Splines auf
  dünnen Gittern hervorragend geeignet sind,
  um glatte Interpolierende für höhere Dimensionalitäten zu erhalten.
  Die Dissertation ist in zwei Hauptteile gegliedert:
  Der erste Teil leitet neue B-Spline-Basen auf dünnen Gittern her und
  untersucht ihre Implikationen bezüglich Theorie und Algorithmen.
  Der zweite Teil behandelt drei Realwelt-Anwendungen aus der Optimierung:
  Topologieoptimierung, biomechanische Kontinuumsmechanik und
  Modelle der dynamischen Portfolio-Wahl in der Finanzmathematik.
  Die Ergebnisse zeigen, dass die Optimierungsprobleme dieser
  Anwendungen durch hierarchische B-Splines auf dünnen Gittern
  genau und effizient gelöst werden können.
\end{otherlanguage}

\cleardoublepage

    \addchap{Preface}

Before I start, I want to thank my advisor Dirk Pflüger.
His ideas and his valuable input have driven me in
my time as PhD student.
I also thank Stephen Roberts and Martin Radetzki
for their readiness to examine my thesis.

Similarly, I thank
Peter Schober
(Prof.\ Dr.\ Raimond Maurer, Goethe University Frankfurt),
Daniel Hübner
(Prof.\ Dr.\ Michael Stingl, FAU Erlangen-Nürnberg),
Michael Sprenger
(Prof.\ Oliver Röhrle, PhD, SimTech/University of Stuttgart), and
Stefan Zimmer
(University of Stuttgart) for
the collaborations and for the enlightening discussions.
I am also grateful for the exciting time with the whole group of
SSE (Simulation Software Engineering) and
SGS (Simulation of Large Systems),
for which I want to thank all past and current PhD students and postdocs of
Dirk Pflüger and Miriam Mehl.
I thank Benjamin, Carolin, Gregor, Henriette, Malte, Michael, Peter,
Ralf, and Theresa
for thoroughly proofreading parts of drafts of this thesis.

Probably the most important role for the success of my PhD thesis
has played my family.
Without their moral support and distraction from the daily work,
I doubt that this thesis would have been possible.

Likewise, I am very grateful for the financial support from
the \foreignlanguage{ngerman}{Juniorprofessurenprogramm} of the
\foreignlanguage{ngerman}{Baden-Württemberg Stiftung}.
I thank the SimTech Cluster of Excellence for supporting
my three-month research stay in Canberra, Australia.

In addition, I want to thank the open-source community for making it possible to
write this thesis in an aesthetically sophisticated manner.
The list of software that was used to write this thesis includes
\LaTeX, Lua\LaTeX, Bib\LaTeX,
\scalebox{0.9}{\KOMAScript}, Ti\emph{k}Z, Python, Matplotlib,
and many more.

\label{page:preface}
Now, I wish that you, dear reader, obtain as much insight as possible
while reading the remaining
\pagedifference{page:preface}{LastPage} pages of this thesis.

Enjoy!

\vspace{1em}

\noindent
Stuttgart, \thedate

\vspace{0.5em}

\noindent
\includegraphics[scale=0.8]{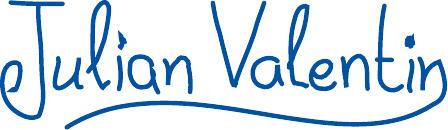}

\cleardoublepage

  }
  
  \setdictum[0.57\textwidth]{%
  There is a fine line between wrong and visionary.
  Unfortunately, you have to be a visionary to see it\dots%
}{%
  Sheldon Cooper (The Big Bang Theory)%
}

\chapter{Introduction}
\label{chap:10introduction}

\initial{0.1em}{B}{efore simulation became available}
as a widespread tool,
knowledge in science and engineering could only advance through
theoretical or experimental considerations.
Nowadays, processes can be simulated
that would be too complicated or even impossible
to be studied theoretically or experimentally,
justifying that simulation is widely viewed as the
``third pillar in knowledge acquisition''
besides theory and experiments \cite{Bungartz14Modeling}.

However, simulations cannot be performed without constructing
a suitable model beforehand (e.g., based on first principles).
Such a model often depends on a number of uncertain or unknown parameters.
Simulations can only represent the real-world circumstances well
if the parameters are well-chosen.
The problem of determining appropriate values for these parameters,
given experimental data, is known as the \term{inverse problem.}
Unfortunately, the solution process of such an inverse problem is non-trivial:
Inverse problems are equivalent to optimization problems of the form
\begin{equation}
  \label{eq:generalOptimizationProblem}
  \min \objfun(\*x),\quad
  \*x \in \real^d\;\;\text{s.t.}\;\;
  \ineqconfun(\*x) \le \*0,
\end{equation}
where $\objfun(\*x)$ gives, for instance,
a measure for the error between the simulation for the model with
the parameter $\*x$ and experimental real-world data and
$\ineqconfun$ constrains the set of feasible parameters.
Since simulations are often time-consuming,
exhaustive search fails
if the dimensionality $d$ is already moderately large ($d > 4$):
Full grid approaches sample each dimension of the domain
independently and construct the Cartesian product of the univariate
samples.
The number of resulting full grid points
grows exponentially in the dimensionality $d$,
which is known as the \term{curse of dimensionality} \cite{Bellman61Adaptive}.

Inverse problems are a motivating example for optimization problems
of the form \eqref{eq:generalOptimizationProblem}, which we will consider in
this thesis as our main application.
The curse affects not only optimization,
but also various tasks such as interpolation, quadrature, and regression,
which play a vital role in computational science.

\paragraph{Sparse grid surrogates}

The surrogate-based approach we pursue in this thesis is simple and
yet powerful:
Instead of directly optimizing the expensive objective function $\objfun$,
we replace it with a surrogate $\sgintp$ that can be evaluated cheaply.
We choose interpolation as our method to construct the surrogates,
although other methods such as quasi-interpolation
\cite{Hoellig13Approximation} or regression \cite{Pflueger10Spatially} exist.
Again, conventional full grid interpolation schemes are afflicted by the
curse of dimensionality, which rules them out for our purposes.

This is where \term{sparse grids} come into play.
In their simplest form, sparse grids give an \term{a priori} selection of
full grid points and corresponding basis functions such that
the exponential dependency of the grid size on the dimensionality
is removed, while not deteriorating the $\Ltwo$ interpolation error too much
\cite{Bungartz04Sparse}.
However, sparse grids can also be employed spatially adaptively,
where grid points are refined \term{a posteriori} according to suitable
refinement criteria.
This is of particular interest for the scope of this thesis,
as spatial adaptivity enables us to increase the
accuracy in regions of interest,
simultaneously keeping the number of grid points at an acceptable level.

\paragraph{B-splines}

Conventional basis functions for sparse grids
(most common are piecewise linear functions)
are not continuously differentiable.
This poses problems for gradient-based optimization algorithms,
which use the gradient $\gradient{\*x}{f}$ of the
objective function $\objfun$ to update the search direction.
Employing finite differences as a remedy
is time-consuming and introduces new error sources.

Previous studies
\multicite{Pflueger10Spatially,Sickel11Spline,Valentin14Hierarchische}
suggest that \term{hierarchical B-splines}
as sparse grid basis functions may significantly improve results.
B-splines of degree~$p$ are $(p-1)$ times continuously differentiable
piecewise polynomials of degree $p$ that form a basis of the space of splines.
As the derivatives of B-splines can be evaluated fast and explicitly,
the convergence of gradient-based optimization techniques is
greatly accelerated.
Additionally, the higher order of B-splines increases the accuracy
of surrogates obtained by interpolation
when compared to piecewise linear bases.

\paragraph{Main goals}

So far, there is no work that brings sparse grids and B-splines together,
thoroughly examining the theoretical implications on algorithms and
assessing the practical performance in real-world applications.
This thesis addresses this very intersection of theory and practice of
\term{B-splines for sparse grids.}
The main goals of the thesis are
\begin{itemize}
  \item
  to establish a consistent notational and theoretical
  framework for sparse grids with general basis functions,
  
  \item
  to construct new algorithmically efficient
  B-spline-based basis function types for sparse grids,
  
  \item
  to study the algorithmic properties of the new bases and
  to formulate suitable new algorithms,
  while proving their formal correctness, and
  
  \item
  to apply the new bases and algorithms to different real-world
  optimization scenarios.
\end{itemize}
These goals set the agenda for the outline of the rest of the thesis,
which is described in the following.

\paragraph{Outline}

We start in \cref{chap:20sparseGrids} by defining sparse grids
for arbitrary tensor product basis functions.
The advantage of introducing the notation independently
of the type of basis functions is that
different hierarchical B-spline bases can be substituted easily
for the following chapters.
We first define sparse grids with points on the boundary of the domain
and then study options for the boundary treatment
(as opposed to some literature
\multicite{Bungartz04Sparse,Pflueger10Spatially}).

In \cref{chap:30BSplines}, we define
the standard hierarchical B-spline basis for sparse grids.
In addition, we construct various new hierarchical
B-spline basis types, such as
non-uniform B-splines (e.g., Clenshaw--Curtis B-splines) and
modified B-splines.
A mismatch of dimensions between the uniform spline space and the
hierarchical B-spline space implies that, surprisingly,
polynomials cannot be replicated by the
standard hierarchical B-spline basis.
Hence, we have to incorporate specific boundary conditions
\term{(not-a-knot conditions)} into the hierarchical B-splines,
which we explain in the second half of \cref{chap:30BSplines}.

The new hierarchical B-spline bases call for novel algorithmic approaches
to solve numerical tasks such as
hierarchization (interpolation) on sparse grids.
In \cref{chap:40algorithms}, we show new algorithms for
spatially adaptive sparse grids with the example problem of hierarchization
based on existing algorithms,
which work for B-splines only in specific special cases.
In the course of \cref{chap:40algorithms}, we construct several new
hierarchical B-spline basis types to enable the
applicability of the new algorithms.
As mentioned above, we prove the formal correctness of every algorithm
that we repeat from the literature or develop from scratch.

\Cref{chap:50optimization} shows how to apply B-splines on sparse grids
to gradient-based optimization problems.
We briefly discuss different optimization scenarios and how to solve them
with various gradient-free and gradient-based optimization techniques.
Numerical results are given for a number of test scenarios as well
as for an example application from fuzzy arithmetic.

Three real-world applications follow in
\cref{chap:60topoOpt,chap:70muscle,chap:80finance}.
In these chapters,
the theoretical knowledge gained in the first half of the thesis
is applied to the solution of the three real-world optimization problems:

First, in \cref{chap:60topoOpt},
we study topology optimization via a homogenized two-scale approach.
For this application, the key ingredient is an interpolation scheme
that preserves both the positive definiteness of the interpolated tensors and
their explicit differentiability.

Second, in \cref{chap:70muscle},
we consider a biomechanical application in which the interpolated data values
are the result of very expensive continuum-mechanical calculations.
The optimization problems posed in this chapter ask for
muscle activation levels such that a specific joint angle is attained,
which is a recurring problem in medicine and robotics.
B-spline surrogates on sparse grids
decrease the necessary computing time significantly.

Third, in \cref{chap:80finance},
we examine dynamic portfolio choice models.
This financial application features some peculiarities that have
to be considered when solving the corresponding optimization problems.
For instance, it is necessary to evaluate interpolants outside their domain
and calculate integrals due to random factors such as stock returns.

\Cref{chap:90conclusion} concludes the thesis by
summarizing its results and giving an overview of possible future work.
In the appendix, one can find supplementary information such as
technical proofs that are too verbose to be included in the main text.

\paragraph{Original contribution}

This thesis is written to be largely self-contained.
Therefore, it is necessary that some introductory definitions and
results are repeated from the literature,
which is properly attributed in the respective chapters.
In addition, some new results have already been published.
Whenever a publication is co-authored by collaborators,
the original contribution of the author of this thesis is highlighted
at the beginning of the respective chapters or sections.

\paragraph{Notation}

The notation of this thesis should be intuitive and suggestive.
It is designed to be as natural as possible (i.e., not distracting)
and as unambiguous as necessary.
One example is that vectors are written in bold face, which leads to
very similar formulas for the univariate and the multivariate cases.
For instance, $\clint{0, 1}$ and $\sum_{l=0}^n$
become
$\clint{\*0, \*1} = \clint{0, 1}^d$ and
$\sum_{\*l=\*0}^{\*n} = \sum_{l_1=0}^{n_1} \dotsb \sum_{l_d=0}^{n_d}$,
respectively.
This and other necessary notation is introduced in the text when needed.
If a symbol or an abbreviation is unclear,
it is likely explained in the glossary at the beginning of the thesis.

\cleardoublepage

  \setdictum[0.6\textwidth]{%
  We combine two sparse grid approximations
  and call it ``deep sparse grids,\hspace{-0.1em}''
  because everything is deep today!%
}{%
  In a talk at the 5th Workshop on\\Sparse Grids and Applications%
}

\longchapter{%
  Sparse Grids with Arbitrary Tensor Product Bases%
}{%
  Sparse Grids with Arbitrary\texorpdfstring{\\}{ }Tensor Product Bases%
}{%
  Sparse Grids with Arbitrary Tensor Product Bases%
}
\label{chap:20sparseGrids}

\initial[lhang=0.06]{0.1em}{S}{parse grids are a versatile tool}
in computational mathematics and scientific computing.
As already mentioned in \cref{chap:10introduction},
their motivation is to ease the curse of dimensionality,
which states that the number of full grid points
grows exponentially in the dimensionality $d$ of the underlying domain.
Their general formulation and the possibility to employ sparse grids in a
regular, dimensionally adaptive, or spatially adaptive fashion
opens a broad field of theoretical and practical applications
to sparse grids.

Sparse grids have been known for at least half a century,
albeit not under this name.
A paper by Smolyak \cite{Smolyak63Quadrature} is usually regarded
as the first modern treatment of sparse grids in the form
of the combination technique \cite{Garcke13Sparse}.
Additionally, there are close connections to
hyperbolic crosses \cite{Temljakov82Approximation}
and to Boolean interpolation operators
\multicite{Delvos82Dvariate,Delvos89Boolean}.
The term \term{sparse grids} was coined by Zenger in 1991
\cite{Zenger91Sparse}.
Some important subsequent work for hierarchical bases was done by
Bungartz and Griebel
\multicite{%
  Bungartz92Duenne,%
  Griebel92Combination,%
  Bungartz98Finite,%
  Bungartz04Sparse%
}.
Since then, sparse grids have been applied to various fields,
for instance,
data mining
\multicite{Garcke01Data,Pandey08Regression,Pflueger10Spatially},
interpolation
\cite{Sickel11Spline},
quadrature
\cite{Gerstner98Numerical},
density estimation
\multicite{Griebel10Finite,Peherstorfer14Density},
PDEs
\multicite{Balder94Adaptive,Bungartz98Finite,Nobile16Adaptive}, and
optimization
\multicite{Ferenczi05Globale,Donahue09Robust,Valentin16Hierarchical}.
Various software toolboxes for sparse grids have been developed
\multicite{Klimke05Algorithm,Pflueger10Spatially,Stoyanov18User}.
For a general introduction to sparse grids,
see the tutorial by Garcke \cite{Garcke13Sparse} or
the more extensive survey by Bungartz and Griebel
\cite{Bungartz04Sparse}.

This chapter provides a consistent notational framework
for the definition of sparse grids with general basis functions.
The reason not to employ specific bases such as the common hat functions
or B-splines of higher degrees is two-fold:
First, we will define various new ``flavors'' of B-splines,
which is easier if the choice of basis is left open.
Second, most of the statements and theorems that we will make in this
thesis will hold for general basis functions
(in some cases with additional assumptions)
and not just for B-splines.

Besides the derivation of sparse grids with
coarser boundary points in \cref{sec:241coarseBoundary},
this section is mostly
a repetition of the definition of sparse grids with general basis functions.
Our notation and presentation will roughly follow
\cite{Pflueger10Spatially} and \cite{Garcke13Sparse}.
A more detailed introduction to sparse grids can be found in
\cite{Bungartz04Sparse}.
Original contributions of the thesis in this chapter
are the formalization of the hierarchical splitting for
arbitrary basis functions in \cref{sec:21nodalSpaces,sec:22hierSubspaces} and
the definition of sparse grids with coarse boundary in
\cref{sec:241coarseBoundary}.

\section{Nodal Basis and Nodal Space}
\label{sec:21nodalSpaces}

\minitoc{62mm}{4}

\mbox{}\vspace{-14mm}

\disableornamentsfornextheadingtrue
\subsection{Univariate Case}
\label{sec:211nodalUV}

\paragraph{Grid and basis functions}

In this thesis, we consider univariate functions
that are defined on the unit interval $\clint{0, 1}$.
\usenotation{l}
We discretize this domain by splitting it into $2^l$ equally-sized segments,
where $l \in \natz$ is the \term{level.}
\usenotation{i}
The resulting $2^l + 1$ \term{grid points} $\gp{l,i}$ are given by
\begin{equation}
  \gp{l,i} \ceq i \cdot \ms{l},\quad
  i = 0, \dotsc, 2^l,
\end{equation}
where $i$ is the \term{index} and $\ms{l} \ceq 2^{-l}$ is the \term{mesh size.}%
\footnote{%
  Note that from a strict formal perspective,
  this equation defined $\gp{l,i}$ only for $i = 0, \dotsc, 2^l$,
  but we will later need $\gp{l,i}$ also for $i < 0$ or $i > 2^l$.
  The convention in this thesis is that all definitions are
  implicitly generalized whenever needed.%
}
Every grid point is associated with a \term{basis function}
\begin{equation}
  \basis{l,i}\colon \clint{0, 1} \to \real.
\end{equation}
We assume $\basis{l,i}$ to be arbitrary,
satisfying required assumptions when needed and stated.
However, it helps for both the theory and the intuition to have a
specific example of basis functions in mind.
\usenotation{zzzz1}
The so-called \term{hat functions} (linear B-splines)
are the most common choice for $\basis{l,i}$:
\begin{equation}
  \label{eq:hatFunctionUV}
  \bspl{l,i}{1}(x)
  \ceq \max(1 - \abs{\tfrac{x}{\ms{l}} - i}, 0).
\end{equation}
Here and in the following,
the superscript ``1'' stands for the degree of the linear B-spline and
is not to be read as an exponent.
We generalize this notation to B-splines $\bspl{l,i}{p}$ of
arbitrary degrees $p$ in \cref{chap:30BSplines}.

\paragraph{Nodal space}

\begin{figure}
  \subcaptionbox{%
    Basis functions $\bspl{l,i}{1}$ ($i = 0, \dotsc, 2^l$)
    and grid points $\gp{l,i}$ \emph{(dots).}%
  }[72mm]{%
    \includegraphics{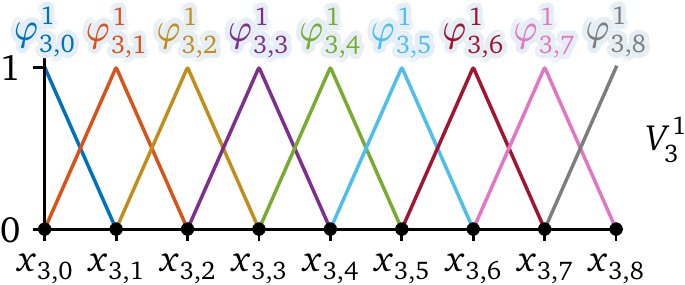}%
  }%
  \hfill%
  \subcaptionbox{%
    Piecewise linear interpolant $\fgintp{l}$
    of some function data $\objfun(\gp{l,i})$
    as a weighted sum of the nodal hat functions.%
  }[72mm]{%
    \includegraphics{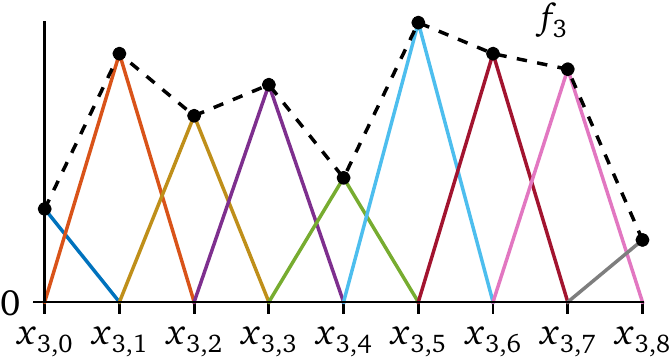}%
  }%
  \caption[%
    Univariate nodal hat functions%
  ]{%
    Univariate nodal hat functions of level $l = 3$.%
  }%
  \label{fig:nodalHat}%
\end{figure}

The \term{nodal space} $\ns{l}$ of level $l$
is defined as the linear span of all basis functions
$\basis{l,i}$:
\begin{equation}
  \ns{l} \ceq \spn\{\basis{l,i} \mid i = 0, \dotsc, 2^l\}.
\end{equation}
We assume that the functions $\basis{l,i}$ form a basis of $\ns{l}$, i.e.,
they are linearly independent.
Consequently, every linear combination of these functions is unique.
This ensures that for every objective function $\objfun\colon \clint{0, 1} \to \real$,
there is a unique function $\fgintp{l}\colon \clint{0, 1} \to \real$ such that
\begin{equation}
  \label{eq:interpFullGridUV}
  \fgintp{l}
  = \sum_{i=0}^{2^l} \interpcoeff{l,i} \basis{l,i},\quad
  \falarge{i = 0, \dotsc, 2^l}{\fgintp{l}(\gp{l,i}) = \objfun(\gp{l,i})},
\end{equation}
for some $\interpcoeff{l,i} \in \real$.
In this case, $\fgintp{l}$ is called \term{interpolant} of $\objfun$ in $\ns{l}$.
The nodal space $\nsbspl{l}{1}$ is defined analogously to $\ns{l}$
as the span of the hat functions $\bspl{l,i}{1}$.
It is the space of all linear splines,
that is, the space of all continuous functions on $\clint{0, 1}$ that are
piecewise linear polynomials on $\clint{\gp{l,i}, \gp{l,i+1}}$ for
$i = 0, \dotsc, 2^l - 1$ \cite{Hoellig13Approximation}.
The nodal hat function basis of level~$l = 3$
and a linear combination are shown in \cref{fig:nodalHat}.

\subsection{Multivariate Case}
\label{sec:212nodalMV}

\paragraph{Cartesian and tensor products}

\begin{SCfigure}
  \includegraphics{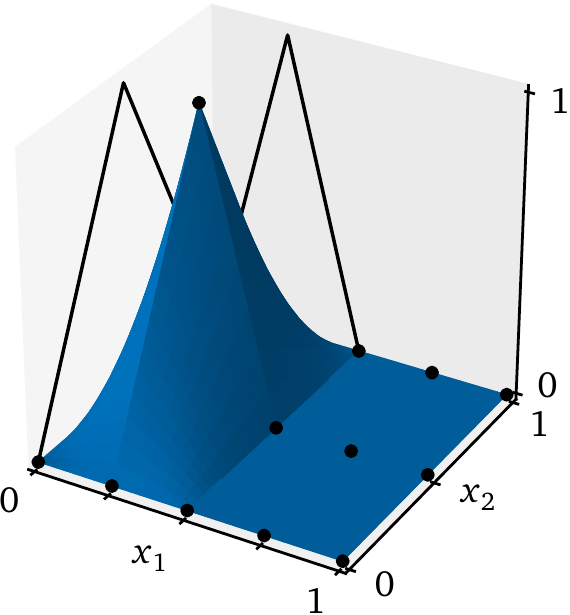}%
  \caption[%
    Bivariate nodal hat function%
  ]{%
    Bivariate nodal hat function of level $\*l = (2, 1)$ and
    index $i = (1, 1)$ as the tensor product of two univariate
    nodal hat functions.%
  }%
  \label{fig:nodalHat2D}%
\end{SCfigure}

\usenotation{d}
For the multivariate case with $d \in \nat$ dimensions,
we employ a tensor product approach,
for which we replace all indices, points, and functions with
multi-indices, Cartesian products, and tensor products, respectively.
\usenotation{@0}
\usenotation{@1}
Therefore, the domain is now $\clint{\*0, \*1} \ceq \clint{0, 1}^d$,
which can be partitioned into
$\prod_{t=1}^d 2^{l_t} = 2^{\normone{\vec{l}}}$ equally-sized hyper-rectangles,
where $\*l = (l_1, \dotsc, l_d) \in \natz^d$ is the $d$-dimensional level
and $\normone{\vec{l}} \ceq \sum_{t=1}^d \abs{l_t}$ is the level sum.
The corners of the hyper-rectangles are given by the grid points
\begin{equation}
  \label{eq:gridPointMultivariate}
  \gp{\*l,\*i} \ceq \*i \cdot \ms{\*l},\quad
  \*i = \*0, \dotsc, \*2^{\*l}.
\end{equation}
Relations and operations with vectors in bold face
are to be read coordinate-wise in this thesis, unless stated otherwise.
Bold-faced numbers like $\*0$ are defined to be the vector $(0, \dotsc, 0)$
in which every entry is equal to that number.
This is to allow a somewhat intuitive and suggestive notation.
For example, \eqref{eq:gridPointMultivariate} is equivalent to
the much longer formula
\begin{equation}
  \gp{\*l,\*i}
  \ceq (i_1 \ms{l_1},\; \dotsc,\; i_d \ms{l_d}),\quad
  i_t = 0, \dotsc, 2^{l_t},\quad
  t = 1, \dotsc, d,
\end{equation}
with the $d$-dimensional mesh size
$\ms{\*l} \ceq \*2^{-\*l} = (\ms{l_1}, \dotsc, \ms{l_d})$.
Again, every grid point is associated with a basis function that is defined
as the tensor product of the univariate functions:%
\footnote{%
  Note that,
  although \cref{eq:tensorProduct} does not cover it,
  one could employ basis functions of different types in
  each dimension, for example B-splines of different degrees.
  All remaining considerations in this thesis
  regarding tensor product basis functions are independent
  of whether we use the same function type or
  different types in each dimension.%
}
\begin{equation}
  \label{eq:tensorProduct}
  \basis{\*l,\*i}\colon \clint{\*0, \*1} \to \real,\quad
  \basis{\*l,\*i}(\*x)
  \ceq \prod_{t=1}^d \basis{l_t,i_t}(x_t).
\end{equation}
\cref{fig:nodalHat2D} shows an example of a bivariate nodal hat function
$\bspl{\*l,\*i}{1}$.

\vspace*{\fill}
\pagebreak

\paragraph{Multivariate nodal space}

The multivariate nodal space $\ns{\*l}$ is defined analogously to
the univariate case:
\begin{equation}
  \ns{\*l}
  \ceq \spn\{\basis{\*l,\*i} \mid \*i = \*0, \dotsc, \*2^{\*l}\}.
\end{equation}
In the case of hat functions $\bspl{\*l,\*i}{1}$,
the nodal space $\nsbspl{\*l}{1}$ is the $d$-linear spline space
\cite{Hoellig13Approximation}, i.e.,
the space of all continuous functions
on $\clint{\*0, \*1}$ that are piecewise $d$-linear polynomials on
all hyper-rectangles
\begin{equation}
  \clint{\gp{\*l,\*i}, \gp{\*l,\*i+\*1}}
  \ceq \clint{\gp{l_1,i_1}, \gp{l_1,i_1+1}} \times \dotsb \times
  \clint{\gp{l_d,i_d}, \gp{l_d,i_d+1}},\quad
  \*i = \*0, \dotsc, \*2^\*l - \*1.
\end{equation}
Analogously to \eqref{eq:interpFullGridUV},
we can interpolate objective functions $\objfun\colon \clint{\*0, \*1} \to \real$
in the nodal space $\ns{\*l}$ with $\fgintp{\*l}\colon \clint{\*0, \*1} \to \real$ satisfying
\begin{equation}
  \label{eq:interpFullGridMV}
  \fgintp{\*l}
  = \sum_{\*i=\*0}^{\*2^\*l} \interpcoeff{\*l,\*i} \basis{\*l,\*i},\quad
  \falarge{\*i = \*0, \dotsc, \*2^\*l}{\fgintp{\*l}(\gp{\*l,\*i}) = \objfun(\gp{\*l,\*i})},
\end{equation}
where $\interpcoeff{\*l,\*i} \in \real$ and
the sum is over all $\*i = \*0, \dotsc, \*2^\*l$
(i.e., $i_t = 0, \dotsc, 2^{l_t}$, $t = 1, \dotsc, d$).
To ensure that the coefficients $\interpcoeff{\*l,\*i}$
exist for every objective function $\objfun$ and are uniquely determined by
the values at the grid points
\begin{equation}
\fgset{\*l}
\ceq \{\gp{\*l,\*i} \mid \*i = \*0, \dotsc, \*2^{\*l}\},
\end{equation}
we prove the following statement:

\vspace*{\fill}
\pagebreak

\begin{lemma}[linear independence of tensor products]
  \label{lemma:tensorProductLinearIndependence}
  The functions $\basis{\*l,\*i}$ ($\*i = \*0, \dotsc, \*2^\*l$)
  form a basis of $\ns{\*l}$, if the univariate functions
  $\basis{l_t,i_t}$ ($i_t = 0, \dotsc, 2^{l_t}$)
  form a basis of the univariate nodal space $\ns{l_t}$
  for $t = 1, \dotsc, d$.
\end{lemma}
\begin{proof}
  Assume that $\interpcoeff{\*l,\*i} \in \real$ are chosen in \eqref{eq:interpFullGridMV}
  such that $\fgintp{\*l} \equiv 0$.
  Then for all $\*i' = \*0, \dotsc, \*2^\*l$,
  we can evaluate \eqref{eq:interpFullGridMV} at $\gp{\*l,\*i'}$ to obtain
  \begin{equation}
    \sum_{i_1=0}^{2^{l_1}}
    \paren*{
      \sum_{i_2=0}^{2^{l_2}} \dotsb \paren*{
        \sum_{i_d=0}^{2^{l_d}}
        \interpcoeff{\*l,\*i} \basis{l_d,i_d}(\gp{l_d,i_d'})
      } \dotsb \basis{l_2,i_2}(\gp{l_2,i_2'})
    } \basis{l_1,i_1}(\gp{l_1,i_1'})
    = 0.
  \end{equation}
  We apply the univariate linear independence ($x_1$ direction) to infer
  that the sum over $i_2$ must vanish for all $i_1 = 0, \dotsc, 2^{l_1}$.
  Repeating this argument for all dimensions, we have
  $\interpcoeff{\*l,\*i} = 0$ for all~$\*i = \*0, \dotsc, \*2^\*l$,
  implying the linear independence of the functions $\basis{\*l,\*i}$.
\end{proof}

\usenotation{n10}
A common choice for the level $\*l$ is $n \cdot \*1$ for some $n \in \natz$.
\usenotation{Vnd}
In this case, we replace ``$\*l$'' in the subscripts with ``$n{,}d$''
(for example, $\ns{n,d} \ceq \ns{n \cdot \*1}$).
For the hat function basis $\bspl{\*l,\*i}{1}$,
it can be shown that the $\Ltwo$ interpolation error of the interpolant
$\fgintp{n,d} \in \ns{n,d}$ is given by
\begin{equation}
  \normLtwo{\objfun - \fgintp{n,d}} = \landauO{\ms{n}^2},
\end{equation}
i.e., the order of the interpolation error is quadratic in the mesh size
\multicite{Hoellig13Approximation,Bungartz04Sparse}.

\section{Hierarchical Basis and Hierarchical Subspace}
\label{sec:22hierSubspaces}

\minitoc{67mm}{5}

\noindent
The dimension of the nodal space $\ns{\*l}$ is given by
\begin{equation}
  \label{eq:dimensionFG}
  \dim \ns{\*l}
  = \setsize{\fgset{\*l}}
  = \prod_{t=1}^d (2^{l_t} + 1).
\end{equation}
If we choose the same level $n \in \natz$ in all dimensions,
then the dimension of $\ns{n,d}$ and the
number of grid points grow at least as fast as
$2^{nd} = (\ms{n}^{-1})^d$.
This exponential dependency between $\dim \ns{n,d}$ and $d$ is known as the
\term{curse of dimensionality} \cite{Bellman61Adaptive}.
The curse makes interpolation on $\ns{\*l}$ computationally infeasible
for dimensionalities $d > 4$,
as we would have to calculate and store
$\dim(\ns{\*l})$-many coefficients $\interpcoeff{\*l,\*i}$.%

\subsection{Hierarchical Splitting in the Univariate Case}
\label{sec:221hierUV}

\paragraph{Hierarchical subspaces}

\begin{SCfigure}
  \includegraphics{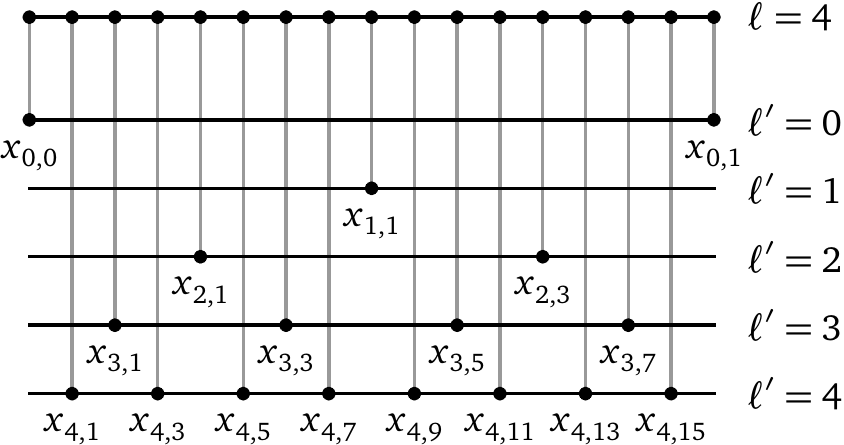}%
  \caption[%
    Decomposition of the set of univariate grid points%
  ]{%
    The set of grid points $\fgset{l}$ of level $l = 4$ \emph{(top)}
    decomposes into hierarchical grids of level $l' \le l$,
    whose grid points $\gp{l',i'}$ have odd indices $i' \in \hiset{l'}$
    ($\gp{0,0}$ being the only exception).%
  }%
  \label{fig:pointSplittingUniform}%
\end{SCfigure}

\begin{figure}
  \subcaptionbox{%
    Basis functions $\bspl{l',i'}{1}$ ($l' \le l$, $i' \in \hiset{l'}$)
    and grid points $\gp{l',i'}$ \emph{(dots).}
    The domain is the unit interval $\clint{0, 1}$.%
  }[72mm]{%
    \includegraphics{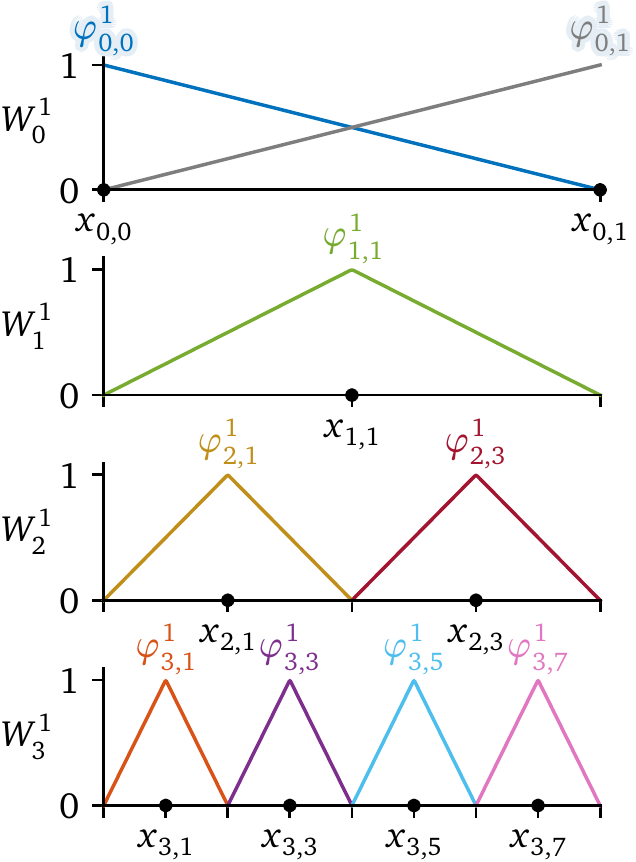}%
  }%
  \hfill%
  \subcaptionbox{%
    Piecewise linear interpolant $\fgintp{l}$
    of some function data $\objfun(\gp{\*l,\*i})$
    as a linear combination of hierarchical hat functions \emph{(stacked).}
    The two boundary functions are combined to a single function
    \emph{\textcolor{C0}{(blue)}} for simplicity.%
  }[72mm]{%
    \includegraphics{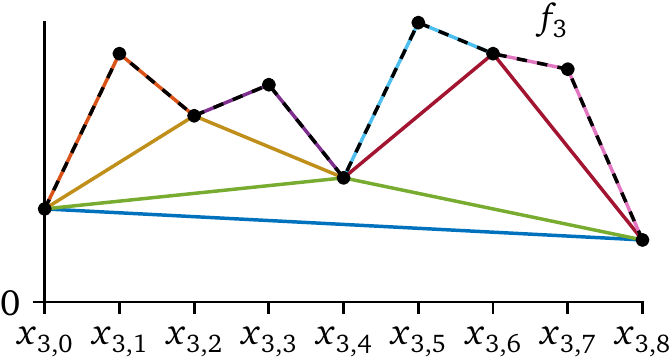}%
  }%
  \caption[%
    Univariate hierarchical hat functions%
  ]{%
    Univariate hierarchical hat functions up to level $l = 3$.%
  }%
  \label{fig:hierarchicalHat}%
\end{figure}

In order to reduce the computational effort,
we first split $\ns{\*l}$ into smaller subspaces and then identify
subspaces that we can omit at the cost of a slightly larger error.
In the univariate case, the key observation is that a grid point of a
level $l$ can be written as a grid point of a higher level~$l'$:
\begin{equation}
  \label{eq:rewriteGridPoint}
  \gp{l,i} = \gp{l',i'},\quad
  l' \ge l,\quad
  i' = 2^{l'-l} i.
\end{equation}
Conversely, this implies that every grid point $\gp{l,i}$ of level $l \ge 1$
and index $i \ge 1$ can be uniquely written
as a grid point of a coarser level $l'$ (or $l' = l$) and an odd index $i'$:
\begin{equation}
  \gp{l,i} = \gp{l',i'},\quad
  l' = l - \bracket*{\log_2(\xor(i, i-1) + 1) - 1},\quad
  i' = 2^{l'-l} i,
\end{equation}
where $\xor$ is the bitwise ``exclusive or'' function.
The term in square brackets is the exponent of the
highest power of two that divides $i$.
The two boundary points zero and one are obtained by
inserting an additional level $l' = 0$ with indices $i' \in \{0, 1\}$.
As shown in \cref{fig:pointSplittingUniform},
this implies that $\fgset{l}$ decomposes into
\begin{equation}
  \fgset{l}
  = \bigdotcup_{l'=0}^l \{\gp{l',i'} \mid i' \in \hiset{l'}\},\quad
  \hiset{l'} \ceq
  \begin{cases}
    \{i' = 0, \dotsc, 2^{l'} \mid \text{$i'$ odd}\},&l' > 0,\\
    \{0, 1\},&l' = 0,
  \end{cases}
\end{equation}
where $\dotcup$ indicates the disjoint union.
We call the spaces spanned by the basis functions that correspond to the
index sets $\hiset{l'}$ \term{hierarchical subspaces} $\hs{l'}$:
\begin{equation}
  \hs{l'}
  \ceq \spn\{\basis{l',i'} \mid i' \in \hiset{l'}\}.
\end{equation}
The corresponding basis functions
$\basis{l',i'}$, $l' = 0, \dotsc, l$, $i' \in \hiset{l'}$,
are called \term{hierarchical basis functions.}
The hierarchical hat function basis is shown in \cref{fig:hierarchicalHat}.

\paragraph{Hierarchical splitting}

For the hat function basis $\bspl{l,i}{1}$ and other basis types,
we can prove that the corresponding nodal space
decomposes into the direct sum of all
hierarchical subspaces of coarser levels or the same level, i.e.,
\begin{equation}
  \label{eq:hierSplittingUV}
  \ns{l}
  \overset{?}{=} \bigoplus_{l'=0}^l \hs{l'}.
\end{equation}
We call this relation \term{hierarchical splitting.}
Here, the direct sum $\oplus$ is
the vector space sum that additionally indicates
that the dimension of the sum $\sum_{l'=0}^l \hs{l'}$ is the sum
of the dimensions of the summands $\hs{l'}$
(analogously to
$\setsize{\fgset{l}}
= \sum_{l'=0}^l \setsize{\{\gp{l',i'} \mid i' \in \hiset{l'}\}}$,
where $\fgset{l}$ is the disjoint union of the sets
$\{\gp{l',i'} \mid i' \in \hiset{l'}\}$).
In general, \eqref{eq:hierSplittingUV} may not be true,
depending on the type of basis functions.
The following lemma provides a characterization
that can be used to prove \eqref{eq:hierSplittingUV} for hat functions.

\vspace*{\fill}
\pagebreak

\begin{lemma}[univariate hierarchical splitting characterization]
  \label{lemma:hierSplittingUV}
  \Cref{eq:hierSplittingUV} is equivalent to the satisfaction of
  both of the following conditions:
  \begin{itemize}
    \item
    The hierarchical subspaces $\hs{l'}$ ($l' \le l$)
    are subspaces of $\ns{l}$.
    
    \item
    The hierarchical functions
    $\basis{l',i'}$ ($l' \le l$, $i' \in \hiset{l'}$)
    are linearly independent.
  \end{itemize}
\end{lemma}
\begin{proof}
  The first condition is equivalent to $\sum_{l'=0}^l \hs{l'} \subset \ns{l}$.
  The second condition is equivalent to
  $\dim \sum_{l'=0}^l \hs{l'} = \sum_{l'=0}^l \dim \hs{l'}$,
  i.e., to the directness of the sum.
  Therefore, the logical conjunction of both is equivalent to
  $\bigoplus_{l'=0}^l \hs{l'} \subset \ns{l}$.
  If the sum is direct,
  the dimension of the sum is equal to $2 + \sum_{l'=1}^l 2^{l'-1} = 2^l + 1$
  (due to $\dim \hs{l'} = \setsize{\hiset{l'}} = 2^{l'-1}$ for $l' > 0$ and
  $\dim \hs{l'} = 2$ for $l' = 0$),
  which is also the dimension of $\ns{l}$.
  The only subspace of $\ns{l}$ that has the same
  dimension as $\ns{l}$ is $\ns{l}$ itself,
  so we infer $\bigoplus_{l'=0}^l \hs{l'} = \ns{l}$.
\end{proof}
\begin{corollary}[univariate hierarchical splitting for hat functions]
  \label{cor:hierSplittingHatUV}
  The hierarchical splitting \eqref{eq:hierSplittingUV}
  holds for the hat function basis.
\end{corollary}
\begin{proof}
  The first condition of \cref{lemma:hierSplittingUV}
  is satisfied as piecewise linear splines of level~$l'$
  are also piecewise linear splines of higher levels $l \ge l'$.
  We can prove the linear independence for the second condition by induction
  over $l$:
  If a linear combination of $\bspl{l',i'}{1}$
  ($l' \le l$, $i' \in \hiset{l'}$)
  vanishes everywhere, then the coefficients of level $l$ must be zero,
  as otherwise the basis functions $\bspl{l,i'}{1}$ ($i' \in \hiset{l}$) would
  introduce kinks at $\gp{l,i'}$, which the zero function does not have.
  This means that we have a zero linear combination of $\bspl{l',i'}{1}$ for
  $l' \le l - 1$, $i' \in \hiset{l'}$,
  and by the induction hypothesis, the other coefficients also vanish.
\end{proof}

\subsection{Hierarchical Splitting in the Multivariate Case}
\label{sec:222hierMV}

Multivariate hierarchical subspaces are defined analogously
to the univariate case:
\begin{equation}
  \hs{\*l}
  \ceq \spn\{\basis{\*l,\*i} \mid \*i \in \hiset{\*l}\},\quad
  \hiset{\*l}
  \ceq \hiset{l_1} \times \dotsb \times \hiset{l_d},\quad
  \*l \in \natz^d.
\end{equation}
The univariate hierarchical splitting \eqref{eq:hierSplittingUV}
can now be generalized to
\begin{equation}
  \label{eq:hierSplittingMV}
  \ns{\*l}
  \overset{?}{=} \bigoplus_{\*l'=\*0}^\*l \hs{\*l'}.
\end{equation}
Again, this relation does not hold in general.
We use a multivariate counterpart of \thmref{lemma:hierSplittingUV}
to prove that \eqref{eq:hierSplittingMV} holds if
the corresponding univariate relation \eqref{eq:hierSplittingUV}
holds for all dimensions:

\begin{lemma}[multivariate hierarchical splitting characterization]
  \label{lemma:hierSplittingMV}
  \Cref{eq:hierSplittingMV} is equivalent to the satisfaction of
  both of the following conditions:
  \begin{itemize}
    \item
    The hierarchical subspaces $\hs{\*l'}$ ($\*l' \le \*l$)
    are subspaces of $\ns{\*l}$.
    
    \item
    The basis functions
    $\basis{\*l',\*i'}$ ($\*l' \le \*l$, $\*i' \in \hiset{\*l'}$)
    are linearly independent.
  \end{itemize}
\end{lemma}

\vspace*{0pt plus 0.3fill}

\begin{proof}
  If the sum is direct, then its dimension is given by
  \begin{equation}
    \hspace*{-5mm}
    \dim \sum_{\*l'=\*0}^\*l \hs{\*l'}
    = \sum_{l_1'=0}^{l_1} \dotsb \sum_{l_d'=0}^{l_d}
    \prod_{t=1}^d \dim \hs{l_t'}
    = \prod_{t=1}^d \sum_{l_t'=0}^{l_t} \dim \hs{l_t'}
    = \prod_{t=1}^d (2^{l_t} + 1)
    = \dim \ns{\*l}
    \hspace*{-5mm}
  \end{equation}
  using \eqref{eq:dimensionFG}.
  The rest is analogous to the proof of \cref{lemma:hierSplittingUV}.
\end{proof}

\vspace*{0pt plus 1fill}

\begin{proposition}[from univariate to multivariate splitting]
  \label{prop:splittingUVToMV}
  If univariate splitting \eqref{eq:hierSplittingUV}
  holds for every dimension,
  then the multivariate splitting \eqref{eq:hierSplittingMV} holds as well.
\end{proposition}

\vspace*{0pt plus 0.3fill}

\begin{proof}
  We check the two conditions of \cref{lemma:hierSplittingMV}
  given the two univariate conditions of \cref{lemma:hierSplittingUV}:
  \begin{enumerate}
    \item
    The hierarchical basis functions $\basis{\*l',\*i'}$
    of $\hs{\*l'}$ ($\*l' \le \*l$, $\*i' \in \hiset{\*l'}$)
    are tensor products of functions $\basis{l_t',i_t'}$.
    According to the first condition of \cref{lemma:hierSplittingUV},
    each $\basis{l_t',i_t'}$ can be written as a linear combination of
    the nodal basis $\basis{l_t,i_t}$ ($i_t = 0, \dotsc, 2^{l_t}$).
    We can expand the tensor product to a linear combination
    of tensor products of the univariate nodal basis functions.
    Therefore, $\basis{\*l',\*i'}$ is a linear combination of
    multivariate nodal functions, i.e., $\basis{\*l',\*i'} \in \ns{\*l}$.
    As this is true for all $\*i' \in \hiset{\*l'}$, we obtain
    $\hs{\*l'} \subset \ns{\*l}$.
    
    \item
    The linear independence of the hierarchical functions $\basis{\*l',\*i'}$
    ($\*l' \le \*l$, $\*i' \in \hiset{\*l'}$) can be shown completely
    analogously to the proof of
    \thmref{lemma:tensorProductLinearIndependence}.
  \end{enumerate}
  According to \cref{lemma:hierSplittingMV},
  the multivariate splitting \eqref{eq:hierSplittingMV} holds.
\end{proof}

\vspace*{0pt plus 1fill}

A direct consequence of \cref{prop:splittingUVToMV} is that
the hierarchical splitting holds for the hierarchical hat function basis.

\pagebreak

\begin{corollary}[multivariate hierarchical splitting for hat functions]
  \label{cor:hierSplittingHatMV}
  The multivariate hierarchical splitting \eqref{eq:hierSplittingMV}
  holds for the hat function basis.
\end{corollary}
\begin{proof}
  Follows directly by applying
  \thmref{cor:hierSplittingHatUV} to \cref{prop:splittingUVToMV}.
\end{proof}

\section{Sparse Grids}
\label{sec:23sparseGrids}

\minitoc{67mm}{4}

\noindent
The idea of sparse grids is to use the
hierarchical splitting \eqref{eq:hierSplittingMV}
to keep only the most important hierarchical subspaces,
omitting the remaining ones.
There are three main ``flavors'' of sparse grids:
regular, dimensionally adaptive, and spatially adaptive.

\subsection{Regular Sparse Grids}
\label{sec:231regularSG}

\paragraph{Hierarchical contributions}

To assess the importance of a subspace, we consider again the
interpolant $\fgintp{\*l} \in \ns{\*l}$ of a function $\objfun\colon \clint{\*0, \*1} \to \real$.
According to the splitting \eqref{eq:hierSplittingMV}, the interpolant can
be written as
\begin{equation}
  \label{eq:interpHierFullGrid}
  \fgintp{\*l}
  = \sum_{\*l'=\*0}^\*l \sum_{\*i' \in \hiset{\*l'}}
  \surplus{\*l',\*i'} \basis{\*l',\*i'},\quad
  \falarge{\*i = \*0, \dotsc, \*2^\*l}{\fgintp{\*l}(\gp{\*l,\*i}) = \objfun(\gp{\*l,\*i})}.
\end{equation}
The coefficients $\surplus{\*l',\*i'}$ with respect to the hierarchical basis
$\basis{\*l',\*i'}$ are the \term{hierarchical surpluses.}
When using the hat function basis $\bspl{\*l,\*i}{1}$,
one can prove the following representation
for the corresponding surpluses \multicite{Bungartz04Sparse,Garcke13Sparse}:
\begin{equation}
  \label{eq:surplusIntegral}
  \surplus{\*l',\*i'}
  = (-1)^d 2^{-\normone{\*l'+\*1}}
  \int_\*0^\*1 \bspl{\*l',\*i'}{1}(\*x)\,
  \partialderiv[2d]{\partialdiff x_1^2 \dotsm \partialdiff x_d^2}{\objfun}(\*x)
  \diff{}\*x,
\end{equation}
if $\*l \ge \*1$ and
$\objfun$ is twice continuously differentiable in every dimension simultaneously,
i.e.,
$\partialderiv[2d]{\partialdiff x_1^2 \dotsm \partialdiff x_d^2}{\objfun}$
exists and is continuous.%
\footnote{%
  Again, the notation implies that the integration domain is
  the unit hyper-cube $\clint{\*0, \*1} = \clint{0, 1}^d$.%
}\multiplefootnoteseparator%
\footnote{%
  The statement is even valid for functions in the Sobolev space
  $H_\mathrm{mix}^2(\clint{\*0, \*1})$ with dominating mixed derivative,
  as its proof mainly relies on integration by parts
  \multicite{Bungartz04Sparse,Garcke13Sparse}.%
}
Consequently, the contribution of the summand of level $\*l$
can be estimated by
\begin{equation}
  \label{eq:componentEstimation}
  \normLtwoscaled{
    \sum_{\*i' \in \hiset{\*l'}}
    \surplus{\*l',\*i'} \bspl{\*l',\*i'}{1}
  }
  \le 3^{-d} \cdot 2^{-2 \normone{\*l}} \cdot
  \normLtwoscaled{
    \partialderiv[2d]{\partialdiff x_1^2 \dotsm \partialdiff x_d^2}{\objfun}
  }
\end{equation}
for the hat function surpluses $\surplus{\*l',\*i'}$
\multicite{Bungartz04Sparse,Garcke13Sparse}.

\paragraph{Definition of regular sparse grids}

\begin{figure}
  \subcaptionbox{%
    Hierarchical splitting and subspace selection.
    The rectangles indicate the support of the
    bivariate hat basis functions.%
  }[85mm]{%
    \includegraphics{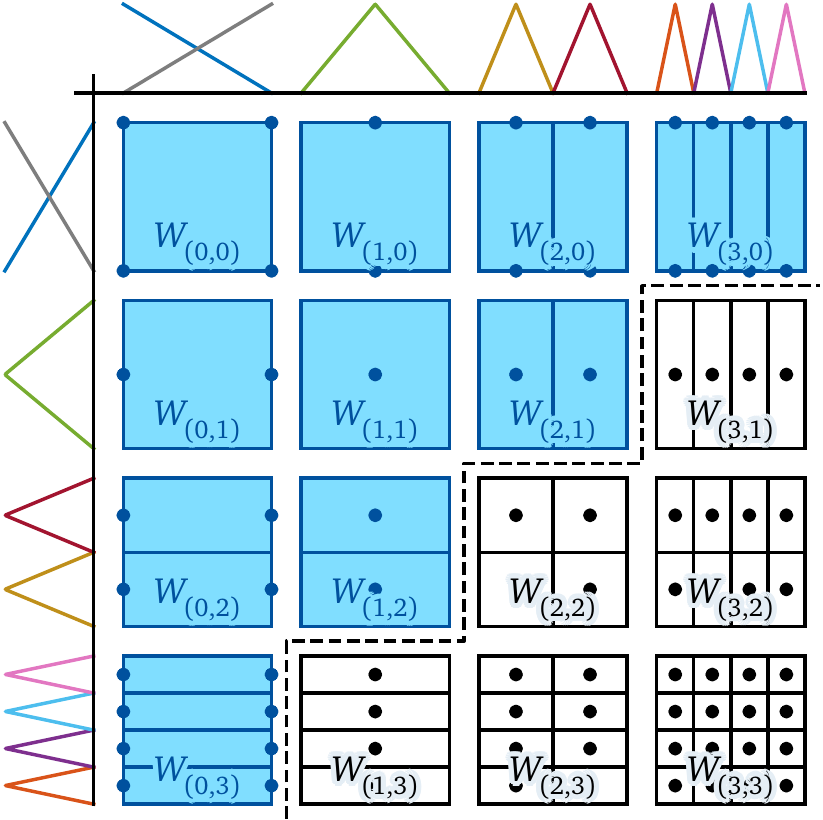}%
  }%
  \hfill%
  \begin{minipage}[b]{59mm}
    \subcaptionbox{%
      Full grid obtained by adding all subspaces of level $\*l \le n \cdot \*1$.%
    }[59mm]{%
      \includegraphics{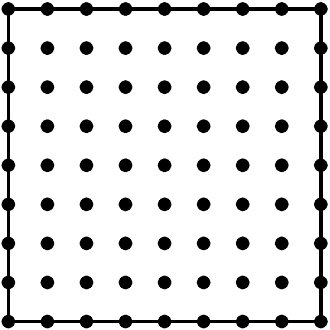}%
    }\\[4mm]%
    \subcaptionbox{%
      Regular sparse grid obtained by adding all subspaces
      whose level $\*l$ satisfies $\normone{\*l} \le n$
      \emph{\textcolor{mittelblau}{(blue)}.}%
    }[59mm]{%
      \includegraphics{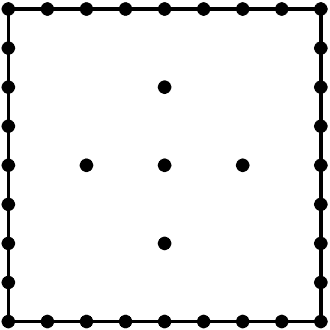}%
    }%
  \end{minipage}%
  \caption[%
    Regular two-dimensional sparse grid%
  ]{%
    Regular sparse grid of level $n = 3$ in two dimensions.%
  }%
  \label{fig:regularSG}%
\end{figure}

Equation \eqref{eq:componentEstimation} motivates to omit those summands
from the sum \eqref{eq:interpHierFullGrid} whose level sum $\normone{\*l}$
exceeds a certain value $n \in \natz$,
as their contribution can be neglected compared to the summands
with coarser level sums.
More formally, the selection of the relevant subspaces can be formulated as a
continuous knapsack problem~\cite{Bungartz04Sparse},
assuming homogeneous boundary conditions.
\usenotation{zzzzs}
This motivates
\begin{equation}
  \label{eq:regularSG}
  \regsgspace{n}{d}
  \ceq \bigoplus_{\normone{\*l} \le n} \hs{\*l},\qquad
  \regsgset{n}{d}
  \ceq \bigdotcup_{\normone{\*l} \le n}
  \{\gp{\*l,\*i} \mid \*i \in \hiset{\*l}\}
\end{equation}
as the definitions for the \term{regular sparse grid space} and
\term{regular sparse grid} of level $n$, respectively.
The functions $\regsgintp{n}{d}$ contained in
$\regsgspace{n}{d}$ have the form
\begin{equation}
  \label{eq:regularSGInterpolant}
  \regsgintp{n}{d}
  = \sum_{\normone{\*l} \le n} \sum_{\*i \in \hiset{\*l}}
  \surplus{\*l,\*i} \basis{\*l,\*i}.
\end{equation}
To better distinguish the different grids,
we call the grids corresponding to the nodal spaces \term{full grids.}
We generalize the definition to arbitrary bases $\basis{\*l,\*i}$,
although sparse grids have been motivated using the hat function
basis $\bspl{\*l,\*i}{1}$
(the estimate \eqref{eq:componentEstimation} does not hold anymore
in the general case).
\Cref{fig:regularSG} shows the construction of a
regular sparse grid in two dimensions.

\paragraph{Grid size and interpolation error}

One can prove that for homogeneous boundary conditions
$\restrictfcn{\objfun}{\bndrydomain{\clint{\*0,\*1}}} \equiv 0$,
the number of required inner grid points
($\gp{\*l,\*i} \in \regsgset{n}{d}$ where $\*l \ge \*1$)
grows like $\landauO{\ms{n}^{-1} (\log_2 \ms{n}^{-1})^{d-1}}$
\multicite{Bungartz04Sparse,Garcke13Sparse}, which is much less than
the corresponding number $\landauO{(\ms{n}^{-1})^d}$ in the full grid case
(see \eqref{eq:dimensionFG}).
The $\Ltwo$ error of the sparse grid interpolant
$\regsgintp{n}{d} \in \regsgspace{n}{d}$ using hat functions
(still assuming homogeneous boundary conditions) decays like
\begin{equation}
  \normLtwo{\objfun - \regsgintp{n}{d}}
  = \landauO{\ms{n}^2 (\log_2 \ms{n}^{-1})^{d-1}},
\end{equation}
which is only slightly worse than the full grid error by the factor of
$(\log_2 \ms{n}^{-1})^{d-1}$
\multicite{Bungartz04Sparse,Garcke13Sparse}.

\subsection{Dimensionally Adaptive Sparse Grids}
\label{sec:232dimensionallyAdaptiveSG}

The idea of dimensional adaptivity is to spend more grid
points along specific dimensions depending on the objective function.
Different criteria for the choice of dimensions exist,
for example the maximal absolute value of the linear hierarchical surpluses.
To incorporate dimensional adaptivity into sparse grids,
one has to generalize the symmetric
choice of subspaces in the definition of regular sparse grids
to allow asymmetric preferences.
Generally, function spaces~$\sgspace$ and grid sets $\sgset$
of \term{dimensionally adaptive sparse grids} have the form
\begin{equation}
  \label{eq:dimensionallyAdaptiveSG}
  \sgspace
  = \bigoplus_{\*l \in \levelset} \hs{\*l},\qquad
  \sgset
  = \bigdotcup_{\*l \in \levelset} \{\gp{\*l,\*i} \mid \*i \in \hiset{\*l}\},
\end{equation}
where $\levelset$ is a \term{downward closed} set, i.e.,
a finite subset $\levelset \subset \natz^d$
for which $\fafa{\*l \in \levelset}{\*l' \le \*l}{\*l' \in \levelset}$.
Regular sparse grids are a special case by setting
$\levelset = \{\*l \in \natz^d \mid \normone{\*l} \le n\}$.

\paragraph{Combination technique}

\begin{SCfigure}
  \includegraphics{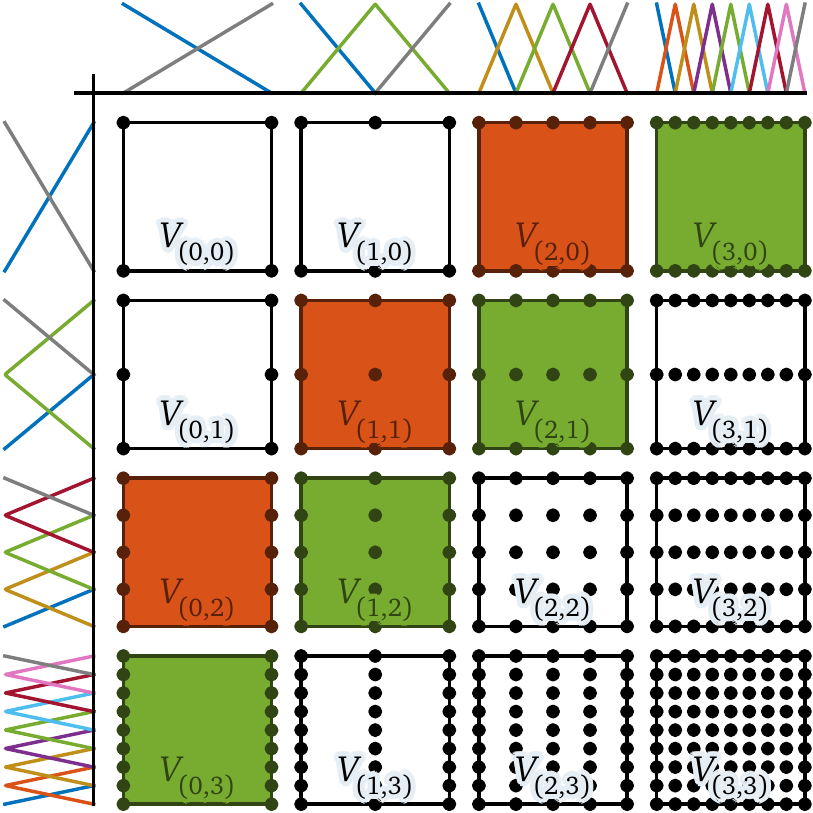}%
  \caption[%
    Sparse grid combination technique%
  ]{%
    The combination technique combines nodal subspaces in a weighted
    sum to form a regular sparse grid space of level $n = 3$ in two dimensions.
    The \textcolor{C1}{red subspaces} ($q = 1$ in \eqref{eq:combiTechnique})
    are subtracted from the sum of the
    \textcolor{C4}{green subspaces} ($q = 0$).%
  }%
  \label{fig:combinationTechnique}%
\end{SCfigure}

The key advantage of dimensionally adaptive sparse grids over
spatially adaptive approaches is the
so-called \term{combination technique.}
For regular sparse grids, one can show that the sparse grid interpolant
$\regsgintp{n}{d}$ can be written as
\begin{equation}
  \label{eq:combiTechnique}
  \regsgintp{n}{d}
  = \sum_{q=0}^{d-1} (-1)^q \binom{d-1}{q} \sum_{\normone{\*l} = n-q}
  \sum_{\*i=\*0}^{\*2^\*l} \interpcoeff{\*l,\*i} \basis{\*l,\*i},
\end{equation}
where the $\interpcoeff{\*l,\*i} \in \real$ ($\*i = \*0, \dotsc, \*2^\*l$)
are the interpolation coefficients on the full grid
$\fgset{\*l}$ of level~$\*l$, i.e.,
$\fa{\*i' = \*0, \dotsc, \*2^\*l}{%
  \sum_{\*i=\*0}^{\*2^\*l} \interpcoeff{\*l,\*i} \basis{\*l,\*i}(\gp{\*l,\*i'})
  = \objfun(\gp{\*l,\*i'})%
}$ \multicite{Smolyak63Quadrature,Zenger91Sparse}.
For general dimensionally adaptive sparse grids, a similar formula exists
\cite{Nobile16Adaptive}.
The combination formula \eqref{eq:combiTechnique} splits the
sparse grid interpolant into a weighted sum of full grid interpolants
(see \cref{fig:combinationTechnique}).
In applications, each grid can be processed in parallel,
drastically speeding up computations like the solution of \pdes{}
\cite{Heene18Massively}.
In addition, existing code working on nodal bases does not have to be
rewritten in terms of implementing hierarchical functions,
which means that the combination technique allows sparse grids to be employed
in existing software in a minimally invasive way.

\subsection{Spatially Adaptive Sparse Grids}
\label{sec:233spatiallyAdaptiveSG}

Dimensional adaptivity does not suffice to resolve local features of the
objective function.
Especially in some applications, it is crucial for the
interpolant to be highly accurate in specific regions of the domain.
For instance in optimization, it is not necessary to have a small global
interpolation error.
Instead, high accuracy near the optima is important.

\begin{figure}
  \subcaptionbox{%
    Hierarchical splitting and grid point selection.
    The rectangles indicate again the support of the
    bivariate hat basis functions.%
  }[85mm]{%
    \includegraphics{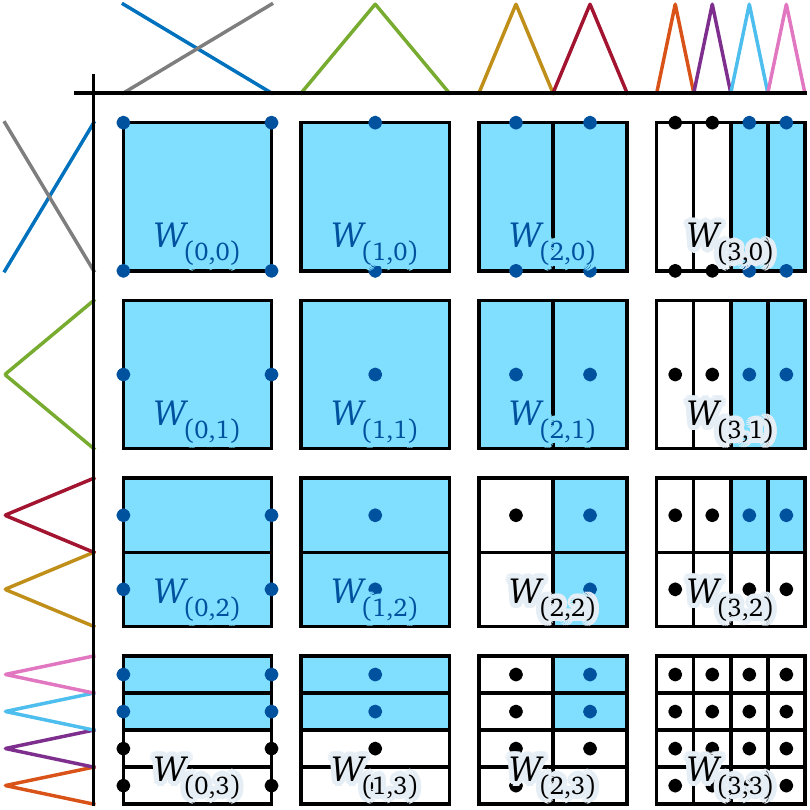}%
    \hspace*{1.411224mm}%
  }%
  \hfill%
  \subcaptionbox{%
    Resulting spatially adaptive sparse grid.%
  }[59mm]{%
    \includegraphics{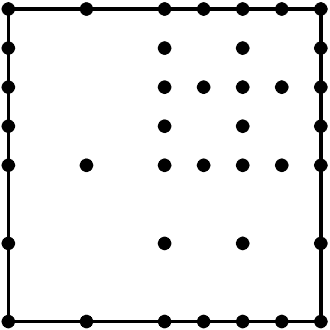}%
  }%
  \caption[%
    Construction of spatially adaptive sparse grids%
  ]{%
    Spatially adaptive sparse grid in two dimensions.
    More grid points were generated in the top right corner,
    which can help to resolve fine oscillations of the objective function.%
  }%
  \label{fig:spatiallyAdaptiveSG}%
\end{figure}

This can be achieved by \term{spatially adaptive sparse grids,}
on which this thesis focuses.
Generally, their function spaces $\sgspace$
and grid sets $\sgset$ have the form
\begin{equation}
  \label{eq:spatiallyAdaptiveSG}
  \sgspace
  = \spn\{\basis{\*l,\*i} \mid (\*l,\*i) \in \liset\},\qquad
  \sgset
  = \{\gp{\*l,\*i} \mid (\*l,\*i) \in \liset\},
\end{equation}
where $\liset$ is a finite set of level-index pairs $(\*l,\*i)$
with $\*l \in \natz^d$ and $\*i \in \hiset{\*l}$.
An example for a spatially adaptive sparse grid is shown in
\cref{fig:spatiallyAdaptiveSG}.

Algorithms for sparse grids often make specific assumptions about $\liset$.
If they are not met, then the algorithms do not produce the correct results.
For example when working with hat functions $\bspl{\*l,\*i}{1}$,
the grid should contain the hierarchical ancestors of every grid point.
Otherwise, the so-called unidirectional principle \cite{Balder94Adaptive},
which is used for instance to efficiently calculate
hierarchical surpluses, does not hold in general.
However, as we will see in \cref{chap:40algorithms},
the unidirectional principle cannot be applied
to B-splines of general degree, even if the hierarchical ancestors exist.
Hence, for most of our considerations, we will not restrict the
choice of $\liset$.

\section{Boundary Treatment}
\label{sec:24boundary}

\minitoc{79mm}{5}

\noindent
One issue of regular sparse grids $\regsgset{n}{d}$
is that the number of grid points still grows very fast
with the level $n$ and the dimensionality $d$ \cite{Pflueger10Spatially}.
This is mainly because the finest mesh size $\ms{n}$ on the
boundary of the domain $\clint{\*0, \*1}$ is finer than
the finest mesh size $\ms{n-d+1}$ that can be found in the interior.
If we define $\interiorregsgset{n}{d}$ as the set of
interior grid points in $\regsgset{n}{d}$,%
\footnote{%
  Note that in the literature (e.g., \cite{Pflueger10Spatially}),
  the regular sparse grid space of level $n$ without boundary points is often
  defined via $\normone{\*l} \le n + d - 1$ to ensure that the finest mesh size
  is given by $\ms{n}$.
  In our notation, this corresponds to $\interiorregsgset{n+d-1}{d}$.%
}
i.e.,
\begin{equation}
  \interiorregsgset{n}{d}
  \ceq \regsgset{n}{d} \cap \opint{\*0, \*1}
  = \{\gp{\*l,\*i} \in \regsgset{n}{d} \mid \*l \ge \*1\},
\end{equation}
then the following relation about the number of grid points
of $\regsgset{n}{d}$ can be shown:

\begin{lemma}[number of regular sparse grid points]
  \label{lemma:numberOfGridPointsBoundary}
  \setlength{\abovedisplayskip}{0pt}%
  \begin{equation}
    \setsize{\regsgset{n}{d}}
    = \sum_{q=0}^d 2^q \binom{d}{q} \setsize{\interiorregsgset{n}{d-q}}
  \end{equation}
\end{lemma}
\begin{proof}
  See \cite{Bungartz04Sparse}.
\end{proof}
Here, we define zero-dimensional grids to contain exactly one grid point
such that $\setsize{\interiorregsgset{n}{0}} = 1$.
The number of interior grid points can be calculated as follows:
\begin{lemma}[number of interior regular sparse grid points]
  \label{lemma:numberOfGridPointsInterior}
  \setlength{\abovedisplayskip}{0pt}%
  \begin{equation}
    \setsize{\interiorregsgset{n}{d}}
    = \sum_{q=0}^{n-d} 2^q \binom{d-1+q}{d-1}
  \end{equation}
\end{lemma}
\begin{proof}
  See \cite{Bungartz04Sparse}.
\end{proof}

\begin{figure}
  \raisebox{-0.5\height}{\includegraphics{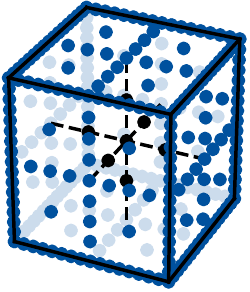}}%
  \raisebox{-0.5\height-0.5mm}{$\;\;=\;\;$}%
  \raisebox{-0.5\height}{\includegraphics{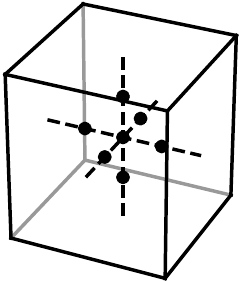}}%
  \raisebox{-0.5\height}{$\;\;\dotcup\;\;$}%
  \raisebox{-0.5\height}{\includegraphics{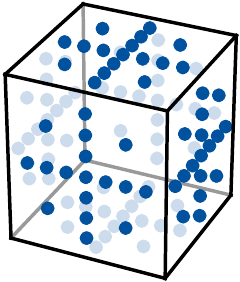}}%
  \raisebox{-0.5\height}{$\;\;\dotcup\;\;$}%
  \raisebox{-0.5\height}{\includegraphics{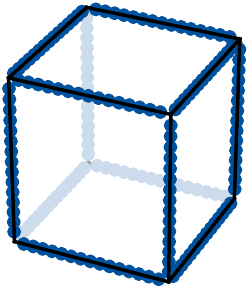}}%
  \raisebox{-0.5\height}{$\;\;\dotcup\;\;$}%
  \raisebox{-0.5\height}{\includegraphics{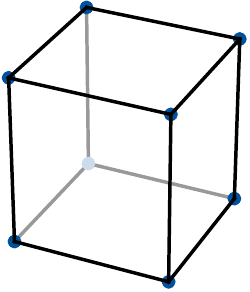}}%
  \caption[%
    Decomposition of a sparse grid into lower-dimensional sparse sub-grids%
  ]{%
    Decomposition of the three-dimensional sparse grid $\regsgset{n}{d}$
    ($n = 4$, $d = 3$) into lower-dimensional sparse sub-grids.
    The main axes (axis-parallel lines through $0.5 \cdot \*1$, \emph{dashed})
    serve as a visual aid.%
  }%
  \label{fig:sgDecompose}%
\end{figure}

Intuitively, \cref{lemma:numberOfGridPointsBoundary} splits the sparse grid
$\regsgset{n}{d}$ into lower-dimensional sparse grids
$\interiorregsgset{n}{d-q}$ of the same level, but without boundary points.
The factor $2^q \binom{d}{q}$ is the number of $(d-q)$-dimensional faces
of the $d$-dimensional unit hyper-cube.
In the three-dimensional example of \cref{fig:sgDecompose},
the unit cube $\clint{0, 1}^3$ decomposes into
\begin{itemize}
  \item
  $2^0 \binom{3}{0} = 1$ interior cube $\opint{0, 1}^3$,
  
  \item
  $2^1 \binom{3}{1} = 6$ sides (two-dimensional faces)
  like $\opint{0, 1}^2 \times \{0\}$,
  
  \item
  $2^2 \binom{3}{2} = 12$ edges (one-dimensional faces)
  like $\opint{0, 1} \times \{(0, 0)\}$, and
  
  \item
  $2^3 \binom{3}{3} = 8$ corners (zero-dimensional faces)
  like $(0, 0, 0)$.
\end{itemize}
On each of these $(d-q)$-dimensional faces,
the sparse grid $\regsgset{n}{d}$ contains
the interior of a sparse grid of level $n$ and dimensionality $d - q$,
the size of which grows like $\landauO{2^n n^{d-q-1}}$.
As the number of boundary faces increases exponentially
with the dimensionality $d$,
the size of $\regsgset{n}{d}$ quickly exhausts the available
computational memory.
To deal with this issue, there are mainly two solutions,
which are described below.

\subsection{Sparse Grids with Coarser Boundaries}
\label{sec:241coarseBoundary}

\paragraph{Inserting boundary points at higher levels}

The first solution is to insert the boundary level functions and grid points
at a higher level than at level zero.
A popular choice is the insertion at level one, which corresponds to
\begin{equation}
  \label{eq:sparseGridB1}
  \coarseregsgset{n}{d}{1}
  \ceq \bigdotcup_{\*l \in \coarselevelset{n}{d}{1}}
  \{\gp{\*l,\*i} \mid \*i \in \hiset{\*l}\},\quad
  \coarselevelset{n}{d}{1}
  \ceq \{\*l \in \natz^d \mid \normone{\vecmax(\*l, \*1)} \le n\},
\end{equation}
where $\vec{\max}$ is to be read coordinate-wise as usual.
This choice is equivalent to treating zero-level components as level one
in the subspace selection.
This ensures that the finest mesh sizes in the interior of
$\clint{\*0, \*1}$ and on its boundary coincide to be $\ms{n-d+1}$,
which reduces the number of grid points on the boundary significantly.

Another solution that can be found in the literature of sparse grids with
hat functions \cite{Baar15Gradient}
is to start with the ``constant one'' function on level zero with
corresponding grid point $0.5$,
then employ the two boundary functions and points on level one,
and finally proceed as usual for the higher levels $\ge 2$.
Apart from a constant shift of the resulting sparse grid levels,
this is equivalent to inserting the boundary functions and points at level two.
This solution leads to even less grid points than the previous approach,
as now the mesh size is finer in the interior of the domain than on the
boundary.
However, for very high dimensionalities this might still lead to
computationally infeasible sparse grids.

\paragraph{Regular sparse grids with coarse boundary}

\begin{figure}
  \subcaptionbox{%
    $d = 2$, $b = 0$%
  }[35mm]{%
    \includegraphics{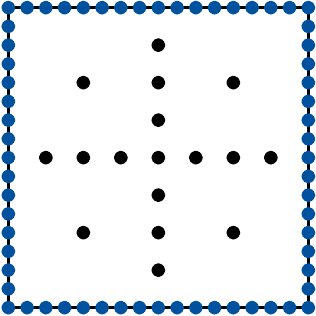}%
  }%
  \hfill%
  \subcaptionbox{%
    $d = 2$, $b = 1$%
  }[35mm]{%
    \includegraphics{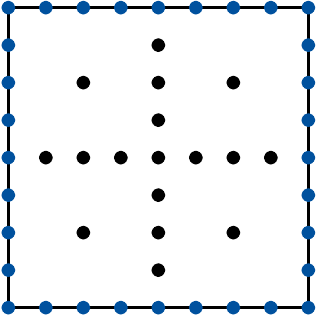}%
  }%
  \hfill%
  \subcaptionbox{%
    $d = 2$, $b = 2$%
  }[35mm]{%
    \includegraphics{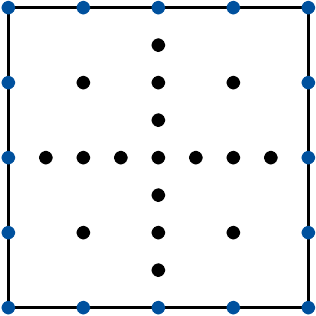}%
  }%
  \hfill%
  \subcaptionbox{%
    $d = 2$, $b = 3$%
  }[35mm]{%
    \includegraphics{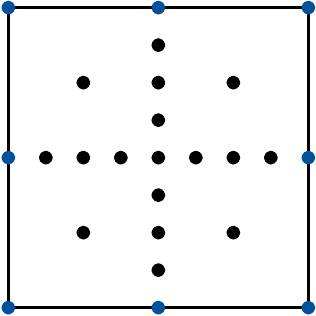}%
  }\\[2mm]%
  \subcaptionbox{%
    $d = 3$, $b = 0$%
  }[35mm]{%
    \includegraphics{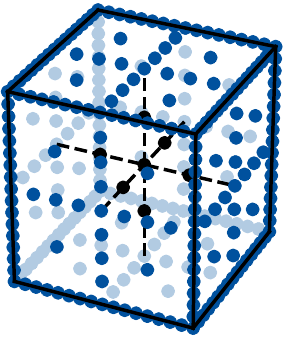}%
  }%
  \hfill%
  \subcaptionbox{%
    $d = 3$, $b = 1$%
  }[35mm]{%
    \includegraphics{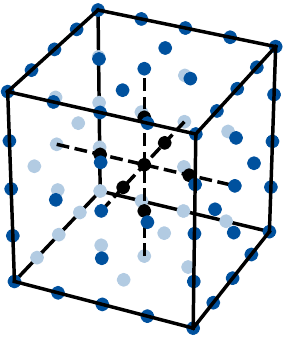}%
  }%
  \hfill%
  \subcaptionbox{%
    $d = 3$, $b = 2$%
  }[35mm]{%
    \includegraphics{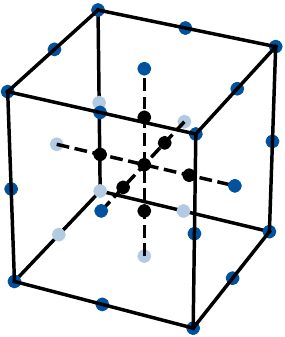}%
  }%
  \hfill%
  \subcaptionbox{%
    $d = 3$, $b = 3$%
  }[35mm]{%
    \includegraphics{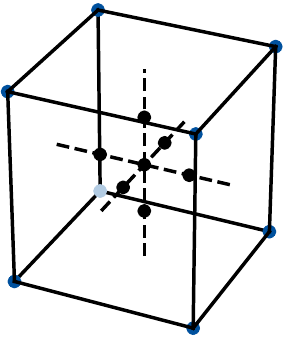}%
  }%
  \caption[%
    Comparison of regular sparse grids with coarse boundary%
  ]{%
    Sparse grids $\coarseregsgset{n}{d}{b}$ of level $n = 4$
    in two and three dimensions for different values of the
    boundary parameter $b$.
    For constant $d$ and $n$,
    the points in the interior of $\clint{\*0, \*1}$
    (black) are the same,
    while the points on the boundary of $\clint{\*0, \*1}$
    \emph{\textcolor{mittelblau}{(blue)}} become coarser
    for increasing values of $b$.
    The main axes (axis-parallel lines through $0.5 \cdot \*1$, \emph{dashed})
    serve as a visual aid.%
  }%
  \label{fig:coarseBoundary}%
\end{figure}

\usenotation{zzzzsb}
We generalize these two solutions to the definition of a
sparse grid $\coarseregsgset{n}{d}{b}$ that is equivalent to inserting
the boundary functions and points at an arbitrary level $b \in \nat$:
\begin{definition}[regular sparse grid with coarse boundary]
  \label{def:coarseBoundary}
  The regular sparse grid of level $n \in \nat$,
  dimensionality $d \le n$, and boundary parameter $b \in \nat$ is defined as
  \begin{subequations}
    \begin{align}
      \label{eq:coarseBoundary1}
      \coarseregsgset{n}{d}{b}
      &\ceq \bigdotcup_{\*l \in \coarselevelset{n}{d}{b}}
      \{\gp{\*l,\*i} \mid \*i \in \hiset{\*l}\},\\
      \label{eq:coarseBoundary2}
      \begin{split}
        \coarselevelset{n}{d}{b}
        &\ceq \{\*l \in \nat^d \mid \normone{\*l} \le n\}\\
        &\qquad {} \dotcup \paren*{
          \{\*l \in \natz^d \setminus \nat^d \mid
          \normone{\vecmax(\*l, \*1)} \le n-b+1\} \cup \{\*0\}
        }.
      \end{split}
    \end{align}
  \end{subequations}
  For convenience, we define
  $\coarseregsgset{n}{d}{0} \ceq \regsgset{n}{d}$.
\end{definition}
The definition is motivated by partitioning the levels $\*l \in \natz^d$
into interior levels ($\*l \in \nat^d$)
and boundary levels ($\*l \in \natz^d \setminus \nat^d$).
By including the levels of the interior grid $\interiorregsgset{n}{d}$,
the mesh size in the interior is the same as before ($\ms{n-d+1}$).
Like in \eqref{eq:sparseGridB1}, we treat boundary levels as level one,
but we subtract $b - 1$ from the upper bound to ensure the correct
mesh size $\ms{n-d-b+2}$ on the boundary.
We append $\*0$ to the level set to ensure that at least the $2^d$ corner
points are included in the resulting sparse grid.
Note that this definition is consistent with \eqref{eq:sparseGridB1} as
$\coarselevelset{n}{d}{b}
= \{\*l \in \natz^d \mid \normone{\vecmax(\*l, \*1)} \le n\}$
for $b = 1$.
Examples of~$\coarseregsgset{n}{d}{b}$ are shown
in \cref{fig:coarseBoundary}.
The flip book animation in the bottom right corner of the
odd-numbered pages of this thesis
visualizes $\coarseregsgset{n}{d}{b}$ for $n = 4$, $d = 3$, and $b = 1$.

The number of grid points of $\coarseregsgset{n}{d}{b}$
can be calculated as follows:
\begin{restatable}[number of regular sparse grid points with coarse boundary]{%
  proposition%
}{%
  propGridSizeCoarseBoundary%
}
  \label{prop:gridSizeCoarseBoundary}
  \setlength{\abovedisplayskip}{0pt}%
  \begin{equation}
    \setsize{\coarseregsgset{n}{d}{b}}
    = \setsize{\interiorregsgset{n}{d}} +
    \sum_{q=1}^d 2^q \binom{d}{q}
    \setsize{\interiorregsgset{n-q-b+1}{d-q}},\quad
    b \in \nat
  \end{equation}
\end{restatable}
\begin{proof}
  See \cref{sec:a111proofGridSizeCoarseBoundary}.
\end{proof}
As can be seen in \cref{tbl:coarseBoundary3D} for three dimensions and
in \cref{tbl:coarseBoundary10D} for ten dimensions,
the number of grid points decreases drastically for increasing values
of $b$, especially when compared with
$\regsgset{n}{d} = \coarseregsgset{n}{d}{0}$.

\begin{table}
  \newcommand*{\myheader}{$\setsize{\interiorregsgset{n}{d}}$}%
  \setnumberoftableheaderrows{2}%
  \begin{tabular}{%
    >{\kern\tabcolsep}=l<{\kern7mm}+r<{\kern7mm}+r+r+r+r+r+r<{\kern\tabcolsep}%
  }
    \toprulec
    \headerrow
    &&\multicolumn{6}{c}{%
      $\setsize{\coarseregsgset{n}{d}{b}}/\setsize{\interiorregsgset{n}{d}}$%
    }\\
    \headerrow
    &          \myheader&     $b = 0$&    $b = 1$&    $b = 2$&    $b = 3$&    $b = 4$&    $b = 5$\\
    \midrulec
    $n = 3$&     \num{1}& \num{123.0}& \num{27.0}& \num{9.00}& \num{9.00}& \num{9.00}& \num{9.00}\\
    $n = 4$&     \num{7}&  \num{42.4}& \num{11.6}& \num{4.71}& \num{2.14}& \num{2.14}& \num{2.14}\\
    $n = 5$&    \num{31}&  \num{22.7}& \num{7.3}&  \num{3.39}& \num{1.84}& \num{1.26}& \num{1.26}\\
    $n = 6$&   \num{111}&  \num{14.9}& \num{5.3}&  \num{2.75}& \num{1.67}& \num{1.23}& \num{1.07}\\
    $n = 7$&   \num{351}&  \num{10.9}& \num{4.3}&  \num{2.37}& \num{1.55}& \num{1.21}& \num{1.07}\\
    $n = 8$&  \num{1023}&   \num{8.5}& \num{3.6}&  \num{2.13}& \num{1.47}& \num{1.19}& \num{1.07}\\
    $n = 9$&  \num{2815}&   \num{7.0}& \num{3.2}&  \num{1.96}& \num{1.41}& \num{1.17}& \num{1.07}\\
    $n = 10$& \num{7423}&   \num{6.0}& \num{2.9}&  \num{1.83}& \num{1.36}& \num{1.16}& \num{1.06}\\
    \bottomrulec
  \end{tabular}
  \caption[%
    Comparison of regular sparse grid sizes with coarse boundary
    ($d = 3$)%
  ]{%
    For $d = 3$:
    Grid size of the interior grid
    \vspace{-0.33em}%
    $\interiorregsgset{n}{d}$ \emph{(second column)}
    and ratios
    $\setsize{\coarseregsgset{n}{d}{b}}/\setsize{\interiorregsgset{n}{d}}$
    \emph{(beginning with the third column)} of the sizes of
    the grid $\coarseregsgset{n}{d}{b}$ with boundary points
    to the size of the interior grid of the same level.
    The table starts with the first level $n = 3$ for which
    the interior grid $\interiorregsgset{n}{d}$ is not empty.%
  }%
  \label{tbl:coarseBoundary3D}%
\end{table}

\begin{table}
  \newcommand*{\myheader}{$\setsize{\interiorregsgset{n}{d}}$}%
  \setnumberoftableheaderrows{2}%
  \begin{tabular}{%
    >{\kern\tabcolsep}=l<{\kern7mm}+r<{\kern7mm}+r+r+r+r+r+r<{\kern\tabcolsep}%
  }
    \toprulec
    \headerrow
    &&\multicolumn{6}{c}{%
      $\setsize{\coarseregsgset{n}{d}{b}}/\setsize{\interiorregsgset{n}{d}}$%
    }\\
    \headerrow
    &             \myheader&     $b = 0$&     $b = 1$&    $b = 2$&      $b = 3$&      $b = 4$&      $b = 5$\\
    \midrulec
    $n = 10$&       \num{1}& \num{3.3e8}& \num{59049}& \num{1025}& \num{1025.0}& \num{1025.0}& \num{1025.0}\\
    $n = 11$&      \num{21}& \num{4.3e7}& \num{21558}& \num{2813}&   \num{49.8}&   \num{49.8}&   \num{49.8}\\
    $n = 12$&     \num{241}& \num{1.0e7}& \num{10046}& \num{1879}&  \num{246.0}&    \num{5.2}&    \num{5.2}\\
    $n = 13$&    \num{2001}& \num{3.4e6}&  \num{5407}& \num{1211}&  \num{227.2}&   \num{30.5}&    \num{1.5}\\
    $n = 14$&   \num{13441}& \num{1.3e6}&  \num{3213}&  \num{806}&  \num{181.1}&   \num{34.7}&    \num{5.4}\\
    $n = 15$&   \num{77505}& \num{6.2e5}&  \num{2054}&  \num{558}&  \num{140.6}&   \num{32.2}&    \num{6.8}\\
    $n = 16$&  \num{397825}& \num{3.1e5}&  \num{1390}&  \num{401}&  \num{109.5}&   \num{28.2}&    \num{7.1}\\
    $n = 17$& \num{1862145}& \num{1.7e5}&   \num{984}&  \num{298}&   \num{86.5}&   \num{24.2}&    \num{6.8}\\
    \bottomrulec
  \end{tabular}
  \caption[%
    Comparison of regular sparse grid sizes with coarse boundary
    ($d = 10$)%
  ]{%
    For $d = 10$:
    Grid size of the interior grid
    \vspace{-0.33em}%
    $\interiorregsgset{n}{d}$ \emph{(second column)}
    and ratios
    $\setsize{\coarseregsgset{n}{d}{b}}/\setsize{\interiorregsgset{n}{d}}$
    \emph{(beginning with the third column)} of the sizes of
    the grid $\coarseregsgset{n}{d}{b}$ with boundary points
    to the size of the interior grid of the same level.
    The table starts with the first level $n = 10$ for which
    the interior grid $\interiorregsgset{n}{d}$ is not empty.%
  }%
  \label{tbl:coarseBoundary10D}%
\end{table}

\Cref{alg:coarseBoundary} shows how to generate the necessary set of
hierarchical levels.
Its correctness can be formally proven with the following invariant:
\begin{restatable}[invariant of SG generation with coarse boundary]{%
  proposition%
}{%
  propInvariantCoarseBoundary%
}
  \label{prop:invariantCoarseBoundary}
  After iteration $t$ of \cref{alg:coarseBoundary}
  ($t = 1, \dotsc, d$), it holds
  \begin{equation}
    \label{eq:coarseInvariant}
    \begin{split}
      \levelset^{(t)}
      &= \{\*l \in \nat^t \mid \normone{\*l} \le n - d + t\}\\
      &\hphantom{{}={}} {} \dotcup \paren*{
        \{\*l \in \natz^t \setminus \nat^t \mid
        \normone{\vecmax(\*l, \*1)} \le n-d+t-b+1\} \cup \{\*0\}
      }.
    \end{split}
  \end{equation}
\end{restatable}
\begin{proof}
  See \cref{sec:a112proofInvariantCoarseBoundary}.
\end{proof}
\begin{shortcorollary}[correctness of SG generation with coarse boundary]
  \label{cor:algCoarseBoundaryCorrectness}
  \Cref{alg:coarseBoundary} is correct.
\end{shortcorollary}
\begin{proof}
  Follows immediately from \cref{prop:invariantCoarseBoundary}
  by setting $t = d$,
  as then \eqref{eq:coarseInvariant} becomes
  \eqref{eq:coarseBoundary2} from \thmref{def:coarseBoundary}.
\end{proof}

\begin{algorithm}
  \begin{algorithmic}[1]
    \Function{$\coarselevelset{n}{d}{b} = \texttt{computeSGCoarseBoundary}$}{%
      $n$, $d$, $b$%
    }
      \State{$\levelset^{(1)} \gets \{0, 1, \dotsc, n - d + 1\}$}
      \Comment{one-dimensional grid}%
      \label{line:algCoarseBoundary6}
      \For{$t = 2, \dotsc, d$}
        \State{$\levelset^{(t)} \gets \emptyset$}
        \Comment{$t$-dimensional grid}%
        \For{$\*l \in \levelset^{(t-1)}$}
          \If{%
            $\normone{\vecmax(\*l, \*1)} \le n - d + t - b$ or
            $\*l = \*0$%
          }%
          \label{line:algCoarseBoundary1}
            \State{$\levelset^{(t)} \gets \levelset^{(t)} \cup \{(\*l, 0)\}$}
            \Comment{%
              add corners (with $(\*l, 0) \ceq (l_1, \dotsc, l_{t-1}, 0)$)%
            }%
            \label{line:algCoarseBoundary5}
          \EndIf{}
          \If{$\*l \in \nat^{t-1}$}
            \State{$l^\ast \gets n - d + t - \normone{\*l}$}%
            \Comment{add interior points}%
            \label{line:algCoarseBoundary2}
          \Else{}
            \State{%
              $l^\ast \gets n - d + t - b + 1 -
              \normone{\vecmax(\*l, \*1)}$%
            }%
            \Comment{add boundary points}%
            \label{line:algCoarseBoundary3}
          \EndIf{}
          \State{%
            $\levelset^{(t)} \gets \levelset^{(t)} \cup
            \{(\*l, l_t) \mid l_t = 1, \dotsc, l^\ast\}$%
          }
          \Comment{%
            with $(\*l, l_t) \ceq (l_1, \dotsc, l_{t-1}, l_t)$%
          }%
          \label{line:algCoarseBoundary4}
        \EndFor{}
      \EndFor{}
      \State{$\coarselevelset{n}{d}{b} \gets \levelset^{(d)}$}
    \EndFunction{}
  \end{algorithmic}
  \caption[%
    Generation of regular sparse grids with coarse boundary%
  ]{%
    Generation of the sparse grid $\coarseregsgset{n}{d}{b}$
    with coarse boundary.
    Inputs are the level $n \in \nat$, the dimensionality $d \le n$, and
    the boundary parameter $b \in \nat$.
    Output is the level set $\coarselevelset{n}{d}{b}$
    that corresponds to $\coarseregsgset{n}{d}{b}$.%
  }%
  \label{alg:coarseBoundary}%
\end{algorithm}

\paragraph{Hierarchization and other algorithms}

An important implication of the regular sparse grids
$\coarseregsgset{n}{d}{b}$ as defined in \cref{def:coarseBoundary}
is that, in general,
the unidirectional principle cannot be directly applied anymore.
For example, this is relevant when calculating hierarchical surpluses
for the hat function basis.
As we mostly deal with B-splines, for which the unidirectional
principle cannot be applied even on regular sparse grids,
this issue is not important in the scope of this thesis.

However, it is possible to calculate the hierarchical surpluses
of hat functions on $\coarseregsgset{n}{d}{b}$ in a three-step algorithm.
First, we compute the surpluses of the boundary grid
$\coarseregsgset{n}{d}{b} \setminus \interiorregsgset{n}{d}$.
Second, we subtract the values of the resulting ``boundary interpolant'' at
the inner grid points
$\interiorregsgset{n}{d}$.
Third, we calculate the surpluses of the inner grid points
as usual with the unidirectional principle.
As the corresponding ``inner interpolant'' vanishes
on the boundary, this does not influence the interpolated values in the
first step.

\subsection{Sparse Grids Without Boundary Points and Modified Bases}
\label{sec:242modified}

\paragraph{Omitting boundary points}

The second solution to reduce the number of grid points on the boundary
is to omit the boundary points and the basis functions altogether.
For the hat function basis $\bspl{\*l,\*i}{1}$,
this is a feasible option if the objective
function $\objfun\colon \clint{\*0, \*1} \to \real$
satisfies homogeneous boundary conditions
$\restrictfcn{\objfun}{\bndrydomain{\clint{\*0, \*1}}} \equiv 0$,
as $\bspl{\*l,\*i}{1}$ vanishes on the boundary if and only if
$\*l \ge \*1$, i.e., if the basis function corresponds to an inner grid point.
Consequently, the surpluses corresponding to boundary points vanish
for a grid with boundary points and homogeneous boundary conditions,
implying that these points can be removed from the grid.

\paragraph{Modified linear basis}

\begin{SCfigure}
  \includegraphics{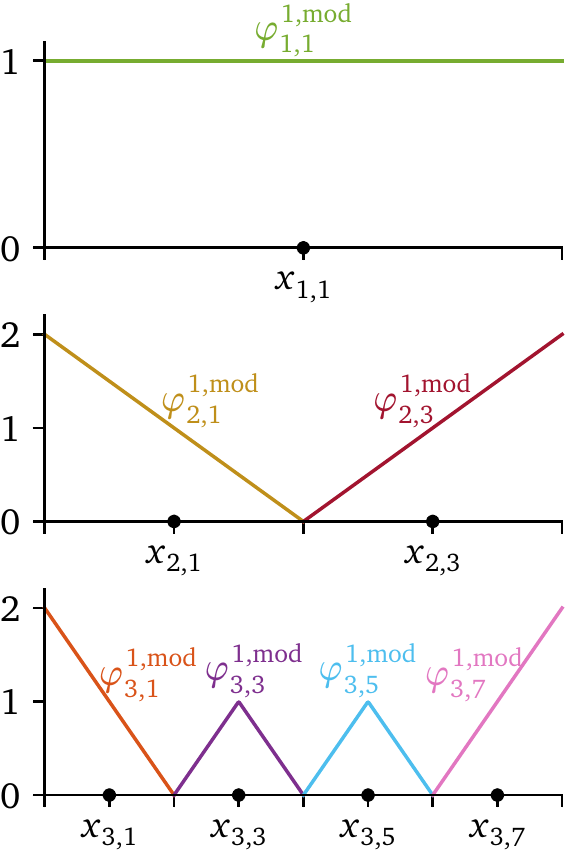}%
  \caption[%
    Modified hierarchical hat functions%
  ]{%
    Modified hierarchical hat functions $\bspl[\modified]{l',i'}{1}$
    ($l' \le l$, $i' \in \hiset{l'}$) up to level $l = 3$.%
  }%
  \label{fig:modifiedHat}%
\end{SCfigure}

Of course, this approach is not viable for functions with non-zero
boundary values or general hierarchical bases,
making it necessary to change the basis.
For hat functions, Pflüger modified the leftmost and rightmost
univariate basis function of each level (with indices $i = 1$ and
$i = 2^l - 1$ respectively) such that the modified functions
extrapolate the inner values linearly towards the boundary
\cite{Pflueger10Spatially}.
The basis function on level one is replaced by the
``constant one'' function.
All other basis functions remain unchanged.
\usenotation{zzzzmod}
The resulting \term{modified hat functions}
$\bspl[\modified]{l,i}{1}\colon \clint{0, 1} \to \real$
are shown in \cref{fig:modifiedHat} and defined as follows:
\begin{equation}
  \bspl[\modified]{l,i}{1}(x)
  \ceq
  \begin{cases}
    1,&
    l = 1,\quad i = 1,\\
    \max(2 - \tfrac{x}{\ms{l}}, 0),&
    l \ge 2,\quad i = 1,\\
    \bspl{l,i}{1}(x),&
    l \ge 2,\quad i \in \hiset{l} \setminus \{1, 2^l - 1\},\\
    \bspl[\modified]{l,1}{1}(1 - x),&
    l \ge 2,\quad i = 2^l - 1.
  \end{cases}
\end{equation}
The modified linear basis provides ``reasonable'' boundary values
without the need to insert basis functions and grid points on the boundary.
For other bases such as B-splines, similar modifications are possible,
which we will discuss when we introduce the corresponding unmodified functions
(see \cref{chap:30BSplines,chap:40algorithms}).

\cleardoublepage

  \setdictum{%
  B-splines are not enough!%
}{%
  In a talk at the 2017 SIAM Conference on\\
  Computational Science and Engineering%
}

\chapter{Hierarchical B-Splines}
\label{chap:30BSplines}

\initial{-0.05em}{P}{iecewise linear ``hat functions''} $\bspl{l,i}{1}$,
which served in the previous chapter as the motivation
to define sparse grids for arbitrary tensor product basis functions,
are not continuously differentiable.
This has two implications.
First, the approximation order of hat functions
is lower than the order of other basis function types
such as higher-degree splines \cite{Sickel11Spline}
or the piecewise polynomial basis by Bungartz \cite{Bungartz98Finite}.
Second, we cannot compute globally continuous gradients of the
interpolant of a smooth objective function
if we use non-smooth basis functions%
\footnote{%
  \dots{} which include the hat function basis as well as
  the piecewise polynomials by Bungartz.%
}.
However, the availability of continuous gradients is
crucial in gradient-based optimization,
which is highly important in
the scope of simulation technology (see \cref{chap:10introduction}) and
which we target in this thesis.
In this chapter, we define the hierarchical and
higher-order B-spline basis
as a generalization of the well-known hat functions
to obtain both higher-order approximations
and continuous derivatives.

The first to study B-splines was Schoenberg in 1946
\cite{Schoenberg46Contributions},
but he claimed that they had already been known to Laplace
\cite{Boor76Splines}.
It was also Schoenberg who coined the term ``B-splines,''
which is short for ``basis splines'' \cite{Schoenberg67Spline}.
De~Boor pioneered B-splines, developed basic algorithms, and
proved fundamental theoretical results \cite{Boor72Calculating}.
Research and industry recognized the possibilities of B-splines when
the \fem emerged in the 1960s.
The \fem remains one of the driving forces behind the research of
B-splines \cite{Hoellig03Finite} as the recent rise of \iga shows
\multicite{Cottrell09Isogeometric,Hoellig12Finite}.
Researchers have also applied B-splines to
geometric modeling with \nurbs
\multicite{Cohen01Geometric,Hoellig13Approximation},
financial mathematics \cite{Pflueger10Spatially},
molecular and atomic physics
\multicite{Bachau01Applications,McCurdy04Implementation},
and numerous other scientific and industrial areas
\multicite{Valentin12Spline,Martin17WEB}.
Theoretical and practical aspects of B-splines on sparse grids
have also been studied before
\multicite{%
  Pandey08Regression,%
  Pflueger10Spatially,%
  Sickel11Spline,%
  Valentin16Hierarchical%
}.

In this chapter, we define hierarchical B-splines on sparse grids.
The chapter is divided into two sections:
First, we define hierarchical B-splines for both
uniform and non-uniform knot sequences in \cref{sec:31standardBSplines}.
Second, we learn in \cref{sec:32notAKnot} that the boundary behavior
of the standard uniform B-spline basis is problematic.
Incorporating not-a-knot boundary conditions into the B-spline basis
mitigates the problems caused by the boundary behavior.

\Cref{sec:311uniform} is a repetition of the definition
of nodal B-splines \multicite{Hoellig03Finite,Hoellig13Approximation} and
hierarchical B-splines \multicite{Pflueger10Spatially,Valentin14Hierarchische}.
Original contributions of the thesis in this chapter are the proof of
the linear independence of hierarchical B-splines in
\cref{sec:312proofHierarchicalSplitting}
(improved version of \cite{Valentin14Hierarchische},
published in \cite{Valentin16Hierarchical}),
the modified hierarchical Clenshaw-Curtis B-splines in
\cref{sec:314nonUniform}, and
the hierarchical not-a-knot B-spline basis in \cref{sec:32notAKnot}.

\section{Uniform and Non-Uniform Hierarchical B-Splines}
\label{sec:31standardBSplines}

\minitoc{83mm}{7}

\noindent
In this section, we mainly follow the presentation of
\multicite{Pflueger10Spatially,Valentin14Hierarchische,Valentin16Hierarchical}
to define hierarchical B-splines
starting from the well-known nodal B-spline basis
\multicite{Hoellig03Finite,Hoellig13Approximation,Quak16About}.
Note that thanks to the groundwork laid in \cref{chap:20sparseGrids},
especially \thmref{lemma:tensorProductLinearIndependence} and
\thmref{prop:splittingUVToMV},
it suffices to study the univariate case
of all bases that will be defined in the rest of this thesis.
The multivariate case is treated canonically by tensor products.

\subsection{Uniform Hierarchical B-Splines}
\label{sec:311uniform}

\paragraph{Cardinal B-splines}

The \term{cardinal B-spline}
$\cardbspl{p}\colon \real \to \real$ of \term{degree} $p \in \natz$
is defined by
\begin{equation}
  \label{eq:cardinalBSpline}
  \cardbspl{p}(x)
  \ceq
  \begin{cases}
    \displaystyle\int_0^1 \cardbspl{p-1}(x - y) \diff{}y,&p \ge 1,\\
    \charfun{\hopint{0, 1}}(x),&p = 0,
  \end{cases}
\end{equation}
where $\charfun{\hopint{0, 1}}$ is the characteristic function of
the half-open unit interval $\hopint{0, 1}$
(see \cite{Hoellig13Approximation}).
The cardinal B-spline $\cardbspl{p}$ has the following properties
\cite{Hoellig03Finite},
which are shown in \cref{fig:cardinalBSplineProps}:

\begin{figure}
  \includegraphics{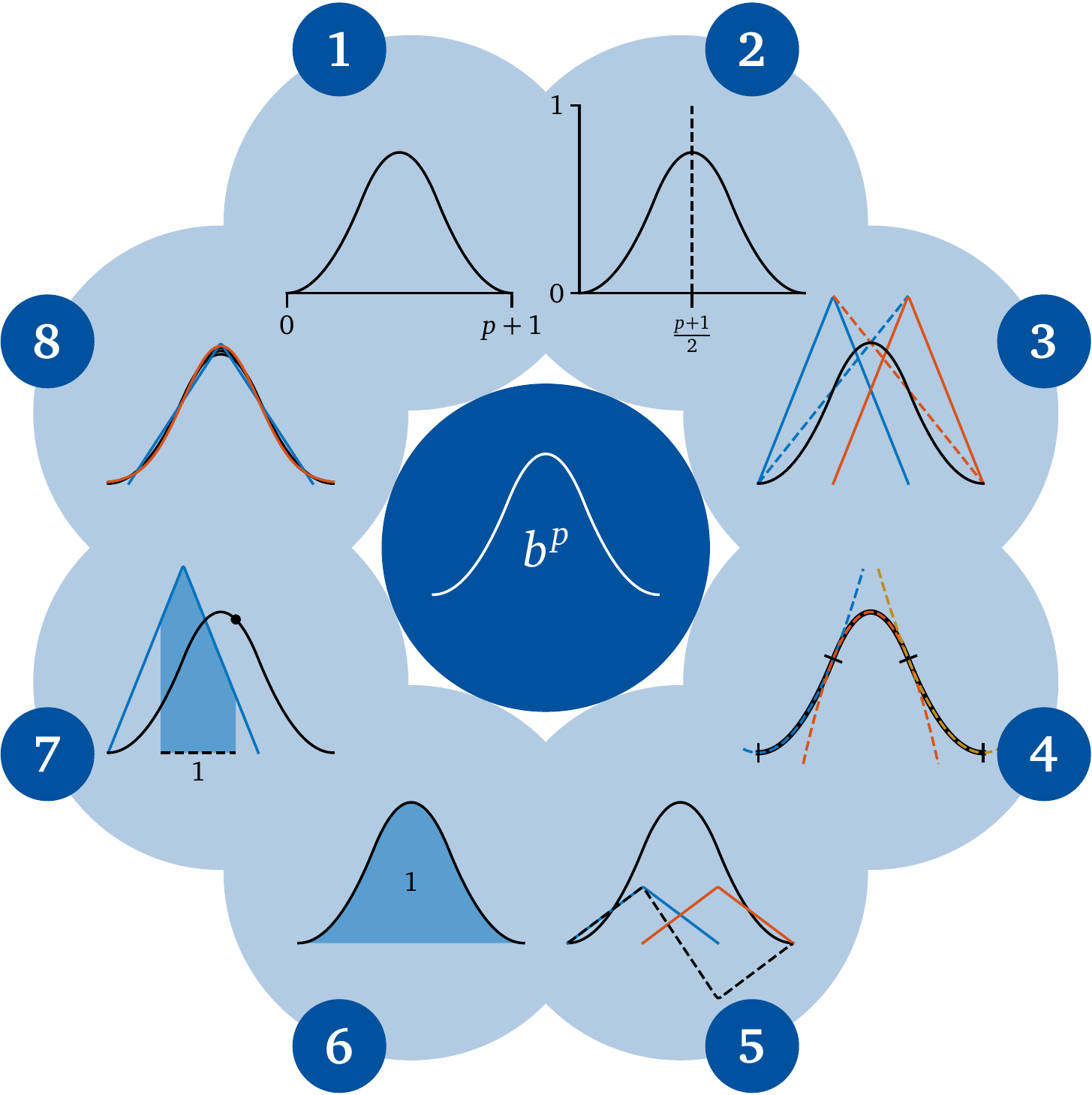}%
  \caption[%
    Properties of cardinal B-splines%
  ]{%
    Eight properties of cardinal B-splines using the quadratic case
    $p = 2$ as an example.\\
    \lefthphantom{1.}{5.}\,
    $\cardbspl{p}$ is compactly supported on $\clint{0, p+1}$.\\
    \lefthphantom{2.}{5.}\,
    $\cardbspl{p}$ is symmetric and $0 \le \cardbspl{p} \le 1$.\\
    \lefthphantom{3.}{5.}\,
    $\cardbspl{p}$ is a weighted combination of
    $\cardbspl{p-1}$ \emph{\textcolor{C0}{(blue)}} and
    $\cardbspl{p-1}({\cdot} - 1)$ \emph{\textcolor{C1}{(red)}.}\\
    \lefthphantom{4.}{5.}\,
    $\cardbspl{p}$ is a piecewise polynomial of degree $p$.\\
    \lefthphantom{5.}{5.}\,
    $\deriv{x}{\cardbspl{p}}$ \emph{(dashed)}
    is the difference of
    $\cardbspl{p-1}$ \emph{\textcolor{C0}{(blue)}} and
    $\cardbspl{p-1}({\cdot} - 1)$ \emph{\textcolor{C1}{(red)}.}\\
    \lefthphantom{6.}{5.}\,
    $\cardbspl{p}$ has unit integral.\\
    \lefthphantom{7.}{5.}\,
    $\cardbspl{p}$ is the convolution of
    $\cardbspl{p-1}$ \emph{\textcolor{C0}{(blue)}} and $\cardbspl{0}$.\\
    \lefthphantom{8.}{5.}\,
    Hat function \emph{\textcolor{C0}{(blue)}} and
    Gaussian function \emph{\textcolor{C1}{(red)}}
    are special cases of $\cardbspl{p}$.%
  }%
  \label{fig:cardinalBSplineProps}%
\end{figure}

\begin{enumerate}
  \item
  \emph{Compact support:}
  The support of $\cardbspl{p}$ is given by $\supp \cardbspl{p} = \clint{0, p + 1}$.
  
  \item
  \emph{Bounds and symmetry:}
  The cardinal B-spline $\cardbspl{p}$ is non-negative and bounded from above by one.
  It is symmetric with respect to $x = \tfrac{p+1}{2}$, i.e.,
  $\cardbspl{p}(x) = \cardbspl{p}(p + 1 - x)$.
  
  \item
  \emph{Recursion:}
  The cardinal B-spline $\cardbspl{p}$ ($p \ge 1$)
  satisfies the following recurrence relation
  (which can be used as an alternative definition):
  \begin{equation}
    \cardbspl{p}(x)
    = \frac{x}{p} \cardbspl{p-1}(x) + \frac{p+1-x}{p} \cardbspl{p-1}(x-1).
  \end{equation}
  
  \item
  \emph{Spline:}
  On every \term{knot interval} $\hopint{k, k+1}$ for $k = 0, \dotsc, p$,
  $\cardbspl{p}$ is a polynomial of degree~$p$, i.e.,
  $\cardbspl{p}$ is a spline of degree $p$ (piecewise polynomial).
  
  \item
  \emph{Derivative:}
  At the \term{knots} $k = 0, \dotsc, p + 1$,
  $\cardbspl{p}$ is $(p - 1)$ times continuously differentiable (if $p \ge 1$).
  The derivative can be computed by differentiating
  \eqref{eq:cardinalBSpline}:
  \begin{equation}
    \label{eq:cardinalBSplineDerivative}
    \deriv{x}{\cardbspl{p}}(x)
    = \cardbspl{p-1}(x) - \cardbspl{p-1}(x-1),\quad
    x \in \real.
  \end{equation}
  
  \item
  \emph{Integral:}
  The B-spline $\cardbspl{p}$ has unit integral, i.e.,
  $\int_{\real} \cardbspl{p}(x) \dx = 1$.
  
  \item
  \emph{Convolution:}
  The integral in the definition of $\cardbspl{p}$
  is the convolution $\cardbspl{p-1} \convolution \cardbspl{0}$
  of the B-spline $\cardbspl{p-1}$
  of degree $p - 1$ with the B-spline $\cardbspl{0}$ of degree zero.
  
  \item
  \emph{Generalization:}
  As a special case, $\cardbspl{1}$ is a hat function,
  interpolating the data
  $\{(k, \kronecker{k}{1}) \mid k \in \integer\}$.
  For $p \to \infty$, the normalized cardinal B-splines converge
  pointwise to the standard Gaussian function
  $\cardbspl{\infty}(x) \ceq (2\pi)^{-1/2} \exp(-x^2/2)$ \cite{Unser92Asymptotic}:%
  \footnote{%
    This can also be seen as a consequence of the central limit theorem
    applied to uniformly distributed random variables.
    The pointwise convergence of the probability density functions
    can be proven from the convergence
    in distribution using a converse to Scheffé's theorem
    \cite{Boos85Converse}.%
  }
  \vspace{-0.5em}
  \begin{equation}
    \lim_{p \to \infty}
    c^p \cardbspl{p}(c^p x + \tfrac{p+1}{2})
    = \cardbspl{\infty}(x),\quad
    c^p \ceq \sqrt{\frac{p+1}{12}},\quad
    x \in \real.
  \end{equation}
\end{enumerate}

\begin{figure}
  \includegraphics{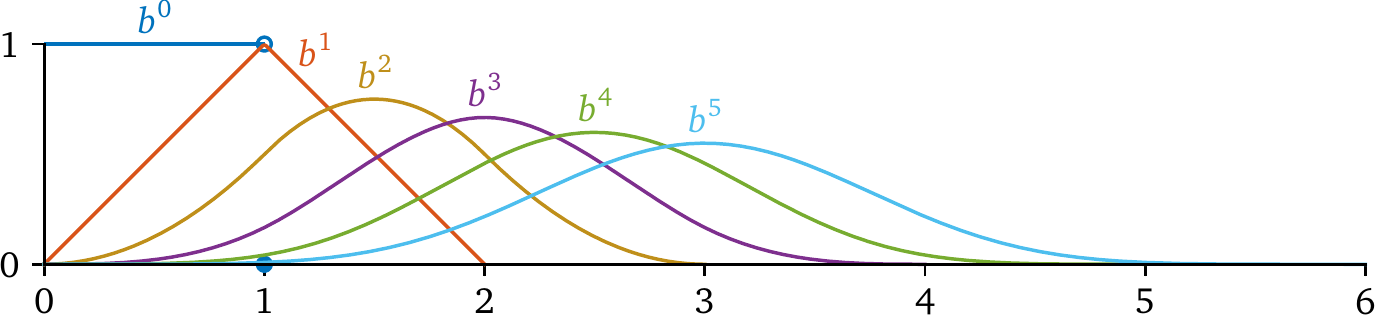}%
  \caption[%
    Cardinal B-splines%
  ]{%
    Cardinal B-splines $\cardbspl{p}$ up to quintic degree $p = 5$.%
  }%
  \label{fig:cardinalBSpline}%
\end{figure}

The cardinal B-splines of the first degrees are shown in
\cref{fig:cardinalBSpline}.
Due to the convolution property,
cardinal B-splines of degree $p \ge 2$ are ``smoothed versions''
of the hat function.
This is shown in the flip book animation in the bottom left corner
of the even-numbered pages of this thesis.

\paragraph{Definition of uniform hierarchical B-splines}

\usenotation{zzzzp}
As for the hat functions in \cref{chap:20sparseGrids},
we can use the cardinal B-spline $\cardbspl{p}$ as a ``parent function'' to
define the uniform hierarchical B-spline
$\bspl{l,i}{p}\colon \clint{0, 1} \to \real$ of level~$l \in \natz$ and index
$i \in \hiset{l}$ via an affine parameter transformation
\cite{Pflueger10Spatially}:
\begin{equation}
  \label{eq:uniformHierarchicalBSplineUV}
  \bspl{l,i}{p}(x)
  \ceq \cardbspl{p}(\tfrac{x}{\ms{l}} + \tfrac{p+1}{2} - i).
\end{equation}
The support of $\bspl{l,i}{p}$ is given
by $\supp \bspl{l,i}{p} = \clint{0, 1} \cap \clint{\gp{l,i-(p+1)/2}, \gp{l,i+(p+1)/2}}$.
The hat function basis $\bspl{l,i}{1}$ defined in
\eqref{eq:hatFunctionUV} is a special case of
\eqref{eq:uniformHierarchicalBSplineUV} for $p = 1$,
which allows us to use the same notation $\bspl{l,i}{p}$ for both.
Note that due to the \term{translation invariance} of $\bspl{l,i}{p}$
(i.e., the basis functions are the same up to scaling and translation),
it suffices to precompute and implement the polynomial pieces of $\cardbspl{p}$
to enable evaluations of all hierarchical B-splines
$\bspl{l,i}{p}$ ($l \in \natz$, $i \in \hiset{l}$).

\begin{figure}
  \subcaptionbox{%
    Nodal B-splines $\bspl{l,i}{p}$ ($i \in \hiset{l}$) and grid points $\gp{l,i}$
    \emph{(dots).}%
  }[67mm]{%
    \includegraphics{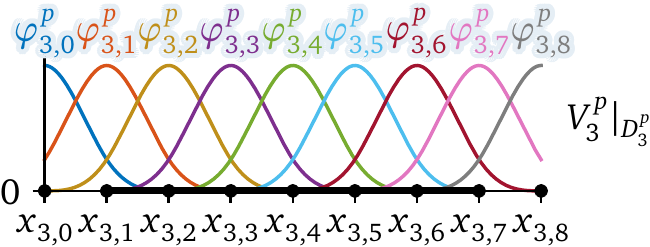}%
  }%
  \hfill%
  \begin{tikzpicture}
    \draw[decorate,decoration={brace,aspect=0.153}] (0,-8.5) -- (0,0);
    \node[anchor=east,inner sep=0mm] at (-0.15,-7.208) {$= \bigoplus$};
  \end{tikzpicture}%
  \hfill%
  \subcaptionbox{%
    Hierarchical B-splines $\bspl{l',i'}{p}$ ($l' \le l$, $i' \in \hiset{l'}$)
    and grid points $\gp{l',i'}$
    \emph{(dots).}%
  }[69mm]{%
    \includegraphics{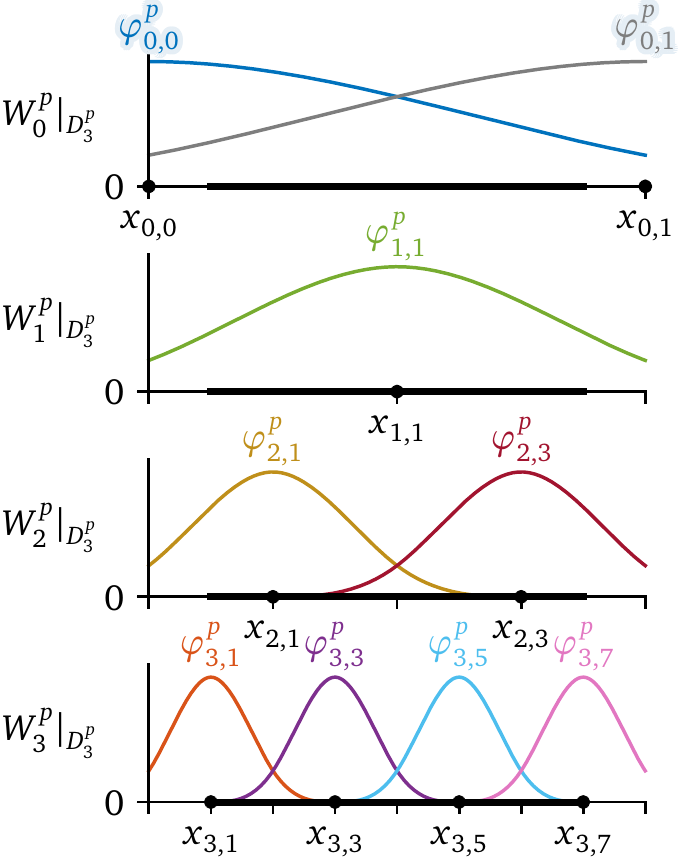}%
  }%
  \caption[%
    Nodal and hierarchical B-splines%
  ]{%
    Univariate nodal and hierarchical cubic B-splines ($p = 3$)
    \vspace{0.1em}%
    up to level $l = 3$.
    When restricting all functions to $\rspldomain{l}{p}$
    \emph{(thick black line),}
    \vspace{-0.2em}%
    the nodal space $\restrictspace{\nsbspl{l}{p}}{\rspldomain{l}{p}}$
    decomposes into the direct sum of the hierarchical subspaces
    $\restrictspace{\hsbspl{l'}{p}}{\rspldomain{l}{p}}$ ($l' \le l$).%
  }%
  \label{fig:hierarchicalBSpline}%
\end{figure}

\paragraph{Odd and even degrees}

In this thesis, we will only allow odd degrees $p = 1, 3, 5, \dotsc$
for hierarchical B-splines.
Many theoretical considerations fail for even degrees.
The basic reason is that for odd degrees, the knots of
$\bspl{l,i}{p}$ coincide with the grid points \cite{Valentin14Hierarchische}
\begin{equation}
  \gp{l,i-(p+1)/2},\quad
  \dotsc,\quad
  \gp{l,i},\quad
  \dotsc,\quad
  \gp{l,i+(p+1)/2}.
\end{equation}
For even degrees $p$, the knots of $\bspl{l,i}{p}$ lie exactly in
the middle between two subsequent grid points:
\begin{equation}
  \gp{l,i-p/2} - \frac{\ms{l}}{2},\quad
  \dotsc,\quad
  \gp{l,i} - \frac{\ms{l}}{2},\quad
  \gp{l,i} + \frac{\ms{l}}{2},\quad
  \dotsc,\quad
  \gp{l,i+p/2} + \frac{\ms{l}}{2}.
\end{equation}
This fact has many adverse implications on the hierarchical basis.
The most crucial implication is
that for even degrees $p$,
the hierarchical splitting \eqref{eq:hierSplittingUV} does not hold.
Furthermore,
we would not be able to define non-uniform hierarchical B-splines as
simple as for odd degrees and
fundamental splines would not be defined at all
(as we will see in \cref{sec:443fundamentalSplines}).
Additionally,
odd degrees include the hat function case~($p = 1$) and the
most commonly applied cubic degree~($p = 3$).
Therefore,
it is reasonable to restrict ourselves to odd degrees
for the rest of the thesis.

\subsection{Non-Uniform B-Splines and Proof of the Hierarchical Splitting}
\label{sec:312proofHierarchicalSplitting}

\paragraph{Non-uniform B-splines and spline space}

With the hierarchical B-splines $\bspl{l,i}{p}$, we can define
the nodal spaces $\nsbspl{l}{p}$ and hierarchical subspaces $\hsbspl{l}{p}$
as in \cref{chap:20sparseGrids}.
However, in order for the hierarchical splitting \eqref{eq:hierSplittingUV}
to be correct, we have to prove that the conditions of
\thmref{lemma:hierSplittingUV} are satisfied.
To investigate how the nodal space $\nsbspl{l}{p}$ looks like,
we introduce the notion of non-uniform B-splines.

\begin{definition}[non-uniform B-splines]
  \label{def:nonUniformBSpline}
  Let $m, p \in \natz$ and $\knotseq = (\knot{0}, \dotsc, \knot{m+p})$ be an
  increasing sequence of real numbers \term{(knot sequence).}
  For $k = 0, \dotsc, m - 1$,
  the \term{(non-uniform) B-spline} $\nonunifbspl{k,\knotseq}{p}$ of degree $p$
  with knots~$\knotseq$ and index $k$ is defined by the
  Cox--de~Boor recurrence
  \multicite{Cox72Numerical,Boor72Calculating,Hoellig13Approximation}
  \begin{equation}
    \nonunifbspl{k,\knotseq}{p}(x)
    \ceq
    \begin{cases}
      \dfrac{x - \knot{k}}{\knot{k+p} - \knot{k}} \nonunifbspl{k,\knotseq}{p-1}(x) +
      \dfrac{\knot{k+p+1} - x}{\knot{k+p+1} - \knot{k+1}}
      \nonunifbspl{k+1,\knotseq}{p-1}(x),&p \ge 1,\\
      \charfun{\hopint{\knot{k}, \knot{k+1}}}(x),&p = 0.
    \end{cases}
    \hspace*{-4mm}
  \end{equation}
\end{definition}
Note that when choosing $\knotseq = (0, 1, \dotsc, p + 1)$ and
$k = 0$, we obtain the cardinal B-spline~$\cardbspl{p}$.
\Cref{def:nonUniformBSpline} can be used to characterize
the nodal space $\nsbspl{l}{p}$:

\begin{proposition}[spline space]
  \label{prop:splineSpace}
  Let $\knotseq = (\knot{0}, \dotsc, \knot{m+p})$ be a knot sequence.
  Then, the B-splines $\nonunifbspl{k,\knotseq}{p}$ ($k = 0, \dotsc, m - 1$)
  form a basis of the \term{spline space}
  \begin{equation}
    \nonunifsplspace{\knotseq}{p}
    \ceq \spn\{\nonunifbspl{k,\knotseq}{p} \mid k = 0, \dotsc, m - 1\}.
  \end{equation}
  $\nonunifsplspace{\knotseq}{p}$ contains exactly those functions that are continuous
  on $\spldomain{\knotseq}{p} \ceq \clint{\knot{p}, \knot{m}}$,
  polynomials of degree $\le p$ on every knot interval
  $\hopint{\knot{k}, \knot{k+1}}$ in
  $\spldomain{\knotseq}{p}$
  ($k = p, \dotsc, m - 1$) and at least $(p - 1)$ times
  continuously differentiable at every knot $\knot{k}$ in the interior of
  $\spldomain{\knotseq}{p}$ ($k = p + 1, \dotsc, m - 1$).
\end{proposition}

\begin{proof}
  See \cite{Hoellig13Approximation}.
\end{proof}

\begin{figure}
  \includegraphics{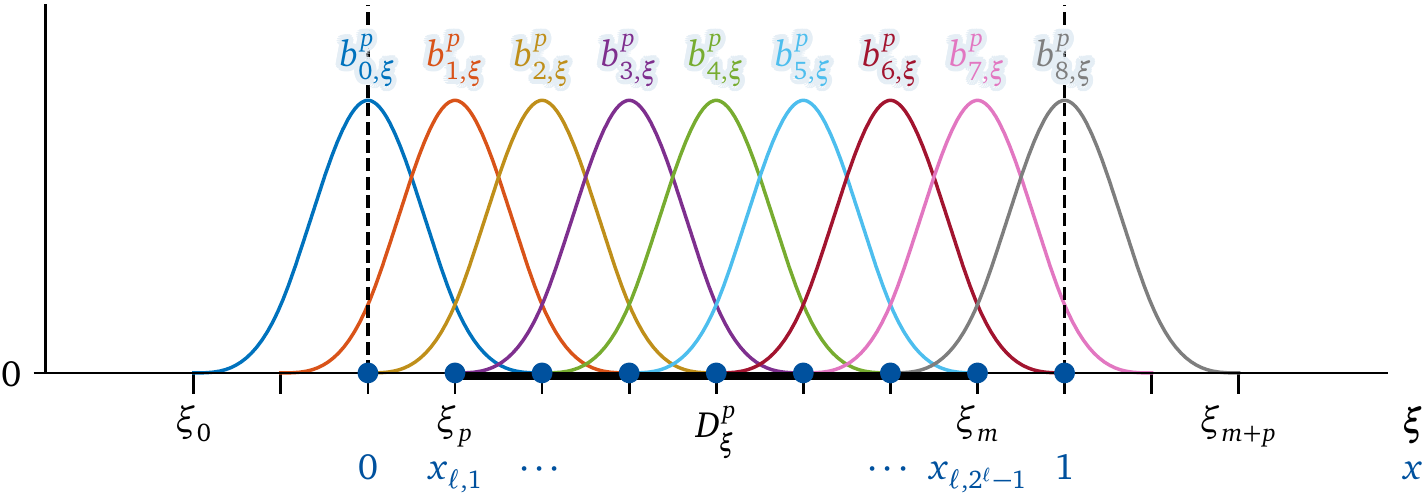}%
  \caption[%
    Uniform nodal B-splines and knot sequence%
  ]{%
    Uniform knot sequence $\nodalknotseq{l}{p}$ \emph{(ticks on horizontal axis)}
    and corresponding nodal cubic B-splines ($p = 3$) of level $l = 3$.
    In the domain $\clint{0, 1}$ \emph{(delimited by dashed lines),}
    the grid points $\fgset{l}$ \emph{\textcolor{mittelblau}{(blue dots)}}
    coincide with the B-spline knots.
    The spline interpolation domain
    $\rspldomain{l}{p} \ceq \spldomain{\nodalknotseq{l}{p}}{p}$
    \emph{(thick line)}
    is only a proper subset of~$\clint{0, 1}$.%
  }%
  \label{fig:splineSpaceUniform}%
\end{figure}

This proposition gives the reason for the letter ``B'' in ``B-splines,''
which stands for ``basis'' (of the space of splines) \cite{Schoenberg67Spline}.
One example of a knot sequence and the corresponding B-splines is
given in \cref{fig:splineSpaceGeneral}.
The key observation is that B-splines of a knot sequence $\knotseq$
do not form a basis of the spline space on the union
$\clint{\knot{0}, \knot{m+p}}$ of the B-spline supports.
Instead, they form a basis of the spline space
on a proper sub-interval $\spldomain{\knotseq}{p}$.
Intuitively, for every point in $\spldomain{\knotseq}{p}$ that is not a knot,
exactly $p + 1$ B-splines must be \term{relevant} (i.e., non-zero)
to uniquely span the spline space,
as on every knot interval, the spline is a polynomial of degree $\le p$
and therefore, there must be $p + 1$ degrees of freedom.
Outside $\spldomain{\knotseq}{p}$, there are too few relevant B-splines
to span the spline space.
This fact,
which is shown in \cref{fig:splineSpaceUniform},
forces us to restrict the nodal space and the hierarchical subspaces
to~$\spldomain{\knotseq}{p}$:

\begin{figure}
  \includegraphics{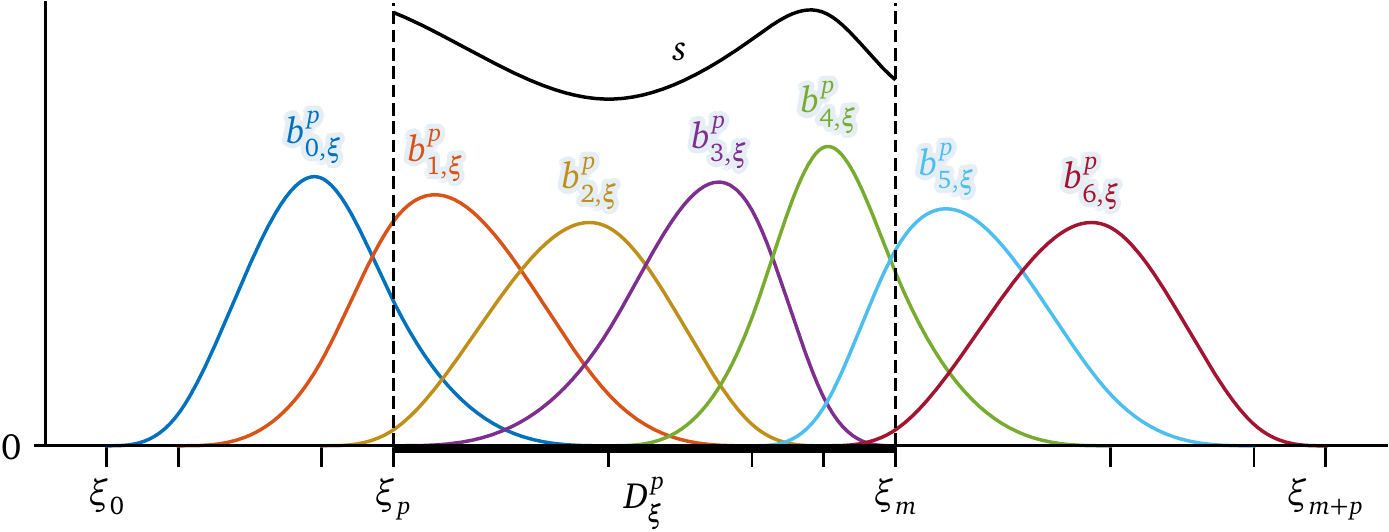}%
  \caption[%
    Non-uniform B-splines with knot sequence and interpolation domain%
  ]{%
    Knot sequence $\knotseq = (\knot{0}, \dotsc, \knot{m+p})$
    \vphantom{$\nonunifbspl{0,\knotseq}{p}$}%
    with the corresponding $m = 7$ non-uniform cubic B-splines
    $\nonunifbspl{k,\knotseq}{p}$ ($k = 0, \dotsc, m - 1$, $p = 3$).
    On $\spldomain{\knotseq}{p}$ \emph{(thick line, delimited by dashed lines),}
    which starts with the last knot interval of the first B-spline
    $\nonunifbspl{0,\knotseq}{p}$
    and ends with the first knot interval of the last B-spline
    $\nonunifbspl{m-1,\knotseq}{p}$,
    the B-splines span the spline space $\nonunifsplspace{\knotseq}{p}$.
    Elements of this space are splines $\spl\colon \spldomain{\knotseq}{p} \to \real$
    \emph{(black line),}
    which are linear combinations
    $\spl = \sum_{k=0}^{m-1} \interpcoeff{k} \nonunifbspl{k,\knotseq}{p}$
    of the B-splines.%
  }%
  \label{fig:splineSpaceGeneral}%
\end{figure}

\begin{corollary}[nodal B-spline space]
  \label{cor:nodalBSplineSpace}
  The restricted nodal B-splines $\restrictfcn{\bspl{l,i}{p}}{\rspldomain{l}{p}}$
  ($i = 0, \dotsc, 2^l$)
  of level $l \in \natz$ are
  a basis of the spline space $\nonunifsplspace{\nodalknotseq{l}{p}}{p}$,
  where
  \begin{subequations}
    \begin{gather}
      \nodalknot{l,k}{p}
      \ceq (k - \tfrac{p+1}{2}) \ms{l},\quad
      k = 0, \dotsc, m + p,\quad
      m \ceq 2^l + 1,\\
      \rspldomain{l}{p} \ceq [\tfrac{p-1}{2} \ms{l},\;
      1 - \tfrac{p-1}{2} \ms{l}],
    \end{gather}
  \end{subequations}
  and consequently
  \begin{equation}
    \restrictspace{\nsbspl{l}{p}}{\rspldomain{l}{p}}
    = \restrictedsplspace{l}{p}
    \ceq \nonunifsplspace{\nodalknotseq{l}{p}}{p}.
  \end{equation}
\end{corollary}

\begin{proof}
  We have $\bspl{l,i}{p} = \nonunifbspl{i,\nodalknotseq{l}{p}}{p}$ for
  $i = 0, \dotsc, m - 1$,
  as the B-splines on both sides have the same knots.
  The assertions now follow from \thmref{prop:splineSpace}.
\end{proof}

\vspace{1em}

Note that $\rspldomain{l}{p}$ might contain only a single point or even be empty,
if $p$ is too large or $l$ is too small.
However, the corresponding B-splines $\bspl{l,i}{p}$ are still linearly
independent on~$\clint{0, 1}$ (see \cite{Hoellig13Approximation}).
Similarly, the corollary also implies that the hierarchical functions
$\bspl{l,i}{p}$ of level $l$ ($i \in \hiset{l}$)
are linearly independent on $\clint{0, 1}$.

\paragraph{Hierarchical splitting for uniform B-splines}

We can use \cref{prop:splineSpace} and \cref{cor:nodalBSplineSpace}
to prove the hierarchical splitting for the uniform B-spline basis
\cite{Valentin16Hierarchical}.

\begin{lemma}[hierarchical B-splines in nodal space]
  \label{lemma:hierBSplineInNodalSpace}
  The restricted hierarchical subspaces
  $\restrictspace{\hsbspl{l'}{p}}{\rspldomain{l}{p}}$ ($l' \le l$) are
  subspaces of the restricted nodal space $\restrictspace{\nsbspl{l}{p}}{\rspldomain{l}{p}}$.
\end{lemma}

\begin{proof}
  Every function $\bspl{l',i'}{p}$ ($i' \in \hiset{l'}$) is continuous on
  $\rspldomain{l}{p}$, a polynomial of degree $\le p$ on every knot interval
  of $\nodalknotseq{l}{p}$ (due to $p$ odd),
  and at the knots themselves at least $(p - 1)$ times continuously
  differentiable.
  \Cref{prop:splineSpace} implies $\bspl{l',i'}{p} \in \restrictedsplspace{l}{p}$
  and from \cref{cor:nodalBSplineSpace}, it follows
  $\bspl{l',i'}{p} \in \restrictspace{\nsbspl{l}{p}}{\rspldomain{l}{p}}$.
  As the functions $\bspl{l',i'}{p}$ ($i' \in \hiset{l'}$) span
  $\restrictspace{\hsbspl{l'}{p}}{\rspldomain{l}{p}}$, we can conclude
  $\restrictspace{\hsbspl{l'}{p}}{\rspldomain{l}{p}} \subset \restrictspace{\nsbspl{l}{p}}{\rspldomain{l}{p}}$.
\end{proof}

\vspace{1em}

It is crucial to note that this lemma does not hold for even $p$,
as the knots of the B-splines of level $l - 1$ are not contained in the
knots of level $l$.
This implies that in general,
$\restrictspace{\hsbspl{l-1}{p}}{\rspldomain{l}{p}}$
is not contained in $\restrictspace{\nsbspl{l}{p}}{\rspldomain{l}{p}}$
for even degrees $p$.
Therefore, the hierarchical splitting equation \eqref{eq:hierSplittingUV}
does not hold.

\begin{restatable}[hierarchical B-splines are linearly independent]{%
  proposition%
}{%
  propHierBSplineLinearlyIndependent%
}
  \label{prop:hierBSplineLinearlyIndependent}
  The hierarchical B-splines
  $\bspl{l',i'}{p}$ ($l' \le l$, $i' \in \hiset{l'}$)
  are linearly independent.
\end{restatable}

\begin{proof}
  See \cref{sec:a121proofHierBSplineLinearlyIndependent}.
\end{proof}

\vspace{1em}

Although we have to restrict all functions and spaces to $\rspldomain{l}{p}$,
\thmref{lemma:hierSplittingUV} is still applicable to prove that
the hierarchical splitting equation \eqref{eq:hierSplittingUV}
is correct for hierarchical B-splines:

\begin{corollary}[hierarchical splitting for uniform B-splines]
  \label{cor:hierSplittingBSpline}
  The hierarchical splitting \eqref{eq:hierSplittingUV}
  holds for the hierarchical B-spline basis
  if restricting all functions to $\rspldomain{l}{p}$:
  \begin{equation}
    \restrictedsplspace{l}{p}
    = \restrictspace{\nsbspl{l}{p}}{\rspldomain{l}{p}}
    = \bigoplus_{l'=0}^l \restrictspace{\hsbspl{l'}{p}}{\rspldomain{l}{p}}.
  \end{equation}
\end{corollary}

\begin{proof}
  Analogously to the proof of
  \hyperlink{lemma:hierSplittingUV}{%
    \cref*{lemma:hierSplittingUV}~%
    (univariate hierarchical splitting characterization)%
  }
  and apply \thmref{cor:nodalBSplineSpace}.
\end{proof}

\vspace{1em}

This corollary is also visualized in \cref{fig:hierarchicalBSpline}.
We can now proceed to define
multivariate nodal spaces $\nsbspl{\*l}{p}$,
multivariate hierarchical subspaces $\hsbspl{\*l}{p}$, and
sparse grid spaces $\regsgspace[p]{n}{d}$ as in \cref{chap:20sparseGrids}.
Note that it is possible to choose different degrees $p_t$ for
different dimensions $t = 1, \dotsc, d$,
since the hierarchical splitting \eqref{eq:hierSplittingMV} does not
require the bases in each dimension to be the same.
Consequently, we can define degree-dimension-adaptive
(so-called \term{$hp$-adaptive}) sparse grids
$\regsgspace[\*p]{n}{d}$ for arbitrary odd degree vectors $\*p$.

In the course of this thesis, we will derive multiple variations
of the standard hierarchical B-spline basis.
We will not repeat formal proofs of the hierarchical splitting equation
\eqref{eq:hierSplittingUV}
(i.e., verifying the two conditions of \cref{lemma:hierSplittingUV})
for each of these bases for the sake of brevity.
The idea of the proof of \cref{prop:hierBSplineLinearlyIndependent},
which is inductively exploiting the smoothness conditions given by
B-splines of previous levels, can be transferred to similar B-spline
bases.

\subsection{Modification}
\label{sec:313modification}

\paragraph{Marsden's identity}

Similar to the piecewise linear case in \cref{sec:242modified},
Pflüger defined modified
hierarchical B-splines to obtain reasonable values on the boundary
without having to place grid points there \cite{Pflueger10Spatially}.
The main motivation is to define basis
functions $\bspl[\modified]{l,i}{p}$ that satisfy natural boundary
conditions, i.e.,
\begin{equation}
  \label{eq:naturalBoundaryConditions}
  \deriv[2]{x}{\bspl[\modified]{l,i}{p}}(x) = 0,\quad
  x \in \bndrydomain{\clint{0, 1}} = \{0, 1\}.
\end{equation}
Originally, this requirement stems from financial problems
\cite{Pflueger10Spatially}.
For the left boundary,
\eqref{eq:naturalBoundaryConditions} can be satisfied by
modifying the left-most function $\bspl{l,1}{p}$ such that
$\bspl[\modified]{l,1}{p}(x) = 2 - \tfrac{x}{\ms{l}}$ is a linear polynomial
when $x$ is ``near'' the boundary.
As in \cite{Pflueger10Spatially},
we append $\bspl{l,1}{p}$ with
B-splines $\bspl{l,i}{p}$ with index $i \le 0$ and
use \term{Marsden's identity} to compute the corresponding coefficients.
The identity enables us to explicitly compute the B-spline coefficients
of an arbitrary polynomial of degree $\le p$ by giving an explicit formula
for the B-spline coefficients of shifted monomials $({\cdot} - y)^p$.
Interestingly, the coefficients are polynomials themselves (in $y$):

\vspace*{0pt plus 0.5fill}

\begin{lemma}[Marsden's identity]
  \label{lemma:marsden}
  Let $p \in \natz$ and
  $\knotseq = (\knot{0}, \dotsc, \knot{m+p})$ be a knot sequence.
  Then, for all $x \in \spldomain{\knotseq}{p}$ and $y \in \real$,
  \begin{equation}
    \label{eq:marsden}
    (x - y)^p
    = \sum_{k=0}^{m-1}
    \;\bracket*{(\knot{k+1} - y) \dotsm (\knot{k+p} - y)}\;
    \nonunifbspl{k,\knotseq}{p}(x).
  \end{equation}
\end{lemma}

\begin{proof}
  See \cite{Hoellig13Approximation}.
\end{proof}

\vspace*{0pt plus 1fill}

\paragraph{Modified hierarchical B-splines}

\begin{figure}
  \subcaptionbox{%
    Construction of the modified
    \vphantom{$\bspl[\modified]{l',i'}{p}$}\vspace{-0.2em}%
    hierarchical cubic B-spline
    $\bspl[\modified]{l',1}{p}$ (\emph{dashed,} $l' \ge 2$)
    \vphantom{$\bspl[\modified]{l',i'}{p}$}\vspace{-0.2em}%
    as a linear combination of neighboring B-splines $\bspl{l',i'}{p}$.%
    \vphantom{$\bspl[\modified]{l',i'}{p}$}%
  }[72mm]{%
    \includegraphics{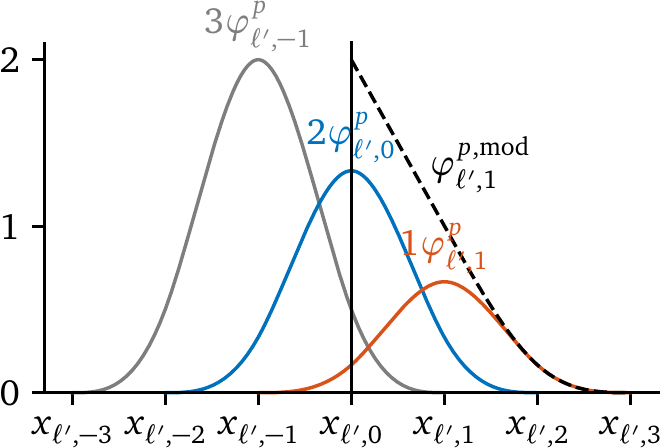}%
  }%
  \hfill%
  \subcaptionbox{%
    Modified hierarchical cubic B-splines
    \vphantom{$\bspl[\modified]{l',i'}{p}$}\vspace{-0.2em}%
    $\bspl[\modified]{l',i'}{p}$
    ($l' = 1, \dotsc, l$, $i' \in \hiset{l'}$) and
    \vphantom{$\bspl[\modified]{l',i'}{p}$}\vspace{-0.2em}%
    grid points $\gp{l',i'}$
    \emph{(dots).}%
    \vphantom{$\bspl[\modified]{l',i'}{p}$}%
  }[72mm]{%
    \includegraphics{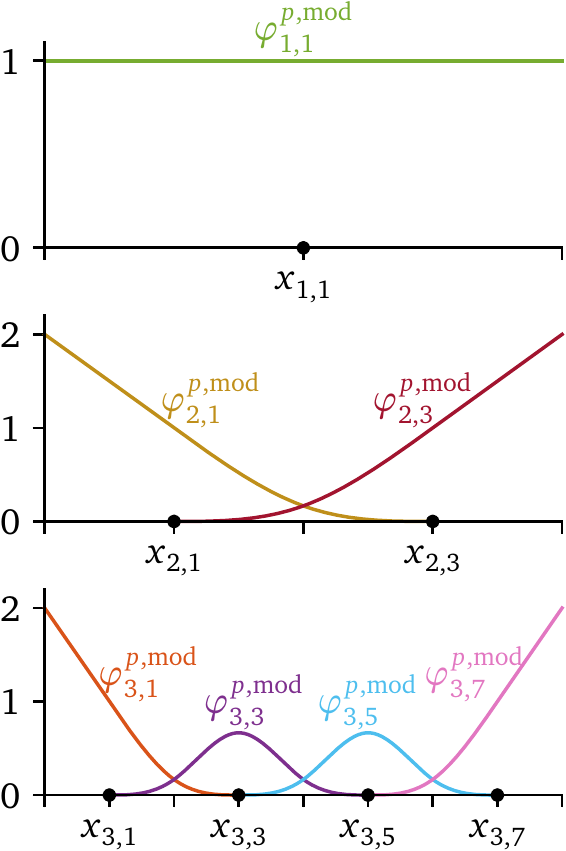}%
  }%
  \caption[%
    Modified hierarchical B-splines%
  ]{%
    Construction of modified hierarchical cubic B-splines ($p = 3$) and
    the resulting basis up to level $l = 3$.%
  }%
  \label{fig:modifiedBSpline}%
\end{figure}

By extending the sum in Marsden's identity to $i \in \integer$ and
comparing the coefficients of $y^p$ and $y^{p-1}$ of both sides,
we obtain representations for the first two monomials:
\begin{equation}
  1
  = \sum_{i \in \integer} \bspl{l,i}{p}(x),\quad
  x
  = \sum_{i \in \integer} \gp{l,i} \bspl{l,i}{p}(x),\quad
  x \in \real.
\end{equation}
This can be used to derive a definition for $\bspl[\modified]{l,i}{p}$:
\begin{equation}
  \label{eq:modifiedBSplineConstruction}
  2 - \frac{x}{\ms{l}}
  = 2 \sum_{i \in \integer} \bspl{l,i}{p}(x)
  - \frac{1}{\ms{l}} \sum_{i \in \integer} \gp{l,i} \bspl{l,i}{p}(x)
  = \sum_{i \in \integer} (2 - i) \bspl{l,i}{p}(x),\quad
  x \in \real.
\end{equation}
Note that only the B-splines with $i \ge 1 - \tfrac{p+1}{2}$
are relevant for the unit interval
as all other B-splines vanish in $\clint{0, 1}$.
Pflüger omits summands with $i > 1$ as he only wants to modify
$\bspl{l,1}{p}$ left of its grid point $\gp{l,1}$ \cite{Pflueger10Spatially}.
The right-most function $\bspl[\modified]{l,2^l-1}{p}$ can be derived
analogously by mirroring $\bspl[\modified]{l,1}{p}$ at $x = \tfrac{1}{2}$.
\pagebreak%
For~$l = 1$, again the ``constant one'' function is taken for the definition
of modified hierarchical B-splines (see \cref{fig:modifiedBSpline}):
{%
  \setlength{\abovedisplayskip}{9pt}%
  \setlength{\belowdisplayskip}{9pt}%
  \begin{equation}
    \label{eq:modifiedBSplineDefinition}
    \bspl[\modified]{l,i}{p}(x)
    \ceq
    \begin{cases}
      1,&
      l = 1,\quad i = 1,\\
      \displaystyle\sum_{i'=1-(p+1)/2}^1\!\!\!\! (2 - i') \bspl{l,i'}{p}(x),&
      l \ge 2,\quad i = 1,\\
      \bspl{l,i}{p}(x),&
      l \ge 2,\quad i \in \hiset{l} \setminus \{1, 2^l - 1\},\\
      \bspl[\modified]{l,1}{p}(1 - x),&
      l \ge 2,\quad i = 2^l - 1.
    \end{cases}
  \end{equation}%
}%
By \thmref{prop:splineSpace},
this definition implies that
$\bspl{l,1}{p}(x) = 2 - \tfrac{x}{\ms{l}}$ ($l \ge 2$)
is only valid for $x \in \clint{0, \tfrac{5-p}{2} \ms{l}}$, i.e.,
the second derivative at $x = 0$ vanishes only for $p \le 5$.
For higher degrees, it is non-zero, albeit very small
in its absolute value.
To enforce natural boundary conditions
for higher than quintic degrees,
the upper bound of $i'$ in the sum in \eqref{eq:modifiedBSplineDefinition}
must be extended to $\tfrac{p+1}{2}$.
In addition, more than just the left-most and the right-most inner
hierarchical B-spline must be modified for $p \ge 5$,
as the size of the supports of $\bspl{l,i}{p}$ increases
for growing $p$.

\subsection{Non-Uniform Hierarchical B-Splines}
\label{sec:314nonUniform}

Sparse grid spaces and their corresponding grid point sets,
as we have defined them in \cref{chap:20sparseGrids},
are completely independent of the actual location of the grid points
$\gp{l,i}$.
Therefore, it is possible to use different distributions for the grid points
than the standard equidistant choice of $\gp{l,i} = i \cdot \ms{l}$,
if the basis functions are altered accordingly
\cite{Valentin14Hierarchische}.
\usenotation{zzzzcc}
The so-called Chebyshev points $\ccgp{l,i}$ and the
resulting Clenshaw--Curtis grids and B-splines will serve as an example.

\paragraph{Chebyshev points and Clenshaw--Curtis grids}

\begin{SCfigure}
  \includegraphics{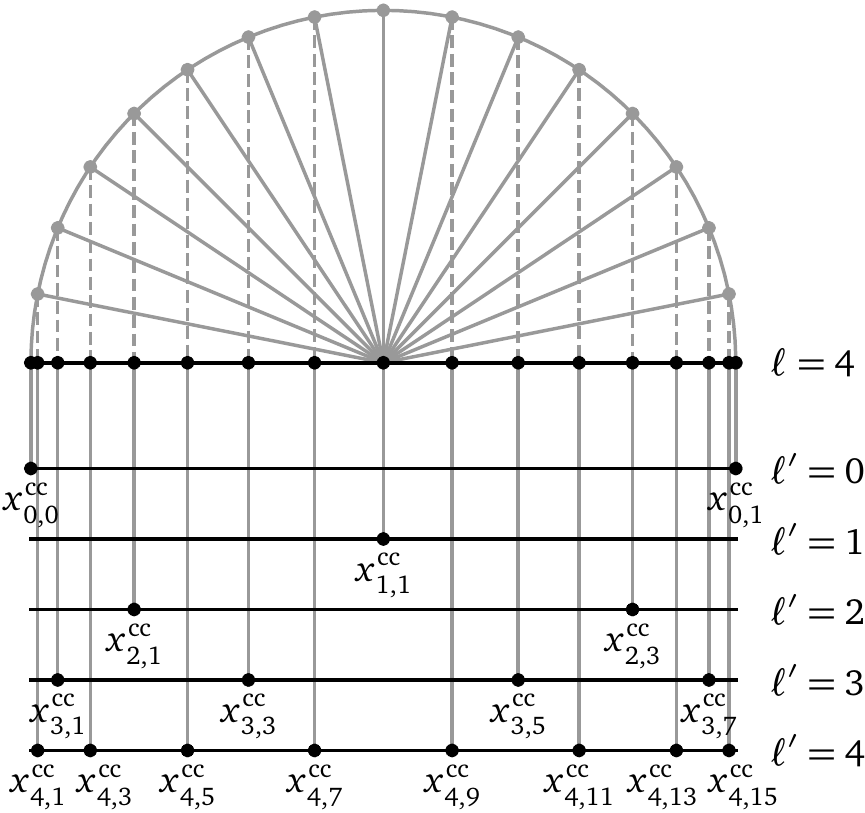}%
  \caption[%
    Decomposition of the set of univariate Clenshaw--Curtis grid points%
  ]{%
    The set of Chebyshev points $\fgset[\cc]{l}$ of level
    $l = 4$ \emph{(top)}
    decomposes into hierarchical grids of level $l' \le l$
    (compare with \cref{fig:pointSplittingUniform}).
    The Chebyshev points are constructed as
    the orthogonal projection of the
    endpoints of $2^l$ equally-sized segments
    of a semicircle onto its diameter \emph{\textcolor{C8}{(gray, top)}.}%
  }%
  \label{fig:pointSplittingChebyshev}%
\end{SCfigure}

The \term{Chebyshev points} $\ccgp{l,i}$ of level $l \in \natz$
are defined as
\begin{equation}
  \ccgp{l,i}
  \ceq \frac{1 - \cos(\pi \gp{l,i})}{2},\quad
  i = 0, \dotsc, 2^l,
\end{equation}
see \cite{Xu16Chebyshev} for example.
Chebyshev points are the locations of the extrema%
\footnote{%
  The literature sometimes uses the name ``Chebyshev points'' for
  the roots of $T_{2^l}$, which are closely connected with the extrema.
  One way to distinguish them is to call the extrema
  ``Chebyshev--Lobatto points'' and the roots
  ``Chebyshev--Gauss points'' \cite{Xu16Chebyshev}.%
}
of the Chebyshev polynomials $T_{2^l}$, defined as
$T_{2^l}\colon \clint{0, 1} \to \real$,
$T_{2^l}(x) \ceq \cos(2^l \arccos(2x - 1))$.
They are geometrically obtained
by dividing a semicircle into $2^l$ equally-sized
segments and subsequently orthogonally projecting the
segment endpoints onto the diameter
(see \cref{fig:pointSplittingChebyshev}).
Analytically, they are determined by minimizing
the interpolation error term in polynomial interpolation.
The most practical use of Chebyshev points is in
polynomial interpolation to avoid Runge's phenomenon and in
numerical integration (quadrature), resulting in the
so-called Clenshaw--Curtis quadrature rules.

In some settings, it may be beneficial to use full or sparse grids consisting
of Chebyshev points, which are then called \term{Clenshaw--Curtis grids,}
instead of uniform grids.
Besides the already mentioned advantages for polynomials and
quadrature, Clenshaw--Curtis grids can help to reduce interpolation
errors in a neighborhood of the boundary of the domain due to the increased
grid point density near the boundary
(at the cost of a lower grid point density in the interior).
If we want to use Chebyshev points as grid points for sparse grids,
we have to employ an appropriate basis to ensure that interpolation
is still possible.

\paragraph{Hierarchical Clenshaw--Curtis B-splines}

\begin{figure}
  \subcaptionbox{%
    Hierarchical cubic Clenshaw--Curtis B-splines
    $\bspl[\cc]{l',i'}{p}$
    ($l' \le l$, $i' \in \hiset{l'}$, $p = 3$) and
    \vspace{-0.2em}%
    modified Clenshaw--Curtis B-splines
     $\bspl[\cc,\modified]{l',i'}{p}$
    \emph{(dashed).}%
  }[72mm]{%
    \includegraphics{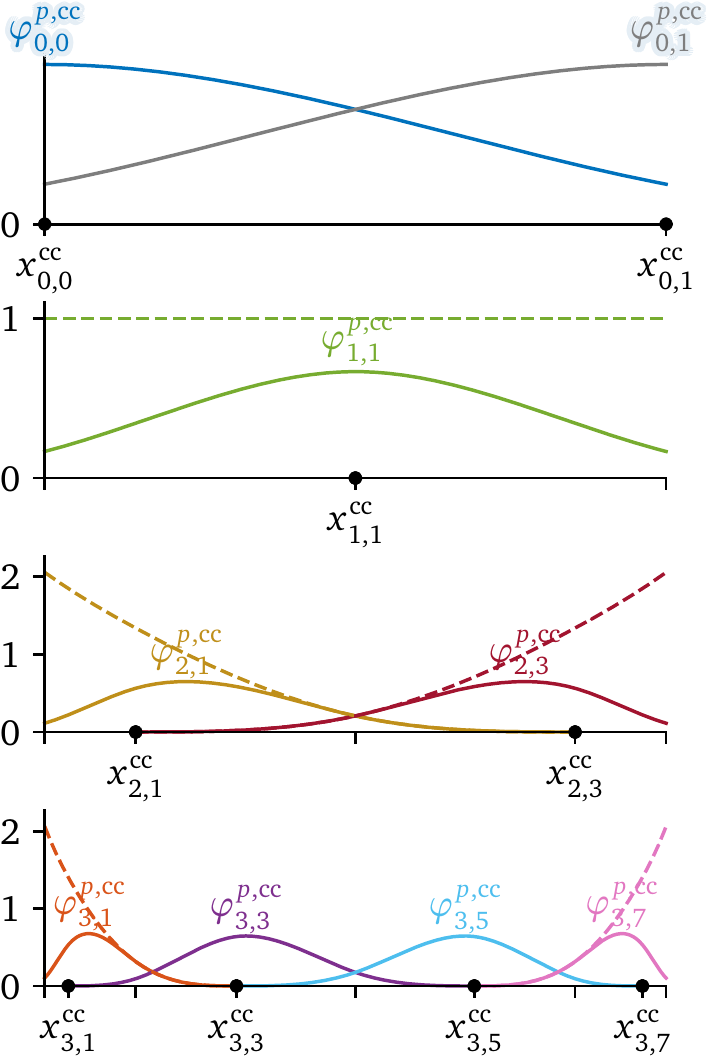}%
  }%
  \hfill%
  \subcaptionbox{%
    Sparse Clenshaw--Curtis grid $\regsgset[\cc]{n}{d}$
    of level $n = 4$ in $d = 2$ dimensions.%
  }[72mm]{%
    \includegraphics{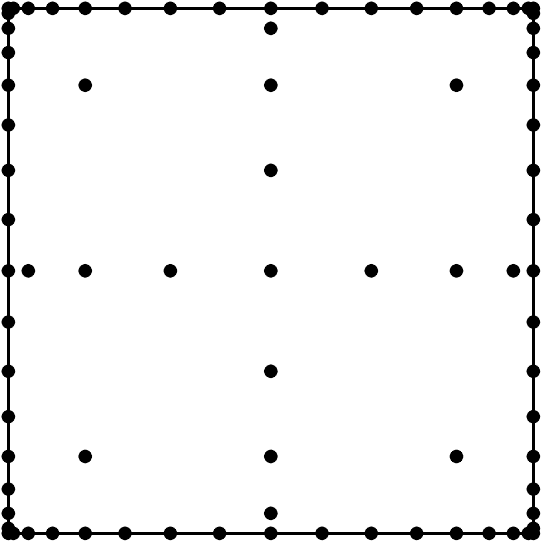}%
  }%
  \caption[%
    Clenshaw--Curtis B-splines and sparse grids%
  ]{%
    Clenshaw--Curtis B-splines and sparse grids.%
  }%
  \label{fig:clenshawCurtis}%
\end{figure}

The \term{hierarchical Clenshaw--Curtis B-spline}
$\bspl[\cc]{l,i}{p}$ of level $l \in \natz$ and index
$i \in \hiset{l}$ is defined as \cite{Valentin14Hierarchische}
\begin{equation}
  \bspl[\cc]{l,i}{p}
  \ceq \nonunifbspl{i,\nodalknotseq[\cc]{l}{p}}{p},
\end{equation}
where
\begin{subequations}
  \label{eq:clenshawCurtisBSplineKnots}
  \begin{gather}
    \nodalknot[\cc]{l,k}{p}
    \ceq
    \begin{cases}
      (k - \frac{p+1}{2}) \cdot \ccgp{l,1},&
      k = 0,\, \dotsc,\, \frac{p-1}{2},\\
      \ccgp{l,k-(p+1)/2},&
      k = \frac{p+1}{2},\, \dotsc,\, 2^l + \frac{p+1}{2},\\
      1 + (k - 2^l - \frac{p+1}{2}) \cdot \ccgp{l,1},&
      k = 2^l + \frac{p+3}{2},\, \dotsc,\, 2^l + p + 1,
    \end{cases}\\
    k = 0, \dotsc, m + p,\quad
    m \ceq 2^l + 1.
  \end{gather}
\end{subequations}
For the construction of the knot sequence $\nodalknotseq[\cc]{l}{p}$,
the Chebyshev points $\ccgp{l,i}$
are equidistantly extended onto the real line $\real$.
As for the standard hierarchical B-spline basis,
it is now straightforward to define nodal spaces
$\nsbspl[\cc]{l}{p}$
and hierarchical subspaces $\hsbspl[\cc]{l}{p}$ as well as
sparse grid spaces $\regsgspace[p,\cc]{n}{d}$ and
grid point sets $\regsgset[\cc]{n}{d}$
using tensor products of Clenshaw--Curtis B-splines
and Cartesian products of Chebyshev point sets.
The one-dimensional cubic Clenshaw--Curtis basis and a
sparse Clenshaw--Curtis grid are shown in \cref{fig:clenshawCurtis}.

Note that Clenshaw--Curtis B-splines
$\bspl[\cc]{l,i}{p}$ are not translation-invariant,
in contrast to uniform B-splines $\bspl{l,i}{p}$.
As a result, both implementation effort and computation time for evaluation
are significantly higher for Clenshaw--Curtis B-splines than
for uniform B-splines,
as the Clenshaw--Curtis B-splines cannot be precomputed.

\paragraph{Modification and generalizations}

We define \term{modified hierarchical Clenshaw--Curtis B-splines}
$\bspl[\cc,\modified]{l,i}{p}$ using the
same method as in \cref{eq:modifiedBSplineDefinition}.
Here, the second derivative of $\bspl[\cc,\modified]{l,1}{p}$
does not vanish at $x = 0$ even for degrees $p \le 5$,
as the formula \eqref{eq:modifiedBSplineConstruction}
derived from \thmref{lemma:marsden} assumes uniform knots.
However, since most of the B-spline knots in the summation formula
lie outside $\clint{0, 1}$ and are thus uniform according
to \cref{eq:clenshawCurtisBSplineKnots},
the absolute deviation of the second derivative from zero is small.
Again, to enforce natural boundary conditions,
we can recompute the coefficients
of the components $\bspl[\cc]{l,i'}{p}$
dynamically with Marsden's identity using the correct Chebyshev knots
in \cref{eq:marsden}.

Note that our framework permits arbitrary grid point distributions
$\gp{l,i}^\ast$,
as long as two requirements are met:
First, their number should grow exponentially
(i.e., there are $2^l + 1$ points $\gp{l,i}^\ast$ in each level $l \in \natz$),
and second, they should be nested
(i.e., \cref{eq:rewriteGridPoint} holds).
Appropriate non-uniform B-spline bases can be defined analogously
to Clenshaw--Curtis B-splines.

\section{Boundary Behavior of Hierarchical B-Splines}
\label{sec:32notAKnot}

\minitoc{87mm}{8}

\noindent
As we have seen in the last section (see \cref{cor:hierSplittingBSpline}),
the hierarchical splitting equation \eqref{eq:hierSplittingUV}
only holds when restricting the function spaces to
$\rspldomain{l}{p} =
\clint{\tfrac{p-1}{2} \ms{l},\; 1 - \tfrac{p-1}{2} \ms{l}}$,
which is a proper subset of the domain $\clint{0, 1}$ if $p > 1$.
As we will see, the implications of this fact on the approximation quality
of the hierarchical B-spline basis are severe.
In this section, we study the underlying reasons of the restriction, and
we introduce a new hierarchical B-spline basis
that does not suffer from this issue.

\subsection{Approximation Quality of Uniform Hierarchical B-Splines}
\label{sec:321approximation}

\paragraph{Interpolation of polynomials}

\begin{figure}
  \subcaptionbox{%
    Objective function $\objfun$ \emph{\textcolor{C0}{(blue)},}
    interpolant $\fgintp{l}$ \emph{\textcolor{C1}{(red, dashed)},}
    grid points \emph{(dots),} and
    spline interpolation domain $\rspldomain{l}{p}$
    \emph{(thick line).}%
  }[72mm]{%
    \includegraphics{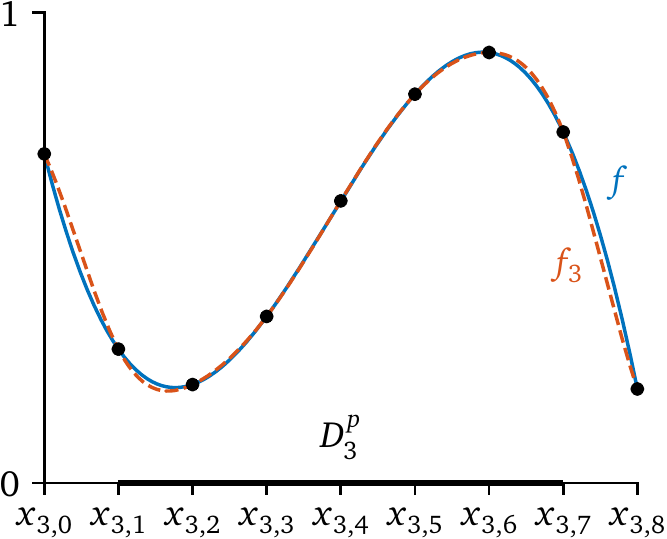}%
  }%
  \hfill%
  \subcaptionbox{%
    Pointwise relative error
    $\abs{(\objfun - \fgintp{l})/\objfun}$ on a logarithmic scale.%
  }[72mm]{%
    \includegraphics{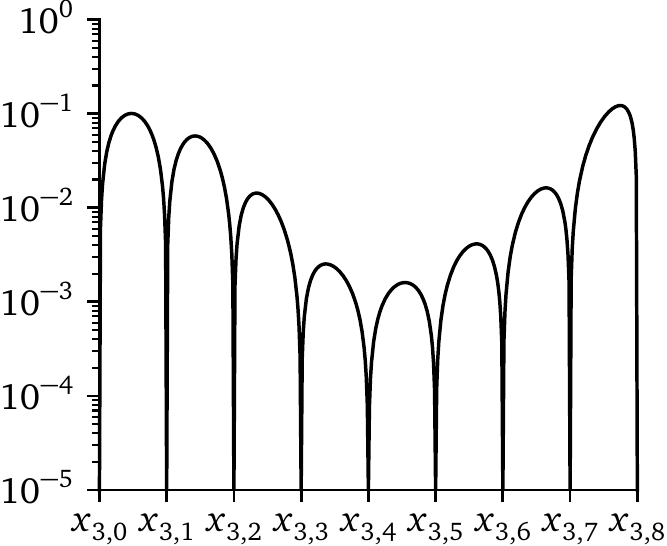}%
  }%
  \caption[%
    Issues when interpolating with uniform hierarchical B-splines%
  ]{%
    Hierarchical cubic B-splines $\bspl{l',i'}{p}$
    ($l' \le l$, $i' \in \hiset{l'}$, $p = 3$)
    fail to replicate a cubic polynomial $\objfun$
    (here: $\objfun(x) \ceq -10.2 x^3 + 14.7 x^2 - 5x + 0.7$)
    \vphantom{$\bspl{l',i'}{p}$}%
    when interpolating on the grid of level $l = 3$.%
    \vphantom{$\bspl{l',i'}{p}$}%
  }%
  \label{fig:nakInterpolation}%
\end{figure}

Splines are a piecewise generalization of polynomials.
Approximation spaces spanned by splines of degree $p$ should
at least contain all polynomials of degree $\le p$.
Unfortunately, this statement is not true for uniform B-splines
$\bspl{l,i}{p}$ as we have defined them in the last section.
A counterexample is given in \cref{fig:nakInterpolation},
in which a cubic polynomial $\objfun$ is interpolated with
hierarchical cubic B-splines.
We can clearly see deviations of the interpolant from the polynomial
near the boundary,
where the pointwise relative error exceeds \SI{10}{\percent}.
The oscillations are even visible in the interior of the
spline interpolation domain $\rspldomain{l}{p}$.
Obviously, this phenomenon impairs the approximation quality
for other non-polynomial functions as well.

This issue can be explained as follows:
According to \thmref{cor:hierSplittingBSpline},
we have $\restrictedsplspace{l}{p}
= \bigoplus_{l'=0}^l \restrictspace{\hsbspl{l'}{p}}{\rspldomain{l}{p}}$
with $p = 3$.
Since cubic polynomials are also cubic splines,
it follows $\objfun \in \restrictedsplspace{l}{p}$ and hence
$\objfun \in \bigoplus_{l'=0}^l \restrictspace{\hsbspl{l'}{p}}{\rspldomain{l}{p}}$.
This means that there is a linear combination of hierarchical B-splines
$\bspl{l',i'}{p}$ ($l' \le l$, $i' \in \hiset{l'}$)
that replicates $\objfun$ on the whole domain $\rspldomain{l}{p}$
(not be confused with $\fgintp{l}$ in \cref{fig:nakInterpolation},
which does not replicate $\objfun$ exactly on $\rspldomain{l}{p}$).
However, in general, this interpolant is not equal $\objfun$ outside
$\rspldomain{l}{p}$ (i.e., in $\clint{0, 1} \setminus \rspldomain{l}{p}$),
as \thmref{prop:splineSpace} only holds for $\rspldomain{l}{p}$.
In particular, the interpolant evaluated at $x \in \{0, 1\}$ is not
equal to $\objfun(x)$.
If we now force the additional interpolation conditions in
$\gp{l,0} = 0$ and $\gp{l,2^l} = 1$,
the resulting interpolant $\fgintp{l}$ cannot be the same as the previous
interpolant,
which is why $\objfun$ and $\fgintp{l}$ differ inside $\rspldomain{l}{p}$.

\paragraph{Schoenberg--Whitney conditions}

Formally, the unique existence of an interpolating spline is
described by the \term{Schoenberg--Whitney conditions:}

\begin{proposition}[Schoenberg--Whitney conditions]
  \label{prop:schoenbergWhitneyConditions}
  Let $\knotseq = (\knot{0}, \dotsc, \knot{m+p})$ be a knot sequence
  and $t_0, \dotsc, t_{m-1}$ a sequence of interpolation points with
  $t_0 < \dotsb < t_{m-1}$ and
  $\knot{p} \le t_0 < t_{m-1} \le \knot{m}$.
  Then, there exists a unique interpolating spline
  $\spl = \sum_{k=0}^{m-1} \interpcoeff{k} \nonunifbspl{k,\knotseq}{p}$
  for arbitrary data if and only if
  \begin{equation}
    \knot{k} < t_k < \knot{k+p+1},\quad
    k = 0, \dotsc, m - 1.
  \end{equation}
\end{proposition}

\vspace*{0pt plus 0.3fill}

\begin{proof}
  See \cite{Hoellig13Approximation}.
\end{proof}

\vspace*{0pt plus 1.0fill}

The Schoenberg--Whitney conditions require that the interpolation points
are contained in $\rspldomain{l}{p}$,
which is not the case for $p = 3$ (see \cref{fig:nakInterpolation}),
as $\rspldomain{l}{p}$ does not contain the points $x = 0$ and $x = 1$.
For general degree $p$, the first $\tfrac{p-1}{2}$ and the
last $\tfrac{p-1}{2}$ grid points of level $l$ are
missing from $\rspldomain{l}{p}$,
thus violating the Schoenberg--Whitney conditions.
One possible remedy would be to move these interpolation points inside
$\rspldomain{l}{p}$ without changing the corresponding basis functions
(i.e., the knots stay the same) \cite{Hoellig13Approximation}.
For instance in the cubic case, we could move
$x = 0$ to $x = 1.5 \ms{l}$ and
$x = 1$ to $x = 1 - 1.5 \ms{l}$.
However, with this approach, we would not be able to interpolate
boundary values.
In addition, the condition of the interpolation problem will most likely
worsen if we place interpolation points near the ends of the supports
of the corresponding basis functions.

\pagebreak

\paragraph{Mismatch of dimensions}

To find a solution for this issue,
let $\wholesplspace{l}{p}$ denote the space of all splines of degree $p$
on the grid of level $l$, i.e., the space $\nonunifsplspace{\knotseq}{p}$ with
\begin{equation}
  \label{eq:fullGridKnots}
  \knot{k} \ceq (k - p) \ms{l},\quad
  k = 0, \dotsc, m + p,\quad
  m \ceq 2^l + p.
\end{equation}
We have $\spldomain{\knotseq}{p} = \clint{0, 1}$
for this choice of $\knotseq$.
Hence, the grid points $\gp{l,i}$ ($i = 0, \dotsc, 2^l$) satisfy
the Schoenberg--Whitney conditions for the uniform B-spline basis.
Clearly, the sum $\bigoplus_{l'=0}^l \hsbspl{l'}{p}$
is a subspace of $\wholesplspace{l}{p}$,
but it cannot be equal to $\wholesplspace{l}{p}$ due to
\begin{equation}
  \dim \bigoplus_{l'=0}^l \hsbspl{l'}{p}
  = 2^l + 1
  < 2^l + p
  = m
  = \dim \wholesplspace{l}{p},\quad
  p > 1,
\end{equation}
by \thmref{prop:splineSpace}.
There are too few nodal (and hierarchical) basis functions to
span the whole spline space $\wholesplspace{l}{p}$.

\paragraph{Restriction to spline subspaces}

The key idea is now to impose additional $p - 1$ boundary conditions
on the basis functions to restrict $\wholesplspace{l}{p}$ to a reasonable subspace
with the correct dimension $(2^l + p) - (p - 1) = 2^l + 1$.
``Reasonable'' means that besides this dimension constraint,
two requirements should be met:
First, the Schoenberg-Whitney conditions should be satisfied for
the new subspace and the grid of level $l$.
Second, the new subspace should contain all polynomials of degree $\le p$,
eliminating the issue discussed in \cref{fig:nakInterpolation}.

\subsection{Hierarchical Not-A-Knot B-Splines}
\label{sec:322NAKBSplines}

\paragraph{Not-a-knot conditions}

A suitable subspace can be obtained by incorporating the
so-called \term{not-a-knot boundary conditions} into $\wholesplspace{l}{p}$.
For the cubic case $p = 3$
(for which we need two conditions),
the not-a-knot conditions demand that
for all splines $\spl$ in the subspace,
$\deriv[3]{x}{s}$ is continuous at the first and at the last
interior knot $\gp{l,1} = \ms{l}$ and $\gp{l,2^l-1} = 1 - \ms{l}$
\cite{Hoellig13Approximation}.
This means that
$\gp{l,1}$ and $\gp{l,2^l-1}$ are effectively removed from the
knot sequence, as this is equivalent to requiring that
$\spl$ is a cubic polynomial on $\clint{0, \gp{l,2}}$ and
$\clint{\gp{l,2^l-2}, 1}$ (hence ``not-a-knot'').
For general degree $p$ (for which we need $(p - 1)$ conditions),
we require that the $p$-th derivative $\deriv[p]{x}{s}$
is continuous at the first $\tfrac{p-1}{2}$ and the last $\tfrac{p-1}{2}$
inner grid points
\begin{equation}
  \label{eq:removedNAKKnots}
  \gp{l,i},\quad
  i \in \{1, \dotsc, \tfrac{p-1}{2}\} \cup
  \{2^l - \tfrac{p-1}{2}, \dotsc, 2^l - 1\}.
\end{equation}
This is equivalent to removing these knots
from the knot sequence $\knotseq$,
or, alternatively, to requiring that $\spl$ is a polynomial
of degree $\le p$ on $\clint{0, \gp{l,(p+1)/2}}$ and on
$\clint{\gp{l,2^l-(p+1)/2}, 1}$.

\begin{figure}
  \includegraphics{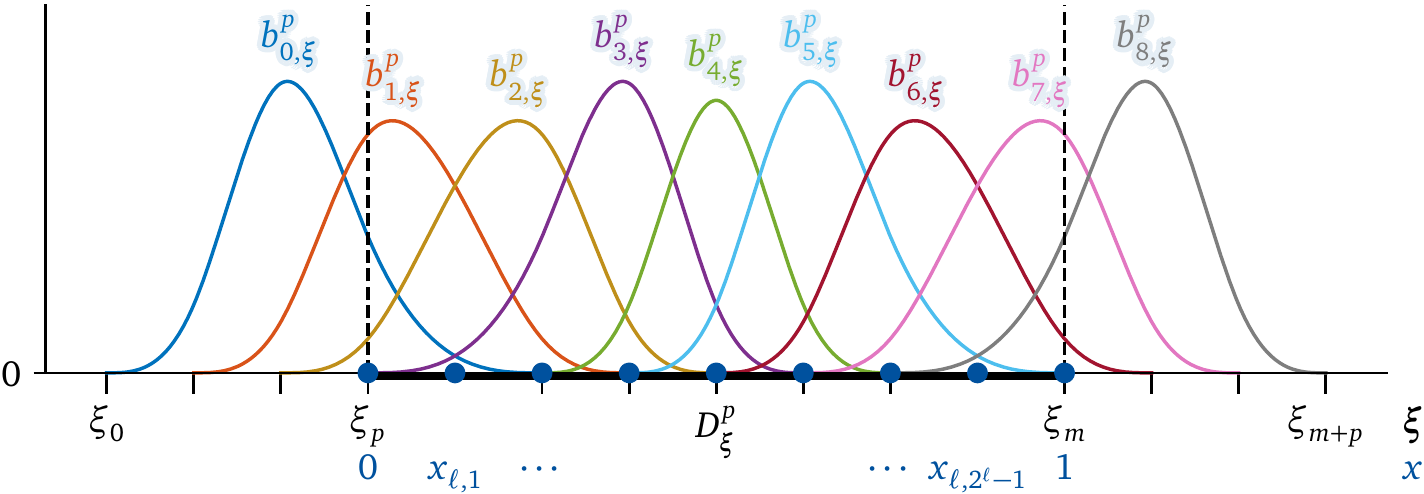}%
  \caption[%
    Nodal not-a-knot B-splines and knot sequence%
  ]{%
    Not-a-knot knot sequence $\nodalknotseq[\nak]{l}{p}$
    \emph{(ticks on horizontal axis)}
    and nodal cubic not-a-knot B-splines ($p = 3$)
    of level $l = 3$.
    In the domain $\clint{0, 1}$ \emph{(delimited by dashed lines),}
    the first $\tfrac{p-1}{2}$ and the last $\tfrac{p-1}{2}$ interior grid points
    of the set of grid points $\fgset{l}$
    \emph{\textcolor{mittelblau}{(blue dots)}}
    have been removed from the set of knots.
    The spline interpolation domain
    $\spldomain{\nodalknotseq[\nak]{l}{p}}{p}$
    \emph{(thick line)}
    equals the whole domain $\clint{0, 1}$.%
  }%
  \label{fig:splineSpaceNotAKnot}%
\end{figure}

\usenotation{zzzznak}
The knot sequence $\nodalknotseq[\nak]{l}{p}$
with not-a-knot boundary conditions is defined as follows:
\begin{subequations}
  \begin{gather}
    \nodalknotseq[\nak]{l}{p}
    \ceq (\nodalknot[\nak]{l,0}{p}, \dotsc,
    \nodalknot[\nak]{l,m+p}{p}),\quad
    m \ceq 2^l + 1,\\
    \nodalknot[\nak]{l,k}{p}
    \ceq
    \begin{cases}
      \gp{l,k-p},&
      k = 0, \dotsc, p,\\
      \gp{l,k-(p+1)/2},&
      k = p + 1, \dotsc, 2^l,\\
      \gp{l,k-1},&
      k = 2^l + 1, \dotsc, 2^l + p + 1.
    \end{cases}
  \end{gather}
\end{subequations}
This knot sequence $\nodalknotseq[\nak]{l}{p}$
can be obtained by removing the knots
given in \eqref{eq:removedNAKKnots} from the
knot sequence \eqref{eq:fullGridKnots} for the full grid of level $l$.
We show $\nodalknotseq[\nak]{l}{p}$ and the corresponding B-spline functions
in \cref{fig:splineSpaceNotAKnot}.
The resulting spline space
\begin{equation}
  \naksplspace{l}{p}
  \ceq \nonunifsplspace{\nodalknotseq[\nak]{l}{p}}{p}
\end{equation}
is a subspace
of $\wholesplspace{l}{p}$ with the desired dimensionality:
\begin{equation}
  \dim \bigoplus_{l'=0}^l \hsbspl{l'}{p}
  = 2^l + 1
  = \dim \naksplspace{l}{p}.
\end{equation}
The space $\naksplspace{l}{p}$ satisfies our two requirements:
\vspace{0.1em}%
First, the spline interpolation domain
$\spldomain{\nodalknotseq[\nak]{l}{p}}{p}
= \clint{\nodalknot[\nak]{l,p}{p}, \nodalknot[\nak]{l,m}{p}}$
of $\naksplspace{l}{p}$ equals the whole unit interval $\clint{0, 1}$.
\vspace{-0.3em}%
This means that the Schoenberg--Whitney conditions are satisfied
\vspace{0.1em}%
for $\naksplspace{l}{p}$, since all interpolation points
(grid points) are contained in
$\spldomain{\nodalknotseq[\nak]{l}{p}}{p} = \clint{0, 1}$.
\vspace{-0.3em}%
Second, $\naksplspace{l}{p}$ still contains all polynomials of
degree $\le p$, as we have only removed knots compared to
$\wholesplspace{l}{p}$.

\vspace*{\fill}

However, $\bigoplus_{l'=0}^l \hsbspl{l'}{p}$ is not a subspace of
$\naksplspace{l}{p}$ anymore,
since the hierarchical basis functions $\bspl{l',i'}{p}$ are not
not-a-knot splines (due to the additional knots that we removed from
$\naksplspace{l}{p}$).
For this reason,
we have to incorporate the not-a-knot boundary conditions into the
hierarchical basis.

\vspace*{\fill}

Before defining the new hierarchical basis functions,
we make two additional observations.
First, $\nodalknotseq[\nak]{l}{p}$ coincides with the
uniform knot sequence $\nodalknotseq{l}{p}$ of \thmref{cor:nodalBSplineSpace}
for the piecewise linear case of $p = 1$.
This is intuitively clear:
For this case,
we do not need to remove any knots as the hierarchical splitting already
holds for the full domain by \cref{cor:hierSplittingHatMV}.
Second, the removal of the knots in \eqref{eq:removedNAKKnots}
is only possible if $p + 1 \le 2^l$,
which is equivalent to $l \ge \ceil{\log_2(p+1)}$.
For coarser levels,
there are not enough interior knots that could be removed.
Without any special treatment,
we would not be able to obtain enough basis functions to span the spline space.

\pagebreak

\paragraph{Definition of hierarchical not-a-knot B-splines}

\begin{figure}
  \subcaptionbox{%
    Nodal not-a-knot B-splines
    $\bspl[\nak]{l,i}{p}$ ($i \in \hiset{l}$)
    and grid points $\gp{l,i}$ \emph{(dots).}%
  }[67.5mm]{%
    \includegraphics{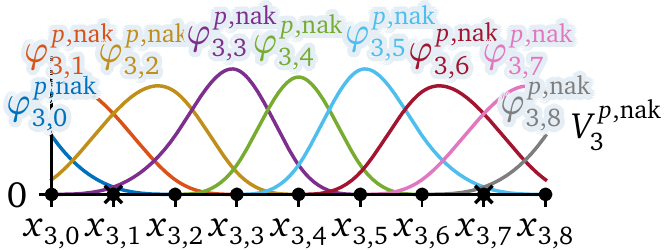}%
  }%
  \hfill%
  \begin{tikzpicture}
    \draw[decorate,decoration={brace,aspect=0.145}] (0,-8.5) -- (0,0);
    \node[anchor=east,inner sep=0mm] at (-0.15,-7.265) {$= \bigoplus$};
  \end{tikzpicture}%
  \hfill%
  \subcaptionbox{%
    Hierarchical not-a-knot B-splines
    $\bspl[\nak]{l',i'}{p}$ ($l' \le l$, $i' \in \hiset{l'}$)
    and grid points $\gp{l',i'}$ \emph{(dots).}%
  }[71mm]{%
    \includegraphics{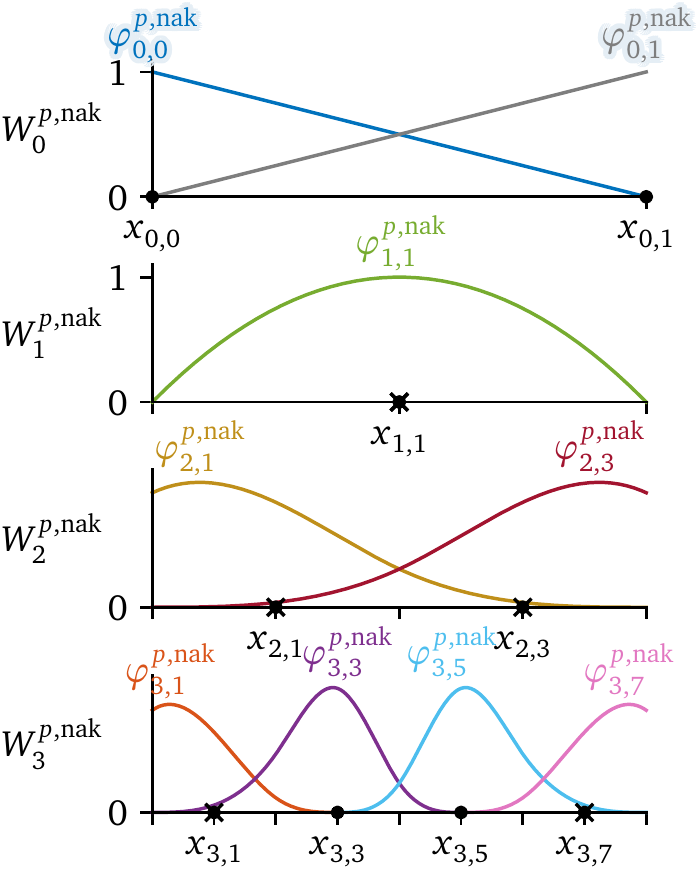}%
  }%
  \caption[%
    Nodal and hierarchical not-a-knot B-splines%
  ]{%
    Univariate nodal and hierarchical cubic not-a-knot B-splines ($p = 3$)
    up to level $l = 3$.
    The nodal space $\nsbspl[\nak]{l}{p}$,
    which coincides with the not-a-knot spline space $\naksplspace{l}{p}$,
    decomposes into the direct sum
    of the hierarchical subspaces $\hsbspl[\nak]{l'}{p}$ ($l' \le l$).
    The knots of each level $l'$ are given by removing the
    first $\tfrac{p-1}{2}$ and last $\tfrac{p-1}{2}$
    inner points \emph{(crosses)}
    from the set of grid points $\gp{l',i'}$
    ($i' = 0, \dotsc, 2^{l'}$).%
  }%
  \label{fig:notAKnotBSpline}%
\end{figure}

For the definition of \term{hierarchical not-a-knot B-splines}
$\bspl[\nak]{l,i}{p}$ based on \thmref{def:nonUniformBSpline},
we use global Lagrange polynomials for the coarser levels:
\begin{subequations}
  \label{eq:hierarchicalNotAKnotBSpline}
  \begin{gather}
    \bspl[\nak]{l,i}{p}
    \ceq
    \begin{cases}
      \lagrangepoly{l,i},&
      l < \ceil{\log_2(p+1)},\\
      \nonunifbspl{i,\nodalknotseq[\nak]{l}{p}}{p},&
      l \ge \ceil{\log_2(p+1)},
    \end{cases}\quad
    l \in \natz,\quad
    i \in \hiset{l},\\
    \lagrangepoly{l,i}\colon \clint{0, 1} \to \real,\quad
    \lagrangepoly{l,i}(x)
    \ceq \!\!\prod_{\substack{i'=0,\dotsc,2^l\\i'\not=i}}
    \frac{x - \gp{l,i'}}{\gp{l,i} - \gp{l,i'}}.
  \end{gather}
\end{subequations}
The hierarchical not-a-knot B-spline basis is shown for the
cubic case $p = 3$ in \cref{fig:notAKnotBSpline}.
The function $\lagrangepoly{l,i}$ is the $i$-th
\term{Lagrange polynomial} of level $l$, that is,
the unique polynomial of degree $\le 2^l$ that interpolates the data
$\{(\gp{l,i'}, \kronecker{i}{i'}) \mid i' = 0, \dotsc, 2^l\}$.
Since its degree $\deg \lagrangepoly{l,i}$ is bounded by $2^l$,
we have $\deg \lagrangepoly{l,i} < p + 1$,
as the Lagrange polynomials are employed only when
$l < \ceil{\log_2(p+1)}$.
Due to $2^l$ even (when $l \ge 1$) and $p$ odd,
we can conclude from $\deg \lagrangepoly{l,i} \le 2^l \le p$ that actually
$\deg \lagrangepoly{l,i} \le 2^l < p$
(for $p > 1$; for $p = 1$, the case $l = 0$ is the exception).

\vspace*{\fill}

The motivation for using Lagrange polynomials for coarse levels
is that they form a basis of the polynomial space
and that they can be implemented and calculated quickly.
However, the specific choice of basis functions for the levels
$l < \ceil{\log_2(p + 1)}$ is arbitrary,
as long as these functions are linearly independent
(of each other and of the ``true'' not-a-knot B-splines
$\bspl[\nak]{l,i}{p}$, $l \ge \ceil{\log_2(p+1)}$)
and contained in the space $\naksplspace{l}{p}$.

\pagebreak

\paragraph{Implementation}

Note that in each level $l \ge \ceil{\log_2(p+1)}$,
only the first $\tfrac{p+1}{2}$
(indices $i = 1, 3, \dotsc, p$)
and the last $\tfrac{p+1}{2}$
(indices $i = 2^l - p, 2^l - p + 2, \dotsc, 2^l - 1$)
hierarchical basis functions $\bspl[\nak]{l,i}{p}$
differ from $\bspl{l,i}{p}$,
i.e., we have
\begin{equation}
  \bspl[\nak]{l,i}{p} = \bspl{l,i}{p},\quad
  i = p + 2,\; p + 4,\; \dotsc,\; 2^l - p - 2.
\end{equation}

\vspace*{\fill}

\noindent
This means that we can reuse uniform B-spline code
for the inner functions.
Due to $\bspl[\nak]{l,i}{p}(x) = \bspl[\nak]{l,2^l-i}{p}(1-x)$
(because of the symmetry of $\nodalknotseq[\nak]{l}{p}$),
we only have to reimplement $\tfrac{p+1}{2}$ not-a-knot B-splines per level $l$.
As $\bspl[\nak]{l,i}{p}$ and $\bspl[\nak]{l+1,i}{p}$
use the same knots up to an affine transformation for $l$ large enough
($l \ge 3$ suffices for $p = 3$),
only a number of special cases for coarse levels $l$ must be implemented.
In other words, the not-a-knot approach is ``minimally invasive''
with respect to an implementation that already uses uniform B-splines.

\pagebreak

\paragraph{Hierarchical splitting}

The main benefit of the hierarchical not-a-knot B-spline basis
is the validity of the hierarchical splitting.
As usual, we define $\nsbspl[\nak]{l}{p}$ and $\hsbspl[\nak]{l}{p}$
as the nodal and the hierarchical not-a-knot subspace of level $l$,
respectively.

\begin{proposition}[%
  univariate hierarchical splitting for not-a-knot B-splines%
]
  \label{prop:hierSplittingNAKBSplineUV}
  The hierarchical splitting \eqref{eq:hierSplittingUV}
  holds for the hierarchical not-a-knot B-spline basis:
  \begin{equation}
    \naksplspace{l}{p}
    = \nsbspl[\nak]{l}{p}
    = \bigoplus_{l'=0}^l \hsbspl[\nak]{l'}{p},
  \end{equation}
  where for $l < \ceil{\log_2(p+1)}$, $\naksplspace{l}{p}$
  is defined as the space $\polyspace{2^l}$ of polynomials of degree
  $\le 2^l$ on $\clint{0, 1}$.
  (For $l \ge \ceil{\log_2(p+1)}$,
  $\naksplspace{l}{p}$ is the not-a-knot spline space.)
\end{proposition}

\begin{proof}
  For $l < \ceil{\log_2(p+1)}$, all
  three spaces coincide with $\polyspace{2^l}$ and nothing is to prove.
  
  For $l \ge \ceil{\log_2(p+1)}$,
  we check the two conditions of \thmref{lemma:hierSplittingUV}.
  First, the hierarchical subspace $\hsbspl[\nak]{l'}{p}$ ($l' \le l$)
  is a subspace of $\naksplspace{l}{p} = \nsbspl[\nak]{l}{p}$.
  This is a conclusion of \thmref{prop:splineSpace}, as
  every function $\bspl[\nak]{l',i'}{p}$ ($i' \in \hiset{l'}$)
  is continuous on $\clint{0, 1}$, a polynomial on every knot interval of
  $\nodalknotseq[\nak]{l}{p}$, and at the knots themselves
  at least $(p - 1)$ times continuously differentiable.
  
  Second, the hierarchical functions $\bspl[\nak]{l',i'}{p}$
  ($l' \le l$, $i' \in \hiset{l'}$) are linearly independent.
  This can be shown similarly to the proof of
  \thmref{prop:hierBSplineLinearlyIndependent}.
  The linear independence of the Lagrange polynomials
  can be checked by inserting grid points into a zero linear combination.
  The linear combination collapses and only one term remains,
  the coefficient corresponding to the grid point.
  Hence, all coefficients must vanish.
\end{proof}

\vspace*{1em}

\begin{corollary}[%
  multivariate hierarchical splitting for not-a-knot B-splines%
]
  \label{cor:hierSplittingNAKBSplineMV}
  It holds
  \begin{equation}
    \naksplspace{\*l}{\*p}
    = \nsbspl[\nak]{\*l}{\*p}
    = \bigoplus_{\*l'=\*0}^\*l \hsbspl[\nak]{\*l'}{\*p},
  \end{equation}
  where $\naksplspace{\*l}{\*p}$ is the
  tensor product space of $\naksplspace{l_t}{p_t}$
  ($t = 1, \dotsc, d$) as defined in \cref{prop:hierSplittingNAKBSplineUV}.
\end{corollary}

\begin{proof}
  Follows directly from \thmref{prop:splittingUVToMV}.
\end{proof}

\paragraph{Sparse grids with not-a-knot B-splines}

Regular sparse grid spaces using the new hierarchical not-a-knot basis
are defined analogously to the uniform case, i.e.,
\begin{equation}
  \label{eq:sparseGridRegularNAK}
  \regsgspace[\*p,\nak]{n}{d}
  \ceq \bigoplus_{\normone{\*l} \le n} \hsbspl[\nak]{\*l}{\*p}.
\end{equation}
If the level $n$ is large enough, then $\regsgspace[\*p,\nak]{n}{d}$
contains the space $\polyspace{\*p}$ of all $d$-variate polynomials of
coordinate degree $\le \*p$ on $\clint{\*0, \*1}$
(i.e., functions $\objfun\colon \clint{0, 1} \to \real$ of the form
$\objfun(\*x) \ceq \sum_{\*q=\*0}^\*p \interpcoeff{\*q} \prod_{t=1}^d x_t^{q_t}$
with $\interpcoeff{\*q} \in \real$).
This means that in contrast to the uniform B-spline basis,
hierarchical not-a-knot B-splines on sparse grids are able to
replicate global polynomials on $\clint{\*0, \*1}$:

\begin{shortcorollary}[%
  sparse grid with not-a-knot B-splines contains polynomials%
]
  \label{cor:sparseGridRegularNAKPolynomials}
  If $n \ge \normone{\ceil{\veclog_\*2(\*p + \*1)}}$,
  then $\polyspace{\*p} \subset \regsgspace[\*p,\nak]{n}{d}$.
\end{shortcorollary}

\begin{proof}
  Let $\*l \ceq \ceil{\veclog_\*2(\*p + \*1)}$ and $n \ge \normone{\*l}$.
  By \cref{cor:hierSplittingNAKBSplineMV}, we have
  $\bigoplus_{\*l'=\*0}^{\*l} \hsbspl[\nak]{\*l'}{\*p} = \naksplspace{\*l}{\*p}$.
  In addition, all $\*l' \in \natz^d$ with $\*l' \le \*l$ satisfy
  $\normone{\*l'} \le n$ and thus,
  $\bigoplus_{\*l'=\*0}^{\*l} \hsbspl[\nak]{\*l'}{\*p} \subset
  \regsgspace[\*p,\nak]{n}{d}$ by \eqref{eq:sparseGridRegularNAK}.
  We conclude
  $\polyspace{\*p} \subset \naksplspace{\*l}{\*p} \subset
  \regsgspace[\*p,\nak]{n}{d}$, which is the asserted claim.
\end{proof}

\subsection{Modified and Non-Uniform Hierarchical Not-A-Knot B-Splines}
\label{sec:323modifiedNAKBSplines}

\paragraph{Modified hierarchical not-a-knot B-splines}

As for uniform and Clenshaw--Curtis B-splines
(\cref{sec:31standardBSplines}),
it is possible to define a modified version of the
hierarchical not-a-knot B-spline basis to obtain
``reasonable'' boundary values without boundary points.
However, we cannot use \thmref{lemma:marsden} similarly to
\eqref{eq:modifiedBSplineConstruction}:
Due to the removal of knots, there is only a single
not-a-knot B-spline $\bspl[\nak]{l,0}{p}$ left of
$\bspl[\nak]{l,1}{p}$.
B-splines $\bspl[\nak]{l,i}{p}$ with index $i < 0$
would vanish on $\clint{0, 1}$.

\begin{figure}
  \subcaptionbox{%
    $\bspl[\nak]{l',i'}{p}$,
    $\bspl[\nak,\modified]{l',i'}{p}$
    \emph{(dashed),} and $\gp{l',i'}$ \emph{(dots).}%
    \label{fig:modifiedNotAKnotBSpline}%
  }[73mm]{%
    \includegraphics{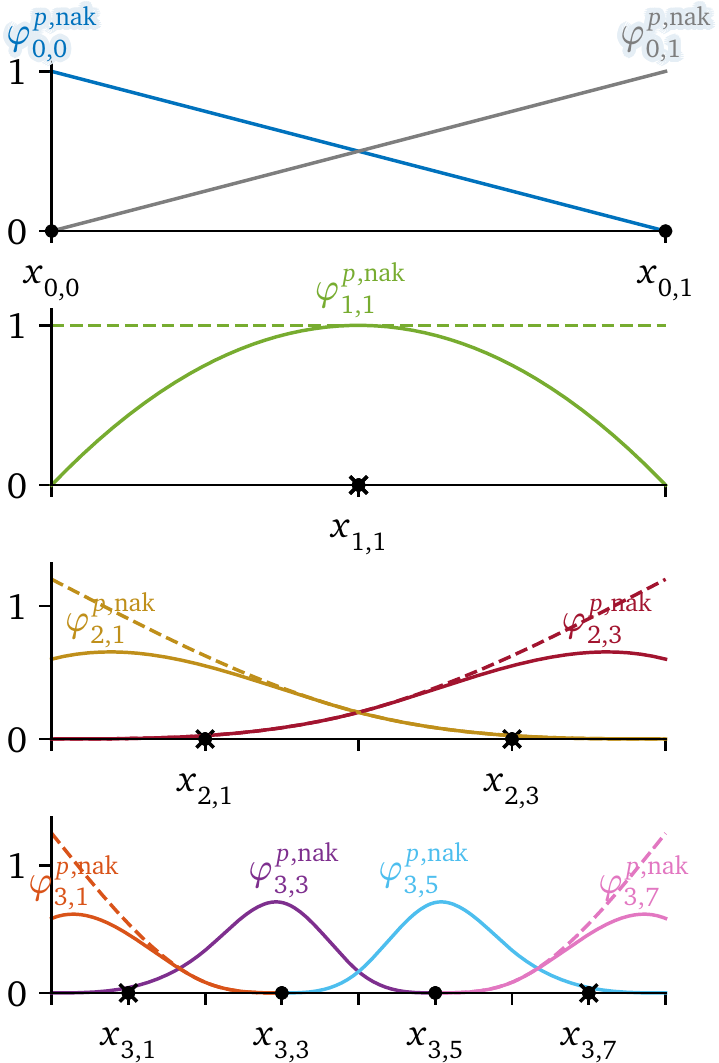}%
  }%
  \hfill
  \subcaptionbox{%
    $\bspl[\cc,\nak]{l',i'}{p}$,
    $\bspl[\cc,\nak,\modified]{l',i'}{p}$
    \emph{(dashed),} and $\ccgp{l',i'}$ \emph{(dots).}%
    \label{fig:clenshawCurtisNotAKnotBSpline}%
  }[76mm]{%
    \includegraphics{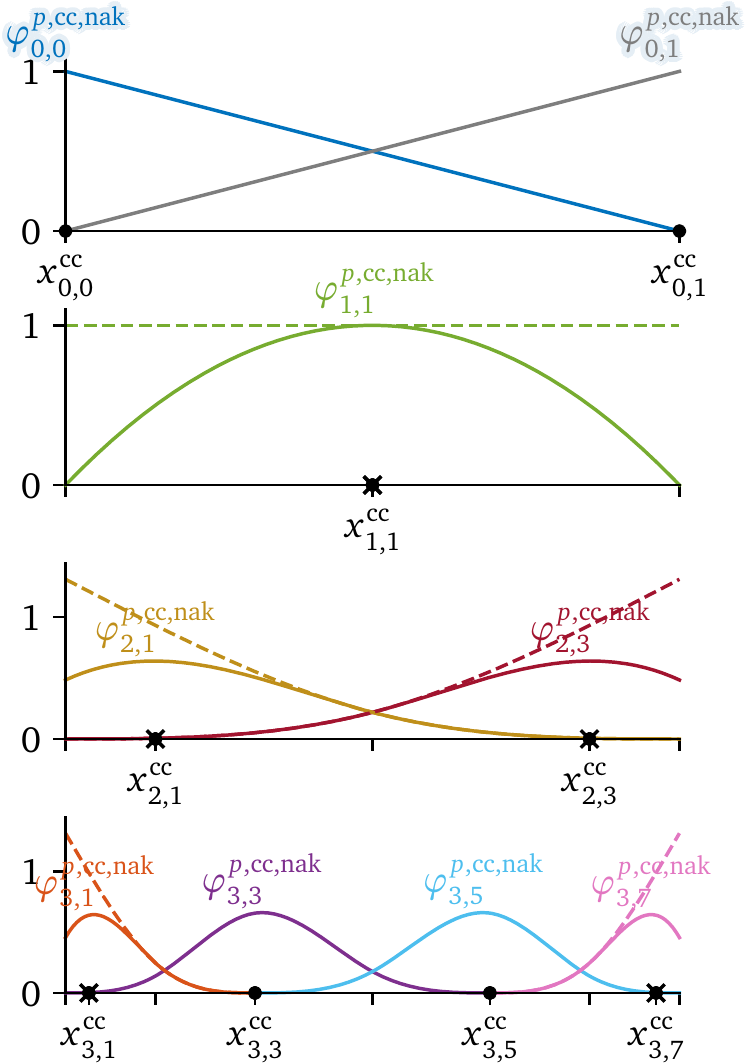}%
  }%
  \caption[%
    Comparison of hierarchical not-a-knot B-splines%
  ]{%
    Comparison of uniform \emph{(left)} and
    Clenshaw--Curtis \emph{(right)} hierarchical cubic not-a-knot
    B-splines $\bspl[\nak]{l',i'}{p}$ and
    $\bspl[\cc,\nak]{l',i'}{p}$
    ($l ' \le l$, $i' \in \hiset{l'}$, $p = 3$) up to level $l = 3$
    together with the respective modified versions
    $\bspl[\nak,\modified]{l',i'}{p}$ and
    $\bspl[\cc,\nak,\modified]{l',i'}{p}$
    \emph{(dashed).}
    The knots of each level $l'$ are given by removing the
    first $\tfrac{p-1}{2}$ and last $\tfrac{p-1}{2}$
    inner points \emph{(crosses)}
    from the set of grid points $\gp{l',i'}$ or
    $\ccgp{l',i'}$
    ($i' = 0, \dotsc, 2^{l'}$), respectively.%
  }%
  \label{fig:uniformAndClenshawCurtisNotAKnotBSpline}%
\end{figure}

While we are therefore not able to construct modified functions
whose second derivative vanishes in a neighborhood of $x = 0$,
we can use $\bspl[\nak]{l,0}{p}$ to let the
second derivative vanish in $x = 0$ itself:
\begin{equation}
  \label{eq:modifiedNotAKnotBSpline}
  \bspl[\nak,\modified]{l,i}{p}(x)
  \ceq
  \begin{cases}
    1,&
    l = 1,\quad i = 1,\\
    \bspl[\nak]{l,1}{p}(x)
    - \dfrac{\deriv[2]{x}{\bspl[\nak]{l,1}{p}}(0)}%
    {\deriv[2]{x}{\bspl[\nak]{l,0}{p}}(0)}
    \bspl[\nak]{l,0}{p}(x),&
    l \ge 2,\quad i = 1,\\
    \bspl[\nak]{l,i}{p}(x),&
    l \ge 2,\quad i \in \hiset{l} \setminus \{1, 2^l - 1\},\\
    \bspl[\nak,\modified]{l,1}{p}(1 - x),&
    l \ge 2,\quad i = 2^l - 1.
  \end{cases}
  \hspace*{-2mm}
\end{equation}
The resulting modified hierarchical not-a-knot B-spline basis
$\bspl[\nak,\modified]{l,i}{p}$ is shown with dashed lines
in \cref{fig:modifiedNotAKnotBSpline}.
As before, we have to implement $\bspl[\nak,\modified]{l,1}{p}$
only for a single level $l$, as modified functions of higher levels
are the same up to an affine parameter transformation.
Note again that for $p \ge 5$, we would have to modify additional
interior B-splines as the interior of their support then extends to the
boundary of $\clint{0, 1}$.
We refrain from doing so to keep the definition
\eqref{eq:modifiedNotAKnotBSpline} simple.

\paragraph{Non-uniform hierarchical not-a-knot B-splines}

The not-a-knot construction is completely independent of the
distribution of the grid points at hand.
Consequently, we can define hierarchical not-a-knot B-splines
for non-uniform distributions.
For instance, to define not-a-knot B-splines for
Chebyshev points (see \cref{sec:314nonUniform}),
we first specify the knot sequence as
\begin{subequations}
  \begin{gather}
    \nodalknotseq[\cc,\nak]{l}{p}
    \ceq (\nodalknot[\cc,\nak]{l,0}{p}, \dotsc,
    \nodalknot[\cc,\nak]{l,m+p}{p}),\quad
    m \ceq 2^l + 1,\\
    \nodalknot[\cc,\nak]{l,k}{p}
    \ceq
    \begin{cases}
      \ccgp{l,k-p},&
      k = 0, \dotsc, p,\\
      \ccgp{l,k-(p+1)/2},&
      k = p + 1, \dotsc, 2^l,\\
      \ccgp{l,k-1},&
      k = 2^l + 1, \dotsc, 2^l + p + 1,
    \end{cases}
  \end{gather}
\end{subequations}
and then define hierarchical not-a-knot Clenshaw--Curtis B-splines as
\begin{subequations}
  \begin{gather}
    \bspl[\cc,\nak]{l,i}{p}
    \ceq
    \begin{cases}
      \lagrangepoly[\cc]{l,i},&
      l < \ceil{\log_2(p+1)},\\
      \nonunifbspl{i,\nodalknotseq[\cc,\nak]{l}{p}}{p},&
      l \ge \ceil{\log_2(p+1)},
    \end{cases}\quad
    l \in \natz,\quad
    i \in \hiset{l},\\
    \lagrangepoly[\cc]{l,i}\colon \clint{0, 1} \to \real,\quad
    \lagrangepoly[\cc]{l,i}(x)
    \ceq \!\!\prod_{\substack{i'=0,\dotsc,2^l\\i'\not=i}}
    \frac{x - \ccgp{l,i'}}%
    {\ccgp{l,i} - \ccgp{l,i'}}.
  \end{gather}
\end{subequations}

This definition can even be combined with the modification
of hierarchical not-a-knot B-splines as discussed above.
We can use exactly the same approach as in
\eqref{eq:modifiedNotAKnotBSpline}, if we replace the
not-a-knot basis functions with their non-uniform not-a-knot counterparts
(not-a-knot Clenshaw--Curtis B-splines in the above example).
The hierarchical not-a-knot Clenshaw--Curtis B-spline basis of
cubic degree and the corresponding modified functions are shown in
\cref{fig:clenshawCurtisNotAKnotBSpline}.

\subsection{Other Approaches to Incorporate Boundary Conditions}
\label{sec:324naturalBoundary}

Not-a-knot boundary conditions are not the only approach
to obtain a subspace of $\wholesplspace{l}{p}$ with the right dimension $2^l - 1$.
\term{Natural boundary conditions} are another possibility,
which we want to discuss briefly.
In the cubic case, for which they are usually formulated
\cite{Hoellig13Approximation},
these boundary conditions require that the
second derivatives $\deriv[2]{x}{\basis{l,i}}$ of the
basis functions vanish at the boundary $x \in \{0, 1\}$.
To obtain the necessary number of $p - 1$ constraints also
for higher degrees $p$,
we require that all derivatives
$\deriv[q]{x}{\basis{l,i}}$ of order
$q = 2, \dotsc, \tfrac{p+1}{2}$ vanish at $x \in \{0, 1\}$.

Consequently, we can define hierarchical natural B-splines as
\begin{subequations}
  \begin{gather}
    \bspl[\ntrl]{l,i}{p}(x)
    \ceq
    \begin{cases}
      \lagrangepoly{0,i}(x),&
      l = 0,\\
      \bspl{l,i}{p} +
      \sum_{j \in J_i^{p,\ntrl}} c_{l,i,j} \bspl{l,j}{p},&
      l \ge 1,
    \end{cases}\quad
    l \in \natz,\quad
    i \in \hiset{l},\\
    J_i^{p,\ntrl}
    \ceq \{i-\tfrac{p-1}{2}, \dotsc, i-1\} \cup
    \{i+1, \dotsc, i+\tfrac{p-1}{2}\},
  \end{gather}
\end{subequations}
where the coefficients $c_{l,i,j} \in \real$ are chosen such that
the natural boundary conditions are satisfied:
\begin{equation}
  \deriv[q]{x}{\bspl[\ntrl]{l,i}{p}}(x)
  = 0,\quad
  l \ge 1,\quad
  i \in \hiset{l},\quad
  q = 2, \dotsc, \tfrac{p+1}{2},\quad
  x \in \{0, 1\}.
\end{equation}
The first half of the coefficients $c_{l,i,j} \in \real$
($j < i$) vanishes if the interior of the support of $\bspl{l,i}{p}$
does not contain $x = 0$
(i.e., $i \ge \tfrac{p+1}{2}$).
The second half of the coefficients~($j > i$) vanishes analogously
if $1 \notin \interiorsupp \bspl{l,i}{p} \iff i \le 2^l - \tfrac{p+1}{2}$.
This means that
only the first $\floor{\tfrac{p+1}{4}}$ and the last $\floor{\tfrac{p+1}{4}}$
hierarchical functions have to be altered in each level.

\begin{SCfigure}
  \includegraphics{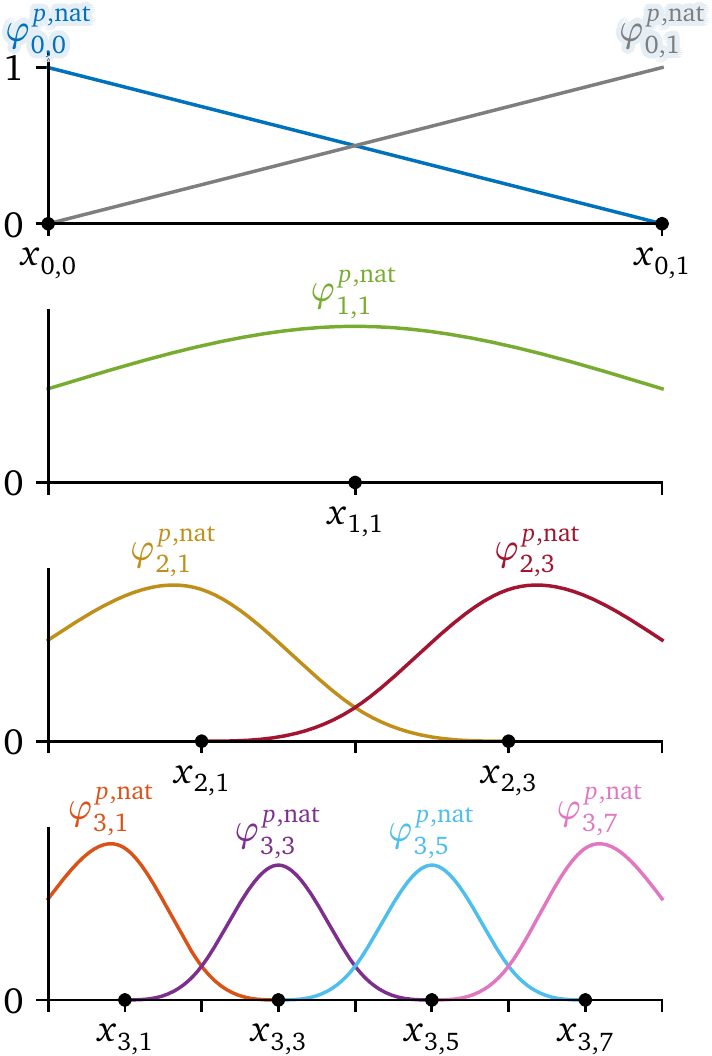}%
  \caption[%
    Hierarchical natural B-splines%
  ]{%
    Hierarchical cubic natural B-splines
    \vspace{-0.1em}%
    $\bspl[\ntrl]{l',i'}{p}$
    ($l' \le l$, $i' \in \hiset{l'}$, $p = 3$) and
    \vspace{0.05em}%
    grid points $\gp{l',i'}$ \emph{(dots)} up to level $l = 3$.%
  }%
  \label{fig:naturalBSpline}%
\end{SCfigure}

\Cref{fig:naturalBSpline} shows the hierarchical natural spline basis.
The main disadvantage of natural boundary conditions is that
we are not able to replicate arbitrary polynomials exactly on $\clint{0, 1}$
with this approach.
Only polynomials that satisfy natural boundary conditions themselves
(linear polynomials for example)
can be replicated exactly.
For this reason, we do not further consider this basis in the
rest of the thesis.

\cleardoublepage

  \setdictum{%
  Who are you? How did you get in my house?%
}{%
  Donald E. Knuth about one-based array indices in algorithms
  (according to \texttt{xkcd}\footnotemark)%
}

\longchapter{%
  Algorithms for B-Splines on Sparse Grids%
}{%
  Algorithms for B-Splines on\texorpdfstring{\\}{ }Sparse Grids%
}{%
  Algorithms for B-Splines on Sparse Grids%
}
\footnotetext{\url{https://xkcd.com/163/}}
\label{chap:40algorithms}

\initial{0.6em}{L}{ittle is known}
about the algorithmic challenges
that hierarchical bases of B-spline type (or even of general type) pose
on sparse grids.
In general, we are not able to directly apply the sparse grid algorithms
that were designed for hat functions $\bspl{l,i}{1}$.
Hence, we have to generalize these algorithms
to higher-order B-splines $\bspl{l,i}{p}$
or even to arbitrary tensor product basis functions.
The main problem is the larger support
of higher-order B-splines when compared to degree $p = 1$.
The larger support introduces new dependencies between
values of grid points that cannot be resolved with conventional algorithms.

This chapter gives an overview of six algorithms for B-splines
on sparse grids.
Two of these algorithms are already known
\multicite{Griebel92Combination,Balder94Adaptive},
while the remaining four are new.
Correctness results are given for every algorithm.
We use hierarchization as the exemplary problem for our algorithms,
but the ideas of the algorithms can be generalized to any linear operator.
Furthermore, most algorithms are not tailored to B-splines $\bspl{l,i}{p}$,
but applicable to general tensor product basis functions $\basis{l,i}$.
Whether an algorithmic approach is feasible for the sparse grid at hand
or not depends on the grid's type:
full grid, dimensionally adaptive sparse grid, or
spatially adaptive sparse grid.
The more assumptions the grid satisfies, the faster and
easier the corresponding algorithms will be.

\Cref{sec:41problem} explains hierarchization as our example problem
and defines the notation used in this chapter.
The remaining four sections treat the three different types of grids:
First, \cref{sec:42fullGrids} deals with full grids to formalize and repeat
the well-known \up.
Second, \cref{sec:43dimAdaptive} focuses on algorithms for
dimensionally adaptive sparse grids.
Third, \cref{sec:44spatAdaptiveBFS} and \cref{sec:45spatAdaptiveUP}
treat arbitrary (spatially adaptive) sparse grids,
which is the most interesting case for us.
\Cref{sec:44spatAdaptiveBFS} employs \bfs for hierarchization,
while \cref{sec:45spatAdaptiveUP} uses the \up\punctfix{.}

This original chapter is the main theoretical contribution of this thesis.
Although the unidirectional principle in \cref{sec:42fullGrids} and
the combination technique in \cref{sec:43dimAdaptive} are well-known,
the presentation with formal proofs of correctness is new.
Parts of the chapter have already been published in scientific papers,
namely \cref{sec:44spatAdaptiveBFS} \cite{Valentin18Fundamental}.
The weakly fundamental splines (\cref{sec:454wfs}) and the
Hermite hierarchization method (\cref{sec:455hermiteHierarchization})
are based on an idea by Dr.\ Stefan Zimmer (University of Stuttgart, Germany).

\section{The Hierarchization Problem}
\label{sec:41problem}

\begin{figure}
  \subcaptionbox{%
    The objective function $\objfun$ is sampled at the grid points
    $\gp{l,i} \in \sgset$ to obtain function values $\fcnval{l,i}$,
    which form the input vector $\vlinin$ for the linear operator
    $\linop = \intpmatinv$.%
  }[46mm]{%
    \includegraphics{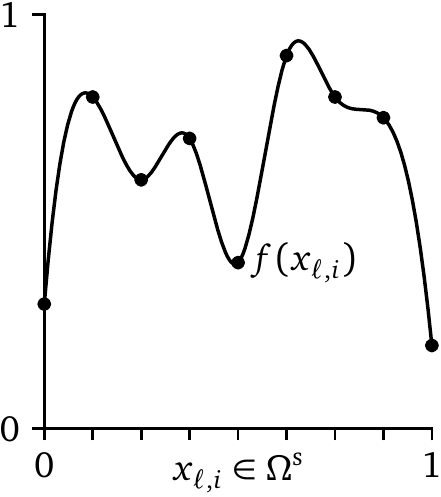}%
  }%
  \hfill%
  \subcaptionbox{%
    The linear operator $\linop$ is applied to $\vlinin$ to obtain
    the output vector $\vlinout$, which contains the hierarchical surpluses
    $\surplus{l,i}$ ($(l,i) \in \liset$).%
  }[46mm]{%
    \includegraphics{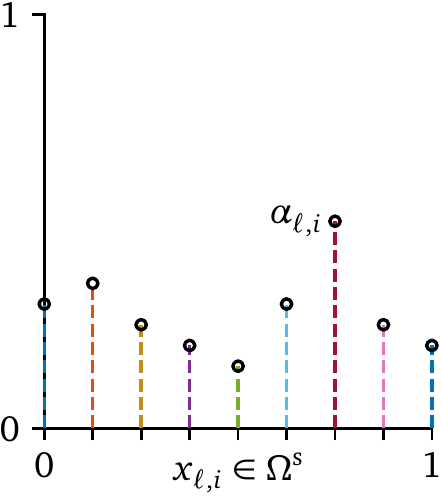}%
  }%
  \hfill%
  \subcaptionbox{%
    The interpolant $\sgintp$ \emph{(black dashed line)}
    is evaluated at $x \in \clint{0, 1}$
    by adding contributions \emph{(black dotted lines)}
    of weighted basis functions $\surplus{l,i} \basis{l,i}$
    \emph{(colored),} obtaining $\sgintp(x)$ \emph{(cross).}%
  }[46mm]{%
    \includegraphics{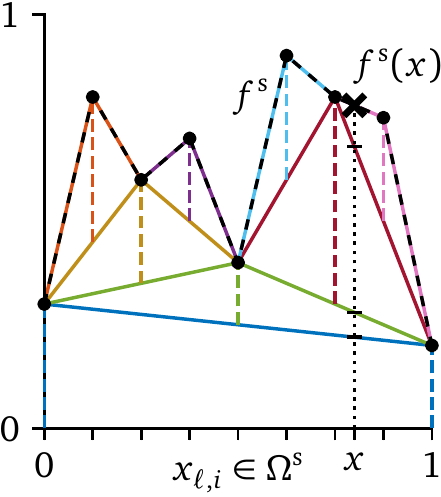}%
  }%
  \caption[%
    Hierarchization of function values and evaluation of interpolant%
  ]{%
    Hierarchization of function values $\fcnval{l,i}$ \emph{(left)}
    to obtain hierarchical surpluses $\surplus{l,i}$ \emph{(center)} and
    evaluation of the resulting interpolant $\sgintp$ \emph{(right),}
    using a univariate grid and the piecewise linear basis as an example.%
  }%
  \label{fig:hierarchization}%
\end{figure}

Let $\sgset \subset \clint{0, 1}^d$ be a general (sparse) grid that
may be spatially adaptive, i.e.,
of the form $\sgset = \{\gp{\*l,\*i} \mid (\*l, \*i) \in \liset\}$,
where $\liset$ is a set of level-index pairs $(\*l, \*i)$ with $\*l \in \natz^d$
and $\*i \in \hiset{\*l}$ such that
$\ngp \ceq \setsize{\sgset} = \setsize{\liset} < \infty$
(see \cref{sec:233spatiallyAdaptiveSG}).
The \term{hierarchization problem} is finding
\term{hierarchical surpluses}
$(\surplus{\*l',\*i'})_{(\*l',\*i') \in \liset} \in \real^{\ngp}$ such that
\begin{equation}
  \label{eq:hierarchizationProblem}
  \largesum{(\*l', \*i') \in \liset} \surplus{\*l',\*i'}
  \basis{\*l',\*i'}(\gp{\*l,\*i}) = \fcnval{\*l,\*i}
  \quad\text{for all}\quad
  (\*l, \*i) \in \liset,
\end{equation}
where $\basis{\*l',\*i'}$ are arbitrary tensor product basis functions and
$(\fcnval{\*l,\*i})_{(\*l,\*i) \in \liset} \in \real^{\ngp}$ is a set of
function values $\objfun(\gp{\*l,\*i})$
at the grid points $\gp{\*l,\*i}$.
This then defines the interpolant $\sgintp$ as
\begin{equation}
  \label{eq:hierarchizationInterpolant}
  \sgintp\colon \clint{\*0, \*1} \to \real,\quad
  \sgintp \ceq
  \largesum{(\*l', \*i') \in \liset} \surplus{\*l',\*i'}
  \basis{\*l',\*i'},
\end{equation}
which interpolates $\objfun$ at the grid points $\gp{\*l,\*i}$ of $\sgset$.
\Cref{fig:hierarchization} shows the process of hierarchizing given
function values and evaluating the resulting interpolant.

We explicitly allow $\basis{\*l',\*i'}$ to be nodal basis functions,
in which case $\*l'$ is constant and
$\sgset$ is a full grid.
Strictly speaking, the problem is then an \term{interpolation problem}
and the $\surplus{\*l',\*i'}$ are \term{interpolation coefficients.}
However, we still apply the terms
``hierarchization'' and ``hierarchical surpluses'' in this case
to keep the terminology consistent.

\paragraph{Hierarchization as a linear operator}

The example of hierarchization can be generalized
to arbitrary linear operators
\begin{equation}
  \linop\colon \real^{\ngp} \to \real^{\ngp},\quad
  \vlinin \mapsto \vlinout = \linop[\vlinin],
\end{equation}
where $\linop$ depends on the grid $\sgset$ at hand.
Input $\vlinin$ and output $\vlinout$ are scalar-valued data%
\begin{equation}
  \vlinin = (\linin{\*l,\*i})_{(\*l,\*i) \in \liset} \in \real^{\ngp},\quad
  \vlinout = (\linout{\*l,\*i})_{(\*l,\*i) \in \liset} \in \real^{\ngp},
\end{equation}
which give one scalar per grid point $\gp{\*l,\*i} \in \sgset$.
For the case of hierarchization,
$\linop$ is the inverse of the \term{interpolation matrix}
$\intpmat \in \real^{\ngp \times \ngp}$:
\begin{subequations}
  \label{eq:hierarchizationSLE}
  \begin{equation}
    \linop = \intpmatinv,\quad
    \intpmat = (\basis{\*l',\*i'}(\gp{\*l,\*i}))_%
    {(\*l,\*i),(\*l',\*i') \in \liset},\quad
    \vlinin = (\fcnval{\*l,\*i})_{(\*l,\*i) \in \liset},\quad
    \vlinout = (\surplus{\*l',\*i'})_{(\*l',\*i') \in \liset}.
    \hspace*{-1mm}
  \end{equation}
  This means that we can determine the $\surplus{\*l',\*i'}$ by solving
  the $\ngp \times \ngp$ system of linear equations
  \begin{equation}
    \vlinout = \linop[\vlinin]
    \quad\iff\quad
    \intpmat \cdot (\surplus{\*l',\*i'})_{(\*l',\*i') \in \liset}
    = (\fcnval{\*l,\*i})_{(\*l,\*i) \in \liset}.
  \end{equation}
\end{subequations}

\paragraph{Complexity of B-spline hierarchization}

\begin{SCfigure}
  \includegraphics{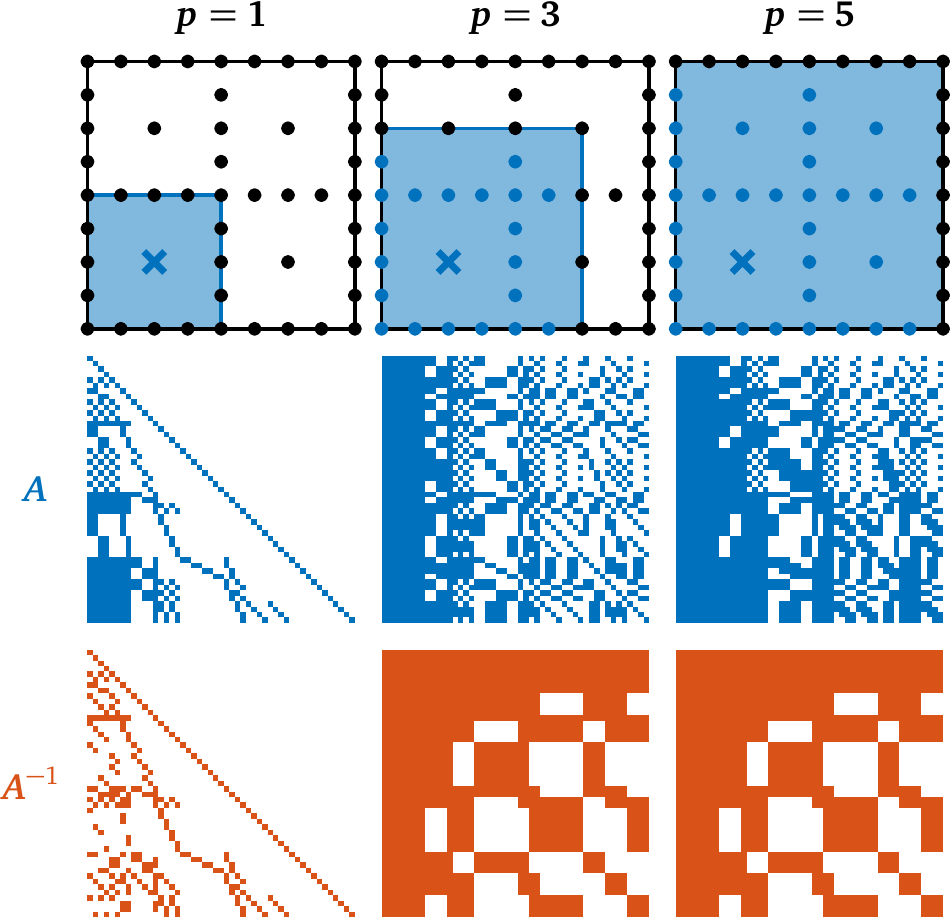}%
  \caption[%
    Density pattern of hierarchization matrices and of their inverses%
  ]{%
    Density pattern
    of the hierarchization matrix $\intpmat$
    \emph{(middle row, \textcolor{C0}{blue})} and
    of its inverse $\intpmatinv$
    \emph{(bottom row, \textcolor{C1}{red})}
    for the regular sparse grid $\coarseregsgset{n}{d}{1}$
    with $n = 4$ and $d = 2$ \emph{(top row)}
    and uniform hierarchical B-splines $\bspl{\*l,\*i}{\*p}$
    for degrees $\*p \in \{\*1, \*3, \*5\}$.
    The \textcolor{C0}{blue areas} in the top row
    show the extent of the support of one specific basis function
    $\bspl{\*l',\*i'}{\*p}$ with $\*l' = (2, 2)$ and $\*i' = (1, 1)$
    (\emph{cross:} corresponding grid point $\gp{\*l',\*i'}$).
    The \textcolor{C0}{blue points} are the grid points at which
    $\bspl{\*l',\*i'}{\*p}$ is non-zero.%
  }%
  \label{fig:matrixDensityPattern}%
\end{SCfigure}

As noted in \cite{Valentin18Fundamental},
hierarchization on sparse grids with hierarchical B-splines
$\bspl{\*l,\*i}{\*p}$ of degree $\*p$
as basis functions $\basis{\*l,\*i}$ is a tedious task.
The corresponding linear system \eqref{eq:hierarchizationSLE} is in general
non-symmetric
(i.e., $\bspl{\*l',\*i'}{\*p}(\gp{\*l,\*i}) \not=
\bspl{\*l,\*i}{\*p}(\gp{\*l',\*i'})$) and densely populated.
This is because the matrix entry in the $(\*l,\*i)$-th row and
$(\*l',\*i')$-th column vanishes if and only if
\begin{equation}
  \gp{\*l,\*i} \notin \interiorsupp \bspl{\*l',\*i'}{\*p}
  \iff
  \exlarge{t = 1, \dotsc, d}{
    \gp{l_t,i_t} \notin
    \opintscaled{
      \gp{l'_t,i'_t} - \tfrac{p_t+1}{2} \ms{l'_t},\,
      \gp{l'_t,i'_t} + \tfrac{p_t+1}{2} \ms{l'_t}
    }
  },
\end{equation}
where $\interiorsupp$ is the interior of the support
\cite{Valentin18Fundamental}.
For coarse levels $\*l'$, the mesh size $\ms{l'_t}$ is large in
every dimension $t$, which implies that $\interiorsupp \bspl{\*l',\*i'}{\*p}$
contains most of the grid points.
In contrast to the hat function case ($\*p = \*1$),
the value of $\surplus{\*l',\*i'}$ depends not only on
$\fcnval{\*l,\*i}$ and the data of the $3^d - 1$ neighboring grid points
on the boundary of $\supp \bspl{\*l',\*i'}{\*1}$,
but potentially on the data of the whole grid.
This is shown in \cref{fig:matrixDensityPattern}:
There are at most $3^d = 9$ non-zero entries in each row of $\intpmatinv$
for $\*p = \*1$ and $d = 2$.
As soon as the B-spline degree is increased,
both $\intpmat$ and $\intpmatinv$ become significantly denser.

This prohibits the use of the \up\punctfix{,}
which we will discuss in the next section,
on sparse grids with hierarchical B-splines.
Consequently, we have to solve the linear system
\eqref{eq:hierarchizationSLE}, which is significantly more time-consuming,
as it takes between $\landauOmega{\ngp^2 d}$ and $\landauO{\ngp^2 (N+d)}$ time
via Gaussian elimination.%
\footnote{%
  $\landauOmega{\ngp^2 d}$ for assembling $\intpmat$ and
  $\landauO{\ngp^3}$ for solving the system.
}
In addition, if we use an explicit solver for the linear system,
we additionally have to store an $\ngp \times \ngp$ matrix in memory.
However, a grid of size $\ngp = \num{116000}$ already exceeds the memory
of a \SI{128}{\gibi\byte} supercomputer node,
if we explicitly store the full matrix in double precision.
In comparison, for the hat function basis,
the \up only requires $\landauO{\ngp d}$ time and $\landauO{\ngp}$ memory.

\paragraph{Notation}

We do not need the hierarchical level-index information $(\*l, \*i)$ in
$\sgset$, $\liset$, $\vlinin$, and $\vlinout$
for most of the considerations in this chapter.
In these cases, we assume that in each dimension $t$, the level-index pairs
$(l_t, i_t)$ ($l_t \in \natz$, $i_t \in \hiset{l_t}$)
are continuously enumerated by a single index $k_t = k_t(l_t, i_t) \in \natz$.
We identify $(\*l, \*i)$ with a single index $\*k$,
whose $t$-th component is given by $k_t(l_t, i_t)$.
Hence,
we regard $\liset$ as a subset $\liset \ceq \{\*k \mid \*x_\*k \in \sgset\}$
of $\natz^d$.
We will switch between the notations whenever appropriate.
All statements that are formulated in the $\*k$ notation are
valid for both the nodal and the hierarchical basis.

\usenotation{kT00}
In the following, $k_t$ denotes the $t$-th component of a $d$-vector $\*k$
as usual.
\usenotation{kT10}
With $\*k_{-t}$, we denote the $(d-1)$-vector that is obtained from $\*k$
by omitting the $t$-th component,
i.e., $\*k_{-t} \ceq (k_1, \dotsc, k_{t-1}, k_{t+1}, \dotsc, k_d)$.
\usenotation{kT20}
For a $j$-tuple $T = (t_1, \dotsc, t_j) \in \{1, \dotsc, d\}^j$,
we define $\*k_T$ to be the $j$-vector $(k_{t_1}, \dotsc, k_{t_j})$
that only contains the entries of the dimensions listed in $T$.
\usenotation{kT30}
Accordingly, $\*k_{-T}$ is defined as the $(d-j)$-vector
that contains the entries of the remaining dimensions
(sorted by the dimension $t$).
We define $\*k_{\range{a}{b}} \ceq (k_a, k_{a+1}, \dotsc, k_b)$
as an indexing shortcut ($a \le b$).

\section{Hierarchization on Full Grids (Unidirectional Principle)}
\label{sec:42fullGrids}

\begin{figure}
  \includegraphics{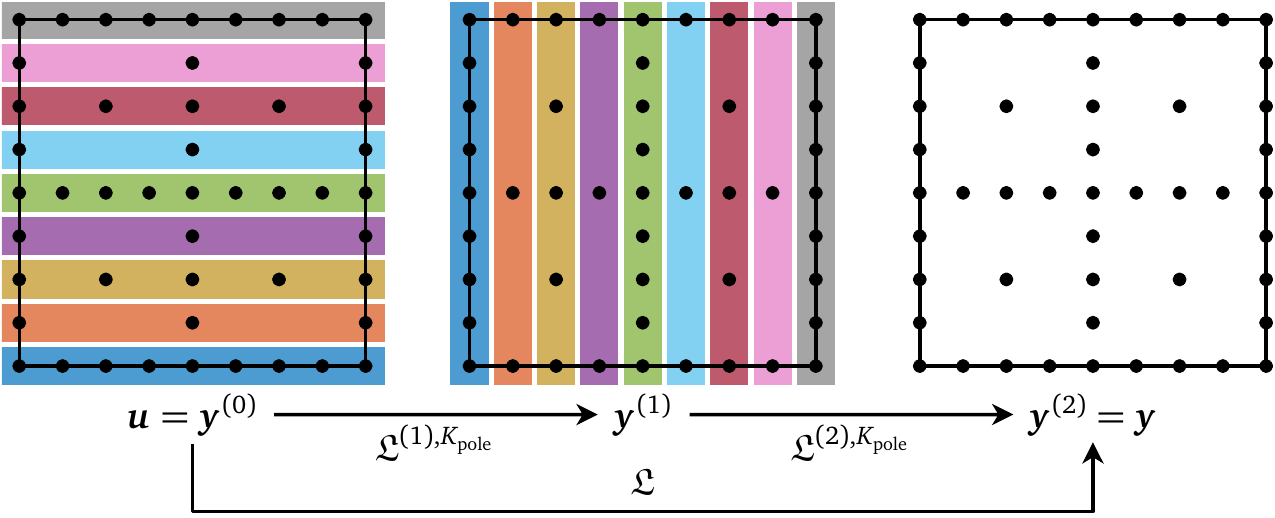}%
  \caption[%
    Unidirectional principle%
  ]{%
    Application of a linear operator $\linop$
    on two-dimensional sparse grid data with the unidirectional principle.
    \emph{Left:}
    The univariate operator $\upopuv{1}{\lisetpole}$ is applied on
    the input data $\vlinin$
    along poles $\lisetpole$ of the first dimension $x_1$.
    \vspace{-0.2em}%
    \emph{Center:}
    The univariate operator $\upopuv{2}{\lisetpole}$ is applied on the
    resulting intermediate data $\vlinout[(1)]$
    along poles $\lisetpole$ of the second dimension $x_2$.
    \emph{Right:}
    Final values $\vlinout = \linop[\vlinin]$ after the application
    on both dimensions.
    All grid points of the same color are part of the same pole $\lisetpole$
    (equivalence classes of $\samepole{t}$ in
    \cref{alg:unidirectionalPrinciple}).%
  }%
  \label{fig:unidirectionalPrinciple}%
\end{figure}

If $\sgset$ is a full grid $\fgset{\*l}$
(see \cref{sec:21nodalSpaces}),
the well-known \up
can be used to apply $\linop$ to input data $\vlinin$.
As shown in \cref{fig:unidirectionalPrinciple} for a sparse grid,
the idea of the \up
is to apply the corresponding one-dimensional operators on the
one-dimensional subgrids (the \term{poles}) of $\sgset$,
which is repeated for all dimensions.
In this section, we first formulate the \up for
general linear operators $\linop$ and then prove its correctness for
the case $\linop = \intpmatinv$ of hierarchization.
The correctness for arbitrary tensor product operators
will follow from \cref{sec:45spatAdaptiveUP}.

\paragraph{Unidirectional principle and its correctness}

We state the \up in \cref{alg:unidirectionalPrinciple}.
The algorithm is given a permutation $(t_1, \dotsc, t_d)$ of $(1, \dotsc, d)$
that specifies the order of dimensions in which the \up should be applied.
We denote with $\upopuv{t_j}{\lisetpole}$ the one-dimensional version of $\linop$
applied in dimension $t_j$ ($j = 1, \dotsc, d$) on the pole $\lisetpole$.
Formally, a pole is an equivalence class of the
\term{pole equivalence relation} $\samepole{t_j}$ on $\liset$:
\begin{equation}
  \label{eq:poleEquivalenceRelation}
  \*k' \samepole{t_j} \*k'' \iff \*k'_{-t_j} = \*k''_{-t_j},\quad
  \*k', \*k'' \in \liset.
\end{equation}
We prove the correctness of the \up with the following invariant:

\begin{algorithm}
  \begin{algorithmic}[1]
    \Function{$\vlinout = \texttt{unidirectionalPrinciple}$}{%
      $\vlinin$, $\liset$, $(t_1, \dotsc, t_d)$%
    }
      \State{$\vlinout[(0)] \gets \vlinin$}
      \For{$j = 1, \dotsc, d$}
        \For{$\lisetpole \in \eqclasses{\liset}{\samepole{t_j}}$}
          \State{%
            $(\linout[(j)]{\*k})_{\*k \in \lisetpole} \gets
            \upopuv{t_j}{\lisetpole}
            \bracket*{(\linout[(j-1)]{\*k})_{\*k \in \lisetpole}}$%
          }
          \Comment{apply univariate operator on pole}%
          \label{line:algUnidirectionalPrinciple1}
        \EndFor{}
      \EndFor{}
      \State{$\vlinout \gets \vlinout[(d)]$}
    \EndFunction{}
  \end{algorithmic}
  \caption[%
    Unidirectional principle%
  ]{%
    Application of a tensor product operator $\linop$ with
    the unidirectional principle.
    Inputs are
    the vector $\vlinin = (\linin{\*k})_{\*k \in \liset}$ of input data,
    the set $\liset$ of grid indices, and
    the permutation $(t_1, \dotsc, t_d)$ specifying the order in which
    the one-dimensional operators $\upopuv{t_j}{\lisetpole}$ should be applied.
    The output is the vector $\vlinout = (\linout{\*k})_{\*k \in \liset}$
    of output data.%
  }%
  \label{alg:unidirectionalPrinciple}%
\end{algorithm}

\begin{proposition}[invariant of unidirectional principle for hierarchization]
  \label{prop:invariantUnidirectionalPrinciple}
  Let $\linop$ be the hierarchization operator on a full grid,
  i.e.,
  $\linop = \intpmatinv$,
  $\vlinin = (\fcnval{\*k})_{\*k \in \liset}$,
  $\vlinout = (\surplus{\*k})_{\*k \in \liset}$,
  $\upopuv{t_j}{\lisetpole}$ is the univariate interpolation operator
  $\intpmatuvinv{t_j}$, and
  $\liset = \{\*0, \dotsc, \*2^\*l\}$
  corresponds to a full grid $\fgset{\*l}$ of level $\*l$.
  After iteration $j$ of \cref{alg:unidirectionalPrinciple}
  ($j = 1, \dotsc, d$), it holds for $T \ceq (t_1, \dotsc, t_j)$
  \begin{equation}
    \sum_{\*k_T=\*0}^{\*2^{\*l_T}}
    \linout[(j)]{(\*k_T,\*k'_{-T})} \basis{\*k_T}(\gp{\*k'_T})
    = \fcnval{\*k'},\quad
    \*k' = \*0, \dotsc, \*2^\*l,
  \end{equation}
  where $(\*k_T,\*k'_{-T})$ is shorthand for $\*k''$
  with $k''_t \ceq k_t$ if $t \in T$ and $k''_t \ceq k'_t$ if $t \notin T$.
\end{proposition}

\begin{proof}
  We prove the assertion by induction over $j = 1, \dotsc, d$.
  We set $T' \ceq (t_1, \dotsc, t_{j-1})$,
  $T \ceq (t_1, \dotsc, t_{j-1}, t_j)$,
  and we exploit the tensor product structure of the basis
  to write the \lhs of the assertion for $j$
  and arbitrary $\*k' = \*0, \dotsc, \*2^\*l$ as
  \begin{equation}
    \label{eq:proofPropInvariantUnidirectionalPrinciple1}
    \sum_{\*k_T=\*0}^{\*2^{\*l_T}}
    \linout[(j)]{(\*k_T,\*k'_{-T})} \basis{\*k_T}(\gp{\*k'_T})
    = \sum_{\*k_{T'}=\*0}^{\*2^{\*l_{T'}}}
    \basis{\*k_{T'}}(\gp{\*k'_{T'}})
    \sum_{k_{t_j}=0}^{2^{l_{t_j}}}
    \linout[(j)]{(\*k_T,\*k'_{-T})} \basis{k_{t_j}}(\gp{k'_{t_j}}).
  \end{equation}
  If we choose the equivalence class
  $\lisetpole \ceq \eqclass{(\*k_T,\*k'_{-T})}{\samepole{t_j}}$
  ($\*k_T$ arbitrary),
  then the inner sum over $k_{t_j}$ equals
  \begin{equation}
    \label{eq:proofPropInvariantUnidirectionalPrinciple2}
    \sum_{\*k'' \in \lisetpole}
    \linout[(j)]{\*k''} \basis{k''_{t_j}}(\gp{k'_{t_j}})
    = \paren*{
      (\upopuv{t_j}{\lisetpole})^{-1}
      \bracket*{(\linout[(j)]{\*k''})_{\*k'' \in \lisetpole}}
    }_{k'_{t_j}}
    = \linout[(j-1)]{(\*k_{T'},\*k'_{-T'})}
  \end{equation}
  by \cref{line:algUnidirectionalPrinciple1} of
  \cref{alg:unidirectionalPrinciple}.
  We can conclude that the \lhs
  \cref{eq:proofPropInvariantUnidirectionalPrinciple1} equals
  \begin{equation}
    \sum_{\*k_{T'}=\*0}^{\*2^{\*l_{T'}}}
    \linout[(j-1)]{(\*k_{T'},\*k'_{-T'})}
    \basis{\*k_{T'}}(\gp{\*k'_{T'}}),
  \end{equation}
  which, by the induction hypothesis, equals $\fcnval{\*k'}$ as desired
  (if $j > 1$).
  The same reasoning for
  \eqref{eq:proofPropInvariantUnidirectionalPrinciple2} can be used
  to establish the base case for $j = 1$.
\end{proof}

\begin{shortcorollary}[%
  correctness of unidirectional principle for hierarchization%
]
  \label{cor:algUnidirectionalPrincipleCorrectness}
  \Cref{alg:unidirectionalPrinciple}
  is correct for hierarchization on full grids.
\end{shortcorollary}

\begin{proof}
  We apply \cref{prop:invariantUnidirectionalPrinciple} for $j = d$ to obtain
  $\sum_{\*k=\*0}^{\*2^\*l}
  \linout[(j)]{\*k} \basis{\*k}(\gp{\*k'})
  = \fcnval{\*k'}$
  for all $\*k' = \*0, \dotsc, \*2^\*l$, i.e.,
  the $\linout[(j)]{\*k}$ are the correct interpolation coefficients
  according to \eqref{eq:hierarchizationProblem}.
\end{proof}

\paragraph{Complexity}

We compare the complexity of the \up for hierarchization with
the direct solution of the system \eqref{eq:hierarchizationSLE} of
linear equations.
If we assume that $d$ is constant and that
$\linop$ and $\upopuv{t_j}{\lisetpole}$ apply Gaussian elimination to
solve the multivariate and univariate systems, respectively,
then directly solving \eqref{eq:hierarchizationSLE} takes
$\landauO{\ngp^2 (\ngp + d)}$ time and
$\landauO{\ngp^2}$ memory.
In contrast, the \up only requires
$\landauO{\ngp \sum_t \ngp_t^2}$ time%
\footnote{%
  There are $\ngp/\ngp_t$ poles in the
  $t$-th iteration of \cref{alg:unidirectionalPrinciple}.
  Each pole requires the solution of an $\ngp_t \times \ngp_t$ linear system,
  which takes $\landauO{\ngp_t^3}$ time.%
}
and $\landauO{\max\{N_1^2, \dotsc, N_d^2, N\}}$ memory,
where $\ngp_t$ is the grid size
$\setsize{\{k_t \mid \*k \in \liset\}}$ in dimension $t = 1, \dotsc, d$.
The dependency from the univariate grid sizes $\ngp_t$ instead of $\ngp$
makes the \up significantly less computationally expensive.
As already mentioned,
the \up is even more efficient in the piecewise linear case,
where the univariate interpolation operators can be applied
in-place.
Hence, it only needs $\landauO{Nd}$ time and
$\landauO{N}$ memory in this case.

\section{Hierarchization on Dimensionally Adaptive Sparse Grids}
\label{sec:43dimAdaptive}

\minitoc{77mm}{7}

\noindent
Dimensionally adaptive sparse grids,
which are sums of different hierarchical subspaces
as described in \cref{sec:232dimensionallyAdaptiveSG},
have the advantage over general spatially adaptive sparse grids
that algorithms can be formulated and applied more easily.
In this section, we describe two methods:
first, the well-known combination technique, which was
already mentioned in \cref{sec:232dimensionallyAdaptiveSG},
and second, a new algorithm based on residual interpolation.

\subsection{The Combination Technique and Its Combinatorial Proof}
\label{sec:431combiTechniqueProof}

The combination technique was one of the first methods that
were developed by Griebel et al. in \cite{Griebel92Combination}
(for two and three dimensions)
after the term ``sparse grids'' was coined in 1991 \cite{Zenger91Sparse}.
However, the combination technique predates the development of sparse grids
by decades, as it was already mentioned by Smolyak in 1963
\multicite{Smolyak63Quadrature,Hegland07Combination}.
Delvos developed and proved the standard combination formula in the
framework of Boolean interpolation operators in 1982
\multicite{Delvos82Dvariate,Delvos89Boolean}.

\paragraph{Formal description and outline of a combinatorial proof}

In the following, we give a formal description of the
sparse grid combination technique, and we outline a new combinatorial proof
of its correctness.
While we discuss a high-level explanation of the proofs in this section,
the proofs themselves can be found in \cref{sec:a131proofCombiTechnique},
since most of them are rather technical.
For simplicity,
we formulate the combination technique and its proof for regular
sparse grids (see \cref{sec:231regularSG}).
However, the main ideas of the chain of proofs are also applicable
to dimensionally adaptive sparse grids
(see \cref{sec:232dimensionallyAdaptiveSG}).

\begin{theorem}[sparse grid combination technique]
  \label{thm:combiTechnique}
  Let $\liset \ceq \{(\*l, \*i) \mid
  \normone{\*l} \le n,\; \*i \in \hiset{\*l}\}$
  correspond to the regular sparse grid
  $\regsgset{n}{d}$ and let $(\fcnval{\*l,\*i})_{(\*l,\*i) \in \liset}$
  be given function values on $\regsgset{n}{d}$.
  If we define
  \begin{itemize}
    \item
    the combined sparse grid interpolant $\regsgintp[\ct]{n}{d}$ via
    \eqref{eq:combiTechnique}, i.e.,
    \begin{equation}
      \regsgintp[\ct]{n}{d}
      = \sum_{q=0}^{d-1} (-1)^q \binom{d-1}{q} \sum_{\normone{\*l'} = n-q}
      \fgintp{\*l'},
    \end{equation}
    where $\fgintp{\*l'} \in \ns{\*l'}$ is the full grid interpolant
    of $\objfun$ with level $\*l'$, and
    
    \item
    the hierarchical sparse grid interpolant $\regsgintp{n}{d}$
    via \eqref{eq:hierarchizationProblem} and
    \eqref{eq:hierarchizationInterpolant}
  \end{itemize}
  and if we assume that the hierarchical splitting equation
  \eqref{eq:hierSplittingMV} holds,
  then the combined and the hierarchical sparse grid interpolants coincide:
  \begin{equation}
    \regsgintp[\ct]{n}{d}
    = \regsgintp{n}{d}.
  \end{equation}
\end{theorem}

\begin{proof}[Proof (sketch)]
  Let $\gp{\*l,\*i} \in \regsgset{n}{d}$ be an arbitrary
  point of the regular sparse grid.
  First, we split the inner sum of $\regsgintp[\ct]{n}{d}(\gp{\*l,\*i})$
  into levels $\*l'$ whose full grid sets $\fgset{\*l'}$
  contain $\gp{\*l,\*i}$ and levels whose full grid sets
  do not contain $\gp{\*l,\*i}$:
  \begin{equation}
    \regsgintp[\ct]{n}{d}(\gp{\*l,\*i})
    = \sum_{q=0}^{d-1} (-1)^q \binom{d-1}{q} \cdot \paren*{
      \sum_{\substack{\normone{\*l'} = n - q\\\fgset{\*l'} \ni \gp{\*l,\*i}}}
      \fgintp{\*l'}(\gp{\*l,\*i}) +
      \sum_{\substack{\normone{\*l'} = n - q\\\fgset{\*l'} \notni \gp{\*l,\*i}}}
      \fgintp{\*l'}(\gp{\*l,\*i})
    }.
  \end{equation}
  The summands $\fgintp{\*l'}(\gp{\*l,\*i})$ of the first inner sum
  each equal $\fcnval{\*l,\*i}$ due to the full grid interpolation
  property \eqref{eq:interpFullGridMV}.
  Therefore, the first inner sum is equal to the product of
  $\fcnval{\*l,\*i}$ with the number of summands:
  \begin{equation}
    \label{eq:combiTechniqueSplitSum}
    \begin{split}
      \regsgintp[\ct]{n}{d}(\gp{\*l,\*i})
      &= \fcnval{\*l,\*i} \cdot \sum_{q=0}^{d-1} (-1)^q \binom{d-1}{q} \cdot
      \setsize{
        \{\*l' \mid \normone{\*l'} = n - q,\; \fgset{\*l'} \ni \gp{\*l,\*i}\}
      }\\
      &\qquad{} {} + \sum_{q=0}^{d-1} (-1)^q \binom{d-1}{q} \cdot
      \sum_{\substack{\normone{\*l'} = n - q\\\fgset{\*l'} \notni \gp{\*l,\*i}}}
      \fgintp{\*l'}(\gp{\*l,\*i}).
    \end{split}
  \end{equation}
  After this sketch of proof,
  we will prove that the first of the two summands in
  \cref{eq:combiTechniqueSplitSum}
  equals one (see \cref{prop:combiTechniqueOne})
  and that the second of the two summands
  equals zero (see \cref{prop:combiTechniqueZero}).
  Consequently, we infer
  \begin{equation}
    \regsgintp[\ct]{n}{d}(\gp{\*l,\*i})
    = \fcnval{\*l,\*i},
  \end{equation}
  i.e., $\regsgintp[\ct]{n}{d}$ interpolates $\objfun$ at $\regsgset{n}{d}$.
  Note that $\regsgintp[\ct]{n}{d}$ is contained in $\regsgspace{n}{d}$,
  if the hierarchical splitting equation \eqref{eq:hierSplittingMV} holds,
  as
  \begin{equation}
    \fgintp{\*l'} \in
    \ns{\*l'}
    = \bigoplus_{\*l''=\*0}^{\*l'} \hs{\*l''}
    \subset \regsgspace{n}{d},\quad
    \normone{\*l'} \le n,
  \end{equation}
  due to $\normone{\*l''} \le \normone{\*l'} \le n$
  for $\*l'' \le \*l'$, i.e.,
  $\hs{\*l''} \subset \regsgspace{n}{d}$ for
  $\*l'' = \*0, \dotsc, \*l'$.%
  \footnote{%
    This argumentation can be straightforwardly adapted
    for general dimensionally adaptive sparse grids
    with downward closed level sets as mentioned in
    \cref{sec:232dimensionallyAdaptiveSG}.%
  }
  As both $\regsgintp[\ct]{n}{d}$ and $\regsgintp{n}{d}$
  are contained in $\regsgspace{n}{d}$ and
  interpolate $\objfun$ on $\regsgset{n}{d}$, they coincide
  due to the uniqueness of sparse grid interpolation
  (linear independence of the hierarchical basis).
\end{proof}

\paragraph{Inclusion-exclusion principle}

It remains to prove that the first sum in \eqref{eq:combiTechniqueSplitSum}
is indeed one and that the second sum vanishes.
The first statement is a direct consequence of the
\term{inclusion-exclusion principle} \cite{Hegland07Combination}.
In its simplest form, the idea of the principle is that the cardinality
of the union of two finite subsets $A$ and $B$ of some set is given by
\begin{equation}
  \setsize{A \cup B}
  = \setsize{A} + \setsize{B} - \setsize{A \cap B},
\end{equation}
i.e., we first count (include)
the elements in $A$ and then in $B$,
but as we have counted the elements of $A \cap B$ twice,
we have to subtract (exclude) its cardinality again.

The setting is similar for the combination technique.
If we add all grids in \cref{fig:combinationTechnique}
on the green diagonal, then every point whose index is not odd
will be counted multiple times.
By subtracting the number of occurrences of the points on the
red diagonal,
the result of the ``weighted counting'' is exactly one for every point.
The following proposition, whose proof is of purely combinatorial nature,
generalizes this argument to higher dimensions:

\begin{restatable}[inclusion-exclusion principle]{%
  proposition%
}{%
  propCombiTechniqueOne%
}
  \label{prop:combiTechniqueOne}
  For every $\gp{\*l,\*i} \in \regsgset{n}{d}$, we have
  \begin{equation}
    \label{eq:combiTechniqueOne}
    \sum_{q=0}^{d-1} (-1)^q \binom{d-1}{q} \cdot
    \setsize{
      \{\*l' \mid \normone{\*l'} = n - q,\; \fgset{\*l'} \ni \gp{\*l,\*i}\}
    }
    = 1.
  \end{equation}
\end{restatable}

\begin{proof}
  See \cref{sec:a131proofCombiTechnique}.
\end{proof}

\paragraph{Canceling out function values}

The second statement about the vanishing second sum in
\eqref{eq:combiTechniqueSplitSum} is much harder to prove.
It says that at every grid point $\gp{\*l,\*i}$,
the contributions $\fgintp{\*l'}$ of levels $\*l'$
that do not contain that point cancel out,
which may seem quite surprising.
The key observation is as follows:
The values of $\fgintp{\*l'}, \fgintp{\*l''}$ for two levels
$\*l', \*l''$ are the same at $\gp{\*l,\*i}$,
if all non-equal entries $l'_t, l''_t$ of the levels are
each greater or equal to $l_t$.

For a higher-level explanation,
note that the statement $l'_t \ge l_t$ is equivalent to
$\fgset{l'_t} \ni \gp{l_t,i_t}$.
Both $\fgintp{\*l'}, \fgintp{\*l''}$ interpolate at
$\gp{\*l,\*i}$ when projected onto the $t$-th dimension,
so their contribution to $\fgintp{\*l'}(\gp{\*l,\*i})$ and
$\fgintp{\*l''}(\gp{\*l,\*i})$ must be the same.
Although there may be dimensions $t$ for which
$\fgset{l'_t} \notni \gp{l_t,i_t}$,
these dimensions do not matter if $l'_t = l''_t$,
as the univariate restrictions of $\fgintp{\*l'}, \fgintp{\*l''}$
interpolate the same data, and they are evaluated at the same point
$\gp{l_t,i_t}$.

One can formalize these considerations by defining an
equivalence relation on the set of levels such that the values of
$\fgintp{\*l'}$ at $\gp{\*l,\*i}$ are constant
on the equivalence classes.

\begin{definition}[%
  equivalence relation for the proof of the combination technique%
]
  \label{def:combiTechniqueEquivalenceRelation}
  Let $\gp{\*l,\*i} \in \regsgset{n}{d}$ be fixed and
  \begin{equation}
    \label{eq:combiTechniqueSpecialLevelSet}
    L
    \ceq \{\*l' \mid \ex{q=0,\dotsc,d-1}{
      \normone{\*l'} = n - q,\; \fgset{\*l'} \notni \gp{\*l,\*i}
    }\}
  \end{equation}
  be the set of levels that do not contain $\gp{\*l,\*i}$.
  We define a relation $\eq$ on $L$ as follows:
  For $\*l', \*l'' \in L$, we set $\*l' \eq \*l''$ if and only if
  \begin{equation}
    \falarge{t \notin T_{\*l',\*l''}}{\min\{l'_t, l''_t\} \ge l_t},\quad
    T_{\*l',\*l''}
    \ceq \{t \mid l'_t = l''_t < l_t\}.
  \end{equation}
\end{definition}

\begin{restatable}[identical values in equivalence classes]{%
  shortlemma%
}{%
  lemmaCombiTechniqueIdenticalValues%
}
  \label{lemma:combiTechniqueIdenticalValues}
  Let $\*l', \*l'' \in L$ with $\*l' \eq \*l''$.
  Then, $\fgintp{\*l'}(\gp{\*l,\*i})
  = \fgintp{\*l''}(\gp{\*l,\*i})$.
\end{restatable}

\begin{proof}
  See \cref{sec:a131proofCombiTechnique}.
\end{proof}

By exploiting the tensor product structure of the basis functions,
the proof shows an even stronger version, which is shown
in \cref{fig:combiTechniqueProof}:
The components $\fgintp{\*l'}$ and $\fgintp{\*l''}$ are equal
on the $m$-dimensional affine subspace through $\gp{\*l,\*i}$
parallel to the $m$ coordinates in $T_{\*l',\*l''}$
(where $m \ceq \setsize{T_{\*l',\*l''}}$).
The lemma allows us to group summands in the second sum of
\eqref{eq:combiTechniqueSplitSum} by function values.
Hence, it remains to count the number of levels in
each equivalence class of $\eq$.
Therefore, we need a characterization of the equivalence classes:

\begin{SCfigure}
  \includegraphics{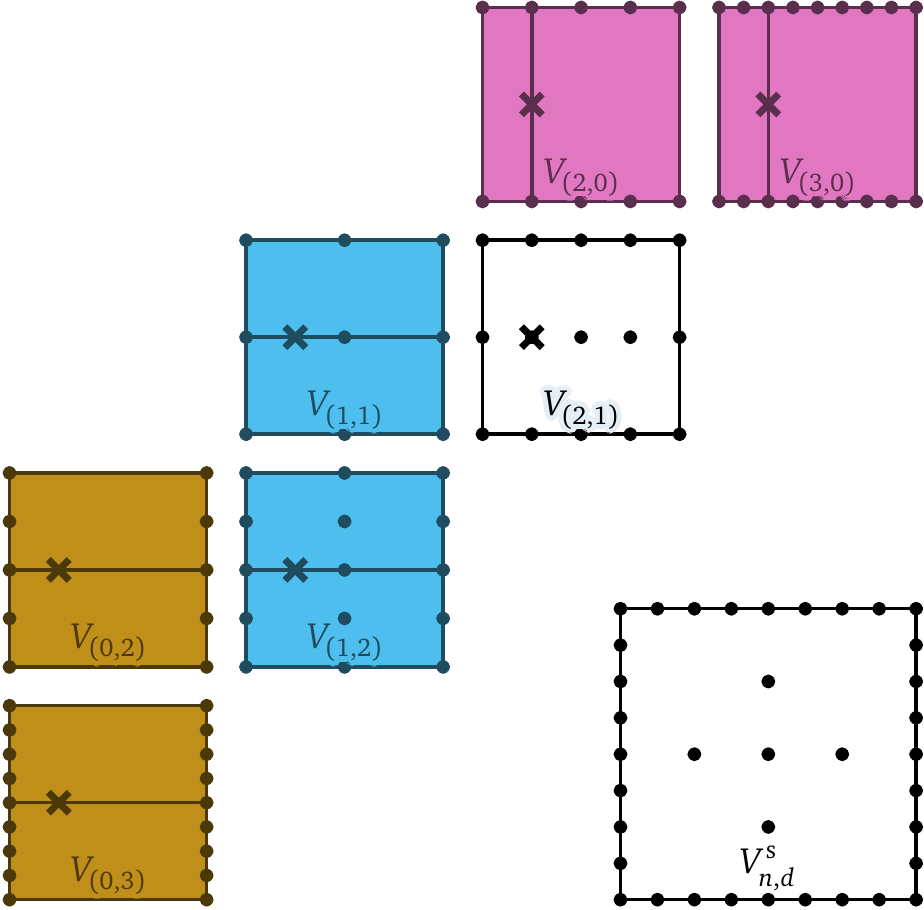}%
  \caption[%
    Canceling out function values in the proof of the combination technique%
  ]{%
    Nodal subspaces $\ns{\*l}$ contributing to the combination
    technique solution for the two-dimensional regular sparse grid
    $\regsgspace{n}{d}$ of level $n = 3$ \emph{(bottom right).}
    After picking a point $\gp{\*l,\*i} \in \regsgset{n}{d}$
    (\emph{cross,} here $\*l = (2, 1)$, $\*i = (1, 1)$),
    the set $L$ of levels whose grids do not contain $\gp{\*l,\*i}$
    \emph{(colored subspaces)}
    decompose into three disjoint equivalence classes
    \emph{(colors)} given by the relation $\eq$.
    In every equivalence class $L_0 \in \eqclasses{L}{\eq}$,
    the interpolants $\fgintp{\*l'}$ ($\*l' \in L_0$)
    equal on an affine subspace
    \emph{(dark lines),} which contains $\gp{\*l,\*i}$.
    Due to the combination coefficients,
    the contribution to the combined solution
    vanishes per equivalence class.%
  }%
  \label{fig:combiTechniqueProof}%
\end{SCfigure}

\begin{restatable}[characterization of equivalence classes]{%
  lemma%
}{%
  lemmaCombiTechniqueCharacterization%
}
  \label{lemma:combiTechniqueCharacterization}
  Let $L_0 \in \eqclasses{L}{\eq}$ be an equivalence class of $\eq$.
  If we define
  \begin{equation}
    T_{L_0}
    \ceq \{t \mid \exfa{l^\ast_t < l_t}{\*l' \in L_0}{l'_t = l^\ast_t}\}
  \end{equation}
  as the set of dimensions $t$ in which all levels in $L_0$
  have the same entry $l^\ast_t < l_t$, then
  \begin{equation}
    L_0
    = \{\*l' \in L \mid
    \fa{t \in T_{L_0}}{l'_t = l^\ast_t},\;
    \fa{t \notin T_{L_0}}{l'_t \ge l_t}\}.
  \end{equation}
\end{restatable}

\begin{proof}
  See \cref{sec:a131proofCombiTechnique}.
\end{proof}

The lemma states that every equivalence class $L_0$ is exactly the
set of the levels whose entries are equal
and smaller than $l_t$ in some dimensions
(which are contained in $T_{L_0}$)
and whose entries are greater or equal than $l_t$ in all other dimensions.
While this statement may seem intuitively correct,
the proof is rather technical.
Finally, we are now able to show that the second sum in
\eqref{eq:combiTechniqueSplitSum} vanishes:

\begin{restatable}[function value cancellation]{%
  proposition%
}{%
  propCombiTechniqueZero%
}
  \label{prop:combiTechniqueZero}
  For every $\gp{\*l,\*i} \in \regsgset{n}{d}$, we have
  \begin{equation}
    \sum_{q=0}^{d-1} (-1)^q \binom{d-1}{q} \cdot
    \sum_{\substack{\normone{\*l'} = n - q\\\fgset{\*l'} \notni \gp{\*l,\*i}}}
    \fgintp{\*l'}(\gp{\*l,\*i})
    = 0.
  \end{equation}
\end{restatable}

\begin{proof}
  See \cref{sec:a131proofCombiTechnique}.
\end{proof}

The proof essentially first counts the number of possible levels in
an equivalence class and then applies known combinatorial identities
to prove that the sum must vanish.
This proves \thmref{thm:combiTechnique}.

\subsection{Hierarchization with the Combination Technique}
\label{sec:432hierarchizationCombiTechnique}

It is straightforward to hierarchize function values
$\fcnval{\*l,\*i}$ on dimensionally adaptive sparse grids
with the combination technique.
The resulting hierarchization algorithm is given as \cref{alg:combiTechnique}.
In \cref{line:algCombiTechnique1},
the hierarchical surpluses corresponding to the full grid
interpolant $\fgintp{\*l'} \in \ns{\*l'}$ have to be computed
(see \eqref{eq:interpFullGridMV}).
As shown in \cref{sec:42fullGrids}, we can easily calculate these
surpluses with the unidirectional principle in
\cref{alg:unidirectionalPrinciple}.
The surpluses are then combined with the same combination formula
as in \thmref{thm:combiTechnique}.
Note that it is imperative to employ the hierarchical basis functions
$\basis{\*l,\*i}$ with $\*l = \*0, \dotsc, \*l'$ and $\*i \in I_{\*l}$
and not the nodal basis,
i.e., $\basis{\*l',\*i'}$ with $\*i' = \*0, \dotsc, \*2^{\*l'}$.

\begin{algorithm}
  \begin{algorithmic}[1]
    \Function{$\vlinout = \texttt{combinationTechnique}$}{%
      $\vlinin$, $n$, $d$%
    }
      \For{$q = 0, \dotsc, d - 1$}
        \For{$\*l' \in \natz^d$ with $\normone{\*l'} = n - q$}
          \State{%
            Let $(\surplus[(\*l')]{\*l,\*i})_{
              \*l = \*0, \dotsc, \*l'\!,\, \*i \in \hiset{\*l}
            }$ be such that
            $\sum_{\*l=\*0}^{\*l'} \sum_{\*i \in \hiset{\*l}}
            \surplus[(\*l')]{\*l,\*i} \basis{\*l,\*i} \equiv
            \fgintp{\*l'}$%
          }
          \label{line:algCombiTechnique1}
          \State{%
            $\surplus[(\*l')]{\*l,\*i} \gets 0$
            for all $(\*l,\*i) \in \liset$
            with $\lnot(\*l \le \*l')$%
          }
          \Comment{extend surpluses}%
          \label{line:algCombiTechnique2}
        \EndFor{}
      \EndFor{}
      \State{%
        $\linout{\*l,\*i}
        = \sum_{q=0}^{d-1} (-1)^q \binom{d-1}{q}
        \sum_{\normone{\*l'} = n-q} \surplus[(\*l')]{\*l,\*i}$
        for all $(\*l, \*i) \in \liset$%
      }
      \Comment{combine surpluses}%
    \EndFunction{}
  \end{algorithmic}
  \caption[%
    Hierarchization with the combination technique%
  ]{%
    Application of the hierarchization operator $\linop = \intpmatinv$
    with the combination technique.
    For simplicity,
    the algorithm is described for regular sparse grids,
    but it can be generalized to arbitrary dimensionally adaptive sparse grids.
    Inputs are
    the vector $\vlinin = (\linin{\*l,\*i})_{(\*l,\*i) \in \liset}$
    of input data (function values $\fcnval{\*l,\*i}$ at the grid points),
    the level $n$, and the dimensionality $d$ of the regular sparse grid,
    where $\liset$ is the set of all feasible level-index pairs $(\*l,\*i)$,
    i.e., $\normone{\*l} \le n$, $\*i \in \hiset{\*l}$.
    The output is the vector
    $\vlinout = (\linout{\*l,\*i})_{(\*l,\*i) \in \liset}$
    of output data (hierarchical surpluses $\surplus{\*l,\*i}$).%
  }%
  \label{alg:combiTechnique}%
\end{algorithm}

\paragraph{Correctness}

Of course, the proof of the correctness of \cref{alg:combiTechnique}
relies on the correctness of the combination technique
(see \cref{thm:combiTechnique}).
If determining the combination coefficients correctly
\cite{Nobile16Adaptive}, the algorithm can even be applied to
all dimensionally adaptive sparse grids.
The proof of the following proposition can be generalized accordingly.

\begin{proposition}[correctness of combination technique]
  \label{prop:correctnessAlgCombiTechnique}
  \Cref{alg:combiTechnique}
  is correct for hierarchization on regular sparse grids.
\end{proposition}

\begin{proof}
  According to \cref{line:algCombiTechnique1} of \cref{alg:combiTechnique},
  the full grid interpolants $\fgintp{\*l'}$ can be written as
  \begin{equation}
    \fgintp{\*l'}
    = \sum_{\normone{\*l} \le n} \sum_{\*i \in \hiset{\*l}}
    \surplus[(\*l')]{\*l,\*i} \basis{\*l,\*i}
  \end{equation}
  where the surpluses have been extended
  from $\*l = \*0, \dotsc, \*l'$ to all $\*l$ with $\normone{\*l} \le n$
  by zero in \cref{line:algCombiTechnique2}.
  \Cref{thm:combiTechnique} now allows to write the hierarchical
  interpolant $\regsgintp{n}{d}$ in terms of the full grid components:
  \begin{subequations}
    \begin{align}
      \regsgintp{n}{d}
      = \regsgintp[\ct]{n}{d}
      &= \sum_{q=0}^{d-1} (-1)^q \binom{d-1}{q} \sum_{\normone{\*l'} = n-q}
      \fgintp{\*l'}\\
      &= \sum_{q=0}^{d-1} (-1)^q \binom{d-1}{q} \sum_{\normone{\*l'} = n-q}
      \sum_{\normone{\*l} \le n} \sum_{\*i \in \hiset{\*l}}
      \surplus[(\*l')]{\*l,\*i} \basis{\*l,\*i}.\\
      \label{eq:propCorrectnessAlgCombiTechnique1}
      &= \sum_{\normone{\*l} \le n} \sum_{\*i \in \hiset{\*l}}
      \underbrace{
        \paren*{
          \sum_{q=0}^{d-1} (-1)^q \binom{d-1}{q} \sum_{\normone{\*l'} = n-q}
          \surplus[(\*l')]{\*l,\*i}
        }
      }_{= \linout{\*l,\*i}}
      \basis{\*l,\*i},
    \end{align}
  \end{subequations}
  where $\linout{\*l,\*i}$ is the $(\*l,\*i)$-th entry of the output vector
  of \cref{alg:combiTechnique}.
  Note that the hierarchical interpolant $\regsgintp{n}{d}$
  can be written as
  $\regsgintp{n}{d} = \sum_{\normone{\*l} \le n} \sum_{\*i \in \hiset{\*l}}
  \surplus{\*l,\*i} \basis{\*l,\*i}$
  (see \eqref{eq:regularSGInterpolant}),
  where the surpluses $\surplus{\*l,\*i}$ are unique due to the
  linear independence of the hierarchical basis.
  As \eqref{eq:propCorrectnessAlgCombiTechnique1}
  equals $\regsgintp{n}{d}$ and has the same form,
  the coefficients $\linout{\*l,\*i}$
  must coincide with the surpluses $\surplus{\*l,\*i}$.
\end{proof}

\subsection{Hierarchization with Residual Interpolation}
\label{sec:433residualInterpolation}

Another method to hierarchize function values on
dimensionally adaptive sparse grids is the
\term{method of residual interpolation.}
The advantage over the combination technique is that
it only needs to operate on so-called \term{active nodal spaces.}
In contrast, the combination technique needs to perform computations
on additional non-active nodal subspaces
(for the regular sparse grid case:
summands with $q \ge 1$ in \eqref{eq:combiTechnique}).

\paragraph{Active nodal spaces}

\Cref{alg:residualInterpolation} describes the procedure of
the method of residual interpolation,
given the function values $\vlinin$ corresponding to the grid points
and the levels $L$ contained in the sparse grid
(see \eqref{eq:dimensionallyAdaptiveSG}).
The list $\*l^{(1)}, \dotsc, \*l^{(m)}$ of active nodal spaces
in \cref{line:algResidualInterpolation1} is determined by the condition
\begin{equation}
  \label{eq:activeNodalSpaces}
  \bigcup_{j=1}^m \{\*l \in \natz^d \mid \*l \le \*l^{(j)}\} = L,\quad
  \falarge{j_1 \not= j_2}{\lnot(\*l^{(j_1)} \le \*l^{(j_2)})}.
\end{equation}
This means that the corresponding sparse grid $\sgset$
is the (non-disjoint) union of the full grid sets $\fgset{\*l^{(j)}}$
($j = 1, \dotsc, m$)
and no full grid set is contained in another, i.e.,
no full grid set can be omitted without
removing points from the union $\sgset$.

\begin{algorithm}
  \begin{algorithmic}[1]
    \Function{$\vlinout = \texttt{residualInterpolation}$}{%
      $\vlinin$, $\levelset$%
    }
      \State{%
        $r^{(0)}(\gp{\*l,\*i}) \gets \fcnval{\*l,\*i}$
        for all $(\*l,\*i) \in \liset$%
      }
      \State{%
        Compute list $\*l^{(1)}, \dotsc, \*l^{(m)}$
        of active nodal spaces from $L$ (see \eqref{eq:activeNodalSpaces})%
      }
      \label{line:algResidualInterpolation1}
      \State{%
        Sort $\*l^{(1)}, \dotsc, \*l^{(m)}$ by decreasing level sum%
      }
      \For{$j = 1, \dotsc, m$}
        \State{%
          Let $r_{\*l^{(j)}}^{(j-1)} \in \ns{\*l^{(j)}}$ be the
          interpolant of $r^{(j-1)}$ on $\fgset{\*l^{(j)}}$%
        }
        \label{line:algResidualInterpolation3}
        \State{%
          Let $(\surplus[(j)]{\*l,\*i})_{(\*l,\*i) \in \liset}$ be such that
          $\sum_{\*l=\*0}^{\*l^{(j)}} \sum_{\*i \in \hiset{\*l}}
          \surplus[(j)]{\*l,\*i} \basis{\*l,\*i}
          \equiv r_{\*l^{(j)}}^{(j-1)}$%
        }
        \Comment{interpolation}%
        \label{line:algResidualInterpolation2}
        \State{%
          $r^{(j)}(\gp{\*l,\*i}) \gets
          r^{(j-1)}(\gp{\*l,\*i}) - r_{\*l^{(j)}}^{(j-1)}(\gp{\*l,\*i})$
          for all $(\*l,\*i) \in \liset$%
        }
        \Comment{new residuals}%
        \label{line:algResidualInterpolation4}
      \EndFor{}
      \State{%
        $\vlinout \gets \sum_{j=1}^{m} \vsurplus^{(j)}$
        (where $\surplus[(j)]{\*l,\*i} = 0$,
        $(\*l,\*i) \in \liset$,
        if $\lnot(\*l \le \*l^{(j)})$)%
      }
      \Comment{combine surpluses}%
    \EndFunction{}
  \end{algorithmic}
  \caption[%
    Hierarchization with residual interpolation%
  ]{%
    Application of the hierarchization operator $\linop = \intpmatinv$
    with residual interpolation
    for dimensionally adaptive sparse grids.
    Inputs are
    the vector $\vlinin = (\linin{\*l,\*i})_{(\*l,\*i) \in \liset}$
    of input data (function values $\fcnval{\*l,\*i}$ at the grid points) and
    the set $\levelset$ of levels that are part of
    the sparse grid (see \eqref{eq:dimensionallyAdaptiveSG}),
    where $\liset$ is the set of all feasible level-index pairs $(\*l,\*i)$,
    i.e., $\*l \in \levelset$, $\*i \in \hiset{\*l}$.
    The output is the vector
    $\vlinout = (\linout{\*l,\*i})_{(\*l,\*i) \in \liset}$
    of output data (hierarchical surpluses $\surplus{\*l,\*i}$).%
  }%
  \label{alg:residualInterpolation}%
\end{algorithm}

\paragraph{Correctness}

The principle of \Cref{alg:residualInterpolation} is maintaining
a vector $(r^{(j)}(\gp{\*l,\*i}))_{(\*l,\*i) \in \liset}$ of residuals
and interpolating the residual data subsequently on the active nodal spaces.
Again, note that it is necessary to compute the coefficients
$\surplus[(j)]{\*l,\*i}$ in the hierarchical basis, despite interpolating
on the full grid $\fgset{\*l^{(j)}}$.
In \cref{chap:a10proofs}, we prove that the algorithm satisfies
the following invariant, which can be used to show its correctness:

\begin{restatable}[invariant of residual interpolation]{%
  proposition%
}{%
  propInvariantResidualInterpolation%
}
  \label{prop:invariantResidualInterpolation}
  For $j = 1, \dotsc, m$, it holds
  \begin{subequations}
    \label{eq:propInvariantResidualInterpolationStatements}
    \begin{alignat}{4}
      \label{eq:propInvariantResidualInterpolation1}
      r_{\*l^{(j)}}^{(j-1)}(\gp{\*l,\*i})
      &= 0,\quad
      &&\*l \le \*l^{(j')},\;\;
      &&\*i \in \hiset{\*l},\quad
      &j'
      &= 1, \dotsc, j - 1,\\
      \label{eq:propInvariantResidualInterpolation2}
      r^{(j)}(\gp{\*l,\*i})
      &= 0,\quad
      &&\*l \le \*l^{(j')},\;\;
      &&\*i \in \hiset{\*l},\quad
      &j'
      &= 1, \dotsc, j,\\
      \label{eq:propInvariantResidualInterpolation3}
      r^{(j)}(\gp{\*l,\*i})
      &= \fcnval{\*l,\*i} - f^{\sparse,(j)}(\gp{\*l,\*i}),\quad
      &&\*l \in L,\;\;
      &&\*i \in \hiset{\*l},&&
    \end{alignat}
  \end{subequations}%
  \setlength{\abovedisplayskip}{0pt}%
  \begin{equation}
    \label{eq:propInvariantResidualInterpolation4}
    \text{where}\quad
    f^{\sparse,(j)}
    \ceq \sum_{\*l' \in \levelset} \sum_{\*i' \in \hiset{\*l'}}
    \paren*{\sum_{j'=1}^{j} \surplus[(j')]{\*l',\*i'}} \basis{\*l',\*i'}.
  \end{equation}
\end{restatable}

\begin{proof}
  See \cref{sec:a132proofResidualInterpolation}.
\end{proof}

\begin{corollary}[correctness of residual interpolation]
  \label{cor:algResidualInterpolationCorrectness}
  \Cref{alg:residualInterpolation} is correct for hierarchization
  on dimensionally adaptive sparse grids.
\end{corollary}

\begin{proof}
  Let $\*l \in \levelset$ and $\*i \in \hiset{\*l}$.
  By construction of the active nodal spaces,
  there exists some $j' \in \{1, \dotsc, m\}$ such that $\*l \le \*l^{(j')}$.
  By \cref{prop:invariantResidualInterpolation}, we obtain
  for $j = m$
  \begin{subequations}
    \label{eq:proofCorAlgResidualInterpolationCorrectness1}
    \begin{align}
      \sum_{\*l' \in L} \sum_{\*i' \in \hiset{\*l'}}
      \smash{
        \underbrace{
          \paren*{\sum_{j''=1}^{m} \surplus[(j'')]{\*l',\*i'}}
        }_{= \linout{\*l',\*i'}}
      }
      \basis{\*l',\*i'}(\gp{\*l,\*i})
      &\quad\;
      \mathclap{\overset{\eqref{eq:propInvariantResidualInterpolation4}}{=}}
      \quad\;
      f^{\sparse,(m)}(\gp{\*l,\*i})
      \overset{\eqref{eq:propInvariantResidualInterpolation3}}{=}
      \fcnval{\*l,\*i} - r^{(m)}(\gp{\*l,\*i})\\
      &\quad\;
      \mathclap{\overset{\eqref{eq:propInvariantResidualInterpolation2}}{=}}
      \quad\;
      \fcnval{\*l,\*i}.
    \end{align}
  \end{subequations}
  As the hierarchical interpolant $\sgintp$
  (see \eqref{eq:hierarchizationInterpolant})
  has the same form
  $\sum_{\*l' \in \levelset} \sum_{\*i' \in \hiset{\*l'}}
  \surplus{\*l',\*i'} \basis{\*l',\*i'}$ as the \lhs of
  \eqref{eq:proofCorAlgResidualInterpolationCorrectness1}
  with unique surpluses $\surplus{\*l',\*i'}$ such that the function values
  are interpolated (see \eqref{eq:hierarchizationProblem}),
  the coefficients $\linout{\*l',\*i'}$
  (output of \cref{alg:residualInterpolation})
  coincide with the surpluses $\surplus{\*l',\*i'}$.
\end{proof}

\vspace{1em}

\Cref{prop:invariantResidualInterpolation} shows that
$r^{(j)}(\gp{\*l,\*i})$ is the residual of the
interpolant $f^{\sparse,(j)}$ of iteration~$j$
to the objective function $\objfun$ at the grid points $\gp{\*l,\*i}$
(\cref{eq:propInvariantResidualInterpolation3}).
After interpolating $r^{(j-1)}$ on the grid $\fgset{\*l^{(j)}}$
to obtain the function $r_{\*l^{(j)}}^{(j-1)}$
and subtracting the resulting values from the old residual values,
the new residual values $r^{(j)}(\gp{\*l,\*i})$ vanish
not only on the grid $\{(\*l^{(j)}, \*i) \mid \*i \in \hiset{\*l^{(j)}}\}$,
but also on all previous grids
$\{(\*l^{(j')}, \*i) \mid \*i \in \hiset{\*l^{(j')}}\}$, $j' \le j$
(\cref{eq:propInvariantResidualInterpolation2}).
The proof of \cref{prop:invariantResidualInterpolation}
shows this by exploiting the auxiliary statement of
\cref{eq:propInvariantResidualInterpolation1}
and the tensor product structure of the hierarchical basis.

\begin{SCfigure}
  \includegraphics{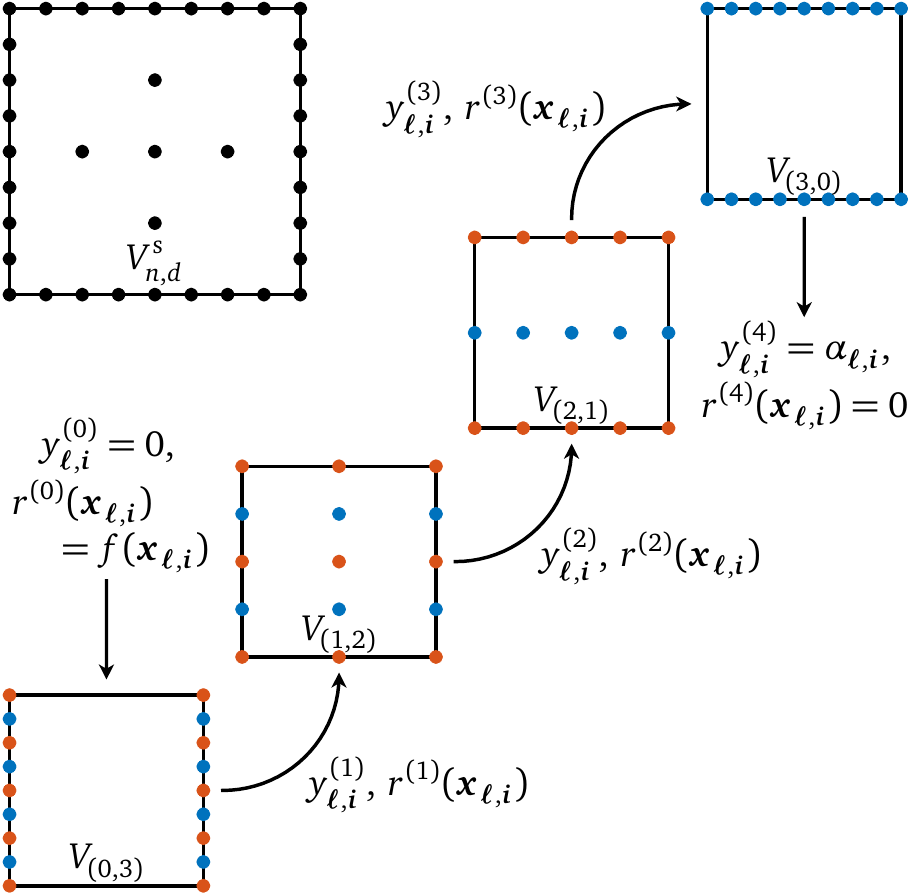}%
  \caption[%
    Hierarchization with residual interpolation%
  ]{%
    Hierarchization of function value data on the
    two-dimensional regular sparse grid
    $\regsgspace{n}{d}$ of level $n = 3$ \emph{(top left)}
    using the method of residual interpolation.
    In this figure, we use
    $\linout[(j)]{\*l,\*i} \ceq \sum_{j'=1}^j \surplus[(j')]{\*l,\*i}$
    as an abbreviation.
    The order of the nodal spaces (here: bottom left to top right)
    does not matter.
    The data $\linout[(j)]{\*l,\*i}$
    corresponding to \textcolor{C0}{blue grid points}
    will not be modified in the remaining iterations
    and, therefore, already equals the correct surpluses $\surplus{\*l,\*i}$.
    The data corresponding to \textcolor{C1}{red grid points}
    will be modified as the grid points appear in one of the remaining
    nodal grids.%
  }%
  \label{fig:residualInterpolation}%
\end{SCfigure}

An example for the application of \cref{alg:residualInterpolation}
on a two-dimensional sparse grid can be seen in
\cref{fig:residualInterpolation}.
Note that
$\surplus[(j)]{\*l,\*i} \not= 0$ can only be true if $\*l \le \*l^{(j)}$.
Therefore, if $(\*l, \*i)$ is not contained in one of the
grids that are processed in one of the remaining iterations $j+1, \dotsc, m$,
then $\linout[(j)]{\*l,\*i}$ is already equal to the correct surplus
$\surplus{\*l,\*i}$,
where $\linout[(j)]{\*l,\*i} \ceq \sum_{j'=1}^j \surplus[(j')]{\*l,\*i}$
denotes the intermediate result obtained after $j$ iterations.

\longsection{%
  \texorpdfstring{%
    Hierarchization on Spatially Adaptive Sparse Grids\\
    with Breadth-First Search%
  }{%
    Hierarchization on Spatially Adaptive Sparse Grids
    with Breadth-First Search%
  }%
}{%
  Hierarchization on Spatially Adaptive Sparse Grids
  with Breadth-First Search%
}{%
  Hierarchization with Breadth-First Search%
}
\label{sec:44spatAdaptiveBFS}

\minitoc{87mm}{8}

\noindent
Unfortunately, we cannot apply the algorithms presented in the last
sections to spatially adaptive sparse grids with
hierarchical B-splines.
The reason is that the algorithms relied on the final interpolant $\sgintp$
being a linear combination of full grid solutions $\fgintp{\*l}$,
which is only possible for dimensionally adaptive sparse grids.
Consequently, the problem of hierarchization becomes significantly
harder if we operate on spatially adaptive sparse grids.
An exception is the case of piecewise linear basis functions ($p = 1$),
where we are still able to apply the \up,
as we will show in \cref{sec:45spatAdaptiveUP}.
In this section, we study one approach to hierarchize on
spatially adaptive sparse grids,
namely transforming the hierarchical basis to so-called fundamental bases
to enable a \bfs algorithm for hierarchization.

The approach in this section has already been published
\cite{Valentin18Fundamental}.
Again, note that while B-splines are our target application,
the considerations in this chapter are fully independent of the choice of
basis functions $\basis{\*l,\*i}$,
as long as they have tensor product structure.
Although we do not state it explicitly, it is possible to employ
different types of basis functions $\basis{l_t,i_t}$ in different
dimensions, e.g., B-splines of different degrees $p_t$
to enable $p$-adaptivity.

\fillsubsectionornament
\subsection{Hierarchization with Breadth-First Search on Fundamental Bases}
\label{sec:441BFSFundamentalBases}

\paragraph{Fundamental property}

\begin{SCfigure}
  \includegraphics{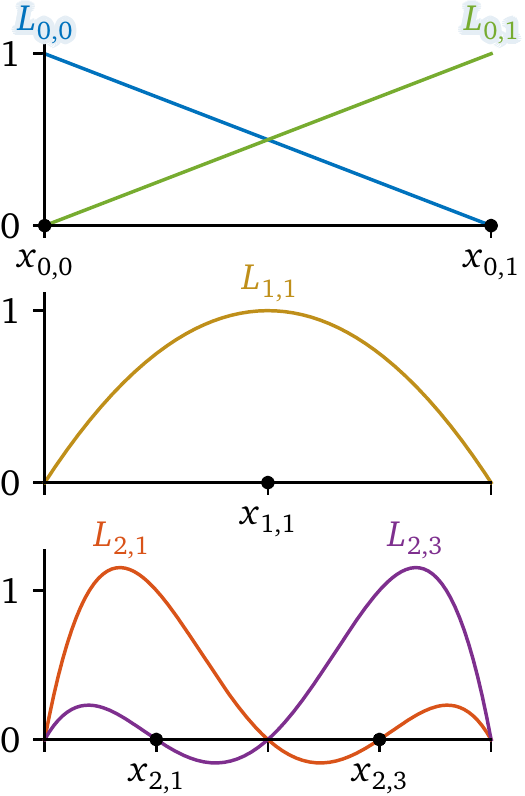}\quad%
  \includegraphics{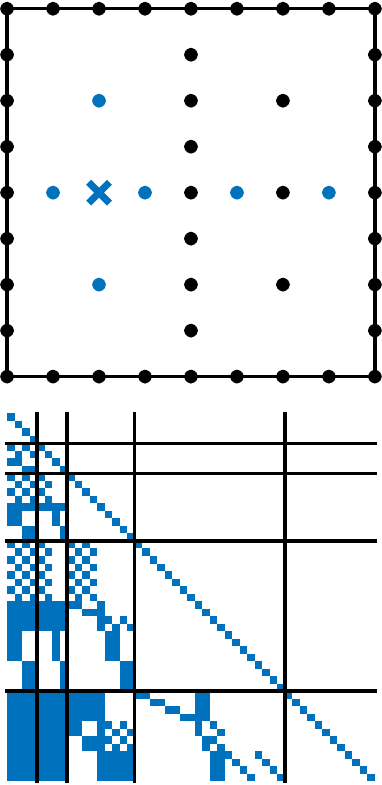}%
  \caption[%
    Fundamental property with Lagrange polynomials as fundamental basis%
  ]{%
    Fundamental property with Lagrange polynomials.\\
    \emph{Left:}
    Univariate Lagrange polynomials up to level $l = 2$.\\
    \emph{Top right:}
    Regular sparse grid $\coarseregsgset{n}{d}{1}$
    ($n = 4$, $d = 2$).
    The fundamental basis function $\fundbasis{\*l',\*i'}$ corresponding
    to the marked grid point \emph{(cross)} does not vanish
    at the \textcolor{C0}{blue points} $\gp{\*l,\*i}$
    (which satisfy \eqref{eq:fundamentalPropertyImplicationMV}).\\
    \emph{Bottom right:}
    Corresponding density pattern of $\intpmat$
    when sorting rows and columns by increasing level sum
    $\normone{\*l} = 0, \dotsc, n$ \emph{(black bars).}%
  }%
  \label{fig:fundamentalProperty}%
\end{SCfigure}

As already discussed in \cref{sec:41problem},
the main cause of the difficulty of the hierarchization with B-splines
$\bspl{\*l,\*i}{p}$ is their overlapping support
(which they need for their approximation order).
Thus, high-level B-splines $\bspl{\*l',\*i'}{p}$ do not vanish
at all coarse-level grid points $\gp{\*l,\*i}$, $\*l < \*l'$.
In the univariate case,
\pagebreak%
the idea is to transform the B-spline basis to obtain
new basis functions $\fundbasis{l',i'}\colon \clint{0, 1} \to \real$
($l' \in \natz$, $i' \in \hiset{l'}$) that satisfy
\begin{subequations}
  \label{eq:fundamentalProperty}
  \begin{alignat}{3}
    \label{eq:fundamentalProperty1}
    \fundbasis{l',i'}(\centerhphantom{\gp{l,i}}{\gp{l',i}})
    &= 0,\qquad
    &&l < l',\quad
    &&i \in \hiset{l},\\
    \label{eq:fundamentalProperty2}
    \fundbasis{l',i'}(\gp{l',i})
    &= \kronecker{i}{i'},\qquad
    &&&&i \in \hiset{l'}.
  \end{alignat}
\end{subequations}
We call \eqref{eq:fundamentalProperty} \term{fundamental property}
and functions $\fundbasis{l',i'}$
that fulfill this property \term{fundamental basis functions.}
The first \cref{eq:fundamentalProperty1} ensures that
basis functions of level $l'$ vanish at
grid points of coarser levels $l < l'$.
The second \cref{eq:fundamentalProperty2} requires the
basis functions $\fundbasis{l',i'}$
to additionally vanish at all grid points of the same level $l'$
with different index $i \not= i'$.
An example for fundamental basis functions are
the piecewise linear B-splines $\bspl{l',i'}{1}$
or the Lagrange polynomials $\lagrangepoly{l',i'}$
(see \cref{fig:fundamentalProperty}, left).
The statement that $\fundbasis{l',i'}(\gp{l',i'})$ should equal one
is not an additional restriction, if
the value $\fundbasis{l',i'}(\gp{l',i'})$ is non-zero,
since we can just replace $\fundbasis{l',i'}$ with
$\fundbasis{l',i'}/\fundbasis{l',i'}(\gp{l',i'})$ to obtain
$\fundbasis{l',i'}(\gp{l',i'}) = 1$.

\paragraph{Multivariate case}

For the multivariate case of $d \in \nat$ dimensions,
we define as usual tensor product versions
$\fundbasis{\*l',\*i'}$ of the univariate fundamental bases
$\fundbasis{l'_t,i'_t}$ ($t = 1, \dotsc, d$).
\Cref{eq:fundamentalProperty} then implies
\begin{equation}
  \label{eq:fundamentalPropertyImplicationMV}
  \fundbasis{\*l',\*i'}(\gp{\*l,\*i}) \not= 0
  \implies
  \falarge{t = 1, \dotsc, d}{
    \bracket*{(l'_t < l_t) \lor ((l'_t, i'_t) = (l_t, i_t))}
  },\quad
  (\*l, \*i), (\*l',\*i') \in \liset.
  \hspace*{-1mm}
\end{equation}
This means that every basis function
$\fundbasis{\*l',\*i'}$ can only be non-zero
at the grid points $\gp{\*l,\*i}$ that, in every dimension $t$,
have a strictly higher level $l_t$ or
the same level-index pair $(l_t, i_t)$ as the basis function.
We show an example for this relation in
\cref{fig:fundamentalProperty} (top right).

\paragraph{Triangular interpolation matrix}

The main motivation for enforcing the fundamental property
is the fact that it results in the hierarchization matrix $\intpmat$
being triangular, if the rows and columns are arranged
in the order of monotonously increasing level sum:
We assume that $k = k(\*l, \*i) \in \{1, \dotsc, N\}$ is a single
continuously enumerated index of the level-index pairs $(\*l, \*i) \in \liset$
(where $N = \setsize{\liset}$) such that
\begin{equation}
  k(\*l, \*i) \le k(\*l', \*i') \implies \normone{\*l} \le \normone{\*l'},\quad
  (\*l, \*i), (\*l', \*i') \in K,
\end{equation}
i.e., we sort the row indices $k = k(\*l, \*i)$ and
the column indices $k' = k(\*l', \*i')$ of $\intpmat$
by level sum $\normone{\cdot}$.
Consequently,
$\intpmat = (\intpmatentry{k}{k'})_{k=1,\dotsc,N,\, k'=1,\dotsc,N}$
is in lower block-triangular form:
\begin{subequations}
  \label{eq:fundamentalTriangularMV}
  \begin{alignat}{2}
    \intpmatentry{k}{k'}
    &= \fundbasis{k'}(\vgp{k})
    = 0,\quad
    &\normone{\*l}
    &< \normone{\*l'},\\
    \intertext{%
      as $\normone{\*l} < \normone{\*l'} \implies \ex{t}{l_t < l'_t}$
      and using \eqref{eq:fundamentalProperty1}.
      Additionally, the diagonal blocks are unit matrices due to%
    }
    \intpmatentry{k}{k'}
    &= \fundbasis{k'}(\vgp{k})
    \overset{(\ast)}{=} \kronecker{(\*l,\*i)}{(\*l',\*i')},\quad
    &\normone{\*l}
    &= \normone{\*l'},
  \end{alignat}
\end{subequations}
since $\normone{\*l} = \normone{\*l'}$ implies that
either $\bracket*{\ex{t}{l_t < l'_t}}$ or $\*l = \*l'$.%
\footnote{%
  Note that as specified in the list of symbols at the beginning
  of this thesis, the Kronecker delta $\kronecker{X}{Y}$
  is defined for arbitrary objects $X$ and $Y$
  that can be compared with ``$=$.''%
}
In the former case, both sides of $(\ast)$ vanish
(according to \eqref{eq:fundamentalProperty1}),
and in the latter case,
both sides equal $\kronecker{\*i}{\*i'}$
(according to \eqref{eq:fundamentalProperty2}).
Hence, $\intpmat$ is a lower-triangular matrix.
This is visualized for a two-dimensional example in
\cref{fig:fundamentalProperty} (bottom right).

\paragraph{Forward substitution}

The triangular structure of $\intpmat$ implies that
we can determine the surpluses $\surplus{\*l,\*i}$
via forward substitution:

\begin{lemma}[forward substitution]
  \label{lemma:forwardSubstitution}
  The hierarchical surpluses $\surplus{\*l,\*i}$, which
  are determined by \eqref{eq:hierarchizationSLE}
  with respect to $\fundbasis{\*l,\*i}$, satisfy
  \begin{equation}
    \surplus{\*l,\*i}
    = \fcnval{\*l,\*i} -
    \largesum{\substack{(\*l',\*i')\in\liset\\\normone{\*l'} < \normone{\*l}}}
    \surplus{\*l',\*i'} \fundbasis{\*l',\*i'}(\gp{\*l,\*i}),\quad
    (\*l, \*i) \in \liset.
  \end{equation}
\end{lemma}

\begin{proof}
  The linear system \eqref{eq:hierarchizationSLE} is given by
  \begin{subequations}
    \setlength{\abovedisplayskip}{9pt}%
    \setlength{\belowdisplayskip}{9pt}%
    \begin{align}
      \fcnval{\*l,\*i}
      &= \largesum{(\*l',\*i') \in \liset}
      \surplus{\*l',\*i'} \fundbasis{\*l',\*i'}(\gp{\*l,\*i}),\quad
      (\*l, \*i) \in \liset.\\
      \intertext{%
        According to \eqref{eq:fundamentalTriangularMV},
        all summands with $\normone{\*l'} > \normone{\*l}$ vanish
        and from the summands with $\normone{\*l'} = \normone{\*l}$,
        only the $(\*l, \*i)$-th summand remains
        with $\fundbasis{\*l,\*i}(\gp{\*l,\*i}) = 1$:%
      }
      \cdots
      &= \surplus{\*l,\*i} +
      \largesum{\substack{(\*l',\*i')\in\liset\\\normone{\*l'} < \normone{\*l}}}
      \surplus{\*l',\*i'} \fundbasis{\*l',\*i'}(\gp{\*l,\*i}).\\[-6.4em]\notag
    \end{align}
  \end{subequations}
\end{proof}

\vspace{1em}

\paragraph{Breadth-first search}

\begin{SCfigure}
  \includegraphics{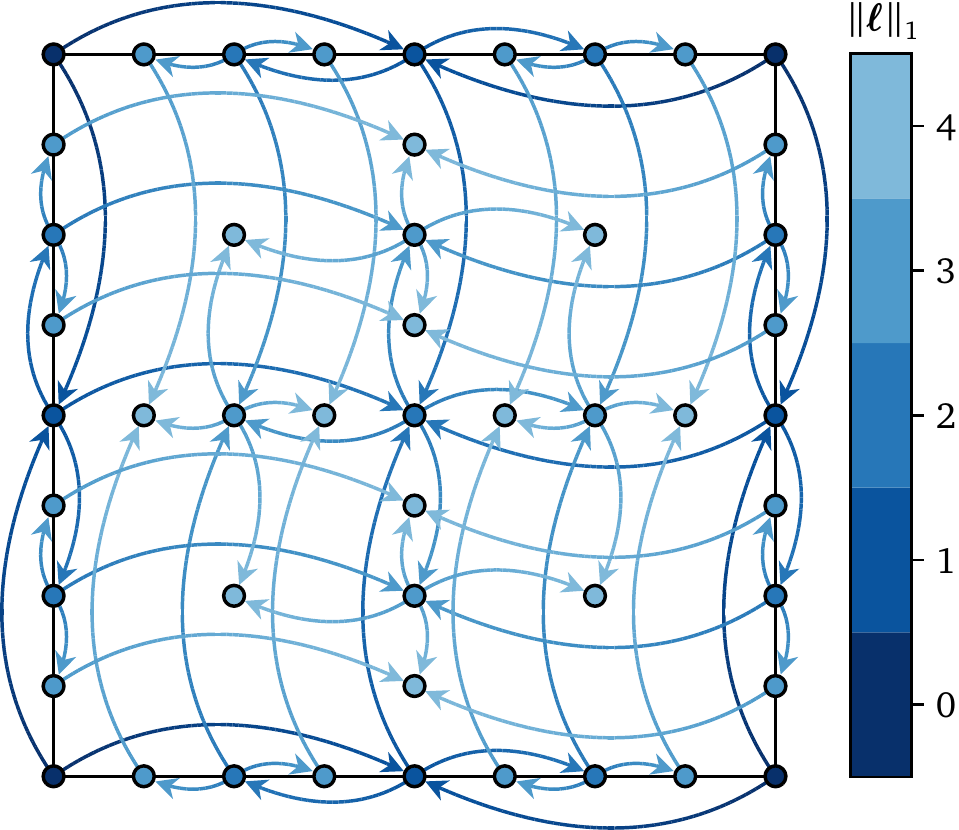}%
  \caption[%
    Sparse grid as directed acyclic graph%
  ]{%
    Ancestor relationships \emph{(arrows)} in a
    regular sparse grid $\coarseregsgset{n}{d}{1}$ \emph{(points)}
    of level $n = 4$ and dimensionality $d = 2$.
    The color indicates the level sum $\normone{\*l}$.
    Breadth-first search, as implemented in \cref{alg:BFS},
    visits all grid points $\gp{\*l,\*i}$ with level sum
    $\normone{\*l} = 0$ first,
    then those with $\normone{\*l} = 1$, and so on.%
  }%
  \label{fig:DAG}%
\end{SCfigure}

Exploiting this lemma, we formulate
a hierarchization algorithm (see \cref{alg:BFS})
that applies forward substitution by \bfs in the \dagr of
the spatially adaptive sparse grid $\sgset$.
The nodes of the \dagr are the level-index pairs $(\*l, \*i) \in \liset$.
An edge connects $(\*l, \*i)$ to $(\*l', \*i')$,
if $(\*l, \*i)$ is a direct ancestor of $(\*l', \*i')$, i.e., if
{%
  \setlength{\abovedisplayskip}{9pt}%
  \setlength{\belowdisplayskip}{9pt}%
  \begin{equation}
    \label{eq:directAncestor}
    \exlarge{t = 1, \dotsc, d}{
      \*l'_{-t} = \*l^{}_{-t},\;
      \*i'_{-t} = \*i^{}_{-t},\;
      l'_t = l^{}_t + 1,\;
      i'_t \in
      \begin{cases}
        \{1\},&l_t = 0,\\
        \{2i_t - 1, 2i_t + 1\},&l_t > 0.
      \end{cases}
    }
  \end{equation}%
}%
An example for the resulting \dagr for a regular sparse grid
in two dimensions is shown in \cref{fig:DAG}.
\begin{algorithm}
  \begin{algorithmic}[1]
    \Function{$\vlinout = \texttt{breadthFirstSearch}$}{%
      $\vlinin$, $\liset$%
    }
      \State{$\vlinout \gets \vlinin$}
      \label{line:algBFS3}
      \State{$K_\mathrm{p} \gets \{\*0\} \times \{0, 1\}^d$}
      \Comment{set of processed points}%
      \State{%
        $Q \gets \text{FIFO queue initialized with contents of } K_\mathrm{p}$%
      }
      \Comment{points to be processed}%
      \While{$Q \not= \emptyset$}
        \State{$(\*l', \*i') \gets Q.\pop()$}
        \Comment{obtain next point}%
        \For{%
          $\{(\*l, \*i) \in \liset \setminus \{(\*l', \*i')\} \mid
          \fa{t=1,\dotsc,d}{(l'_t < l_t) \lor ((l'_t, i'_t) = (l_t, i_t))}\}$%
        }
        \label{line:algBFS1}
          \State{%
            $\linout{\*l,\*i} \gets \linout{\*l,\*i} -
            \linout{\*l',\*i'} \fundbasis{\*l',\*i'}(\gp{\*l,\*i})$%
          }
          \Comment{%
            update surpluses according to \cref{lemma:forwardSubstitution}%
          }%
          \label{line:algBFS2}
        \EndFor{}
        \For{%
          $\{(\*l, \*i) \in \liset \setminus \liset_\mathrm{p} \mid
          \text{$(\*l, \*i)$ direct child of $(\*l', \*i')$}\}$%
        }
          \State{$Q.\push((\*l, \*i))$}
          \Comment{add children to queue}%
          \State{%
            $\liset_\mathrm{p} \gets
            \liset_\mathrm{p} \cup \{(\*l, \*i)\}$%
          }
          \Comment{mark as processed}%
        \EndFor{}
      \EndWhile{}
    \EndFunction{}
  \end{algorithmic}
  \caption[%
    Hierarchization with breadth-first search (BFS)%
  ]{%
    Hierarchization with breadth-first search
    on spatially adaptive sparse grids with fundamental bases.
    Inputs are
    the vector $\vlinin = (\linin{\*l,\*i})_{(\*l,\*i) \in \liset}$
    of input data (function values $\fcnval{\*l,\*i}$ at the grid points) and
    the set $\liset$ of level-index pairs of the
    sparse grid (see \eqref{eq:spatiallyAdaptiveSG}).
    The output is the vector
    $\vlinout = (\linout{\*l,\*i})_{(\*l,\*i) \in \liset}$
    of output data (hierarchical surpluses $\surplus{\*l,\*i}$).%
  }%
  \label{alg:BFS}%
\end{algorithm}
We make two assumptions for \cref{alg:BFS}.
First, $\liset$ should contain at least all $2^d$ corners of the
domain $\clint{\*0, \*1}$:
\begin{subequations}
  \label{eq:BFSAssumptions}
  \begin{equation}
    \label{eq:BFSAssumption1}
    \liset \supset \{\*0\} \times \{0, 1\}^d
    = \{(\*0, \*i) \mid \*i \in \{0, 1\}^d\}.
  \end{equation}
  Second, all grid points should be reachable from the corners:
  \begin{equation}
    \label{eq:bfsAssumption2}
    \begin{split}
      \forall_{(\*l', \*i') \in \liset}
      \exists_{m \in \natz}
      \exists_{
        (\*l^{(0)}, \*i^{(0)}), \dotsc, (\*l^{(m)}, \*i^{(m)}) \in \liset
      }\;\;
      &\bigl[(\*l^{(0)}, \*i^{(0)}) \to \dotsb \to (\*l^{(m)}, \*i^{(m)}),\\
      &\hspace{3mm} \*l^{(0)} = \*0,\;
      (\*l^{(m)}, \*i^{(m)}) = (\*l', \*i')\bigr],
    \end{split}
  \end{equation}
\end{subequations}
where ``$\to$'' is the direct ancestor relation \eqref{eq:directAncestor}.
One can also use a different initial set than the corners
of $\clint{\*0, \*1}$, e.g., when working with sparse grids
without boundary points.
In general, there are three requirements on the initial set:
First, all grid points should be reachable from this set.
Second, the grid points in the set are sorted by increasing level sum
(if the set contains grid points with different level sums).
Third, the surpluses corresponding to the initial grid points
need to be precalculated correctly
before the \texttt{\algorithmicwhile} loop in \cref{alg:BFS}.

\paragraph{Correctness}

The correctness of \cref{alg:BFS} can be shown with the following invariant:

\begin{restatable}[invariant of breadth-first-search hierarchization]{%
  proposition%
}{%
  propInvariantBFS%
}
  \label{prop:invariantBFS}
  Under the assumption \eqref{eq:BFSAssumptions},
  it holds after \pop{}\,ping all grid points with level sum $< q$
  from the queue $Q$ in \cref{alg:BFS}:
  \begin{equation}
    \label{eq:propInvariantBFS}
    \linout{\*l,\*i}
    = \fcnval{\*l,\*i} -
    \largesum{\substack{(\*l',\*i')\in\liset\\\normone{\*l'} < q}}
    \linout{\*l',\*i'} \fundbasis{\*l',\*i'}(\gp{\*l,\*i}),\quad
    (\*l, \*i) \in \liset,\;\;
    \normone{\*l} = q.
    \hspace*{-6mm}
  \end{equation}
\end{restatable}

\vspace{-0.8em}

\begin{proof}
  See \cref{sec:a133proofBFS}.
\end{proof}

\vspace{0.8em}

\begin{shortcorollary}[correctness of breadth-first-search hierarchization]
  \label{cor:algBFSCorrectness}
  \Cref{alg:BFS} is correct.
\end{shortcorollary}

\begin{proof}
  As noted in the proof of \cref{prop:invariantBFS},
  the result $\linout{\*l,\*i}$ after \pop{}ping all grid points with level sum
  $q$ as stated in \cref{prop:invariantBFS} is also the final result
  of the algorithm:
  \begin{equation}
    \linout{\*l,\*i}
    = \fcnval{\*l,\*i} -
    \largesum{\substack{(\*l',\*i')\in\liset\\\normone{\*l'} < \normone{\*l}}}
    \linout{\*l',\*i'} \fundbasis{\*l',\*i'}(\gp{\*l,\*i}),\quad
    (\*l, \*i) \in \liset.
  \end{equation}
  By \cref{lemma:forwardSubstitution},
  the correct hierarchical surpluses $\surplus{\*l,\*i}$
  satisfy the same relation.
  Inductively, $\linout{\*l,\*i}$ and $\surplus{\*l,\*i}$ must coincide.
\end{proof}

\paragraph{Complexity}

The \bfs algorithm in \cref{alg:BFS} is not as efficient as the \up:
It still needs to perform $\landauO{N^2 d}$ many univariate
basis evaluations (compared to $\landauO{Nd}$ for the \up).
However, it only needs linear space $\landauO{N}$ similar to the
\up\punctfix{.}
This is a significant advantage over directly solving the system
\eqref{eq:hierarchizationSLE} of linear equations, which typically
needs quadratic space $\landauO{N^2}$.

\subsection{Constructing Fundamental Bases}
\label{sec:442constructingFundamentalBases}

Unfortunately, the hierarchical B-splines $\bspl{l',i'}{p}$ do not
satisfy the fundamental property \eqref{eq:fundamentalProperty}.
We now focus on the construction of univariate fundamental bases
$\fundbasis{l',i'}$ starting from an
arbitrary hierarchical basis $\basis{l',i'}$.
To this end, we study two transformations
$\basis{l',i'} \mapsto \fundbasis{l',i'}$.
As usual, the multivariate case is treated with the
tensor product approach.

\paragraph{Hierarchical fundamental transformation (HFT)}

\begin{SCfigure}
  \includegraphics{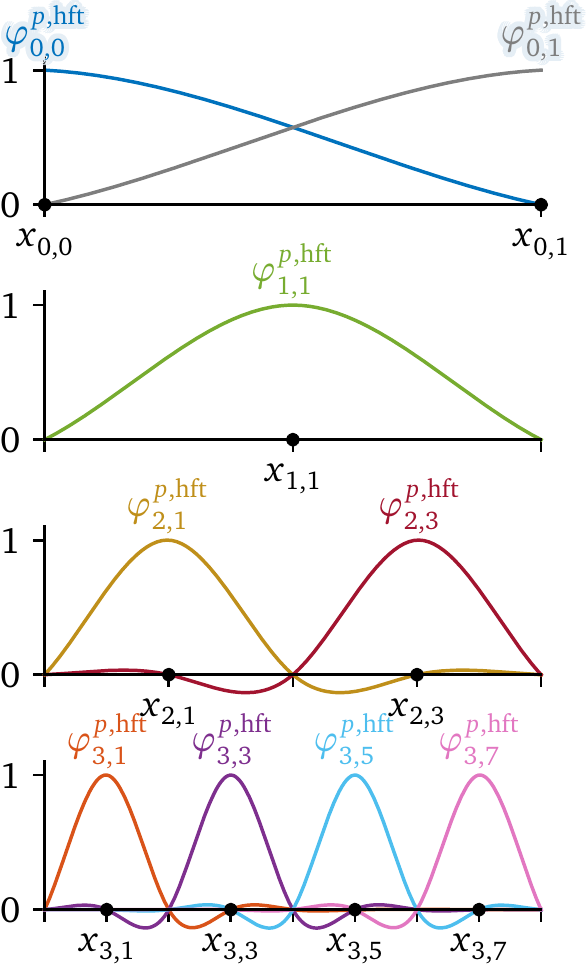}%
  \caption[%
    Hierarchical fundamental transformation on hierarchical B-splines%
  ]{%
    Resulting basis functions $\bspl[\hft]{l',i'}{p}$
    ($l' \le l$, $i' \in \hiset{l'}$)
    after applying the hierarchical fundamental transformation
    to hierarchical cubic B-splines ($p = 3$) and
    grid points $\gp{l',i'}$ \emph{(dots)} up to level $l = 3$.%
  }%
  \label{fig:hftBSpline}%
\end{SCfigure}

The canonical way to find a fundamental basis $\fundbasis{l',i'}$ is to
use a linear combination $\basis[\hft]{l',i'}$
of coarser basis functions $\basis{l'',i''}$ as an ansatz
and require that the fundamental property \eqref{eq:fundamentalProperty}
is fulfilled:
\begin{equation}
  \label{eq:hierarchicalFundamentalTransformation}
  \basis[\hft]{l',i'}
  \ceq \sum_{l''=0}^{l'} \sum_{i'' \in \hiset{l''}}
  c_{l'',i''}^{l',i'} \basis{l'',i''}
  \quad\text{such that}\quad
  \fafalarge{l = 0, \dotsc, l'}{i \in \hiset{l}}{
    \basis[\hft]{l',i'}(\gp{l,i}) = \kronecker{(l,i)}{(l',i')}.
  }
\end{equation}
This means that
$\basis[\hft]{l',i'}$ ($l' \in \natz$, $i' = 0, \dotsc, 2^{l'}$)
interpolates the data
$\{(\gp{l',i}, \kronecker{i}{i'}) \mid i = 0, \dotsc, 2^{l'}\}$.
The coefficients $c_{l'',i''}^{l',i'} \in \real$
are, in general, different for each basis function $\basis[\hft]{l',i'}$.
This complicates precomputation and storage of the $2^{l'} + 1$ coefficients,
as they have to be determined by solving a system of linear equations.
In addition, the transformation $\basis{l',i'} \mapsto \basis[\hft]{l',i'}$
does not preserve the locality of the support of the basis functions.
Consequently, $\basis[\hft]{l',i'}$ may be globally supported,
which means that we have to evaluate up to $2^{l'} + 1$ basis functions
$\basis{l'',i''}$ when evaluating $\basis[\hft]{l',i'}$ at a single point
$x \in \clint{0, 1}$.
The global support of the resulting transformed basis
for uniform hierarchical B-splines (which are locally supported)
can be seen in \cref{fig:hftBSpline}.

We call the transformation $\basis{l',i'} \mapsto \basis[\hft]{l',i'}$
\term{\hftr.}
The following proposition shows that this is only a change of basis,
as the spanned sparse grid space remains unchanged.
While the proposition is formulated for regular sparse grids,
a similar statement can be proven for the dimensionally adaptive case.

\begin{proposition}[%
  spanned sparse grid space for the HFT%
]
  \label{prop:hftSparseGridSpace}
  If $\liset \ceq \{(\*l, \*i) \mid
  \normone{\*l'} \le n,\; \*i' \in \hiset{\*l'}\}$
  is the set of level-index pairs for the regular sparse grid of level $n$
  and dimensionality $d$, then
  \begin{equation}
    \regsgspace{n}{d}
    = \regsgspace[\hft]{n}{d}
    \ceq \spn\{\basis[\hft]{\*l',\*i'} \mid (\*l', \*i') \in \liset\}.
  \end{equation}
\end{proposition}

\begin{proof}
  We have $\regsgspace{n}{d} \supset \regsgspace[\hft]{n}{d}$ as
  $\basis[\hft]{\*l',\*i'} \in \regsgspace{n}{d}$
  for all $(\*l', \*i') \in \liset$:
  \begin{equation}
    \basis[\hft]{\*l',\*i'}
    = \prod_{t=1}^d \, \sum_{l''_t=0}^{l'_t} \, \sum_{i''_t \in \hiset{l''_t}}
    c_{l''_t,i''_t}^{l'_t,i'_t} \basis{l''_t,i''_t}
    = \sum_{\*l''=\*0}^{\*l'} \, \sum_{\*i'' \in \hiset{\*l''}}
    c_{\*l'',\*i''}^{\*l',\*i'} \basis{\*l'',\*i''}
    \in \regsgspace{n}{d},\quad
    c_{\*l'',\*i''}^{\*l',\*i'}
    \ceq \prod_{t=1}^d c_{l''_t,i''_t}^{l'_t,i'_t}.
  \end{equation}
  
  To prove that $\regsgspace{n}{d} \subset \regsgspace[\hft]{n}{d}$,
  we show that the dimension of $\regsgspace[\hft]{n}{d}$
  matches $\dim \regsgspace{n}{d} = \setsize{\liset}$.
  It suffices to show that
  the functions $\basis[\hft]{\*l',\*i'}$ ($(\*l', \*i') \in \liset$),
  are linearly independent.
  Let $\surplus{\*l',\*i'} \in \real$ be with
  $\sum_{(\*l', \*i') \in \liset}
  \surplus{\*l',\*i'} \basis[\hft]{\*l',\*i'} \equiv 0$.
  By evaluating at $\gp{\*l,\*i}$ ($(\*l, \*i) \in \liset$), we obtain
  \begin{equation}
    \largesum{(\*l', \*i') \in \liset}
    \surplus{\*l',\*i'} \basis[\hft]{\*l',\*i'}(\gp{\*l,\*i}) = 0,\quad
    (\*l, \*i) \in \liset.
  \end{equation}
  This is a lower triangular system according to
  \eqref{eq:fundamentalTriangularMV},
  which implies $\surplus{\*l',\*i'} = 0$ for all $(\*l, \*i) \in \liset$.
  Hence, the functions $\basis[\hft]{\*l',\*i'}$ ($(\*l', \*i') \in \liset$)
  are linearly independent.
\end{proof}

\paragraph{Translation-invariant fundamental transformation (TIFT)}

Another disadvantage of the \hftr is that it does not preserve
the so-called \term{translation invariance} of the original basis.
A~basis $\basis{l,i}$ ($l \in \natz$, $i = 0, \dotsc, 2^l$)
is translation-invariant, if there is a \term{parent function}
$\parentfcn\colon \real \to \real$
such that
\begin{equation}
  \label{eq:translationInvariance}
  \basis{l,i}(x)
  = \parentfcn(\tfrac{x}{\ms{l}} - i),\quad
  l \in \natz,\;\;
  i = 0, \dotsc, 2^l,\;\;
  x \in \clint{0, 1}.
\end{equation}
The fact that the \hftr does not preserve translation invariance means
that for each basis function $\basis[\hft]{l',i'}$, we have to
calculate its individual $2^{l'} + 1$ coefficients $c_{l'',i''}^{l',i'}$.

To solve this problem,
we use a similar ansatz as for the \hftr,
but we replace the hierarchical basis functions $\basis{l'',i''}$
($l'' = 0, \dotsc, l'$,\, $i'' \in \hiset{l''}$)
with nodal basis functions $\basis{l',i''}$ and
allow general integer indices $i'' \in \integer$:
\begin{equation}
  \label{eq:tiftDefinition}
  \basis[\tift]{l',i'}
  \ceq \sum_{i'' \in \integer} c_{i''}^{l',i'} \basis{l',i''}
  \quad\text{such that}\quad
  \falarge{i \in \hiset{l'}}{
    \basis[\tift]{l',i'}(\gp{l',i}) = \kronecker{i}{i'},
  }
\end{equation}
where $l' \in \natz$, $i' = 0, \dotsc, 2^{l'}$, and
$c_{i''}^{l',i'} \in \real$.
We have to make three assumptions for \eqref{eq:tiftDefinition} to make sense:

\begin{itemize}
  \item
  The functions $\basis{l',i''}$ have to be defined for integer indices
  $i'' \in \integer$, i.e.,
  the functions $\basis{l',i''}\colon \clint{0, 1} \to \real$
  must also exist for $i'' < 0$ or $i'' > 2^{l'}$.
  This is the case for translation-invariant bases
  $\basis{l',i''}$ (such as B-splines $\bspl{l',i''}{p}$),
  as they can be generalized to $i'' \in \integer$
  via \cref{eq:translationInvariance}.
  
  \item
  The set
  \begin{equation}
    \relindexset{l'}
    \ceq \braced*{
      i'' \in \integer \bigm\vert
      \restrictfcn{\basis{l',i''}}{\clint{0, 1}} \not\equiv 0
    },\quad
    l' \in \natz,
  \end{equation}
  of relevant indices should be finite,
  so that in each point $x \in \clint{0, 1}$
  only a finite number of basis functions $\basis{l',i''}$ of level $l'$
  is non-zero.
  This means that the series in \eqref{eq:tiftDefinition} collapses to a
  finite sum over $i'' \in \relindexset{l'}$.
  The condition is met for compactly supported and translation-invariant
  basis functions such as B-splines $\bspl{l',i''}{p}$.
  For $d \in \nat$ dimensions and $\*l' \in \natz^d$, we define
  $\relindexset{\*l'} \ceq
  \relindexset{l'_1} \times \dotsb \times \relindexset{l'_d}$.
  
  \item
  The coefficients $c_{i''}^{l',i'}$, such that
  \eqref{eq:tiftDefinition} holds, exist and are uniquely determined.
\end{itemize}

\vspace{1em}

\noindent
Let $\basis{l',i''}$ be translation-invariant and
let $l' \in \natz$ and $i' = 0, \dotsc, 2^{l'}$ be arbitrary.
Then we have
\begin{equation}
  \label{eq:tiftParentFunctionDerivation}
  \basis[\tift]{l',i'}(x)
  = \sum_{i'' \in \integer} c_{i''}^{l',i'} \basis{l',i''}(x)
  \!\overset{\eqref{eq:translationInvariance}}{=}\!
  \sum_{i'' \in \integer} c_{i''}^{l',i'}
  \parentfcn(\tfrac{x}{\ms{l'}} - i'')
  = \sum_{i'' \in \integer}
  c_{i'+i''}^{l',i'} \parentfcn((\tfrac{x}{\ms{l'}} - i') - i''),
\end{equation}
where the function defined by
$\sum_{i'' \in \integer} c_{i'+i''}^{l',i'} \parentfcn({\cdot} - i'')$
satisfies
\begin{subequations}
  \label{eq:tiftParentFunctionConditions}
  \begin{align}
    \sum_{i'' \in \integer} c_{i'+i''}^{l',i'} \parentfcn(i - i'')
    &= \sum_{i'' \in \integer} c_{i'+i''}^{l',i'}
    \parentfcn((\tfrac{\gp{l',i+i'}}{\ms{l'}} - i') - i'')
    \!\overset{\eqref{eq:tiftParentFunctionDerivation}}{=}\!
    \basis[\tift]{l',i'}(\gp{l',i+i'})
    \!\overset{\eqref{eq:tiftDefinition}}{=}\!
    \kronecker{i'}{i+i'}\\
    &= \kronecker{i}{0},\quad
    i \in \integer.
  \end{align}
\end{subequations}
Due to the translation invariance and the third assumption in the list above,
there is only a single set of coefficients of $\parentfcn({\cdot} - i'')$
such that \eqref{eq:tiftParentFunctionConditions} holds.
This means that $c_{i'+i''}^{l',i'}$ does not depend on $l'$ or $i'$,
if $\basis{l',i'}$ is translation-invariant.
Consequently, if we set $c_{i''} \ceq c_{i'+i''}^{l',i'}$
for some arbitrary $l' \in \natz$ and $i' = 0, \dotsc, 2^{l'}$,
then $\parentfcn[\tift]\colon \real \to \real$ defined by
\begin{equation}
  \label{eq:tiftParentFunctionDefinition}
  \parentfcn[\tift](x)
  \ceq \sum_{i'' \in \integer} c_{i''} \parentfcn(x - i''),\quad
  x \in \clint{0, 1},
\end{equation}
is a parent function of $\basis[\tift]{l',i'}$ that satisfies
\begin{equation}
  \label{eq:tiftParentFunctionInterpolation}
  \falarge{i \in \integer}{
    \parentfcn[\tift](i) = \kronecker{i}{0}.
  }
\end{equation}
The fact that $\parentfcn[\tift]$ is the parent function of
$\basis[\tift]{l',i'}$ easily follows from
\eqref{eq:tiftParentFunctionDerivation}
as the \rhs is exactly $\parentfcn[\tift](\tfrac{x}{\ms{l'}} - i')$,
as required by \eqref{eq:translationInvariance}.
This shows that the transformation $\basis{l',i'} \mapsto \basis[\tift]{l',i'}$
preserves translation-invariance.
Therefore, we call the transformation
$\basis{l',i'} \mapsto \basis[\tift]{l',i'}$ \term{\tiftr.}

\vspace{2em}

In contrast to the \hftr, the \tiftr is only a change of basis
if we consider the extended nodal spaces that also include
basis functions with indices outside $\{0, \dotsc, 2^{l'}\}$.
This is the statement of the following proposition
(generalized to the $d$-variate case).
Note that although the proposition involves basis functions
$\basis{\*l',\*i'}$ and $\basis[\tift]{\*l',\*i'}$
outside the domain $\clint{0, 1}$
(in the sense that $\gp{\*l',\*i'} \notin \clint{0, 1}$),
we still restrict all functions to $\clint{0, 1}$.
We cannot formulate an equivalent version of \thmref{prop:hftSparseGridSpace},
as it might be that, in one dimension,
$\basis{l',i''}$ ($i'' < 0$ or $i'' > 2^{l'}$)
is not contained in the nodal space $\ns{l'}$.

\pagebreak

\begin{proposition}[%
  spanned nodal space for the TIFT%
]
  \label{prop:tiftNodalSpace}
  We have
  \begin{equation}
    \hspace*{-10mm}
    \spn\{\basis{\*l',\*i'} \mid \*i' \in \relindexset{\*l'}\}
    =: \ns[\ext]{\*l'}
    = \ns[\tift,\ext]{\*l'}
    \ceq \spn\{\basis[\tift]{\*l',\*i'} \mid \*i' \in \relindexset{\*l'}\},\quad
    \*l' \in \natz^d.
    \hspace*{-20mm}
  \end{equation}
\end{proposition}

\vspace{1em}

\begin{proof}
  We have $\ns[\ext]{\*l'} \supset \ns[\tift,\ext]{\*l'}$ as
  $\basis[\tift]{\*l',\*i'} \in \ns[\ext]{\*l'}$
  for all $\*i' \in \relindexset{\*l'}$:
  \begin{equation}
    \basis[\tift]{\*l',\*i'}
    = \prod_{t=1}^d \sum_{i''_t \in \relindexset{l'_t}}
    c_{i''_t}^{l'_t,i'_t} \basis{l'_t,i''_t}
    = \sum_{\*i'' \in \relindexset{\*l'}}
    c_{\*i''}^{\*l',\*i'} \basis{\*l',\*i''}
    \in \ns[\ext]{\*l'},\quad
    c_{\*i''}^{\*l',\*i'}
    \ceq \prod_{t=1}^d c_{i''_t}^{l'_t,i'_t}.
  \end{equation}
  
  To prove that
  $\ns[\ext]{\*l'} \subset \ns[\tift,\ext]{\*l'}$,
  we show that the dimensions of the two spaces match.
  As before, it suffices to show that
  the functions $\basis[\tift]{\*l',\*i'}$ ($\*i' \in \relindexset{\*l'}$)
  are linearly independent.
  Let $\surplus{\*l',\*i'} \in \real$ be with
  $\sum_{\*i' \in \relindexset{\*l'}}
  \surplus{\*l',\*i'} \basis[\tift]{\*l',\*i'} \equiv 0$.
  By evaluating at $\gp{\*l',\*i}$ ($\*i \in \relindexset{\*l'}$), we obtain
  \begin{equation}
    0
    = \sum_{\*i' \in \relindexset{\*l'}} \surplus{\*l',\*i'}
    \underbrace{\basis[\tift]{\*l',\*i'}(\gp{\*l',\*i})}_{=\kronecker{\*i}{\*i'}}
    = \surplus{\*l',\*i},\quad
    \*i \in \relindexset{\*l'},
  \end{equation}
  i.e., all coefficients $\surplus{\*l',\*i'}$ must vanish.%
  \footnote{%
    Note that we have to allow evaluations outside the
    domain $\clint{0, 1}$ for this step.
    However, this is feasible for proving the linear independence
    of $\basis[\tift]{\*l',\*i'}$, since we can just restrict the functions
    after showing that the extended nodal spaces equal.%
  }
  Hence, the functions $\basis[\tift]{\*l',\*i'}$ ($\*i' \in \relindexset{\*l'}$)
  are linearly independent.
\end{proof}

\fillsubsectionornament
\subsection{Hierarchical Fundamental Splines}
\label{sec:443fundamentalSplines}

\paragraph{Definition}

We now apply the translation-invariant fundamental transformation
to hierarchical B-splines $\bspl{l,i}{p}$ of degree $p$.
The parent function $\parentbspl{p}\colon \real \to \real$ of B-splines
and the set $\relindexset[p]{l}$ of relevant indices for level $l \in \natz$
are given by
\begin{equation}
  \parentbspl{p}(x)
  \ceq \cardbspl{p}(x + \tfrac{p+1}{2}),\qquad
  \relindexset[p]{l}
  \ceq \{-\tfrac{p-1}{2},\;
  -\tfrac{p-1}{2} + 1,\;
  \dotsc,\;
  2^l + \tfrac{p-1}{2}\},
\end{equation}
respectively.
According to \eqref{eq:tiftParentFunctionInterpolation},
the coefficients $\fundsplcoeff{k}{p} \in \real$ of
the transformed parent function $\parentfcn[p,\tift]$
in \cref{eq:tiftParentFunctionDefinition} are determined by a
\pagebreak%
bi-infinite-dimensional system of linear equations:
\begin{equation}
  \label{eq:fundamentalSplineSLE}
  \newcommand*{\vala}{\cardbspl{p}(\tfrac{p+1}{2} - 2)}
  \newcommand*{\valb}{\cardbspl{p}(\tfrac{p+1}{2} - 1)}
  \newcommand*{\valc}{\cardbspl{p}(\tfrac{p+1}{2})}
  \newcommand*{\vald}{\cardbspl{p}(\tfrac{p+1}{2} + 1)}
  \newcommand*{\vale}{\cardbspl{p}(\tfrac{p+1}{2} + 2)}
  \newcommand*{\ddotsr}{\rotatebox{16}{$\ddots$}}
  \newcommand*{\raiseentry}[1]{\raisebox{0.3em}{$#1$}}
  \paren*{
    \begin{array}{%
      @{}>{\raggedright\arraybackslash$}m{14mm}<{$}@{}%
      *{3}{%
        @{}>{\centering\arraybackslash$}m{23mm}<{$}@{}%
      }%
      @{}>{\raggedleft\arraybackslash$}m{14mm}<{$}@{}%
    }
      \ddotsr&\ddotsr&\ddotsr&       &       \\
      \ddotsr&\valc  &\valb  &\vala  &       \\
      \ddotsr&\vald  &\valc  &\valb  &\ddotsr\\
             &\vale  &\vald  &\valc  &\ddotsr\\
             &       &\ddotsr&\ddotsr&\ddotsr
    \end{array}
  }
  \cdot
  \begin{pmatrix}
    \vphantom{\ddotsr}\raiseentry{\vdots}\\
    \vphantom{\ddotsr}\raiseentry{\fundsplcoeff{-1}{p}}\\
    \vphantom{\ddotsr}\raiseentry{\fundsplcoeff{0}{p}}\\
    \vphantom{\ddotsr}\raiseentry{\fundsplcoeff{1}{p}}\\
    \vphantom{\ddotsr}\raiseentry{\vdots}
  \end{pmatrix}
  =
  \begin{pmatrix}
    \vphantom{\ddotsr}\raiseentry{\vdots}\\
    \vphantom{\ddotsr}\raiseentry{0}\\
    \vphantom{\ddotsr}\raiseentry{1}\\
    \vphantom{\ddotsr}\raiseentry{0}\\
    \vphantom{\ddotsr}\raiseentry{\vdots}
  \end{pmatrix}.
\end{equation}
As in each row only $p$ entries are non-zero,
the system matrix is a symmetric banded Toeplitz matrix%
\footnote{%
  The entries $a_{k,j}$ of a Toeplitz matrix $\mat{A}$
  solely depend on $k - j$, i.e.,
  $a_{k,j} = c_{k-j}$ for some vector $\*c$.%
}.
One can show that the linear system \eqref{eq:fundamentalSplineSLE}
is uniquely solvable:

\begin{theorem}[unique existence of fundamental spline coefficients]
  \label{thm:fundamentalSplineExistence}
  \usenotation{zzzzfs}
  The system \eqref{eq:fundamentalSplineSLE} has a unique solution
  $(\fundsplcoeff{k}{p})_{k \in \integer}$ and the corresponding parent function
  $\parentfundspl{p}\colon \real \to \real$ defined by
  $\parentfundspl{p}(x) \ceq
  \sum_{k \in \integer} \fundsplcoeff{k}{p} \parentbspl{p}(x - k)$
  satisfies
  \begin{equation}
    \exfalarge{\beta_p, \gamma_p \in \posreal}{x \in \real}{
      \abs{\parentfundspl{p}(x)}
      \le \beta_p \cdot (\gamma_p)^{-\abs{x}}.
    }
  \end{equation}
\end{theorem}

\begin{proof}
  See Theorems 1 and 2 in \cite{Schoenberg72Cardinal}.
\end{proof}

\begin{figure}
  \subcaptionbox{%
    $p = 1$%
  }[72mm]{%
    \includegraphics{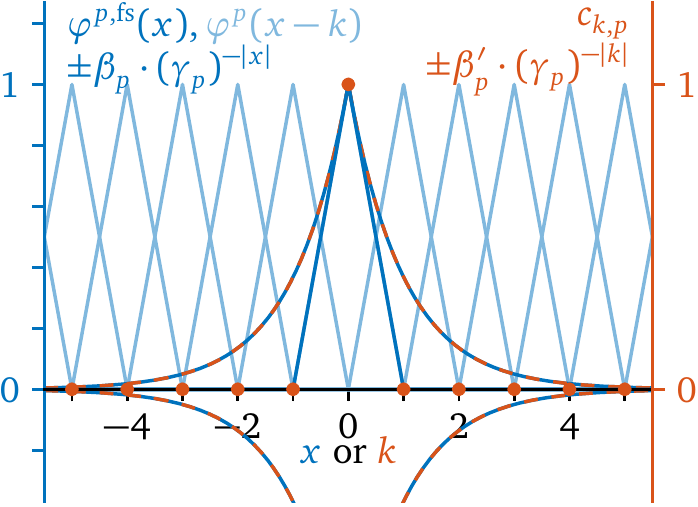}%
  }%
  \hfill%
  \subcaptionbox{%
    $p = 3$%
  }[72mm]{%
    \includegraphics{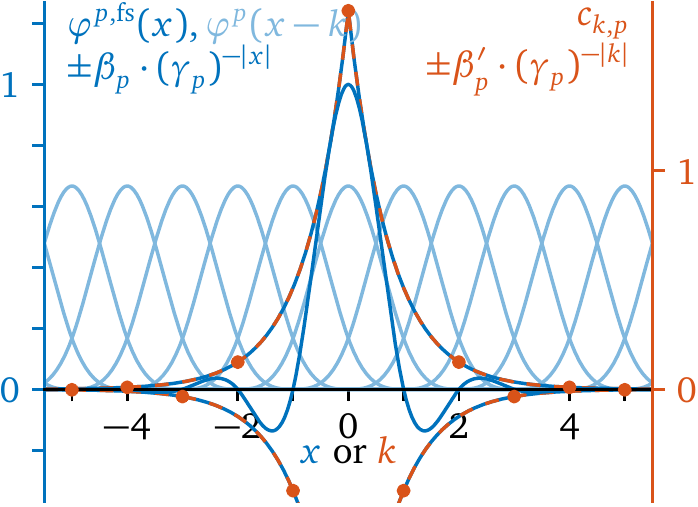}%
  }\\[4mm]%
  \subcaptionbox{%
    $p = 5$%
  }[72mm]{%
    \includegraphics{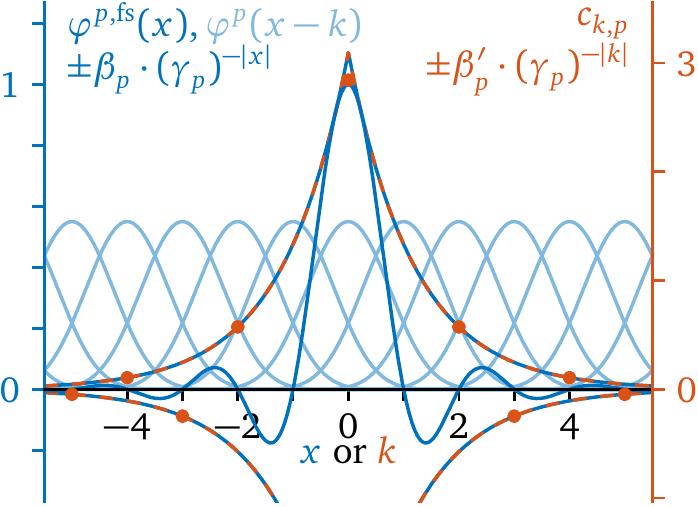}%
  }%
  \hfill%
  \subcaptionbox{%
    $p = 7$%
  }[72mm]{%
    \includegraphics{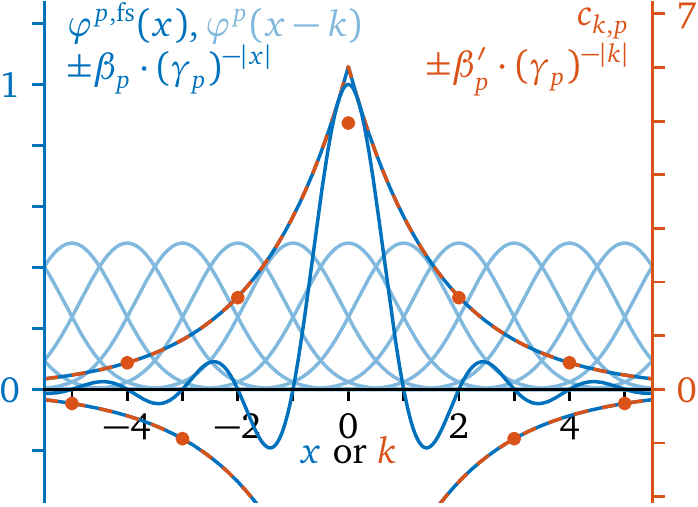}%
  }%
  \caption[%
    Fundamental splines and their B-spline coefficients%
  ]{%
    The fundamental spline $\parentfundspl{p}$ \emph{\textcolor{C0}{(blue)}}
    is a linear combination of cardinal B-splines $\cardbspl{p}({\cdot} - k)$,
    $k \in \integer$ \emph{\textcolor{C0!50}{(light blue)},}
    \vspace{-0.1em}%
    with coefficients $\fundsplcoeff{k}{p}$ \emph{\textcolor{C1}{(red points)}.}
    The absolute values of the fundamental spline $\parentfundspl{p}$ and
    its coefficients are bounded by a multiple of $(\gamma_p)^{-\abs{k}}$
    \emph{(\textcolor{C0}{blue}-\textcolor{C1}{red}-dashed line).}
    The axis for $\fundsplcoeff{k}{p}$ (on the right side) is scaled such that
    both bounding functions are on top of each other.%
  }%
  \label{fig:fundamentalSpline}%
\end{figure}

The function $\parentfundspl{p}$ from \cref{thm:fundamentalSplineExistence}
is well-known as the \term{fundamental spline} of degree~$p$
\multicite{Schoenberg72Cardinal,Schoenberg73Cardinal}.
Applications of fundamental splines are interpolation and
the definition of spline wavelets \cite{Chui92Introduction}.
The fundamental splines $\parentfundspl{p}$ of low degrees $p$ and
their bounding functions $\beta_p \cdot (\gamma_p)^{-\abs{x}}$
are plotted in \cref{fig:fundamentalSpline}.

\paragraph{Definition of hierarchical fundamental splines}

\begin{SCfigure}
  \includegraphics{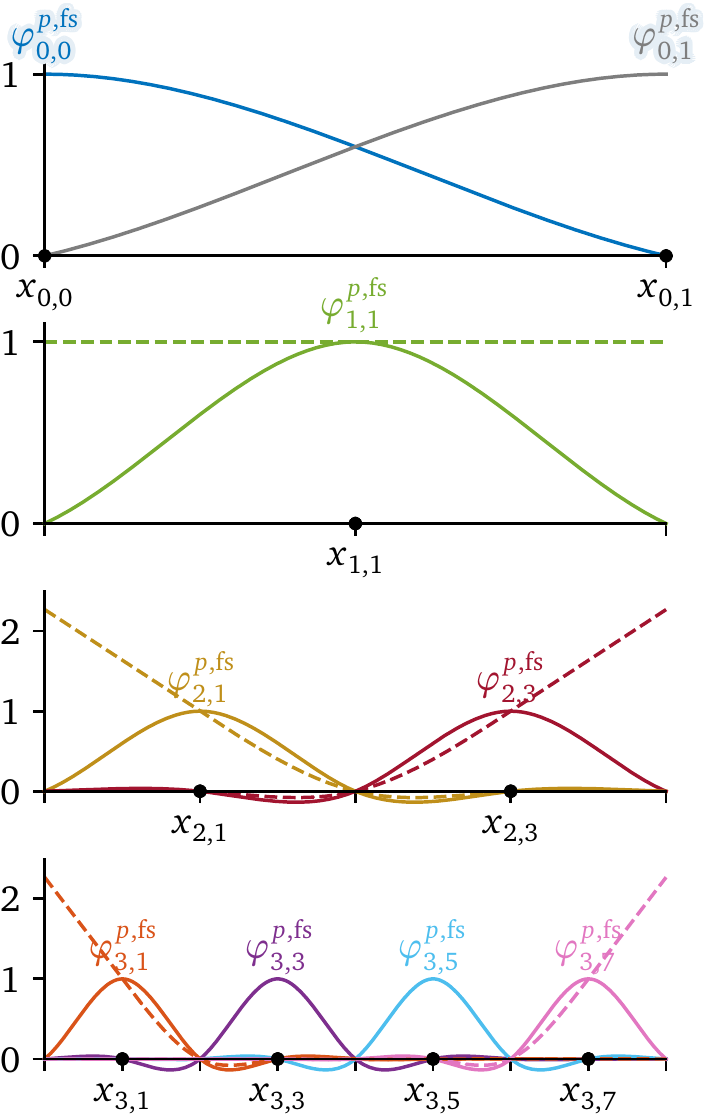}%
  \caption[%
    Hierarchical fundamental splines%
  ]{%
    Hierarchical cubic fundamental splines
    $\bspl[\fs]{l',i'}{p}$
    ($l' \le l$, $i' \in \hiset{l'}$, $p = 3$),
    their modified versions $\bspl[\fs,\modified]{l',i'}{p}$
    \emph{(dashed),} and
    grid points $\gp{l',i'}$ \emph{(dots)} up to level $l = 3$.%
  }%
  \label{fig:hierarchicalFundamentalSpline}%
\end{SCfigure}

The fundamental spline $\parentfundspl{p}$ defines
hierarchical fundamental spline functions
$\bspl[\fs]{l,i}{p}\colon \clint{0, 1} \to \real$ via
\cref{eq:translationInvariance}, i.e.,
\begin{equation}
  \bspl[\fs]{l,i}{p}(x)
  \ceq \parentfundspl{p}(\tfrac{x}{\ms{l}} - i),\quad
  l \in \natz,\;\;
  i = 0, \dotsc, 2^l,\;\;
  x \in \clint{0, 1}.
\end{equation}
The hierarchical cubic fundamental spline basis is depicted in
\cref{fig:hierarchicalFundamentalSpline}.
As usual, we define $d$-variate hierarchical fundamental splines
as tensor products of their univariate counterparts.
According to \thmref{prop:splineSpace} and \thmref{prop:tiftNodalSpace},
the common extended nodal space
$\nsbspl[\ext]{\*l}{\*p} = \nsbspl[\fs,\ext]{\*l}{\*p}$
is equal to the spline space $\wholesplspace{\*l}{\*p}$
defined by the Cartesian product of
knot sequences of the form \eqref{eq:fullGridKnots},
i.e., the space of all splines of degree $\*p$ on the full grid of level $\*l$.

The B-spline coefficients $(\fundsplcoeff{k}{p})_{k \in \integer}$ of the
fundamental spline $\parentfundspl{p}$ decay with the same rate
as the fundamental spline itself due to the stability of the B-spline basis
\cite{Hoellig13Approximation}, i.e.,
\begin{equation}
  \label{eq:fundamentalSplineCoefficientsDecay}
  \abs{\fundsplcoeff{k}{p}}
  \le \beta'_p \cdot (\gamma_p)^{-\abs{k}},\quad
  k \in \integer,
\end{equation}
for some $\beta'_p > 0$ independent of $k$,
which is also shown in \cref{fig:fundamentalSpline}.
For $p > 1$,
there is a surprising relationship between the optimal (i.e., largest)
decay rate $\gamma_p$ and the polynomial
$\sum_{k=1}^p \cardbspl{p}(k) x^{k-1}$,
whose coefficients are the values of the
cardinal B-spline $\cardbspl{p}$ at its inner knots:
The decay rate is given by the absolute value of the largest root smaller
than $-1$ of said polynomial
(see \multicite{Chui92Introduction,Schoenberg73Cardinal}).

\vspace*{\fill}

Due to \eqref{eq:fundamentalSplineCoefficientsDecay},
we may solve the system \eqref{eq:fundamentalSplineSLE}
of linear equations approximately,
if we symmetrically truncate the linear system to
\pagebreak
$2\fundsplcutoff{p} - 1$ rows and columns
and set $\fundsplcoeff{k}{p} \ceq 0$ for all $\abs{k} \ge \fundsplcutoff{p}$,
where $\fundsplcutoff{p} \in \nat$ is a truncation index.
Note that we only have to perform $p + 1$ cardinal B-spline evaluations
to evaluate $\parentfundspl{p}$ once.
In \cref{tbl:fundamentalSplineDecay}, we list the decay rates $\gamma_p$,
the factors $\beta_p$ and $\beta'_p$, and the truncation indices
$\fundsplcutoff{p}$ for different $p$.

\begin{table}
  \setnumberoftableheaderrows{1}%
  \begin{tabular}{%
    >{\kern\tabcolsep}=l<{\kern5mm}*{8}{+c}<{\kern\tabcolsep}%
  }
    \toprulec
    \headerrow
    $p$&$1$&$3$&$5$&$7$&$9$&$11$&$13$&$15$\\
    \midrulec
    $\gamma_p$&$2.718$&$3.732$&$2.322$&$1.868$&$1.645$&$1.512$&$1.425$&$1.363$\\
    $\beta_p$&$1$&$1.241$&$1.104$&$1.058$&$1.037$&$1.026$&$1.019$&$1.014$\\
    $\beta'_p$&$1$&$1.732$&$3.095$&$6.016$&$12.27$&$25.82$&$55.56$&$121.6$\\
    $\fundsplcutoff{p}$&$1$&$18$&$29$&$40$&$52$&$64$&$77$&$90$\\
    \bottomrulec
  \end{tabular}%
  \caption[%
    Decay rates of fundamental splines%
  ]{%
    Optimal decay rates $\gamma_p$ and corresponding factors
    $\beta_p$ and $\beta'_p$ for the bound functions of
    the fundamental spline $\parentfundspl{p}$ and
    its coefficients $\fundsplcoeff{k}{p}$, i.e.,
    $\fa{x \in \real}{\abs{\parentfundspl{p}(x)} \le \beta_p (\gamma_p)^{-\abs{x}}}$
    and
    $\fa{k \in \integer}{\abs{\fundsplcoeff{k}{p}} \le \beta'_p (\gamma_p)^{-\abs{k}}}$
    (approximated values).
    The truncation indices $\fundsplcutoff{p}$ are the smallest numbers such that
    $\fa{\abs{k} \ge \fundsplcutoff{p}}{\abs{\fundsplcoeff{k}{p}} < 10^{-10}}$.%
  }%
  \label{tbl:fundamentalSplineDecay}%
\end{table}

\subsection{Modified Hierarchical Fundamental Splines}
\label{sec:444modifiedFundamentalSplines}

\begin{figure}
  \subcaptionbox{%
    $p = 3$%
  }[48mm]{%
    \includegraphics{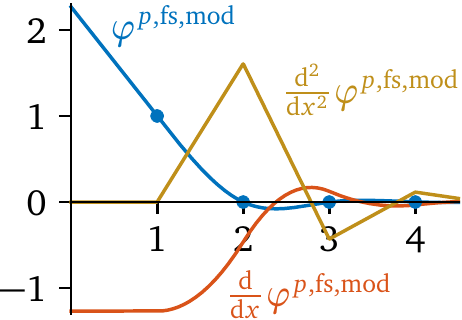}%
  }%
  \hfill%
  \subcaptionbox{%
    $p = 5$%
  }[48mm]{%
    \includegraphics{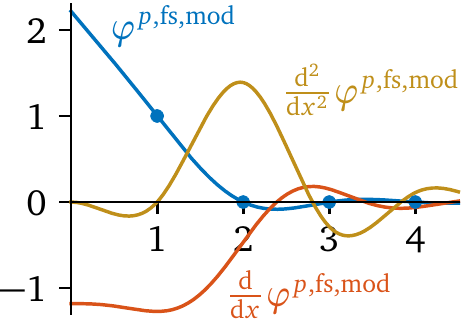}%
  }%
  \hfill%
  \subcaptionbox{%
    $p = 7$%
  }[48mm]{%
    \includegraphics{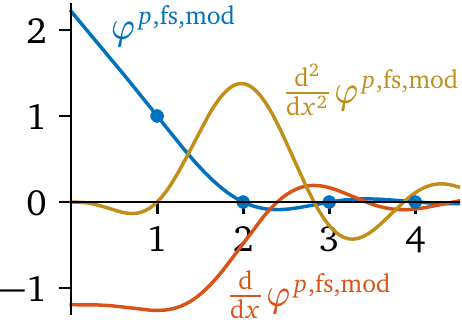}%
  }%
  \caption[%
    Modified fundamental spline and its derivatives%
  ]{%
    Modified fundamental spline $\parentfundspl[\modified]{p}$
    \emph{\textcolor{C0}{(blue)}} together with its
    first \emph{\textcolor{C1}{(red)}} and
    second \emph{\textcolor{C2}{(brown)}} derivatives
    and the function value interpolation conditions from
    \cref{eq:modifiedFundamentalSplineConditions}
    \emph{\textcolor{C0}{(blue dots)}.}
    For $p = 3$, the second derivative vanishes on
    $\clint{0, 1}$.
    For higher degrees $p > 3$, the second derivative is close to zero
    on this interval, vanishing at $x = 0$.%
  }%
  \label{fig:modifiedFundamentalSpline}%
\end{figure}

Similar to the B-spline bases introduced in \cref{chap:30BSplines},
it is possible to define a modified version of the
hierarchical fundamental spline basis to obtain reasonable
boundary values when working with sparse grids without boundary points.
The definition of the modified fundamental spline
$\bspl[\fs,\modified]{l,i}{p}\colon \clint{0, 1} \to \real$ of
level $l \in \nat$, index $i \in \hiset{l}$, and degree $p$ is
defined as follows:
\begin{equation}
  \bspl[\fs,\modified]{l,i}{p}(x)
  \ceq
  \begin{cases}
    1,&
    l = 1,\quad i = 1,\\
    \parentfundspl[\modified]{p}(\tfrac{x}{\ms{l}}),&
    l \ge 2,\quad i = 1,\\
    \bspl[\fs]{l,i}{p}(x),&
    l \ge 2,\quad i \in \hiset{l} \setminus \{1, 2^l - 1\},\\
    \bspl[\fs,\modified]{l,1}{p}(1 - x),&
    l \ge 2,\quad i = 2^l - 1,
  \end{cases}
\end{equation}
where $\parentfundspl[\modified]{p}$ is a linear combination
\begin{equation}
  \parentfundspl[\modified]{p}\colon \nonnegreal \to \real,\quad
  \parentfundspl[\modified]{p}(x)
  \ceq \largesum[\infty]{k=1-(p+1)/2}
  \fundsplcoeff[\modified]{k}{p} \parentbspl{p}(x - k),
\end{equation}
whose coefficients $\fundsplcoeff[\modified]{k}{p} \in \real$
are chosen such that
\begin{subequations}
  \label{eq:modifiedFundamentalSplineConditions}
  \begin{alignat}{2}
    \parentfundspl[\modified]{p}(i)
    &= \kronecker{i}{1},\quad
    &&i \in \nat,\\
    \tderiv[2]{x}{\parentfundspl[\modified]{p}}(1)
    &= 0,\quad&&\\
    \tderiv[q]{x}{\parentfundspl[\modified]{p}}(0)
    &= 0,\quad
    &&q = 2, 3, \dotsc, \tfrac{p+1}{2},
  \end{alignat}
\end{subequations}
if $p > 1$.
For $p = 1$, we define $\parentfundspl[\modified]{p}
\ceq \bspl[\modified]{2,1}{p}(\tfrac{{\cdot}}{4})$.
Since the modification coefficients $\fundsplcoeff[\modified]{k}{p}$
experience the same decay as the coefficients $\fundsplcoeff{k}{p}$ of the
fundamental spline,
we can also approximate $\fundsplcoeff[\modified]{k}{p}$ by solving
a truncated system of linear equations.
The resulting function $\parentfundspl[\modified]{p}$
is shown in \cref{fig:modifiedFundamentalSpline}.
The corresponding hierarchical
basis $\bspl[\fs,\modified]{l,i}{p}$ is included in
\cref{fig:hierarchicalFundamentalSpline} (dashed lines).

\vspace*{\fill}

The conditions stated in \eqref{eq:modifiedFundamentalSplineConditions}
are motivated by the case $p = 3$ of cubic fundamental splines.
The first relevant cardinal B-spline is the one with
index $k = 1 - \tfrac{p+1}{2}$ ($k = -1$ in the cubic case),
as the B-splines with indices $\le -\tfrac{p+1}{2}$ vanish on $\nonnegreal$.
The modified function $\parentfundspl[\modified]{p}$
should satisfy the fundamental property \eqref{eq:tiftDefinition}
at all positive integer points $k \in \nat$.
In contrast to the standard fundamental spline $\parentfundspl{p}$,
we do not enforce the fundamental property at $k = 0$,
as our aim is to obtain non-zero boundary values.
This leaves us exactly two degrees of freedom in the cubic case,
namely $k = -1$ and $k = 0$.
We use these to let $\parentfundspl[\modified]{p}$ extrapolate
linearly on $\clint{0, 1}$, as we did for uniform hierarchical
B-splines (see \cref{sec:313modification}).
\pagebreak
Therefore, in the cubic case, we set the second derivative
$\tderiv[2]{x}{\parentfundspl[\modified]{p}}$ to zero
at $x = 0$ and at $x = 1$.
This suffices since
$\tderiv[2]{x}{\parentfundspl[\modified]{p}}$ is piecewise
linear for $p = 3$.
For higher degrees $p > 3$,
we use the additional degrees of freedom (in total $\tfrac{p+1}{2}$)
to increase the multiplicity of the root of
$\tderiv[2]{x}{\parentfundspl[\modified]{p}}$ at $x = 0$.
This ensures that $\parentfundspl[\modified]{p}$ is
``as linear as possible'' near $x = 0$.
Note that we cannot maintain
$\tderiv[2]{x}{\parentfundspl[\modified]{p}} \equiv 0$
on $\clint{0, 1}$ for higher degrees $p > 3$,
since this would require $p - 1$ conditions,
and we only have $\tfrac{p+1}{2}$ degrees of freedom left,
after taking the fundamental conditions into account.

\subsection{Fundamental Not-A-Knot Splines}
\label{sec:445fundamentalNotAKnotSplines}

The hierarchical fundamental spline basis suffers from the same problem
as the uniform hierarchical B-spline basis.
As explained in \cref{sec:321approximation}, there is a mismatch
of dimensions of the nodal B-spline space $\nsbspl{l}{p}$ of level $l$
when compared with the spline space $\wholesplspace{l}{p}$ on the grid
$\{\gp{l,i} \mid i = 0, \dotsc, 2^l\}$ of level $l$.
This issue also affects the fundamental spline basis.

\paragraph{Definition of fundamental not-a-knot splines}

\begin{SCfigure}
  \includegraphics{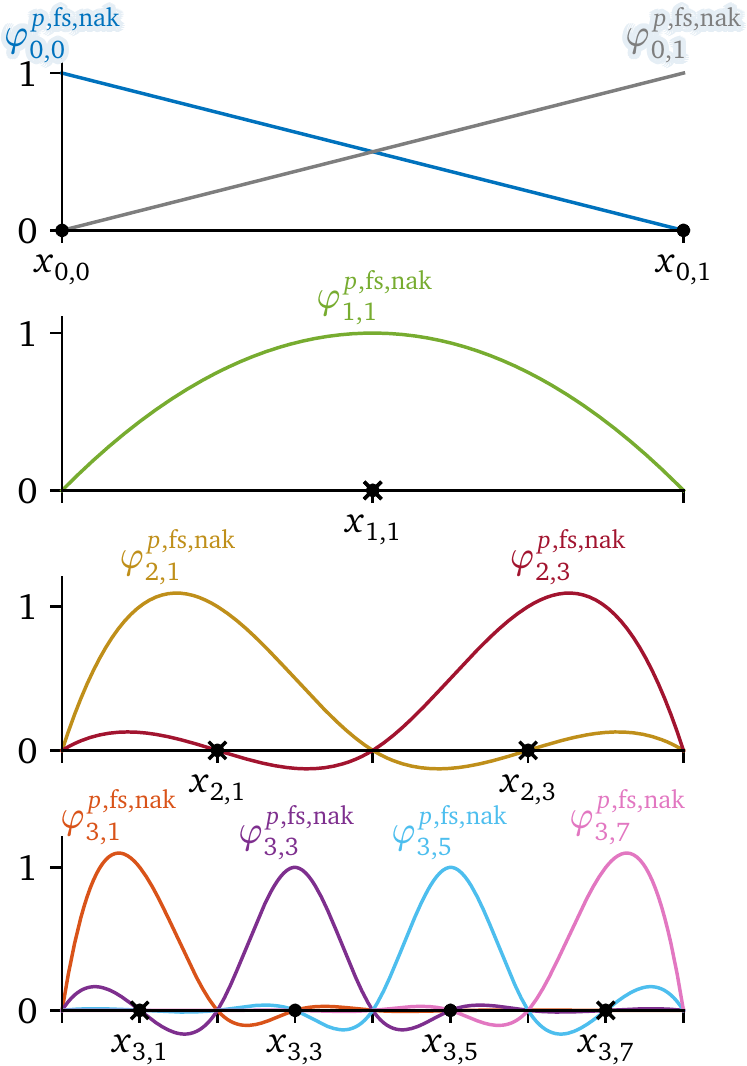}%
  \caption[%
    Hierarchical fundamental not-a-knot splines%
  ]{%
    Hierarchical cubic fundamental not-a-knot splines
    $\bspl[\fs,\nak]{l',i'}{p}$
    ($l' \le l$, $i' \in \hiset{l'}$, $p = 3$),
    grid points $\gp{l',i'}$ \emph{(dots),} and
    removed knots \emph{(crosses)} up to level $l = 3$.%
  }%
  \label{fig:hierarchicalFundamentalNotAKnotSpline}%
\end{SCfigure}

It is possible to combine the idea of fundamental splines
with the not-a-knot approach from \cref{sec:32notAKnot}.
We define hierarchical fundamental not-a-knot splines
$\bspl[\fs,\nak]{l',i'}{p}\colon \clint{0, 1} \to \real$ as
linear combinations of nodal not-a-knot B-splines of the same level,
where the coefficients are chosen such that the
fundamental property \eqref{eq:fundamentalProperty} is satisfied:
\begin{equation}
  \label{eq:fundamentalNotAKnotSplines}
  \bspl[\fs,\nak]{l',i'}{p}
  \ceq \sum_{i''=0}^{2^{l'}}
  \fundsplcoeff[l',i',\fs]{i''}{p} \bspl[\nak]{l',i''}{p}
  \quad\text{such that}\quad
  \falarge{i = 0, \dotsc, 2^{l'}}{
    \bspl[\fs,\nak]{l',i'}{p}(\gp{l',i}) = \kronecker{i}{i'}
  },
\end{equation}
where $l' \in \natz$, $i' = 0, \dotsc, 2^{l'}$, and
$\fundsplcoeff[l',i',\fs]{i''}{p} \in \real$.
This approach is similar to the \hftr in
\cref{sec:442constructingFundamentalBases},
see \cref{eq:hierarchicalFundamentalTransformation}.
We show the hierarchical fundamental not-a-knot spline basis
of cubic degree in \cref{fig:hierarchicalFundamentalNotAKnotSpline}.

The fundamental not-a-knot splines $\bspl[\fs,\nak]{l',i'}{p}$
of level $l' < \ceil{\log_2(p+1)}$ equal the Lagrange polynomials
$\lagrangepoly{l',i'}$ ($i' = 0, \dotsc, 2^{l'}$),
This is because the $i'$-th summand $\bspl[\nak]{l',i'}{p}$
of \eqref{eq:fundamentalNotAKnotSplines} equals $\lagrangepoly{l',i'}$ and
as $\lagrangepoly{l',i'}$ already fulfills the
fundamental interpolation conditions given in
\eqref{eq:fundamentalNotAKnotSplines}
(see \cref{eq:hierarchicalNotAKnotBSpline}),
we obtain $\fundsplcoeff[l',i',\fs]{i''}{p} = \kronecker{i'}{i''}$, i.e.,
$\bspl[\fs,\nak]{l',i'}{p} = \lagrangepoly{l',i'}$.

\paragraph{Implementation}

We make two remarks with respect to the efficient implementation
of hierarchical fundamental not-a-knot splines.
First,
\cref{eq:fundamentalNotAKnotSplines} requires the solution of a system
of linear equations with dimension $2^{l'} + 1$,
which grows exponentially in the level~$l'$.
However, as the coefficients decay roughly in the same order
as the fundamental spline coefficients $\fundsplcoeff{k}{p}$ in
\cref{eq:fundamentalSplineCoefficientsDecay},
we can solve a truncated system of linear equations instead.

Second, the fundamental not-a-knot spline basis
$\bspl[\fs,\nak]{l',i'}{p}$ is not translation-invariant anymore.
This means that theoretically, we have to compute the
$\fundsplcoeff[l',i',\fs]{i''}{p}$ individually for each basis function
$\bspl[\fs,\nak]{l',i'}{p}$.
Nevertheless, when truncating the linear system for a fixed level $l'$,
almost all the inner basis functions $\bspl[\fs,\nak]{l',i'}{p}$
will be identical to hierarchical fundamental splines
$\bspl[\fs]{l',i'}{p}$, if the distance of the region with the removed knots
to the grid point $\gp{l',i'}$ is large enough
(if the removed knots are outside the truncated support of
$\bspl[\fs,\nak]{l',i'}{p}$).
For different levels $l'$, the fundamental not-a-knot splines
$\bspl[\fs,\nak]{l',i'}{p}$ are the same up to scaling
(if the level $l'$ is high enough).

Consequently, an efficient implementation only has to implement
$\bspl[\fs,\nak]{l',i'}{p}$ for some special cases for coarse levels.
The other basis functions can then be derived via an affine
parameter transformation.

\longsection{%
  \texorpdfstring{%
    Hierarchization on Spatially Adaptive Sparse Grids\\
    with the Unidirectional Principle%
  }{%
    Hierarchization on Spatially Adaptive Sparse Grids
    with the Unidirectional Principle%
  }%
}{%
  Hierarchization on Spatially Adaptive Sparse Grids
  with the Unidirectional Principle%
}{%
  Hierarchization with the Unidirectional Principle%
}
\label{sec:45spatAdaptiveUP}

\minitoc{90mm}{9}

\noindent
In this final section of the chapter,
we further decrease the computational complexity
for the application of the linear operator $\linop$
on spatially adaptive sparse grids
from quadratic to linear time with two algorithms based on the \up.

\subsection{%
  Iteratively Applying the Unidirectional Principle with Iterative Refinement%
}
\label{sec:451iterativeRefinement}

The first algorithm can be applied if two requirements are met:
\begin{itemize}
  \item
  The inverse $\linop^{-1}$ is known and can be efficiently applied.
  
  \item
  There is an operator $\linop'$
  that is ``sufficiently close'' to $\linop$ and can be efficiently applied.
\end{itemize}
For hierarchization with B-splines on sparse grids,
we choose $\linop$ to be the hierarchization
operator given in \cref{eq:hierarchizationSLE} and
$\linop'$ to be the \up directly applied on the
sparse grid.
Both of the assumptions are then satisfied,
as $\linop^{-1}$ is known
(interpolation matrix $\intpmat$ of basis function evaluations)
and $\linop^{-1}$ and $\linop'$ can be applied fast.
The \up $\linop'$ generally produces wrong
results for hierarchical B-splines due to missing coupling points.
However, especially for low B-spline degrees,
$\linop'$ does not deviate too much from the true operator $\linop$.
Below, we will specify a sufficient criterion for the ``closeness.''

\paragraph{Iterative refinement}

Under the two assumptions above, we can apply the procedure given in
\cref{alg:iterativeRefinement}.
The algorithm is equivalent to the well-known method of
\term{iterative refinement,} which has been developed to
stabilize the numerical solution of a linear system
influenced by rounding errors \cite{Higham02Accuracy}.
The operator $\linop'$ acts like a preconditioner,
which is why it is required to be close to $\linop$.
Note that the algorithm is similar to the repeated application
of the method of residual interpolation
(see \cref{sec:433residualInterpolation}) on the whole sparse grid.

\begin{algorithm}
  \begin{algorithmic}[1]
    \Function{$\vlinout = \texttt{iterativeRefinement}$}{%
      $\vlinin$, $\vlinout[(0)]$%
    }
      \State{$\*r^{(0)} \gets \vlinin - \linop^{-1}\vlinout[(0)]$}
      \Comment{initial residual}%
      \For{$m = 0, 1, 2, \dotsc$}
        \State{$\vlinout[(m+1)] \gets \vlinout[(m)] + \linop' \*r^{(m)}$}
        \Comment{update solution}%
        \State{$\*r^{(m+1)} \gets \*r^{(m)} - \linop^{-1} \linop' \*r^{(m)}$}
        \Comment{update residual}%
      \EndFor{}
      \vspace{-1mm}
      \State{$\vlinout \gets \text{last computed } \vlinout[(m)]$}
    \EndFunction{}
  \end{algorithmic}
  \caption[%
    Iterative refinement%
  ]{%
    Application of a tensor product operator $\linop$
    on spatially adaptive sparse grids with iterative refinement,
    where $\linop'$ is an approximation of $\linop$.
    Inputs are the vector $\vlinin = (\linin{\*l,\*i})_{(\*l,\*i) \in \liset}$
    of input data (function values $\fcnval{\*l,\*i}$ at the grid points) and
    an initial solution $\vlinout[(0)]$.
    The output is the vector
    $\vlinout = (\linout{\*l,\*i})_{(\*l,\*i) \in \liset}$
    of output data (hierarchical surpluses $\surplus{\*l,\*i}$).%
  }%
  \label{alg:iterativeRefinement}%
\end{algorithm}

The loop in \cref{alg:iterativeRefinement} has to be terminated
after some iterations.
The following simple lemma allows to use a stopping criterion based on the
size of the residual $\*r^{(m)}$ to the true solution,
which we denote with $\vlinout[\ast] \ceq \linop \vlinin$.

\begin{shortlemma}[equivalent convergence for iterative refinement]
  \label{lemma:iterativeRefinementEquivalent}
  In \cref{alg:iterativeRefinement}, we have
  $\vlinout[(m)] \to \vlinout[\ast] \iff \*r^{(m)} \to \*0$ for
  $m \to \infty$.
\end{shortlemma}

\vspace{-0.5em}

\begin{proof}
  It suffices to prove $\linop \*r^{(m)} = \vlinout[\ast] - \vlinout[(m)]$
  for $m \in \nat$ by induction.
  For $m = 0$, we have
  $\linop \*r^{(0)}
  = \linop \vlinin - \linop \linop^{-1} \vlinout[(0)]
  = \vlinout[\ast] - \vlinout[(0)]$.
  For $m \to m+1$, it holds
  $\linop \*r^{(m+1)}
  = \linop \*r^{(m)} - \linop \linop^{-1} \linop' \*r^{(m)}
  = (\vlinout[\ast] - \vlinout[(m)]) - \linop' \*r^{(m)}
  = \vlinout[\ast] - \vlinout[(m+1)]$.
\end{proof}

\vspace{0.5em}

\noindent
Next, we give a sufficient condition for the
convergence of \cref{alg:iterativeRefinement} to the true solution.

\vspace{0.5em}

\begin{proposition}[%
  sufficient condition for the convergence of
  {\hyperref[alg:iterativeRefinement]{Alg.\ \ref*{alg:iterativeRefinement}}}%
]
  \label{prop:iterativeRefinementSufficient}
  If we have $\limsup_{m \to \infty}
  \sqrt[m]{\norm{(\idop - \linop^{-1} \linop')^m}} < 1$
  with an arbitrary operator matrix norm $\norm{\cdot}$ and the
  identity operator $\idop$,
  then $\vlinout[(m)] \to \vlinout[\ast]$ for $m \to \infty$
  in \cref{alg:iterativeRefinement}
  for every initial solution $\vlinout[(0)]$.
\end{proposition}

\vspace{-0.5em}

\begin{proof}
  A short proof by induction shows that
  \begin{equation}
    \label{eq:proofPropIterativeRefinementSufficient}
    \vlinout[(m)]
    = \vlinout[(0)] + \linop'
    \sum_{m'=0}^{m-1} (\idop - \linop^{-1} \linop')^{m'} \*r^{(0)},
  \end{equation}
  where $(\idop - \linop^{-1} \linop')^{m'} \*r^{(0)} = \*r^{(m')}$.
  For $m \to \infty$ and with the assumption on
  $\norm{(\idop - \linop^{-1} \linop')^m}$,
  the sum converges to the Neumann series
  $\sum_{m'=0}^\infty (\idop - \linop^{-1} \linop')^{m'}
  = (\idop - (\idop - \linop^{-1} \linop'))^{-1} = (\linop')^{-1} \linop$
  (see, e.g., \cite{Werner11Funktionalanalysis}).
  In this case, we infer that the limit of $\vlinout[(m)]$ is given by
  \begin{equation}
    \vlinout[(0)] + \linop' (\linop')^{-1} \linop \*r^{(0)}
    = \vlinout[(0)] + \linop \vlinin - \linop \linop^{-1} \vlinout[(0)]
    = \linop \vlinin
    = \vlinout[\ast],
  \end{equation}
  as claimed.
\end{proof}

The sufficient condition given in \cref{prop:iterativeRefinementSufficient}
is quite strong, as it can be shown that $\limsup_{m \to \infty}
\sqrt[m]{\norm{(\idop - \linop^{-1} \linop')^m}} \le 1$ is necessary for
convergence.
Unfortunately, in the case of hierarchization with B-splines,
numerical experiments show that
this condition is only met for low dimensionalities $d$ and low
B-spline degrees $p$.
\Cref{alg:iterativeRefinement} generally diverges
for higher dimensionalities or higher degrees.

\subsection{Duality of the Unidirectional Principle}
\label{sec:452duality}

To motivate the second algorithm that we present in this section,
we study why we cannot directly apply the \up
(as introduced in \cref{alg:unidirectionalPrinciple})
on spatially adaptive sparse grids.
As before, we denote with $\liset$ the level-index set of
the spatially adaptive sparse grid (see \cref{sec:41problem}).

The \up\punctfix{,} as stated in
\cref{alg:unidirectionalPrinciple} for full grids,
subsequently applies one-dimensional operators
$\upopuv{t_j}{\lisetpole}\colon \real^{\setsize{\lisetpole}} \to
\real^{\setsize{\lisetpole}}$
on the poles $\lisetpole$ of the sparse grid at hand,
iterating over a permutation $t_1, \dotsc, t_d$
of the dimensions $1, \dotsc, d$.
We recall the pole equivalence relation $\samepole{t_j}$
from \cref{eq:poleEquivalenceRelation}:
Two points $\*k', \*k'' \in \liset$ are $\samepole{t_j}$-equivalent,
if $\*k'$ is contained in the pole through $\*k''$
with respect to the $t_j$-th dimension, i.e.,
\begin{equation}
  \*k' \samepole{t_j} \*k'' \iff \*k'_{-t_j} = \*k''_{-t_j},\quad
  \*k', \*k'' \in \liset.
\end{equation}

\paragraph{Operators for the unidirectional principle}

The combined application of all one-di\-men\-sional operators
$\upopuv{t_j}{\lisetpole}$
($\lisetpole \in \eqclasses{\liset}{\samepole{t_j}}$)
of the $j$-th iteration of \cref{alg:unidirectionalPrinciple}
is equivalent to a single application of the following operator
$\upop{t_j}\colon \real^{\setsize{\liset}} \to \real^{\setsize{\liset}}$:
\begin{equation}
  \label{eq:upopEntries}
  (\upop{t_j})_{\*k'',\*k'}
  \ceq
  \begin{cases}
    (\upopuv{t_j}{\lisetpole})_{k''_{t_j},k'_{t_j}},&
    \ex{\lisetpole \in \eqclasses{\liset}{\samepole{t_j}}}{
      \*k', \*k'' \in \lisetpole
    },\\
    0,&\*k' \not\samepole{t_j} \*k'',
  \end{cases}
\end{equation}
where $(\upop{t_j})_{\*k'',\*k'}$ denotes the entry of row $\*k''$
and column $\*k'$ of the matrix corresponding to $\upop{t_j}$
(similar for $(\upopuv{t_j}{\lisetpole})_{k''_{t_j},k'_{t_j}}$).
The reason for this equivalence is that
the poles $\lisetpole$ are pairwise disjoint equivalence classes.
Consequently, every point $\*k$ is only acted upon by a single
one-dimensional operator $\upopuv{t_j}{\lisetpole}$,
namely the one with $\lisetpole = \eqclass{\*k}{\samepole{t_j}}$.
This leads to the block-diagonal structure
of $\upop{t_j}$ given in \eqref{eq:upopEntries},
if the rows of the matrix of $\upop{t_j}$
are grouped by poles $\lisetpole$ and the columns are arranged accordingly.

\paragraph{Correctness and duality of the unidirectional principle}

For the remaining considerations, we assume that
the operators $\linop$ and $\upopuv{t_j}{\lisetpole}$ are invertible.
In this case, $\upop{t_j}$ is also invertible and
$\upopinv{t_j}$ is given by the block-diagonal matrix composed of
the inverses of the blocks $\upopuv{t_j}{\lisetpole}$ of $\upop{t_j}$.
This is satisfied by dehierarchization operators $\intpmat$ due to the
linear independence of the hierarchical basis functions.

We are now able to describe the whole \up of
\cref{alg:unidirectionalPrinciple} as the operator
$\upop{t_1,\dotsc,t_d}\colon \real^{\setsize{\liset}} \to
\real^{\setsize{\liset}}$ given by
\begin{equation}
  \label{eq:upopProduct}
  \upop{t_1,\dotsc,t_d}
  \ceq \upop{t_d} \dotsm \upop{t_1}.
\end{equation}
The right-most operator is $\upop{t_1}$, since it is applied first.
We say that the \up is \term{correct} for $\linop$ and
$(t_1, \dotsc, t_d)$, if
\begin{equation}
  \upop{t_1,\dotsc,t_d}
  \overset{?}{=} \linop.
\end{equation}
This relation is not satisfied in general,
especially for B-spline hierarchization with the operator
$\linop = \intpmatinv$.
However, for operators like these, whose inverse
$\linopinv = \intpmat$ can be described and applied much easier,
we can make use of the so-called \term{duality of the \up:}

\begin{lemma}[duality of the unidirectional principle]
  \label{lemma:dualityUnidirectionalPrinciple}
  Let the operators $\linop$ and $\upopuv{t_j}{\lisetpole}$ be invertible
  for all poles $\lisetpole$ in $\liset$.
  Then the \up is correct for $\linop$ and $(t_1, \dotsc, t_d)$
  if and only if the \up is correct for $\linopinv$ and $(t_d, \dotsc, t_1)$.
\end{lemma}

\begin{proof}
  The correctness of the \up for $\linop$ and $(t_1, \dotsc, t_d)$
  is by definition equivalent to
  \begin{equation}
    \upop{t_d} \dotsm \upop{t_1} = \linop.
  \end{equation}
  By inverting both sides, we obtain the definition of the
  correctness of the \up for $\linopinv$ and $(t_d, \dotsc, t_1)$.
\end{proof}

This duality means that in order to establish the correctness of $\linop$
for some arbitrary permutation $(t_1, \dotsc, t_d)$ of $1, \dotsc, d$,
it suffices to establish the \up's correctness for the
inverse operator $\linopinv$ and the reverse permutation $(t_d, \dotsc, t_1)$.
This is especially of interest for our main application,
the hierarchization operator $\linop = \intpmatinv$ for B-splines.

\subsection{Chains and Equivalent Correctness Conditions}
\label{sec:453chains}

We first define the notion of a chain between two grid points
$\*k'$ and $\*k''$.

\begin{definition}[chain]
  \label{def:chain}
  Let $\*k', \*k'' \in \liset$ and
  $(t_1, \dotsc, t_j)$ be a permutation of $j$ of the
  dimensions $1, \dotsc, d$.
  We define the \term{chain} from $\*k'$ to $\*k''$ with respect to
  $(t_1, \dotsc, t_j)$ as the sequence
  $(\chain{0}, \dotsc, \chain{j})$, where
  \begin{equation}
    \chain{j'}_{T_{j'}}
    \ceq \*k''_{T_{j'}},\quad
    \chain{j'}_{-T_{j'}}
    \ceq \*k'_{-T_{j'}},\quad
    T_{j'}
    \ceq (t_1, \dotsc, t_{j'}),\quad
    j' = 0, \dotsc, j,
    \hspace*{-10mm}
  \end{equation}
  if $\chain{j} = \*k''$ and
  $\chain{j'} \in \liset$ for all $j' = 0, \dotsc, j$.
\end{definition}

\begin{figure}
  \subcaptionbox{%
    $d = 2$, $(t_1, t_2) = (2, 1)$%
  }[47mm]{%
    \includegraphics{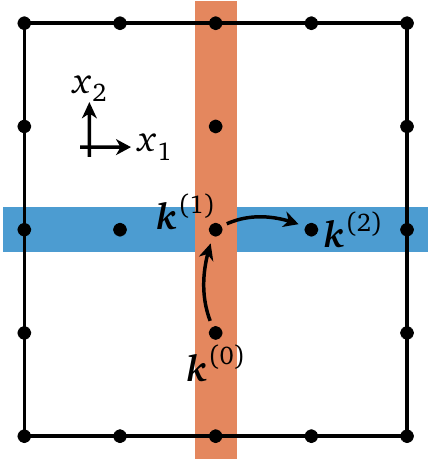}%
  }%
  \hfill%
  \subcaptionbox{%
    $d = 2$, $(t_1, t_2) = (1, 2)$%
    \label{fig:chainDefinition2}%
  }[47mm]{%
    \includegraphics{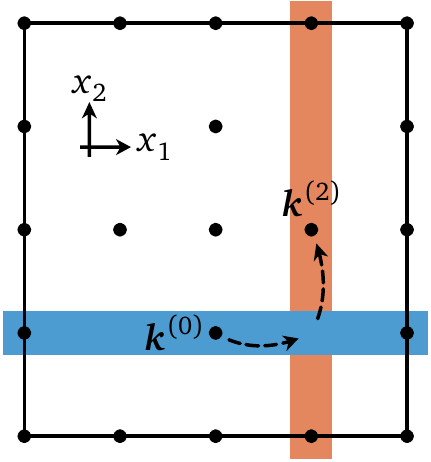}%
  }%
  \hfill%
  \subcaptionbox{%
    $d = 3$, $(t_1, t_2, t_3) = (2, 3, 1)$%
  }[50mm]{%
    \includegraphics{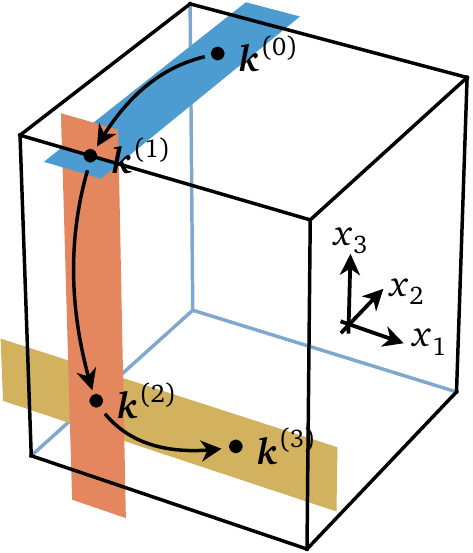}%
  }%
  \caption[%
    Examples for the definition of chains%
  ]{%
    Examples for chains in two and three dimensions.
    \emph{Left:} A chain from $\chain{0}$ to $\chain{2}$
    with respect to $(t_1, t_2) = (2, 1)$ in a two-dimensional sparse grid.
    \emph{Center:}
    With respect to the reverse permutation
    $(t_1, t_2) = (1, 2)$ of the dimensions,
    there is no chain from $\chain{0}$ to $\chain{2}$,
    because the corresponding chain point $\chain{1}$ is missing in the grid.
    \emph{Right:} A chain in three dimensions.%
  }%
  \label{fig:chainDefinition}%
\end{figure}

This definition is equivalent to
$\chain{j'-1} \samepole{t_{j'}} \chain{j'}$ for $j' = 1, \dotsc, j$.
\Cref{fig:chainDefinition} shows examples of chains in
two and three dimensions.
As it is shown in \cref{fig:chainDefinition2},
the order $(t_1, \dotsc, t_j)$ of the dimensions is important
for whether the grid contains the chain from $\*k'$ to $\*k''$.
The grid must contain all intermediate points, otherwise
it is not a chain.

We now show two lemmas.
First, we prove that $(\upop{t_1,\dotsc,t_j})_{\*k'',\*k'} \not= 0$
is sufficient for the existence of a chain from $\*k'$ to $\*k''$:

\begin{restatable}[sufficient condition for chain existence]{%
  lemma%
}{%
  lemmaChainExistenceSufficient%
}
  \label{lemma:chainExistenceSufficient}
  If $(\upop{t_1,\dotsc,t_j})_{\*k'',\*k'} \not= 0$
  for some $j = 0, \dotsc, d$,
  then the grid $\liset$ contains the chain from $\*k'$ to $\*k''$
  with respect to $(t_1, \dotsc, t_j)$.
\end{restatable}

\vspace{-0.5em}

\begin{proof}
  See \cref{sec:a134proofCorrectnessUnidirectionalPrincipleSASG}.
\end{proof}

\vspace{0.5em}

Second, we show that the equality of
$(\upop{t_1,\dotsc,t_j})_{\chain{j},\*k'}$ and the product of
the one-dimensional operators is necessary for the
existence of a chain from $\*k'$ to $\*k''$:

\begin{restatable}[necessary condition for chain existence]{%
  lemma%
}{%
  lemmaChainExistenceNecessary%
}
  \label{lemma:chainExistenceNecessary}
  If the grid $\liset$ contains the chain $(\chain{0}, \dotsc, \chain{j})$
  from $\*k'$ to $\*k''$ with respect to $(t_1, \dotsc, t_j)$
  for some $j = 0, \dotsc, d$, then
  \begin{equation}
    \label{eq:lemmaChainExistenceNecessary}
    (\upop{t_1,\dotsc,t_j})_{\chain{j},\*k'}
    = (\upopuv{t_1}{\eqclass{\chain{1}}{\samepole{t_1}}})_{k''_{t_1},k'_{t_1}}
    \dotsm
    (\upopuv{t_j}{\eqclass{\chain{j}}{\samepole{t_j}}})_{k''_{t_j},k'_{t_j}}.
  \end{equation}
\end{restatable}

\vspace{-0.5em}

\begin{proof}
  See \cref{sec:a134proofCorrectnessUnidirectionalPrincipleSASG}.
\end{proof}

\vspace{0.5em}

These two lemmas can be used to prove the following characterization
of the correctness of the \up\punctfix{.}
Here, we need an additional assumption on the structure of the
operator $\linop$, which we call \term{tensor product structure:}

\begin{restatable}[characterization of the correctness of the UP]{%
  proposition%
}{%
  propCorrectnessUPCharacterization%
}
  \label{prop:correctnessUPCharacterization}
  Let $\linop$ have tensor product structure:
  For all $\*k', \*k'' \in \liset$ with the chain
  $(\chain{0}, \dotsc, \chain{d})$ from $\*k'$ to $\*k''$
  with respect to $(t_1, \dotsc, t_d)$,
  we assume that
  \begin{equation}
    \label{eq:tensorProductOperator}
    (\linop)_{\*k'',\*k'}
    = \prod_{j=1}^d
    (\upopuv{t_j}{\eqclass{\chain{j}}{\samepole{t_j}}})_{k''_{t_j},k'_{t_j}}.
  \end{equation}
  Then the \up is correct for $\linop$ and $(t_1, \dotsc, t_d)$
  if and only if the grid $\liset$ contains the chain from $\*k'$ to $\*k''$
  with respect to $(t_1, \dotsc, t_d)$ for all $\*k', \*k'' \in \liset$
  for which $(\linop)_{\*k'',\*k'} \not= 0$.
\end{restatable}

\vspace{-0.5em}

\begin{proof}
  See \cref{sec:a134proofCorrectnessUnidirectionalPrincipleSASG}.
\end{proof}

\vspace{0.5em}

When applied to the hierarchization operator,
the combination of \cref{prop:correctnessUPCharacterization} with
\thmref{lemma:dualityUnidirectionalPrinciple} can be summarized in
the following corollary:

\begin{corollary}[%
  equivalent statements for correctness of UP for hierarchization%
]
  \label{cor:equivalentCorrectnessUPHierarchization}
  The following statements are equivalent:
  \begin{itemize}
    \item
    The \up is correct for $\intpmatinv$ and $(t_1, \dotsc, t_d)$.
    
    \item
    The \up is correct for $\intpmat$ and $(t_d, \dotsc, t_1)$.
    
    \item
    The grid $\liset$ contains the chain from $\*k'$ to $\*k''$
    with respect to $(t_d, \dotsc, t_1)$ for all $\*k', \*k'' \in \liset$
    for which $\basis{\*k'}(\gp{\*k''}) \not= 0$.
  \end{itemize}
\end{corollary}

\begin{proof}
  The corollary is a direct consequence of
  \cref{lemma:dualityUnidirectionalPrinciple} and
  \cref{prop:correctnessUPCharacterization},
  applied to the dehierarchization operator $\linop = \intpmat$.
  
  The assumption of \cref{lemma:dualityUnidirectionalPrinciple}
  is satisfied:
  The operators $\upopuv{t_j}{\lisetpole}$ are invertible
  for all poles $\lisetpole$ in $\liset$
  due to the uniqueness of univariate interpolants
  (linear independence of the basis functions).
  Similarly, $\linop$ is invertible
  due to the uniqueness of multivariate interpolants.
  In addition, the assumption of \cref{prop:correctnessUPCharacterization}
  is satisfied, since
  \begin{equation}
    (\linop)_{\*k'',\*k'}
    = (\intpmat)_{\*k'',\*k'}
    = \prod_{j=1}^d \basis{k'_{t_j}}(\gp{k''_{t_j}})
    = \prod_{j=1}^d
    (\upopuv{t_j}{\eqclass{\chain{j}}{\samepole{t_j}}})_{k''_{t_j},k'_{t_j}}
  \end{equation}
  due to the tensor product basis functions.
\end{proof}

\paragraph{Inserting chain points}

This means that we can establish the correctness of the \up
for the hierarchization operator $\linop = \intpmatinv$,
if we insert all missing chain points that are specified by
\cref{prop:correctnessUPCharacterization} into the grid.

\begin{figure}
  \subcaptionbox{%
    Original grid ($N = 85$)%
    \label{fig:chainInsertionBSpline1}%
  }[48mm]{%
    \includegraphics{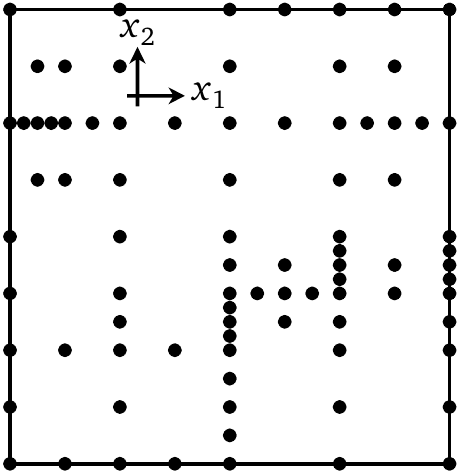}%
  }%
  \hfill%
  \subcaptionbox{%
    Chain points for $p = 1$\\
    \rlap{\hspace*{10.5mm}($N = 121$)}%
    \label{fig:chainInsertionBSpline2}%
  }[48mm]{%
    \includegraphics{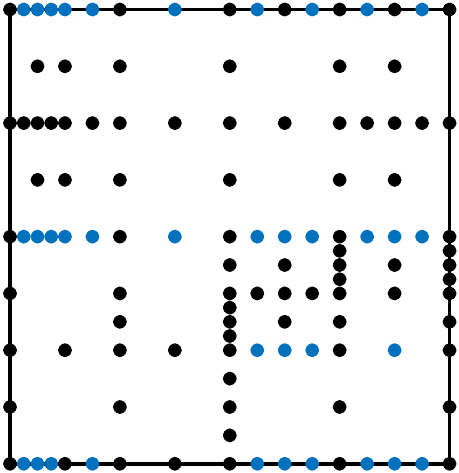}%
  }%
  \hfill%
  \subcaptionbox{%
    Chain points for $p = 3$\\
    \rlap{\hspace*{10.5mm}($N = 289 = 17 \times 17$)}%
    \label{fig:chainInsertionBSpline3}%
  }[48mm]{%
    \includegraphics{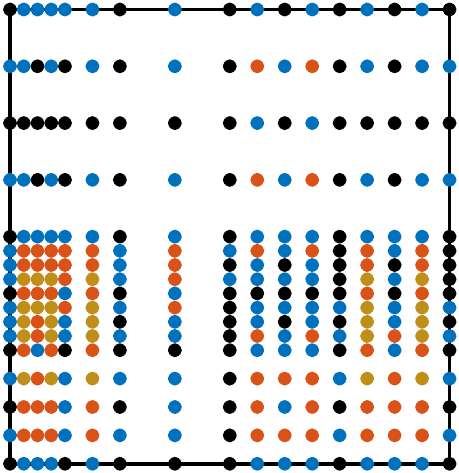}%
  }%
  \caption[%
    Chain points for hierarchical B-splines on a sparse grid%
  ]{%
    Necessary chain points for the correctness of the unidirectional principle
    with respect to $(t_1, t_2) = (1, 2)$
    for hierarchical B-splines $\bspl{l,i}{p}$ on a
    two-dimensional spatially adaptive sparse grid.
    The colors indicate the recursion depth in which the
    chain points have been inserted.
    Black points are contained in the original grid
    (``zero-order points'').
    \textcolor{C0}{Blue points} are part of chains
    between original grid points (``first-order chain points'').
    \textcolor{C1}{Red points} are second-order chain points,
    i.e., they are part of chains
    from $\*k'$ to $\*k''$ where $\*k'$ and $\*k''$ are
    original grid points or first-order chain points
    and at least one of them is a first-order chain point.
    Analogously,
    \textcolor{C2}{brown points} are third-order chain points.
    $N$ is the number of points in the final grid.%
  }%
  \label{fig:chainInsertionBSpline}%
\end{figure}

We take the case $p = 1$ of piecewise linear
standard B-splines $\bspl{\*l,\*i}{1}$ as an example.
We assume that we iteratively generated a spatially adaptive sparse grid
such that all grid points are reachable from the corners of $\clint{\*0, \*1}$
in the sense of \cref{eq:bfsAssumption2}.
If we want to ensure the correctness of the \up for all possible permutations
$(t_1, \dotsc, t_d)$ of the dimensions $(1, \dotsc, d)$,
then the existence of the necessary chains in
\cref{cor:equivalentCorrectnessUPHierarchization} is equivalent to the
requirement that the grid should contain
the hierarchical ancestors of every grid point in every direction:
\begin{subequations}
  \begin{alignat}{4}
    \fafa{(\*l',\*i') \in \liset}{\{t = 1, \dotsc, d \mid l'_t > 1\}}{
      (\*l,\*i) \in \liset
    },\quad
    &&&\*l \ceq \*l' - \stdbasis{t},\quad
    &&i_t \ceq 2 \floor{\tfrac{i'_t}{4}} + 1,\quad
    &&\*i_{-t} = \*i'_{-t},\\
    \fafa{(\*l',\*i') \in \liset}{\{t = 1, \dotsc, d \mid l'_t = 1\}}{
      (\*l,\*i) \in \liset
    },\quad
    &&&\*l \ceq \*l' - \stdbasis{t},\quad
    &&i_t \ceq 0,\quad
    &&\*i_{-t} = \*i'_{-t},
  \end{alignat}
\end{subequations}
where $\stdbasis{t}$ is the $t$-th standard basis vector.
This is a standard assumption on spatially adaptive sparse grids with
piecewise linear basis functions \cite{Pflueger10Spatially}.
However, we only have to satisfy the conditions of
\cref{cor:equivalentCorrectnessUPHierarchization} for a single permutation
$(t_1, \dotsc, t_d)$ of the dimensions
in order to hierarchize with the \up\punctfix{.}
\Cref{fig:chainInsertionBSpline} shows the necessary ancestor chain points
(colored points in \cref{fig:chainInsertionBSpline2})
for an example of a two-dimensional spatially adaptive sparse grid
(\cref{fig:chainInsertionBSpline1}).

Unfortunately, we have to insert these points recursively,
e.g., the inserted points may generate new chains,
for which other missing points have to be inserted and so on
(``higher-order chain points'' in \cref{fig:chainInsertionBSpline}).
Therefore, the number of points to be inserted may be large.
The worst case is that the final grid is a full grid, i.e.,
the Cartesian product of the union of the poles in the different dimensions:
\begin{equation}
  \paren*{\bigcup_{\*k \in \liset} \eqclass{\*k}{\samepole{1}}}
  \times \dotsb \times
  \paren*{\bigcup_{\*k \in \liset} \eqclass{\*k}{\samepole{d}}},
\end{equation}
i.e., we fully lose the advantage of sparse grids,
whose purpose is to ease the curse of dimensionality.
For the standard hierarchical B-spline basis $\bspl{l,i}{p}$,
this worst case often occurs as there are many non-zero entries
in the corresponding interpolation matrices $\intpmat$
(see \cref{sec:41problem} and \cref{fig:chainInsertionBSpline3}).

\subsection{Hierarchical Weakly Fundamental Splines}
\label{sec:454wfs}

\paragraph{Motivation}

In order to reduce the number of chain points to be inserted,
we have to use other spline bases such that
the resulting interpolation matrices $\intpmat$ have more zero entries.
The hierarchical fundamental splines
as introduced in \cref{sec:443fundamentalSplines} are one possibility.
However, they are globally supported, which implies a number
of disadvantages concerning the algorithms and the implementations.
The most significant disadvantage is that although
we can use \bfs for the univariate hierarchization operators,
the time complexity for the univariate hierarchization is still quadratic.
We search for a locally supported spline basis for which
the univariate hierarchization can be done in linear time.

To meet these goals, we have to relax the fundamental property
to a weaker version, which results in the so-called
\term{weakly fundamental property.}
A univariate hierarchical basis
$\wfundbasis{l',i'}\colon \clint{0, 1} \to \real$
is called \term{weakly fundamental,} if
\begin{equation}
  \label{eq:weaklyFundamentalProperty}
  \wfundbasis{l',i'}(\gp{l,i}) = 0,\quad
  l < l',\;\;
  i \in \hiset{l}.
\end{equation}
This is exactly the first condition \eqref{eq:fundamentalProperty1}
of the fundamental property \eqref{eq:fundamentalProperty}.
We drop the requirement that the basis functions
should vanish at the other grid points of the same level.
The relation \eqref{eq:fundamentalPropertyImplicationMV} from the
fundamental case becomes
\begin{equation}
  \label{eq:weaklyFundamentalPropertyImplicationMV}
  \wfundbasis{\*l',\*i'}(\gp{\*l,\*i})
  \not= 0
  \implies
  \*l' \le \*l,
\end{equation}
i.e., every basis function $\wfundbasis{\*l',\*i'}$
can only be non-zero at grid points $\gp{\*l,\*i}$ with
higher or equal level $\*l$.

\paragraph{Definition of hierarchical weakly fundamental splines}

\begin{SCfigure}
  \includegraphics{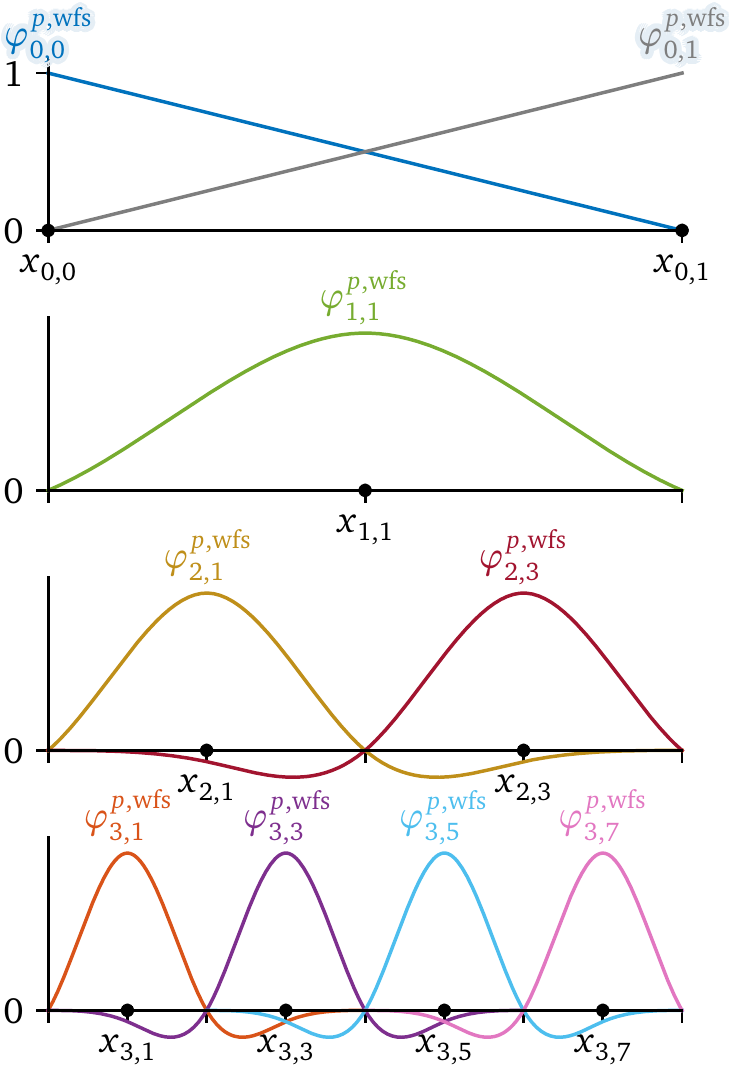}%
  \caption[%
    Hierarchical weakly fundamental splines%
  ]{%
    Hierarchical cubic weakly
    \vspace{-0.1em}%
    fundamental splines
    $\bspl[\wfs]{l',i'}{p}$
    ($l' \le l$, $i' \in \hiset{l'}$, $p = 3$) and
    grid points $\gp{l',i'}$ \emph{(dots)} up to level $l = 3$.%
  }%
  \label{fig:hierarchicalWeaklyFundamentalSpline}%
\end{SCfigure}

We construct the \term{weakly fundamental spline parent function}
$\parentwfundspl{p}\colon \real \to \real$
by forming a linear combination of as few neighboring
uniform B-splines as possible such that $\parentwfundspl{p}$
satisfies the weakly fundamental property
\eqref{eq:weaklyFundamentalProperty}:
\begin{subequations}
  \begin{gather}
    \label{eq:weaklyFundamentalSplineParent}
    \parentwfundspl{p}(x)
    \ceq \largesum[(p-1)/2]{k=-(p-1)/2}
    \wfundsplcoeff{k}{p} \parentbspl{p}(x - k)
    \quad\text{such that}\\
    \wfundsplcoeff{0}{p} = 1,\quad
    \parentwfundspl{p}(k') = 0,\;\;
    k' = -p + 2,\; -p + 4,\; \dotsc,\; p - 2.
  \end{gather}
\end{subequations}
\usenotation{zzzzwfs}
\term{Hierarchical weakly fundamental splines}
$\bspl[\wfs]{l,i}{p}\colon \clint{0, 1} \to \real$
are now defined canonically via an affine parameter transformation:
\begin{equation}
  \bspl[\wfs]{l,i}{p}(x)
  \ceq \parentwfundspl{p}(\tfrac{x}{\ms{l}} - i),\quad
  l \ge 1.
\end{equation}
For $l = 0$, we define $\bspl[\wfs]{l,i}{p}$ to be the
linear Lagrange polynomial of level zero:%
\footnote{%
  This will simplify the description of the
  Hermite hierarchization algorithm in \cref{sec:455hermiteHierarchization}.%
}
\begin{equation}
  \bspl[\wfs]{0,i}{p}
  \ceq \lagrangepoly{0,i},\quad
  i = 0, 1.
\end{equation}
The hierarchical weakly fundamental spline basis is shown in
\cref{fig:hierarchicalWeaklyFundamentalSpline}.
Note that these basis functions are translation-invariant by construction
(starting with level $l \ge 1$).
As the weakly fundamental parent spline $\parentwfundspl{p}$
vanishes at all odd integers and as the
support of $\bspl[\wfs]{l,i}{p}$ is local
($\supp \bspl[\wfs]{l,i}{p}
= \clint{\gp{l,i-p}, \gp{l,i+p}} \cap \clint{0, 1}$),
this implies that the weakly fundamental property
\eqref{eq:weaklyFundamentalProperty} is fulfilled.

\paragraph{Chain points for weakly fundamental splines}

\begin{SCfigure}
  \includegraphics{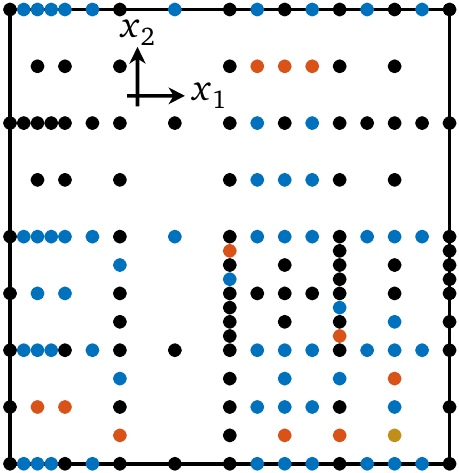}%
  \caption[%
    Chain points for hierarchical weakly fundamental splines on a
    sparse grid%
  ]{%
    Necessary chain points for the correctness of the unidirectional principle
    with respect to $(t_1, t_2) = (1, 2)$
    for hierarchical cubic weakly fundamental splines
    $\bspl[\wfs]{l,i}{p}$ ($p = 3$)
    on the same two-dimensional spatially adaptive sparse grid
    as in \cref{fig:chainInsertionBSpline1}.
    The colors indicate the recursion depth in which the
    chain points have been inserted
    (see caption of \cref{fig:chainInsertionBSpline}).
    The number of points in the final grid is $N = 157$.%
  }%
  \label{fig:chainInsertionWeaklyFundamentalSpline}%
\end{SCfigure}

The first advantage of the
weakly fundamental spline basis $\bspl[\wfs]{l,i}{p}$
over standard uniform B-splines $\bspl{l,i}{p}$ is that
the condition $\basis{\*k'}(\gp{\*k''}) \not= 0$ in
\cref{cor:equivalentCorrectnessUPHierarchization} is
satisfied for much fewer $\*k', \*k''$.
Consequently, fewer chain grid points have to be inserted to
ensure the correctness of the \up for hierarchization.
\Cref{fig:chainInsertionWeaklyFundamentalSpline} shows the inserted points
for the same grid as in \cref{fig:chainInsertionBSpline}.

In the special case of regular sparse grids $\regsgset{n}{d}$,
we do not have to insert any grid points for the correctness of the
\up\punctfix{.}
We can verify this statement with
\thmref{cor:equivalentCorrectnessUPHierarchization}:
Let $(\*l',\*i')$ and $(\*l'',\*i'')$ with
$\normone{\*l'}, \normone{\*l''} \le n$ and
$\*i' \in \hiset{\*l'}$, $\*i'' \in \hiset{\*l''}$,
such that $\bspl[\wfs]{\*l',\*i'}{p}(\gp{\*l'',\*i''}) \not= 0$.
Furthermore, let $(\chain[\*l]{0}, \chain[\*i]{0}), \dotsc,
(\chain[\*l]{d}, \chain[\*i]{d})$ be the chain
from $\*k'$ to $\*k''$ with respect to $t_1, \dotsc, t_d$.
Note that $\chain[\*l]{j} \le \vecmax\{\*l', \*l''\}$ due to the
definition of chain points (\cref{def:chain}).
Therefore, we have for $j = 0, \dotsc, d$
by \eqref{eq:weaklyFundamentalPropertyImplicationMV}:
\begin{equation}
\*l' \le \*l''
\implies
\chain[\*l]{j} \le \vecmax\{\*l', \*l''\} \le \*l''
\implies
\normone{\chain[\*l]{j}} \le \normone{\*l''} \le n.
\end{equation}
Hence, $\regsgset{n}{d}$ contains the grid points corresponding to
$(\chain[\*l]{j}, \chain[\*i]{j})$ for all $j = 0, \dotsc, d$.
Consequently, the conditions of
\cref{cor:equivalentCorrectnessUPHierarchization} are satisfied without
inserting any additional chain points.
This statement is even valid for arbitrary
dimensionally adaptive sparse grids.

\subsection{Hermite Hierarchization}
\label{sec:455hermiteHierarchization}

\paragraph{Hermite interpolation}

The second advantage of the weakly fundamental spline basis
is that due to the reduced coupling,
the univariate hierarchization operators can be applied easier
than for standard uniform B-splines.
This results in the formulation of the so-called
\term{Hermite hierarchization} algorithm.
We first recall higher-order Hermite interpolation:

\begin{lemma}[higher-order Hermite interpolation]
  \label{lemma:hermiteInterpolation}
  Let $p \in \nat$ be odd and $a, b \in \real$ with $a < b$.
  Furthermore, let
  $\deriv[q]{x}{\objfun}(a) \in \real$ and
  $\deriv[q]{x}{\objfun}(b) \in \real$ be given data
  for $q = 0, \dotsc, \frac{p-1}{2}$.
  Then there is a unique polynomial $\spl \in \polyspace{p}$ such that
  \begin{equation}
    \deriv[q]{x}{\objfun}(a)
    = \deriv[q]{x}{\spl}(a),\quad
    \deriv[q]{x}{\objfun}(b)
    = \deriv[q]{x}{\spl}(b),\quad
    q = 0, \dotsc, \frac{p-1}{2}.
    \hspace*{-10mm}
  \end{equation}
\end{lemma}

\begin{proof}
  See \cite{Freund07Stoer}.
\end{proof}

\paragraph{Hermite hierarchization algorithm}

The interpolating polynomial $\spl$ and its derivatives can be
efficiently evaluated using Hermite basis functions
(generalized Lagrange polynomials \cite{Freund07Stoer}).
With Hermite interpolation, we formulate
\cref{alg:hermiteHierarchization}
for the hierarchization with hierarchical weakly fundamental splines.
While we formulate \cref{alg:hermiteHierarchization}
only for regular univariate grids and weakly fundamental splines,
a slightly reformulated version of the algorithm
also correctly operates on spatially adaptive univariate grids
(with the assumption that the grids contain the parents of their grid points)
and other weakly fundamental bases that are
piecewise polynomials of degree $\le p$.

\begin{algorithm}
  \begin{algorithmic}[1]
    \Function{$\vlinout = \texttt{hermiteHierarchization1D}$}{%
      $\vlinin$, $n$%
    }
      \For{$i = 0, 1$}
      \Comment{set values for level $0$}%
        \State{%
          $\linout{0,i} \gets \fcnval{0,i}$%
        }
        \label{line:algHermiteHierarchization1}
        \State{%
          $\deriv[q]{x}{\fgintp{0}}(\gp{0,i})
          \gets \kronecker{q}{0} \cdot \fcnval{0,i} +
          \kronecker{q}{1} \cdot (\fcnval{0,1} - \fcnval{0,0})$
          for all $q = 0, \dotsc, \frac{p-1}{2}$%
        }
        \label{line:algHermiteHierarchization3}
      \EndFor{}
      \For{$l = 1, \dotsc, n$}
        \For{$i \in \hiset{l}$}
          \State{%
            $\fgintp{l-1}(\gp{l,i}) \gets \text{Hermite interpolation of}$
            $\deriv[q]{x}{\fgintp{l-1}}(\gp{l,i\pm1})$
            ($q = 0, \dotsc, \frac{p-1}{2}$)%
          }
          \label{line:algHermiteHierarchization4}
          \State{%
            $r^{(l)}(\gp{l,i})
            \gets \fcnval{l,i} - \fgintp{l-1}(\gp{l,i})$%
          }
          \label{line:algHermiteHierarchization5}
          \Comment{residual to be interpolated}%
        \EndFor{}
        \State{%
          Let $r^{(l)}_l$ be of the form
          $\sum_{i' \in \hiset{l}} \linout{l,i'} \bspl[\wfs]{l,i'}{p}$%
        }
        \Comment{contribution of level $l$}%
        \label{line:algHermiteHierarchization6}
        \State{%
          Choose $(\linout{l,i'})_{i' \in \hiset{l}}$ such that
          $r^{(l)}_l(\gp{l,i}) = r^{(l)}(\gp{l,i})$ for all $i \in \hiset{l}$%
        }
        \label{line:algHermiteHierarchization7}
        \For{$i = 0, \dotsc, 2^l$}
        \Comment{for all points (current level and ancestors)}%
          \For{$q = 0, \dotsc, \frac{p-1}{2}$}
            \State{%
              $\deriv[q]{x}{\fgintp{l}}(\gp{l,i})
              \gets \deriv[q]{x}{\fgintp{l-1}}(\gp{l,i}) +
              \deriv[q]{x}{r^{(l)}_l}(\gp{l,i})$%
            }
            \Comment{update values}%
            \label{line:algHermiteHierarchization8}
          \EndFor{}
        \EndFor{}
      \EndFor{}
    \EndFunction{}
  \end{algorithmic}
  \caption[%
    Hermite hierarchization%
  ]{%
    Hermite hierarchization on one-dimensional regular grids.
    Inputs are
    the vector $\vlinin = (\linin{l,i})_{(l,i) \in \liset}$
    of input data (function values $\fcnval{l,i}$ at the grid points) and
    the level $n$ of the regular grid,
    where $\liset = \{(l, i) \mid l = 0, \dotsc, n,\; i \in \hiset{l}\}$.
    The output is the vector
    $\vlinout = (\linout{l,i})_{(l,i) \in \liset}$
    of output data (hierarchical surpluses $\surplus{l,i}$).%
  }%
  \label{alg:hermiteHierarchization}%
\end{algorithm}

\vspace*{\fill}

\begin{SCfigure}
  \includegraphics{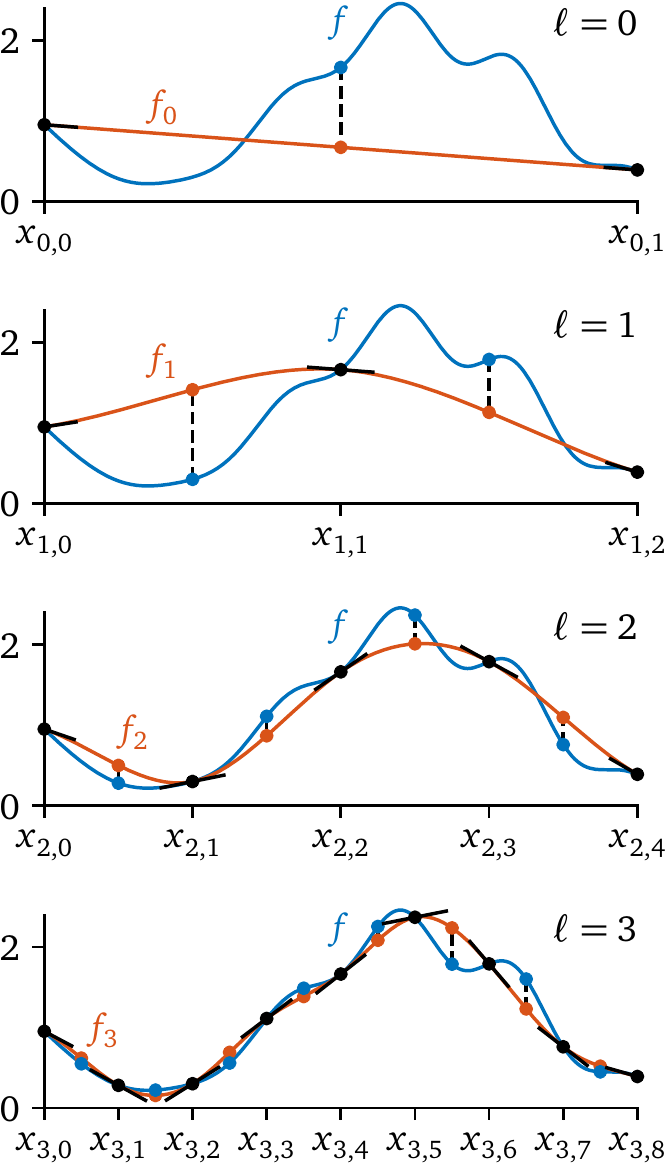}%
  \caption[%
    Hermite hierarchization%
  ]{%
    Hermite hierarchization on a regular grid in one dimension
    with cubic weakly fundamental splines $\bspl[\wfs]{l,i}{p}$ ($p = 3$).
    The interpolants $\fgintp{l}$ \emph{\textcolor{C1}{(red)}}
    of the objective function \emph{\textcolor{C0}{(blue)}}
    are computed level by level.
    For each level $l$,
    \vspace{-0.1em}%
    the values $\fgintp{l}(\gp{l,i})$ and the derivatives
    $\deriv{x}{\fgintp{l}}(\gp{l,i})$ of the
    current interpolant $\fgintp{l}$ at the
    grid points $\gp{l,i}$ ($i = 0, \dotsc, 2^l$) are saved
    \emph{(black dots and bars).}
    The values and derivatives are used for the Hermite interpolation
    of the residual $\objfun - \fgintp{l}$.
    The interpolated residual is then added to the current interpolant
    such that the sum vanishes in the grid points of the next level $l + 1$
    \emph{%
      (black dashed lines between \textcolor{C1}{red} and
      \textcolor{C0}{blue} dots).%
    }
    Due to the weakly fundamental property, the previously
    interpolated values of $\objfun$ remain unchanged.%
  }%
  \label{fig:hermiteHierarchization}%
\end{SCfigure}

The idea of \cref{alg:hermiteHierarchization},
which is also illustrated in \cref{fig:hermiteHierarchization},
is to hierarchize the function value data level by level,
which is only possible because of the weakly fundamental property
\eqref{eq:weaklyFundamentalProperty}.
For each level $l$, we calculate surpluses
$\surplus{l,i} = \linout{l,i}$, while keeping track of
the values and derivatives
$\deriv[q]{x}{\fgintp{l}}(\gp{l,i})$ of the
``current'' interpolant $\fgintp{l}$ (up to level $l$).
Hermite interpolation is used to determine the ``delta''
to the interpolant of the next level.
Note that in \cref{line:algHermiteHierarchization8},
we have to evaluate the derivatives of
$\deriv[q]{x}{\fgintp{l-1}}(\gp{l,i})$ of the Hermite interpolant
determined in \cref{line:algHermiteHierarchization4}.
This is not an issue since in an implementation
one would typically simultaneously evaluate the
Hermite interpolant and its derivatives.

\vspace*{\fill}

For hierarchical weakly fundamental splines,
the complexity of the $l$-th iteration of \cref{alg:hermiteHierarchization}
is linear in the number of grid points of level $l$, i.e., $\landauO{2^l}$.
The reason for this is the bandedness (with bandwidth $\landauO{p}$) of the
system of linear equations corresponding to the interpolation problem of
\cref{line:algHermiteHierarchization6,line:algHermiteHierarchization7},
which means that the interpolation problem can be solved in
linear time and memory.
In total, the complexity of \cref{alg:hermiteHierarchization} is
given by $\landauO{\sum_{l=0}^n 2^l} = \landauO{2^n}$, i.e.,
the time and memory required by \cref{alg:hermiteHierarchization}
is only linear in the number of grid points.

\pagebreak

\paragraph{Correctness}

We prove the correctness of Hermite hierarchization
with the following invariant.

\begin{restatable}[invariant of Hermite hierarchization]{%
  proposition%
}{%
  propInvariantHermiteHierarchization%
}
  \label{prop:invariantHermiteHierarchization}
  In \cref{alg:hermiteHierarchization}, it holds
  for $l = 0, \dotsc, n$ and $i = 0, \dotsc, 2^l$
  \begin{equation}
    \label{eq:propInvariantHermiteHierarchization}
    \deriv[q]{x}{\fgintp{l}}(\gp{l,i})
    = \sum_{l'=0}^l \sum_{i' \in \hiset{l'}}
    \linout{l',i'} \deriv[q]{x}{\bspl[\wfs]{l',i'}{p}}(\gp{l,i}),\quad
    q = 0, \dotsc, \frac{p-1}{2}.
    \hspace*{-6mm}
  \end{equation}
\end{restatable}

\begin{proof}
  See \cref{sec:a135proofHermiteHierarchization}.
\end{proof}

\begin{restatable}[correctness of Hermite hierarchization]{%
  shortcorollary%
}{%
  corAlgHermiteHierarchizationCorrectness%
}
  \label{cor:algHermiteHierarchizationCorrectness}
  \Cref{alg:hermiteHierarchization} is correct.
\end{restatable}

\begin{proof}
  See \cref{sec:a135proofHermiteHierarchization}.
\end{proof}

\subsection{Hierarchical Weakly Fundamental Not-A-Knot Splines}
\label{sec:456wfsNotAKnot}

Finally, as for fundamental splines,
it is possible to combine the weakly fundamental basis
with the not-a-knot idea from \cref{sec:32notAKnot} to construct
hierarchical weakly fundamental not-a-knot spline functions
$\bspl[\wfs,\nak]{l',i'}{p}$.
The approach is similar to the fundamental not-a-knot splines
in \cref{sec:445fundamentalNotAKnotSplines}
(see \cref{eq:fundamentalNotAKnotSplines}):
Instead of combining uniform B-splines as in
\eqref{eq:weaklyFundamentalSplineParent},
we combine not-a-knot B-splines such that the
weakly fundamental property is satisfied.

\begin{SCfigure}
  \includegraphics{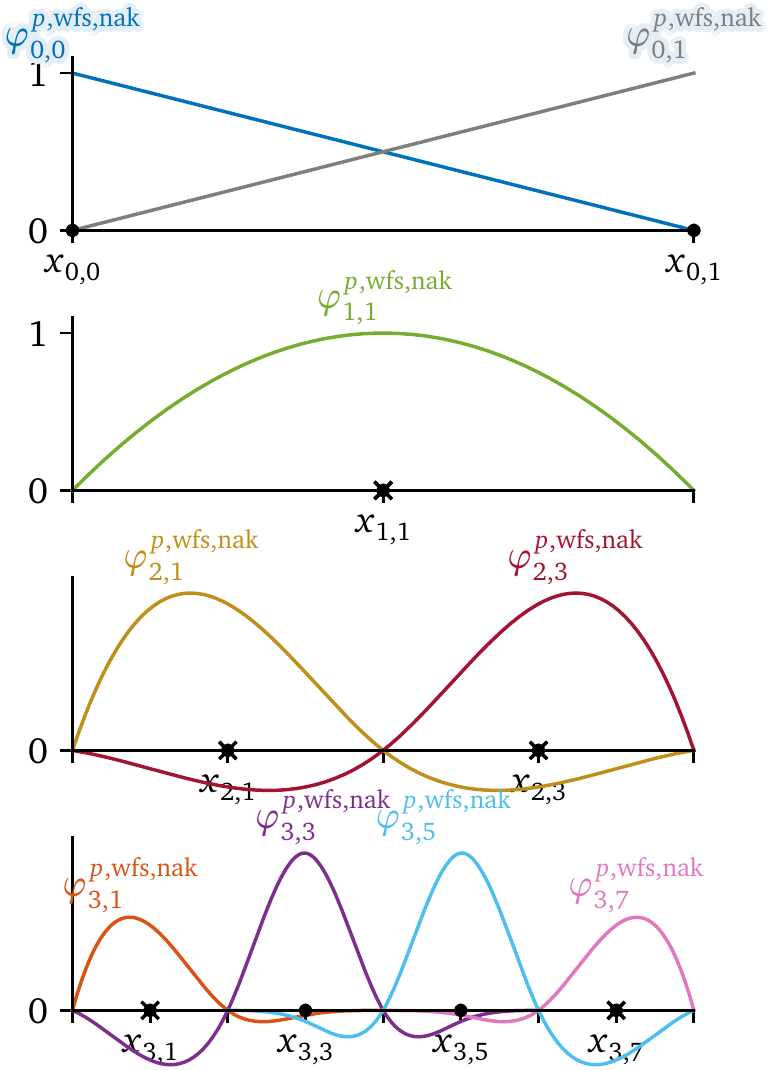}%
  \caption[%
    Hierarchical weakly fundamental not-a-knot splines%
  ]{%
    Hierarchical cubic weakly fundamental not-a-knot splines
    $\bspl[\wfs,\nak]{l',i'}{p}$
    ($l' \le l$, $i' \in \hiset{l'}$, $p = 3$),
    grid points $\gp{l',i'}$ \emph{(dots),} and
    removed knots \emph{(crosses)} up to level $l = 3$.%
  }%
  \label{fig:hierarchicalWeaklyFundamentalNotAKnotSpline}%
\end{SCfigure}

However, the exact construction is somewhat complicated,
as one has to carefully consider which conditions to enforce
with which basis functions.
There are some special cases, if the index of the basis function
$\bspl[\wfs,\nak]{l',i'}{p}$ is near the boundary
(near $i' = 0$ or near $i' = 2^{l'}$).
Nevertheless, there are only finitely many special cases;
for higher levels $l'$, one can just scale the basis functions
of coarser levels.
In the scope of this thesis,
it suffices to show the resulting basis functions for
the cubic case ($p = 3$) in
\cref{fig:hierarchicalWeaklyFundamentalNotAKnotSpline},
instead of rigorously stating the technical formulas.

\cleardoublepage

  \setdictum{%
  Premature optimization is the root of all evil.%
}{%
  Donald E.\ Knuth \cite{Knuth74Structured}%
}

\longchapter{%
  Gradient-Based Optimization with B-Splines on Sparse Grids%
}{%
  Gradient-Based Optimization with B-Splines on Sparse Grids%
}{%
  Gradient-Based Optimization%
}
\label{chap:50optimization}

\initial{0em}{I}{n this chapter,}
we apply the hierarchical B-spline bases derived in
\cref{chap:30BSplines,chap:40algorithms} to optimization,
which is a major task in simulation technology,
for instance in inverse problems (see \cref{chap:10introduction}).
We pursue a three-step surrogate-based optimization approach:
First, we sample the objective function at specific sparse grid points
to retrieve objective function values.
Second, by interpolating these values with hierarchical bases,
we obtain a surrogate for the objective function.
Third and finally, we discard the original objective function and apply
already existing optimization methods to the surrogate.

\vspace*{\fill}

One of the key advantages of hierarchical B-splines
over common hierarchical bases for
sparse grids is their continuous differentiability.
The derivatives of B-spline surrogates on sparse grids are not only continuous,
but also explicitly known, and they can be evaluated fast.
This gives the opportunity to employ gradient-based optimization methods,
which usually converge significantly faster than gradient-free alternatives.

\vspace*{\fill}

The outline of this chapter is as follows:
We start in \cref{sec:51algorithms}
with a compact overview of textbook optimization algorithms,
\pagebreak%
which comprises gradient-free and gradient-based optimization algorithms
for unconstrained problems as well as algorithms for constrained problems.
In \cref{sec:52method}, we present the main method that
conflates the various optimization algorithms and
the generation of spatially adaptive sparse grids with the
Novak--Ritter criterion to create a single method for the
optimization of sparse grid surrogates.
\Cref{sec:53testProblems} continues with a small array of test problems
for unconstrained and constrained optimization.
In \cref{sec:54results}, we apply the presented method
to the test problems, studying the influence of the different
hierarchical B-splines on optimality gaps and conducting
numerical experiments.
Finally, in \cref{sec:55fuzzy}, we examine the fuzzy extension principle
as an example application of the optimization of
hierarchical B-spline surrogates on sparse grids.

Parts of this chapter have already been published in previous work
\cite{Valentin14Hierarchische}, especially
the overview of optimization algorithms (\cref{sec:51algorithms})
and the methodology of the optimization of sparse grid surrogates
(\cref{sec:52method}).
However, as the previous work included other basis functions as well,
this chapter represents the first comprehensive study
that focuses on the application of hierarchical B-splines to optimization.

\section{Overview of Optimization Algorithms}
\label{sec:51algorithms}

\paragraph{Problem setting}

\minitoc{80mm}{7}

\usenotation{zzzzopt}
Generally, \term{unconstrained optimization problems} have the form
\begin{equation}
  \label{eq:unconstrainedOptimization}
  \xopt = \vecargmin \objfun(\*x),\quad
  \*x \in \real^d,
\end{equation}
where $\objfun\colon \real^d \to \real$ is the \term{objective function.}
\term{Constrained optimization problems} are given by
\begin{equation}
  \label{eq:constrainedOptimization}
  \xopt = \vecargmin \objfun(\*x),\quad
  \*x \in \real^d\;\;\text{s.t.}\;\;
  \ineqconfun(\*x) \le \*0,
\end{equation}
where
$\ineqconfun\colon \real^d \to \real^{m_{\ineqconfun}}$
($m_{\ineqconfun} \in \nat$)
is the \term{(inequality) constraint function.}
This formulation also contains optimization problems
with equality constraints $\eqconfun(\*x) = \*0$
by setting $\ineqconfun(\*x) \ceq (\eqconfun(\*x), -\eqconfun(\*x))$.
Equality constraints can also be solved by incorporating them
into the unconstrained solver (e.g., see \cite{Boyd04Convex}
for an equality-constrained Newton method).

As sparse grid surrogates $\objfun = \sgintp$ are only defined on the
unit hyper-cube,
the choice of $\*x$ has to be restricted to $\clint{\*0, \*1}$.
In the case of \eqref{eq:unconstrainedOptimization},
this results in a \term{box-constrained optimization problem.}
A simple method for applying unconstrained optimization algorithms
to box-constrained problems is extending $\sgintp$ to $\real^d$ by
$\sgintp(\*x) \ceq +\infty$ for all $\*x \in \real^d \setminus \clint{\*0, \*1}$.
However, more sophisticated
approaches are also available \cite{More87Optimization}.

\paragraph{Black-box optimization methods}

Problems of the form \eqref{eq:unconstrainedOptimization} or
\eqref{eq:constrainedOptimization} are \term{black-box optimization problems,}
where we cannot gain any insight into the structure or algebraic
properties of $\objfun$.
Black-box optimization methods perform a series of evaluations
$\objfun(\*x_k)$,
choosing the next evaluation point $\*x_{k+1}$
based on the previous function values $\objfun(\*x_0), \dotsc, \objfun(\*x_k)$.
Gradient-based methods differ from gradient-free approaches
in such a way that they also take values
of the gradient $\gradient{\*x}{\objfun}(\*x_k)$,
of the Hessian $\hessian{\*x}{\objfun}(\*x_k)$, or
of even higher-order derivatives into account.

A vast range of optimization methods exists in literature.
Some methods are better suited for specific optimization problems
than others.
However, according to the ``no-free-lunch theorem''
and under some assumptions \cite{Wolpert97No},
all methods perform equally well (or equally badly) in the mean of all possible
optimization problems.

\paragraph{Local and global optima}

Most optimization methods depend on an initial point $\*x_0$ and
only find local optima,
where \eqref{eq:unconstrainedOptimization} or
\eqref{eq:constrainedOptimization} only holds for $\*x$
in a neighborhood of $\xopt$.
One can globalize local methods to increase the probability
of finding a global optimum with a Monte Carlo multi-start approach:
The local method is repeated with different pseudo-random initial points
and the best local optimum is chosen as the result.

In the following,
we give a brief survey of a small selection of optimization methods
(see
\cref{tbl:optimizationMethod},
\cref{fig:optimizationMethodGradientFree}, and
\cref{fig:optimizationMethodGradientBased}),
highlighting the key ingredients for each method.

\begin{table}
  \setnumberoftableheaderrows{1}%
  \begin{tabular}{%
    >{\kern\tabcolsep}=l+l<{\kern5mm}*{5}{+c}<{\kern\tabcolsep}%
  }
    \toprulec
    \headerrow
    Method&                 Type&                    C&    D&  S&    I\\
    \midrulec
    Nelder--Mead&           Simplex heuristic&       \no&  0&  \no&  \yes\\
    Differential evolution& Evolutionary&            \no&  0&  \yes& \yes\\
    CMA-ES&                 Evolutionary&            \no&  0&  \yes& \yes\\
    Simulated annealing&    Temperature heuristic&   \no&  0&  \yes& \no\\
    PSO&                    Swarm heuristic&         \no&  0&  \yes& \no\\
    GP-LCB&                 Bayesian&                \no&  0&  \yes& \no\\
    \midrulec
    Gradient descent&       Descent&                 \no&  1&  \no&  \yes\\
    NLCG&                   Descent&                 \no&  1&  \no&  \yes\\
    Newton&                 Newton&                  \no&  2&  \no&  \yes\\
    BFGS&                   Quasi-Newton&            \no&  1&  \no&  \yes\\
    Rprop&                  Heuristic&               \no&  1&  \no&  \yes\\
    Levenberg--Marquardt&   Least sq., trust-region& \no&  1&  \no&  \yes\\
    \midrulec
    Log-barrier&            Interior-point&          \yes& 0+& --&   \yes\\
    Squared penalty&        Penalty&                 \yes& 0+& --&   \yes\\
    Augmented Lagrangian&   Penalty&                 \yes& 0+& --&   \yes\\
    SQP&                    Quadratic subproblems&   \yes& 2&  --&   \no\\
    \bottomrulec
  \end{tabular}
  \caption[Selection of optimization methods]{%
    Selection of optimization methods.
    The columns show
    if constrained problems are supported (C),
    the order of required derivatives (D),
    if the algorithm is stochastic (S), and
    if the algorithm has been implemented in \sgpp (I).%
  }%
  \label{tbl:optimizationMethod}%
\end{table}

\subsection{Gradient-Free Unconstrained Optimization Methods}
\label{sec:511gradientFree}

\newcommand*{\optImage}[1]{%
  \raisebox{-0.5\height}{\includegraphics{optimizationMethod_#1}}%
}

\begin{figure}
  \optImage{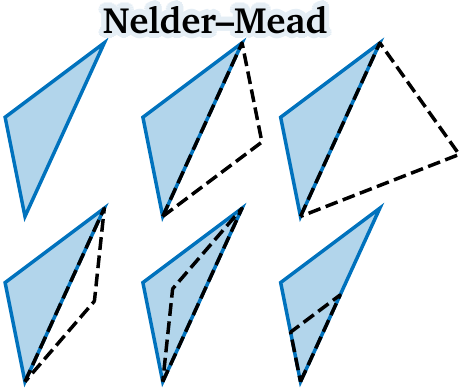}\quad\optImage{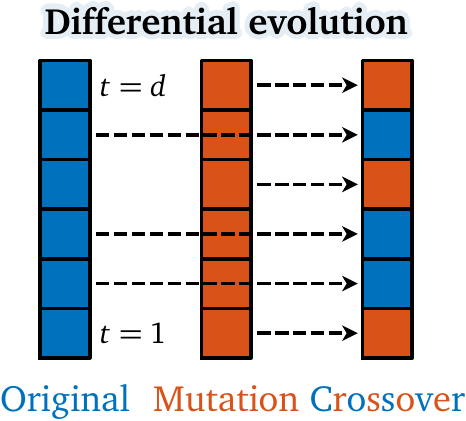}\quad\optImage{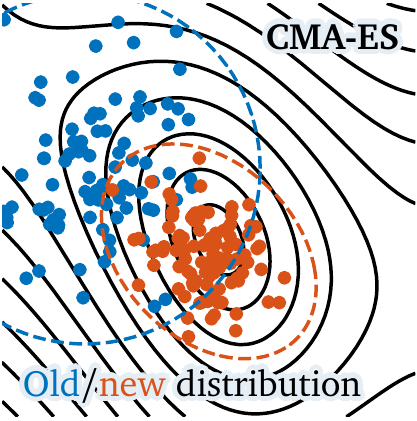}%
  \\[3mm]%
  \optImage{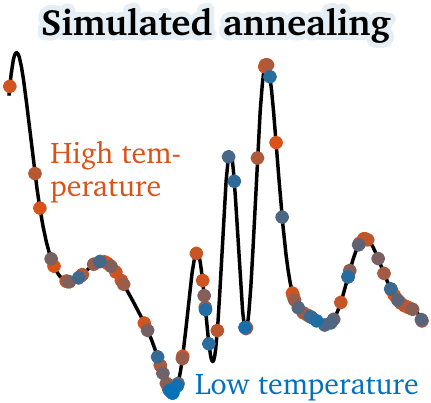}\qquad\optImage{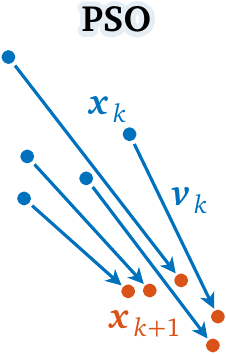}\qquad\optImage{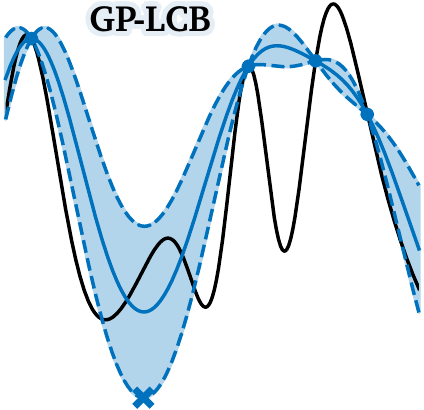}%
  \caption[Ideas of various gradient-free optimization methods]{%
    Sketch of the ideas of gradient-free optimization methods.%
  }%
  \label{fig:optimizationMethodGradientFree}%
\end{figure}

\paragraph{Nelder--Mead}

The \term{Nelder--Mead method}
\multicite{Nelder65Simplex,Gao12Implementing,Valentin14Hierarchische}
maintains a list of $d + 1$ vertices of a $d$-dimensional simplex,
sorted by ascending function value.
In each iteration,
the method performs one of the operations
\term{reflection,}
\term{expansion,}
\term{outer contraction,}
\term{inner contraction,} and
\term{shrinking}
on the vertices.
Typically, convergence can be detected by
checking the size of the simplex,
as the simplex tends to contract around local minima.
However, there are counterexamples where the method converges to
a non-critical point for an only bivariate objective function
that is strictly convex and twice continuously differentiable
\cite{McKinnon98Convergence}.

\paragraph{Differential evolution}

The method of \term{differential evolution}
\multicite{%
  Storn97Differential,%
  Zielinski09Optimizing,%
  Valentin14Hierarchische%
}
is an evolutionary meta-heuristic algorithm.
Being similar to genetic algorithms,
the method maintains a \term{population} of $m$ points
that is iteratively updated according to pseudo-random \term{mutations,}
which are weighted sums of the points of the previous generation.
The mutated vector is \term{crossed over} with the original vector
entry by entry.
The resulting \term{offspring} are only accepted if they lead to
an improvement in terms of objective function value.

\paragraph{CMA-ES}

\term{CMA-ES (covariance matrix adaption, evolution strategy)}
\cite{Hansen03Reducing}
is an evolutionary algorithm that addresses the issue
that simple evolution strategies do not prefer a search direction
due to the lack of gradients \cite{Toussaint15Introduction}.
The name of the algorithm stems from the fact
that it keeps track of the \term{covariance matrix} of the
Gaussian search distribution.
After $m$ points have been sampled from the current distribution,
the mean of the distribution for the next iteration
is calculated as the weighted mean of the $k$ best samples and
the covariance matrix is adapted accordingly.
An advantage of the method is that if the population is large enough,
local minima are smoothed out \cite{Toussaint15Introduction}.

\paragraph{Simulated annealing}

\term{Simulated annealing}
\multicite{Laarhoven87Simulated,Press07Numerical,Kiranyaz14Multidimensional}
imitates the cooling of a solid by randomly drawing samples
from a proposal distribution and calculating an \term{acceptance probability}
that depends on the function value improvement as well as
on a \term{temperature} $T$.
This temperature is slowly decreased in the course of the algorithm.
Simulated annealing is closely connected to the
Metropolis--Hastings algorithm for drawing random samples of arbitrary
probability distributions.
If run long enough,
simulated annealing is guaranteed to find the global
optimum \cite{Toussaint15Introduction}.

\paragraph{Particle swarm optimization (PSO)}

The method of \term{particle swarm optimization (PSO)}
\multicite{Kennedy95Particle,Zielinski09Optimizing,Kiranyaz14Multidimensional}
can be seen as another evolutionary algorithm
that stems from swarm intelligence.
For each \term{particle} of the population,
not only the \term{position} $\*x_k$ is stored,
but also the current \term{velocity} $\*v_k$,
the best known position in a neighborhood of $\*x_k$
(which may be the whole swarm), and
the best known position of the $k$-th particle.
The next velocity $\*v_{k+1}$ is computed as
a pseudo-randomly weighted sum of $\*v_k$,
the vector from $\*x_k$ to the best neighborhood position, and
the vector from $\*x_k$ to the best own position.

\paragraph{GP-LCB}

\term{GP-LCB (Gaussian process, lower confidence bound)}
\multicite{Srinivas10Gaussian,Toussaint15Introduction} is an example
for a \term{Bayesian optimization} strategy.
The objective function is treated as a stochastic process.
A \term{prior distribution} is updated according to the previous function
evaluations to calculate the \term{posterior distribution.}
The posterior distribution is used to form the \term{acquisition function,}
which in turn determines the point at which the objective
function is evaluated next.
The GP-LCB method is obtained by choosing
\term{Gaussian processes} for the family of stochastic processes and
\term{lower confidence bounds} (which are the difference of the mean
and a multiple of the standard deviation) for the acquisition function.

\subsection{Gradient-Based Unconstrained Optimization Methods}
\label{sec:512gradientBasedUnconstrained}

\begin{figure}
  \optImage{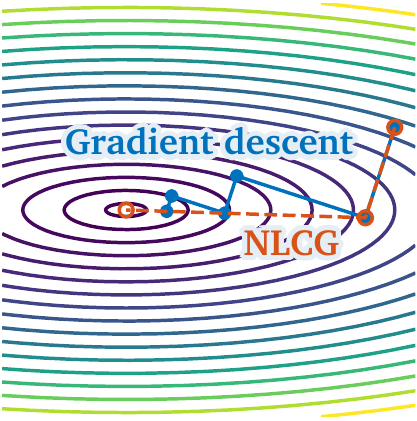}\qquad\optImage{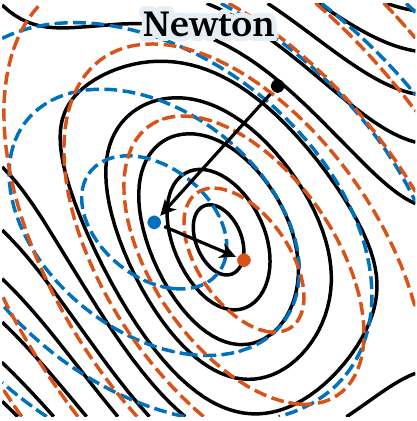}\qquad\optImage{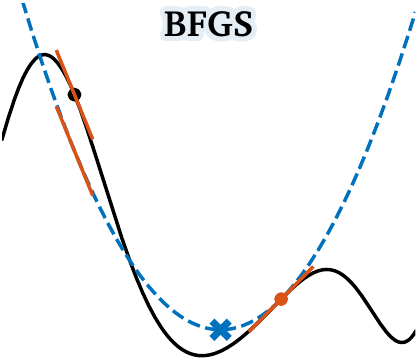}%
  \\[3mm]%
  \optImage{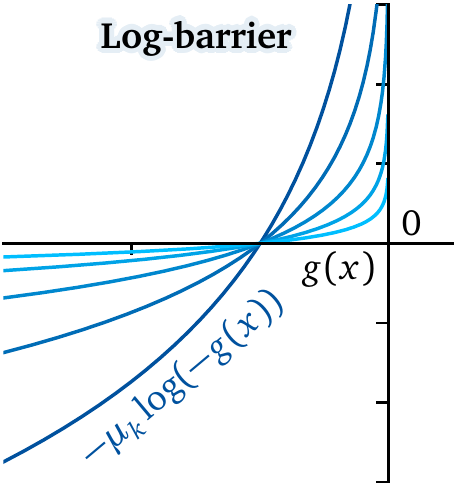}\quad\optImage{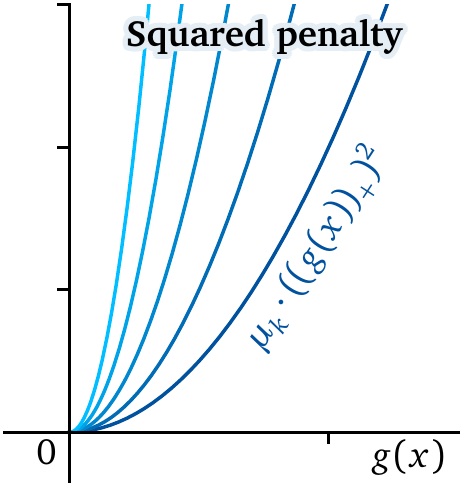}\quad\optImage{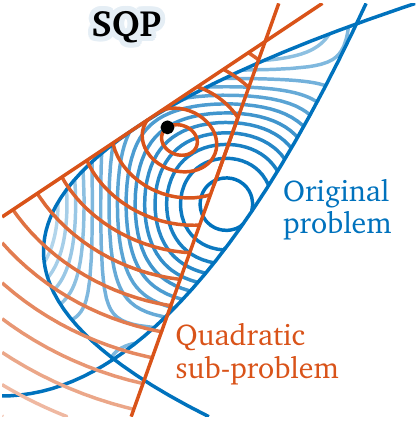}%
  \caption[Ideas of various gradient-based optimization methods]{%
    Sketch of the ideas of gradient-based or constrained
    optimization methods.%
  }%
  \label{fig:optimizationMethodGradientBased}%
\end{figure}

Most gradient-based optimization algorithms determine in
each iteration $k$ a unit \term{search direction} $\*d_k \in \real^d$
($\norm[2]{\*d_k} = 1$) to update the current iterate $\*x_k$:
\begin{equation}
  \*x_k
  \to \*x_{k+1}
  \ceq \*x_k + \delta_k \*d_k,\qquad
  \delta_k
  \ceq \argmin_{\delta \in \posreal} \objfun(\*x_k + \delta \*d_k),
\end{equation}
where $\delta_k \in \posreal$ is the \term{step size.}
The algorithms essentially differ in the
choice of the search direction $\*d_k$,
which should be oriented like the negative gradient
($\innerprod[2]{\*d_k}{\gradient{\*x}{\objfun}(\*x_k)} < 0$).
The step size $\delta_k$ can then be determined independently of the
algorithm via \term{line search,}
for instance, the \term{Armijo line search algorithm}
\multicite{Nocedal99Numerical,Ulbrich12Nichtlineare,Valentin14Hierarchische},
which uses a heuristic acceptance criterion
to find $\delta_k$ with a good enough improvement.

\paragraph{Gradient descent}

\term{Gradient descent}
\multicite{%
  Ulbrich12Nichtlineare,%
  Valentin14Hierarchische,%
  Toussaint15Introduction%
}
chooses $\*d_k
\propto -\gradient{\*x}{\!\objfun}(\*x_k)$ (i.e., normalized).
The method suffers from slow convergence,
if the Hessian $\hessian{\*x}{\objfun}$ is ill-conditioned:
One can show that for strictly convex quadratic functions,
the error $\objfun(\*x_k) - \objfun(\xopt)$
can decrease in each iteration only by the factor of
$(\tfrac{\lambda^{\max} - \lambda^{\min}}{\lambda^{\max} + \lambda^{\min}})^2$,
where $\lambda^{\min}$ and $\lambda^{\max}$ are the minimum and maximum
eigenvalue of $\hessian{\*x}{\objfun}$, respectively
\cite{Ulbrich12Nichtlineare}.
If the condition number
$\tfrac{\lambda^{\max}}{\lambda^{\min}}$
of $\hessian{\*x}{\objfun}$ is large,
then this factor will be very close to one.

\paragraph{NLCG}

A possible remedy for this issue is the method of
\term{non-linear conjugate gradients (NLCG)}
\multicite{%
  Nocedal99Numerical,%
  Valentin14Hierarchische,%
  Toussaint15Introduction%
}.
It is equivalent to the CG method for
solving symmetric positive definite
linear systems $\mat{A} \*x = \*b$,
if we optimize the strictly convex quadratic function
$\objfun(\*x) \ceq \frac{1}{2} \tr{\*x} \mat{A} \*x - \tr{\*b} \*x$
\multicite{Reinhardt13Nichtlineare,Valentin14Hierarchische}, i.e.,
it finds the optimum after only $d$ steps for strictly convex
quadratic functions.
The NLCG method quickly converges even for non-convex objective functions,
as due to the Taylor theorem,
three times continuously differentiable functions
with positive definite Hessian are ``similar'' to a
strictly convex quadratic function in a neighborhood of $\xopt$
\cite{Valentin14Hierarchische}.

\paragraph{Newton}

The \term{Newton method}
\multicite{%
  Ulbrich12Nichtlineare,%
  Valentin14Hierarchische,%
  Toussaint15Introduction%
}
replaces the objective function with the second-order Taylor approximation
given by
$\objfun(\*x_k + \*d_k)
\!\approx\! \objfun(\*x_k) +
\tr{(\gradient{\*x}{\objfun}(\*x_k))} \*d_k \,+
\frac{1}{2} \tr{(\*d_k)} (\hessian{\*x}{\objfun}(\*x_k)) \*d_k$
and determines the search direction such that $\*x_k + \*d_k$ is
the minimum of the approximation, i.e.,
$\*d_k \propto
-(\hessian{\*x}{\objfun}(\*x_k))^{-1} \gradient{\*x}{\objfun}(\*x_k)$.
Despite converging for strictly convex quadratic functions in a single step,
the Hessian must not be ill-conditioned for the Newton method as well,
as we have to solve a linear system with the matrix
$\hessian{\*x}{\objfun}(\*x_k)$.
Hence, often a \term{regularization/damping term} $\lambda \eye$
for some $\lambda > 0$ is added to the Hessian.

\paragraph{BFGS}

The Newton method has the disadvantage that it needs to evaluate the
Hessian $\hessian{\*x}{\objfun}$,
which may be unavailable or too expensive.
\term{Quasi-Newton methods} such as the method of
\term{BFGS (Broyden, Fletcher, Goldfarb, Shanno)}
\multicite{%
  Nocedal99Numerical,%
  Ulbrich12Nichtlineare,%
  Toussaint15Introduction%
}
approximate the Hessian by a solution of the secant equation
$\hessian{\*x}{\objfun}(\*x_k) (\*x_k - \*x_{k-1}) \approx
\gradient{\*x}{\objfun}(\*x_k) - \gradient{\*x}{\objfun}(\*x_{k-1})$.
As the solution is not unique for $d > 1$,
Quasi-Newton methods differ in which solution to choose.
The BFGS method performs a simple rank-one update.

\paragraph{Rprop}

\term{Rprop (resilient propagation)}
\multicite{Riedmiller93Direct,Toussaint15Introduction}
considers the gradient entries $(\gradient{\*x}{\objfun}(\*x_k))_t$
of each dimension $t = 1, \dotsc, d$ separately
and updates the entries $x_{k,t}$ of $\*x_k$
according to the sign of the respective gradient entry,
while adapting the step size dimension-wise.
Although the algorithm is independent of the exact direction
of $\gradient{\*x}{\objfun}(\*x_k)$,
it was found to often work robustly in machine learning scenarios
\cite{Toussaint15Introduction}.

\paragraph{Levenberg--Marquardt}

The \term{Levenberg--Marquardt method}
\multicite{Nocedal99Numerical,Freund07Stoer,Toussaint15Introduction}
can only solve \term{non-linear least-squares problems,} i.e.,
the objective function must be of the form
$
  \objfun(\*x)
  = \norm[2]{\*\phi(\*x)}^2
  = \sum_{i=1}^{m_{\*\phi}} \abs{\phi_i(\*x)}^2
$
for some function $\*\phi\colon \real^d \to \real^{m_{\*\phi}}$.
It is an improvement over the \term{Gauss--Newton method}
(which is in turn a slight modification of the Newton method)
and can be obtained by replacing the line search in the
Gauss--Newton method with a \term{trust-region approach.}

\subsection{Constrained Optimization Methods}
\label{sec:513gradientBasedConstrained}

Methods for constrained optimization usually
solve a series of unconstrained \term{auxiliary problems} with an arbitrary
unconstrained optimization method.
The auxiliary function to be minimized is
the sum of the objective function and \term{penalty terms,}
which penalize if the current point $\*x_k$ is near the boundary
of the feasible domain or even outside.
The penalty terms slowly increase to enforce
the feasibility of the final result.
Constrained optimization methods can roughly be divided
into \term{interior-point or barrier methods,}
where $\*x_k$ always stays in the feasible domain,
and \term{penalty methods,}
where intermediate solutions $\*x_k$ may be infeasible,
in which case the penalty term is applied.

At least for the interior-point methods,
a feasible initial solution $\*x_0$ is required.
We can find an initial solution by solving another auxiliary problem
\cite{Toussaint15Introduction}, for instance
\begin{equation}
  \min_{(\*x, s) \in \real^{d+1}} s
  \quad\text{s.t.}\quad
  s \ge 0,\;\;
  \ineqconfun(\*x) \le s \cdot \*1_{m_{\ineqconfun}},
\end{equation}
where $\*1_{m_{\ineqconfun}} \in \real^{m_{\ineqconfun}}$
is the all-one vector.
An initial solution for this problem can be explicitly given
(for example, $\*x_0 = \*0$ and
$s_0 = \max(\max(\ineqconfun(\*x_0)), 0)$).

\paragraph{Log-barrier}

The \term{log-barrier method}
\multicite{Boyd04Convex,Reinhardt13Nichtlineare,Toussaint15Introduction}
is an interior-point method that adds
a logarithmic \term{barrier function term} to the objective function
near the boundary.
The method solves
$\min\, [\objfun(\*x) - \mu_k \sum_{i=1}^{m_{\ineqconfun}} \log(-g_i(\*x))]$
for some decreasing $\mu_k \in \posreal$.

\paragraph{Squared penalty}

The \term{squared penalty method}
\multicite{%
  Polak71Computational,%
  Ulbrich12Nichtlineare,%
  Toussaint15Introduction%
}
replaces the constrained problem with the penalized problem
$\min\, [\objfun(\*x) + \mu_k \norm[2]{\nonnegpart{\ineqconfun(\*x)}}^2]$,
where $\mu_k \in \posreal$ is a penalty parameter and
$\nonnegpart{\cdot} \ceq \vecmax(\cdot, \*0)$ denotes the non-negative part.
With increasing $\mu_k$, the constraint violation of the solution of
the penalized problem decreases, although it may happen
that it never vanishes.

\paragraph{Augmented Lagrangian}

The method of the \term{augmented Lagrangian}
\multicite{Reinhardt13Nichtlineare,Toussaint15Introduction}
considers the auxiliary problem
\begin{equation}
  \min_{\*x \in \real^d} \bracket*{
    \objfun(\*x) + \mu_k \sum_{i=1}^{m_{\ineqconfun}} [\lambda_{k,i} > 0]
    ((\ineqconfun[i](\*x))_{+})^2 + \tr{\*\lambda_k} \ineqconfun(\*x)
  },
\end{equation}
where $[\lambda_{k,i} > 0] \in \{0, 1\}$ is defined as one
if and only if $\lambda_{k,i} > 0$, and
$\*\lambda_k \in \nonnegreal^{m_{\ineqconfun}}$ is an estimate of the
\term{Lagrangian multipliers.}
They are updated according to the penalty of the previous iteration,
generating a ``virtual gradient'' that drastically decreases
the necessary magnitude of the penalty parameter $\mu_k$
to achieve feasibility of the solution \cite{Toussaint15Introduction}.

\paragraph{Sequential quadratic programming (SQP)}

\term{Sequential quadratic programming (SQP) methods}
\multicite{%
  Ulbrich12Nichtlineare,%
  Reinhardt13Nichtlineare,%
  Toussaint15Introduction%
}
are one of the most powerful method classes for constrained optimization.
They are motivated by the \term{Karush--Kuhn--Tucker (KKT) conditions,}
which are necessary to hold in any optimal point
(similarly to critical points in unconstrained optimization).
The Newton method can be employed to solve the KKT conditions when
written as a non-linear system of equations.
The linear system of the resulting \emph{Newton--Lagrange method}
is equivalent to the KKT conditions of a
\term{quadratic programming (QP) problem,}
for which objective and constraint functions have the
form $\objfun(\*x) = \frac{1}{2} \tr{\*x} \mat{Q} \*x + \tr{\*d} \*x$ and
$\ineqconfun(\*x) = \mat{A} \*x - \*b$, respectively.

\section{Optimization of Surrogates on Sparse Grids}
\label{sec:52method}

\minitoc{72mm}{5}

\noindent
The methods presented in the last section can be combined to a
``meta-method'' for surrogate optimization.
The surrogates are constructed as interpolants on spatially adaptive
sparse grids, which we explain in the following.

\subsection{Novak--Ritter Adaptivity Criterion}
\label{sec:521novakRitter}

The classic surplus-based refinement strategy for
spatially adaptive sparse grids is not tailored to optimization,
as this refinement strategy aims to minimize the overall $\Ltwo$ error.
However, in optimization, it is reasonable to generate more
points in regions where we suspect the global minimum to be
to increase the interpolant's accuracy in these regions.
Hence, we employ an adaptivity criterion proposed by
Novak and Ritter \cite{Novak96Global} for hyperbolic cross points.
The Novak--Ritter criterion has also been applied to sparse grids
\multicite{Ferenczi05Globale,Valentin14Hierarchische,Valentin16Hierarchical}.

\paragraph{$m$-th order children}

As usual, the criterion works iteratively:
Starting with an initial regular sparse grid of a very coarse level,
the criterion selects a specific point $\gp{\*l,\*i}$ in each iteration
and inserts all its children into the grid.
This process is repeated until a given number $\ngpMax$ of grid points is
reached,
since we evaluate $\objfun$ at every grid point once, and we assume that
function evaluations dominate the overall complexity.
The difference to common refinement criteria is that
a point may be selected multiple times, in which case
\term{higher-order children} are inserted.
The $m$-th order children $\gp{\*l',\*i'}$ of a grid point $\gp{\*l,\*i}$
satisfy
\begin{equation}
  \label{eq:indirectChild}
  \*l'_{-t} = \*l^{}_{-t},\;\,
  \*i'_{-t} = \*i^{}_{-t},\;\,
  l'_t = l^{}_t + m,\;\,
  i'_t \in
  \begin{cases}
    \{1\},&(l_t = 0) \land (i_t = 0),\\
    \{2^m - 1\},&(l_t = 0) \land (i_t = 1),\\
    \{2^m i_t - 1,\, 2^m i_t + 1\},&\hphantom{(}l_t > 0,
  \end{cases}
\end{equation}
where $m \in \nat$ and $t \in \{1, \dotsc, d\}$
(cf.\ \cref{eq:directAncestor} for $m = 1$).
The order is chosen individually for each child point to be inserted
as the lowest number $m$ such that $\gp{\*l',\*i'}$ does not yet exist
in the grid.

\paragraph{Criterion}

The Novak--Ritter refinement criterion \cite{Novak96Global}
refines the grid point $\gp{\*l,\*i}$ that minimizes the product%
\footnote{%
  Compared to \cite{Novak96Global},
  we added one in the base of each factor to avoid ambiguities
  for $0^0$.
  In addition, we swapped $\gamma$ with $1-\gamma$
  to make $\gamma$ more consistent with its name as adaptivity.%
}
\begin{equation}
  (r_{\*l,\*i} + 1)^\gamma \cdot
  (\normone{\*l} + d_{\*l,\*i} + 1)^{1 - \gamma}.
\end{equation}
Here, $r_{\*l,\*i} \ceq \setsize{
  \{(\*l',\*i') \in \liset \mid
  \objfun(\gp{\*l',\*i'}) \le \objfun(\gp{\*l,\*i})\}
} \in \{1, \dotsc, \setsize{\liset}\}$ is the \term{rank} of $\gp{\*l,\*i}$
(where $\liset$ is the current set of level-index pairs of the grid), i.e.,
the place of the function value at $\gp{\*l,\*i}$
in the ascending order of the function values at all points
of the current grid.
We denote the \term{degree} $d_{\*l,\*i} \in \natz$ of $\gp{\*l,\*i}$
as the number of previous refinements of this point.
Finally, $\gamma \in \clint{0, 1}$ is the \term{adaptivity parameter.}
We have to choose a suitable compromise between exploration ($\gamma = 0$)
and exploitation ($\gamma = 1$).
The best choice of course depends on the objective function $\objfun$ at hand,
but for our purposes, we choose a priori a value of $\gamma = 0.15$.
However, it may be an option to adapt the value of $\gamma$ automatically
during the grid generation phase.

\subsection{Global Optimization of Sparse Grid Surrogates}
\label{sec:522method}

\paragraph{Global, local, and globalized optimization methods}

In \cref{sec:51algorithms}, we presented various optimization methods
for the unconstrained case,
divided into global gradient-free methods such as differential evolution and
local gradient-based methods, for example, gradient descent.
A subset of these methods has been implemented in \sgpp{}
\cite{Pflueger10Spatially}, see \cref{tbl:optimizationMethod}.
The gradient-based methods need an initial point, and
they may get stuck in local minima.
Hence, we additionally implemented globalized versions
of the gradient-based methods
via a multi-start Monte Carlo approach with $m \ceq \min(10d, 100)$
uniformly distributed pseudo-random initial points.%
\footnote{%
  We split the number of permitted function evaluations evenly
  among the $m$ parallel calls.%
}
This means there are three types of methods:

\begin{enumerate}[label=T\arabic*.,ref=T\arabic*,leftmargin=2.7em]
  \item
  \label{item:globalMethods}
  Global gradient-free methods listed as implemented in
  \cref{tbl:optimizationMethod}
  
  \item
  \label{item:localMethods}
  Local gradient-based methods listed as implemented in
  \cref{tbl:optimizationMethod}%
  \footnote{%
     Excluding Levenberg--Marquardt, which is only applicable
     to least-squares problems.%
  }
  
  \item
  \label{item:globalizedMethods}
  Globalized versions of the methods of type \ref{item:localMethods}
\end{enumerate}

\paragraph{Unconstrained optimization of sparse grid surrogates}

Given the objective function $\objfun\colon \clint{\*0, \*1} \to \real$,
the maximal number $\ngpMax \in \nat$ of evaluations of $f$, and
the adaptivity parameter $\gamma \in \clint{0, 1}$,
we determine an approximation $\xoptappr \in \clint{\*0, \*1}$
of the global minimum $\xopt$ of $\objfun$ as follows:

\begin{enumerate}
  \item
  Generate a spatially adaptive sparse grid $\sgset$
  with the Novak--Ritter refinement criterion
  for $\objfun$, $\ngpMax$, and $\gamma$.
  
  \item
  Determine the sparse grid interpolant $\sgintp$ of $\objfun$
  by solving the linear system \eqref{eq:hierarchizationProblem}.
  
  \item
  Optimize the interpolant:
  First, find the best grid point
  $\*x^{(0)} \ceq \vecargmin_{\gp{\*l,\*i} \in \sgset} \objfun(\*x_{\*l,\*i})$.
  Second, apply the local methods of type \ref{item:localMethods}
  to the interpolant $\sgintp$ with $\*x^{(0)}$ as initial point.
  Let $\*x^{(1)}$ be the resulting point with minimal objective function value.
  Third, we apply the global and globalized methods
  of types \ref{item:globalMethods} and \ref{item:globalizedMethods}
  to the interpolant $\sgintp$.
  Again, let $\*x^{(2)}$ be the point with
  minimal $\objfun$ value.
  Finally, determine the point of $\{\*x^{(0)}, \*x^{(1)}, \*x^{(2)}\}$
  with minimal $\objfun$ value and return it as $\xoptappr$.
\end{enumerate}

\noindent
Note that the third step requires a fixed number of additional
evaluations of the objective function,
which can be neglected compared to $\ngpMax$.
By default, we use the cubic modified hierarchical not-a-knot B-spline basis
$\bspl[\nak,\modified]{\*l,\*i}{p}$ ($p = 3$)
for the construction of the sparse grid surrogate.
However, we could apply any of the hierarchical (B-)spline bases presented in
\cref{chap:30BSplines,chap:40algorithms}.

\paragraph{Comparison methods}

We use two comparison methods.
First, we apply the gradient-free methods
(type \ref{item:globalMethods}) to the sparse grid interpolant
using modified piecewise linear hierarchical basis functions
(i.e., $p = 1$) on the same sparse grid as the cubic B-splines.
We cannot employ gradient-based optimization as the objective function
should be continuously differentiable and
discontinuous derivatives are usually numerically problematic
for gradient-based optimization methods
(see, e.g., \cite{Huebner14Mehrdimensionale}).
Second, we apply the gradient-free methods
(type \ref{item:globalMethods}) directly to the objective function.
We cannot use the gradient-based methods here as the gradient of the
objective function is assumed to be unknown.
For both of the comparison methods,
we make sure that the objective function is evaluated at most $\ngpMax$ times
by splitting the $\ngpMax$ evaluations
evenly among all employed optimization methods.

\paragraph{Constrained optimization}

For optimization problems with constraints,
we proceed exactly as for unconstrained optimization,
except that for optimizing the interpolant, we use the
constrained optimization algorithms implemented in \sgpp
as listed in \cref{tbl:optimizationMethod}.
We only replace the objective function $\objfun$ with a sparse grid
surrogate $\sgintp$, and we assume that the constraint function
$\ineqconfun$ can be evaluated fast.
However, it would also be possible to replace $\ineqconfun$
with a sparse grid surrogate.
In this case, it cannot be guaranteed that the resulting optimal point
$\xoptappr$ is feasible, i.e., we could have
$\lnot(\ineqconfun(\xoptappr) \le \*0)$.

\section{Test Problems}
\label{sec:53testProblems}

It is impossible to assess the capability of optimization methods
for every possible optimization problem.
The most widespread approach in literature
is the selection of a subset of specific problems
with different characteristics \term{(test problems)} and
the application of the methods to only these problems,
in the hope that the methods perform similarly in
actual application settings.

\paragraph{Trivial test functions}

When testing methods that involve sparse grid interpolation,
one has to consider that the function to be interpolated
does not satisfy a specific \term{trivial property.}
A test function $\objfun\colon \clint{\*0, \*1} \to \real$ is trivial if
$\objfun$ is a sum of tensor products of which at
least one factor is a linear polynomial, i.e., if
$\objfun$ is of the form
{%
  \setlength{\abovedisplayskip}{9pt}%
  \setlength{\belowdisplayskip}{9pt}%
  \begin{equation}
    \objfun(\*x) \equiv \sum_{q=1}^m \prod_{t=1}^d \objfun_{q,t}(x_t),\quad
    m \in \natz,\;\;
    \objfun_{q,t}\colon \clint{0, 1} \to \real,\;\;
    \faex{q = 1, \dotsc, m}{t \in \{1, \dotsc, d\}}{
      \objfun_{q,t} \in \polyspace{1},
    }
  \end{equation}%
}%
where $\polyspace{1}$ is the space of univariate polynomials
up to linear degree.
This is already fulfilled if the summands of $\objfun(\*x)$
do not depend on all coordinates $x_t$ of $\*x$.
One can show that for hat functions on sparse grids,
the hierarchical surpluses $\surplus{\*l,\*i}$ for trivial functions
vanish if $\*l \ge \*1$.
This means that trivial functions can be well-approximated by hat functions
on sparse grids just with boundary points, without placing any points
in the interior.
As this would distort our results,
we avoid trivial test functions in the following,
which include popular functions such as the
Branin01, Rosenbrock, and Schwefel26 functions.

\paragraph{Selection of test problems}

\begin{figure}
  \subcaptionbox{%
    Bra02%
  }[71mm]{%
    \includegraphics{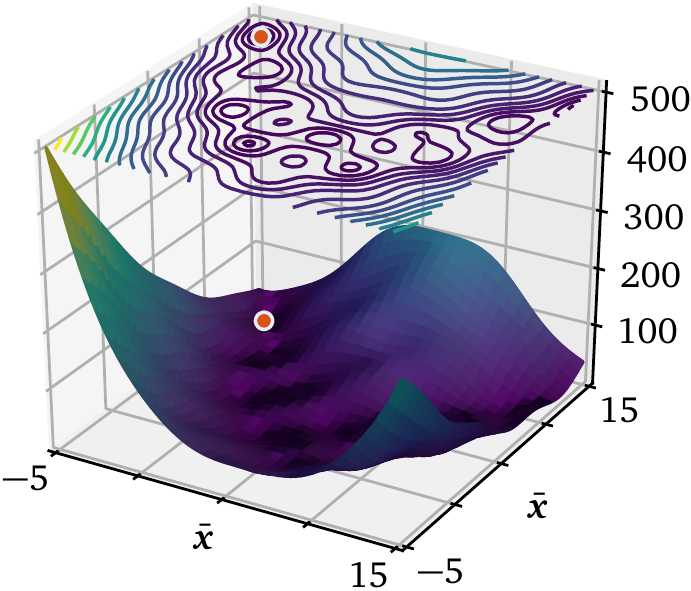}%
  }%
  \hfill%
  \subcaptionbox{%
    GoP%
  }[76mm]{%
    \includegraphics{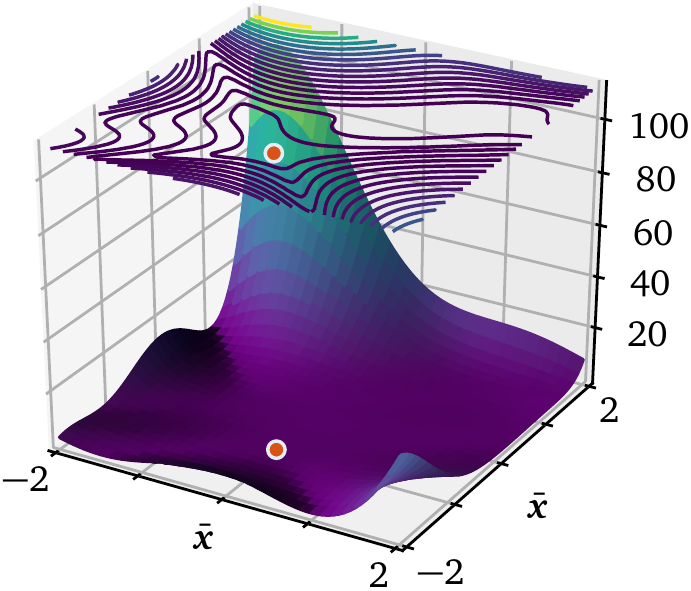}%
  }\\[2.5mm]%
  \subcaptionbox{%
    Sch06%
  }[71mm]{%
    \includegraphics{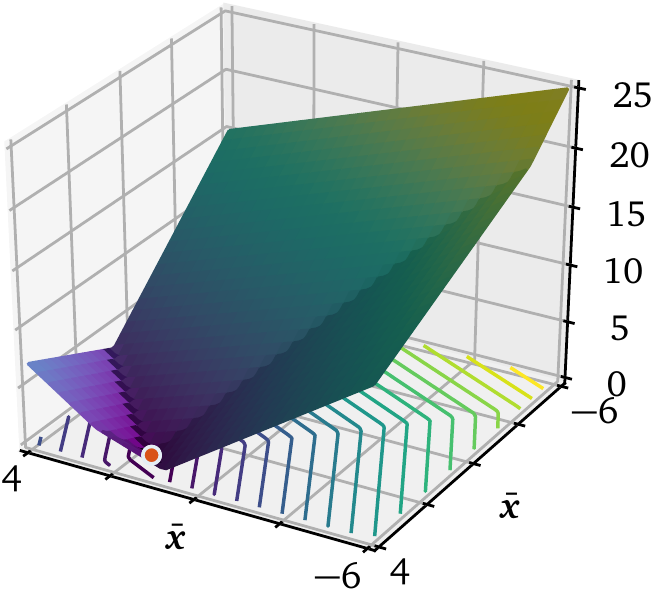}%
  }%
  \hfill%
  \subcaptionbox{%
    Ack for $d = 2$%
  }[76mm]{%
    \includegraphics{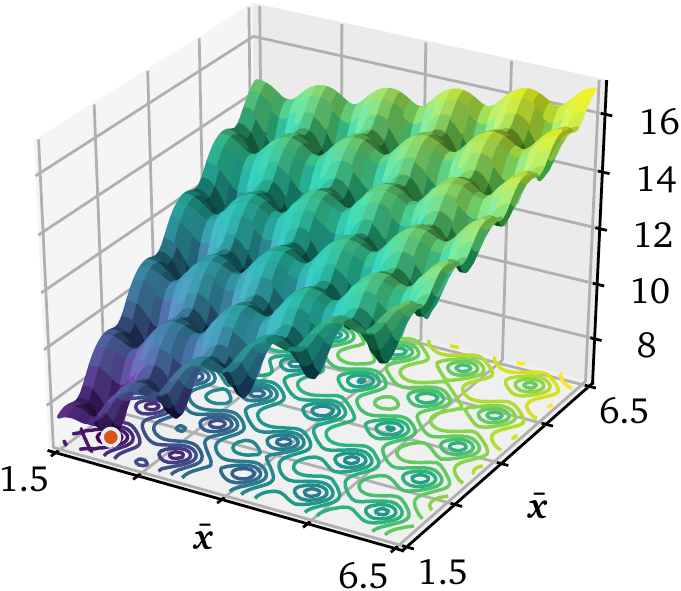}%
  }\\[2.5mm]%
  \subcaptionbox{%
    Alp02 for $d = 2$%
  }[71mm]{%
    \includegraphics{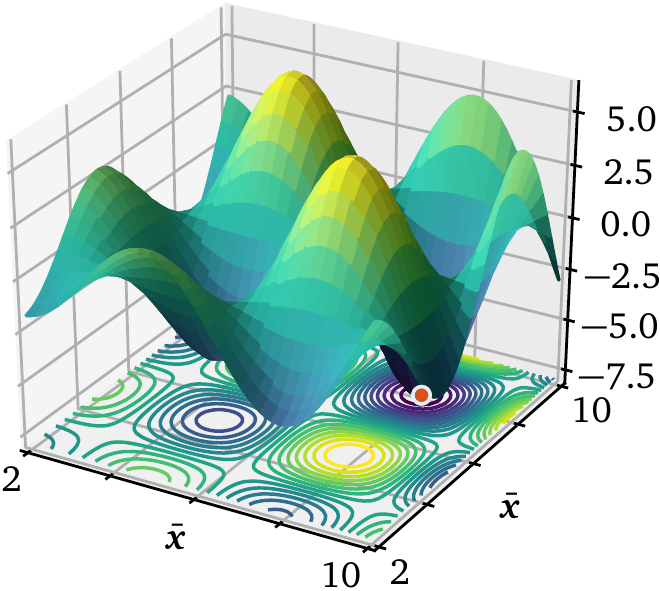}%
  }%
  \hfill%
  \subcaptionbox{%
    Sch22 for $d = 2$%
  }[76mm]{%
    \includegraphics{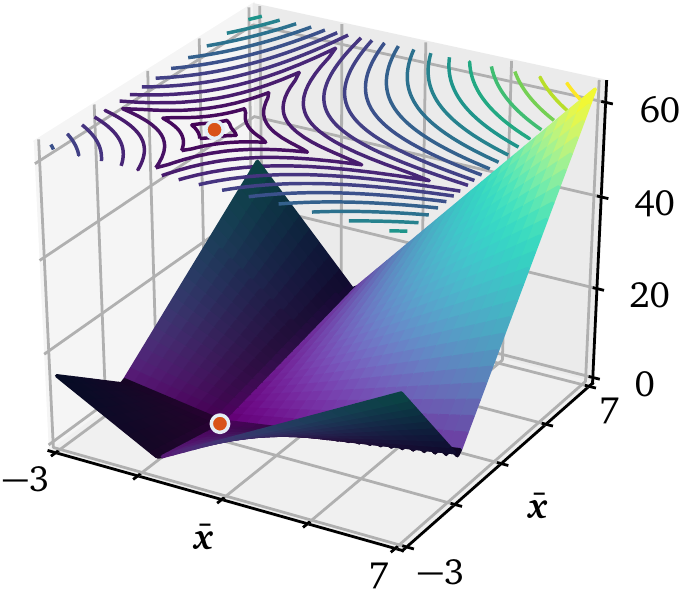}%
  }%
  \caption[%
    Unconstrained test problems%
  ]{%
    Bivariate test functions $\objfunscaled$ in unconstrained optimization.
    The \textcolor{C1}{red dot} indicates the location of the
    global minimum.%
  }%
  \label{fig:unconstrainedOptimizationProblem}%
\end{figure}

\begin{figure}
  \subcaptionbox{%
    G08%
  }[72mm]{%
    \includegraphics{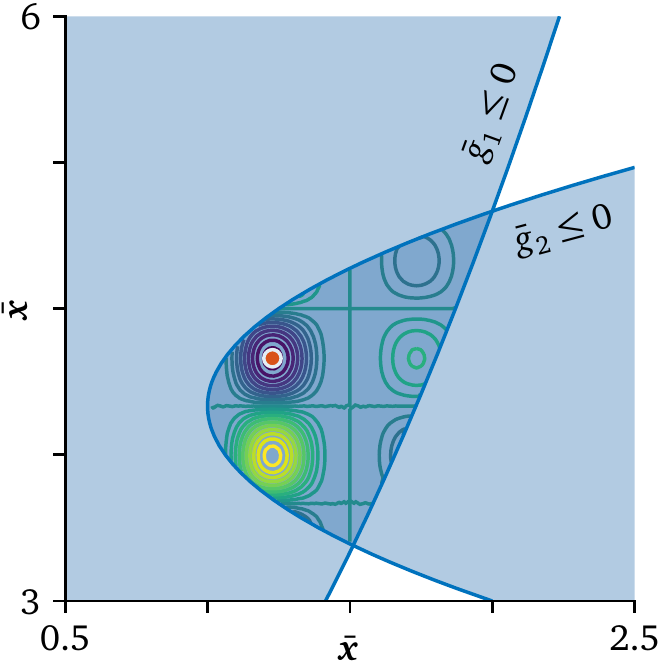}%
  }%
  \hfill%
  \subcaptionbox{%
    G04Sq (bivariate projection over $\xscaledentry{3}$ and $\xscaledentry{5}$
    onto $\xscaledentry{t} = \xoptscaledentry{t}$ for $t = 1, 2, 4$)%
  }[72mm]{%
    \includegraphics{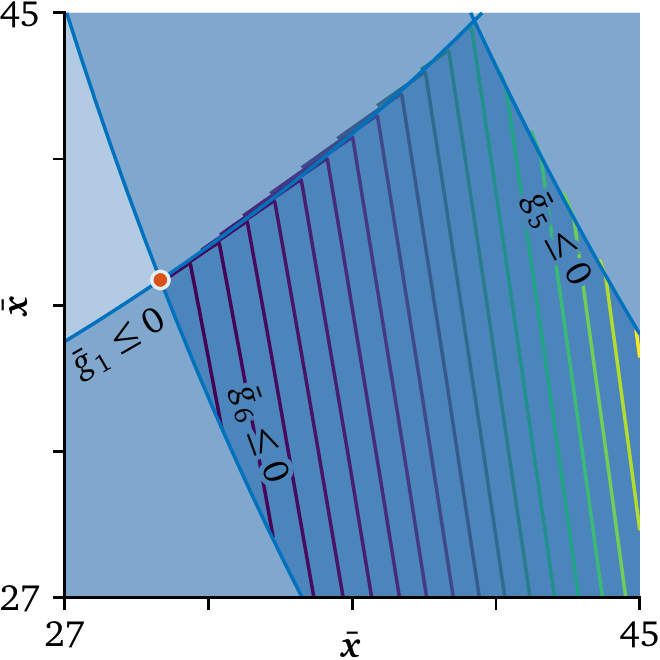}%
  }%
  \caption[%
    Constrained test problems%
  ]{%
    Test problems in constrained optimization.
    The \textcolor{C0}{blue areas} denote the inequality constraints and
    the \textcolor{C1}{red dot} indicates the location of the
    global minimum.%
  }%
  \label{fig:constrainedOptimizationProblem}%
\end{figure}

In the following, we select six unconstrained test problems
and two constrained test problems, which are listed in
\cref{tbl:optimizationProblem} and plotted in
\cref{fig:unconstrainedOptimizationProblem,fig:constrainedOptimizationProblem}.
The definitions of the problems are given in \cref{chap:a20testProblems}.
For the unconstrained case and the standard hierarchical
B-spline basis, a more exhaustive list of test functions has been
studied previously \cite{Valentin14Hierarchische}.
Gavana \cite{Gavana13Global} and Runarsson/Yao \cite{Runarsson00Stochastic}
provide a good overview of unconstrained and constrained test problems,
respectively.

\begin{table}
  \setnumberoftableheaderrows{1}%
  \begin{tabular}{%
    >{\kern\tabcolsep}=l<{\kern5mm}+l<{\kern5mm}*{5}{+c}%
    <{\kern5mm}+l<{\kern\tabcolsep}%
  }
    \toprulec
    \headerrow
    Name&           Abbr.& $d$& $m_{\ineqconfun}$& C&    CD&   MM&   Reference\\
    \midrulec
    Branin02&       Bra02& 2&   0&                 \yes& \yes& \yes& \cite{Munteanu98Global}\\
    GoldsteinPrice& GoP&   2&   0&                 \yes& \yes& \yes& \cite{Goldstein71Descent}\\
    Schwefel06&     Sch06& 2&   0&                 \yes& \no&  \no&  \cite{Schwefel77Numerische}\\
    Ackley&         Ack&   $d$& 0&                 \yes& \yes& \yes& \cite{Ackley87Connectionist}\\
    Alpine02&       Alp02& $d$& 0&                 \yes& \yes& \yes& \cite{Clerc99Swarm}\\
    Schwefel22&     Sch22& $d$& 0&                 \yes& \no&  \no&  \cite{Schwefel77Numerische}\\
    \midrulec
    G08&            G08&   2&   2&                 \yes& \yes& \yes& \cite{Schoenauer93Constrained}\\
    G04Squared&     G04Sq& 5&   6&                 \yes& \yes& \no&  \cite{Colville68Comparative}\\
    \bottomrulec
  \end{tabular}
  \caption[Selection of test problems in optimization]{%
    Unconstrained \emph{(top)} and constrained \emph{(bottom)} test problems.
    The columns state the full and abbreviated names,
    the dimensionality $d$ of the objective function $\objfun$,
    the number $m_{\ineqconfun}$ of constraints,
    whether $\objfun$ is continuous in the domain
    $\clint{\*0, \*1}$ (C),
    whether $\objfun$ is continuously differentiable in the domain
    $\clint{\*0, \*1}$ (CD),
    whether $\objfun$ is multi-modal (MM, i.e.,
    whether there are multiple local minima), and
    a reference to the original literature that defines the problem.%
  }%
  \label{tbl:optimizationProblem}%
\end{table}

For each test problem, we state unscaled versions of objective functions
$\objfunscaled\colon \clint{\*a, \*b} \to \real$,
$\xscaled \mapsto \objfunscaled(\xscaled)$
(and the unscaled constraint function $\ineqconfunscaled$, if present).
The actual objective function $\objfun\colon \clint{\*0, \*1} \to \real$
can be obtained by $\objfun(\*x) \ceq \objfunscaled(\xscaled)$
with the affine parameter transformation
$x_t = \tfrac{\xscaledentry{t} - a_t}{b_t - a_t}$, $t = 1, \dotsc, d$
(similarly for the constraint function).

\vspace*{1em}

The parameter domains of some test problems have been slightly translated
compared to the literature
to avoid that the minima are located exactly at or close to
the center of the domain.
In these cases, sparse grids would be in advantage as
they tend to place more points near the center of the domain
(especially for high dimensionalities).

\fillsectionornament
\section{Numerical Results}
\label{sec:54results}

\minitoc{85mm}{5}

\noindent
The following numerical experiments can be roughly divided into two parts.
First, we study interpolation errors for the test functions
to assess the effects of the hierarchical B-spline bases introduced in
\cref{chap:20sparseGrids,chap:30BSplines} on interpolation.
Second, we consider the optimality gaps $\objfun(\xoptappr) - \objfun(\xopt)$
of the calculated approximations $\xoptappr$
of the point $\xopt$ at which the objective function $\objfun$
is minimal.

\pagebreak

The results have been computed with the sparse grid toolbox \sgpp
\cite{Pflueger10Spatially},%
\footnote{%
  \url{http://sgpp.sparsegrids.org/}%
}
which has been extended in the scope of this thesis.
The new code has been written in such a way that
it is scalable and efficient, while still being maintainable and
portable \cite{Pflueger16Scalability}.

\subsection{Interpolation Error and Decay of Surpluses}
\label{sec:541interpolation}

\paragraph{Interpolation error for different test functions}

\begin{figure}
  \includegraphics{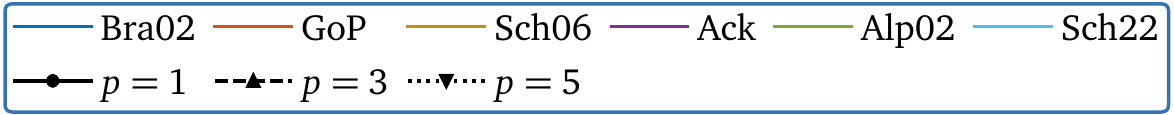}\\[2mm]%
  \includegraphics{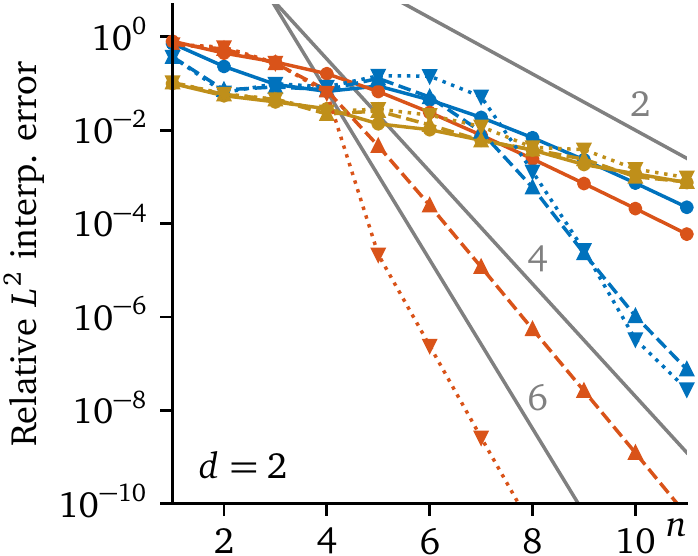}%
  \hfill%
  \includegraphics{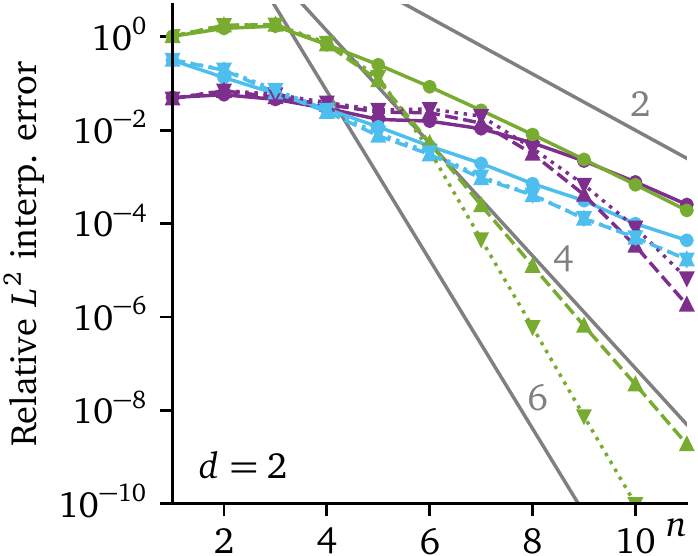}%
  \\[2mm]%
  \includegraphics{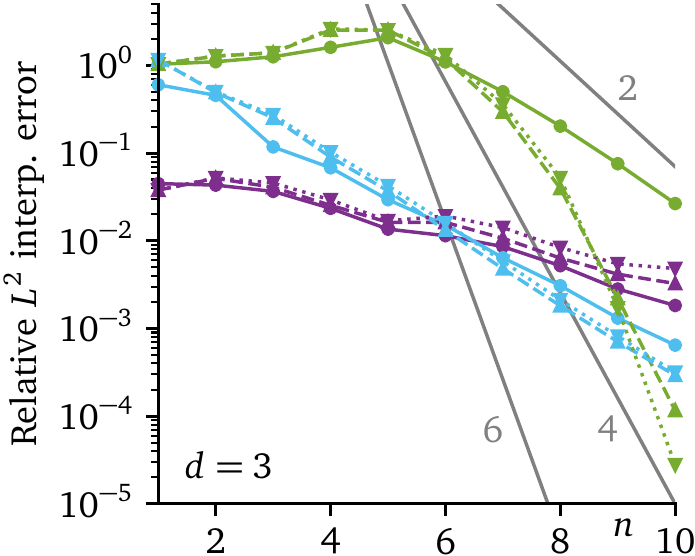}%
  \hfill%
  \includegraphics{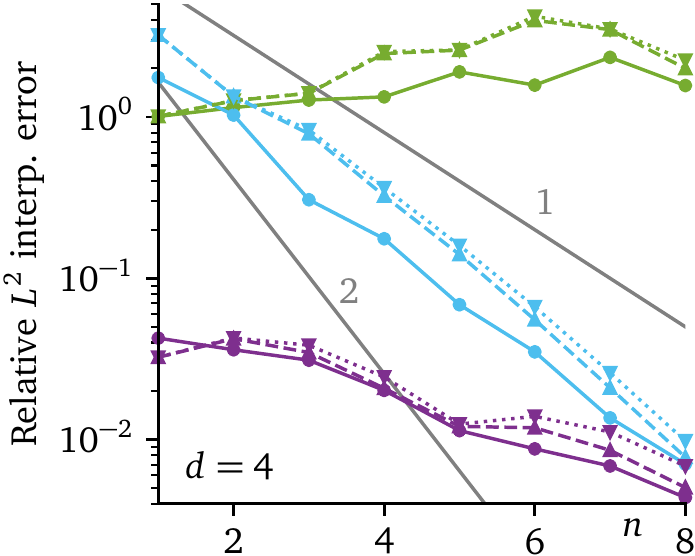}%
  \caption[Relative interpolation error for different test functions]{%
    Relative $\Ltwo$ interpolation error
    \vspace{-0.05em}%
    $\normLtwo{\objfun - \sgintp}/\normLtwo{\objfun}$
    for different test functions $\objfun$ \emph{(colors)}
    using hierarchical not-a-knot B-splines
    $\bspl[\nak]{\*l,\*i}{p}$ of different degrees $p$
    \emph{(line styles/markers)} on
    regular sparse grids $\regsgset{n}{d}$ of different levels $n$.%
  }%
  \label{fig:resultsInterpolationErrorTestFunctions}%
\end{figure}

\Cref{fig:resultsInterpolationErrorTestFunctions} shows the
relative $\Ltwo$ interpolation error
$\tfrac{\normLtwo{\objfun - \sgintp}}{\normLtwo{\objfun}}$
of sparse grid interpolants $\sgintp$ to
different objective functions $\objfun$
(approximated via Monte Carlo quadrature using
$10^4$ uniformly pseudo-random samples).
The interpolation is performed on regular sparse grids of increasing levels
using hierarchical not-a-knot B-splines $\bspl[\nak]{\*l,\*i}{p}$
of degree $p = 1, 3, 5$.
As a visual aid, the plots include gray lines that indicate different
orders of convergence.

It is already known that---%
if the objective function is sufficiently smooth---%
the $\Ltwo$ error of spline interpolants of degree $p$ on
$d$-dimensional regular sparse grids of level $n$
asymptotically behaves like
$\landauO{\ms{n}^{p+1} (\log_2 \ms{n}^{-1})^{d-1}}
= \landauO{2^{-(p+1)n} n^{d-1}}$ for $n \to \infty$ \cite{Sickel11Spline}.
We can numerically verify this fact easily with
\cref{fig:resultsInterpolationErrorTestFunctions},
in which we obtain the asserted orders of convergence
\pagebreak%
for the bivariate functions that are continuously differentiable.
For the functions Sch06 and Sch22, which have a non-differentiable kink,
only linear convergence can be achieved regardless of the B-spline degree.

\vspace*{\fill}

The region where the asymptotic behavior dominates largely depends
on the objective function at hand.
Functions like Bra02 and Ack with many small oscillations
require more interpolation points than ``smoother'' functions like
GoP and Alp02.
This is also the case for all functions in higher dimensionalities,
as more interpolation points are necessary to sufficiently explore the domain
(curse of dimensionality).
In \cref{fig:resultsInterpolationErrorTestFunctions}, this can already be seen
for $d \ge 3$.
This is not a consequence of employing higher-order B-splines for
the hierarchical basis.
However, it seems that higher-order B-splines lead to a slight increase
of the interpolation error in the preasymptotic range.

\pagebreak

\paragraph{Interpolation error for different basis functions}

\begin{SCfigure}
  \begin{minipage}{98mm}%
    \hspace*{4mm}%
    \includegraphics{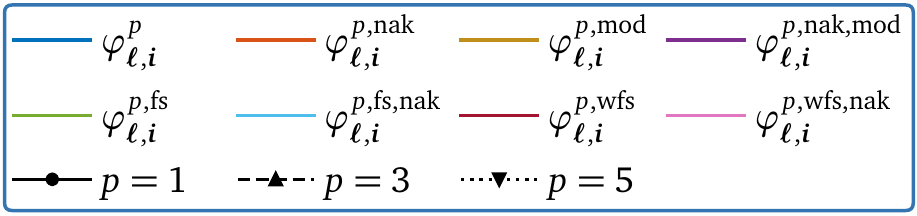}\\[2mm]%
    \includegraphics{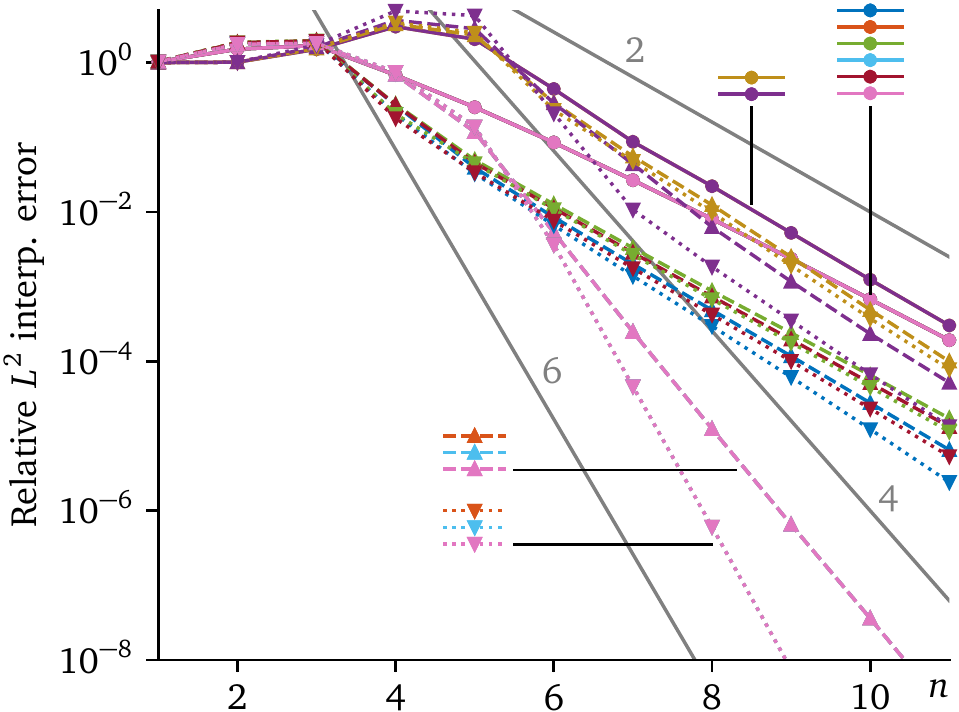}%
  \end{minipage}%
  \caption[Relative interpolation error for different basis functions]{%
    Relative $\Ltwo$ interpolation error
    $\normLtwo{\objfun - \sgintp}/\normLtwo{\objfun}$
    for the bivariate Alp02 function ($d = 2$)
    using different hierarchical basis functions
    $\basis{\*l,\*i}$ \emph{(colors)}
    of different degrees $p$ \emph{(line styles/markers)} and
    regular sparse grids $\regsgset{n}{d}$ of different levels $n$.\\
    The basis functions shown here involve
    standard \emph{(no superscript),}
    not-a-knot ($\mathrm{nak}$),
    modified ($\mathrm{mod}$),
    fundamental ($\mathrm{fs}$), and
    weakly fundamental ($\mathrm{wfs}$)
    splines as well as the combinations
    introduced in \cref{chap:20sparseGrids,chap:30BSplines}.%
  }%
  \label{fig:resultsInterpolationErrorBasisFunctions}%
\end{SCfigure}

In \cref{fig:resultsInterpolationErrorBasisFunctions},
we fix the objective function and study the influence of the choice
of hierarchical basis functions on the interpolation error.
Shown are eight types of hierarchical B-spline bases as introduced in
\cref{chap:20sparseGrids,chap:30BSplines} for the degrees $p = 1, 3, 5$.
Note that some lines exactly overlap, which is indicated in
the figure.

For $p = 1$, the non-modified bases and the modified bases coincide.
For higher degrees, the modified bases show worse results than
the corresponding non-modified versions for the same level $n$.
However, modified bases need significantly less grid points
(no boundary points),
which means that a direct comparison based on the sparse grid level $n$
is somewhat skewed.
In addition, we see that the not-a-knot bases coincide exactly for $p > 1$,
as they span the same space for regular and
dimensionally adaptive sparse grids.
Only with the not-a-knot boundary conditions, we obtain the true
theoretical order of convergence, which is $p + 1$ for degree $p$.
Otherwise, only quadratic convergence can be achieved regardless of $p$,
albeit with a smaller constant (offset).

\vspace*{-0.5em}

\paragraph{Pointwise interpolation error}

\begin{figure}
  \includegraphics{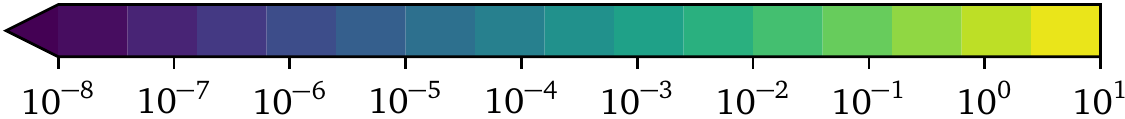}\\[2mm]%
  \subcaptionbox{%
    $\bspl{\*l,\*i}{p}\vphantom{\bspl[\nak]{\*l,\*i}{p}}$%
  }[48mm]{%
    \includegraphics{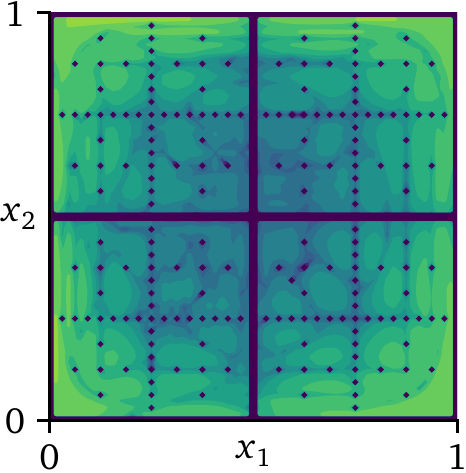}%
  }%
  \hfill%
  \subcaptionbox{%
    $\bspl[\nak]{\*l,\*i}{p}$%
  }[48mm]{%
    \includegraphics{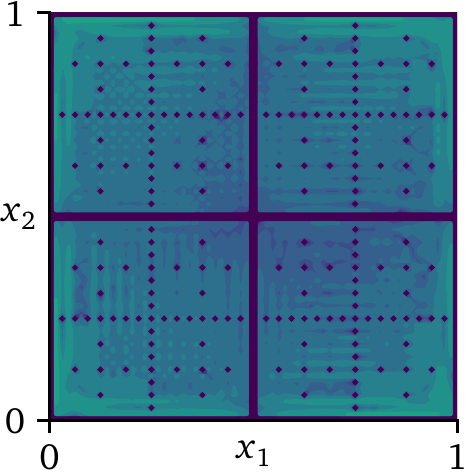}%
  }%
  \hfill%
  \subcaptionbox{%
    $\bspl[\modified]{\*l,\*i}{p}$%
  }[48mm]{%
    \includegraphics{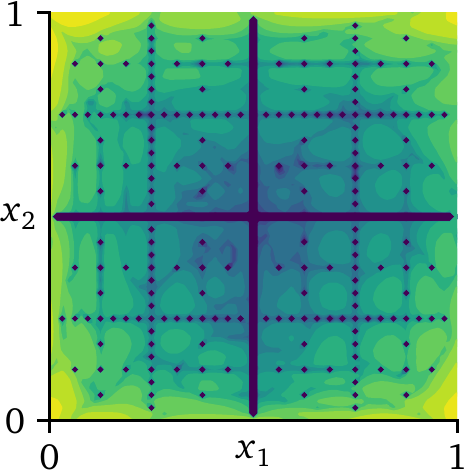}%
  }%
  \caption[Pointwise interpolation error for the GoP function]{%
    Pointwise interpolation error
    $\abs{\objfun(\*x) - \sgintp(\*x)}$ on a logarithmic scale
    for the bivariate GoP function ($d = 2$)
    using different hierarchical basis functions
    $\basis{\*l,\*i}$ \emph{(left, center, right)} of degree $p = 3$ on
    the regular sparse grid $\regsgset{n}{d}$ of level $n = 7$.%
  }%
  \label{fig:resultsInterpolationErrorPointwise}%
\end{figure}

The importance of not-a-knot boundary conditions is also evident
from plots of the pointwise interpolation error as in
\cref{fig:resultsInterpolationErrorPointwise}.
The interpolation error grows for the standard hierarchical
B-spline basis $\bspl{\*l,\*i}{p}$
as we move towards the boundary of the domain $\clint{\*0, \*1}$,
before dropping to zero or near-zero values at or near boundary grid points.
With not-a-knot B-splines $\bspl[\nak]{\*l,\*i}{p}$,
the interpolation error is uniformly low.
For comparison, modified B-splines $\bspl[\modified]{\*l,\*i}{p}$
incur even worse issues near the boundary, since the corresponding sparse grids
do not contain boundary points.

\vspace*{-0.3em}

\paragraph{Decay of surpluses}

In the piecewise linear case ($p = 1$),
the hierarchical surpluses $\surplus{\*l,\*i}$
can be represented as the $\Ltwo$ inner product of
the corresponding hat function $\bspl{\*l,\*i}{1}$ with the
second mixed derivative
$\partialderiv[2d]{\partialdiff x_1^2 \dotsm \partialdiff x_d^2}{\objfun}$
of the objective function $\objfun$,
if $\*l \ge \*1$ and if this derivative exists and is continuous
(see \cref{eq:surplusIntegral}).
Consequently, one can prove that
$\abs{\surplus{\*l,\*i}} \le 2^{-d} 2^{-2\normone{\*l}}
\normLinftyscaled{
  \partialderiv[2d]{\partialdiff x_1^2 \dotsm \partialdiff x_d^2}{\objfun}
}$ \cite{Bungartz04Sparse},
i.e., the absolute values of the hierarchical surpluses
decay  in quadratic order with the level sum $\normone{\*l}$.
This relation can be used to estimate the convergent range
of the corresponding interpolation error (\cref{%
  fig:resultsInterpolationErrorTestFunctions,%
  fig:resultsInterpolationErrorBasisFunctions%
}).
A generalization of this estimate to higher B-spline degrees $p > 1$
is not straightforward, as the surpluses $\surplus{\*l,\*i}$
then also depend on function values $\objfun(\gp{\*l',\*i'})$ at
grid points of higher levels $\*l' \ge \*l$.

\begin{figure}
  \includegraphics{resultsInterpolationLegend_1}\\[2mm]%
  \includegraphics{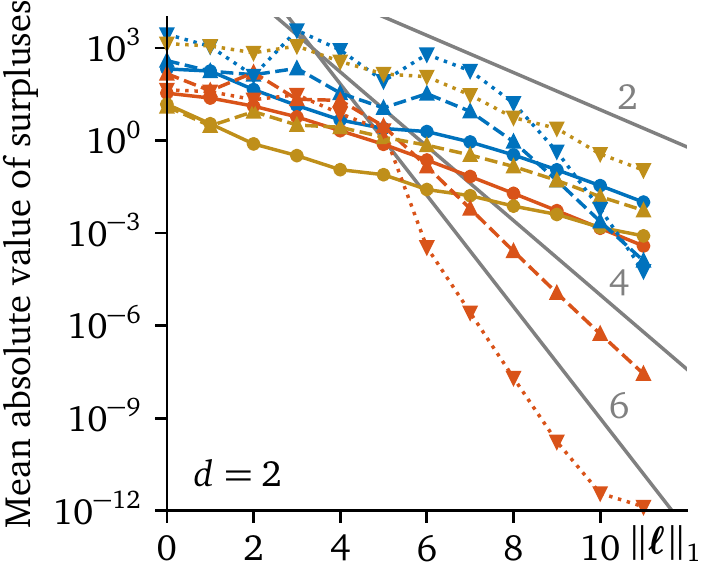}%
  \hfill%
  \includegraphics{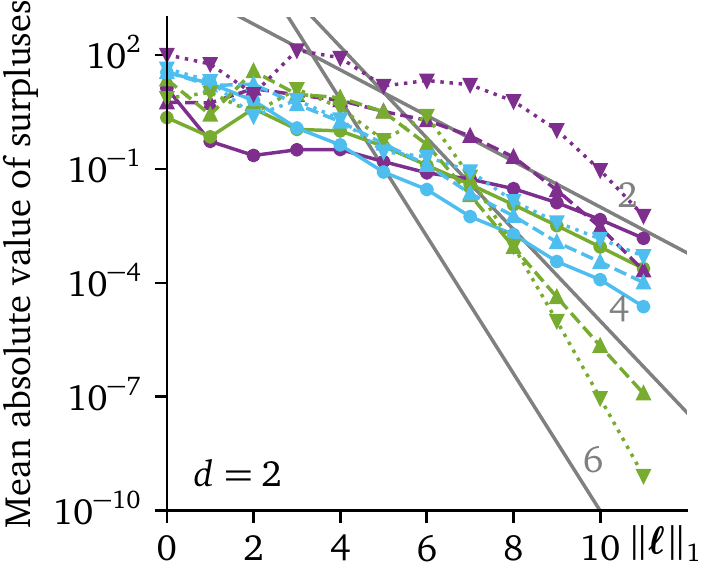}%
  \caption[Decay of surpluses for different test functions]{%
    Mean absolute value of surpluses by level sum $\normone{\*l}$
    \vspace{-0.05em}%
    for different test functions $\objfun$ \emph{(colors)}
    using hierarchical not-a-knot B-splines
    $\bspl[\nak]{\*l,\*i}{p}$ of different degrees $p$
    \emph{(line styles/markers)} on
    the regular sparse grid $\regsgset{n}{d}$ of level $n = 11$.%
  }%
  \label{fig:resultsDecaySurpluses}%
\end{figure}

The decay of surpluses can be seen in \cref{fig:resultsDecaySurpluses},
which shows the mean absolute value of surpluses corresponding to
grid points grouped by their level sum $\normone{\*l}$.
Due to the dependency of coarse-level surpluses on high-level grid points
for $p > 1$,
we have to fix the level $n$ of the regular sparse grid for this analysis.
\Cref{fig:resultsDecaySurpluses} suggests that the
absolute value of the surpluses decays with order $p + 1$ for
B-spline degree $p$, although no theoretical evidence
is known to support this claim.
Higher B-spline degrees seem to imply that
$\abs{\surplus{\*l,\*i}}$ generally increases, if $\normone{\*l}$ is
in the preasymptotic range.

\subsection{Complexity of Hierarchization}
\label{sec:543complexity}

In \cref{chap:40algorithms}, we introduced a number of new
hierarchical spline bases with the aim to reduce the complexity
of algorithms with the key example of hierarchization.
In the following, we study the suitability of the new bases
to achieve this goal \cite{Valentin18Fundamental}.

\paragraph{Complexity of fundamental splines}

\begin{figure}
  \includegraphics{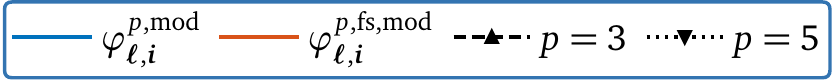}\\[2mm]%
  \subcaptionbox{%
    Computation time%
  }[72mm]{%
    \includegraphics{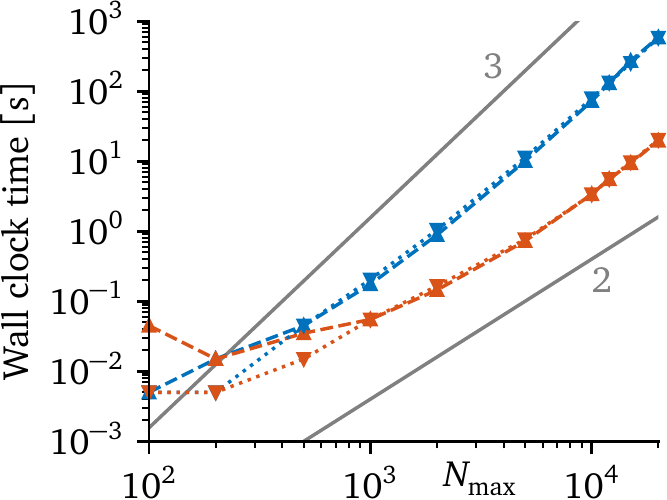}%
  }%
  \hfill%
  \subcaptionbox{%
    Memory consumption%
  }[72mm]{%
    \includegraphics{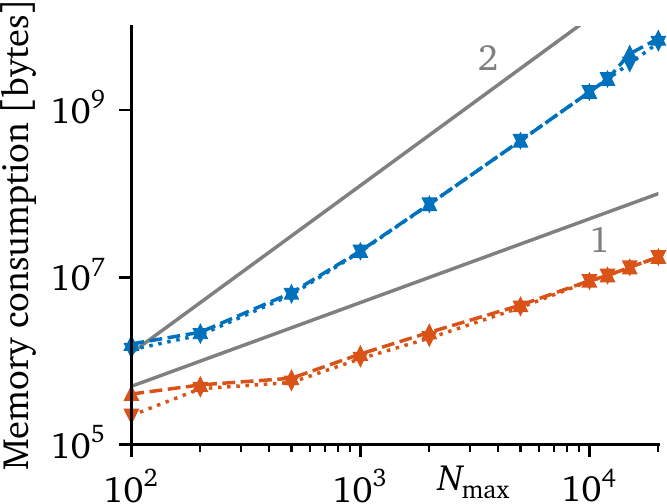}%
  }%
  \caption[Complexity of fundamental splines]{%
    Computation time and memory consumption of hierarchization
    with modified hierarchical B-splines
    $\bspl[\modified]{\*l,\*i}{p}$ \emph{\textcolor{C0}{(blue)}} and
    \vspace*{-0.3em}%
    modified hierarchical fundamental splines
    $\bspl[\fs,\modified]{\*l,\*i}{p}$ \emph{\textcolor{C1}{(red)}}
    of degrees $p = 3$ and $p = 5$
    on the spatially adaptive sparse grids generated by the criterion of
    Novak--Ritter for the optimization of Ack with $d = 4$ and
    $\ngpMax$ grid points.
    Adapted from \cite{Valentin18Fundamental}.%
  }%
  \label{fig:complexityFundamental}%
\end{figure}

\Cref{fig:complexityFundamental} compares the hierarchization complexity of
modified hierarchical B-splines $\bspl[\modified]{\*l,\*i}{p}$ with the new
modified hierarchical fundamental splines $\bspl[\fs,\modified]{\*l,\*i}{p}$
as measured on a laptop with Intel Core i5-4300U.
For the modified hierarchical B-spline basis,
we solve a linear system of size $\ngp \times \ngp$,
for which Gaussian elimination takes
$\landauTheta{\ngp^3}$ time and $\landauTheta{\ngp^2}$ memory for
$\ngp \to \infty$ (where $\ngp$ is the number of sparse grid points and
$d$ is assumed to be constant).
More sophisticated methods to solve linear systems
are not able to significantly
reduce this complexity without any further assumptions on the system matrix
$\intpmat$ (e.g., symmetry, positive definiteness, or bandedness).
As $\ngp$ grows,
the space needed to store an $\ngp \times \ngp$ matrix quickly exceeds
the available memory.

For the modified hierarchical fundamental splines,
we can use the \bfs algorithm presented in \cref{sec:44spatAdaptiveBFS}.
\bfs works in quadratic time $\landauO{\ngp^2}$, but more importantly,
it works in linear space $\landauO{\ngp}$.
Both can be seen very well in \cref{fig:complexityFundamental}:
The computation time drops from cubic to quadratic complexity for fundamental splines
and the consumed memory is reduced from quadratic to linear complexity.

\paragraph{Complexity of weakly fundamental splines}

\begin{figure}
  \includegraphics{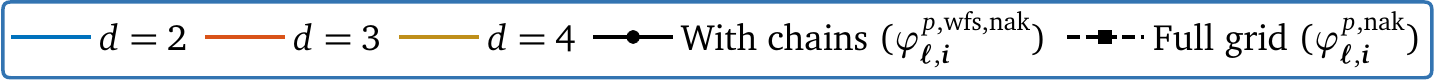}\\[2mm]%
  \subcaptionbox{%
    $\gamma = 0.05$%
  }[48mm]{%
    \includegraphics{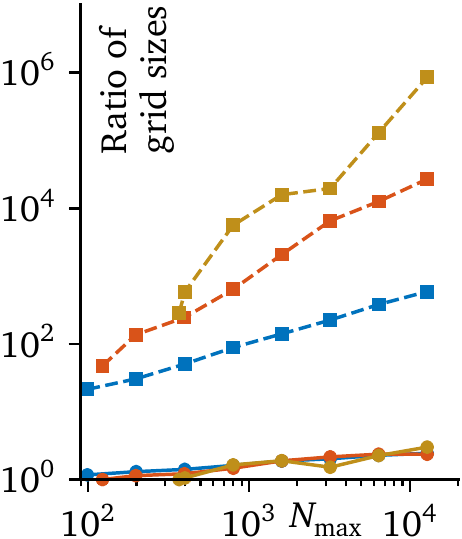}%
  }%
  \hfill%
  \subcaptionbox{%
    $\gamma = 0.15$%
  }[48mm]{%
    \includegraphics{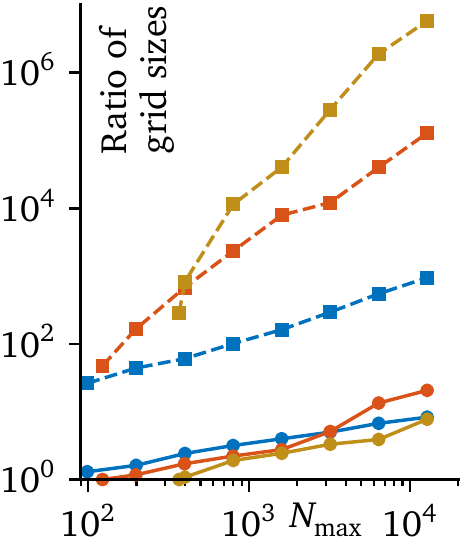}%
  }%
  \hfill%
  \subcaptionbox{%
    $\gamma = 0.25$%
  }[48mm]{%
    \includegraphics{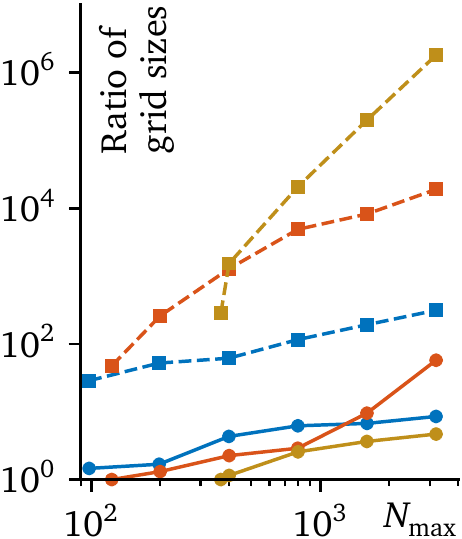}%
  }%
  \caption[Complexity of weakly fundamental splines]{%
    Total number of grid points after inserting all missing chains
    for cubic weakly fundamental not-a-knot splines
    $\bspl[\wfs,\nak]{\*l,\*i}{p}$ \emph{($p = 3$, solid lines),}
    and after inserting all missing full grid points \emph{(dashed).}
    Shown are the ratios of the resulting grid sizes to the
    initial grid sizes before inserting points.
    The initial grids are the spatially adaptive sparse grids generated
    by the criterion of Novak--Ritter for the optimization of Ack with
    different dimensionalities \emph{(colors)},
    different adaptivity parameters $\gamma$ \emph{(left, center, right),}
    and $\ngpMax$ grid points.%
  }%
  \label{fig:complexityWeaklyFundamental}%
\end{figure}

It is not straightforward to include weakly fundamental splines
in \cref{fig:complexityFundamental} as we have to insert missing
chain points to apply the unidirectional principle
(see \cref{sec:45spatAdaptiveUP}).
As this increases the number of necessary evaluations of $\objfun$,
a comparison of computation times with standard B-splines would be skewed.
Instead, we study in \cref{fig:complexityWeaklyFundamental}
the number of grid points that have to be inserted
to ensure the correctness of the unidirectional principle.
As we have seen in \cref{sec:453chains}
(cf.\ \cref{fig:chainInsertionBSpline}),
inserting all chains needed for standard hierarchical B-splines
often results in a full grid, which suffers from the curse of dimensionality.
This can be seen in \cref{fig:complexityWeaklyFundamental}
for hierarchical not-a-knot B-splines (dashed lines).
By inserting all full grid points,
the number of grid points increases by several orders of magnitude:
If the initial grid has $\ngpMax = \num{10000}$ points,
then the grid size increases roughly by the factor $10^2$ for $d = 2$,
$10^4$ for $d = 3$, and $10^6$ for $d = 4$,
resulting in computationally infeasible
grids with $10^6$, $10^8$, and $10^{10}$ points, respectively.
If we instead only insert the missing chain points needed for the
hierarchical weakly fundamental not-a-knot basis
(solid lines, cf.\ \cref{fig:chainInsertionWeaklyFundamentalSpline}),
then the number of grid points increases only slightly.
For grids that have a low adaptivity (which correspond
to low adaptivity parameters $\gamma$ in the Novak--Ritter criterion,
see \cref{sec:521novakRitter}), the grid size only increases by the
factor of two.
For highly-adaptive grids
(corresponding to large $\gamma$),
the number of necessary chain grid points increases significantly.

\subsection{Optimality Gap}
\label{sec:542optimization}

\paragraph{Optimality gaps and displacements}

With the method described in \cref{sec:52method},
we find approximations $\xoptappr$ of the
global minimum $\xopt$ of some objective function $\objfun$
using optimization of a B-spline surrogate $\sgintp$ of $\objfun$
on sparse grids.
Obviously, the more accurate the sparse grid surrogate is,
the better the approximation $\xoptappr$ will be.
In the following plots,
we show the optimality gaps $\objfun(\xoptappr) - \objfun(\xopt)$
in terms of function values.%
\footnote{%
  In order to calculate the optimality gap,
  it is crucial to determine $\objfun(\xopt)$ as exact as possible.
  Otherwise, the optimality gap might either not converge to zero
  or it might even become negative.%
}
The results are sensitive to even small displacements
of the objective function, i.e.,
the results may change for the
function $\*x \mapsto \objfun(\*x - \*a)$
instead of $\*x \mapsto \objfun(\*x)$ for $\*x \in \clint{\*0, \*1}$
and some small $\*a \in \real^d$.%
\footnote{%
  By using the formulas in \cref{chap:a20testProblems},
  all test functions $\objfun$ in \cref{sec:53testProblems}
  can be extended such that they can be evaluated at $\*x - \*a$
  for all $\*x \in \clint{\*0, \*1}$, if $\*a \in \real^d$ is small enough.
  Note that we set $a_t$ to zero if a non-zero displacement in
  the $t$-th component would change the location of the global minimum.%
}
Therefore, the optimization for each of the
data points for \cref{%
  fig:resultsOptimizationUnconstrainedTestFunctions,%
  fig:resultsOptimizationConstrainedTestFunctions%
}
was repeated five times with replacements $\*a$
whose entries $a_t$ were independent and identically distributed Gaussian
pseudo-random numbers with zero mean and a standard deviation of $0.01$.
The optimality gaps shown in the figures of this section were computed
as the mean of the five runs to increase confidence in the results.

\paragraph{Unconstrained optimization}

\begin{figure}
  \includegraphics{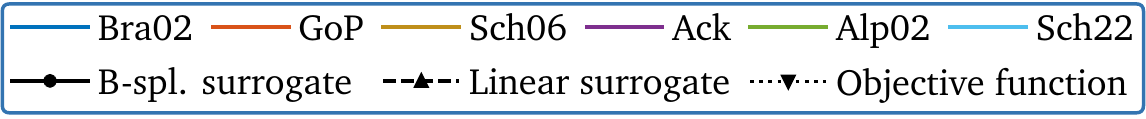}\\[2mm]%
  \includegraphics{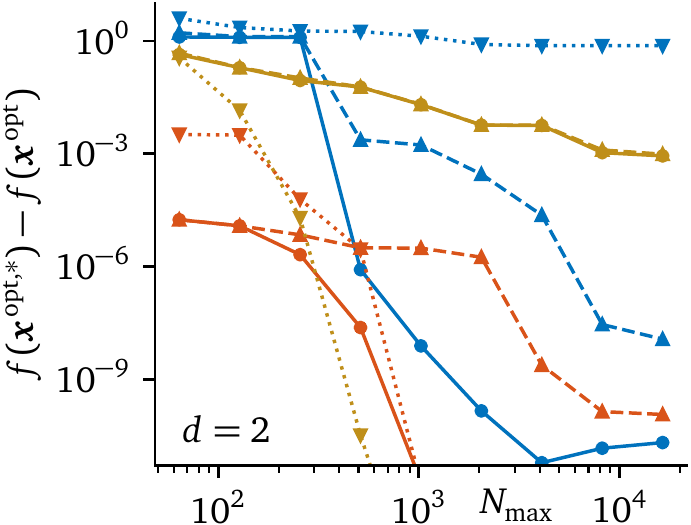}%
  \hfill%
  \includegraphics{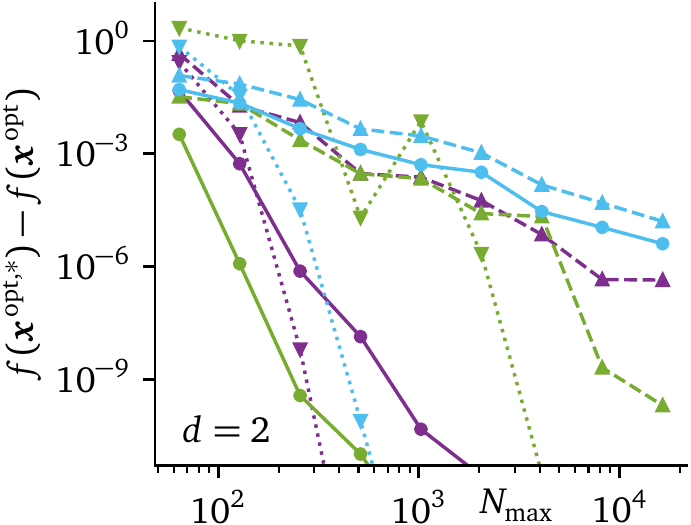}%
  \\[2mm]%
  \includegraphics{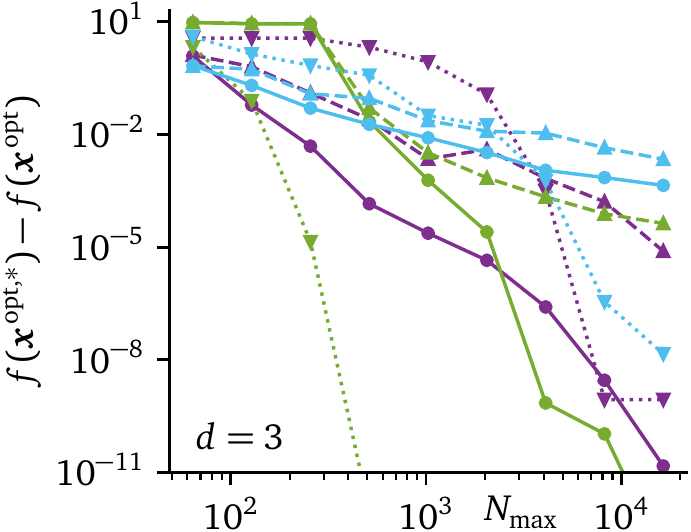}%
  \hfill%
  \includegraphics{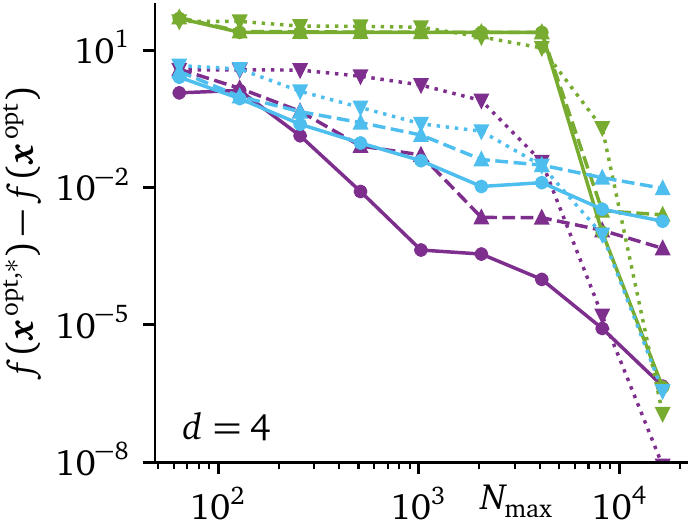}%
  \\[2mm]%
  \includegraphics{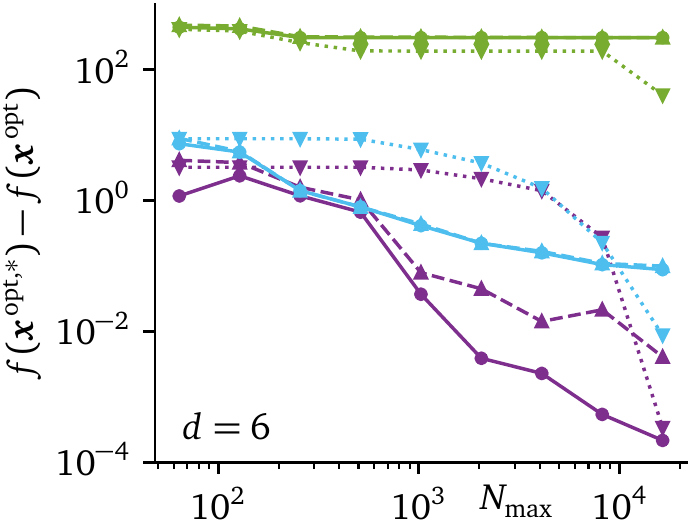}%
  \hfill%
  \includegraphics{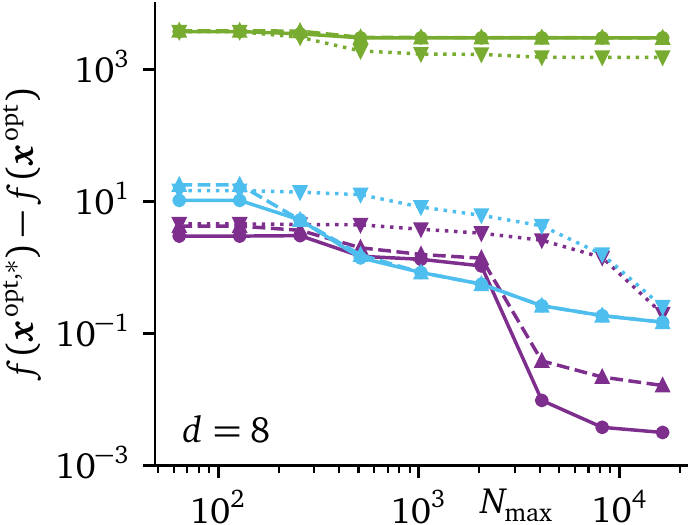}%
  \caption[Optimality gaps for different objective functions (unconstrained)]{%
    Optimality gaps $\objfun(\xoptappr) - \objfun(\xopt)$ between
    the function value at the approximated optimum $\xoptappr$ and
    the minimal function value at the actual optimum $\xopt$
    over the number $\ngpMax$ of objective function evaluations
    for different unconstrained
    objective functions $\objfun$ \emph{(colors).}
    Shown are the optimization results of the B-spline surrogate
    \emph{(solid lines),}
    the optimization results of the piecewise linear surrogate
    \emph{(dashed),} and
    the optimization results of the actual objective function
    \emph{(dotted)} as described in \cref{sec:52method}.%
  }%
  \label{fig:resultsOptimizationUnconstrainedTestFunctions}%
\end{figure}

\Cref{fig:resultsOptimizationUnconstrainedTestFunctions}
shows the optimality gaps for different test functions $\objfun$
over the number $\ngpMax$ of allowed evaluations of $\objfun$.
For the continuously differentiable functions
Bra02, GoP\punctfix{,} Ack, and Alp02
in $d = 2$ dimensions (top row),
the optimization of the corresponding cubic B-spline surrogates (solid lines)
performs significantly better than using piecewise linear basis functions
(dashed lines).
The reason is two-fold:
First, by using higher-order basis functions, the surrogates are more accurate
in general as seen in the discussion of the interpolation error in
\cref{sec:541interpolation}.
Second, the availability of surrogate gradients accelerates the
convergence of the employed optimization methods.
For some test functions, B-splines give better results than even
the direct optimization of the objective function (dotted lines).

For the test functions Sch06 and Sch22 with discontinuous derivatives,
the advantage of higher-order B-splines is not as evident (Sch22) or
does not even exist (Sch06).
However, in low dimensions, i.e., $d \le 4$, B-splines
achieve a slight advantage compared to the piecewise linear basis
for the Sch22 function.
In higher dimensionalities, i.e., $d \ge 6$ (bottom row),
convergence visibly slows down for all methods shown in
\cref{fig:resultsOptimizationUnconstrainedTestFunctions},
although for some objective functions, B-splines are still able
to perform better than the comparison methods
(most notably for the Ack function).

\paragraph{Constrained optimization}

\begin{figure}
  \includegraphics{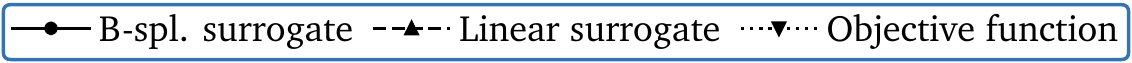}\\[2mm]%
  \subcaptionbox{%
    G08%
  }[72mm]{%
    \includegraphics{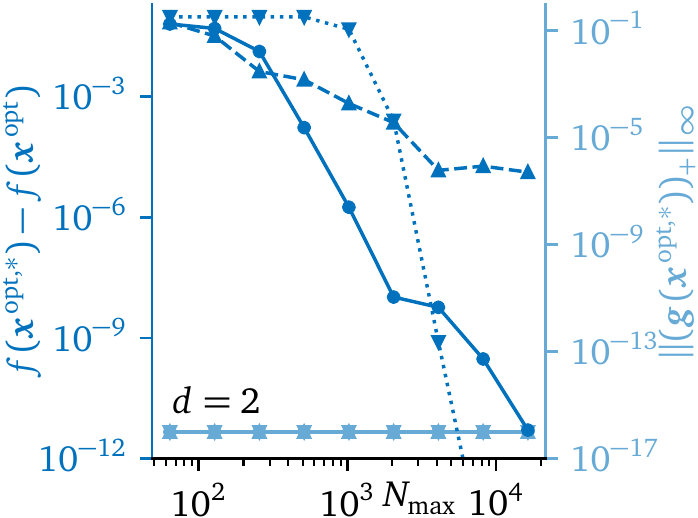}%
  }%
  \hfill%
  \subcaptionbox{%
    G04Sq%
  }[72mm]{%
    \includegraphics{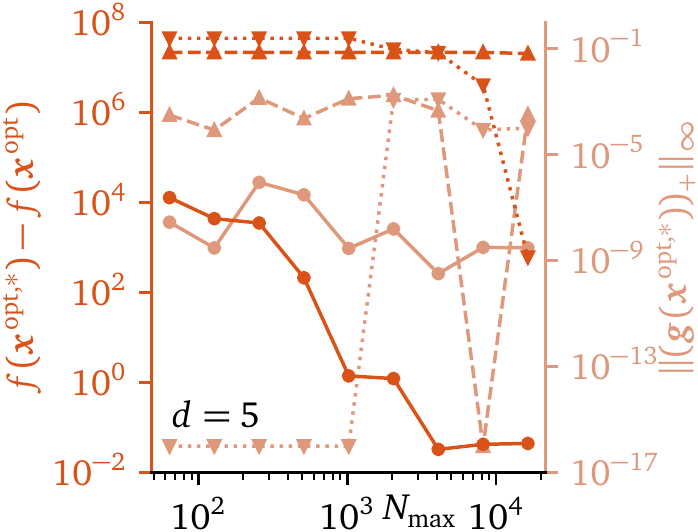}%
  }%
  \caption[Optimality gaps for different objective functions (constrained)]{%
    Optimality gaps $\objfun(\xoptappr) - \objfun(\xopt)$ between
    the function value at the approximated optimum $\xoptappr$ and
    the minimal function value at the actual optimum $\xopt$
    over the number $\ngpMax$ of objective function evaluations
    for different constrained optimization problems
    \emph{(dark, left vertical axes).}
    In addition to the lines of
    \cref{fig:resultsOptimizationUnconstrainedTestFunctions},
    the constraint violation
    $\norm[\infty]{\nonnegpart{\ineqconfun(\xoptappr)}}$
    at the approximated optimum $\xoptappr$ is plotted
    \emph{(light, right vertical axes).}%
  }%
  \label{fig:resultsOptimizationConstrainedTestFunctions}%
\end{figure}

\Cref{fig:resultsOptimizationConstrainedTestFunctions}
shows the result for the two constrained optimization problems.
The objective function value $\objfun(\xoptappr)$
at the approximated optimum $\xoptappr$ should not only
be as small as possible, but $\xoptappr$ should also be feasible, i.e.,
$\ineqconfun(\xoptappr) \le \*0$.
Hence, we also plot the maximal violation
$\norm[\infty]{\nonnegpart{\ineqconfun(\xoptappr)}}$
of the constraints in the respective optimal points $\xoptappr$.

For the bivariate G08 problem, the hierarchical B-splines surrogates
perform better than the direct gradient-free optimization of the problem
for $\ngpMax \le 3500$ objective function evaluations
and better than the piecewise linear surrogate for $\ngpMax \ge 300$
objective function evaluations.
All calculated points are feasible.%
\footnote{%
  For plotting reasons,
  \cref{fig:resultsOptimizationConstrainedTestFunctions} shows
  $\max(\norm[\infty]{\nonnegpart{\ineqconfun(\xoptappr)}}, 10^{-16})$
  instead of the true constraint violation.%
}

The range of the objective function of the five-dimensional G04Sq problem
is larger than the range of G08.
This results in generally higher optimality gaps
$\objfun(\xoptappr) - \objfun(\xopt)$ as we do not normalize
with respect to the range.
B-splines achieve good approximations $\xoptappr$
of $\xopt$ already for $\ngpMax = 1000$ with an optimality gap of
around one.
Both comparison methods show optimality gaps that are
seven orders of magnitude higher.

Additionally, the corresponding values of constraint violation
are between $10^{-10}$ and $10^{-6}$, i.e.,
the constraints are numerically met.
In contrast, the optimizers struggle more for
the comparison methods (optimization of the linear surrogate and
of the objective function) to meet the constraints,
as the values of the constraint violation partly exceed $10^{-3}$.
The availability of gradients seems to allow the constrained optimization
methods to better enforce the feasibility of the
resulting points $\xoptappr$.

Note that while the results look already promising for the B-spline surrogate
method, these results could still be improved upon.
The Novak--Ritter criterion used to generate the
spatially adaptive sparse grids does not take the constraints into account.
Consequently, many sparse grid points are created outside the feasible domain.
By modifying the criterion to prefer points that are in a neighborhood of
the feasible domain, the quality of the interpolant
close to potential optima should increase.

\section{Example Application: Fuzzy Extension Principle}
\label{sec:55fuzzy}

\minitoc[-1mm]{84mm}{4}

\noindent
To conclude this chapter, we consider the fuzzy extension principle
as an example application of optimization of B-spline sparse grid surrogates.

\paragraph{Aleatoric and epistemic uncertainties}

Classical uncertainty quantification (UQ) distinguishes between
aleatoric and epistemic uncertainties \cite{Walz16Fuzzy}.
Aleatoric uncertainties result from the variability of inputs or
model components and from the ``intrinsic randomness''
of quantities.
They are best described by probability theory, giving exact probabilities.
Epistemic uncertainties arise from subjectivity,
simplifying modeling assumptions, and incomplete knowledge.
These uncertainties are better captured by fuzzy theory,
which is more imprecise than the ``exact'' stochastic assumptions
of probabilities \cite{Walz16Fuzzy}.

\paragraph{Uncertainty quantification with fuzzy uncertainties}

In uncertainty quantification, the key question is as follows:
Given a model and uncertain input parameters for the model,
how uncertain is the model output?
While there are many approaches available
for probabilistic uncertainties,
it is not straightforward to solve this task
for fuzzy uncertainties.
Fortunately, Zadeh proposed in 1975
the \term{fuzzy extension principle} \cite{Zadeh75Concept},
which addresses this very question.

\paragraph{Sparse grids and B-splines for fuzzy uncertainties}

As we explain in this section,
the fuzzy extension principle requires the solution of numerous
optimization problems that involve the original objective function
$\objfun$.
This predestines the replacement of $\objfun$ with sparse grid surrogates,
as explained in the beginning of the chapter.
Previous work by Klimke \cite{Klimke06Uncertainty} already
studied this approach for piecewise linear functions on uniform sparse grids
and for global polynomials on sparse Clenshaw--Curtis grids.
We assess the suitability of interpolation with higher-order
hierarchical B-splines on sparse grids for the fuzzy extension principle.
It should be mentioned that there is also work
directly incorporating (non-hierarchical) B-splines
into the framework of fuzzy theory for modeling uncertain surfaces
\multicite{Anile00Modeling,Zakaria14Fuzzy}.

\subsection{Fuzzy Sets and Fuzzy Intervals}
\label{sec:551fuzzySets}

\begin{figure}
  \subcaptionbox{%
    Non-convex fuzzy set \emph{\textcolor{C0}{(blue)}}
    and $\alpha$-cut \emph{\textcolor{C1}{(red)}.}%
  }[46mm]{%
    \includegraphics{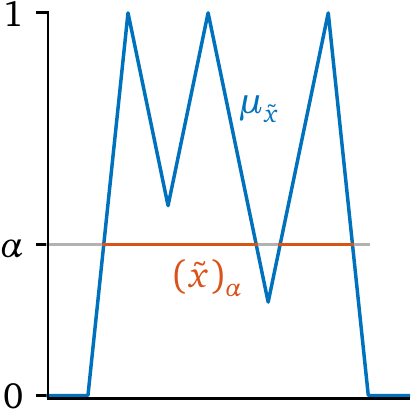}%
  }%
  \hfill%
  \subcaptionbox{%
    Fuzzy interval.%
  }[43mm]{%
    \includegraphics{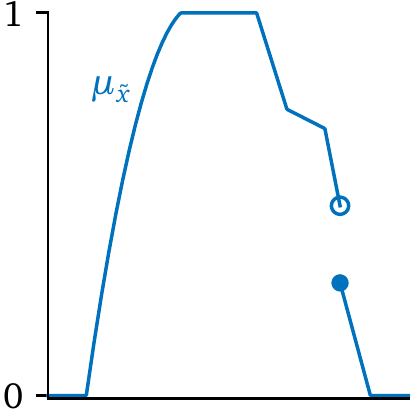}%
  }%
  \hfill%
  \subcaptionbox{%
    Common types of fuzzy numbers and intervals (see text).%
  }[56mm]{%
    \includegraphics{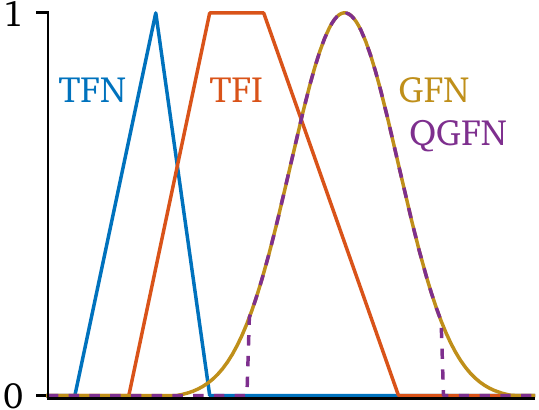}%
  }%
  \caption[%
    Examples of fuzzy sets and $\alpha$-cuts%
  ]{%
    Examples of membership functions of fuzzy sets and $\alpha$-cuts.%
  }%
  \label{fig:fuzzySet}%
\end{figure}

In the following, we repeat very briefly the necessary
definitions of basic fuzzy theory.
Examples for the definitions are shown in \cref{fig:fuzzySet}.
A more in-depth introduction can be found in
\multicite{Hanss05Applied,Klimke06Uncertainty,Walz16Fuzzy}.

\paragraph{Fuzzy sets}

Let $X \subset \real$ be a closed interval on the real line
and $\memfun{x}\colon X \to \clint{0, 1}$ be a function.
We call the graph $\fuzzy{x} \ceq \{(x, \memfun{x}(x) \mid x \in X\}$
of $\memfun{x}$ a \term{fuzzy set} with
\term{membership function} $\memfun{x}$.
Fuzzy sets generalize ordinary subsets of $X$,
which can be obtained by requiring $\memfun{x}(X) \subset \{0, 1\}$.
In this case, the fuzzy set is called \term{crisp} and
$\fuzzy{x}$ can be identified with the ordinary set
$\{x \in X \mid \memfun{x}(x) = 1\}$.
A fuzzy set $\fuzzy{x}$ is \term{normalized}
if $\max_{x \in X} \memfun{x}(x) = 1$.
A \term{convex} fuzzy set $\fuzzy{x}$ satisfies
$\min(\memfun{x}(a), \memfun{x}(c)) \le \memfun{x}(b)$ for all $a, b, c \in X$
with $a \le b \le c$.

\paragraph{Fuzzy intervals and $\alpha$-cuts}

A convex and normalized fuzzy set $\fuzzy{x}$ with
piecewise continuous membership function $\memfun{x}$ is called
\term{fuzzy interval.}
If $\{x \in X \mid \memfun{x}(x) = 1\} = \{a\}$ for some $a \in X$,
then the fuzzy interval $\fuzzy{x}$ is called \term{fuzzy number.}

For $\alpha \in \clint{0, 1}$, the $\alpha$-cut of $\fuzzy{x}$ is
defined as $\acut{x}{\alpha} \ceq \{x \in X \mid \memfun{x}(x) \ge \alpha\}$
for $\alpha > 0$ and $\acut{x}{0} \ceq \supp \memfun{x}$ for $\alpha = 0$.
The $\alpha$-cuts of fuzzy intervals $\fuzzy{x}$ are always
nested closed intervals, i.e.,
$\acut{x}{\alpha} = [a, b]$ for some $a \le b$ and
$\acut{x}{\alpha_1} \supset \acut{x}{\alpha_2}$ for $\alpha_1 \le \alpha_2$.

\paragraph{Common types of fuzzy numbers and intervals}

There are various types of fuzzy numbers and intervals
\cite{Klimke06Uncertainty}.
Most common are
\term{triangular fuzzy numbers} (TFNs, i.e., linear B-splines),
\term{trapezoidal fuzzy intervals}
(TFIs, where a plateau of height one is inserted at the peak, i.e.,
sums of two neighboring linear B-splines), and
\term{Gaussian fuzzy numbers} (GFNs) with membership function
$\memfun{x}(x) = \exp(-\frac{(x - \mu)^2}{(2\sigma)^2})$.
As the support of Gaussian fuzzy numbers is unbounded,
\term{quasi-Gaussian fuzzy numbers} (QGFNs) truncate the support
to a fixed multiple of the standard deviation $\sigma$
\cite{Klimke06Uncertainty}.
However, it would be more natural to directly employ B-splines of
degree $p > 1$ (normalized adequately), since they generalize
triangular fuzzy numbers and their limit with respect to $p$
is a Gaussian fuzzy number.

\subsection{Fuzzy Extension Principle}
\label{sec:552fuzzyExtensionPrinciple}

Let $\objfun\colon \clint{\*0, \*1} \to \real$ be an objective function,
whose values $y = \objfun(\*x)$ represent the results of the
simulation of a model with input parameters $(x_1, \dotsc, x_d) = \*x$.
If the input parameters are uncertain and
given as fuzzy sets $\fuzzy{x}_1, \dotsc, \fuzzy{x}_d$,
what is the resulting uncertain outcome
``$\fuzzy{y} \ceq \objfun(\fuzzy{x}_1, \dotsc, \fuzzy{x}_d)$''?
Note that there is no definite answer to this question,
as ``$\objfun(\fuzzy{x}_1, \dotsc, \fuzzy{x}_d)$'' is not well-defined.
The fuzzy extension principle, suggested by Zadeh \cite{Zadeh75Concept},
provides one possible definition.

\paragraph{Alternative fuzzy extension principle}

We use an alternative formulation of the fuzzy extension principle,
which is stated in \cite{Klimke06Uncertainty}.
The original formulation is computationally more complex,
as it requires the solution of equality-constrained optimization problems
and one needs to know the range of $\objfun$, which might not be given.
The two formulations are equivalent,
if $\fuzzy{x}_1, \dotsc, \fuzzy{x}_d$ are (compactly supported)
fuzzy intervals and $\objfun$ is continuous \cite{Buckley90Using},
which we assume in the following.

\begin{figure}
  \includegraphics{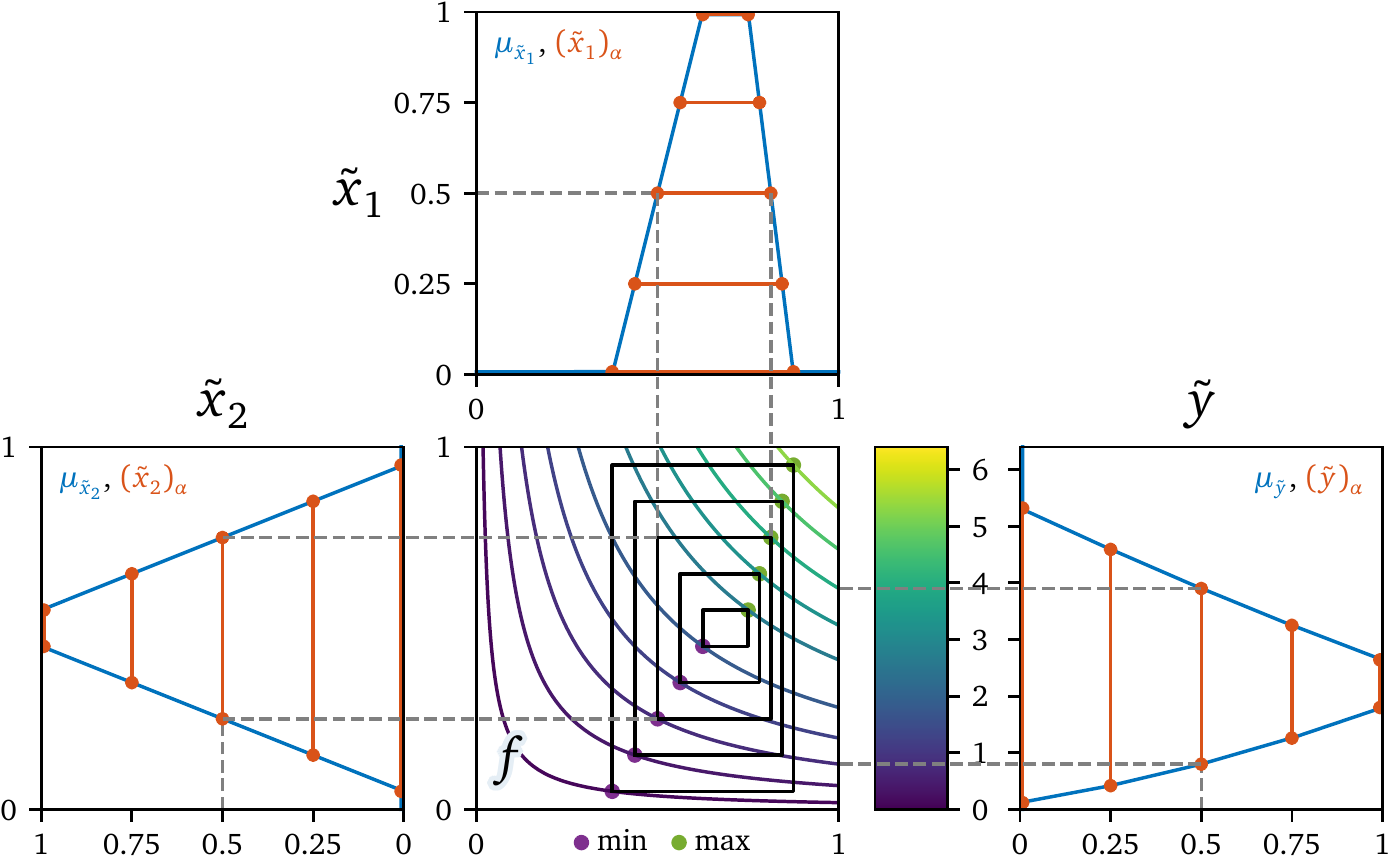}%
  \caption[%
    Alternative fuzzy extension principle%
  ]{%
    Example of the application of the
    alternative fuzzy extension principle to the bivariate objective function
    $\objfun(\*x) = 6.4 x_1 x_2$ \emph{(bottom)}
    and triangular fuzzy input intervals
    $\fuzzy{x}_1$ and $\fuzzy{x}_2$ \emph{(top and left)}
    to obtain the fuzzy output interval $\fuzzy{y}$ \emph{(right).}
    Adapted from \cite{Klimke06Uncertainty}.%
  }%
  \label{fig:fuzzyExtensionPrinciple}%
\end{figure}

The alternative fuzzy extension principle defines
``$\fuzzy{y} = \objfun(\fuzzy{x}_1, \dotsc, \fuzzy{x}_d)$'' as the fuzzy set
$\fuzzy{y}$ with
\begin{subequations}
  \label{eq:alternativeFuzzyExtensionPrinciple}
  \begin{alignat}{2}
    \memfun{y}(y)
    &\ceq \sup\{\alpha \in \clint{0, 1} \mid y \in \acut{y}{\alpha}\},\quad
    &&y \in \real,\\
    \acut{y}{\alpha}
    &\ceq \bracket*{
      \min_{\*x \in \Omega_\alpha} \objfun(\*x),\;
      \max_{\*x \in \Omega_\alpha} \objfun(\*x)
    },\quad
    &&\alpha \in \clint{0, 1},\\
    \Omega_\alpha
    &\ceq \acut[1]{x}{\alpha} \times \dotsb \times \acut[d]{x}{\alpha},\quad
    &&\alpha \in \clint{0, 1}.
  \end{alignat}
\end{subequations}
This definition is visualized in \cref{fig:fuzzyExtensionPrinciple}.
The first equation defines $\fuzzy{y}$ via its $\alpha$-cuts,
which are given in the second equation as the closed interval
between the minimal and the maximal value of $\objfun$ on some
hyper-rectangular domain $\Omega_\alpha$.
The third equation specifies this domain $\Omega_\alpha$ as the
Cartesian product of the univariate $\alpha$-cuts.
Hence, we only have to solve box-constrained optimization problems,
as opposed to the general equality-constrained problems
in the original formulation of the fuzzy extension principle.

\paragraph{Implementation}

The implementation of the alternative fuzzy extension principle
is straightforward and shown in \cref{alg:alternativeFuzzyExtensionPrinciple}
\cite{Klimke06Uncertainty}.
The range $\clint{0, 1}$ of $\alpha$
is discretized into $m + 1$ uniformly spaced values $\alpha_j$
(where $m \in \nat$).
For each of these values $\alpha_j$, we compute the corresponding
$\alpha_j$-cut of $\fuzzy{y}$ by solving the two box-constrained
optimization problems of \cref{eq:alternativeFuzzyExtensionPrinciple}.
The fuzzy output interval $\fuzzy{y}$ can then be approximated by
interpolating the interval bounds of the $\alpha_j$-cuts of $\fuzzy{y}$.

\begin{algorithm}
  \begin{algorithmic}[1]
    \Function{$\fuzzy{y} = \texttt{alternativeFuzzyExtensionPrinciple}$}{%
      $m$, $\fuzzy{x}_1$, \dots, $\fuzzy{x}_d$%
    }
      \For{$j = 0, \dotsc, m$}
        \State{$\alpha_j \gets j/m$}
        \ForOneLine{$t = 1, \dotsc, d$}{Compute $\acut[t]{x}{\alpha_j} = \clint{a_{j,t}, b_{j,t}}$}
        \State{%
          $\Omega_{\alpha_j} \gets \acut[1]{x}{\alpha_j} \times \dotsb \times
          \acut[d]{x}{\alpha_j} = \clint{\*a_j, \*b_j}$%
        }
        \State{%
          Solve $\min_{\*x \in \Omega_{\alpha_j}} \objfun(\*x)$ and
          $\max_{\*x \in \Omega_{\alpha_j}} \objfun(\*x)$%
        }
        \State{%
          $\acut{y}{\alpha_j} = [c_j, d_j] \gets \clint{
            \min_{\*x \in \Omega_{\alpha_j}} \objfun(\*x),
            \max_{\*x \in \Omega_{\alpha_j}} \objfun(\*x)
          }$
        }
      \EndFor{}\vspace{-2mm}
      \State{%
        $D \gets
        \{(c_j, \alpha_j) \mid j = 0,\, 1,\, \dotsc,\, m\} \cup
        \{(d_j, \alpha_j) \mid j = m,\, m - 1,\, \dotsc,\, 0\}$%
      }
      \State{%
        $\memfun{y} \gets \text{Piecewise linear interpolant of $D$}$%
      }%
      \Comment{extend to $X$ by zero}%
    \EndFunction{}
  \end{algorithmic}
  \caption[Alternative fuzzy extension principle]{%
    Alternative fuzzy extension principle.
    Inputs are the number of $\alpha$ segments to use as discretization and
    the $d$ fuzzy intervals $\fuzzy{x}_1, \dotsc, \fuzzy{x}_d$
    (we have to be able to determine $\alpha$-cuts
    of these fuzzy input intervals).
    The output is an approximation to the output $\fuzzy{y}$
    of the alternative fuzzy extension principle
    (given by an approximation of its membership function $\memfun{y}$).%
  }%
  \label{alg:alternativeFuzzyExtensionPrinciple}%
\end{algorithm}

\subsection{Using B-Splines on Sparse Grids to Propagate Fuzzy Uncertainties}
\label{sec:553fuzzyBSplines}

Following Klimke's approach \cite{Klimke06Uncertainty},
we replace the objective function $\objfun$ in
\cref{alg:alternativeFuzzyExtensionPrinciple}
with a sparse grid surrogate $\sgintp$.
The solution of the optimization problems
$\min_{\*x \in \Omega_{\alpha_j}} \sgintp(\*x)$ and
$\max_{\*x \in \Omega_{\alpha_j}} \sgintp(\*x)$ with respect to the
surrogate $\sgintp$ instead of the true objective function $\objfun$
takes significantly less time, if evaluations of the objective function
are expensive.

However, Klimke used piecewise linear functions as the hierarchical basis on
uniform sparse grids and global polynomials on sparse Clenshaw--Curtis grids.
The drawbacks of each of the bases are evident:
First, piecewise linear surrogates are not continuously differentiable and
can thus not be optimized well with gradient-based optimization methods.
Second, global polynomials are only suitable for
Clenshaw--Curtis grids (Chebyshev-distributed points)
due to Runge's phenomenon,
unnecessarily restricting the choice of grid points.
Hierarchical B-splines of degree $p$ are $(p - 1)$ times
continuously differentiable and defined for arbitrary point
distributions, eliminating both drawbacks simultaneously.

\paragraph{Methodology}

Given a sparse grid $\sgset$, which may be regular or spatially adaptive,
we compute three solutions of the alternative fuzzy extension principle
as follows:

\begin{itemize}
  \item
  First,
  we replace $\objfun$ in \cref{alg:alternativeFuzzyExtensionPrinciple}
  with the sparse grid interpolant $\sgintp[p]$
  on $\sgset$ using modified hierarchical not-a-knot B-splines
  $\bspl[\nak,\modified]{l,i}{p}$ of cubic degree ($p = 3$).
  For solving the optimization problems over $\sgintp[p]$ in
  \cref{alg:alternativeFuzzyExtensionPrinciple},
  we use the globalized version of the method of gradient descent
  as described in \cref{sec:522method} using 100 initial points.
  The resulting fuzzy output interval is denoted by $\fuzzy[\sparse,p]{y}$.
  
  \item
  Second,
  we replace $\objfun$ in \cref{alg:alternativeFuzzyExtensionPrinciple}
  with the sparse grid interpolant $\sgintp[1]$
  on $\sgset$ using modified piecewise linear basis functions.
  For solving the optimization problems over $\sgintp[1]$ in
  \cref{alg:alternativeFuzzyExtensionPrinciple},
  we use a multi-start version of the Nelder--Mead method
  as described in \cref{sec:522method}
  using 100 initial simplices%
  \footnote{%
    The Nelder--Mead method does not require an initial point,
    but an initial simplex.
    The method is a hybrid between global and local optimization.
    If the initial simplex is chosen badly, Nelder--Mead may get stuck
    in local minima.
    Hence, we restart the algorithm for different initial simplices.%
  }.
  The resulting fuzzy output interval corresponds to Klimke's method and
  is denoted by $\fuzzy[\sparse,1]{y}$.
  
  \item
  Third,
  for comparison, we solve
  \cref{alg:alternativeFuzzyExtensionPrinciple} for the
  actual objective function $\objfun$.
  For solving the optimization problems over $\objfun$,
  we use a multi-start version of the Nelder--Mead method
  as described in \cref{sec:522method}
  using 1000 initial simplices
  and \num{2000000} allowed evaluations of $\objfun$.
  The resulting fuzzy output interval is denoted by $\fuzzy[\reference]{y}$
  \term{(reference solution).}
\end{itemize}

\noindent
In the following, we fix the number of $\alpha$ segments
in \cref{alg:alternativeFuzzyExtensionPrinciple} as $m = 100$.
As fuzzy input intervals $\fuzzy{x}_t$, $t = 1, \dotsc, d$, we use
the trapezoidal fuzzy interval with $0$-cut $\clint{0.125, 0.625}$
and $1$-cut $\clint{0.25, 0.375}$ if $t$ is odd and
the quasi-Gaussian fuzzy number with mean $0.5$, standard deviation $0.125$,
and $0$-cut $\clint{0.125, 0.875}$ if $t$ is even.

\paragraph{Convergence of fuzzy intervals on regular sparse grids}

\begin{SCfigure}
  \includegraphics{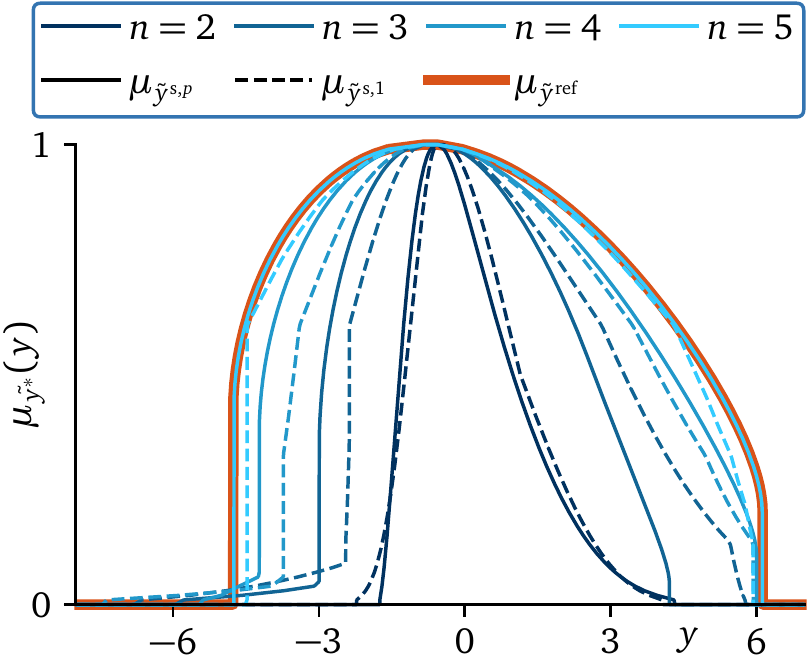}%
  \caption[Convergence of fuzzy output intervals]{%
    Convergence of the membership functions of the fuzzy output intervals
    $\fuzzy[\sparse,p]{y}$
    (\emph{solid lines,} modified hierarchical cubic not-a-knot B-splines)
    and $\fuzzy[\sparse,1]{y}$
    (\emph{dashed,} modified hierarchical hat functions)
    to the reference solution $\fuzzy[\reference]{y}$
    \emph{\textcolor{C1}{(red)}} for the bivariate Alp02 function using
    regular sparse grids of level $n = 2, \dotsc, 5$.%
  }%
  \label{fig:resultsFuzzyPropagation}%
\end{SCfigure}

As an example,
\cref{fig:resultsFuzzyPropagation} shows the convergence of the
fuzzy output intervals $\fuzzy[\sparse,p]{y}$ and $\fuzzy[\sparse,1]{y}$
obtained by the interpolation of the
bivariate Alp02 function on regular sparse grids $\sgset = \regsgset{n}{d}$
to the reference solution $\fuzzy[\reference]{y}$.
Already for $n = 4$, the B-spline approximation is better than the
piecewise linear approximation.
For $n = 5$, no difference is visible anymore between $\fuzzy[\sparse,p]{y}$
and $\fuzzy[\reference]{y}$, while $\fuzzy[\sparse,1]{y}$ still clearly
deviates from $\fuzzy[\reference]{y}$.

\begin{figure}
  \includegraphics{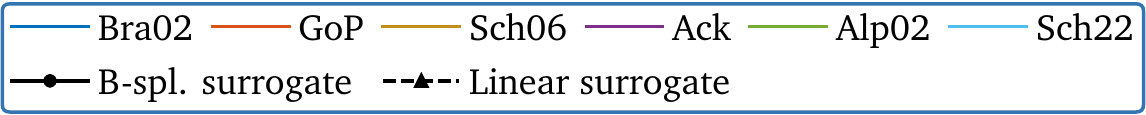}\\[2mm]%
  \includegraphics{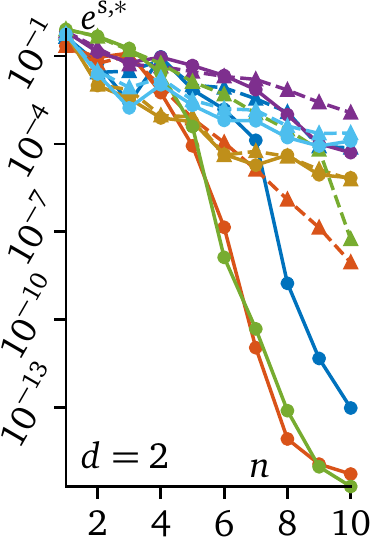}%
  \hfill%
  \includegraphics{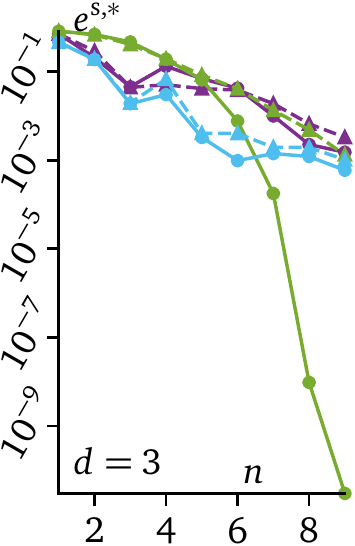}%
  \hfill%
  \includegraphics{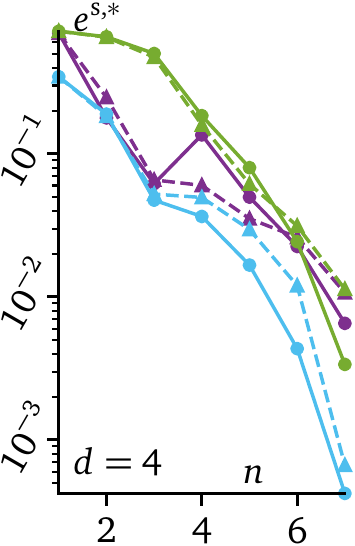}%
  \hfill%
  \includegraphics{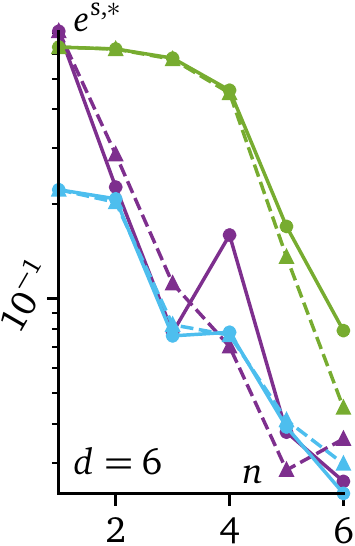}%
  \caption[Fuzzy errors for regular sparse grids]{%
    Fuzzy errors
    $e^{\sparse,\ast}
    \ceq \normLtwo{\memfun[\reference]{y} - \memfun[\sparse,\ast]{y}}/
    \normLtwo{\memfun[\reference]{y}}$
    for regular sparse grids $\regsgset{n}{d}$
    and different objective functions $\objfun$ \emph{(colors)}
    over the level $n$ of the sparse grid.%
  }%
  \label{fig:resultsFuzzyRegular}%
\end{figure}

In \cref{fig:resultsFuzzyRegular}, we study the convergence of the
relative $\Ltwo$ errors
\begin{equation}
  e^{\sparse,\ast}
  \ceq \frac{
    \normLtwo{\memfun[\reference]{y} - \memfun[\sparse,\ast]{y}}
  }{
    \normLtwo{\memfun[\reference]{y}}
  },\quad
  \ast \in \{1, p\},
\end{equation}
of the membership functions (``fuzzy errors'').
The Alp02 ($d = 2$) errors
that correspond to \cref{fig:resultsFuzzyPropagation}
are shown in green in the left-most plot of \cref{fig:resultsFuzzyRegular}.
For the Bra02, GoP\punctfix{,} and Alp02 functions,
B-spline surrogates achieve
dramatic improvements over the hat function surrogates in the bivariate case.
For the bivariate Ack function, B-splines yield an error that
is still an order of magnitude smaller than the error of hat functions.
Just little or even no improvement can be seen
for the functions Sch06 and Sch22 with discontinuous derivatives or
higher dimensionalities $d \ge 4$.

\paragraph{Fuzzy Novak--Ritter method}

We want to employ spatial adaptivity to improve the results
of regular sparse grids.
To this end, we modify the Novak--Ritter criterion
to create a grid generation method that is tailored for the
fuzzy extension principle, resulting in \cref{alg:fuzzyNovakRitterMethod}.
Its main idea is to generate more points near the optima
of the fuzzy extension principle
(\cref{alg:alternativeFuzzyExtensionPrinciple}) than in other regions
of $\clint{\*0, \*1}$.
Therefore, we apply the Novak--Ritter criterion twice to
every $\alpha$ level $\alpha_j = \tfrac{j}{m}$ ($j = 0, \dotsc, m$),
once for the minimum and once for the maximum.
For all $\alpha_j$, the points to be refined are collected in a set.
If a point is selected multiple times for different $\alpha_j$,
it is refined only once.
In addition, we enlarge the search domain $\Omega_{\alpha_j}$
by \SI{10}{\percent}, since the minimum or the maximum might be
near the boundary $\Omega_{\alpha_j}$ and
since the points to be inserted might not be close to
the points to be refined.
We ensure that the size of $\Omega_{\alpha_j}$ is at least $0.05$
in every coordinate direction.
The remaining experiments use $\gamma = 0.1$ as adaptivity.

\begin{algorithm}
  \begin{algorithmic}[1]
    \Function{$\liset = \texttt{fuzzyNovakRitterMethod}$}{%
      $\objfun$, $\gamma$, $m$, $\liset$, $\fuzzy{x}_1$, \dots, $\fuzzy{x}_d$%
    }
      \ForOneLine{$(\*l, \*i) \in \liset$}{$d_{\*l,\*i} \gets 0$}
      \Comment{degrees (number of refinements)}%
      \While{$\setsize{\liset} < \ngpMax$}
        \State{$R \gets \emptyset$}
        \For{$j = 0, \dotsc, m$}
          \State{$\alpha_j \gets j/m$}
          \For{$t = 1, \dotsc, d$}
            \State{$\clint{a_{j,t}, b_{j,t}} \gets \acut[t]{x}{\alpha_j}$}
            \Comment{determine $\alpha_j$-cut}%
            \If{$b_{j,t} - a_{j,t} < 0.05$}
            \Comment{ensure minimal size of $0.05$}%
              \State{%
                $(a_{j,t}, b_{j,t}) \gets
                ((a_{j,t} + b_{j,t})/2 - 0.025,
                (a_{j,t} + b_{j,t})/2 + 0.025)$%
              }\vspace{-1mm}
            \EndIf{}
            \State{%
              $(a_{j,t}, b_{j,t}) \gets
              (a_{j,t} - 0.05 (b_{j,t} - a_{j,t}),
              b_{j,t} + 0.05 (b_{j,t} - a_{j,t}))$%
            }
            \Comment{enlarge by \SI{10}{\percent}}%
            \vspace{-1mm}
          \EndFor{}
          \State{%
            $\liset_j \gets \{(\*l, \*i) \in \liset \mid
            \gp{\*l,\*i} \in \clint{\*a_j, \*b_j} \cap \clint{\*0, \*1}\}$%
          }
          \Comment{%
            set of feasible points
            ($\clint{\*a_j, \*b_j} = \Omega_{\alpha_j}$)%
          }%
          \ForOneLine{$(\*l, \*i) \in \liset_j$}{%
            $r_{\*l,\*i} \gets \setsize{
              \{(\*l', \*i') \in \liset_j \mid
              \objfun(\gp{\*l',\*i'}) \le \objfun(\gp{\*l,\*i})\}
            }$%
          }
          \Comment{ranks}%
          \State{%
            $(\*l^\ast, \*i^\ast) \gets
            \vecargmin_{(\*l,\*i) \in \liset_j} \bracket*{
              (r_{\*l,\*i} + 1)^\gamma
              (\normone{\*l} + d_{\*l,\*i} + 1)^{1 - \gamma}
            }$%
          }
          \Comment{for minimum}%
          \vspace{-0.7mm}
          \State{%
            $(\*l^{\ast\ast}, \*i^{\ast\ast}) \gets
            \vecargmin_{(\*l,\*i) \in \liset_j} \bracket*{
              (\setsize{\liset_j} - r_{\*l,\*i} + 2)^\gamma
              (\normone{\*l} + d_{\*l,\*i} + 1)^{1 - \gamma}
            }$%
          }
          \Comment{for maximum}%
          \State{%
            $R \gets R \cup \{(\*l^\ast, \*i^\ast),
            (\*l^{\ast\ast}, \*i^{\ast\ast})\}$%
          }
        \EndFor{}
        \State{Refine all points in $\liset$ that are in $R$}
        \ForOneLine{$(\*l, \*i) \in R$}{$d_{\*l,\*i} \gets d_{\*l,\*i} + 1$}
      \EndWhile{}
    \EndFunction{}
  \end{algorithmic}
  \caption[Fuzzy Novak--Ritter method]{%
    Fuzzy Novak--Ritter method to generate spatially adaptive sparse grids
    for the fuzzy extension principle.
    Inputs are
    the objective function $\objfun$,
    the adaptivity parameter $\gamma \in \clint{0, 1}$,
    the number of $\alpha$ segments,
    the initial sparse grid $\liset$ as a set of level-index pairs, and
    the $d$ fuzzy intervals $\fuzzy{x}_1, \dotsc, \fuzzy{x}_d$.
    The output is the spatially adaptive sparse grid $\liset$.%
  }%
  \label{alg:fuzzyNovakRitterMethod}%
\end{algorithm}

\paragraph{Convergence of fuzzy intervals on spatially adaptive sparse grids}

\begin{figure}
  \includegraphics{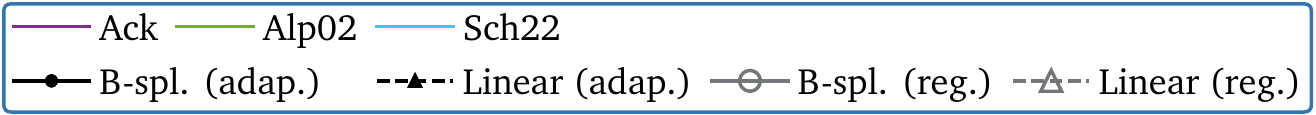}\\[2mm]%
  \includegraphics{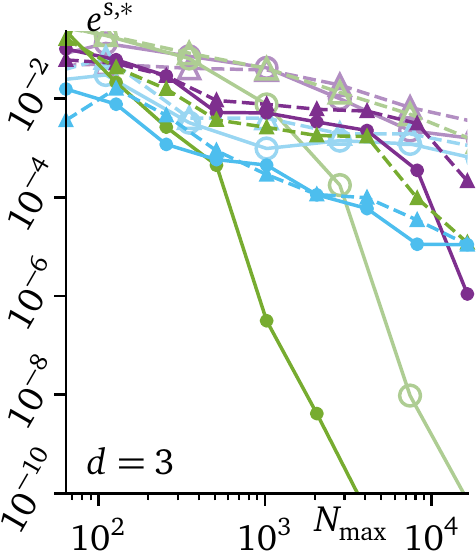}%
  \hfill%
  \includegraphics{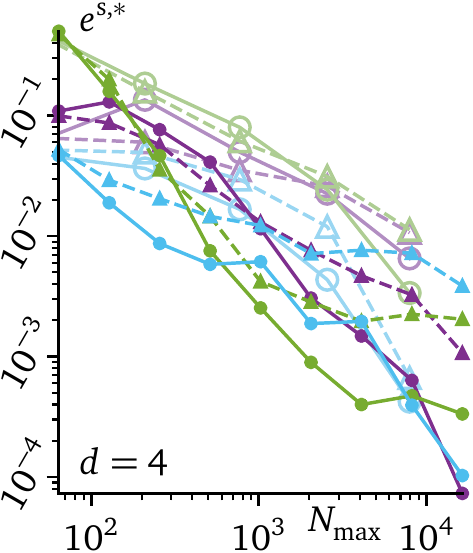}%
  \hfill%
  \includegraphics{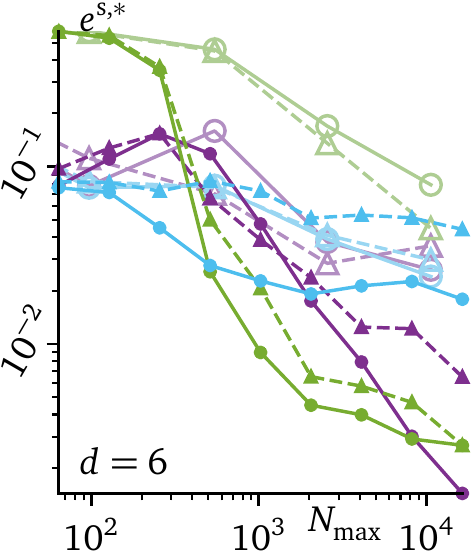}%
  \caption[Fuzzy errors for spatially adaptive sparse grids]{%
    Fuzzy errors
    $e^{\sparse,\ast}
    \ceq \normLtwo{\memfun[\reference]{y} - \memfun[\sparse,\ast]{y}}/
    \normLtwo{\memfun[\reference]{y}}$
    for spatially adaptive sparse grids $\sgset$ \emph{(solid markers)}
    and different objective functions $\objfun$ \emph{(colors)}
    over the number $\ngpMax$ of objective function evaluations.
    For comparison, the results of \cref{fig:resultsFuzzyRegular}
    for regular sparse grids are repeated \emph{(hollow markers).}%
  }%
  \label{fig:resultsFuzzyAdaptive}%
\end{figure}

As we can see in \cref{fig:resultsFuzzyAdaptive},
the spatially adaptive sparse grids generated by the fuzzy Novak--Ritter
method improve results significantly
for both cubic B-spline and piecewise linear surrogates.
However, the performance of the B-spline surrogates benefits more
from the spatial adaptivity.
Even for higher-dimensional settings such as $d = 6$,
the spatial adaptivity helps to decrease the errors by one order of magnitude.
For instance, for the Ack function in six variables,
we can achieve an error of \SI{2.6}{\percent}
with a budget of \num{10000} objective function evaluations (grid points)
on regular sparse grids.
With the same budget and with spatial adaptivity, the error drops below
\SI{0.25}{\percent}.
Conversely, to achieve the same error as in the regular case
(\SI{2.6}{\percent}),
only \SI{1600} evaluations are needed for spatially adaptive grids.

\cleardoublepage

  \setdictum{%
  Money. A social life. A shave.\\
  A Ph.D.\ student needs not such things.%
}{%
  Mike Slackenerny
  (PHD Comics\footnotemark)%
}

\longchapter{%
  Application 1: Topology Optimization%
}{%
  Application 1:\texorpdfstring{\\}{ }Topology Optimization%
}{%
  Application 1 -- Topology Optimization%
}
\footnotetext{\url{http://phdcomics.com/comics/archive.php?comicid=40}}
\label{chap:60topoOpt}

\initial{0em}{N}{ow, we want to investigate}
the first real-world application,
which is the field of topology optimization.
The classical and widely-used method in engineering is shape optimization,
where the shape of a component
(parameterized by $\*x \in \real^d$) has to be
determined such that some objective function value $\objfun(\*x)$ is optimal,
i.e., minimal or maximal.
For example, a bridge over a valley can be built in the shape of a
parabolic arc.
The task of shape optimization is then to choose the coefficients of the
parabola such that the bridge's stability is maximized,
possibly with the constraint that the volume occupied by the bridge
does not exceed a certain value (to save construction costs) or
that the size of the resulting passage meets some size requirements
(e.g., at least \SI{20}{\meter} wide and \SI{6}{\meter} tall).

However, the framework of shape optimization unnecessarily prescribes the
topology of the shapes in the search space \cite{Allaire16Towards}.
In the bridge example, it may well be that a viaduct-type bridge with
three arcs instead of one is more stable while occupying less volume.
We are not able to find such a bridge with shape optimization
in the previous example,
as we have restricted the search space to single-arc bridges.
This issue is resolved by the more sophisticated
framework of topology optimization, where the topology%
\footnote{%
  Two objects are considered ``topologically different''
  if their numbers of ``holes'' differ.
  This stems from the fact that in the field of mathematical topology,
  the \term{genus} (i.e., the number of holes)
  of a topological space is a \term{topological invariant,} i.e.,
  the genus is invariant under homeomorphism.
  If the genera of two topological spaces differ, then they cannot be
  homeomorphic and are thus considered topologically different.%
}
is not given by the user,
but chosen by the optimization algorithm (in a hopefully optimal way),
rendering topology optimization a key area of simulation technology.

Recently, B-splines have been used for
shape optimization \cite{Martin16Formoptimierung} and
topology optimization \multicite{Qian13Topology,Zhang17Topology}.
Sparse grids have been employed for
topology optimization \cite{Huebner14Mehrdimensionale} as well.
In this chapter, we want to combine these two numerical tools,
which have been used in isolation until now,
to perform topology optimization using B-splines on sparse grids.
The two most common approaches for topology optimization are
the \term{level-set method} and
the \term{homogenization method} \cite{Allaire16Towards}.
The level-set method describes the boundary of the object
as the zero level set $\psi^{-1}(0)$ of a function
$\psi\colon \objdomain \to \real$ \term{(level-set function)}
and uses a \pde to iteratively transport this function and,
consequently, the object's boundary \cite{Allaire04Topology}.
However, we want to focus on the second method:
the method of homogenization.

This chapter is structured as follows:
\Cref{sec:61homogenization} explains the homogenization method.
In \cref{sec:62tensors}, we discuss the details of applying B-splines on
sparse grids to this method.
We set up different micro-cell models and scenarios in \cref{sec:63models},
before reviewing numerical results in \cref{sec:64results}.
The results in this chapter have been obtained in collaboration with
Prof.\ Dr.\ Michael Stingl and Daniel Hübner
(both FAU Erlangen-Nürnberg, Germany).
The author of this thesis contributed the parts related to
interpolation and sparse grids, while the collaborators at FAU
studied the engineering and application parts of the joint project
(for example, they provided optimization scenarios and
assessed the quality of the results).

\section{Homogenization and the Two-Scale Approach}
\label{sec:61homogenization}

\minitoc{70mm}{4}

\noindent
We roughly follow the presentation given in
\multicite{%
  Huebner14Mehrdimensionale,%
  Valentin14Hierarchische,%
  Valentin16Hierarchical%
}.
The necessary notation is summarized in
\cref{tbl:glossaryTopologyOptimization}.

\begin{table}
  \setnumberoftableheaderrows{0}%
  \newcommand*{\pnst}[1]{\printnotationsymbol{#1}&\printnotationtext{#1}}%
  \begin{tabular}{%
    >{\kern\tabcolsep}=l<{\kern-1.5mm}+l<{\kern2.9mm}+l<{\kern-1.5mm}+l%
    <{\kern2.9mm}+l<{\kern-1.25mm}+l<{\kern\tabcolsep}%
  }
    \toprulec
    \pnst{\objdomain}&       \pnst{\force}&        \pnst{\densglobal}\\
    \pnst{\dimobjdomain}&    \pnst{\displacement}& \pnst{\denscell}\\
    $d$&\#micro-cell param.& \pnst{\compliance}&   \pnst{\densub}\\
    \pnst{\mcp}&             \pnst{\vol}&          \pnst{\etensor}\\
    &&                       \pnst{\voldens}&      \pnst{\cholfactor}\\
    \bottomrulec
  \end{tabular}%
  \caption[Glossary for topology optimization]{%
    Glossary of the notation for topology optimization.%
  }%
  \label{tbl:glossaryTopologyOptimization}%
\end{table}

\subsection{Homogenization}
\label{sec:611homogenization}

\paragraph{Density function}

Let $\objdomain \subset \real^{\dimobjdomain}$ be the object domain.%
\footnote{%
  We use tildes to denote variables and quantities
  that correspond to the object domain $\objdomain$
  (e.g., $\tilde{\*x}$ is a point in $\objdomain$).
  In contrast, variables without a tilde will correspond
  to the sparse grid domain $\clint{\*0, \*1} = \clint{0, 1}^d$
  (e.g., $\gp{\*l,\*i} \in \clint{\*0, \*1}$ will be a sparse grid point).%
}
Usually, we assume $\dimobjdomain = 2$ or $\dimobjdomain = 3$,
although the method can be generalized to
arbitrary dimensionalities $\dimobjdomain \in \nat$.
Shapes and topologies are described by \term{density functions}
$\densglobal\colon \objdomain \to \clint{0, 1}$.
The function values $\densglobal(\tilde{\*x}) \in \clint{0, 1}$
tell if $\tilde{\*x}$ is contained in the object (value of one) or
not (value of zero).
The \term{homogenization} approach also allows values between
zero and one, giving the physical density of the material in $\tilde{\*x}$.

\paragraph{Optimization of compliance values}

\begin{SCfigure}
  \includegraphics{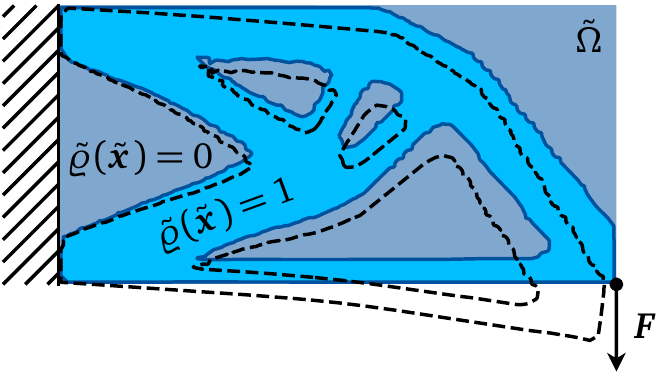}%
  \caption[%
    Example scenario for topology optimization%
  ]{%
    Example scenario for topology optimization.
    An object \emph{\textcolor{hellblau}{(light blue)}}
    is fixed on the left side
    of the object domain $\objdomain$
    \emph{\textcolor{mittelblau!50}{(darker blue)}}
    and deformed by a force $\force$, resulting in a displaced object
    \emph{(dashed).}
    The density function $\densglobal(\tilde{\*x})$ is one inside the object
    and zero outside.%
  }%
  \label{fig:topoOptExample}%
\end{SCfigure}

Furthermore, for every density function $\densglobal$,
let $\compliance(\densglobal)$ be an objective function value.
In our setting, which is shown in \cref{fig:topoOptExample},
we exert a force $\force$ on the object,
measure the resulting deformation, and
compute the \term{compliance} (i.e., the inverse of the stiffness) as
the objective function value~$\compliance(\densglobal)$:
\begin{equation}
  \compliance(\densglobal)
  = \int_{\objdomain} \tr{\force} \displacement_{\densglobal}(\tilde{\*x})
  \diff\tilde{\*x},
\end{equation}
where the \term{displacement function}
$\displacement_{\densglobal}\colon \objdomain \to \real^{\dimobjdomain}$
depends on the density \cite{Huebner14Mehrdimensionale}.
We want to find the density function
with the minimal compliance value:
\begin{equation}
  \label{eq:topoOptProblemContinuous}
  \min_{\densglobal}\, \compliance(\densglobal).
\end{equation}
If we do not impose additional conditions,
then there are often uninteresting trivial solutions.
For example, choosing $\densglobal :\equiv 1$
(i.e., filling the entire domain $\objdomain$ with material)
usually leads to the topology with the
highest stiffness and, thus, the smallest displacement and compliance value.
Therefore, we introduce the following volume constraint:
\begin{equation}
  \frac{\voldens{\densglobal}{\objdomain}}{\vol{\objdomain}} \le \densub,\quad
  \voldens{\densglobal}{\objdomain}
  \ceq \int_{\objdomain} \densglobal(\tilde{\*x}) \diff\tilde{\*x},\quad
  \vol{\objdomain}
  \ceq \voldens{1}{\objdomain},
\end{equation}
where $\vol{\objdomain} = \int_{\objdomain} 1 \diff\tilde{\*x}$
is the volume of the object domain and
$\densub \in \clint{0, 1}$ is an upper bound on the volume fraction.

\subsection{Two-Scale Approach}
\label{sec:612twoScale}

\paragraph{Discretization and two-scale approach}

Of course, we cannot solve the problem \eqref{eq:topoOptProblemContinuous}
numerically,
as there are infinitely many density functions $\densglobal$.
For simplicity, we assume that $\objdomain$ is some hyper-rectangle
$\clint{\tilde{\*a}, \tilde{\*b}}
= \clint{\tilde{a}_1, \tilde{b}_1} \times \dotsb \times
\clint{\tilde{a}_{\dimobjdomain}, \tilde{b}_{\dimobjdomain}}$;
if it is not, we replace $\objdomain$ with its bounding box.
The object domain $\objdomain$ can then be split into
$M_1 \times \dotsb \times M_{\dimobjdomain}$
equally-sized and axis-aligned sub-hyper-rectangles,
which we call \term{macro-cells}
(where $M_1, \dotsc, M_{\dimobjdomain} \in \nat$).

\begin{figure}
  \includegraphics{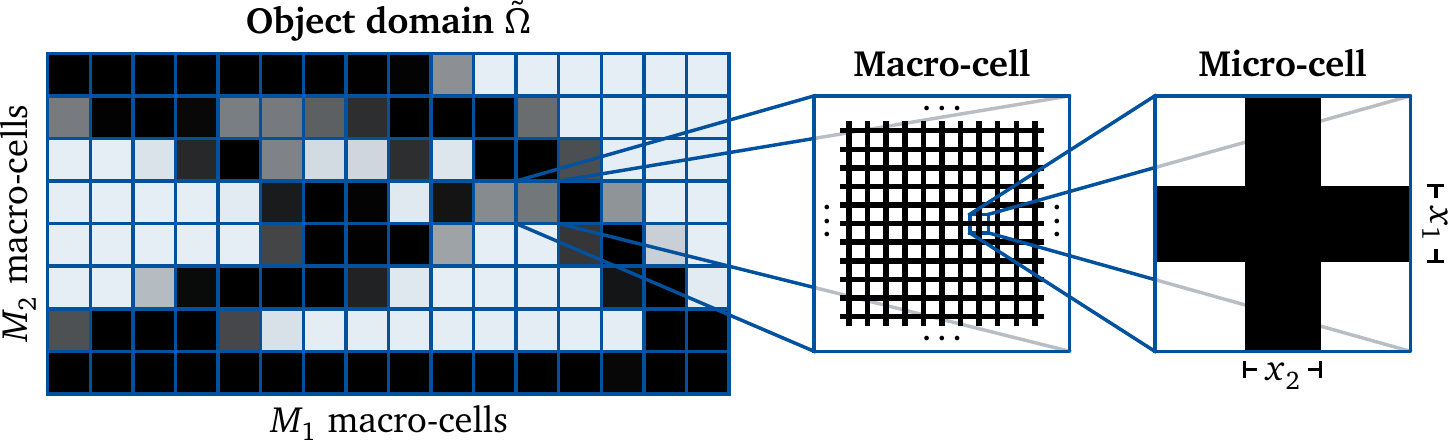}%
  \caption[%
    Two-scale approach for topology optimization%
  ]{%
    Two-scale approach to discretize the homogenized topology
    optimization problem in two dimensions ($\dimobjdomain = 2$).
    \emph{Left:} The object domain $\objdomain$ is
    subdivided into $M_1 \times M_2$ macro-cells,
    each with its own density \emph{(gray squares).}
    \emph{Center:} Every macro-cell is the repetition of infinitesimally small
    periodic micro-cells.
    \emph{Right:} The shape of the structure in every micro-cell is
    described by a micro-cell model with $d$ parameters $x_1, \dotsc, x_d$.
    Here, the micro-cell model is a cross with two parameters
    that represent the thickness of each crossbar.%
  }%
  \label{fig:twoScale}%
\end{figure}

In the \term{two-scale approach,}
we assume the material of the macro-cells to be
repetitions of infinitesimally small periodic structures
(i.e., identical for each macro-cell),
called \term{micro-cells.}
These micro-cells have a specific shape, which is parameterized by
$d$ \term{micro-cell parameters} $x_1, \dotsc, x_d$,
normalized to values in the unit interval $\clint{0, 1}$.
For instance, in two dimensions,
this shape may be an axis-aligned cross
with thicknesses $x_1$ and $x_2$, as shown in \cref{fig:twoScale}.
The choice of a suitable \term{micro-cell model}
(parametrization of the micro-cells)
depends on the optimization scenario and has to be done a priori.

\paragraph{Elasticity tensors}

Note that while the shape of all micro-cells in one macro-cell is identical,
the micro-cell parameters corresponding to different macro-cells differ
in general.
This enables varying densities in different regions of $\objdomain$.
We denote the micro-cell parameters corresponding to the $q$-th macro-cell
with $\mcp{q} = (\mcpentry{1}{q}, \dotsc, \mcpentry{d}{q}) \in
\clint{\*0, \*1} = \clint{0, 1}^d$,
where $q = 1, \dotsc, M$ and
$M \ceq M_1 \dotsm M_{\dimobjdomain}$ is the number of macro-cells.
With linear elasticity,
one can compute so-called \term{elasticity tensors} $\etensor^{(q)}$,
which encode information about the material properties
of the different macro-cells.
The elasticity tensors can be written as symmetric matrices
in $\real^{3 \times 3}$ (for $\dimobjdomain = 2$) or
in $\real^{6 \times 6}$ (for $\dimobjdomain = 3$).%
\footnote{%
  In general, the elasticity tensor is a fourth-order tensor in
  $\real^{\dimobjdomain \times \dimobjdomain \times \dimobjdomain \times \dimobjdomain}$.
  One can reduce the size of the tensor by exploiting various symmetries
  \cite{Huebner14Mehrdimensionale}
  to obtain $6$ or $21$ stiffness coefficients
  in two or three dimensions, respectively.
  These coefficients can then be expressed as a symmetric matrix.%
}
To simplify the following considerations,
we assume that $\dimobjdomain = 3$, i.e.,
$\etensor^{(q)} \in \real^{6 \times 6}$.
The elasticity tensors are usually computed as the solution of a \fem problem
\term{(micro-problem).}
Once all $\etensor^{(q)}$ are known,
we can compute the compliance value
by solving another \fem problem \term{(macro-problem),}
see \cite{Allaire04Topology} and \cite{Huebner14Mehrdimensionale}.

\vspace{-0.5em}

\paragraph{Discretized optimization problem}

The new optimization problem emerging from the
two-scale discretization process has the form
\begin{subequations}
  \label{eq:topoOptProblemDiscrete}
  \setlength{\belowdisplayskip}{5pt}%
  \begin{gather}
    \min J(\mcp{1}, \dotsc, \mcp{M}),\quad
    \mcp{1}, \dotsc, \mcp{M} \in \clint{\*0, \*1}
    \quad\text{s.t.}\quad
    \densmean(\mcp{1}, \dotsc, \mcp{M}) \le \densub,\\
    \densmean(\mcp{1}, \dotsc, \mcp{M})
    \ceq \frac{1}{M} \sum_{q=1}^M \denscell^{(q)}(\mcp{q}).
  \end{gather}
\end{subequations}
Here, $\denscell^{(q)}(\mcp{q}) \in \clint{0, 1}$ is the
density of the $q$-th macro-cell with micro-cell parameter $\mcp{q}$
(i.e., the fraction of material volume of one micro-cell
with respect to its total volume)
and $\densmean(\mcp{1}, \dotsc, \mcp{M}) \in \clint{0, 1}$
is the resulting total mean density.
This discretized optimization problem can now be implemented and
solved numerically.

\section{Approximating Elasticity Tensors}
\label{sec:62tensors}

\paragraph{Optimization process}

\minitoc{75mm}{7}

During the process of solving \cref{eq:topoOptProblemDiscrete},
optimization algorithms typically
evaluate the objective function $\compliance(\mcp{1}, \dotsc, \mcp{M})$
iteratively at different \term{micro-cell parameter combinations}
$(\mcp{1}, \dotsc, \mcp{M}) \in (\real^d)^M$.
Every evaluation of $\compliance$ corresponds to one solution of a
macro-problem.
However, to solve the macro-problem,
the elasticity tensors $\etensor^{(q)}$ of all $M$ macro-cells
need to be known.
Hence, in every optimization iteration, it is necessary to solve
one macro-problem and $M$ micro-cell problems,
all with the \fem.
This naive approach has two major drawbacks, which we explain in
the following.

\subsection{Drawbacks of the Naive Approach}
\label{sec:621drawbacks}

\paragraph{Drawback 1: Computation time}

First, this approach is computationally infeasible
even for simple micro-cell models and optimization scenarios.
The computation of a single elasticity tensor usually takes seconds to
minutes.
All $M$ micro-cell problems per optimization iteration
can be solved in parallel without any communication.
However, $M$ is typically in the range of thousands and
there are thousands or tens of thousands optimization iterations
(the optimization problem is $(d \cdot M)$-dimensional!).
This implies that the overall computation may still take
several days or even weeks to complete.

\paragraph{Drawback 2: Approximation of gradients}

Second, most optimization algorithms require gradients of the
objective function and of the constraints, i.e.,%
\begin{equation}
  \partialderiv{\partialdiff{} x_t}{\etensor^{(q)}}(\mcp{q}),\quad
  \partialderiv{\partialdiff{} x_t}{\denscell^{(q)}}(\mcp{q}),\qquad
  q = 1, \dotsc, M,\quad
  t = 1, \dotsc, d.
\end{equation}
However, in general, both gradients are unavailable and
have to be approximated by finite differences.
This introduces new error sources and
increases the number of elasticity tensors to be evaluated,
further slowing down the solution process.
Additionally, the number of optimization iterations necessary to
achieve convergence might increase
if there are discontinuities in the objective function
or its gradient.
Such discontinuities can already be caused by the inexact solution of the \fem.
If we need Hessians or other higher-order derivatives,
then the issues even worsen.

\subsection{B-Splines on Sparse Grids for Topology Optimization}
\label{sec:622BSplines}

\paragraph{Elasticity tensor function}

As a remedy, we replace the costly elasticity tensors with cheap surrogates.
If we assume that all macro-cells use the same micro-cell model,
the elasticity tensor $\etensor^{(q)}$ of the $q$-th macro-cell
with the parameter $\mcp{q} \in \clint{\*0, \*1}$
can be written as the value $\etensor(\mcp{q})$ of some function
$\etensor\colon \clint{\*0, \*1} \to \real^{6 \times 6}$
(assuming that $\dimobjdomain = 3$) at the point $\mcp{q}$.
In the following,
$\etensor\colon \clint{\*0, \*1} \to \real^m$
gives $m \in \nat$ values from which the symmetric elasticity tensor
can be uniquely reconstructed,
i.e., $m = 6$ for $\dimobjdomain = 2$ and $m = 21$ for $\dimobjdomain = 3$.
The vector-valued/matrix-valued versions of $\etensor$
will be used interchangeably.

\paragraph{Elasticity tensor surrogate}

The idea is to use B-splines on sparse grids to approximate
the elasticity tensor function $\etensor$.
In contrast to the theoretical framework that we established in
\cref{chap:20sparseGrids,chap:30BSplines,chap:40algorithms},
the function to be interpolated is not scalar-valued, but vector-valued.
This means that we have to construct $m$ sparse grid interpolants
$\etensorentryintp{j}$
for the $m$ components $\etensorentry{j}$ of $\etensor$ ($j = 1, \dotsc, m$).
Note that one could generate different spatially adaptive sparse grids for the
different components $\etensorentryintp{j}$.
However, it is not possible to evaluate only specific entries of $\etensor$
without also evaluating all other entries,
which means that we would waste computational resources by selecting only
a subset of the calculated entries.
Therefore, we use the same grid for all components.

Additionally, we approximate the density $\denscell^{(q)}$
of the $q$-th macro-cell with a surrogate $\denscellintp$ using
B-splines on the same sparse grid as for $\etensorentryintp{j}$
for reasons of implementation,
resulting in $m + 1$ sparse grid interpolants in total.
From a theoretical perspective, this is not necessary,
since the density can be explicitly calculated with simple formulas
for most micro-cell models, independently of evaluations of the
elasticity tensor.

\paragraph{Advantages}

Our approach has multiple obvious advantages:
\begin{itemize}
  \item
  The sparse grid interpolant $\etensorintp$ has to be generated only
  once in an \term{offline step} before the optimization algorithm starts.
  During the optimization \term{(online phase),}
  only inexpensive evaluations of $\etensorintp$ are performed,
  saving much computation time.
  
  \item
  Sparse grids ease the curse of dimensionality, which prohibits
  conventional full grid interpolation methods if $d > 4$.
  
  \item
  With spatially adaptive sparse grids and a suitable refinement criterion,
  we can spend more grid points in regions of interest of $\etensor$,
  e.g., regions with large oscillations.
  
  \item
  By using B-splines as basis functions,
  the interpolant $\etensorintp$ will be more accurate
  than with piecewise linear basis functions.
  In addition, we can calculate its derivatives
  $\tpartialderiv{\partialdiff{} x_t}{\etensorintp}(\mcp{q})$
  fast and explicitly,
  accelerating the speed of convergence of the optimizer.
\end{itemize}

\breakpagebeforenextheadingtrue
\subsection{Cholesky Factor Interpolation}
\label{sec:623cholesky}

\paragraph{Positive definiteness of elasticity tensors}

Unfortunately, just replacing elasticity tensors with
B-spline surrogates often does not lead to correct results in practice.
Experiments show that for only for some sparse grids,
the optimization algorithm converges to an optimal point
\cite{Valentin16Hierarchical}.
The optimization algorithm crashes for most spatially adaptive grids,
not being able to find any meaningful optimum.
The root of the problem proves to be that
the interpolated elasticity tensors $\etensorintp(\*x)$ are not
positive definite for specific
micro-cell parameters $\*x \in \clint{\*0, \*1}$.
However, indefinite or even negative definite tensors $\etensorintp$
would mimic unphysical behavior.%
\footnote{%
  In the scalar case, this is analogous to Hooke's law for linear springs,
  where the force $F = kx$ needed to displace the end of a spring
  (fixed at the other end) by $x$ is proportional to $x$.
  The proportionality constant $k$ (which corresponds to the elasticity tensor)
  has to be positive.%
}
Hence, it is imperative for the optimization process that
the interpolated elasticity tensors are \spd.

\paragraph{Positive definiteness of sparse grid interpolants}

\begin{figure}
  \subcaptionbox{%
    Minimal eigenvalue of $\etensorintp(\*x)$%
    \label{fig:cholesky1}%
  }[65mm]{%
    \includegraphics{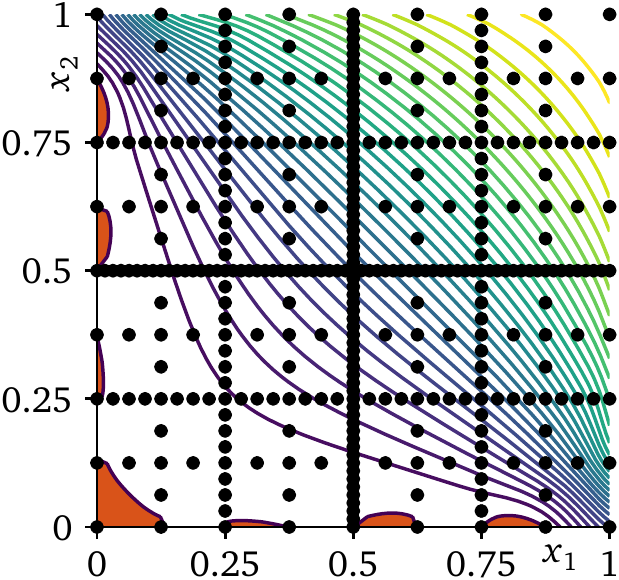}%
  }%
  \hfill%
  \subcaptionbox{%
    Minimal eigenvalue of $\etensorcholintp(\*x)$%
    \label{fig:cholesky2}%
  }[65mm]{%
    \includegraphics{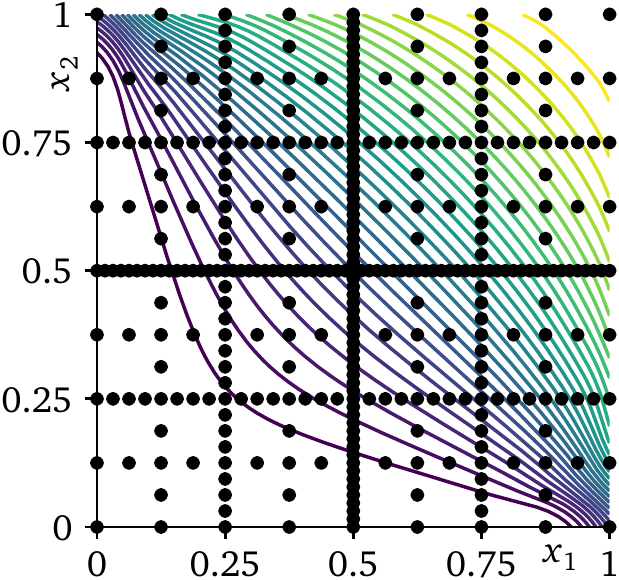}%
  }%
  \hfill\hfill%
  \raisebox{2.2mm}{\includegraphics{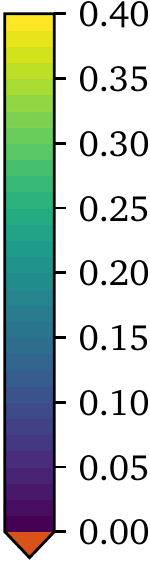}}%
  \caption[%
    Minimal eigenvalue of interpolated elasticity tensors%
  ]{%
    Minimal eigenvalue \emph{(colored contour lines)}
    of elasticity tensor surrogates
    for the 2D cross model ($\dimobjdomain = 2$, $d = 2$)
    \vspace{-0.3em}%
    and cubic hierarchical B-splines $\bspl{\*l,\*i}{p}$ ($p = 3$) on
    the regular sparse grid $\coarseregsgset{n}{d}{b}$ \emph{(dots)}
    with $n = 6$ and $b = 4$.
    \emph{Left:} The minimal eigenvalue of $\etensorintp(\*x)$
    becomes negative in some regions \emph{\textcolor{C1}{(red areas)}}
    of the domain $\clint{\*0, \*1}$,
    indicating that $\etensorintp(\*x)$ is not positive definite.
    \emph{Right:} The minimal eigenvalue of $\etensorcholintp(\*x)$
    is non-negative in the whole domain $\clint{\*0, \*1}$.%
  }%
  \label{fig:cholesky}%
\end{figure}

Interpolation on sparse grids per se does not preserve
positive definiteness.
A counterexample is shown in \cref{fig:cholesky1},
which displays the minimal eigenvalue of the elasticity tensor surrogate
resulting from interpolation on a regular sparse grid.
As the positivity of the diagonal is a necessary condition
for positive definiteness,
small oscillations of the interpolant of some entries
already make the whole elasticity tensor non-positive-definite.

These oscillations are more likely to occur near the boundary of the domain
$\clint{\*0, \*1}$, such that there are larger regions
where the interpolated tensor is not positive definite anymore.
The reason is two-fold:
First, sparse grids without boundary points
are notoriously biased towards the center of the domain,
as they place only few points near the boundary \cite{Pflueger10Spatially}.
This leads to a loss of interpolation accuracy near the boundary
when compared to the center of $\clint{\*0, \*1}$.
Second, both the minimal eigenvalue of $\etensor(\*x)$ and the norm of its
gradient with respect to $\*x$ are small near $x_1 = 0$ or $x_2 = 0$.
Consequently, the ``surface'' of the minimal eigenvalue function
is rather flat in these regions and almost vanishes,
facilitating the existence of negative eigenvalues of
surrogate functions $\etensorintp(\*x)$.

Note that for most micro-cell models,
the optimization algorithm often evaluates the objective function
$\compliance(\mcp{1}, \dotsc, \mcp{M})$ at micro-cell parameter
combinations for which many of the points $\mcp{j}$ are near the boundary
of $\clint{\*0, \*1}$.
This is because many of the macro-cells will either be empty or
fully filled with material, which usually corresponds to micro-cell
parameters near zero or one, respectively.
Thus, $\etensorintp$ is frequently evaluated in the regions of
indefiniteness, which further worsens the issue.

\pagebreak

\paragraph{Positivity-preserving methods}

Even in one dimension, it cannot be guaranteed that the interpolant of
positive data remains positive,
which is a key problem in the estimation of probability densities
\multicite{Pflueger10Spatially,Griebel10Finite,Franzelin16From}.
Just clamping the interpolated values via $\max(\cdot, 0)$ does not help:
In our application, the tensor may still be indefinite;
additionally, the calculated gradients of the interpolants do not match
the actual gradients anymore.
In density estimation, clamping a density-like function changes its
integral, making it necessary to recalculate its normalization constant
\cite{Franzelin17Data}.

One possible workaround is to apply a continuous injective transformation
$T\colon \posreal \to \real$ on the positive values (e.g., $\ln$),
then interpolate the resulting values, and finally
apply the inverse transformation $T^{-1}\colon \real \to \posreal$
on the interpolated values (e.g., $\exp$).%
\footnote{%
  Formally, the inverse function $T^{-1}\colon T(\posreal) \to \posreal$
  is only defined on the image $T(\posreal)$ of $T$,
  which might not be the whole real line.
  However, we assume that $T^{-1}$ can be ``reasonably'' extended to $\real$
  (e.g., $T \ceq \sqrt{\cdot}$ and $T^{-1} = ({\cdot})^2$).%
}
For the piecewise linear hierarchical basis, another approach has been
developed recently \cite{Franzelin17Data},
maintaining the positivity by inserting additional sparse grid points.
In the context of spline approximation,
positivity-preserving approximation schemes based on so-called
quasi-interpolation are known \cite{Hoellig13Approximation}.
For our application, for which we need to preserve positive definiteness,
it is conceivable that one could apply these positivity-preserving
methods in the eigenspace,
interpolating the positive eigenvalues.

\paragraph{Interpolation of Cholesky factors}

Instead, we pursue a different, more canonical
approach based on Cholesky factorization:

\begin{proposition}[Cholesky factorization]
  For every \spd matrix $\etensor \in \real^{6 \times 6}$,
  there is a unique upper triangular matrix
  $\cholfactor \in \real^{6 \times 6}$
  with positive diagonal entries such that
  \begin{equation}
    \etensor
    = \tr{\cholfactor} \cholfactor.
  \end{equation}
\end{proposition}

\begin{proof}
  See \cite{Benoit24Note} or \cite{Freund07Stoer}.
\end{proof}

In one dimension, the Cholesky factorization is equivalent
to the application of a transformation $T$ as above by choosing
$T \ceq \sqrt{\cdot}$ and $T^{-1} = (\cdot)^2$.
Our approach is as follows:
\begin{enumerate}
  \item
  Define $\cholfactor\colon \clint{\*0, \*1} \to \real^{6 \times 6}$
  as the Cholesky factor of
  $\etensor\colon \clint{\*0, \*1} \to \real^{6 \times 6}$, i.e.,
  $\etensor(\*x) = \tr{\cholfactor(\*x)} \cholfactor(\*x)$
  for all $\*x \in \clint{\*0, \*1}$.
  
  \item
  During the grid generation (offline phase),
  evaluate $\etensor(\gp{\*l,\*i})$ at the grid points $\gp{\*l,\*i}$,
  compute the Cholesky factors $\cholfactor(\gp{\*l,\*i})$ of
  $\etensor(\gp{\*l,\*i})$,
  and interpolate them instead of the elasticity tensors
  to obtain an interpolant
  $\cholfactorintp\colon \clint{\*0, \*1} \to \real^{6 \times 6}$.
  
  \item
  During the optimization (online phase),
  every time the value $\etensor(\*x)$ of an elasticity tensor is needed,
  the interpolant $\cholfactorintp(\*x)$ is evaluated and we return
  \begin{equation}
    \etensorcholintp(\*x)
    \ceq \tr{\cholfactorintp(\*x)} \cholfactorintp(\*x).
  \end{equation}
\end{enumerate}

\paragraph{Advantages of Cholesky factor interpolation}

As shown in \cref{fig:cholesky2},
the resulting elasticity tensor surrogate $\etensorcholintp$
is positive semidefinite on the whole domain and
positive definite almost everywhere:
The surrogate $\etensorcholintp(\*x)$ is singular if and only if
$\cholfactorintp(\*x)$ is singular, which is in general
only the case on a negligible null set in $\clint{\*0, \*1}$.

Another advantage of this approach is that not only the
positive definiteness, but also the explicit differentiability
of the surrogate $\etensorcholintp$ is preserved.
The gradient can be computed easily and fast with the product rule:
\begin{equation}
  \label{eq:choleskyFactorDerivative}
  \partialderiv{\partialdiff{} x_t}{\etensorcholintp}(\*x)
  = \tr{\cholfactorintp(\*x)} \cdot
  \partialderiv{\partialdiff{} x_t}{\cholfactorintp}(\*x) +
  \tr{\partialderiv{\partialdiff{} x_t}{\cholfactorintp}(\*x)} \cdot
  \cholfactorintp(\*x),\quad
  t = 1, \dotsc, d,
\end{equation}
where both the sparse grid interpolant $\cholfactorintp(\*x)$ and
its derivative $\partialderiv{\partialdiff{} x_t}{\cholfactorintp}(\*x)$
are known.
As discussed above,
this is key to the applicability of gradient-based optimization.

\section{Micro-Cell Models and Optimization Scenarios}
\label{sec:63models}

\minitoc{62mm}{3}

\noindent
In the following, we present the different micro-cell models
and optimization scenarios for which we perform numerical
experiments in the next section.

\subsection{Micro-Cell Models}
\label{sec:631models}

\begin{figure}
  \subcaptionbox{%
    \lefthphantom{2D cross}{(2D-C, $d = 2$)}\\(2D-C, $d = 2$)%
    \label{fig:microCell_1}%
  }[31mm]{%
    \includegraphics{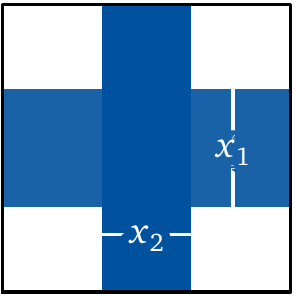}%
  }%
  \hfill%
  \subcaptionbox{%
    2D framed cross\\(2D-FC, $d = 4$)%
    \label{fig:microCell_2}%
  }[31mm]{%
    \includegraphics{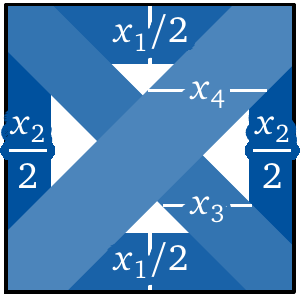}%
  }%
  \hfill%
  \subcaptionbox{%
    2D sheared cross\\\rlap{\hspace*{11mm}{(2D-SC, $d = 3$)}}%
    \label{fig:microCell_3}%
  }[41mm]{%
    \includegraphics{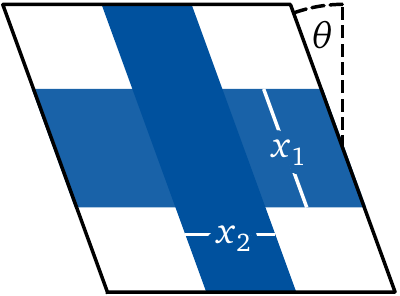}%
  }%
  \hfill%
  \subcaptionbox{%
    2D sheared framed cross (2D-SFC, $d = 5$)%
    \label{fig:microCell_4}%
  }[37.5mm]{%
    \hspace*{-45mm}%
    \rlap{\includegraphics{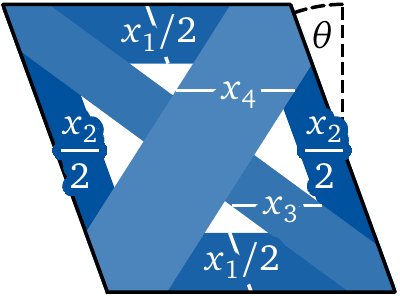}}%
  }\\[2mm]%
  \subcaptionbox{%
    3D cross (3D-C, $d = 3$)%
    \label{fig:microCell_5}%
  }[44mm]{%
    \includegraphics{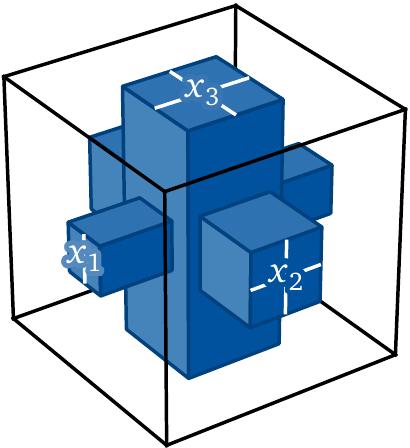}%
  }%
  \quad%
  \subcaptionbox{%
    3D sheared cross (3D-SC, $d = 5$)%
    \label{fig:microCell_6}%
  }[56mm]{%
    \hspace*{5mm}%
    \includegraphics{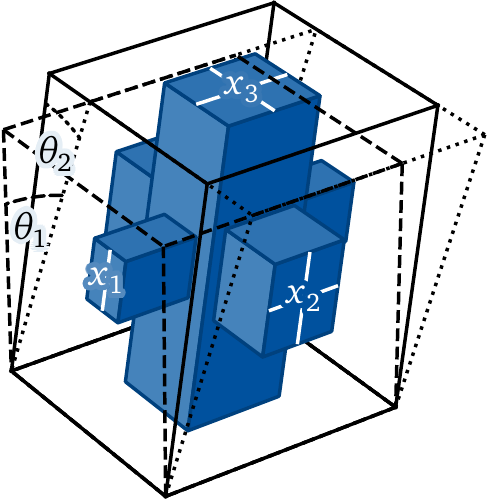}%
  }%
  \caption[Types of micro-cell models]{%
    Types of micro-cell models in two dimensions \emph{(top row)}
    and three dimensions \emph{(bottom row)}.%
  }%
  \label{fig:microCell}%
\end{figure}

We use the various micro-cell models that are depicted in \cref{fig:microCell}.
The models differ in the spatial dimensionality $\dimobjdomain$
and the number $d$ of micro-cell parameters
$\*x \in \clint{\*0, \*1} = \clint{0, 1}^d$.
Note that the presented models are only some examples.
One can easily design complicated micro-cell models
with larger numbers of parameters.

\paragraph{Orthogonal (non-sheared) models in two dimensions}

The basic component of the four two-dimensional models
is a square with a cross (\cref{fig:microCell_1})
of two axis-aligned orthogonal bars,
whose widths are determined by two micro-cell parameters $x_1$ and $x_2$.
The micro-cell parameters are ratios of the bar widths
to the edge lengths of the micro-cell
(although the actual micro-cells are infinitesimally small).
This results in the \emph{cross model.}
For the \term{framed cross model} (\cref{fig:microCell_2}),
we add a diagonal cross with orthogonal bars
of widths $x_3$ and $x_4$ (horizontally measured).
To simplify the boundary treatment,
we shift the contents of the framed cross micro-cell by
\SI{50}{\percent} of the micro-cell's edge lengths in both directions,
such that previous corners of the micro-cell correspond to the new center.

\paragraph{Sheared models in two dimensions}

Both of these models can be extended by shearing.
The idea is to increase the stability of the resulting macro-structure
with respect to forces that act at angles other than
\ang{0} and \ang{90} (cross model) or
\ang{0}, \ang{90}, and \ang{45} (framed cross model).
If we just rotated the crosses in the micro-cells,
then the micro-structure would not be periodic.
Instead, we shear the whole micro-cell in the horizontal direction,
where the shearing angle $\theta$ is an additional micro-cell parameter,
which gives us another degree of freedom.%
\footnote{%
  To be more precise, the angle $\theta$ corresponds to an
  additional micro-cell parameter $x_3$ (sheared cross) or
  $x_5$ (sheared framed cross) that is determined by normalization
  from $\clint{-0.35\pi, 0.35\pi}$, i.e., $\theta/(0.7\pi) + 1/2$.%
}
This results in the \term{sheared cross model} (\cref{fig:microCell_3})
and \term{sheared framed cross model} (\cref{fig:microCell_4})
with three and five micro-parameters each.

\paragraph{Models in three dimensions}

The two-dimensional cross model can be transferred to three
spatial dimensions by just adding another bar in the new dimension.
Each of the three bars has square cross-section with given edge lengths
$x_1$, $x_2$, or $x_3$, respectively,
resulting in the \term{3D cross model} with three micro-cell parameters
(\cref{fig:microCell_5}).
By shearing in the two horizontal directions,
we obtain two new degrees of freedom $\theta_1$ and $\theta_2$
(shearing angles).
The emerging \term{3D sheared cross model} has five micro-cell parameters
(\cref{fig:microCell_6}).

\subsection{Test Scenarios}
\label{sec:632scenarios}

\begin{figure}
  \subcaptionbox{%
    2D cantilever%
    \label{fig:topoOptScenario_1}%
  }[72mm]{%
    \includegraphics{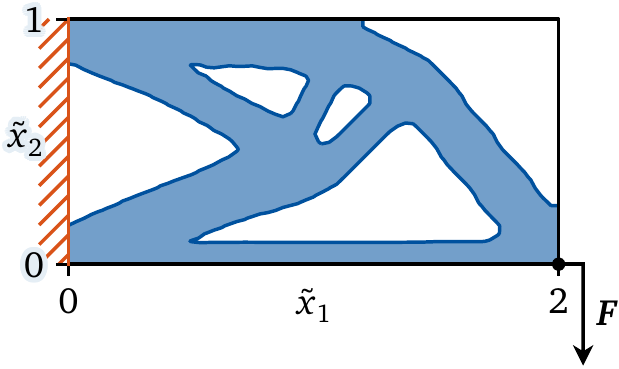}%
  }%
  \hfill%
  \smash{\raisebox{-22mm}{\subcaptionbox{%
    2D L-shape%
    \label{fig:topoOptScenario_2}%
  }[72mm]{%
    \includegraphics{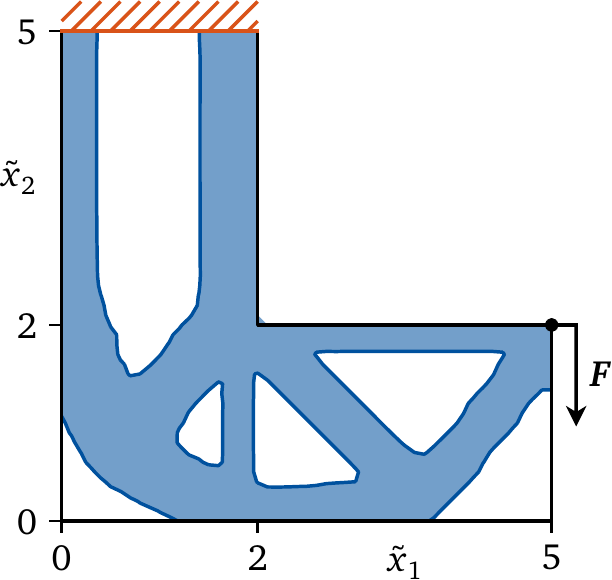}%
  }}}%
  \\[7mm]%
  \subcaptionbox{%
    3D cantilever%
    \label{fig:topoOptScenario_3}%
  }[72mm]{%
    \includegraphics{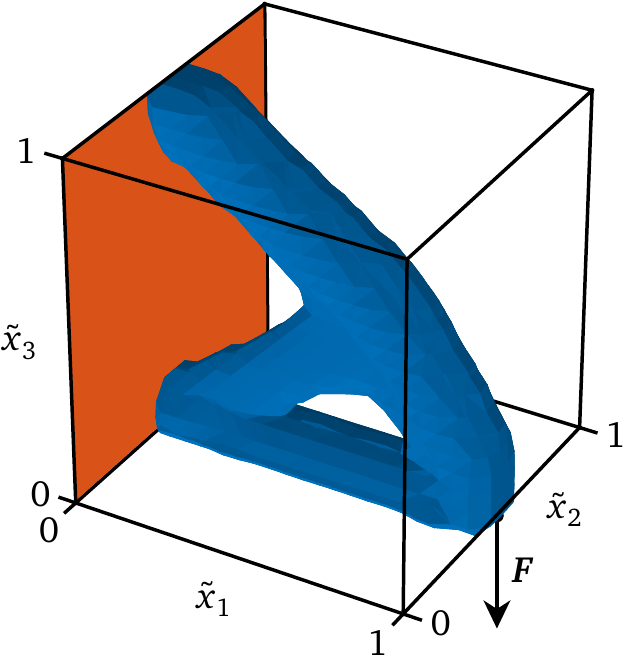}%
  }%
  \hfill%
  \subcaptionbox{%
    3D center-load%
    \label{fig:topoOptScenario_4}%
  }[72mm]{%
    \includegraphics{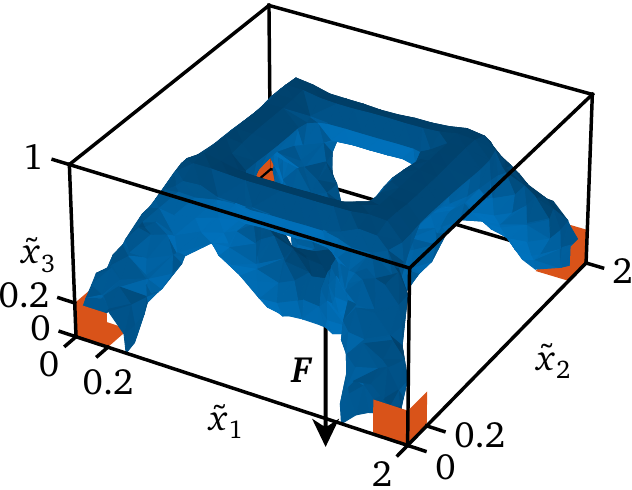}%
  }%
  \caption[Test scenarios in topology optimization]{%
    Test scenarios in topology optimization in
    two and three spatial dimensions.
    Shown are
    the domains $\objdomain$,
    load points,
    locations of homogeneous Dirichlet boundary conditions
    \emph{\textcolor{C1}{(red)},} and
    exemplary optimal structures \emph{\textcolor{C0}{(blue)}.}%
  }%
  \label{fig:topoOptScenario}%
\end{figure}

To benchmark the performance of the new method,
we take a subset of the scenarios given in \cite{Valdez17Topology},
which reviews more than 100 papers on topology optimization
to determine the most common test scenarios in the field.
The geometry and the boundary conditions of the
four scenarios (two for each 2D and 3D)
are given in \cref{fig:topoOptScenario} (dimensions in meters).
In contrast to \cite{Valdez17Topology},
we only use single-point loads
(i.e., not loads applied to line segments, areas, or volumes)
for implementational reasons.
The upper bound on the density (see \cref{sec:611homogenization})
is $\densub = \SI{50}{\percent}$ for the 2D scenarios and
$\densub = \SI{10}{\percent}$ for the 3D scenarios.
As in \cite{Sigmund01Line} and for reasons of simplicity,
we apply a force $\force$ with unit value
(i.e., $\norm[2]{\force} = \SI{1}{\newton}$),
and we use a hypothetical material with
a Young's modulus (stiffness) of \SI{1}{\pascal} and
a Poisson ratio (transversal expansion to axial compression) of $0.3$.

\section{Implementation and Numerical Results}
\label{sec:64results}

\minitoc[-3mm]{73mm}{5}

\noindent
In this final section of the chapter,
we study optimal results of the test scenarios and
analyze interpolation errors and optimization results
for topology optimization with B-spline surrogates on sparse grids.

\subsection{Implementation}
\label{sec:641implementation}

In the following, for simplicity,
we combine the two functions to be interpolated,
i.e., the Cholesky factor
$\cholfactor\colon \clint{\*0, \*1} \to \real^{6 \times 6}$ and
the micro-cell density $\denscell\colon \clint{\*0, \*1} \to \real$,
to one single objective function
$\*\objfun\colon \clint{\*0, \*1} \to \real^{m+1}$,
from which both functions can be recovered.

\paragraph{Overview of offline and online phase}

Our method is divided into an offline phase and an online phase,
both of which are sketched in \cref{fig:topoOptPhases}.
The offline phase consists of
generating the spatially adaptive sparse grid
$\sgset = \{\gp{\*l_k,\*i_k} \mid k = 1, \dotsc, \ngp\}$,
solving the corresponding micro-problems,
computing the Cholesky factors, and
hierarchizing the Cholesky factor entries and micro-cell densities
to obtain the sparse grid interpolant $\vsgintp$.
Each optimization iteration of the online phase consists of
evaluating the interpolant $\vsgintp$
for each micro-cell parameter $\mcp{j}$ ($j = 1, \dotsc, M$),
reconstructing the elasticity tensor $\etensorcholintp$ from
the Cholesky factors $\cholfactorintp$,%
\footnote{%
  In addition, the partial derivatives
  $\partialdiff{} \etensorcholintp/\partialdiff{} x_t$
  ($t = 1, \dotsc, d$)
  are evaluated using \cref{eq:choleskyFactorDerivative}.
  This is necessary to employ gradient-based optimization.%
}
and solving the macro-problem to retrieve the approximated compliance value
$\complianceintp(\mcp{1}, \dotsc, \mcp{M})$.
The superscript in $\complianceintp$ indicates that
we do not use the exact elasticity tensors $\etensor$
to compute the compliance value,
but rather the reconstructed and interpolated tensors
$\etensorcholintp$.

\begin{figure}
  \tikzset{
    myCircle/.style={
      circle,
      fill=mittelblau!30,
      draw=mittelblau,
      inner sep=0.5mm,
    }
  }%
  \subcaptionbox{%
    Offline phase (without the actual grid generation).%
  }[149mm]{%
    \begin{tikzpicture}
      \node[myCircle] (points) at (0mm,0mm) {%
        $
          \begin{matrix}
            \gp{\*l_1,\*i_1},\\
            \dots,\\
            \gp{\*l_{\ngp},\*i_{\ngp}}
          \end{matrix}
        $%
      };
      \node[myCircle] (elasticityTensors) at (43mm,0mm) {%
        $
          \begin{matrix}
            \etensor(\gp{\*l_1,\*i_1}),\\
            \dots,\\
            \etensor(\gp{\*l_{\ngp},\*i_{\ngp}})
          \end{matrix}
        $%
      };
      \node[myCircle] (choleskyFactors) at (80mm,0mm) {%
        $
          \begin{matrix}
            \cholfactor(\gp{\*l_1,\*i_1}),\\
            \dots,\\
            \cholfactor(\gp{\*l_{\ngp},\*i_{\ngp}})
          \end{matrix}
        $%
      };
      \node[myCircle] (choleskyInterpolant) at (118mm,0mm) {%
        $
          \begin{matrix}
            \cholfactorintp\colon \clint{\*0, \*1}\\
            {} \to \real^{6 \times 6}
          \end{matrix}
        $%
      };
      \draw[->,draw=C0] (points) -- node[above] {%
        \footnotesize{}micro-problem%
      } (elasticityTensors);
      \draw[->,draw=C0] (elasticityTensors) -- node[above] {%
        \footnotesize{}%
        $
          \tr{\cholfactor} \cholfactor = \etensor
        $\vphantom{p}%
      } (choleskyFactors);
      \draw[->,draw=C0] (choleskyFactors) -- node[above] {%
        \footnotesize{}interpolate%
      } (choleskyInterpolant);
    \end{tikzpicture}%
  }%
  \\[2mm]%
  \subcaptionbox{%
    Online phase (one iteration of the optimizer).%
  }[149mm]{%
    \begin{tikzpicture}
      \node[myCircle] (points) at (0mm,0mm) {%
        $
          \begin{matrix}
            \mcp{1},\\
            \dots,\\
            \mcp{M}
          \end{matrix}
        $%
      };
      \node[myCircle] (choleskyFactors) at (34mm,0mm) {%
        $
          \begin{matrix}
            \cholfactorintp(\mcp{1}),\\
            \dots,\\
            \cholfactorintp(\mcp{M})
          \end{matrix}
        $%
      };
      \node[myCircle] (elasticityTensors) at (83mm,0mm) {%
        $
          \begin{matrix}
            \etensorcholintp(\mcp{1}),\\
            \dots,\\
            \etensorcholintp(\mcp{M})
          \end{matrix}
        $%
      };
      \node[myCircle] (complianceValue) at (129.5mm,0mm) {%
        $
          \begin{matrix}
            \complianceintp(\mcp{1},\\
            \dotsc,\\
            \mcp{M})
          \end{matrix}
        $%
      };
      \draw[->,draw=C0] (points) -- node[above] {%
        \footnotesize{}evaluate\vphantom{p}%
      } (choleskyFactors);
      \draw[->,draw=C0] (choleskyFactors) -- node[above] {%
        \footnotesize{}%
        $
          \etensorcholintp
          = \tr{(\cholfactorintp)} \cholfactorintp
        $\vphantom{p}%
      } (elasticityTensors);
      \draw[->,draw=C0] (elasticityTensors) -- node[above] {%
        \footnotesize{}macro-problem%
      } (complianceValue);
    \end{tikzpicture}%
  }%
  \caption[Offline and online phase for topology optimization]{%
    Offline and online phase for topology optimization.
    The interpolation of the
    micro-cell density $\denscell$ with $\denscellintp$
    (see \cref{sec:622BSplines}) has been omitted for brevity.%
  }%
  \label{fig:topoOptPhases}%
\end{figure}

\paragraph{Generation of spatially adaptive sparse grids}

We use the classical surplus-based refinement criterion
(see, e.g., \cite{Pflueger10Spatially})
as shown in \cref{alg:topoOptGridGeneration}
to generate the spatially adaptive sparse grids.
The difference to common surrogate settings is that the objective function
$\*f\colon \clint{\*0, \*1} \to \real^{m+1}$ is vector-valued.
As the entries of $\cholfactor$ cannot be evaluated individually,
the adaptivity criterion has to consider all entries at once
to avoid performing unnecessary evaluations.
We use the surpluses in the piecewise linear hierarchical basis,
as their absolute values correlate with the second mixed derivative
of the objective function due to \cref{eq:surplusIntegral}.
The surpluses are combined using the formula
$\beta_k \ceq \tr{\*c} \vabs{\vsurplus_{\*l_k,\*i_k}}$
(with entry-wise absolute value) and the
points with largest $\beta_k$ are refined.

\begin{algorithm}
  \begin{algorithmic}[1]
    \Function{$\sgset = \texttt{offlinePhase}$}{%
      $\*\objfun$, $n$, $b$, $\*c$, $l_{\max}$, $\refinetol$,
      $\ngp_{\mathrm{refine}}$%
    }
      \State{$\sgset \gets \coarseregsgset{n}{d}{b}$}
      \Comment{initial regular sparse grid}%
      \While{\True}
        \State{$\ngp \gets \setsize{\sgset}$}
        \Comment{number of grid points}%
        \State{%
          Let $(\vsurplus_{\*l_{k'},\*i_{k'}})_{k' = 1, \dotsc, \ngp}$
          satisfy $
            \fa{k = 1, \dotsc, \ngp}{
              \sum_{k'=1}^{\ngp} \vsurplus_{\*l_{k'},\*i_{k'}}
              \bspl{\*l_{k'},\*i_{k'}}{1}(\gp{\*l_k,\*i_k})
              = \*\objfun(\gp{\*l_k,\*i_k})
            }
          $%
        }
        \ForOneLine{$k = 1, \dotsc, \ngp$}{%
          $\beta_k \gets \tr{\*c} \vabs{\vsurplus_{\*l_k,\*i_k}}$%
        }
        \Comment{combine surpluses to a scalar value}%
        \State{%
          $
            \liset^\ast \gets \{
              k = 1, \dotsc, \ngp \mid
              \ex{\gp{\*l,\*i} \notin \sgset}{
                \gp{\*l_k,\*i_k} \to \gp{\*l,\*i}
              },\,
              \norm[\infty]{\*l_k} < l_{\max},\,
              \abs{\beta_k} > \refinetol
            \}
          $%
        }
        \IfOneLine{$\liset^\ast = \emptyset$}{\Break{}}
        \Comment{stop when there are no refinable grid points left}%
        \State{%
          Refine $\le \ngp_{\mathrm{refine}}$ of the points
          $\{\gp{\*l_k,\*i_k} \in \sgset \mid k \in \liset^\ast\}$
          with largest $\beta_k$%
        }
      \EndWhile{}
    \EndFunction{}
  \end{algorithmic}
  \caption[%
    Generation of spatially adaptive sparse grids for topology optimization%
  ]{%
    Generation of spatially adaptive sparse grids for topology optimization.
    Inputs are
    the objective function $\*f\colon \clint{\*0, \*1} \to \real^{m+1}$
    (combination of the Cholesky factor of the elasticity tensor and
    the micro-cell density),
    the level $n \ge d$ and boundary parameter $b \in \nat$ of the
    initial regular sparse grid,
    the vector $\*c \in \real^{m+1}$ of coefficients with which the
    absolute values of the entries of the surpluses are combined,
    the maximal level $l_{\max} \in \nat$,
    the refinement threshold $\refinetol \in \posreal$, and
    the number $\ngp_{\mathrm{refine}} \in \nat$ of points to refine
    in each iteration.
    Output is the spatially adaptive sparse grid $\sgset$.%
  }%
  \label{alg:topoOptGridGeneration}%
\end{algorithm}

\paragraph{Parameter bounds}

In the micro-cell models presented in \cref{sec:631models},
extreme micro-cell parameters near zero or one may cause problems
with the resulting elasticity tensors.
For instance, many elasticity tensor entries corresponding to
the 2D cross model are discontinuous near the lines $x_1 = 1$ or $x_2 = 1$
\multicite{Huebner14Mehrdimensionale,Valentin14Hierarchische}.
This is due to the fact that the micro-cell is completely filled with material
on these lines,
independent of the other micro-cell parameter.
Similar issues occur for the other models and the shearing angles.
Hence, we have to restrict the range of the feasible micro-cell parameters,
i.e., the sparse grid points
$\*x = \gp{\*l_k,\*i_k}$ are still defined on the unit hyper-cube
$\clint{\*0, \*1}$,
but the actual micro-cell parameters $\xscaled$ are retrieved by an
affine transformation $\xscaled \ceq \*a + (\*b - \*a) \*x$.
For the models in \cref{sec:631models},
we restrict the bar widths to $\clint{0.01, 0.99}$ and
the shearing angles to $\clint{-0.35\pi, 0.35\pi}$.

\paragraph{Software, algorithms, and domain discretization}

The micro-problems and macro-problems were solved with the
\fem software package CFS++ \cite{Kaltenbacher10Advanced}.%
\footnote{%
  \url{http://www.lse.uni-erlangen.de/cfs/}%
}
The micro-prob\-lems were discretized by dividing the micro-cells into
$128 \times 128 = \num{16384}$ elements (models in two dimensions) or
$16 \times 16 \times 16 = \num{4096}$ elements (models in three dimensions).
The macro-domains $\objdomain$ were discretized using
32 macro-cells per meter in the 2D cantilever scenario
(i.e., $64 \times 32 = \num{2048}$ cells),
20 macro-cells per meter in the 3D cantilever scenario
(i.e., $20 \times 20 \times 20 = \num{8000}$ cells), and
10 macro-cells per meter in the other scenarios
(i.e.,
\num{1600} cells for the 2D L-shape and
\num{4000} cells for the 3D center-load).
The generation of the sparse grids (offline phase) was done via a MATLAB code,
while the evaluation of the interpolants (online phase) was performed
by the sparse grid toolbox \sgpp \cite{Pflueger10Spatially}.%
\footnote{%
  \url{http://sgpp.sparsegrids.org/}%
}
For the solution of the emerging optimization problems,
a sequential quadratic programming method was employed
(see \cref{sec:513gradientBasedConstrained}).

\subsection{Error Sources}
\label{sec:642errorSources}

There are multiple sources that contribute to the numerical error
of our method:

\begin{enumerate}[label=E\arabic*.,ref=E\arabic*,leftmargin=2.7em]
  \item
  \label{item:topoOptErrorMicro}
  Discretization of the micro-problem
  (i.e., the elasticity tensors $\etensor$ are inaccurate)
  
  \item
  \label{item:topoOptErrorInterpolation}
  Sparse grid interpolation
  (i.e., $\etensorintp \not= \etensor$)
  
  \item
  \label{item:topoOptErrorCholesky}
  Reconstruction of elasticity tensors with Cholesky factors
  (i.e., $\etensorcholintp \not= \etensorintp$)
  
  \item
  \label{item:topoOptErrorMacro}
  Discretization of the macro-problem
  (i.e., the compliance $\compliance$ is inaccurate)
  
  \item
  \label{item:topoOptErrorOptimization}
  Optimization
  (i.e., the minimum found by the optimizer is inaccurate or not global)
  
  \item
  \label{item:topoOptErrorRounding}
  Floating-point rounding errors
  (i.e., arithmetical operations are inaccurate)
\end{enumerate}

\noindent
\ref{item:topoOptErrorRounding}-type errors are always present and
will not be analyzed in this chapter.
Errors of type \ref{item:topoOptErrorMicro} and \ref{item:topoOptErrorMacro}
are intrinsic to the homogenization approach
and will not be discussed here either.
The optimization error \ref{item:topoOptErrorOptimization}
has already been discussed in \cref{sec:542optimization}
for explicit test functions.
Therefore, in the remainder of this chapter,
we will focus on the analysis of the errors of types
\ref{item:topoOptErrorInterpolation} and \ref{item:topoOptErrorCholesky},
since the interpolation of Cholesky factors is the
major new contribution to this application.

\subsection{Interpolation Error}
\label{sec:643interpolation}

\paragraph{Spectral interpolation error measure}

For the interpolation error \ref{item:topoOptErrorInterpolation} and
the Cholesky factorization error \ref{item:topoOptErrorCholesky},
we cannot simply take the absolute value of the difference
of the objective function $\*\objfun\colon \clint{\*0, \*1} \to \real^{m+1}$
and its surrogate $\vsgintp$, since both are vector-valued.
As the micro-cell density $\denscell$
is not affected by the Cholesky factorization,
we consider only the elasticity tensor
$\etensor\colon \clint{\*0, \*1} \to \real^{6 \times 6}$ and
its surrogate
$\etensorcholintp\colon \clint{\*0, \*1} \to \real^{6 \times 6}$
obtained by Cholesky factorization.
To retrieve a scalar error measure,
we use the spectral norm
\begin{equation}
  \norm[2]{\etensor(\*x) - \etensorcholintp(\*x)},\quad
  \*x \in \clint{\*0, \*1},
\end{equation}
i.e., the largest absolute eigenvalue of
$\etensor(\*x) - \etensorcholintp(\*x)$.
However, the choice of the norm is arbitrary,
as all matrix norms on $\real^{6 \times 6}$ are equivalent to each other.

\paragraph{Pointwise spectral interpolation error}

\begin{figure}
  \subcaptionbox{%
    $\norm[2]{\etensor(\*x) - \etensorintp(\*x)}$%
    \label{fig:topoOptInterpolationErrorPointwise_1}%
  }[63mm]{%
    \includegraphics{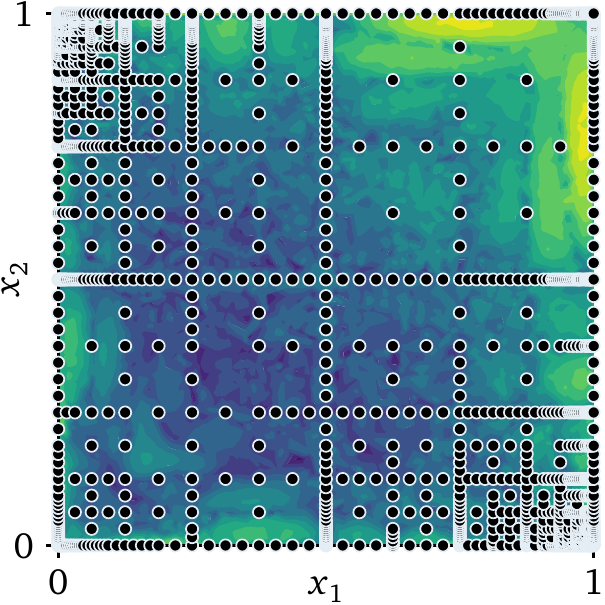}%
  }%
  \hspace{3mm}%
  \subcaptionbox{%
    $\norm[2]{\etensor(\*x) - \etensorcholintp(\*x)}$%
    \label{fig:topoOptInterpolationErrorPointwise_2}%
  }[63mm]{%
    \includegraphics{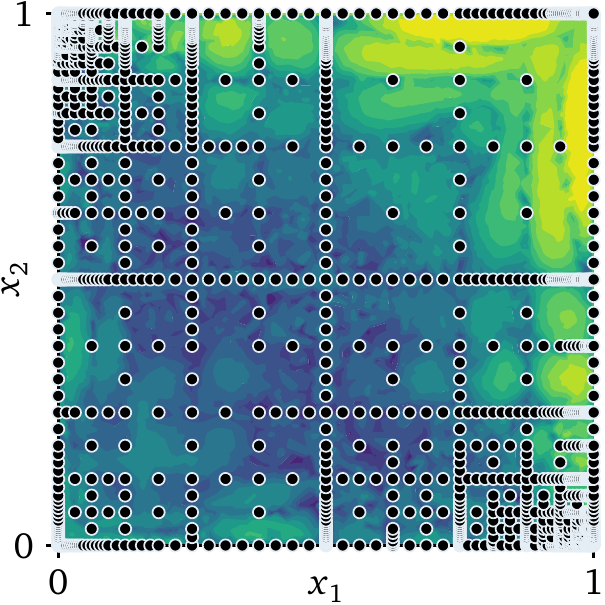}%
  }%
  \hfill%
  \includegraphics{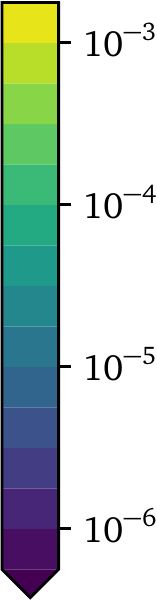}%
  \caption[Pointwise spectral interpolation error for the 2D cross model]{%
    Pointwise spectral interpolation error for the 2D cross model and
    cubic B-splines on
    $\ngp = 1320$ spatially adaptive sparse grid points \emph{(dots)} for
    the direct elasticity tensor interpolation \emph{(left)} and
    the Cholesky factor interpolation \emph{(right).}%
  }%
  \label{fig:topoOptInterpolationErrorPointwise}%
\end{figure}

\Cref{fig:topoOptInterpolationErrorPointwise}
shows the pointwise spectral interpolation error for the 2D cross model
and the corresponding spatially adaptive sparse grid
generated with the refinement algorithm as explained in
\cref{sec:641implementation}.
The above-mentioned discontinuity of elasticity tensor entries
near $x_1 = 1$ or $x_2 = 1$
is most severe near the corners $\*x \in \{(0, 1), (1, 0)\}$
(cf.\ \cref{fig:cholesky}),
as some entries vanish if one of the micro-cell bars has zero width.
Hence, most points are placed near these singularity corners.

The left plot (\cref{fig:topoOptInterpolationErrorPointwise_1})
shows the spectral interpolation error
$\norm[2]{\etensor(\*x) - \etensorintp(\*x)}$
of the direct elasticity tensor interpolant without Cholesky factorization
(i.e., error \ref{item:topoOptErrorInterpolation}).
The maximum error is \num{1.2e-3},
which is attained near the critical lines $x_1 = 1$ or $x_2 = 1$.
Note that the mean error over the whole domain $\clint{\*0, \*1}$
is only \num{4.5e-5}.
In the right plot (\cref{fig:topoOptInterpolationErrorPointwise_2}),
the picture changes slightly when looking at the spectral interpolation error
$\norm[2]{\etensor(\*x) - \etensorcholintp(\*x)}$
of the elasticity tensor resulting from Cholesky factorization
(i.e., errors \ref{item:topoOptErrorInterpolation} and
\ref{item:topoOptErrorCholesky} combined).
The maximum error becomes \num{3.4e-3},
while the mean error increases to \num{1.1e-4}.
We conclude that the Cholesky factorization leads to an increase
of interpolation errors by only less than half an order of magnitude.

\paragraph{Convergence of spectral interpolation error}

\begin{figure}
  \hspace*{5mm}%
  \includegraphics{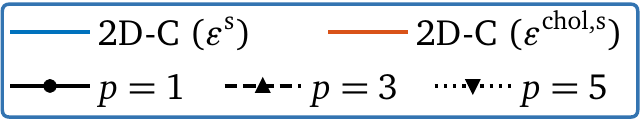}%
  \hfill%
  \raisebox{0.5mm}{\includegraphics{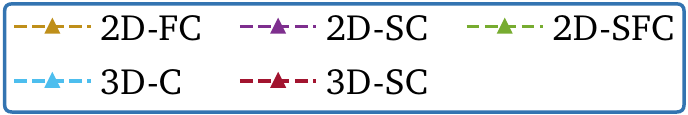}}%
  \\[2mm]%
  \subcaptionbox{%
    $\error^{\sparse}$ \emph{\textcolor{C0}{(blue)}} and
    $\error^{\chol,\sparse}$ \emph{\textcolor{C1}{(red)}}
    for the 2D cross model and different degrees.%
    \label{fig:topoOptInterpolationErrorBasisFunctions_1}%
  }[72mm]{%
    \includegraphics{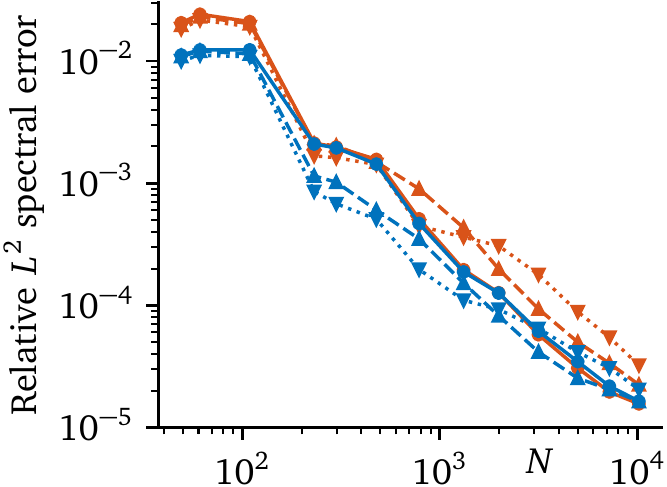}%
  }%
  \hfill%
  \subcaptionbox{%
    $\error^{\chol,\sparse}$
    for the other models and $p = 3$.%
    \label{fig:topoOptInterpolationErrorBasisFunctions_2}%
  }[72mm]{%
    \includegraphics{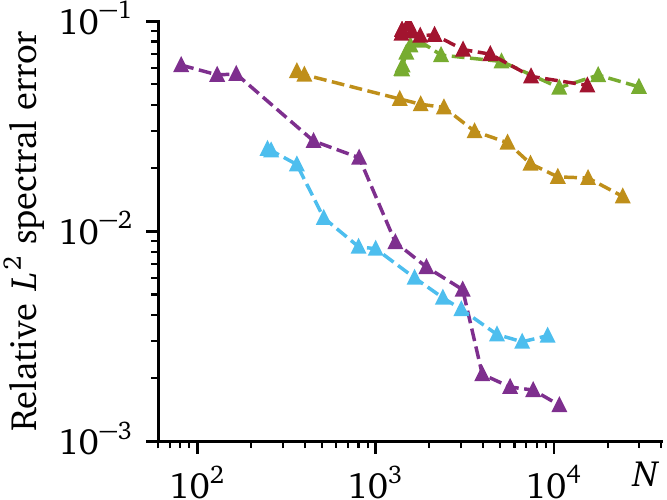}%
  }%
  \caption[Convergence of relative $L^2$ spectral interpolation errors]{%
    Convergence of relative $\Ltwo$ spectral interpolation errors
    over the increasing number $\ngp$ of spatially adaptive grid points
    (i.e., decreasing threshold $\refinetol$)
    for the 2D cross model without or with Cholesky factor interpolation
    and different degrees $p$ \emph{(left)} and
    for the other models and cubic degree \emph{(right)}.%
  }%
  \label{fig:topoOptInterpolationErrorBasisFunctions}%
\end{figure}

\Cref{fig:topoOptInterpolationErrorBasisFunctions_1} shows
the convergence of the relative $\Ltwo$ spectral interpolation errors
\begin{equation}
  \error^{\sparse} \ceq
  \frac{
    \normLtwoscaled{
      \vphantom{\big(}
      \norm[2]{\etensor({\cdot}) - \etensorintp({\cdot})}
    }
  }{
    \normLtwoscaled{
      \vphantom{\big(}
      \norm[2]{\etensor({\cdot})}
    }
  }, \qquad
  \error^{\chol,\sparse} \ceq
  \frac{
    \normLtwoscaled{
      \vphantom{\big(}
      \norm[2]{\etensor({\cdot}) - \etensorcholintp({\cdot})}
    }
  }{
    \normLtwoscaled{
      \vphantom{\big(}
      \norm[2]{\etensor({\cdot})}
    }
  }
\end{equation}
for the 2D cross model, i.e.,
the relative $\Ltwo$ error of the functions depicted in
\cref{fig:topoOptInterpolationErrorPointwise}.
Relative errors of \SI{1}{\permille} are already obtained
for $\ngp = 200$ grid points.
Unfortunately, even for higher B-spline degrees $p > 1$,
the order of convergence is only quadratic
due to the singularities of the elasticity tensor.
This slow convergence does not improve for the other
micro-cell models as shown in
\Cref{fig:topoOptInterpolationErrorBasisFunctions_2}.
In fact, the convergence decelerates even more
as the number of micro-cell parameters increases.
For the 2D sheared cross and 3D cross models with three parameters,
the spatially adaptive sparse grid with $\ngp \approx \num{10000}$ grid points
is able to achieve a relative error of around \SI{3}{\permille}.
However, for the 2D sheared framed cross and 3D sheared cross models
with five parameters, only errors of about \SI{5}{\percent} are reached
for the same grid size.

\subsection{Optimal Compliance Values and Structures}
\label{sec:644optimization}

\paragraph{Optimal compliance values for different micro-cell models}

In the following, we use for each micro-cell model
a specific spatially adaptive sparse grid with around \num{10000} points.
The exact grid sizes and other details about the employed sparse grids
can be found in \cref{tbl:topoOptModels}
(located in \cref{chap:a30topoOptDetails}).
For hierarchical cubic B-splines ($p = 3$),
\cref{tbl:topoOptResultsModels} lists the
compliance values $\compliance(\mcpoptappr{1}, \dotsc, \mcpoptappr{M})$
for each of the four scenarios
and the corresponding possible micro-cell models,
where $(\mcpoptappr{1}, \dotsc, \mcpoptappr{M}) \in (\real^d)^M$
is the micro-cell parameter combination that is returned by the optimizer.%
\footnote{%
  Note that this true compliance value differs from the
  approximated value
  $\complianceintp(\mcpoptappr{1}, \dotsc, \mcpoptappr{M})$,
  which the optimizer reports as the optimal objective function value.%
}
It is obvious that more complicated micro-cell models
lead to lower (better) compliance values,
as they are a generalization of the simple models.
For instance, the 2D cross is a special case of
the 2D framed cross, the 2D sheared cross, and the 2D sheared framed cross.
By choosing the respectively best model for each scenario,
we are able to decrease the compliance value
(and, hence, increase the stability of the resulting structure)
by \SI{9.6}{\percent} in the 2D cantilever scenario,
by \SI{7.7}{\percent} in the 2D L-shape scenario,
by \SI{34}{\percent} in the 3D cantilever scenario, and
by \SI{73}{\percent} in the 3D center-load scenario.
In general, this motivates the usage of more complicated micro-cell models,
which cannot be computationally handled with
conventional full grid interpolation methods.
Consequently, sparse grids or similar methods have to be used.

\begin{table}
  \setnumberoftableheaderrows{1}%
  \begin{tabular}{%
    >{\kern\tabcolsep}=l<{\kern5mm}*{6}{+c}<{\kern\tabcolsep}%
  }
    \toprulec
    \headerrow
    Scenario&       2D-C&   2D-FC&  2D-SC&           2D-SFC& 3D-C&   3D-SC\\
    \midrulec
    2D cantilever&  74.974& 70.816& \textbf{67.809}& 68.602& ---&    ---\\
    2D L-shape&     183.68& 177.51& \textbf{169.60}& 174.55& ---&    ---\\
    \midrulec
    3D cantilever&  ---&    ---&    ---&             ---&    247.60& \textbf{162.59}\\
    3D center-load& ---&    ---&    ---&             ---&    169.27& \textbf{46.171}\\
    \bottomrulec
  \end{tabular}
  \caption[Optimal compliance values for different micro-cell models]{%
    Optimal compliance values for the different scenarios
    and micro-cell models using cubic B-splines
    (spatially adaptive grids with around \num{10000} points).
    The entries highlighted in \textbf{bold face} indicate the best choice
    of micro-cell models for a given scenario.
    More details can be found in \cref{tbl:topoOptResultsDetailed}.%
  }%
  \label{tbl:topoOptResultsModels}%
\end{table}

\paragraph{Corresponding optimal structures}

\begin{figure}
  \subcaptionbox{%
    2D cross%
  }[72mm]{%
    \includegraphics{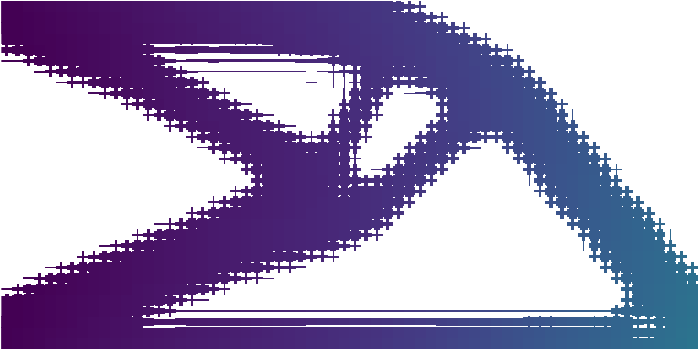}%
  }%
  \hfill%
  \subcaptionbox{%
    2D framed cross%
  }[72mm]{%
    \includegraphics{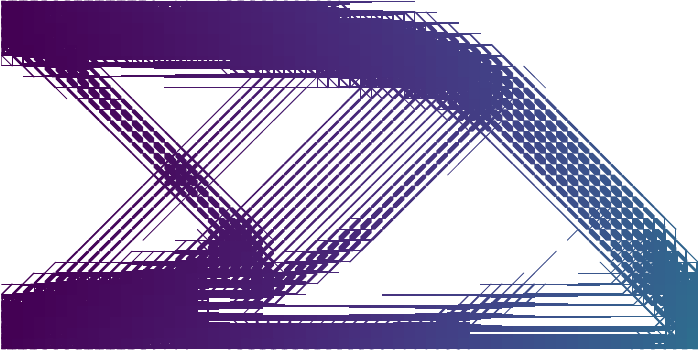}%
  }%
  \\[2mm]%
  \subcaptionbox{%
    2D sheared cross%
  }[72mm]{%
    \includegraphics{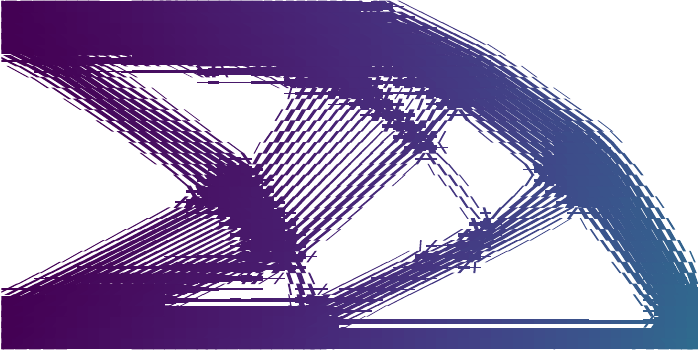}%
  }%
  \hfill%
  \subcaptionbox{%
    2D sheared framed cross%
  }[72mm]{%
    \includegraphics{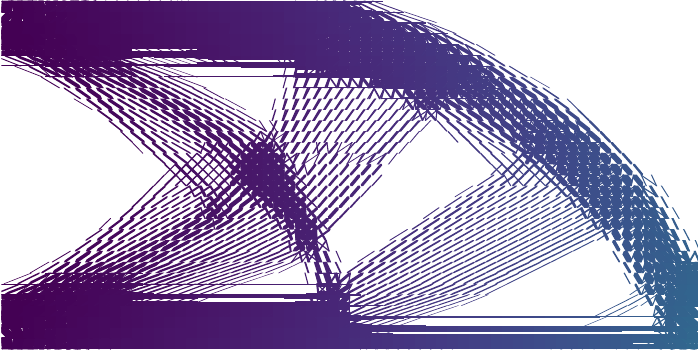}%
  }%
  \caption[Optimal structures in the 2D cantilever scenario]{%
    Topologically optimal structures in the 2D cantilever scenario
    for different micro-cell models using cubic B-splines
    (spatially adaptive grids with around \num{10000} points).
    The colors indicate the length of the displacement,
    where dark regions correspond to weak displacements and
    bright regions to strong displacements.
    The color map is the same as in
    \cref{fig:topoOptStructure2DLShape}.
    Only bars with widths $\ge 0.1$ are shown.
    More details can be found in \cref{tbl:topoOptResultsDetailed}.%
  }%
  \label{fig:topoOptStructure2DCantilever}%
\end{figure}

The corresponding optimal structures are shown in
\cref{fig:topoOptStructure2DCantilever} for the 2D cantilever scenario
and, for reasons of space, in \cref{chap:a30topoOptDetails} in
\cref{fig:topoOptStructure2DLShape,fig:topoOptStructure3D}
for the other three scenarios.
Of course, the periodic micro-cell structures cannot be plotted directly,
as the micro-cells are infinitesimally small.
Therefore, the figures show for each macro-cell only
one single large micro-cell.

Two effects can be seen in the plots of the optimal structures:
First, the simpler models are not able to direct the emerging forces at
arbitrary angles.
For example, the 2D framed cross model strongly prefers
angles of \ang{45}, which results in structures that are not as stable
as they could be.
The 2D sheared cross and 2D sheared framed cross models
are considerably more flexible, allowing
internal forces to act at almost arbitrary angles.
Second, the sheared micro-cell models use the available
material volume more efficiently than the cross model.
This is most striking in the 3D case (see \cref{fig:topoOptStructure3D}),
where it seems that the sheared cross structures
use more volume than the simple cross structures,
although the structures spend exactly the same amount of material volume.
The reason is that for the cross model, both bars
have to be used in order to connect the macro-cell to its neighbors.
For the sheared cross model, a shearing of the vertical bar suffices,
and we save volume by not using the horizontal bar.
Both of these effects explain the significantly lower compliance
values for the sheared micro-cell models.

\paragraph{Comparison to the direct solution}

B-splines on sparse grids lead to a drastic reduction in computation time.
Solving the 2D cantilever scenario with the best-placed sheared cross model
would take 453 days with
exact elasticity tensor evaluations (i.e., without surrogates),
assuming the same number of iterations as for the surrogate tensor case
and sequential computation of the elasticity tensors.
This estimate does not account for approximating the missing derivatives
of the elasticity tensor.
If we incorporate this and use 100 parallel processes,
we still need weeks for the solution.
In contrast, the computation time using our sparse grid surrogates
is a matter of minutes or hours at most,
resulting in speedups of around 200.
This is excluding the time for the offline phase,
which is in the range of hours, but which has to be spent only once,
as the resulting grid can be reused for different scenarios.

\vspace{2em}

\paragraph{Optimality-interpolation gaps}

\begin{figure}
  \includegraphics{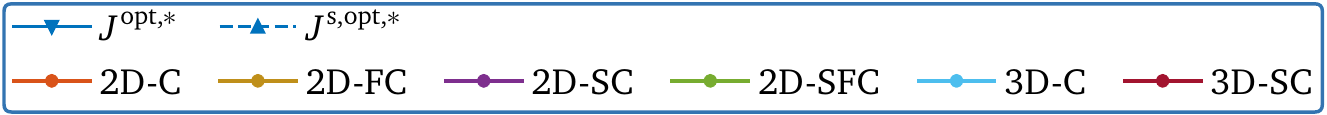}%
  \\[2mm]%
  \includegraphics{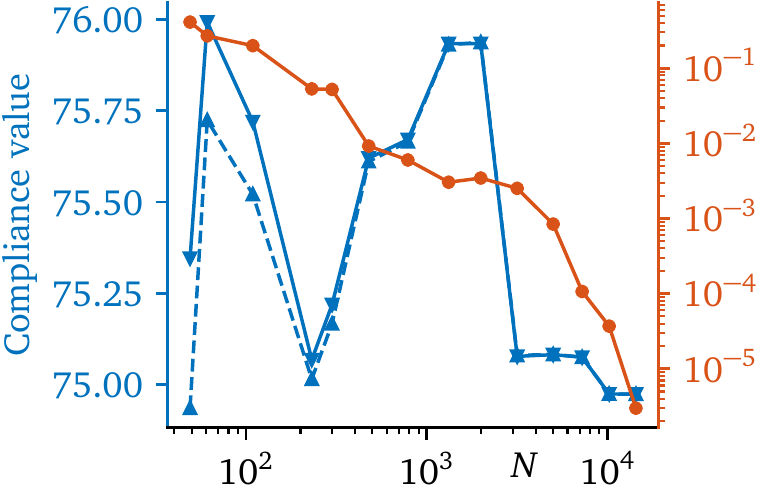}%
  \hfill%
  \includegraphics{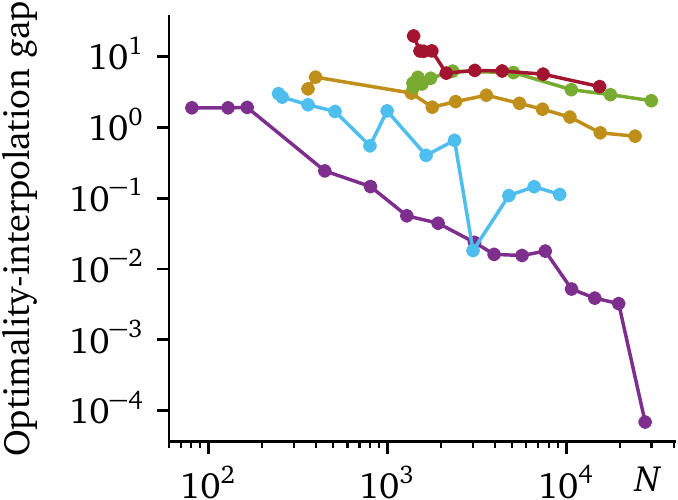}%
  \caption[Convergence of the optimality-interpolation gap]{%
    Convergence of the optimality-interpolation gap
    $\abs{\compliance[\opt,\ast] - \complianceintp[\opt,\ast]}$
    for the 2D/3D cantilever scenario
    and different micro-cell models using cubic B-splines ($p = 3$).
    The left plot additionally shows
    $\compliance[\opt,\ast]
    \ceq \compliance(\mcpoptappr{1}, \dotsc, \mcpoptappr{M})$ and
    $\complianceintp[\opt,\ast]
    \ceq \complianceintp(\mcpoptappr{1}, \dotsc, \mcpoptappr{M})$
    for the 2D cross model.%
  }%
  \label{fig:topoOptOptimalityGap}%
\end{figure}

Ideally, we would measure the true optimality gap
\begin{equation}
  \label{eq:topoOptOptimalityGapTrue}
  \compliance(\mcpoptappr{1}, \dotsc, \mcpoptappr{M}) -
  \compliance(\mcpopt{1}, \dotsc, \mcpopt{M}),
\end{equation}
cf.\ error \ref{item:topoOptErrorOptimization}.
Unfortunately, the true optimum
$(\mcpopt{1}, \dotsc, \mcpopt{M})$ could not be computed:
Apart from the time issue mentioned above,
oscillations in the elasticity tensor evaluation and
errors stemming from types
\ref{item:topoOptErrorMicro} and \ref{item:topoOptErrorRounding}
reliably led to optimizer crashes as it ran into discontinuities,
which are smoothed out when using B-spline surrogates.
However, as in \cref{fig:topoOptOptimalityGap},
we can at least calculate the \term{optimality-interpolation gap}
\begin{equation}
  \label{eq:topoOptOptimalityGapPseudo}
  \abs{
    \compliance(\mcpoptappr{1}, \dotsc, \mcpoptappr{M}) -
    \complianceintp(\mcpoptappr{1}, \dotsc, \mcpoptappr{M})
  }
\end{equation}
between the actual compliance value
and the approximated, reported compliance value.
This gap does not constitute any kind of bound on the true optimality gap;
however, the idea is that
as the interpolation error converges to zero,
the optimality-interpolation gap should converge to zero, too.

\pagebreak

\Cref{fig:topoOptOptimalityGap} (left) shows that for the 2D cross model,
the optimizer reports compliance values that are smaller than in reality
($\complianceintp[\opt,\ast]$ vs. $\compliance[\opt,\ast]$).
However, the difference steadily converges to zero.
This is similar for the other micro-cell model as shown in the
right part of \cref{fig:topoOptOptimalityGap},
although the convergence is much slower due to the
higher number $d$ of micro-cell parameters.

\paragraph{Optimal compliance values for different B-spline degrees}

Finally, to study the effect of the B-spline degree on the
optimization performance,
\cref{tbl:topoOptResultsDegrees} lists the compliance values
for the degrees $p = 1, 3, 5$ and the 2D/3D cross and sheared cross
micro-cell models.
In the two-dimensional scenarios,
higher-order B-splines decrease the compliance value
by up to \SI{9}{\percent}.
In the three-dimensional scenarios,
higher-order B-splines may perform worse than the piecewise linear
functions ($p = 1$).
(However, as indicated in \cref{tbl:topoOptResultsDegrees},
all optimization runs with piecewise linear functions
terminated prematurely due to numerical difficulties with the
discontinuous derivatives.)
It may be suspected that if we used micro-cell models with
less prominent discontinuities (i.e., ``smoother'' elasticity tensors),
the advantage of higher-order B-splines would be more visible.
All in all, the application of topology optimization underlines
that good interpolation (and thus a good quality of the surrogate)
is key to good optimization results.

\begin{table}
  \setnumberoftableheaderrows{2}%
  \begin{tabular}{%
      >{\kern\tabcolsep}=l<{\kern5mm}*{3}{+c}%
      <{\kern5mm}*{3}{+c}<{\kern\tabcolsep}%
    }
    \toprulec
    \headerrow
    &
    \multicolumn{3}{c}{\hspace*{-12pt}2D/3D cross}&
    \multicolumn{3}{c}{\hspace*{-6pt}2D/3D sheared cross}\\
    \headerrow
    Scenario&       $p = 1$&         $p = 3$&         $p = 5$&          $p = 1$&                $p = 3$&         $p = 5$\\
    \midrulec
    2D cantilever&  \emph{82.365}&   \textbf{74.974}& 76.070&           \emph{68.889}&          \textbf{67.809}& 68.018\\
    2D L-shape&     \emph{193.83}&   \textbf{183.68}& 183.70&           \emph{169.85}&          169.60&          \textbf{169.60}\\
    \midrulec
    3D cantilever&  \emph{249.75}&   247.60&          \textbf{247.33}&  \emph{\textbf{148.72}}& 162.59&          152.34\\
    3D center-load& \textbf{162.68}& 169.27&          163.94&           \emph{\textbf{45.713}}& 46.171&          47.074\\
    \bottomrulec
  \end{tabular}
  \caption[Optimal compliance values for different B-spline degrees]{%
    Optimal compliance values for the different scenarios
    and B-spline degrees using the 2D/3D cross micro-cell model \emph{(left)}
    and the 2D/3D sheared cross micro-cell model \emph{(right).}
    The spatially adaptive sparse grids are the same as in
    \cref{tbl:topoOptResultsModels}.
    The entries highlighted in \textbf{bold face} indicate the best choice
    of B-spline degree for a given scenario and micro-cell model.
    Optimization runs of entries marked as \emph{italic}
    terminated prior to success due to numerical difficulties.%
  }%
  \label{tbl:topoOptResultsDegrees}%
\end{table}

\cleardoublepage

  \setdictum{%
  Beware of bugs in the above code;\\
  I have only proved it correct, not tried it.%
}{%
  Donald E.\ Knuth \cite{Knuth77Notes}%
}

\longchapter{%
  Application 2: Musculoskeletal Models%
}{%
  Application 2:\texorpdfstring{\\}{ }Musculoskeletal Models%
}{%
  Application 2 -- Musculoskeletal Models%
}
\label{chap:70muscle}

\initial{0em}{E}{xisting musculoskeletal models of muscle-tendon complexes,}
e.g., of the human upper limb, can mainly be divided into two different types.
\term{Lumped-parameter musculoskeletal models,}
for example Hill-type models based on multi-body simulations
\multicite{Roehrle16Two,Valentin18Gradient},
constitute the most common type.
These models assume that the components of the musculoskeletal system
are rigid.
The mechanics is reduced to point masses
associated with their moment of inertia;
thus, these models can be described by few parameters.

\term{Continuum-mechanical musculoskeletal models} form the other type.
Their advantage is that they are more detailed and, hence, more realistic.
However, their increased complexity leads to higher computational costs.
As an example, we consider an inverse problem (see \cref{chap:10introduction})
that involves a continuum-mechanical simulation of such a
musculoskeletal model,
where we search values of model parameters
such that a specific movement is attained.
Each iteration of the solution process for such an inverse problem
may take hours or even days, depending on the model at hand.

Surrogate methods based on sparse grids help to decrease the complexity
in two ways:
First, the evaluation of surrogates is obviously drastically cheaper
than the solution of the inverse problem.
Second, the particular choice of sparse grids decreases the number
of necessary samples to construct the surrogates.
As for the previous application,
the choice of B-splines as hierarchical basis functions enables
the evaluation of continuously differentiable surrogate gradients.
For the example of inverse problems, this means that
gradient-based optimization methods may be employed,
which significantly accelerates convergence compared to gradient-free methods.

This chapter is split into three sections.
First, in \cref{sec:71model}, we introduce a continuum-mechanical
model of the human upper limb.
Second, in \cref{sec:72methodology}, we list the types of inverse problems
of interest.
Finally, in \cref{sec:73results}, we present numerical results
regarding the solution of these inverse problems.

The results of this chapter are based on a collaboration with
Prof.\ Oliver Röhrle, PhD, and Dr.\ Michael Sprenger
(both SimTech/University of Stuttgart, Germany).%
\footnote{%
  Michael Sprenger left the University of Stuttgart in 2015.%
}
The collaborators contributed the biomechanical model
with its theory, its geometry, and its implementation,
while the author of this thesis contributed the
sparse grid/B-spline methodology and
computed the numerical results.
Note that the results have already been
published in a paper \cite{Valentin18Gradient},
which we will follow closely in this chapter.

\section{Continuum-Mechanical Model of the Upper Limb}
\label{sec:71model}

\minitoc{71mm}{6}

\noindent
In the following, we first discuss the state of the art
in biomechanical modeling.
Then, we address details of the model of the human upper limb.
For convenience, the most relevant symbols are listed in
\cref{tbl:glossaryMusculoskeletal}.

\begin{table}
  \setnumberoftableheaderrows{0}%
  \newcommand*{\pnst}[1]{%
    \printnotationsymbol{#1}%
    \vphantom{\printnotationsymbol{\equielbang}$(\cdot)^{\reference}$}&%
    \printnotationtext{#1}%
  }%
  \begin{tabular}{%
      >{\kern\tabcolsep}=l+l<{\kern5mm}+l+l<{\kern5mm}+l+l<{\kern\tabcolsep}%
    }
    \toprulec
    \pnst{\forceT}& \pnst{\armT}&       \pnst{\actT}\\
    \pnst{\forceB}& \pnst{\armB}&       \pnst{\actB}\\
    \pnst{\forceL}& \pnst{\armL}&       \pnst{\moment}\\
    \pnst{\elbang}& \pnst{\tarelbang}&  \pnst{\equielbang}\\
    \pnst{t}&
    $(\cdot)^{\sparse}$&Sparse grid solution&
    $(\cdot)^{\reference}$&Reference solution\\
    \bottomrulec
  \end{tabular}%
  \caption[Glossary for musculoskeletal models]{%
    Glossary of the notation for musculoskeletal models.%
  }%
  \label{tbl:glossaryMusculoskeletal}%
\end{table}

\subsection{Continuum-Mechanical Musculoskeletal Models}
\label{sec:711models}

\paragraph{%
  Limitations of classical models and
  benefits of continuum-mechanical models%
}

Due to the simplicity of classical lumped-parameter models,
their degree of realism is limited.
Without any modifications,
lumped-parameter models are not able to represent
detailed heterogeneous material characteristics or non-trivial
muscle force paths \cite{Roehrle16Two}.

The exploitation of continuum-mechani\-cal constitutive laws
for musculoskeletal models is a more recent development \cite{Roehrle16Two}.
The resulting models are capable to model spatial quantities
such as complex muscle fiber field architectures,
local activation principles, complex muscle geometries, or contact mechanics
\multicite{Roehrle16Two,Valentin18Gradient}.
Most of the existing work only treats single skeletal muscles in isolation
\multicite{Lemos05Modeling,Sharafi11Strains,Heidlauf14Multiscale}.
The model used in this thesis
(which is the same model as in
\multicite{Sprenger15Continuum,Roehrle16Two,Valentin18Gradient})
aims at studying the interplay of multiple muscles and bones.

\paragraph{Overdetermined antagonistic systems}

Musculoskeletal systems are typically overdetermined \cite{Roehrle16Two}.
This means that the number of muscles that act on a specific joint
is usually larger than the number of the joint's degrees of freedom.
For instance, in a simple model of the human upper limb,
there are two antagonistic muscles
(i.e., muscles that work against each other), namely triceps and biceps,
but only one joint angle at the elbow.
Mathematically speaking, a single muscle suffices to attain
a large range of elbow angles
that are possible with an antagonistic muscle pair.
However, the usage of two muscles enables faster movements and
allows for abrupt changes of direction.

The overdetermination of most musculoskeletal models
implies that additional conditions have to be imposed in order
to obtain unique solutions.
There exist various types of such conditions,
for instance, minimal control effort, minimal control change, and
minimal kinematic energy \cite{Valentin18Gradient}.
The idea behind these conditions is that the human body
tries to minimize the energy effort that is associated with
all types of motion.

\paragraph{Forward and inverse dynamics}

Musculoskeletal simulations are usually based on
either forward dynamics or inverse dynamics \cite{Valentin18Gradient}.
\term{Forward-dynamic approaches} use activation parameters
for the muscles as the input and simulate the corresponding motion
as the output.
This requires that we know the muscle forces
(depending on the activation levels) beforehand.
For example, one can prescribe activation levels of facial muscles
to achieve specific facial expressions \cite{Wu13Modelling}.
In contrast, \term{inverse-dynamic approaches}
use experimental motion data as the input
to estimate the muscle forces as the output \cite{Roehrle16Two}.
With inverse-dynamic simulations,
one can investigate the wrapping of muscles
around the knee joint \cite{Fernandez05Anatomically} or
visualize the motion of skin \cite{Lee09Comprehensive}, for instance.

\subsection{Details of the Human Upper Limb Model}
\label{sec:712details}

As shown in \cref{fig:raisingArm},
our model of the human upper limb consists
of the three bones humerus, ulna, and radius,
the elbow joint with one degree of freedom, and
the antagonistic muscle pair of triceps brachii and biceps brachii
\cite{Valentin18Gradient}.
The bones are rigid bodies and the muscle-tendon complex
is simulated with a continuum-mechanical approach.
This implies that the muscles deform when they contract.

\begin{figure}
  \begin{tikzpicture}[
    dashed/.style={dash pattern=on 2pt off 2pt},
    contour/.style={line width=3pt,draw=mittelblau!10},
  ]
    \newcommand*{\myhspace}{30.34mm}
    \newcommand*{\myscale}{0.14}
    \newcommand*{\mytricepsarmstart}{0.1}
    \newcommand*{\mybicepsstart}{0.2}
    \newcommand*{\myradius}{9mm}
    \newcommand*{\mysmallradius}{2mm}
    
    \foreach \myi in {1,...,5} {
      \node[anchor=north west] at ({(\myi-1)*\myhspace},0mm) {%
        \includegraphics[scale=\myscale]{upperLimb_\myi}%
      };
      \node[anchor=north west] at ({(\myi-1)*\myhspace+11mm},0mm) {%
        \pgfmathparse{10+35*(\myi-1)}%
        $\theta = \ang[
          round-mode=places,round-precision=0,
        ]{\pgfmathresult}$%
      };
    }
    
    \begin{scope}[shift={({(3-1)*\myhspace},0mm)}]
      \coordinate (joint)               at (8mm,-30mm);
      \coordinate (loadForceStart)      at (30.6mm,-35mm);
      \coordinate (loadForceEnd)        at (30.6mm,-45mm);
      \coordinate (tricepsForceStart)   at (5.7mm,-31.9mm);
      \coordinate (tricepsForceEnd)     at (4mm,-15mm);
      \coordinate (tricepsArmStart)     at (
        ${(1-\mytricepsarmstart)}*(joint) +
        \mytricepsarmstart*(loadForceStart)$
      );
      \coordinate (bicepsForceStart)    at (
        ${(1-\mybicepsstart)}*(joint) + \mybicepsstart*(loadForceStart)$
      );
      \coordinate (bicepsForceEnd)      at (11.5mm,-10mm);
      \coordinate (horizontalDashedEnd) at ($(loadForceEnd) + (0mm,3mm)$);
      \coordinate (verticalDashedEnd)   at ($(joint) - (0mm,15mm)$);
      \pgfmathanglebetweenpoints{
        \pgfpointanchor{joint}{center}
      }{
        \pgfpointanchor{loadForceStart}{center}
      }
      \let\myelbowangle\pgfmathresult
      \coordinate (arcStart) at ($(joint) - (0mm,\myradius)$);
      \coordinate (arcEnd) at (
        $(joint) + (
          {\myradius*cos(\myelbowangle)},
          {\myradius*sin(\myelbowangle)}
        )$
      );
      
      \makeatletter
      \newcommand*{\drawleverarm}[5][\@nil]{
        \def\myoptarg{#1}
        \ifx\myoptarg\@nnil
          \def\mypointsize{1pt}
          \def\myfill{C1}
          \def\myoptarg{}
        \else
          \def\mypointsize{2pt}
          \def\myfill{mittelblau!10}
        \fi
        \coordinate (leverArmBase) at ($(#3)!(#2)!(#4)$);
        \draw[dashed,draw=C1,\myoptarg] (#2) -- (leverArmBase);
        \pgfmathanglebetweenpoints{
          \pgfpointanchor{#2}{center}
        }{
          \pgfpointanchor{leverArmBase}{center}
        }
        \let\myangle\pgfmathresult
        \pgfmathparse{\myangle #5 180}
        \let\mystartangle\pgfmathresult
        \pgfmathparse{\myangle #5 90}
        \let\myendangle\pgfmathresult
        \centerarc[draw=C1,\myoptarg](leverArmBase)(
          \mystartangle:\myendangle:\mysmallradius
        );
        \fill[draw=none,fill=\myfill] (
          $(leverArmBase) + (
            {0.62*\mysmallradius*cos((\mystartangle+\myendangle)/2)},
            {0.62*\mysmallradius*sin((\mystartangle+\myendangle)/2)}
          )$
        ) circle (\mypointsize);
      }
      \makeatother
      
      \draw[contour] (joint) -- (loadForceStart);
      \draw[contour] (tricepsArmStart) -- (tricepsForceStart);
      
      \draw[dashed,contour] (joint) -- (verticalDashedEnd);
      \draw[dashed,contour] (
        $(joint)!(horizontalDashedEnd)!(verticalDashedEnd)$
      ) -- (horizontalDashedEnd);
      \draw[dashed,contour] (joint) -- (
        $(tricepsForceStart)!(joint)!(tricepsForceEnd)$
      );
      \draw[dashed,contour] (joint) -- (
        $(bicepsForceStart)!(joint)!(bicepsForceEnd)$
      );
      \centerarc[dashed,contour](joint)(270:\myelbowangle:\myradius);
      
      \drawleverarm[contour]{horizontalDashedEnd}{joint}{verticalDashedEnd}{-}
      \drawleverarm[contour]{joint}{tricepsForceStart}{tricepsForceEnd}{-}
      \drawleverarm[contour]{joint}{bicepsForceStart}{bicepsForceEnd}{+}
      
      \draw[
        ->,contour,>={Stealth[width=9pt,length=9pt]},
      ] (tricepsForceStart) -- (
        $(tricepsForceEnd)!-1.5pt!(tricepsForceStart)$
      );
      \draw[
        ->,contour,>={Stealth[width=9pt,length=9pt]},
      ] (bicepsForceStart) -- (
        $(bicepsForceEnd)!-1.5pt!(bicepsForceStart)$
      );
      \draw[
        ->,contour,>={Stealth[width=9pt,length=9pt]},
      ] (loadForceStart) -- (
        $(loadForceEnd)!-1.5pt!(loadForceStart)$
      );
      
      \draw (joint) -- (loadForceStart);
      \draw (tricepsArmStart) -- (tricepsForceStart);
      
      \draw[dashed] (joint) -- (verticalDashedEnd);
      \centerarc[dashed](joint)(270:\myelbowangle:\myradius);
      
      \drawleverarm{horizontalDashedEnd}{joint}{verticalDashedEnd}{-}
      \drawleverarm{joint}{tricepsForceStart}{tricepsForceEnd}{-}
      \drawleverarm{joint}{bicepsForceStart}{bicepsForceEnd}{+}
      
      \draw[->,draw=C0] (tricepsForceStart) -- (tricepsForceEnd);
      \draw[->,draw=C0] (bicepsForceStart) -- (bicepsForceEnd);
      \draw[->,draw=C0] (loadForceStart) -- (loadForceEnd);
      
      \fill[draw=none,fill=black] (joint) circle (2pt);
      
      \node at (
        $(joint) + (
          {0.62*\myradius*cos((270+\myelbowangle)/2)},
          {0.62*\myradius*sin((270+\myelbowangle)/2)}
        )$
      ) {\contour{mittelblau!10}{$\elbang$}};
      
      \node[anchor=south] at (
        $0.5*(joint) +
        0.5*(tricepsForceStart)!(joint)!(tricepsForceEnd) +
        (-1mm,1mm)$
      ) {\contour{mittelblau!10}{\textcolor{C1}{$\armT$}}};
      \node[anchor=south] at (
        $0.5*(joint) +
        0.5*(bicepsForceStart)!(joint)!(bicepsForceEnd) +
        (1mm,1mm)$
      ) {\contour{mittelblau!10}{\textcolor{C1}{$\armB$}}};
      \node[anchor=north] at (
        $0.5*(horizontalDashedEnd) +
        0.5*(joint)!(horizontalDashedEnd)!(verticalDashedEnd) +
        (0mm,-1mm)$
      ) {\contour{mittelblau!10}{\textcolor{C1}{$\armL$}}};
      
      \node[anchor=east] at (
        $0.2*(tricepsForceStart) + 0.8*(tricepsForceEnd) + (-1mm,0mm)$
      ) {\contour{mittelblau!10}{\textcolor{C0}{$\forceT$}}};
      \node[anchor=west] at (
        $0.2*(bicepsForceStart) + 0.8*(bicepsForceEnd) + (1mm,0mm)$
      ) {\contour{mittelblau!10}{\textcolor{C0}{$\forceB$}}};
      \node[anchor=west] at (
        $0.2*(loadForceStart) + 0.8*(loadForceEnd) + (1mm,0mm)$
      ) {\contour{mittelblau!10}{\textcolor{C0}{$\forceL$}}};
    \end{scope}
  \end{tikzpicture}%
  \caption[Human upper limb model geometry as a raising arm movement]{%
    Human upper limb model geometry shown as a raising arm movement
    for the elbow angles
    $\theta = \ang{10}$, \ang{45}, \ang{80},
    \ang{115}, and \ang{150} \emph{(from left to right).}
    For $\theta = \ang{80}$,
    the contributing forces \emph{\textcolor{C0}{(blue)}} and
    lever arms \emph{\textcolor{C1}{(red)}} are shown.
    Taken and adapted from
    \multicite{Sprenger15Continuum,Valentin18Gradient}.%
  }%
  \label{fig:raisingArm}%
\end{figure}

\paragraph{Overall stress components}

The continuum-mechanical part of the model
is based on the theory of finite elasticity.
When muscles deform, forces act on each infinitesimally small element
of the muscles, which is known as \term{stress.}
Usually, especially in linear elasticity,
stress is measured with the \term{Cauchy stress tensor}
(also called the \term{true stress}) \cite{Soennerlind13Why}.
For non-linear stress-strain relations,
one may use other measures such as the
\term{second Piola--Kirchhoff stress}.
The second Piola--Kirchhoff stress has the additional advantage
that it is defined along the material directions,
in contrast to the Cauchy stress tensor,
which measures the stress in coordinate directions \cite{Soennerlind13Why}.

In \multicite{Sprenger15Continuum,Roehrle16Two,Valentin18Gradient},
the strain energy function is defined such that the
resulting overall second Piola--Kirchhoff stress $\mat{S}_\mathrm{MTC}$ of
the muscle-tendon complex satisfies
\begin{equation}
  \mat{S}_\mathrm{MTC}
  = \mat{S}_\mathrm{iso} + \mat{S}_\mathrm{aniso} - p \mat{C}^{-1},
\end{equation}
where $\mat{S}_\mathrm{iso}$ and $\mat{S}_\mathrm{aniso}$
are the \term{isotropic and anisotropic parts} of the stress, respectively,
$p$ is the \term{hydrostatic pressure} to ensure incompressibility,
and $\mat{C}$ is the \term{right Cauchy--Green deformation tensor.}
The anisotropic part $\mat{S}_\mathrm{aniso}$ is defined as
\begin{equation}
  \mat{S}_\mathrm{aniso}
  \ceq (
    \mat{S}_\mathrm{passive} +
    \act \gamma_\mathrm{M} \mat{S}_\mathrm{active}
  ) (1 - \gamma_\mathrm{ST}),
\end{equation}
cf.\ \cite{Valentin18Gradient}.
Here, $\mat{S}_\mathrm{passive}$ and $\mat{S}_\mathrm{active}$
are the \term{passive and active contributions}
due to the skeletal muscle fibers,
$\act \in \clint{0, 1}$ is the \term{activation parameter} of the
respective muscle-tendon complex, and
$\gamma_\mathrm{M}, \gamma_\mathrm{ST}$ are two \term{material parameters}
with which we can differentiate between the different types of soft tissues
of the muscle-tendon complex: fat, tendon, and muscle
\cite{Valentin18Gradient}.
Isotropic fat tissue can be obtained
by setting $\gamma_\mathrm{ST} \ceq 1$,
passive anisotropic tendon tissue
by setting $\gamma_\mathrm{ST} \ceq 0$ and $\gamma_\mathrm{M} \ceq 0$, and
skeletal muscle tissue
by setting $\gamma_\mathrm{ST} \ceq 0$ and $\gamma_\mathrm{M} \ceq 1$.
A mixture of these pure materials is
achieved by linear interpolation when setting
$\gamma_\mathrm{M}$ and $\gamma_\mathrm{ST}$ to values between zero and one
\cite{Valentin18Gradient}.
More details about the theory part of the model are given in
\multicite{Sprenger15Continuum,Roehrle16Two,Valentin18Gradient}.

\section{Momentum Equilibrium and Elbow Angle Optimization}
\label{sec:72methodology}

\minitoc[-7mm]{70mm}{4}

\noindent
In this section, we give an overview of the methodology of our approach.
We continue following the presentation of \cite{Valentin18Gradient}.

\subsection{From Muscle Forces to Equilibrium Angles}
\label{sec:721equilibrium}

\paragraph{Model inputs and outputs}

In the following, we regard simulations of the
human upper limb model described in \cref{sec:71model} as a black box,
which receives as its input
the elbow angle $\elbang \in \clint{\ang{10}, \ang{150}}$
and the activation parameters
$(\actT, \actB) \in \clint{\*0, \*1} = \clint{0, 1}^2$
of triceps and biceps.%
\footnote{%
  Here and in the following, the subscripts T, B, and L stand for
  triceps, biceps, and load, respectively.%
}
The outputs of the black box simulation are the forces
$\forceT(\elbang, \actT)$ and $\forceB(\elbang, \actB)$
that triceps and biceps exert.
These forces depend on the elbow angle as well as on the respective
activation parameter.
Gravitational forces due to the masses of bones or muscles
are neglected in this context.
However, we allow the specification of an external load $\forceL$,
which is applied to the end of the forearm.
This load may be the weight force of some object
that the arm is supposed to keep in position.

\paragraph{Moments and lever arms}

Each force exerts a \term{moment} (or \term{torque}) on the elbow joint.
The moments are the products of the forces $\forceX$
with the respective lever arms $\armX$
($X \in \{\mathrm{T}, \mathrm{B}, \mathrm{L}\}$).
The lever arms are approximated as in
\multicite{Roehrle16Two,Valentin18Gradient} by using
the tendon-displacement method of \cite{An84Determination}:
\begin{subequations}
  \begin{align}
    \armT(\elbang)
    &\ceq (-0.0009399 \{\elbang\}^2 + 0.1126 \{\elbang\} + 22.21)\;
    \si{\milli\meter},\\
    \armB(\elbang)
    &\ceq (-0.001482 \{\elbang\}^2 + 0.1776 \{\elbang\} + 35.02)\;
    \si{\milli\meter},\\
    \armL(\elbang)
    &\ceq \sin(\elbang) \cdot \SI{282.5}{\milli\meter},
  \end{align}
\end{subequations}
where $\{\elbang\}$ denotes the dimensionless value of $\elbang$
in degrees.
The lever arms are non-negative and the forces are signed, i.e.,
positive forces pull the forearm downwards and
negative forces pull it upwards.
In general, $\forceT, \forceL \ge \SI{0}{\newton}$ and
$\forceB \le \SI{0}{\newton}$.

\paragraph{Total moment and equilibrium elbow angle}

The \term{total moment} of the system is given by the function
\begin{subequations}
  \label{eq:totalMoment}
  \begin{gather}
    \moment_{\forceL,\actT,\actB}\colon
    \clint{\ang{10}, \ang{150}} \to \real,\\
    \moment_{\forceL,\actT,\actB}(\elbang)
    \ceq \forceT(\elbang, \actT) \armT(\elbang) +
    \forceB(\elbang, \actB) \armB(\elbang) +
    \forceL \armL(\elbang),
  \end{gather}
\end{subequations}
cf.\ \cite{Valentin18Gradient}.
The system is in \term{equilibrium}
if the total moment vanishes, i.e.,
$\moment_{\forceL,\actT,\actB}(\elbang) = \SI{0}{\newton\meter}$.
We call the corresponding angle $\elbang$ the
\term{equilibrium elbow angle}
for the load $\forceL$ and the activation parameters $\actT, \actB$.
To find this angle for a given load $\forceL$ and activation parameters
$\actT$ and $\actB$, we first note that
$\moment_{\forceL,\actT,\actB}$ may have zero, exactly one,
or multiple zeros in $\clint{\ang{10}, \ang{150}}$.
Hence, the inverse function evaluated at $\SI{0}{\newton\meter}$
is partially defined depending on the load and the activation parameters:
\begin{equation}
  \label{eq:equilibriumAngle}
  \equielbang{\forceL}\colon \actdomain{\forceL} \!\to
  \clint{\ang{10}, \ang{150}},\quad
  \actdomain{\forceL} \subset \clint{\*0, \*1},\quad
  \equielbang{\forceL}(\actT,\actB)
  \ceq (\moment_{\forceL,\actT,\actB})^{-1}(\SI{0}{\newton\meter}),
\end{equation}
which is well-defined whenever $\moment_{\forceL,\actT,\actB}$
has a unique root.
We approximate $\equielbang{\forceL}(\actT,\actB)$ with the Newton method
\multicite{Roehrle16Two,Valentin18Gradient}:
\begin{equation}
  \label{eq:newtonAngle}
  \elbang^{(j+1)}
  \ceq \elbang^{(j)} -
  \frac{
    \moment_{\forceL,\actT,\actB}(\elbang^{(j)})
  }{
    \partialderiv{\partialdiff{} \elbang}{\moment_{\forceL,\actT,\actB}}
    (\elbang^{(j)})
  },\quad
  j \in \nat,
\end{equation}
with an initial value
$\elbang^{(0)} \in \clint{\ang{10}, \ang{150}}$
and the stopping criterion of
$\abs{\moment_{\forceL,\actT,\actB}(\elbang^{(j)})} <
\SI{e-9}{\newton\meter}$.
We repeat the Newton method for the initial values
$\elbang^{(0)} = \ang{80}, \ang{40}, \ang{120}$
and use the first converged result
(i.e., we check if $\elbang^{(0)} = \ang{80}$ converges;
if not, we proceed with $\elbang^{(0)} = \ang{40}$, and so on).
If all three initial values do not converge,
we conclude that $(\actT, \actB) \notin \actdomain{\forceL}$.

\subsection{Optimization Problems}
\label{sec:722optimization}

\paragraph{General problem}

The general problem in our setting is as follows:
For a given external load $\forceL$ and a target elbow angle $\tarelbang$,
find activation parameters $(\actT, \actB) \in \clint{\*0, \*1}$
such that the target elbow angle is attained in the equilibrium,
i.e., $\equielbang{\forceL}(\actT,\actB) = \tarelbang$.
Example applications of such a scenario are medicine and robotics,
when a specific movement should be carried out.

\paragraph{List of optimization problems}

As discussed in \cref{sec:711models},
musculoskeletal systems with an antagonistic muscle pair
such as our human upper limb model are usually overdetermined.
This means that there are multiple solutions to this general problem.
As a remedy, one may solve one of the following two
optimization problems \cite{Valentin18Gradient}:

\begin{enumerate}[label=O\arabic*.,ref=O\arabic*,leftmargin=2.7em]
  \item
  \label{item:biomech2MinSum}
  For a given external load $\forceL$ and a target angle
  $\tarelbang \in \clint{\ang{10}, \ang{150}}$,
  find the activation parameters $(\actT, \actB) \in \clint{\*0, \*1}$
  such that $\actT + \actB$ is minimized under the constraint
  $\equielbang{\forceL}(\actT, \actB) = \tarelbang$.
  
  \item
  \label{item:biomech2MinDist}
  For a given external load $\forceL(t_2)$ for a time $t_2 > t_1$,
  a target angle $\tarelbang(t_2) \in \clint{\ang{10}, \ang{150}}$,
  and initial activation parameters
  $(\actT(t_1), \actB(t_1)) \in \clint{\*0, \*1}$,
  find new activation parameters
  $(\actT(t_2), \actB(t_2)) \in \clint{\*0, \*1}$ such that
  $(\actT(t_2) - \actT(t_1))^2 + (\actB(t_2) - \actB(t_1))^2$
  is minimized under the constraint
  $\equielbang{\forceL(t_2)}(\actT(t_2), \actB(t_2)) = \tarelbang(t_2)$.
\end{enumerate}

\noindent
The motivation of both problems is that the human body tries to
achieve a given movement with minimal energy effort.

\paragraph{Motivation of problem \ref{item:biomech2MinSum}}

For the first problem \ref{item:biomech2MinSum},
this effort is quantified by $\actT + \actB$,
i.e., the energy effort for each muscle is assumed to be proportional
to its activation parameter.

\paragraph{Motivation of problem \ref{item:biomech2MinDist}}

The second problem \ref{item:biomech2MinDist} is motivated as follows:
Before time $t = t_1$, the musculoskeletal system is in equilibrium for
the external load $\forceL(t_1)$,
activation parameters $\actT(t_1), \actB(t_1)$, and elbow angle
$\tarelbang(t_1) \ceq \equielbang{\forceL(t_1)}(\actT(t_1), \actB(t_1))$, i.e.,
$\moment_{\forceL(t_1),\actT(t_1),\actB(t_1)}(\tarelbang(t_1))
= \SI{0}{\newton\meter}$.
Directly after $t = t_1$,
the external force and/or the target angle is suddenly changed
to $\forceL(t_2)$ and $\tarelbang(t_2)$, respectively.
Consequently, triceps and biceps adapt their activation parameters
such that the musculoskeletal system returns to equilibrium
at some time $t = t_2 > t_1$.
Hence, we have to determine the new activation parameters
$\actT(t_2), \actB(t_2)$ such that
$\moment_{\forceL(t_2),\actT(t_2),\actB(t_2)}(\tarelbang(t_2))
= \SI{0}{\newton\meter}$.
Again, these parameters
$\actT(t_2)$ and $\actB(t_2)$ are not uniquely determined.
Therefore, we want to find the pair of activation parameters
that is closest (in terms of the Euclidean norm) to the initial
activation parameters $\actT(t_1), \actB(t_1)$.

\paragraph{Optimization method}

Problems \ref{item:biomech2MinSum} and \ref{item:biomech2MinDist}
are both constrained optimization problems.
For their solution, we employ the augmented Lagrangian method as
described in \cref{sec:513gradientBasedConstrained}
using an adaptive gradient descent algorithm
for the gradient-based optimization of the penalized objective function
(see \cref{sec:512gradientBasedUnconstrained}).

\subsection{B-Spline Surrogates on Sparse Grids}
\label{sec:723surrogates}

\paragraph{Complexity}

To solve optimization problems
\ref{item:biomech2MinSum} and \ref{item:biomech2MinDist},
the optimization method needs to evaluate the objective
and constraint functions multiple times during the algorithm.
This requires the evaluation of $\equielbang{\forceL}$,
which in turn has to be approximated with the Newton method.
As we see in \cref{eq:newtonAngle},
each iteration of the Newton method needs not only the values of the
muscle forces $\forceT$ and $\forceB$, but also their
partial derivatives with respect to $\elbang$.
These partial derivatives have to be approximated with finite differences.

Unfortunately, simulations of continuum-mechanical models are
computationally expensive.
One evaluation of the muscle force pair $\forceT, \forceB$
requires the solution of a solid mechanics model
with a complex constitutive law, pre-stretch, and contact between
bone and muscles \cite{Valentin18Gradient}.
On average, a single evaluation of $\forceT$ and $\forceB$ takes
about half an hour on current desktop computers.
If we assume that we need four Newton iterations on average,
then a single iteration of the optimization algorithm to solve
problems \ref{item:biomech2MinSum} and \ref{item:biomech2MinDist}
will take four hours to complete
(assuming one evaluation of objective and constraint functions
per optimizer iteration and
two evaluations of the muscle force pair
per Newton iteration to approximate the missing derivative).
Consequently, the whole optimization process takes
two weeks to complete, if the optimizer converges after 100 iterations.

\paragraph{Sparse grid surrogates}

A popular way to reduce complexity is to employ surrogates.
In this case, the idea is to replace the muscle force functions
$\forceT, \forceB$ with surrogates $\forceTintp, \forceBintp$
\cite{Valentin18Gradient}, e.g., by interpolation.
We then automatically obtain a surrogate
\begin{subequations}
  \label{eq:totalMomentSurrogate}
  \begin{gather}
    \momentintp_{\forceL,\actT,\actB}\colon
    \clint{\ang{10}, \ang{150}} \to \real,\\
    \momentintp_{\forceL,\actT,\actB}(\elbang)
    \ceq \forceTintp(\elbang, \actT) \armT(\elbang) +
    \forceBintp(\elbang, \actB) \armB(\elbang) +
    \forceL \armL(\elbang),
  \end{gather}
\end{subequations}
for the total moment (cf.\ \cref{eq:totalMoment}) and,
consequently, a surrogate
\begin{equation}
  \label{eq:equilibriumAngleSurrogate}
  \equielbangintp{\forceL}\colon \actdomainintp{\forceL} \!\to
  \clint{\ang{10}, \ang{150}},\quad
  \actdomainintp{\forceL} \subset \clint{\*0, \*1},\quad
  \equielbangintp{\forceL}(\actT,\actB)
  \ceq (\momentintp_{\forceL,\actT,\actB})^{-1}(\SI{0}{\newton\meter}),
\end{equation}
for the equilibrium elbow angle function (cf.\ \cref{eq:equilibriumAngle}).
Since the surrogates are much cheaper to evaluate,
the computation time is decreased by up to seven orders of magnitude,
as experiments show.

The approach in \cite{Valentin18Gradient} and in this thesis is
to determine surrogates
$\forceXintp\colon \clint{\*0, \*1} \to \real$
($X \in \{\mathrm{T}, \mathrm{B}\}$) by sparse grid interpolation.
Compared to surrogate construction techniques based on full grids,
sparse grids help to reduce the number of samples that
are necessary to build ``reasonably'' accurate surrogates,
especially if the number of dimensions is moderately large
($d \ge 4$, \term{curse of dimensionality}).

The present model only has $d = 2$ dimensions ($\actT$ and $\actB$),
since the model contains only two muscles.
However, as we will see,
already for this low-dimensional problem,
sparse grids outperform conventional full grid interpolation.
The results have to be seen as a proof of concept.
One will be able to handle higher dimensionalities
(i.e., models with a larger number of muscles) similarly with little
or even no adjustments at all.
The low dimensionality of the model in this thesis
enables us to compute and compare against reference solutions,
which would not be possible in a higher-dimensional setting.

\paragraph{Benefiting from B-splines}

As in \cite{Valentin18Gradient},
we use higher-order hierarchical B-splines as basis functions
for the sparse grid surrogates.
This has three advantages when compared with conventional
sparse grid bases such as piecewise linear functions:
First, the partial derivative
$\partialderiv{\partialdiff{} \elbang}{\momentintp}$ needed
for the Newton method in \cref{eq:newtonAngle} is continuous and
explicitly known.
There is no need to approximate the derivative with
finite differences, reducing both error and computation time.
Second, we can use gradient-based optimization methods
for the solution of the optimization problems \ref{item:biomech2MinSum} and
\ref{item:biomech2MinDist},
which involve the equilibrium elbow angle function
$\equielbangintp{\forceL}\colon \clint{\*0, \*1} \to \real$.
With the implicit function theorem \cite{Kudryavtsev95Implicit},
we obtain for the derivative of $\equielbangintp{\forceL}$
\begin{equation}
  \gradient{\actT,\actB}{\equielbangintp{\forceL}}
  = -\,(\gradient{\actT,\actB}{\momentintp}) \cdot
  (\gradient{\elbang}{\momentintp})^{-1}
  = -\,\frac{
    \gradient{\actT,\actB}{\momentintp}
  }{
    \partialderiv{\partialdiff{} \elbang}{\momentintp}
  },
\end{equation}
where $\gradient{\actT,\actB}{}$ is the transposed Jacobian
with respect to $\actT$ and $\actB$.%
\footnote{%
  For example, the first column is the gradient with respect to $\actT$
  and the second column is the gradient with respect to $\actB$.%
}
For B-splines,
both the transposed Jacobian $\gradient{\actT,\actB}{\momentintp}$ and
the partial derivative $\partialderiv{\partialdiff{} \elbang}{\momentintp}$
are continuous, explicitly known, and can be evaluated fast.
Third and finally,
the usage of higher-order B-splines as basis functions
increases the order of convergence of interpolation errors
as shown for test functions in \cref{sec:541interpolation}.
Thus, fewer interpolation points are necessary to construct a surrogate
with the same error as for piecewise linear functions.

\breakpagebeforenextheadingtrue
\section{Implementation and Numerical Results}
\label{sec:73results}

\minitoc[0mm]{69mm}{8}

\parbox{1em}{}
\vspace{-3em}

\disableornamentsfornextheadingtrue
\subsection{Implementation}
\label{sec:731implementation}

\paragraph{Parameters, implementation, and geometry}

Details about implementational aspects of the model can be found in
\multicite{Sprenger15Continuum,Roehrle16Two,Valentin18Gradient},
for instance, values for the material parameters.
The constitutive law has been implemented in the CMISS software package
(an interactive computer program for Continuum Mechanics,
Image analysis, Signal processing and System identification%
\footnote{%
  \url{https://www.cmiss.org/}%
}).
The emerging PDEs are discretized using quadratic finite element basis
functions and the resulting linearized system is solved with CMISS.
The geometry of the human upper limb model is based on
the Visible Human Male's dataset \cite{Spitzer96Visible}.
Again, we refer to \multicite{Sprenger15Continuum,Roehrle16Two} for details
about the geometry.

\subsection{Reference and Sparse Grid Solution}
\label{sec:732solutionTypes}

\paragraph{Reference solution}

\begin{figure}
  \includegraphics{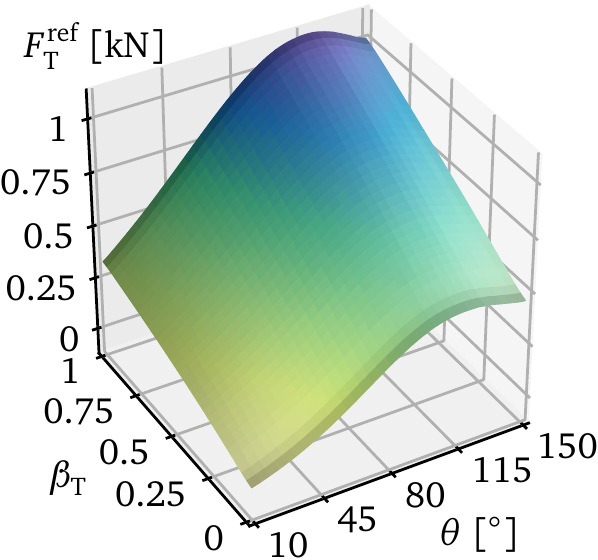}%
  \;\;%
  \includegraphics{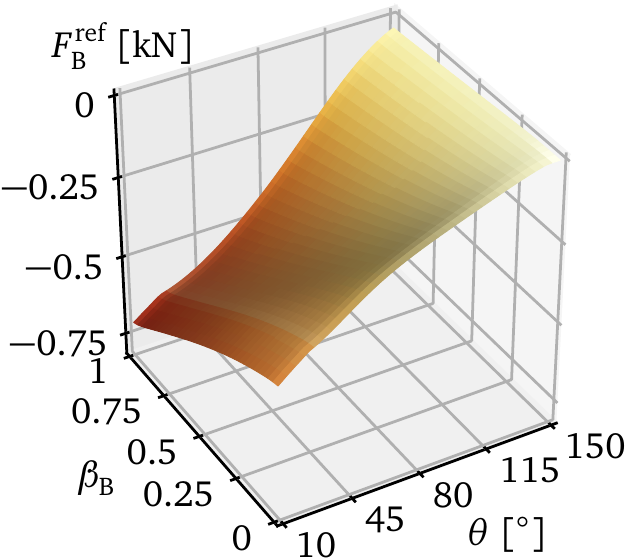}%
  \hfill%
  \rlap{\raisebox{53mm}{\;$\forceXref$ [\si{\kilo\newton}]}}%
  \includegraphics{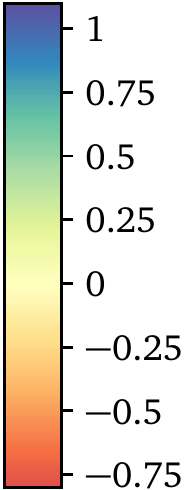}%
  \caption[Reference triceps and biceps forces]{%
    Reference triceps and biceps forces $\forceXref$
    ($X \in \{\mathrm{T}, \mathrm{B}\}$).%
  }%
  \label{fig:biomech2ReferenceForce}%
\end{figure}

\begin{figure}
  \includegraphics{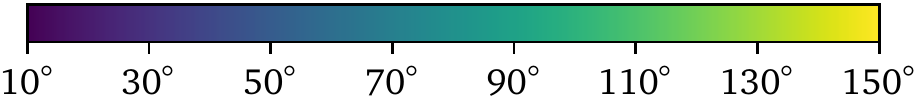}%
  \\[2mm]%
  \subcaptionbox{%
    $\forceL = \SI{22}{\newton}$%
  }[49mm]{%
    \includegraphics{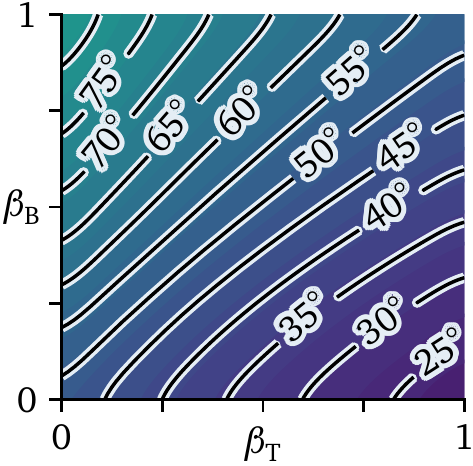}%
  }%
  \hfill%
  \subcaptionbox{%
    $\forceL = \SI{-60}{\newton}$%
  }[49mm]{%
    \includegraphics{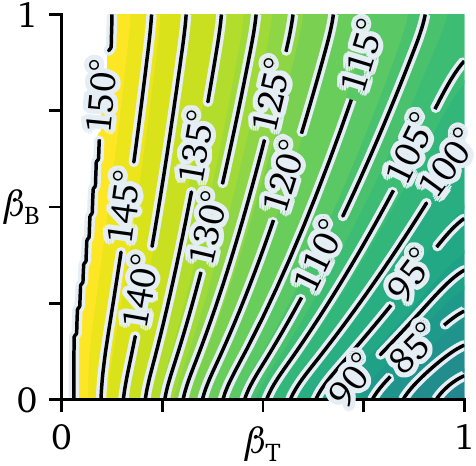}%
  }%
  \hfill%
  \subcaptionbox{%
    $\forceL = \SI{180}{\newton}$%
  }[49mm]{%
    \includegraphics{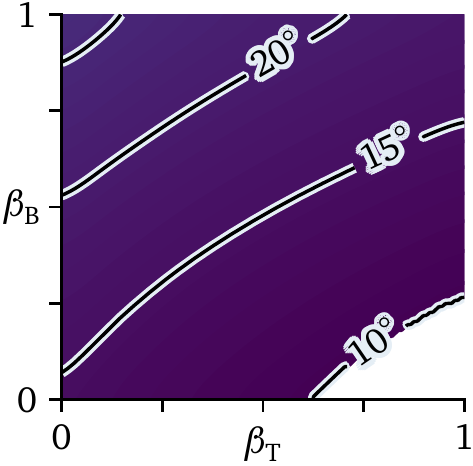}%
  }%
  \caption[Reference equilibrium elbow angle]{%
    Reference equilibrium elbow angle $\equielbangref{\forceL}$
    for different loads $\forceL$.
    The empty areas correspond to activation pairs
    at which $\equielbangref{\forceL}$ is not well-defined
    (see \cref{eq:equilibriumAngle}).%
  }%
  \label{fig:biomech2ReferenceEquilibriumAngle}%
\end{figure}

Since the model is only two-dimensional, we can compute a reference solution
on a full grid.
To this end, we evaluate the exerted muscle forces $\forceT$ and $\forceB$ on
the full grid
\begin{equation}
  \{\ang{10}, \ang{11}, \dotsc, \ang{150}\} \times \{0, 0.1, \dotsc, 1\}
  \ni (\elbang, \actX),\quad
  X \in \{\mathrm{T}, \mathrm{B}\}.
\end{equation}
The resulting \num{1551} grid points
are interpolated with bicubic full grid splines%
\footnote{%
  Computed with the Geometric Tools Engine \cite{Schneider03Geometric},
  see \url{https://www.geometrictools.com/}.%
}
to obtain \term{reference solutions}
$\forceTref, \forceBref\colon
\clint{\ang{10}, \ang{150}} \times \clint{0, 1} \to \real$,
which are shown in \cref{fig:biomech2ReferenceForce}.
Due to the high resolution of the full grid,
we may assume that the reference solutions are accurate enough
to ensure $\forceTref \approx \forceT$ and $\forceBref \approx \forceB$.
We refer to the resulting equilibrium elbow angle
with $\equielbangref{\forceL}$.
It is displayed in
\cref{fig:biomech2ReferenceEquilibriumAngle}
for the loads of $\forceL = \SI{22}{\newton}$,
$\SI{-60}{\newton}$, and $\SI{180}{\newton}$.

\paragraph{Sparse grid solution}

Additionally, we evaluate $\forceT$ and $\forceB$ at the $\ngp = 49$
grid points
\begin{equation}
  \{(\elbang^{(k,\mathrm{unif})}, \actX^{(k,\mathrm{unif})}) \mid
  k = 1, \dotsc, \ngp\}
  \subset \clint{\ang{10}, \ang{150}} \times \clint{0, 1},\quad
  X \in \{\mathrm{T}, \mathrm{B}\},
\end{equation}
of the uniform regular sparse grid $\interiorregsgset{n}{d}$ of
level $n = 5$ in $d = 2$ dimensions
without boundary points (to reduce the number of samples)
and at the sparse Clenshaw--Curtis grid
\begin{equation}
  \{(\elbang^{(k,\cc)}, \actX^{(k,\cc)}) \mid
  k = 1, \dotsc, \ngp\}
  \subset \clint{\ang{10}, \ang{150}} \times \clint{0, 1},\quad
  X \in \{\mathrm{T}, \mathrm{B}\},
\end{equation}
of the same size and level.%
\footnote{%
  The domain $\clint{\ang{10}, \ang{150}} \times \clint{0, 1}$
  is assumed to be implicitly normalized to the unit square
  $\clint{\*0, \*1}$.%
}
These values are interpolated using three
different hierarchical B-spline bases of degree $p = 1$, $3$, and $5$:
modified hierarchical uniform B-splines
$\bspl[\modified]{\*l,\*i}{p}$
(see \cref{sec:313modification}),
modified hierarchical Clenshaw--Curtis B-splines
$\bspl[\cc,\modified]{\*l,\*i}{p}$
(see \cref{sec:314nonUniform}), and
modified hierarchical uniform not-a-knot B-splines
$\bspl[\nak,\modified]{\*l,\*i}{p}$
(see \cref{sec:323modifiedNAKBSplines}).
The implementation was done using the sparse grid toolbox
\sgpp{} \cite{Pflueger10Spatially}.%
\footnote{%
  \url{http://sgpp.sparsegrids.org/}%
}
The corresponding interpolants and resulting quantities
are denoted with the superscripts
``$\sparse,\!p$'', ``$\sparse,\!p,\!\cc$'', or ``$\sparse,\!p,\!\nak$'',
respectively.
A superscript of ``$\sparse$'' without any further specification
means one of the three hierarchical B-spline bases in general.
Note that the equilibrium elbow angle is \emph{not} interpolated
(neither in the full grid nor in the sparse grid case),
but rather obtained by inserting the interpolated muscle forces
into \eqref{eq:totalMomentSurrogate} and \eqref{eq:equilibriumAngleSurrogate}.

\subsection{Errors of Muscle Forces and Equilibrium Angle}
\label{sec:733errors}

\paragraph{Quality of reference interpolants}

Before we turn to the sparse grid interpolants,
we assess the quality of the reference interpolants on the full grid.
For this purpose, we evaluate the full grid interpolants
$\forceTintp, \forceBintp$
at the sparse grid points $(\elbang^{(k)}, \actX^{(k)})$
(which are not a subset of the full grid points!)
and compare the resulting values with the known exact values
$\forceT(\elbang^{(k)}, \actT^{(k)})$ and
$\forceB(\elbang^{(k)}, \actB^{(k)})$
of the muscle forces $\forceT, \forceB$.
We also incorporate the known values at the sparse
Clenshaw--Curtis grid points.
In particular, let $G$ be the union of
$\{(\elbang^{(k,\mathrm{unif})}, \actX^{(k,\mathrm{unif})}) \mid
k = 1, \dotsc, \ngp\}$ and
$\{(\elbang^{(k,\cc)}, \actX^{(k,\cc)}) \mid k = 1, \dotsc, \ngp\}$.
We then approximate the relative $\Ltwo$ interpolation error
of the reference interpolants by
\begin{equation}
  \frac{\normLtwo{\forceX - \forceXref}}{\normLtwo{\forceX}}
  \approx
  \frac{
    \setsize{G}^{-1/2}
    \norm[2]{
      (\forceX(\elbang, \actX) - \forceXref(\elbang, \actX))_
      {(\elbang, \actX) \in G}
    }
  }{
    \setsize{G}^{-1/2}
    \norm[2]{(\forceX(\elbang, \actX))_{(\elbang, \actX) \in G}}
  },\quad
  X \in \{\mathrm{T}, \mathrm{B}\},
\end{equation}
where $\norm[2]{\cdot}$ is the Euclidean norm.%
\footnote{%
  We have $\setsize{G} = 2\ngp - 1$, since sparse grids of
  uniform and Clenshaw--Curtis type only
  share the center point $(\elbang, \actX) = (\ang{80}, 0.5)$,
  if there are no boundary points.%
}
After inserting the known values $\forceX(\elbang, \actX)$ and
$\forceXref(\elbang, \actX)$ ($(\elbang, \actX) \in G$)
on the right-hand side,
we obtain
\begin{equation}
  \frac{\normLtwo{\forceT - \forceTref}}{\normLtwo{\forceT}}
  \approx \SI{2.19}{\permille},\qquad
  \frac{\normLtwo{\forceB - \forceBref}}{\normLtwo{\forceB}}
  \approx \SI{2.06}{\permille}.
\end{equation}
These errors are very small, which justifies our assumption of
$\forceTref \approx \forceT$ and $\forceBref \approx \forceB$.

\paragraph{Error of sparse grid muscle forces}

\Cref{tbl:biomech2ErrorL2_1} contains the relative $\Ltwo$
interpolation errors
$\normLtwo{\forceXref - \forceXintp}/\normLtwo{\forceXref}$
($X \in \{\mathrm{T}, \mathrm{B}\}$) of the sparse grid interpolants
for all hierarchical bases and degrees $p = 1, 3, 5$.
All reported errors are relatively small
due to the smoothness of the original functions
(cf.\ $\forceXref$ in \cref{fig:biomech2ReferenceForce}).
All in all, the modified Clenshaw--Curtis B-splines perform best,
achieving relative $\Ltwo$ errors of below \SI{3.6}{\permille}
in the cubic case.
Surprisingly, the not-a-knot B-splines are the worst choice in our
comparison.
Their corresponding errors exceed \SI{1}{\percent} for the triceps
and $p > 1$.
The possible reasons are two-fold:
First, there might be slight noise in the given muscle force data,
which is visible in \cref{fig:biomech2ReferenceForce},
as there seems to be a kink in $\forceBref$ at $\elbang \approx \ang{25}$.
Second, the employed regular sparse grids might be too coarse
as the higher convergence order of not-a-knot B-splines
only pays off in the asymptotic range (see \cref{sec:541interpolation}).
The same observations hold for the degree $p$,
for which $p = 3$ seems to be the best choice,
as the errors increase again for $p = 5$.

\begin{table}
  \newcommand*{\bi}{$\bspl[\modified]{l,i}{p}$}
  \newcommand*{\bii}{$\bspl[\cc,\modified]{l,i}{p}$}
  \newcommand*{\biii}{$\bspl[\nak,\modified]{l,i}{p}$}
  \subcaptionbox{%
    $\normLtwo{\forceXref - \forceXintp}/\normLtwo{\forceXref}$
    [\si{\permille}] given as triceps/biceps pairs
    ($X \in \{\mathrm{T}, \mathrm{B}\}$).%
    \label{tbl:biomech2ErrorL2_1}%
  }[85.2mm]{%
    \setnumberoftableheaderrows{1}%
    \begin{tabular}{%
      >{\kern\tabcolsep}=l<{\kern2mm}%
      +c<{\kern-1mm}+c<{\kern-1mm}+c<{\kern\tabcolsep}%
    }
      \toprulec
      \headerrow
      $p$&   $1$&                  $3$&                  $5$\\
      \midrulec
      \bi&   $3.60,7.12$&          $3.05,7.00$&          $\mathbf{2.98},7.90$\\
      \bii&  $\mathbf{3.28},4.35$& $3.31,\mathbf{3.56}$& $3.35,3.64$\\
      \biii& $3.60,7.12$&          $3.09,10.0$&          $7.13,24.6$\\
      \bottomrulec
    \end{tabular}%
  }%
  \hfill%
  \subcaptionbox{%
    $\normLtwo{\equielbangref{\forceL} - \equielbangintp{\forceL}}/
    \normLtwo{\equielbangref{\forceL}}$
    [\si{\permille}] for $\forceL = \SI{22}{\newton}$.%
    \label{tbl:biomech2ErrorL2_2}%
  }[59mm]{%
    \setnumberoftableheaderrows{1}%
    \begin{tabular}{%
      >{\kern\tabcolsep}=l<{\kern2mm}%
      +c<{\kern-1mm}+c<{\kern-1mm}+c<{\kern\tabcolsep}%
    }
      \toprulec
      \headerrow
      $p$&   $1$&    $3$&             $5$\\
      \midrulec
      \bi&   $4.15$& $3.74$&          $3.72$\\
      \bii&  $3.42$& $\mathbf{2.83}$& $2.86$\\
      \biii& $4.15$& $4.06$&          $8.28$\\
      \bottomrulec
    \end{tabular}%
  }%
  \caption[Relative $L^2$ errors of forces and equilibrium elbow angle]{%
    Relative $\Ltwo$ errors of triceps/biceps force \emph{(left)} and
    equilibrium elbow angle \emph{(right)}
    for different hierarchical bases $\basis{\*l,\*i}$ and
    B-spline degrees $p$.
    Highlighted entries are the best among those with
    the same hierarchical basis or the same degree
    (similar to Nash equilibria).%
  }%
  \label{tbl:biomech2ErrorL2}%
\end{table}

\vspace{\fill}

\begin{figure}
  \includegraphics{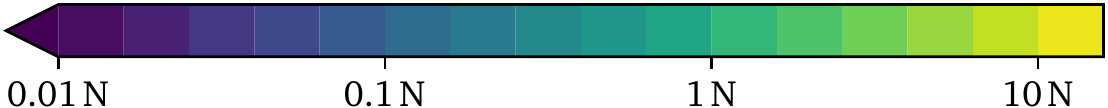}%
  \\[2mm]%
  \subcaptionbox{%
    $\abs{\forceXref - \forceXintp[p]}$ for
    $X = \mathrm{T}$ \emph{(left)} and
    $X = \mathrm{B}$ \emph{(right).}%
  }[74mm]{%
    \includegraphics{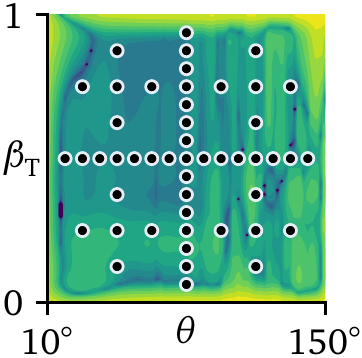}%
    \hfill%
    \includegraphics{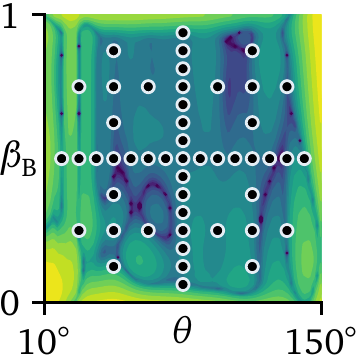}%
  }%
  \hfill%
  \subcaptionbox{%
    $\abs{\forceXref - \forceXintp[p,\cc]}$ for
    $X = \mathrm{T}$ \emph{(left)} and
    $X = \mathrm{B}$ \emph{(right).}%
  }[74mm]{%
    \includegraphics{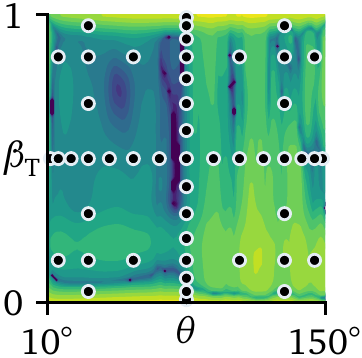}%
    \hfill%
    \includegraphics{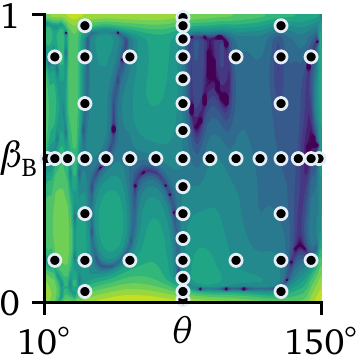}%
  }%
  \caption[Absolute error of muscle forces]{%
    Absolute error of muscle forces $\forceT, \forceB$ for
    modified cubic B-splines ($p = 3$)
    on sparse grids of uniform type \emph{(left two plots)} and
    of Clenshaw--Curtis type \emph{(right two plots)}
    together with the points of the sparse grid \emph{(dots).}%
  }%
  \label{fig:biomech2ErrorForce}%
\end{figure}

\Cref{fig:biomech2ErrorForce} shows the pointwise absolute error
$\abs{\forceXref(\elbang, \actX) - \forceXintp(\elbang, \actX)}$
for the modified B-splines $\bspl[\modified]{l,i}{p}$ and
$\bspl[\cc,\modified]{l,i}{p}$
on uniform and Clenshaw--Curtis grids in the cubic case $p = 3$.
Note that in contrast to usual interpolation settings,
the absolute errors $\abs{\forceXref - \forceXintp}$
shown in \cref{fig:biomech2ErrorForce} do not vanish at the
sparse grid points $(\elbang^{(k)}, \actX^{(k)})$
($X \in \{\mathrm{T}, \mathrm{B}\}$, $k = 1, \dotsc, \ngp$),
since $\forceXintp$ does not interpolate $\forceXref$
at these points.%
\footnote{%
  It would have been possible to construct $\forceXintp$
  as a sparse grid interpolant of $\forceXref$.
  However, building a spline surrogate ($\forceXintp$)
  of another spline surrogate ($\forceXref$) would skew the results.%
}
As it is typical for (modified) sparse grid interpolants,
the error is the largest near the boundary of the domain.
However, the Clenshaw--Curtis points help to decrease the error
due to the higher density of grid points near the boundary.
In the Clenshaw--Curtis case, the maximal errors are
\begin{equation}
  \normLinfty{\forceTref - \forceTintp[p,\cc]}
  \approx \SI{10.6}{\newton},\qquad
  \normLinfty{\forceBref - \forceBintp[p,\cc]}
  \approx \SI{9.51}{\newton},
\end{equation}
where $\normLinfty{\forceXref - \forceXintp[p,\cc]}
\ceq \max_{(\elbang, \actX)}
\abs{\forceXref(\elbang, \actX) - \forceXintp[p,\cc](\elbang, \actX)}$
(since the functions are continuous).
If we restrict the domain to
$\clint{\ang{31}, \ang{129}} \times \clint{0.15, 0.85}$
by omitting \SI{15}{\percent} on each side of the original domain,
then the maximal absolute errors drop to only
\SI{6.73}{\newton} (triceps) and \SI{0.967}{\newton} (biceps),
which is small compared to maximal possible forces of
around \SI{1}{\kilo\newton}.

\pagebreak

\paragraph{Error of the equilibrium elbow angle}

\begin{figure}
  \includegraphics{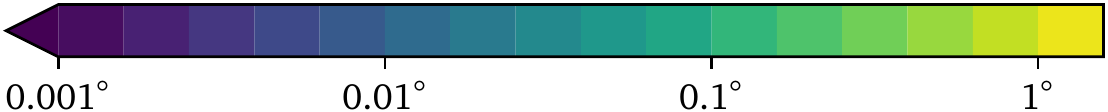}%
  \\[2mm]%
  \subcaptionbox{%
    $\forceL = \SI{22}{\newton}$%
  }[49mm]{%
    \includegraphics{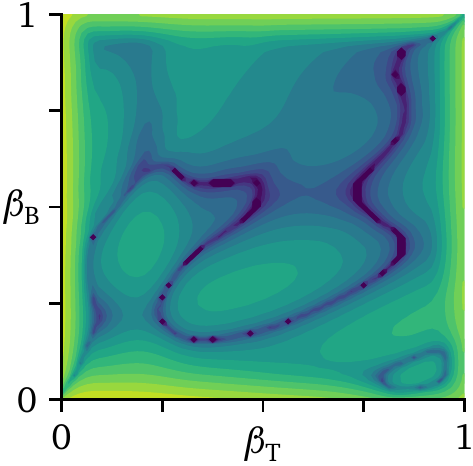}%
  }%
  \hfill%
  \subcaptionbox{%
    $\forceL = \SI{-60}{\newton}$%
  }[49mm]{%
    \includegraphics{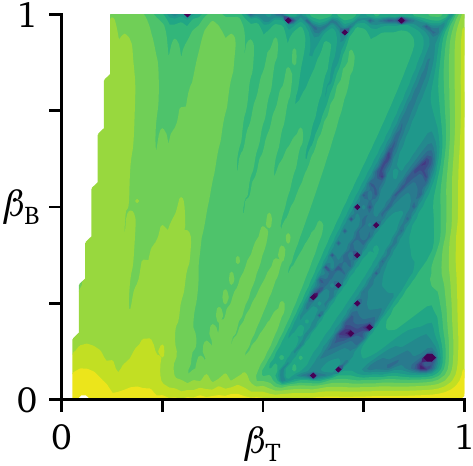}%
  }%
  \hfill%
  \subcaptionbox{%
    $\forceL = \SI{180}{\newton}$%
  }[49mm]{%
    \includegraphics{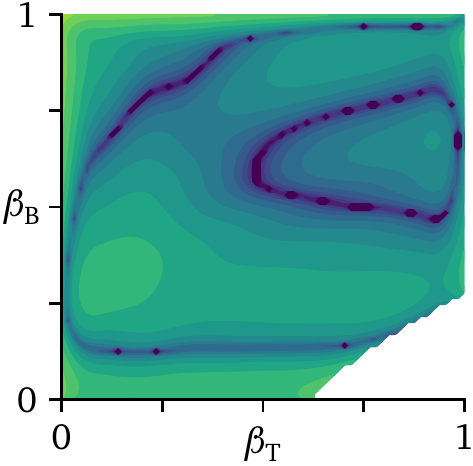}%
  }%
  \caption[Absolute error of the equilibrium elbow angle]{%
    Absolute error
    $\abs{\equielbangref{\forceL} - \equielbangintp[p,\cc]{\forceL}}$
    of the equilibrium elbow angle for
    modified hierarchical cubic Clenshaw--Curtis B-splines ($p = 3$)
    for different loads $\forceL$.
    In the empty areas, at least one of
    $\equielbangref{\forceL}$ and $\equielbangintp[p,\cc]{\forceL}$
    is not well-defined (see \cref{eq:equilibriumAngle}).%
  }%
  \label{fig:biomech2ErrorEquilibriumAngle}%
\end{figure}

The relative $\Ltwo$ errors
$\normLtwo{\equielbangref{\forceL} - \equielbangintp{\forceL}}/
\normLtwo{\equielbangref{\forceL}}$
of the equilibrium elbow angle function are shown in
\cref{tbl:biomech2ErrorL2_1} for the load of
$\forceL = \SI{22}{\newton}$.
Modified cubic Clenshaw--Curtis B-splines achieve the best results.
Therefore, we use this type of hierarchical basis
for the remainder of this chapter.
Pointwise plots of the absolute error
$\abs{\equielbangref{\forceL} - \equielbangintp[p,\cc]{\forceL}}$
are presented in \cref{fig:biomech2ErrorEquilibriumAngle}.
Again, the maximal error is comparatively small:
For $\forceL = \SI{22}{\newton}$, it is only \ang{0.886}.
If we restrict the domain to $\clint{0.15, 0.85}^2$,
then this maximal error drops to \ang{0.103} (or \ang{;6.18;}),
as the areas near the boundary of $\clint{\*0, \*1}$
contribute the most to the error.

\subsection{Test Scenario}
\label{sec:734scenario}

\paragraph{Definition of the test scenario}

In the following, we want to assess the performance
of the sparse grid interpolants for the optimization problems
\ref{item:biomech2MinSum} and \ref{item:biomech2MinDist}.
For this goal, we create a test scenario \cite{Valentin18Gradient}
that simulates a pseudo-dynamic sequence of motions
by varying the load force and/or the target elbow angle
in discrete time steps $t$ as seen in \cref{fig:biomech2ScenarioA_1}.
The test scenario is as follows:
\begin{figure}
  \includegraphics{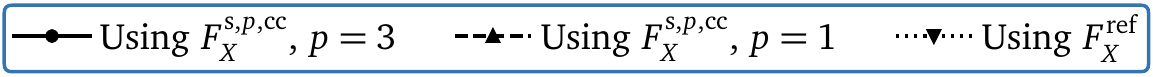}%
  \\[2mm]%
  \subcaptionbox{%
    Load $\forceL$ and target elbow angle $\tarelbang$.%
    \label{fig:biomech2ScenarioA_1}%
  }[73.3mm]{%
    \hspace*{4.5mm}%
    \includegraphics{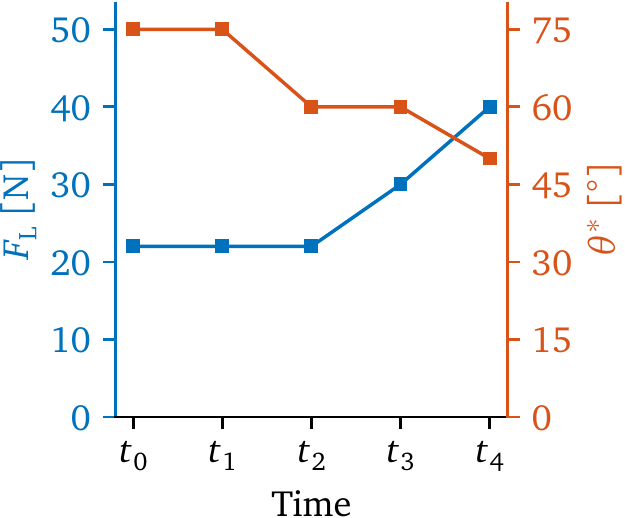}%
    \hspace*{2.4mm}%
  }%
  \hfill%
  \subcaptionbox{%
    Optimal activation parameters $\actT$ and $\actB$.%
    \label{fig:biomech2ScenarioA_2}%
  }[73.3mm]{%
    \hspace*{0.0mm}%
    \includegraphics{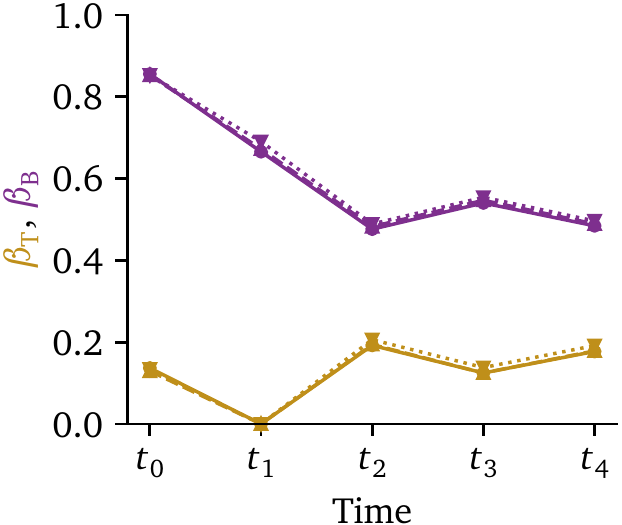}%
    \hspace*{9.8mm}%
  }%
  \\[1mm]%
  \subcaptionbox{%
    Deviation $\abs{\equielbang{\forceL} - \tarelbang}$
    of attained elbow angle to target and
    deviation \smash{$|\momentref|$} of the moment from equilibrium.%
    \label{fig:biomech2ScenarioA_3}%
  }[73.3mm]{%
    \includegraphics{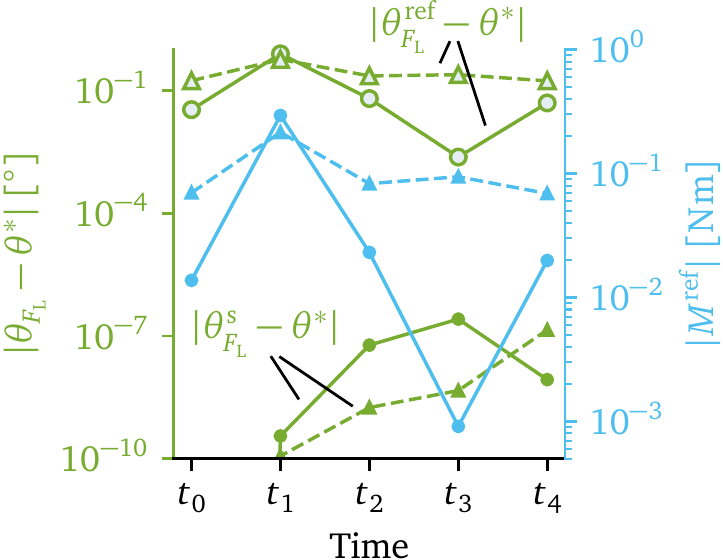}%
  }%
  \hfill%
  \subcaptionbox{%
    Number of evaluations of $\equielbang{\forceL}$
    and number of Newton iterations per evaluation
    of $\equielbang{\forceL}$.%
    \label{fig:biomech2ScenarioA_4}%
  }[73.3mm]{%
    \includegraphics{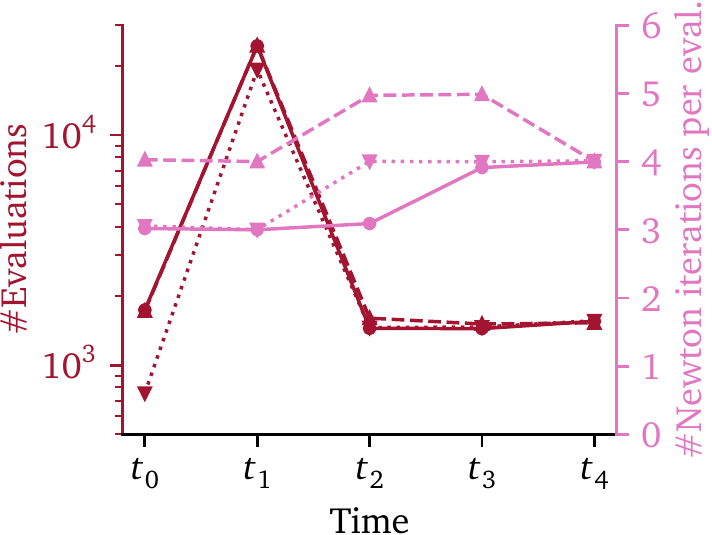}%
  }%
  \caption[Settings and results of the test scenario]{%
    Setting (a) of the test scenario and corresponding results (b, c, d).%
  }%
  \label{fig:biomech2ScenarioA}%
\end{figure}
\begin{enumerate}
  \item
  Find a feasible initial solution for problem \ref{item:biomech2MinSum}
  with $\forceL(t_0) \ceq \SI{22}{\newton}$ and
  $\tarelbang(t_0) \ceq \ang{75}$.
  
  \item
  Apply \ref{item:biomech2MinSum} with $\forceL(t_1) \ceq \SI{22}{\newton}$ and
  $\tarelbang(t_1) \ceq \ang{75}$.
  
  \item
  Apply \ref{item:biomech2MinDist} with $\forceL(t_2) \ceq \SI{22}{\newton}$ and
  $\tarelbang(t_2) \ceq \ang{60}$ (changed target angle).
  
  \item
  Apply \ref{item:biomech2MinDist} with $\forceL(t_3) \ceq \SI{30}{\newton}$ and
  $\tarelbang(t_3) \ceq \ang{60}$ (changed load).
  
  \item
  Apply \ref{item:biomech2MinDist} with $\forceL(t_4) \ceq \SI{40}{\newton}$ and
  $\tarelbang(t_4) \ceq \ang{50}$ (changed load and target angle).
\end{enumerate}
For each of the steps 2 to 5, the activation levels $\actT, \actB$ obtained
in the previous step (i.e., either the feasible initial solution
of step 1 or the optimal solution of steps 2 to 4) are used
as the input of the optimization problem
\ref{item:biomech2MinSum} or \ref{item:biomech2MinDist}.
The feasible initial solution in step 1 is determined as explained
in \cref{sec:513gradientBasedConstrained}.

\paragraph{Solutions of problem \ref{item:biomech2MinSum}}

We note that independently of $\forceL$ and $\tarelbang$,
every solution $(\actT, \actB)$ of problem \ref{item:biomech2MinSum} will be
on the boundary part of the domain $\clint{\*0, \*1}$,
on which at least one activation parameter vanishes, i.e.,
\begin{equation}
  \{(\actT, \actB) \in \clint{\*0, \*1} \mid
  (\actT = 0) \lor (\actB = 0)\}.
\end{equation}
The reason is that the two muscles triceps and biceps are antagonistic
(see \cref{sec:711models}), meaning that they work against each other.
If both $\actT > 0$ and $\actB > 0$, then the body will waste energy,
as the same target elbow angle can be attained by reducing both
$\actT$ and $\actB$ simultaneously, thus requiring less energy.
A visual example for this is \cref{fig:biomech2ReferenceEquilibriumAngle},
where the contour lines generally go from the bottom left
(small $\actT, \actB$) to the top right (large $\actT, \actB$).
This issue may be prevented by either
more complicated musculoskeletal models with more
than two muscles or different optimization problems
such as problem \ref{item:biomech2MinDist},
where the objective function differs.

\paragraph{Plots of optimization results}

\Cref{fig:biomech2ScenarioA_2,fig:biomech2ScenarioA_3,fig:biomech2ScenarioA_4}
show the results of the test scenario using the muscle forces
$\forceXintp[p,\cc]$ obtained by interpolating with
modified hierarchical cubic Clenshaw--Curtis B-splines (solid lines, $p = 3$).
As comparison, we repeat the solution process
with the forces obtained by interpolating with the
corresponding hierarchical piecewise linear basis (dashed lines, $p = 1$) and
with the reference forces $\forceXref$ (dotted lines).
For the piecewise linear basis,
we use exactly the same method as for the cubic case
(Newton method for $\equielbangintp{\forceL}$,
Augmented Lagrangian with adaptive gradient descent for the
solution of problems \ref{item:biomech2MinSum} and \ref{item:biomech2MinDist}),
although the derivatives of the muscle forces are discontinuous.
For the reference forces, we use the fact that the reference surrogates
are full grid spline interpolants, which can be explicitly differentiated.
Without the full grid interpolants,
we would have to approximate the derivatives with finite differences.

\vspace*{\fill}
\pagebreak

\paragraph{Equilibrium elbow angle}

In \cref{fig:biomech2ScenarioA_2}, we see that the activation levels
of all three methods are more or less the same.
However, \cref{fig:biomech2ScenarioA_3} reveals that even these small
differences lead to deviations of the resulting equilibrium elbow angle
to the target angle that differ by up to two orders of magnitude.
The two green lines with filled markers at the bottom of
\cref{fig:biomech2ScenarioA_3} show the error of
the equilibrium elbow angle $\equielbangintp{\forceL}$
using sparse grid interpolation to the desired target angle $\tarelbang$.
Unsurprisingly, this error is very small as
it is minimized by the optimizer as part of the constraint.
The true error, which is obtained by
using the reference equilibrium elbow angle $\equielbangref{\forceL}$,
is in general much larger
(top two green lines in \cref{fig:biomech2ScenarioA_3}
with hollow markers).
We see that the cubic B-splines decrease the error
by up to two orders of magnitude compared to the
piecewise linear basis.
There are two reasons for this:
First, the error of $\equielbang{\forceL}$ is generally smaller
when using higher-order B-splines as we have seen above.
Second, higher-order B-splines are continuously differentiable,
which makes them suitable for gradient-based optimization.
In contrast, the surrogates obtained by piecewise linear interpolation
have kinks, which may complicate finding optimal points
in the augmented Lagrangian and Newton methods.

\paragraph{Number of evaluations and Newton iterations}

This is supported by \cref{fig:biomech2ScenarioA_4},
which shows the number of evaluations of $\equielbang{\forceL}$
during the optimization and the average number of Newton iterations
per evaluation.
While the number of total evaluations is similar for all three methods,
the number of required Newton iterations to achieve convergence
is in general around \SI{50}{\percent} larger for the piecewise linear
basis functions.

\subsection{Spatial Adaptivity}
\label{sec:735adaptivity}

\paragraph{Generation of a spatially adaptive sparse grid}

As mentioned in \cite{Valentin18Gradient}, spatial adaptivity
may be employed to reduce the number of necessary muscle force samples
even further,
especially for more complicated musculoskeletal systems with
more parameters.
To verify this statement, we remove all grid points
$(\elbang^{(k,\cc)}, \actX^{(k,\cc)})$
from the regular sparse Clenshaw--Curtis grid that satisfy
\begin{equation}
  \frac{
    \abs{\alpha_\mathrm{T}^{(k,p,\cc)}}
  }{
    \max_{k'} \abs{\alpha_\mathrm{T}^{(k',p,\cc)}}
  } < \SI{1}{\percent}
  \quad\text{and}\quad
  \frac{
    \abs{\alpha_\mathrm{B}^{(k,p,\cc)}}
  }{
    \max_{k'} \abs{\alpha_\mathrm{B}^{(k',p,\cc)}}
  } < \SI{1}{\percent},
\end{equation}
where $\alpha_X^{(k,p,\cc)}$ ($X \in \{\mathrm{T}, \mathrm{B}\}$)
is the hierarchical surplus of the basis function $\bspl[\cc,\modified]{k}{p}$
corresponding to $(\elbang^{(k,\cc)}, \actX^{(k,\cc)})$.
For higher-dimensional models,
one would of course not sample muscle data on a regular sparse grid
and then coarsen the data by removing points,
but rather use an a posteriori adaptivity criterion to
decide which grid points to refine iteratively.

\paragraph{Comparison with the regular case}

\begin{figure}
  \hspace*{2mm}%
  \raisebox{0.2mm}{\includegraphics{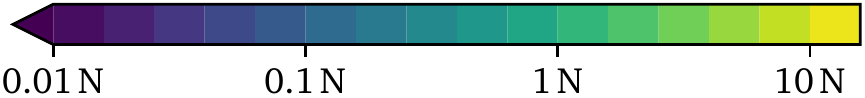}}%
  \hspace*{12mm}%
  \includegraphics{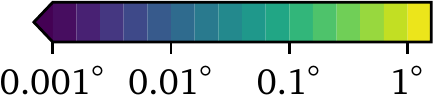}%
  \\[2mm]%
  \subcaptionbox{%
    $\abs{\forceTref - \forceTintp[p,\cc,\mathrm{adap}]}$%
  }[49mm]{%
    \raisebox{1.02mm}{\includegraphics{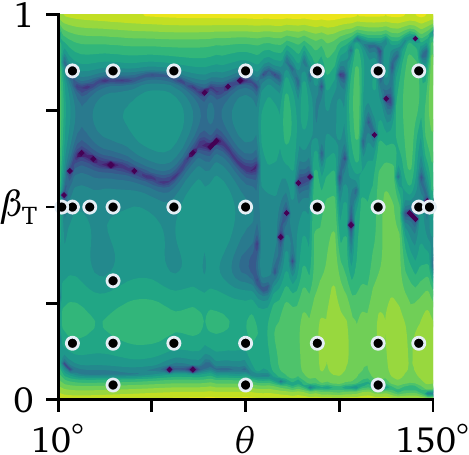}}%
  }%
  \hfill%
  \subcaptionbox{%
    $\abs{\forceBref - \forceBintp[p,\cc,\mathrm{adap}]}$%
  }[49mm]{%
    \raisebox{1.02mm}{\includegraphics{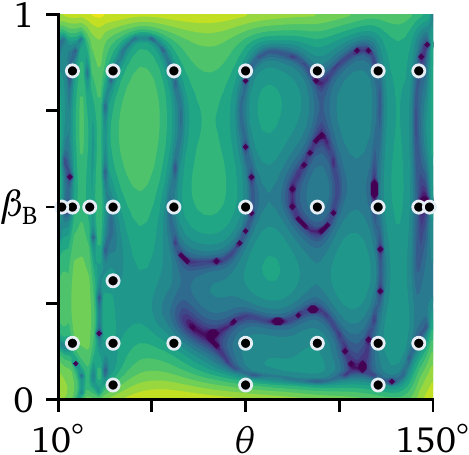}}%
  }%
  \hfill%
  \subcaptionbox{%
    $\abs{
      \equielbangref{\forceL} - \equielbangintp[p,\cc,\mathrm{adap}]{\forceL}
    }$
    for $\forceL = \SI{22}{\newton}$%
  }[49mm]{%
    \includegraphics{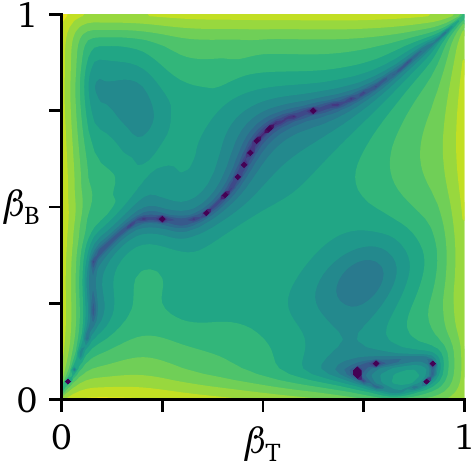}%
  }%
  \caption[%
    Errors of muscle forces and equilibrium angle
    for the spatially adaptive case%
  ]{%
    Errors of muscle forces and equilibrium elbow angle
    for the spatially adaptive case
    (modified hierarchical cubic Clenshaw--Curtis B-splines,
    i.e., $p = 3$) together with the points of the
    spatially adaptive sparse grid \emph{(dots).}%
  }%
  \label{fig:biomech2SpatiallyAdaptive}%
\end{figure}

For the cubic case $p = 3$,
the resulting force interpolants $\forceXintp[p,\cc,\mathrm{adap}]$ together
with the spatially adaptive sparse grid
(which has been coarsened from 49 to 28 points) and
equilibrium elbow angle $\equielbangintp[p,\cc,\mathrm{adap}]{\forceL}$
for $\forceL = \SI{22}{\newton}$ are shown in
\cref{fig:biomech2SpatiallyAdaptive}.
The sparse grid is almost dimensionally adaptive,
as $\forceXref$ seems to be almost linear in the $\actX$ direction
for both $X = \mathrm{T}$ and $X = \mathrm{B}$.
The errors increase slightly:
The relative $\Ltwo$ force errors for $(\mathrm{T}, \mathrm{B})$ increase
from $(\SI{3.31}{\permille}, \SI{3.56}{\permille})$
to $(\SI{3.36}{\permille}, \SI{4.43}{\permille})$,
and the absolute $\Linfty$ errors increase
from $(\SI{10.6}{\newton}, \SI{9.51}{\newton})$
to $(\SI{12.3}{\newton}, \SI{9.57}{\newton})$.
In addition, the relative $\Ltwo$ and absolute $\Linfty$ errors
for $\equielbang{\forceL}$ increase
from \SI{2.83}{\permille} and \ang{0.886}
to \SI{4.12}{\permille} and \ang{1.09}, respectively.
While all these errors are somewhat larger than for the regular sparse grid,
they are still at an acceptable level,
but the number of necessary muscle force evaluations is halved compared
to the regular case.
Additionally, the solution of the test scenario doesn't change significantly
due to the similar errors of $\forceX$ and $\equielbang{\forceL}$.

\cleardoublepage

  \setdictum{%
  The goal is to buy as many iPads as possible during your lifetime.%
}{%
  In a talk at the 5th Workshop on\\Sparse Grids and Applications%
}

\longchapter{%
  Application 3: Dynamic Portfolio Choice Models%
}{%
  Application 3:\texorpdfstring{\\}{ }Dynamic Portfolio Choice Models%
}{%
  Application 3 -- Dynamic Portfolio Choice Models%
}
\label{chap:80finance}

\initial[lhang=0.06]{0.1em}{S}{urrogates based on B-splines on sparse grids}
can also be used for our third application,
which stems from finance.
In this application, we optimize financial decisions of an individual
over their lifetime in discrete time steps or iterations $t = 0, \dotsc, T$
(for example, years $t = 0, \dotsc, 80$, where $20+t$ is the age
of the individual), depending on internal and external factors.
There are three types of variables:

\begin{itemize}
  \item
  \term{State variables} $\state_t$
  such as the individual's wealth $\wealth_t$ and their income
  cannot be controlled directly by the individual.
  Instead, the individual's decisions may influence the value of
  state variables of future iterations.%
  \footnote{%
    The time $t$ can also be regarded as a state variable.%
  }
  
  \item
  \term{Policy variables} $\policy_t$
  such as consumption $\consume_t$ and the amount of stocks to buy or sell
  represent the investment decisions the individual can make in
  each iteration.
  They are subject to specific constraints
  (for instance, you cannot spend more money than you have,
  if you do not allow debts).
  
  \item
  \term{Stochastic variables} $\stochastic_t$
  such as return rates of stocks and inflation
  cannot be controlled by the individual at all.
  Therefore, statements about optimal investment conditions
  are usually made for expected values instead of exact values.
\end{itemize}

\noindent
We discretize the state space with a spatially adaptive sparse grid.
For each state grid point, an optimization problem over the policy
variables has to be solved, where the objective function depends
on the expected value over the stochastic variables.
By using B-splines as hierarchical basis functions,
the accuracy of the interpolants is increased and
the explicitly known gradients enable the usage of
gradient-based optimization methods, thus accelerating convergence.
The process is repeated for each time step,
which is possible due to the Bellman principle,
which implies that the objective functions
occurring at time $t$ depend on the interpolant of the next iteration $t+1$.
Hence, the problem has to be solved backward in time
via a scheme that closely resembles dynamic programming.

The outline of this chapter is as follows:
In \cref{sec:81models}, we formalize the framework of
dynamic portfolio choice models and describe our approach.
Afterwards, we explain in \cref{sec:82algorithms} the necessary algorithms
for implementing the solution of these models.
\Cref{sec:83problem} introduces the transaction costs problem
as an example application of the general framework presented
in \cref{sec:81models}.
Finally, in \cref{sec:84results}, we study numerical results.

This chapter is based on a collaboration with Prof.\ Dr.\ Raimond Maurer
and Peter Schober (both Goethe University Frankfurt, Germany).
In previous work, Peter Schober treated the solution of
dynamic portfolio choice models with piecewise linear basis functions
on spatially adaptive sparse grids \cite{Schober18Solving}.
The original contribution of this thesis is the introduction
of higher-order B-splines for the solution of these problems.
The author of this thesis contributed the methodology of
hierarchical B-splines and large parts of the implementation.
The contributions of the collaborators at Goethe University Frankfurt
are the financial models, the literature review of related work,
and the assessment of the quality of results.

\section{Solving the Bellman Equation}
\label{sec:81models}

\minitoc[-0.5mm]{76mm}{4}

\noindent
In this section, we give a mathematical framework for
dynamic portfolio choice models,
briefly mention related literature, and
explain where B-splines on sparse grids come into play.
\Cref{tbl:glossaryFinance} summarizes the symbols
that are introduced in this chapter.
Rust provides a more detailed introduction to
dynamic portfolio choice models \cite{Rust18Dynamic}.

\begin{table}
  \setnumberoftableheaderrows{0}%
  \newcommand*{\vph}{%
    \vphantom{\printnotationsymbol{\buysell}}%
  }%
  \newcommand*{\pnst}[1]{%
    \printnotationsymbol{#1}\vph&\printnotationtext{#1}%
  }%
  \newcommand*{\pnsta}[1]{%
    \printnotationsymbol{#1}\vph&\multicolumn{3}{l}{\printnotationtext{#1}}%
  }%
  \begin{tabular}{%
    >{\kern\tabcolsep}=l+l<{\kern4.5mm}+l<{\kern-1.5mm}+l%
    <{\kern4.5mm}+l<{\kern-1mm}+l<{\kern\tabcolsep}%
  }
    \toprulec
    \pnst{t}&            \pnst{\wealth}&      \pnst{\utilityfcn}\\
    \pnst{\state}&       \pnst{\consume}&     \pnst{\statefcn}\\
    \pnst{\policy}&      \pnst{\bond}&        \pnst{\valuefcn}\\
    \pnst{\stochastic}&  \pnsta{\cetvalueintp}\\
    \pnst{\riskav}&      \pnst{\stock}&       \pnst{\optpolicyfcn}\\
    \pnst{\patience}&    \pnst{\buysell}&     $\hat{({\cdot})}$&Normalized quantity\\
    \pnst{\bondreturn}&  \pnst{\stockreturn}& \pnst{\wealthratio}\\
    \pnst{\tac}&         \pnsta{\weightedeulererror}\\
    \bottomrulec
  \end{tabular}%
  \caption[Glossary for dynamic portfolio choice models]{%
    Glossary of the notation for dynamic portfolio choice models.%
  }%
  \label{tbl:glossaryFinance}%
\end{table}

\subsection{Bellman Equation}
\label{sec:811bellmanEquation}

\paragraph{Utility maximization}

\begin{SCfigure}
  \includegraphics{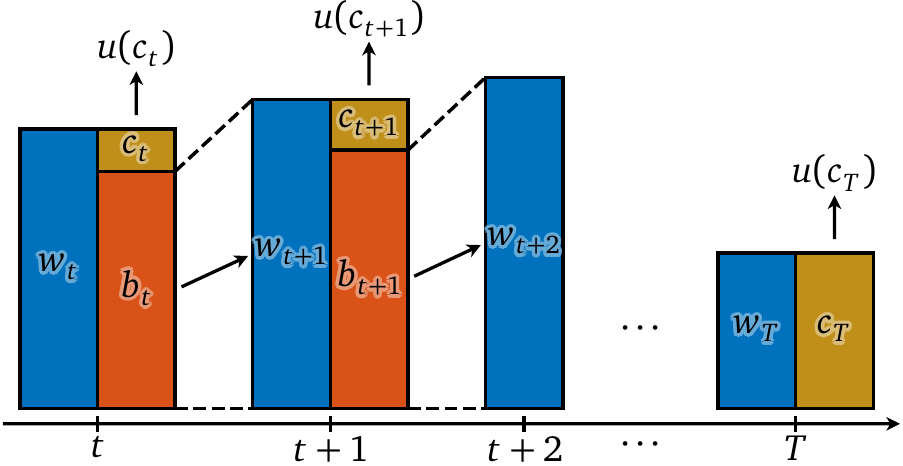}%
  \caption[Example of a dynamic portfolio choice model]{%
    Example of a dynamic portfolio choice model.
    The available wealth $\wealth_t$ is either
    invested into risk-free bonds ($\bond_t$) or consumed ($\consume_t$),
    resulting in utility $\utilityfcn(\consume_t)$.
    In the last time step $T$ \emph{(far right),}
    the optimal solution is to consume the whole wealth,
    if we do not take inheritance into account.%
  }%
  \label{fig:dynamicPortfolioChoice}%
\end{SCfigure}

\usenotation{t}
In the following, dynamic portfolio choice models aim to maximize the expected
\term{discounted time-additive utility} over the lifetime of the individual,
where the terminal utility is derived solemnly from consumption
(i.e., no inheritance motive).
If we neglect stochastic factors, then these models solve
\begin{equation}
  \label{eq:utilityMaximization}
  (\optpolicyfcn_0, \dotsc, \optpolicyfcn_T)
  = \vecargmax_{\policy_0, \dotsc, \policy_T}
  \sum_{t=0}^T \patience^t \utilityfcn(\consume_t(\state_t, \policy_t))
  \quad\text{s.t. specific constraints.}
\end{equation}
Here, $\state_t \in \clint{\*0, \*1} \subset \real^d$ and
$\policy_t \in \real^{m_{\policy}}$
are the state%
\footnote{%
  We assume that each state variable $\stateentry_{t,o}$ ($o = 1, \dotsc, d$)
  is bounded, since the state space will be discretized with sparse grids.
  Without loss of generality, we may then assume that
  $\state_t \in \clint{\*0, \*1}$.
  If some state variables are unbounded in reality,
  then extrapolation is necessary,
  which will be explained in \cref{sec:825interpolation}.%
}
and policy of time $t = 0, \dotsc, T$, respectively.
The constraints ensure that for instance, we do not spend more money
than we actually have.
Starting from a given initial state $\state_0$,
the state $\state_{t+1}$ of time $t+1$ can be computed from
$\state_0$ and $\policy_0, \dotsc, \policy_t$ with a
\term{state transition function} $(\state_t, \policy_t) \mapsto \state_{t+1}$.
As shown in \cref{fig:dynamicPortfolioChoice},
in each time step, a fraction of the available wealth is consumed
\term{(consumption $\consume_t$),}
which can be computed from the state $\state_t$ and the policy $\policy_t$.
The individuals rate the consumption with a \term{utility function}
$\utilityfcn(\consume_t)$.
A common choice for $\utilityfcn$ is the \term{\crra utility}
$\utilityfcn(\consume_t) \ceq c_t^{1-\riskav}/(1-\riskav)$
with the \term{risk aversion} $\riskav \in \real \setminus \{1\}$.
Positive and negative values of $\riskav$
correspond to risk-averse and risk-affine individuals, respectively.
The factor $\patience \in \pohint{0, 1}$ is the \term{patience}
or \term{time discount factor.}

\paragraph{Limitations of naive utility maximization}

When solving the utility maximization problem
in \cref{eq:utilityMaximization}, there are two issues.
First, solving \cref{eq:utilityMaximization} for all times $t$ at once
implies solving a $(T+1) m_{\policy}$-dimensional optimization problem,
which is usually computationally infeasible.
Second, \cref{eq:utilityMaximization} does not take stochastic variables
$\stochastic_t$ such as stock return rates into account.
These variables influence the state transition, i.e.,
$(\state_t, \policy_t, \stochastic_t) \mapsto \state_{t+1}$.
Consequently, $\state_{t+1}$ cannot be computed from $\state_0$ and
$\policy_0, \dotsc, \policy_t$ alone,
which complicates the solution of \cref{eq:utilityMaximization}
even for expected values.

\paragraph{Bellman principle}

To resolve the first issue,
Bellman's principle of optimality \cite{Bellman57Dynamic}
can be applied to problems like
\cref{eq:utilityMaximization} that are said to have
\term{optimal substructure.}
The principle states that the optimal policy for all times $t = 0, \dotsc, T$
is also optimal with respect to $t = 1, \dotsc, T$, i.e.,
{%
  \setlength{\abovedisplayskip}{6pt}%
  \begin{equation}
    \max_{\policy_0, \dotsc, \policy_T}
    \sum_{t=0}^T \patience^t \utilityfcn(\consume_t(\state_t, \policy_t))
    = \max_{\policy_0} \left(
      \utilityfcn(\consume_0(\state_0, \policy_0))
      + \patience \max_{\policy_1, \dotsc, \policy_T}
      \sum_{t=1}^T \patience^{t-1} \utilityfcn(\consume_t(\state_t, \policy_t))
    \right),
  \end{equation}%
}%
where we omitted the constraints for brevity.
The inner maximum problem over $\policy_1, \dotsc, \policy_T$
has the same structure as the problem on the \lhs.
With the \term{value function}
$\valuefcn_t\colon \clint{\*0, \*1} \to \real$,
$\valuefcn_t(\state_t) \ceq
\max_{\policy_t, \dotsc, \policy_T}
\sum_{t'=t}^T \patience^{t'-t}
\utilityfcn(\consume_{t'}(\state_{t'}, \policy_{t'}))$,
this can be rewritten as
{%
  \setlength{\belowdisplayskip}{9pt}%
  \begin{equation}
    \label{eq:simpleBellman}
    \valuefcn_0(\state_0)
    = \max_{\policy_0} \left(
      \utilityfcn(\consume_0(\state_0, \policy_0)) +
      \patience \valuefcn_1(\state_1)
    \right)
    \quad\text{s.t. specific constraints,}
  \end{equation}%
}%
where $\state_1$ is the result of the state transition
starting from $(\state_0, \policy_0)$.

\paragraph{General Bellman equation}

If we formulate \cref{eq:simpleBellman} for arbitrary times $t$ and
consider constraints, state transition, and stochastic variables,
we obtain the \term{Bellman equation:}
\begin{subequations}
  \setlength{\abovedisplayskip}{9pt}%
  \label{eq:generalBellman}
  \begin{gather}
    \valuefcn_t(\state_t)
    = \max_{\policy_t} \left(
      \utilityfcn(\consume_t(\state_t, \policy_t)) +
      \patience \expectation[t]{
        \valuefcn_{t+1}(\statefcn_t(\state_t, \policy_t, \stochastic_t))
      }
    \right),\quad
    t = 0, \dotsc, T,\\[-2mm]
    \policy_t \in \real^{m_{\policy}}\;\;\text{s.t.}\;\;
    \ineqconfun_t(\state_t, \policy_t) \le \*0,
  \end{gather}
\end{subequations}
where $\valuefcn_{T+1} :\equiv 0$ for simplicity,
$\statefcn_t\colon \clint{\*0, \*1} \times \real^{m_{\policy}} \times
\stochdomain \to \clint{\*0, \*1}$,
$(\state_t, \policy_t, \stochastic_t) \mapsto \state_{t+1}$,
is the \term{state transition function,}
$\ineqconfun_t\colon \clint{\*0, \*1} \times \real^{m_{\policy}} \to
\real^{m_{\ineqconfun}}$ is the \term{constraint function,} and
\begin{equation}
  \expectation[t]{
    \valuefcn_{t+1}(\statefcn_t(\state_t, \policy_t, \stochastic_t))
  }
  \ceq \int_{\stochdomain}
  \valuefcn_{t+1}(\statefcn_t(\state_t, \policy_t, \stochastic_t))
  P_{t,\stochastic}(\stochastic_t) \diff{}\stochastic_t
\end{equation}
with the probability density function
$P_{t,\stochastic}\colon \stochdomain \to \nonnegreal$ of $\stochastic_t$.%
\footnote{%
  While the state $\state_t \in \clint{\*0, \*1}$
  is continuous in this thesis,
  \term{Markov-chain discrete states} $\discrstate_t \in \discrstdomain$
  such as alive/dead
  (i.e., $\discrstdomain$ is the Cartesian product of finite sets)
  can be incorporated into \eqref{eq:generalBellman}.
  The objective function of $\valuefcn_t(\state_t, \discrstate_t)$ then equals
  $\utilityfcn(\consume_t(\state_t, \discrstate_t, \policy_t)) +
  \patience \expectation[t]{
    \valuefcn_{t+1}(\statefcn_t(\state_t, \discrstate_t,
    \policy_t, \stochastic_t), \discrstate_{t+1}) \mid \discrstate_t
  }$.%
}
We denote the location of the maximum of \eqref{eq:generalBellman}
as the optimal policy $\optpolicyfcn_t$,
which may be regarded as a function
$\optpolicyfcn_t\colon \clint{\*0, \*1} \to \real^{m_{\policy}}$,
$\state_t \mapsto \optpolicyfcn_t(\state_t)$.

\paragraph{Dynamic programming scheme}

The Bellman equation \eqref{eq:generalBellman} can be solved
backwards in time with a dynamic programming scheme.
Starting from the solution $\valuefcn_T$ and $\optpolicyfcn_T$
of time $T$, which is determined by maximizing the utility
for the terminal time step,
we can determine $\valuefcn_t$ and $\optpolicyfcn_t$
from $\valuefcn_{t+1}$ and $\optpolicyfcn_{t+1}$
for $t = T - 1,\, T - 2,\, \dotsc,\, 0$ with the Bellman equation.
This way, we only have to solve $T+1$ separate $m_{\policy}$-dimensional
optimization problems instead of a single large
$(T+1) m_{\policy}$-dimensional problem.
Often, the terminal solutions $\valuefcn_T$ and $\optpolicyfcn_T$
are explicitly known.
In our case, the optimal terminal solution is
to consume the whole wealth $\wealth_T$
(see \cref{fig:dynamicPortfolioChoice}).

\paragraph{Implementation and interpolation}

For the implementation of \eqref{eq:generalBellman},
we discretize the state space $\clint{\*0, \*1}$ into
$\ngp_t$ grid points $\state_t^{(k)}$, $k = 1, \dotsc, \ngp_t$,
and we tabulate the values of $\valuefcn_t$ and $\optpolicyfcn_t$
at $\state_t^{(k)}$ for all $t = 0, \dotsc, T$ and $k = 1, \dotsc, \ngp_t$.
However, in general, the next state
$\statefcn_t(\state_t^{(k)}, \policy_t, \stochastic_t)$ does not correspond
to a grid point $\state_{t+1}^{(k')}$,
which means that we cannot lookup the value of $\valuefcn_{t+1}$
at $\statefcn_t(\state_t^{(k)}, \policy_t, \stochastic_t)$.
Therefore, we have to interpolate $\valuefcn_{t+1}$
at the grid points, obtaining the interpolant $\valueintp_{t+1}$
as a result:
\begin{equation}
  \label{eq:gridBellman}
  \valueintp_t(\state_t^{(k)})
  = \max_{\policy_t} \left(
    \utilityfcn(\consume_t(\state_t^{(k)}, \policy_t)) +
    \patience \expectation[t]{
      \valueintp_{t+1}(\statefcn_t(\state_t^{(k)}, \policy_t, \stochastic_t))
    }
  \right),\quad
  k = 1, \dotsc, \ngp_t,
\end{equation}
where $\valueintp_{T+1} :\equiv 0$ for simplicity.
As $\valueintp_{t+1}$ on the \rhs is only an approximation to
$\valuefcn_{t+1}$, the values $\valueintp_t(\state_t^{(k)})$
on the \lhs are approximations, too.
Since we are mainly interested in the optimal policy decisions
$\optpolicyfcn_t$, we have to interpolate them as well, i.e.,
\begin{equation}
  \optpolicyintp_t(\state_t^{(k)})
  = \vecargmax_{\policy_t} \left(
    \utilityfcn(\consume_t(\state_t^{(k)}, \policy_t)) +
    \patience \expectation[t]{
      \valueintp_{t+1}(\statefcn_t(\state_t^{(k)}, \policy_t, \stochastic_t))
    }
  \right).
\end{equation}
Note that the employed grids for $\optpolicyintp_t$
may be different from the grids for $\valueintp_t$.

\subsection{Solution with B-Spline Surrogates on Sparse Grids}
\label{sec:812surrogates}

\paragraph{Sparse grids for dynamic models and related work}

As interpolation approaches for $\valuefcn_t$ based on full grids
suffer from the curse of dimensionality,
we want to use interpolation on spatially adaptive sparse grids instead.
Recently, sparse grids have found increasing interest in the
solution of dynamic models in finance
\multicite{Brumm17Using,Judd14Smolyak,Schober18Solving,Winschel10Solving}.
For example in \cite{Brumm17Using},
discrete choices in the value iteration are computed
using piecewise linear basis functions on spatially adaptive sparse grids.
Schober employs spatially adaptive sparse grids
for the interpolation of dynamic portfolio choice models,
but uses piecewise linear basis functions \cite{Schober18Solving}.
Judd et al.\ use global polynomials
on sparse Clenshaw--Curtis grids for the interpolation of
higher-dimensional economic models \cite{Judd14Smolyak}.

\paragraph{B-splines on sparse grids for dynamic portfolio choice models}

The shortcomings of the two approaches of piecewise linear functions
\multicite{Brumm17Using,Schober18Solving} or global polynomials
\cite{Judd14Smolyak} are evident:
Piecewise linear functions are not continuously differentiable,
impeding convergence of interpolation errors
(see \cref{sec:541interpolation}) and prohibiting the use
of gradient-based optimization methods to solve \cref{eq:gridBellman}.
The reason for the latter statement is that gradient-based optimizers require
the derivatives of the objective function of \cref{eq:gridBellman}
with respect to the entries $\policyentry_{t,j}$ of $\policy_t$
($j = 1, \dotsc, m_{\policy}$), i.e.,
\begin{equation}
  \begin{split}
    &\utilityfcn'(\consume_t(\state_t^{(k)}, \policy_t))
    \partialderiv{\partialdiff{}\policyentry_{t,j}}{c_t}\paren{
      \state_t^{(k)}, \policy_t
    }\\
    &{} + \patience \expectation[t]{
      \tr{
        \paren*{
          \gradient{\state_{t+1}}{\valueintp_{t+1}}\paren{
            \statefcn_t(\state_t^{(k)}, \policy_t, \stochastic_t)
          }
        }
      }
      \partialderiv{\partialdiff{}\policyentry_{t,j}}{\statefcn_t}\paren{
        \state_t^{(k)}, \policy_t, \stochastic_t
      }
    },
  \end{split}
\end{equation}
which involves the gradient
$\gradient{\state_{t+1}}{\valueintp_{t+1}}$ of the
value function interpolant $\valueintp_{t+1}$.
Gradient-based optimization methods do not converge fast
if this gradient is discontinuous.
Moreover, piecewise linear basis functions introduce
many additional local minima.
In contrast, global polynomials only work well on Clenshaw--Curtis grids
with Chebyshev-distributed nodes due to Runge's phenomenon.

In the following, we use higher-order B-splines as basis functions for
the interpolation of $\valuefcn_t$ and $\optpolicyfcn_t$.
This method has two advantages:
First, B-splines of degree $p > 1$ are continuously differentiable,
increasing the order of convergence and enabling
gradient-based optimization for solving \cref{eq:gridBellman}.
Second, B-splines are defined for arbitrary knot sequences,
leading to a greater flexibility when compared to global polynomials.

\breakpagebeforenextheadingtrue
\section{Algorithms}
\label{sec:82algorithms}

\minitoc{83mm}{7}

\noindent
This section gives an overview of the algorithms that
we use to implement the solution process after discretization of
the Bellman equation \eqref{eq:gridBellman}.
In the following, we assume that the probability density functions of the
stochastic variables are known.

\subsection{General Structure}
\label{sec:821generalStructure}

The general approach to solve dynamic portfolio choice models is as follows:
\begin{enumerate}
  \item
  Generation of value function interpolants $\valueintp_t$
  
  \item
  Generation of optimal policy interpolants $\optpolicyintp_t$
  
  \item
  Post-processing, e.g., Monte Carlo simulation
\end{enumerate}
The separation of the solution processes for
the value function interpolants $\valueintp_t$
and the optimal policy interpolants $\optpolicyintp_t$
enables the generation of different spatially adaptive sparse grids
for the value function and the optimal policies.
This is useful if the shapes of value function and optimal policies
have different characteristics.

In the following \cref{%
  sec:822solveValueFunction,%
  sec:823optimization,%
  sec:824quadrature,%
  sec:825interpolation,%
  sec:826gridGeneration%
}, we describe the algorithmic details of
\texttt{solveValueFunction} (step 1).
The treatment of the other steps \texttt{solvePolicy} (step 2) and
post-processing (step 3) follows with
\cref{sec:827solvePolicies} and \cref{sec:828postProecessing},
respectively.

We track two interpolants $\valueintp[1]_t$ and $\valueintp[p]_t$
for each $t = 0, \dotsc, T$.
The former interpolates value function data at the grid points
with the hierarchical piecewise linear basis
(used for the surplus-based grid generation),
while the latter interpolates the same data with hierarchical B-splines
of degree $p > 1$.
Each $\valueintp[\ast]_t$ ($\ast \in \{1, p\}$)
additionally stores the grid points $\state_t^{(k)}$
and the optimal policies $\optpolicyintp_t(\state_t^{(k)})$
at the grid points ($k = 1, \dotsc, \ngp_t$).
For simplicity, we do not pass them as separate data
to the algorithms.

\subsection{Solution for the Value Function}
\label{sec:822solveValueFunction}

\paragraph{\texttt{solveValueFunction} algorithm}

\Cref{alg:financeSolveValueFunction} shows \texttt{solveValueFunction},
which generates the value function interpolants
$\valueintp[1]_t$ and $\valueintp[p]_t$ ($t = 0, \dotsc, T$).
The algorithm follows a simple optimize--refine--interpolate scheme,
which is visualized in \cref{fig:structureSolveValueFunction}:
First, the Bellman equation \eqref{eq:gridBellman} is solved
on an initial sparse grid (\texttt{optimize}).
Then, we \texttt{refine} the grid spatially adaptively.
Finally, the resulting grid point data are \texttt{interpolate}d
with hierarchical higher-order B-splines.

\begin{algorithm}
  \begin{algorithmic}[1]
    \Function{%
      $\text{%
        $(\valueintp[p]_t)_{t=0,\dotsc,T}$%
      } = \texttt{solveValueFunction}$%
    }{%
      \hspace*{0mm}%
    }
      \State{$\valueintp[p]_{T+1} \gets \emptyset$}
      \Comment{dummy variable (is not used)}%
      \For{$t = T,\, T - 1,\, \dotsc,\, 0$}
        \State{%
          $\valueintp[1]_t \gets \text{%
            Initial regular sparse grid with no values%
          }$%
        }
        \State{%
          $\valueintp[1]_t \gets
          \texttt{optimize($t$, $\valueintp[1]_t$, $\valueintp[p]_{t+1}$)}$%
        }
        \State{%
          $\valueintp[1]_t \gets
          \texttt{refine($t$, $\valueintp[1]_t$, $\valueintp[p]_{t+1}$)}$%
        }
        \State{%
          $\valueintp[p]_t \gets
          \texttt{interpolate($\valueintp[1]_t$)}$%
        }
      \EndFor{}
    \EndFunction{}
  \end{algorithmic}
  \caption[%
    Generation of value function interpolants (\texttt{solveValueFunction})%
  ]{%
    Generation of value function interpolants.
    The output is the higher-order B-spline interpolant $\valueintp[p]_t$
    for all $t = 0, \dotsc, T$.%
  }%
  \label{alg:financeSolveValueFunction}%
\end{algorithm}

\begin{figure}
  \newcommand*{\getscale}{
    \pgfgettransformentries{%
      \myxscaletmp%
    }{\@tempa}{\@tempa}{\myyscaletmp}{\@tempa}{\@tempa}
    \gdef\myxscale{\myxscaletmp}
  }%
  \newcommand*{\drawcircle}[2][]{
    \node[
      circle,minimum size=60mm,fill=mittelblau!#2,draw=mittelblau,
    ] at (0mm,0mm) (myCircle#2) {};
    \IfStrEq{#1}{noArrow}{}{
      \getscale
      \pgfmathparse{150+20/\myxscale^0.8}
      \let\myin\pgfmathresult
      \draw[->] (myCircle#2) to[
        out=150,in=\myin,looseness=2.3,
      ] node[below,scale=1/\myxscale] {\mytext} (myCircle#2);
    }
  }%
  \newcommand*{\drawtext}[1]{
    \getscale
    \centerarc[
      draw=none,postaction={decorate},decoration={
        text along path,
        text={#1},
        text align=center
      },
    ](0mm,0mm)(180:0:30mm-1em/\myxscale);
  }%
  \newenvironment*{circlescope}{
    \getscale
    \begin{scope}[
      scale=1-0.15/\myxscale,
      shift={(3mm/\myxscale,-3mm/\myxscale)},
      transform shape,
    ]
  }{
    \end{scope}
  }%
  \begin{tikzpicture}
    \draw[dashed] (0mm,9mm) -- (84mm,9mm);
    \draw[dashed] (-4.3mm,-8mm) -- (71.2mm,-48.2mm);
    \begin{scope}[rotate=-130]
      \node[
        circle,
        minimum size=18mm,
        fill=mittelblau!30,
        draw=mittelblau,
        inner sep=0mm,
      ] at (0mm,0mm) (optimize) {\texttt{optimize}};
      \begin{scope}[
        shift={(40mm/2,40mm*sin(60))},
        transform shape,
      ]
        \begin{scope}[
          rotate=130,
          scale=18mm/60mm+0.15,
          transform shape,
        ]
          \getscale
          \drawcircle[noArrow]{30}
          \drawtext{|\ttfamily|refine}
          \begin{scope}[
            scale=1-0.15/\myxscale,
            shift={(3mm/\myxscale,-3mm/\myxscale)},
            transform shape,
          ]
            \def\mytext{}
            \drawcircle{50}
            \getscale
            \node[scale=1/\myxscale] at (0mm,0mm) {\texttt{optimize}};
          \end{scope}
        \end{scope}
      \end{scope}
      \node[
        circle,
        minimum size=18mm,
        fill=mittelblau!30,
        draw=mittelblau,
        inner sep=0mm,
        text width=15mm,
        align=center,
      ] at (40mm,0mm) (interpolate) {\ttfamily{}inter-\\polate};
      \draw[->] (optimize) to[
        bend left=30,
      ] (myCircle30);
      \draw[->] (myCircle30) to[
        bend left=30,
      ] (interpolate);
      \draw[->] (interpolate) to[
        bend left=30,
      ] node[above,rotate=50] {$t \gets (t - 1)$} (optimize);
    \end{scope}
    \begin{scope}[shift={(84mm,-21mm)}]
      \drawcircle[noArrow]{20}
      \drawtext{|\ttfamily|optimize}
      \begin{circlescope}
        \def\mytext{\hspace*{5mm}$\state_t$}
        \drawcircle{30}
        \drawtext{|\ttfamily|optimizeSinglePoint}
        \begin{circlescope}
          \def\mytext{\hspace*{5mm}$\policy_t$}
          \drawcircle{40}
          \drawtext{|\ttfamily|evalObjFcnGrad}
          \begin{circlescope}
          \def\mytext{\hspace*{5mm}$\stochastic_t$}
            \drawcircle{50}
            \drawtext{|\ttfamily|evalQuadPoint}
            \begin{circlescope}
              \def\mytext{\hspace*{4mm}$o$}
              \drawcircle{60}
              \getscale
              \node[
                text width=30mm,align=center,scale=1/\myxscale,
              ] at (0mm,0mm) {%
                \ttfamily{}evalInterp-\\PartDeriv%
              };
            \end{circlescope}
          \end{circlescope}
        \end{circlescope}
      \end{circlescope}
    \end{scope}
  \end{tikzpicture}%
  \caption[Scheme of the generation of value function interpolants]{%
    Scheme of the generation of value function interpolants
    with \texttt{solveValueFunction}
    (\cref{alg:financeSolveValueFunction}, \emph{left}),
    which repeatedly calls the \texttt{optimize} routine
    (\cref{alg:financeOptimize}, \emph{right}),
    which in turn consists of various sub-functions.
    The function \texttt{optimize} iterates over all state grid points
    $\state_t = \state_t^{(k)}$ ($k = 1, \dotsc, \ngp_t$)
    and calls \texttt{optimizeSinglePoint} for each point.
    The optimization method evaluates the objective function and
    its gradient at a sequence of different policy points $\policy_t$
    to find $\optpolicyintp_t(\state_t^{(k)})$.
    This evaluation (denoted by \texttt{evalObjFcnGrad})
    has to compute the expectation in \cref{eq:gridBellman},
    which is done using a quadrature rule.
    For every quadrature point $\stochastic_t = \stochastic_t^{(j)}$
    ($j = 1, \dotsc, m_{\quadweight}$),
    \texttt{evalQuadPoint} computes the corresponding value of
    the expression in the expectation.
    Finally, \texttt{evalInterpPartDeriv} evaluates the interpolant
    $\valueintp[p]_{t+1}$ and its partial derivatives,
    for which we have to loop over the state dimensions $o = 1, \dotsc, d$.%
  }%
  \label{fig:structureSolveValueFunction}%
\end{figure}

At the beginning of every iteration $t$,
the grid of the piecewise linear interpolant is reset
to an initial, possibly regular sparse grid.
It would also be possible to reuse the grid from the
previous iteration $t + 1$.
Nevertheless, the results we then obtain become worse,
likely due to the different characteristics of $\valueintp[1]_t$
for different $t$ (e.g., kinks).

The higher-order B-spline interpolant
$\valueintp[p]_{t+1}$ of the previous iteration $t+1$ is used
for the \rhs of the Bellman equation \eqref{eq:gridBellman},
if $t < T$.
In the first iteration $t = T$,
there is no such interpolant.
However,
the terminal solution $\valuefcn_T$ is usually a known function.

\subsection{Optimization}
\label{sec:823optimization}

\paragraph{\texttt{optimize} algorithm}

The \texttt{optimize} step is given as \cref{alg:financeOptimize}.
The grid of the argument $\valueintp[1]_t$
is some spatially adaptive sparse grid
$\gridset{t}{\sparse}
= \{\state_t^{(k)} \mid k = 1, \dotsc, \ngp_t\}$,
where the function values $\valueintp[1]_t(\state_t^{(k)})$
may already be known for some grid points $\state_t^{(k)}$,
if \texttt{optimize} is called from within \texttt{refine}.
The function \texttt{optimize} computes the missing value function values.
For $t = T$, we assume that the terminal solution
$\valuefcn_T$ can be computed by some function
\texttt{computeKnownTerminalSolution}.%
\footnote{%
  In any case, the terminal solution may be computed as the
  solution of the corresponding single-time optimization problem,
  e.g., $\valuefcn_T(\state_T^{(k)})
  = \max_{\policy_T} \utilityfcn(\consume_T(\state_T^{(k)}, \policy_T))$.%
}
Otherwise, for $t < T$, we solve the Bellman equation
\eqref{eq:gridBellman} by using the higher-order B-spline interpolant
$\valueintp[p]_{t+1}$ of the previous iteration $t + 1$
(\texttt{optimizeSinglePoint}).
The computations for the different $\state_t^{(k)}$ are
independent of each other,
which means that they can be computed in parallel \cite{Horneff16Efficient}.%
\footnote{%
  Such a problem is usually referred to as \term{embarrassingly parallel.}%
}
After generating all missing data,
we update the hierarchical surpluses of the
piecewise linear interpolant $\valueintp[1]_t$
to interpolate the new data at all grid points of $\gridset{t}{\sparse}$.

\begin{algorithm}
  \begin{algorithmic}[1]
    \Function{$\valueintp[1]_t = \texttt{optimize}$}{%
      $t$, $\valueintp[1]_t$, $\valueintp[p]_{t+1}$%
    }
      \State{%
        $(\state_t^{(k)})_{k=1,\dotsc,\ngp_t}
        \gets \text{grid of $\valueintp[1]_t$}$%
      }
      \For{$k = 1, \dotsc, \ngp_t$}
        \If{$\valueintp[1]_t(\state_t^{(k)})$ not previously computed}
          \IfOneLine{$t = T$}{%
            $\valueintp[1]_T(\state_T^{(k)})
            \gets \texttt{computeKnownTerminalSolution($\state_T^{(k)}$)}$%
          }
          \ElseOneLine{%
            $\valueintp[1]_t(\state_t^{(k)})
            \gets \texttt{%
              optimizeSinglePoint(%
                $t$, $\state_t^{(k)}$, $\valueintp[p]_{t+1}$%
              )%
            }$%
          }
        \EndIf{}
      \EndFor{}
      \State{%
        Re-interpolate
        $(\valueintp[1]_t(\state_t^{(k)}))_{k=1,\dotsc,\ngp_t}$
        with piecewise linear functions%
      }
    \EndFunction{}
  \end{algorithmic}
  \caption[Evaluation of the value function (\texttt{optimize})]{%
    Evaluation of the value function at all
    grid points $\state_t^{(k)}$ of $\valueintp[1]_t$
    at which the value function has not been evaluated yet.
    Inputs are
    the time $t$,
    the piecewise linear interpolant $\valueintp[1]_t$
    of the current iteration $t$ (with the underlying sparse grid and
    corresponding function values, possibly unset), and
    the higher-order B-spline interpolant $\valueintp[p]_{t+1}$
    of the previous iteration $t + 1$
    (not used if $t = T$).
    The output is the updated piecewise linear interpolant $\valueintp[1]_t$,
    where all missing function values at grid points have been computed.%
  }%
  \label{alg:financeOptimize}%
\end{algorithm}

\paragraph{Certainty-equivalent transformation}

For utility functions of \crra-type, i.e., of the form
$\utilityfcn(\consume_t) = \consume_t^{1-\riskav}/(1-\riskav)$,
the curvature of the objective function in the Bellman equation
\eqref{eq:gridBellman} can be very high
(depending on the risk aversion parameter $\riskav$),
which may impede convergence of the optimizer.
As a remedy, we transform the value function $\valueintp_t$ with the
\term{certainty-equivalent transformation}
$\valueintp_t \mapsto \cetvalueintp_t
\ceq ((1 - \riskav) \valueintp_t)^{1/(1-\riskav)}$ if $\riskav > 1$.
\Cref{eq:gridBellman} then becomes
$\cetvalueintp_T(\state_T^{(k)}) = \max_{\policy_T}
\consume_T(\state_T^{(k)}, \policy_T)$ for $t = T$ and
\begin{equation}
  \label{eq:gridBellmanCET}
  \cetvalueintp_t(\state_t^{(k)})
  = \max_{\policy_t} \left(
    \left(
      \consume_t(\state_t^{(k)}, \policy_t)^{1-\riskav} +
      \patience \expectation[t]{
        \bigl(
          \cetvalueintp_{t+1}(
            \statefcn_t(\state_t^{(k)}, \policy_t, \stochastic_t)
          )
        \bigr)^{1-\riskav}
      }
    \right)^{1/(1-\riskav)}
  \right)
\end{equation}
for $t < T$, since for $\riskav > 1$,
$(\cdot)^{1/(1-\riskav)}$ is strictly monotonously decreasing and
$(1-\riskav) < 0$.
The notation in the remainder of this section
does not distinguish between $\valueintp[\ast]_t$ and $\cetvalueintp[\ast]_t$
and uses $\valueintp[\ast]_t$ for both if it is not relevant
whether the value function is transformed.

\subsection{Quadrature}
\label{sec:824quadrature}

We need to approximate the expectation in \cref{eq:gridBellmanCET}
by quadrature,
\begin{equation}
  \label{eq:bellmanQuadrature}
  \expectation[t]{
    \bigl(
      \cetvalueintp[p]_{t+1}(
        \statefcn_t(\state_t^{(k)}, \policy_t, \stochastic_t)
      )
    \bigr)^{1-\riskav}
  }
  \approx \sum_{j=1}^{m_{\quadweight}} \quadweight_t^{(j)}
  \bigl(
    \cetvalueintp[p]_{t+1}(
      \statefcn_t(\state_t^{(k)}, \policy_t, \stochastic_t^{(j)})
    )
  \bigr)^{1-\riskav},
\end{equation}
for some weights $\quadweight_t^{(j)} \in \real$ and
nodes $\stochastic_t^{(j)} \in \stochdomain$
($j = 1, \dotsc, m_{\quadweight}$).
Since the stochastic domain $\stochdomain \subset \real^{m_{\stochastic}}$
might be high-dimensional as well,
full grid quadrature rules suffer from the curse of dimensionality.
Therefore, we use sparse grid quadrature rules based
on Gauss--Hermite quadrature
\multicite{Gerstner98Numerical,Horneff16Efficient}.
Note that this sparse grid in the stochastic space $\stochdomain$
is independent of the sparse grid in the state space $\clint{\*0, \*1}$.
However, it would also be feasible to employ Monte Carlo quadrature,
albeit usually far more expensive.

\subsection{Interpolation and Extrapolation}
\label{sec:825interpolation}

\paragraph{Sparse grid interpolation}

As already mentioned,
$\valueintp[1]_t$ is constructed as the sparse grid interpolant
of the grid data $\state_t^{(k)}$ ($k = 1, \dotsc, \ngp_t$)
using the hierarchical piecewise linear basis.
For $\valueintp[p]_t$,
we use cubic hierarchical weakly fundamental not-a-knot splines
(see \cref{sec:454wfs}).
The not-a-knot boundary conditions help to decrease the
interpolation error (see \cref{sec:541interpolation}),
while the weakly fundamental property eases the hierarchization
complexity by enabling us to use the unidirectional principle
(see \cref{sec:45spatAdaptiveUP,sec:543complexity}).

\paragraph{Extrapolation}

Unfortunately, for many dynamic portfolio choice models,
the state transition is not a function
$\statefcn_t\colon \clint{\*0, \*1} \times \real^{m_{\policy}} \times
\stochdomain \to \clint{\*0, \*1}$,
especially if the state space is actually unbounded.
It may then happen that
$\statefcn_t(\state_t^{(k)}, \policy_t, \stochastic_t^{(j)})
\notin \clint{\*0, \*1}$
for some quadrature nodes
$\stochastic_t^{(j)} \in \stochdomain$ in \cref{eq:bellmanQuadrature}.
Hence, we might not be able to evaluate the value function interpolant
$\valueintp[p]_{t+1}(
  \statefcn_t(\state_t^{(k)}, \policy_t, \stochastic_t^{(j)})
)$, as it is only defined on $\clint{\*0, \*1}$.
Scaling of the domain is not an option due to the dynamic nature of
the problem.

Instead, we extend the interpolant $\valueintp[p]_{t+1}$
to $\real^d$ by an extrapolation method based on Taylor approximation.
First, we crop the evaluation point
$\state_{t+1} \in \real^d \setminus \clint{\*0, \*1}$ to a point
$
  \state_{t+1}^\mathrm{in}
  = \statefcn_t(\state_t^{(k)}, \policy_t, \stochastic_t^{(j)})
  \in \clint{\*0, \*1}
$
with $\state_{t+1}^\mathrm{in} \ceq \vecmin(\vecmax(\state_{t+1}, \*0), \*1)$
(component-wise minimum/maximum).
The extrapolation type,
which may be \texttt{constant}, \texttt{linear}, and \texttt{quadratic},
determines the degree of the Taylor approximation:
\begin{equation}
  \begin{split}
    \valueintp[p]_{t+1}(\state_{t+1})
    &\approx \valueintp[p]_{t+1}(\state_{t+1}^\mathrm{in}) +
    \tr{(\gradient{\state_{t+1}}{\valueintp[p]_{t+1}}(\state_{t+1}^\mathrm{in}))}
    (\state_{t+1} - \state_{t+1}^\mathrm{in})\\
    &\qquad + \tr{(\state_{t+1} - \state_{t+1}^\mathrm{in})}
    (\hessian{\state_{t+1}}{\valueintp[p]_{t+1}}(\state_{t+1}^\mathrm{in}))
    (\state_{t+1} - \state_{t+1}^\mathrm{in}),
  \end{split}
\end{equation}
where \texttt{constant} and \texttt{linear}
only use the first summand and first two summands, respectively.
Since hierarchical B-splines enable us to
exactly and efficiently compute
the gradient $\gradient{\state_{t+1}}{\valueintp[p]_{t+1}}$ and
the Hessian $\hessian{\state_{t+1}}{\valueintp[p]_{t+1}}$,
we do not have to approximate the derivatives with finite differences.

\subsection{Grid Generation}
\label{sec:826gridGeneration}

\paragraph{\texttt{refine} algorithm}

\Cref{alg:financeRefine} shows how to generate the spatially adaptive
sparse grid in \texttt{solveValueFunction}
(\cref{alg:financeSolveValueFunction}).
The underlying criterion is the common surplus-based refinement criterion
\cite{Pflueger13Spatially}.
As for the application in topology optimization (see \cref{chap:60topoOpt}),
we use the piecewise linear interpolant for the surplus-based
grid generation,
since the surpluses are easier to compute in the piecewise linear case,
and they are more meaningful
due to the integral representation formula \eqref{eq:surplusIntegral}.
Parameters for \cref{alg:financeRefine} are
the tolerance $\refinetol_t \in \nonnegreal$,
by which the set of grid points to be refined is determined, and
the number $\norefine_t \in \natz$ of refinement iterations.
These parameters may depend on the time $t$,
since it might be beneficial to change the adaptivity of the
grid over time.

\begin{algorithm}
  \begin{algorithmic}[1]
    \Function{$\valueintp[1]_t = \texttt{refine}$}{%
      $t$, $\valueintp[1]_t$, $\valueintp[p]_{t+1}$%
    }
      \For{$j = 1, \dotsc, \norefine_t$}
        \State{%
          $\ngp_t
          \gets \text{number of grid points of $\valueintp[1]_t$}$%
        }
        \ForOneLine{$k = 1, \dotsc, \ngp_t$}{%
          $\surplus[(k)]{t}
          \gets \text{%
            surplus of $\state_t^{(k)}$ in $\valueintp[1]_t$%
          }$%
        }
        \State{%
          $K_\mathrm{refine}
          \gets \{k = 1, \dotsc, \ngp_t \mid
          \abs{\surplus[(k)]{t}} \ge \refinetol_t\}$%
        }
        \IfOneLine{$K_\mathrm{refine} = \emptyset$}{\Break}
        \State{%
          Refine all grid points in
          $\{\state_t^{(k)} \mid k \in K_\mathrm{refine}\}$%
        }
        \State{%
          $\valueintp[1]_t \gets
          \texttt{optimize($t$, $\valueintp[1]_t$, $\valueintp[p]_{t+1}$)}$%
        }
      \EndFor{}
    \EndFunction{}
  \end{algorithmic}
  \caption[Refinement of the value function (\texttt{refine})]{%
    In-place refinement of the value function $\valueintp[1]_t$.
    Inputs are
    the time $t$,
    the piecewise linear interpolant $\valueintp[1]_t$
    of the current iteration $t$, and
    the higher-order B-spline interpolant $\valueintp[p]_{t+1}$
    of the previous iteration $t + 1$
    (not used if $t = T$).
    The output is the updated piecewise linear interpolant $\valueintp[1]_t$
    with the refined sparse grid.%
  }%
  \label{alg:financeRefine}%
\end{algorithm}

\paragraph{Gradient grids}

The classical surplus-refinement criterion focuses on
regions where the mixed second derivative
$\partialderiv[2d]{
  \partialdiff \stateentry_{t,1}^2 \dotsm \partialdiff \stateentry_{t,d}^2
}{
  \valueintp[1]_t
}$
of $\valueintp[1]_t$ has large absolute values, i.e.,
where $\valueintp[1]_t$ has large high-frequency oscillations.
In gradient-based optimization,
it might be advisable to apply this criterion also
to the partial derivatives
$\partialderiv{\partialdiff \stateentry_{t,o}}{\valueintp[1]_t}$
of $\valueintp[1]_t$ ($o = 1, \dotsc, d$),
since the optimizer depends on the accuracy of the gradient.
In this case, we have to track in \cref{alg:financeSolveValueFunction}
additional sparse grid interpolants for every partial derivative
$\partialderiv{\partialdiff \stateentry_{t,o}}{\valueintp[1]_t}$
that is affected by a policy variable.
This possibility is omitted from the algorithms in this section,
as it would unnecessarily complicate their presentation.

\subsection{Solution for Optimal Policies}
\label{sec:827solvePolicies}

\paragraph{\texttt{solvePolicies} algorithm}

After explaining the generation of the value function interpolants
$\valueintp[p]_t$ ($t = 0, \dotsc, T$),
we move on to step 2 of the general structure of our method
(see \cref{sec:821generalStructure}),
which is the generation of optimal policy interpolants.
The corresponding \cref{alg:financeSolvePolicies} is similar to
\texttt{solveValueFunction} (\cref{alg:financeSolveValueFunction}),
except that it operates on the policy instead of
the value function interpolants.
The functions \texttt{optimize}, \texttt{refine}, and \texttt{interpolate}
have been replaced by corresponding policy versions
\texttt{optimizePolicy}, \texttt{refinePolicy}, and \texttt{interpolatePolicy}
that work very much like their value function counterparts.
\texttt{optimizePolicy} only has to generate new values
if the initial regular sparse grid for the policies
is not contained in the grid of $\valueintp[p]_t$.
The policy grid is then refined and interpolated independently
of the value function grid.
The iterations over time are independent of each other,
which means that they can be parallelized.

\begin{algorithm}
  \begin{algorithmic}[1]
    \Function{%
      $\text{%
        $(\optpolicyintp[p]_t)_{t=0,\dotsc,T}$%
      } = \texttt{solvePolicies}$%
    }{%
      $(\valueintp[p]_t)_{t=0,\dotsc,T}$%
    }
      \State{$\valueintp[p]_{T+1} \gets \emptyset$}
      \Comment{dummy variable (is not used)}%
      \For{$t = 0, \dotsc, T$}
        \State{%
          $\optpolicyintp[1]_t \gets \text{%
            Initial regular sparse grid,
            retrieve values from $\valueintp[p]_t$%
          }$%
        }
        \State{%
          $\optpolicyintp[1]_t \gets
          \texttt{%
            optimizePolicy($t$, $\optpolicyintp[1]_t$, $\valueintp[p]_{t+1}$)%
          }$%
        }
        \State{%
          $\optpolicyintp[1]_t \gets
          \texttt{%
            refinePolicy($t$, $\optpolicyintp[1]_t$, $\valueintp[p]_{t+1}$)%
          }$%
        }
        \State{%
          $\optpolicyintp[p]_t \gets
          \texttt{interpolatePolicy($\optpolicyintp[1]_t$)}$%
        }
      \EndFor{}
    \EndFunction{}
  \end{algorithmic}
  \caption[%
    Generation of interpolants for optimal policies (\texttt{solvePolicies})%
  ]{%
    Generation of interpolants for optimal policies.
    The input is the higher-order B-spline interpolant $\valueintp[p]_t$
    of the value function for all $t = 0, \dotsc, T$.
    The output is the higher-order B-spline interpolant $\optpolicyintp[p]_t$
    of the optimal policies for all $t = 0, \dotsc, T$.%
  }%
  \label{alg:financeSolvePolicies}%
\end{algorithm}

\vspace{-0.5em}
\subsection{Post-Processing}
\label{sec:828postProecessing}
\vspace{-0.5em}

\paragraph{Monte Carlo simulation}

There are various ways to assess
whether the resulting optimal policy B-spline interpolants
$(\optpolicyintp[p]_t)_{t=0,\dotsc,T}$
are reasonable.
One possibility is a Monte Carlo simulation,
where we calculate the mean optimal policy
{%
  \setlength{\abovedisplayskip}{9pt}%
  \setlength{\belowdisplayskip}{9pt}%
  \begin{equation}
    \optpolicymean_t
    \ceq \frac{1}{m_\mathrm{MC}} \sum_{j=1}^{m_\mathrm{MC}} \optpolicyfcn_{t,(j)}
  \end{equation}%
}%
for $m_\mathrm{MC} \in \nat$ individuals.
The optimal policies $\optpolicyfcn_{t,(j)}$ of the individuals
($t = 0, \dotsc, T$ and $j = 1, \dotsc, m_\mathrm{MC}$)
are determined by
\begin{subequations}
  \setlength{\abovedisplayskip}{9pt}%
  \setlength{\belowdisplayskip}{9pt}%
  \begin{align}
    \optpolicyfcn_{t,(j)}
    &\ceq \optpolicyintp[p]_t(\state_{t,(j)}),\\
    \state_{t,(j)}
    &\ceq \statefcn_{t-1}(
      \state_{t-1,(j)}, \optpolicyfcn_{t-1,(j)}, \stochastic_{t-1,(j)}
    ),\quad
    t > 0,\qquad
    \state_{0,(j)}
    \sim P_{0,\state},\\
    \stochastic_{t,(j)}
    &\centerhphantom{\sim}{\hspace*{1.6em}} P_{t,\stochastic},
  \end{align}
\end{subequations}
i.e., the initial state $\state_{0,(j)}$ and the
stochastic variables $\stochastic_{t,(j)}$ are samples of random variables.
Monte Carlo simulations enable us to draw macro-economic conclusions,
e.g., the evolution of the amount of consumption of the average individual
over time.

\section{Transaction Costs Problem}
\label{sec:83problem}

\minitoc[-6mm]{70mm}{4}

\vspace{-1.5em}

\paragraph{Description}

In the \term{transaction costs problem,}
the individual can invest their money risk-free in bonds
(with a fixed interest rate similar to a bank account)
or in $m_{\vstock} \in \nat$ different risk-affected stocks
\cite{Schober18Solving}.
Every stock transaction,
i.e., buy $\buy_{t,j}$ or sell $\sell_{t,j}$,
inflicts transaction costs $\tac \buysell_{t,j}$ ($\tac \in \nonnegreal$)
proportional to the amount $\buysell_{t,j}$ bought or sold
($j = 1, \dotsc, m_{\vstock}$).
The individual only wants to invest a fixed
amount $\wealth_0$ in stocks, i.e., we omit the individual's income.

\subsection{Unnormalized Problem}
\label{sec:831unnormalized}

\paragraph{Consumption and state transition}

In the following,
$\stock_{t,j}$ denotes the fraction of the total wealth $\wealth_t$
that is invested in the $j$-th stock.
We combine these \term{stock fractions} $\stock_{t,j}$
in a vector $\vstock_t \ceq (\stock_{t,1}, \dotsc, \stock_{t,m_{\vstock}})$;
similarly, $\vbuysell_t \ceq (\buysell_{t,1}, \dotsc, \buysell_{t,m_{\vstock}})$
combines buy and sell amounts.
Then, the consumption can be computed as a residual variable
(i.e., a variable that can be fully computed from $\state$ and $\policy$
and is thus omitted from $\policy$),
which is given by
\begin{equation}
  \consume_t
  \ceq (1 - \sumfcn(\vstock_t)) \wealth_t - \bond_t -
  (1 + \tac) \sumfcn(\vbuy_t) + (1 - \tac) \sumfcn(\vsell_t),
\end{equation}
where $\sumfcn(\*a) \ceq \tr{\*1} \*a$
is the sum of all entries of $\*a$.
The state transition is computed by adding the returns of bonds and stocks:%
\begin{equation}
  \wealth_{t+1}
  \ceq \bond_t \bondreturn_t +
  \tr{(\vstock_t \wealth_t + \vbuy_t - \vsell_t)} \vstockreturn_t,\qquad
  \vstock_{t+1}
  \ceq \frac{
    (\vstock_t \wealth_t + \vbuy_t - \vsell_t) \compmult \vstockreturn_t
  }{\wealth_{t+1}},
\end{equation}
where $\bondreturn_t \in \real$ is the bond interest rate,
$
  \vstockreturn_t
  = (\stockreturn_{t,1}, \dotsc, \stockreturn_{t,m_{\vstock}})
  \in \real^{m_{\vstock}}
$
is the vector of (stochastic) stock return rates, and
$\compmult$ is component-wise multiplication.

\subsection{Normalization}
\label{sec:832normalized}

\paragraph{State transition}

The above equations can be normalized with respect to the wealth
$\wealth_t$:
By setting
$\normconsume_t  \ceq \consume_t/\wealth_t$,
$\normbond_t     \ceq \bond_t   /\wealth_t$, and
$\vnormbuysell_t \ceq \vbuysell_t/\wealth_t$, we obtain
\begin{subequations}
  \label{eq:normalizedTCPStateTransition}
  \begin{align}
    \normconsume_t
    &\mathrel{\righthphantom{=}{\ceq}}
    (1 - \sumfcn(\vstock_t)) - \normbond_t -
    (1 + \tac) \sumfcn(\vnormbuy_t) + (1 - \tac) \sumfcn(\vnormsell_t),\\
    \wealthratio_{t+1}
    &\ceq \normbond_t \bondreturn_t +
    \tr{(\vstock_t + \vnormbuy_t - \vnormsell_t)} \vstockreturn_t,
    \qquad(= \wealth_{t+1}/\wealth_t)\\
    \vstock_{t+1}
    &\mathrel{\righthphantom{=}{\ceq}}
    \frac{
      (\vstock_t + \vnormbuy_t - \vnormsell_t) \compmult \vstockreturn_t
    }{\wealthratio_{t+1}},
  \end{align}
\end{subequations}
where $\normconsume_t$ and $\wealthratio_{t+1}$ are residual
variables that specify \term{normalized consumption} and
\term{wealth ratio,} respectively.
All in all, the resulting dynamic portfolio choice model has
the following variables:
\begin{itemize}
  \item
  $\centerhphantom{d}{m_{\stochastic}} = m_{\vstock}$
  state variables $\normstate_t$:
  Stock fractions $\stock_{t,1}, \dotsc, \stock_{t,m_{\vstock}}$
  
  \item
  $\centerhphantom{m_{\policy}}{m_{\stochastic}} = 2m_{\vstock} + 1$
  policy variables $\normpolicy_t$:
  Normalized bonds $\normbond_t$,
  normalized buy amounts $\normbuy_{t,1}, \dotsc, \normbuy_{t,m_{\vstock}}$ and
  normalized sell amounts $\normsell_{t,1}, \dotsc, \normsell_{t,m_{\vstock}}$
  
  \item
  $m_{\stochastic} = m_{\vstock}$
  stochastic variables $\stochastic_t$:
  Stock return rates $\stockreturn_{t,1}, \dotsc, \stockreturn_{t,m_{\vstock}}$
\end{itemize}
The state space and policy space constraints are given by
\begin{subequations}
  \label{eq:normalizedTCPConstraints}
  \newcommand*{\centereqline}[1]{%
    \mathclap{\hphantom{\mathrm{(8.99a)}}#1}%
  }%
  \begin{gather}
    \label{eq:normalizedTCPConstraintsShort}
    \centereqline{
      \vstock_t \ge \*0,\qquad
      \sumfcn(\vstock_t) \le 1,\qquad
      \normbond_t \ge 0,\qquad
      \vnormbuysell_t \ge \*0,\qquad
      \vnormsell_t \le \vstock_t,\qquad
      \wealthratio_{t+1} \ge 0,
    }\\
    \label{eq:normalizedTCPConstraintsLong}
    \centereqline{
      \normconsume_{\min} + \normbond_t +
      (1 + \tac) \sumfcn(\vnormbuy_t) - (1 - \tac) \sumfcn(\vnormsell_t)
      \le 1 - \sumfcn(\vstock_t),
    }
  \end{gather}
\end{subequations}
where $\normconsume_{\min} \in \nonnegreal$ is some minimal consumption
that must be maintained.

\paragraph{Bellman equation}

Consequently, the Bellman equation \eqref{eq:gridBellmanCET}
after the certainty-equiva\-lent transformation has to be
normalized as well.
By setting $\normcetvalueintp_t(\state_t^{(k)})
\ceq \cetvalueintp_t(\state_t^{(k)})/\wealth_t$, we obtain
\begin{subequations}
  \begin{align}
    &\hphantom{=}\hspace{0.6em} \normcetvalueintp_t(\state_t^{(k)})
    = \wealth_t^{-1} \cetvalueintp_t(\state_t^{(k)})\\
    &= \max_{\policy_t} \left(
      \left(
        \left(
          \wealth_t^{-1} \consume_t(\state_t^{(k)}, \policy_t)
        \right)^{1-\riskav} +
        \patience \expectation[t]{
          \left(
            \wealth_t^{-1} \cetvalueintp_{t+1}(
              \statefcn_t(\state_t^{(k)}, \policy_t, \stochastic_t)
            )
          \right)^{1-\riskav}
        }
      \right)^{1/(1-\riskav)}
    \right)\\
    \label{eq:normalizedTCPBellmanEquation}
    &= \max_{\normpolicy_t} \left(
      \left(
        \normconsume_t(\state_t^{(k)}, \normpolicy_t)^{1-\riskav} +
        \patience \expectation[t]{
          \bigl(
            \wealthratio_{t+1} \normcetvalueintp_{t+1}(
              \normstatefcn_t(\state_t^{(k)}, \normpolicy_t, \stochastic_t)
            )
          \bigr)^{1-\riskav}
        }
      \right)^{1/(1-\riskav)}
    \right).
  \end{align}
\end{subequations}
This means that compared with
\eqref{eq:gridBellmanCET},
the value function in the expectation has to be multiplied by
the wealth ratio $\wealthratio_{t+1}$ introduced above in
\eqref{eq:normalizedTCPStateTransition}.
Since there is no inheritance, the optimal terminal solution
is to sell all stocks and consume everything:
\begin{equation}
  \normcetvalueintp_t(\state_T^{(k)})
  = 1 - \tac \sumfcn(\vstock_T^{(k)}),\quad
  \normbond_T^{\opt}(\state_T^{(k)})
  = 0,\quad
  \vnormbuy[\opt]_T(\state_T^{(k)})
  = \*0,\quad
  \vnormsell[\opt]_T(\state_T^{(k)})
  = \vstock_T^{(k)}.
\end{equation}

\subsection{State Space Cropping}
\label{sec:833cropping}

\paragraph{Sparse grids on non-rectangular domains}

Unfortunately, the constraint $\sumfcn(\vstock_t) \le 1$
from \cref{eq:normalizedTCPConstraints} limits the feasible state space
region to a proper subset (which is the unit simplex)
of the unit hypercube $\clint{\*0, \*1}$,
which impedes the direct application of sparse grids.
There are three possible remedies:
transforming the unit hypercube to the feasible state space,
applying extrapolation techniques as discussed in
\cref{sec:825interpolation}, or
choosing a model-tailored approach to obtain
function values outside the feasible state space.

\paragraph{Virtual selling of stocks}

We choose the third remedy and \term{virtually sell,}
if $\sumfcn(\vstock_t) > 1$,
as many stocks as needed to meet the constraint $\sumfcn(\vstock_t) \le 1$.
We already might need to sell stocks
even if $\sumfcn(\vstock_t)$ is smaller but close to one
in order to satisfy the minimum consumption requirement
\eqref{eq:normalizedTCPConstraintsLong}.
In detail, we replace $\vstock_t$ by $\normcropfactor \vstock_t$
whenever $\normcropfactor < 1$,
where $\normcropfactor \in \posreal$ is a \term{cropping factor}
that is determined by
\begin{equation}
  \label{eq:virtualSelling}
  \Bigl[
    1 - \tac\, \bigl(
      \sumfcn(\vstock_t) - \sumfcn(\normcropfactor\vstock_t)
    \bigr)
  \Bigr]
  \cdot \bigl(1 - \sumfcn(\normcropfactor\vstock_t)\bigr)
  = \normconsume_{\min}.
\end{equation}
Here, $\bigl(\sumfcn(\vstock_t) - \sumfcn(\normcropfactor\vstock_t)\bigr)$
is the amount of virtually sold stocks.
Hence, the term in square brackets is the fraction of wealth
that is still available after deducting the induced transaction costs.
The product of this term with
$\bigl(1 - \sumfcn(\normcropfactor\vstock_t)\bigr)$
is the fraction of wealth that can be consumed after the virtual selling,
which needs to be at least $\normconsume_{\min}$.
Solving \cref{eq:virtualSelling} for $\normcropfactor$ and
choosing the positive solution, we finally obtain
\begin{equation}
  \newcommand*{\sumX}{\sumfcn(\vstock_t)}
  \newcommand*{\cMin}{\normconsume_{\min}}
  \normcropfactor
  \ceq \frac{
    \tac\, \bigl(1 + \sumX\bigr) - 1 +
    \sqrt{
      \tac^2\, \bigl(1 - \sumX\bigr)^2
      - 2 \tac\, \bigl(2 \cMin - 1 + \sumX\bigr) + 1
    }
  }{
    2 \tac \sumX
  }.
\end{equation}

\subsection{Euler Equation Errors}
\label{sec:834eulerErrors}

\paragraph{Motivation}

Due to the curse of dimensionality,
reasonably accurate full grid reference solutions
of the transaction costs problem can only be computed
if the number $m_{\vstock}$ of stocks is small.
Mainly (but not only) in higher-dimensional settings,
a different means of assessing the
quality of sparse grid solutions is desirable.
We use Euler equation errors to measure the deviation in
the first-order optimality conditions.

\paragraph{Derivation}

In the following, we fix the state $\normstate_t \in \clint{\*0, \*1}$
for which we want to compute the Euler equation error.
We abbreviate
the
value function interpolant
$
  \normcetvalueintp_t
  \ceq \normcetvalueintp_t(\normstate_t)
$,
the state transition function
$
  \normstatefcn_t
  \ceq \normstatefcn_t(\normstate_t, \normpolicy_t, \stochastic_t)
$,
the wealth ratio
$
  \wealthratio_{t+1}
  \ceq \wealthratio_{t+1}(
    \normstate_t, \normpolicy_t, \stochastic_t
  )
$, and
the consumption
$
  \normconsume_t
  \ceq \normconsume_t(\normstate_t, \normpolicy_t)
$.
The Lagrangian of the optimization problem corresponding
to the Bellman equation \eqref{eq:normalizedTCPBellmanEquation}
of the normalized transaction costs problem
with respect to the problem's constraints
\eqref{eq:normalizedTCPConstraints} is given by
{%
  \setlength{\abovedisplayskip}{9pt}%
  \setlength{\belowdisplayskip}{9pt}%
  \begin{equation}
    \begin{split}
      \lagrangian_t(\normstate_t, \normpolicy_t, \*\multiplier)
      &\ceq \left(
        (\normconsume_t)^{1-\riskav} +
        \patience \expectation[t]{
          \bigl(
            \wealthratio_{t+1}\;
            \normcetvalueintp_{t+1}(\normstatefcn_t)
          \bigr)^{1-\riskav}
        }
      \right)^{1/(1-\riskav)}\\
      &\hspace*{6mm} {}
      - \multiplier_1 \normbond_t
      - \tr{\*\multiplier_2} \vnormbuy_t
      - \tr{\*\multiplier_3} \vnormsell_t
      + \tr{\*\multiplier_4}\; (\vnormsell_t - \vstock_t)
      + \multiplier_5\; (\normconsume_{\min} - \normconsume_t)
    \end{split}
  \end{equation}%
}%
with $
\*\multiplier \ceq (
  \multiplier_1,
  \*\multiplier_2,
  \*\multiplier_3,
  \*\multiplier_4,
  \multiplier_5
)$,
$\multiplier_1, \multiplier_5 \in \real$, and
$\*\multiplier_2, \*\multiplier_3, \*\multiplier_4 \in \real^{m_{\vstock}}$.
According to the first-order conditions
\term{(Karush--Kuhn--Tucker (KKT) conditions),}
the partial derivative
$
  \partialderiv{\partialdiff \normbond_t}{\lagrangian_t}(
    \normstate_t, \normpolicy_t, \*\multiplier
  )
$
with respect to $\normbond_t$
vanishes in the exact optimum
$\normpolicy_t = \optnormpolicyfcn_t \ceq \optnormpolicyfcn_t(\normstate_t)$,
i.e.,
{%
  \setlength{\abovedisplayskip}{9pt}%
  \setlength{\belowdisplayskip}{9pt}%
  \begin{equation}
    \label{eq:eulerErrorFirstOrderCondition}
    \partialderiv{\partialdiff \normbond_t}{
      \left(
        (\normconsume_t^{\opt})^{1-\riskav} +
        \patience \expectation[t]{
          \bigl(
            \wealthratio_{t+1}^{\opt}\;
            \normcetvalueintp_{t+1}(\normstatefcn[\opt]_t)
          \bigr)^{1-\riskav}
        }
      \right)^{1/(1-\riskav)}
    }
    - \multiplier_1
    - \multiplier_5
    \partialderiv{\partialdiff \normbond_t}{\normconsume_t^{\opt}}
    = 0,
  \end{equation}%
}%
where
$
  \normstatefcn[\opt]_t
  \ceq \normstatefcn_t(\normstate_t, \optnormpolicyfcn_t, \stochastic_t)
$,
$
  \wealthratio_{t+1}^{\opt}
  \ceq \wealthratio_{t+1}(
    \normstate_t, \optnormpolicyfcn_t, \stochastic_t
  )
$, and
$
  \normconsume_t^{\opt}
  \ceq \normconsume_t(\normstate_t, \optnormpolicyfcn_t)
$.
We now neglect binding constraints,
i.e., we assume that $\multiplier_1 = \multiplier_5 = 0$,
otherwise we cannot compute the error.
After calculating the derivatives,
\cref{eq:eulerErrorFirstOrderCondition} becomes
{%
  \setlength{\abovedisplayskip}{9pt}%
  \setlength{\belowdisplayskip}{9pt}%
  \begin{equation}
    \patience \bondreturn_t
    \cdot \expectationsign[t]\Bigl[
      \bigl(
        \normcetvalueintp_t
        - \tr{(\gradient{\normstate_t}{\normcetvalueintp_t})}
        \normstatefcn[\opt]_t
      \bigr) \cdot
      \bigl(
        \wealthratio_{t+1}^{\opt}\; \normcetvalueintp_t
      \bigr)^{-\riskav}
    \Bigr]
    = (\normconsume_t^{\opt})^{-\riskav}.
  \end{equation}%
}%
This equation can be used as an error measure by substituting
$\optnormpolicyfcn_t$ for the interpolated optimum
$\optnormpolicyintp_t = \optnormpolicyintp_t(\normstate_t)$.
By multiplying the resulting equation by
$
  (\normconsume_t^{\opt,\sparse})^{\riskav}
  \ceq (\normconsume_t(\normstate_t, \optnormpolicyintp_t))^{\riskav}
$, we obtain the
\term{unit-free Euler equation errors $\eulererror_t(\normstate_t)$}
with respect to $\normbond_t$:
\begin{equation}
  \eulererror_t(\normstate_t)
  \ceq \Bigl|
    1 - \Bigl(
      \patience \bondreturn_t (\normconsume_t^{\opt,\sparse})^{\riskav}
      \cdot \expectationsign[t]\Bigl[
        \bigl(
          \normcetvalueintp_t
          - \tr{(\gradient{\normstate_t}{\normcetvalueintp_t})}
          \normstatefcn[\opt,\sparse]_t
        \bigr) \cdot
        \bigl(
          \wealthratio_{t+1}^{\opt,\sparse}\;
          \normcetvalueintp_t
        \bigr)^{-\riskav}
      \Bigr]
    \Bigr)^{-1/\riskav}
  \Bigr|
\end{equation}
with
$
  \normstatefcn[\opt,\sparse]_t
  \ceq \normstatefcn_t(\normstate_t, \optnormpolicyintp_t, \stochastic_t)
$ and
$
  \wealthratio_{t+1}^{\opt,\sparse}
  \ceq \wealthratio_{t+1}(
    \normstate_t, \optnormpolicyintp_t, \stochastic_t
  )
$.

\paragraph{Weighted Euler equation errors}

However, the state space cropping as introduced above
distorts Euler equation errors:
The error $\eulererror_t(\normstate_t)$ does not vanish
even for the exact solution and even inside the feasible state space.
This is because the cropping already occurs for large stock holdings
$\sumfcn(\normstate_t)$ that are less than one,
as stocks have to be sold to maintain
minimum consumption $\normconsume_{\min}$.
Numerical experiments show that due to this issue,
the error attains large values in the region near the hyperplane
$\sumfcn(\normstate_t) = 1$.
Economically, this region is not significant
as such large stock fractions are highly unusual,
which is confirmed by Monte Carlo simulations.
We therefore use the \term{weighted Euler equation error}
\begin{equation}
  \weightedeulererror_t(\normstate_t)
  \ceq \bigl(1 - \sumfcn(\normstate_t)\bigr) \cdot
  \eulererror_t(\normstate_t)
\end{equation}
instead of $\eulererror_t$,
although other strategies exist
such as restricting the state domain where the error is computed or
weighting the error with the probability that a given state
occurs in Monte Carlo simulations.

\section{Implementation and Numerical Results}
\label{sec:84results}

\minitoc[3mm]{65mm}{6}

\parbox{1em}{}
\vspace{-3em}

\disableornamentsfornextheadingtrue
\subsection{Implementation}
\label{sec:841implementation}

\paragraph{Parameter values}

We used
a risk aversion factor of $\riskav \ceq 3.5$,
a patience factor of $\patience \ceq 0.97$,
a transaction cost rate of $\tac \ceq \SI{1}{\percent}$, and
a minimum consumption of $\normconsume_{\min} \ceq 0.001$.
The bond and stock return rates $\bondreturn_t$ and $\vstockreturn_t$
were taken from \cite{Cai10Stable};
the log-normally distributed stock return rates were generalized
from the three-stock case to five stocks via
$\ln \vstockreturn_t \sim \mathcal{N}(\*\mu, \mat{\Sigma})$, where
{%
  \setlength{\abovedisplayskip}{6pt}%
  \setlength{\belowdisplayskip}{6pt}%
  \begin{equation}
    \label{eq:financeStockReturnMeanCovariance}
    \*\mu
    \ceq
    \scalebox{0.92}{$
      \begin{pmatrix*}[l]
        0.0572\\
        0.0638\\
        0.07\\
        0.0764\\
        0.0828
      \end{pmatrix*}
    $},\quad
    \mat{\Sigma}
    \ceq 10^{-2}
    \scalebox{0.92}{$
      \begin{pmatrix*}[l]
        2.56&  0.576&   0.288&   0.176&  0.096\\
        0.576& 3.24&    0.90432& 1.0692& 1.296\\
        0.288& 0.90432& 4&       1.32&   1.68\\
        0.176& 1.0692&  1.32&    4.84&   2.112\\
        0.096& 1.296&   1.68&    2.112&  5.76
      \end{pmatrix*}
    $}.
  \end{equation}%
}%
The models were solved for $T = 6$ time steps;
this number suffices to show all relevant numerical effects and results,
while keeping the computational effort at a reasonable level.
As initial grids, we employed regular sparse grids
$\coarseregsgset{n}{d}{b}$ with $b = 1$
to decrease the number of grid points
(see \cref{sec:241coarseBoundary}).

\vspace*{0.25em}

\paragraph{Software}

The dynamic portfolio choice models were solved using a self-written
MATLAB framework.
The object-oriented framework was designed in such a way that
not only transaction costs problems,
but many other types of dynamic portfolio choice models can be handled.
For instance, the base class \texttt{LifecycleProblem} provides
an interface with abstract functions such as
\texttt{computeTerminalValueFunction} and
\texttt{computeStateTransition}.
The actual functionality implemented in the base class strongly resembles
the algorithms presented in \cref{sec:82algorithms}.
This is not only desirable from a modeling perspective,
but also facilitates future usage by other researchers.
For creating (i.e., hierarchizing) and evaluating sparse grid interpolants,
the sparse grid toolbox \sgpp was used \cite{Pflueger10Spatially}.%
\footnote{%
  \url{http://sgpp.sparsegrids.org/}%
}
The emerging optimization problems were solved using
sequential quadratic programming methods supplied by the
NAG Toolbox for MATLAB.%
\footnote{%
  \url{https://www.nag.com/}%
}
To avoid being stuck in local minima,
we repeated the optimization process for a varying number
of initial multi-start points (in the range of a few dozens).
All computation times were measured on a shared-memory computer
with 4x Intel Xeon E7-8880v3 (72 cores, 144 threads).

\subsection{Error Sources and Error Measure}
\label{sec:842errorSources}

\vspace*{0.25em}

\paragraph{Error sources}

In this application, there are the following error sources:

\begin{enumerate}[
  label=E\arabic*.,
  ref=E\arabic*,
  leftmargin=2.7em,
]
  \item
  \label{item:financeErrorInterpolationValue}
  Interpolation of the value function
  (i.e., $\normcetvalueintp_{t+1} \not= \normcetvaluefcn_{t+1}$)
  
  \item
  \label{item:financeErrorInterpolationPolicy}
  Interpolation of the policy functions
  (i.e., $\optnormpolicyintp_t \not= \optnormpolicyfcn_t$)
  
  \item
  \label{item:financeErrorExtrapolation}
  Extrapolation
  (i.e., $
  \normcetvalueintp_{t+1}(\state_{t+1})
  \not= \normcetvaluefcn_{t+1}(\state_{t+1})
  $)
  
  \item
  \label{item:financeErrorCropping}
  State space cropping
  (i.e., Euler errors do not vanish for exact solution)
  
  \item
  \label{item:financeErrorOptimization}
  Optimization
  (i.e., the minimum found by the optimizer is inaccurate or not global)
  
  \item
  \label{item:financeErrorQuadrature}
  Quadrature
  ($
    \expectation[t]{\cdots}
    \not= \sum_{j=1}^{m_{\quadweight}} \quadweight_t^{(j)}
    \cdot [\cdots](\stochastic_t^{(j)})
  $)
  
  \item
  \label{item:financeErrorRounding}
  Floating-point rounding errors
  (i.e., arithmetical operations are inaccurate)
\end{enumerate}
Due to the dynamic programming scheme,
the combination of all errors accumulates over $t$.
For instance, if the optimization does not find the global optimum
exactly or it only finds a local one for $t + 1$,
the error propagates from the interpolant $\normcetvalueintp_{t+1}$
on the right-hand side of the Bellman equation
\eqref{eq:normalizedTCPBellmanEquation} to $\normcetvalueintp_t$
on the left-hand side, and so on.
If the system does not damp these errors,
the error steadily becomes larger backwards in time $t$.

\paragraph{Error measure}

We use the weighted Euler equation error
$\weightedeulererror_t(\normstate_t)$
to assess the quality of the resulting policies
($\Ltwo$ norm or pointwise).
As the errors generally grow backwards in time,
it suffices to consider $t = 0$.
However, since Euler equation errors can only be evaluated at points in
the simplex
$
  \Omega_\mathrm{simplex}
  \ceq \{\normstate_t \in \clint{\*0, \*1} \mid \sumfcn(\normstate_t) \le 1\}
$, the $\Ltwo$ norm would
quickly converge to zero with growing dimensionality, even if the mean error
stayed constant.
Therefore, we normalize the $\Ltwo$ norm:
\begin{equation}
  \weightedeulererrorLtwo_t
  \ceq \sqrt{d!} \cdot \normLtwo{\weightedeulererror_t}
  = \sqrt{
    \frac{1}{\vol{\Omega_\mathrm{simplex}}}
    \int_{\mathrlap{\Omega_\mathrm{simplex}}\hphantom{\Omega}}
    \weightedeulererror_t(\normstate_t)^2 \diff{}\normstate_t
  },
\end{equation}
where the expression under the root sign is approximated
via Monte Carlo quadrature as the mean of samples
of $\weightedeulererror_t(\normstate_t)^2$.

\subsection{Numerical Results}
\label{sec:843results}

\paragraph{Full grid solution}

\begin{figure}
  \makebox[49mm][r]{%
    \includegraphics{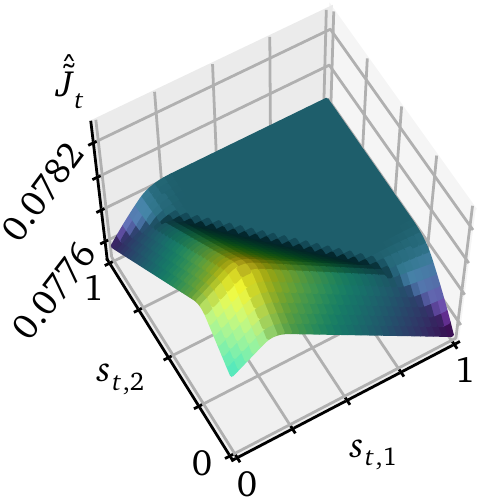}%
  }%
  \hfill%
  \makebox[49mm][r]{%
    \includegraphics{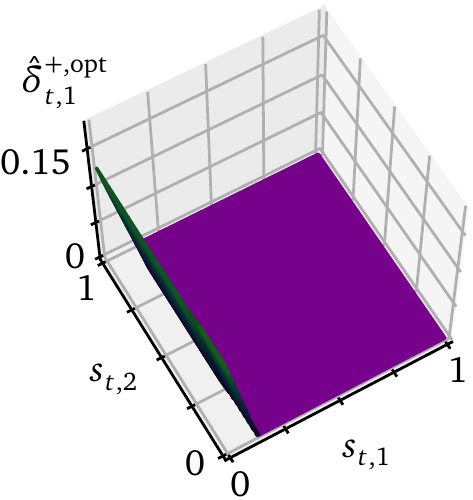}%
  }%
  \hfill%
  \makebox[49mm][r]{%
    \includegraphics{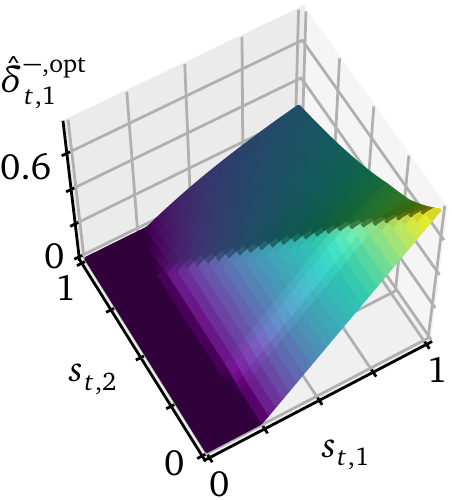}%
  }%
  \\[0mm]%
  \makebox[49mm][r]{%
    \includegraphics{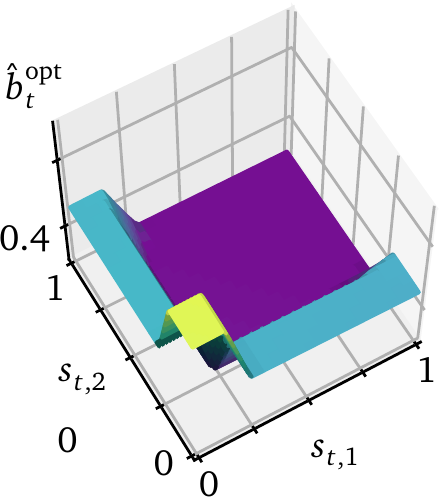}%
  }%
  \hfill%
  \makebox[49mm][r]{%
    \includegraphics{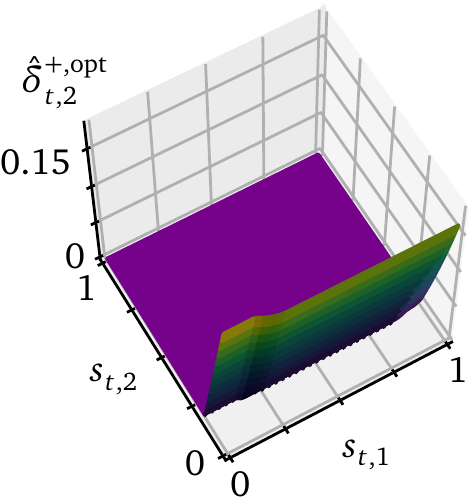}%
  }%
  \hfill%
  \makebox[49mm][r]{%
    \includegraphics{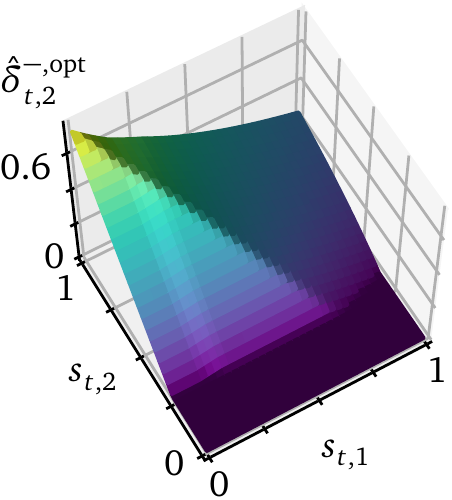}%
  }%
  \caption[Reference solution for the two-dimensional TCP]{%
    Full grid solution for the transaction costs problem
    with $d = 2$ stocks.
    Shown are the value function $\normcetvaluefcn_t$ \emph{(top left)} and the
    optimal policy $\optnormpolicyfcn_t$ for $t = 0$.%
  }%
  \label{fig:financeSolution2DReference}%
\end{figure}

We show in \cref{fig:financeSolution2DReference}
a full grid solution for the case of $d = 2$ stocks,
i.e., $\{\state_t^{(k)} \mid k = 1, \dotsc, \ngp_t\} = \fgset{n,d}$
for some fixed level $n \in \nat$
(here, $n = 7$ and $\ngp_t = (2^7 + 1)^2 = \num{16641}$) and
for all $t = 0, \dotsc, T$.
Obviously, this is only computationally feasible
for low dimensionalities $d$ due to the curse of dimensionality.
The two-dimensional solution of level $n = 7$
took over nine hours to compute.
The solution of the next level is estimated to already take one week.
Hence, full grid solutions can only be computed up to $d = 3$
due to excessive computation time for $d \ge 4$.
This underlines the need for sophisticated
discretization techniques such as sparse grids.

\paragraph{Convergence of the weighted Euler equation error}

\begin{figure}
  \includegraphics{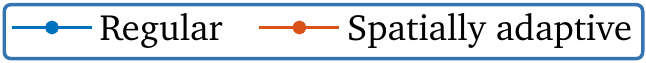}%
  \\[2mm]%
  \subcaptionbox{%
    $d = 1$%
  }[37mm]{%
    \includegraphics{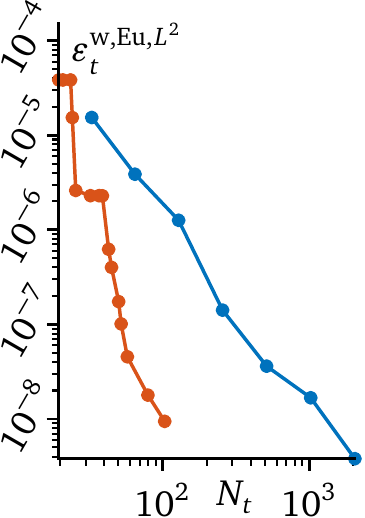}%
  }%
  \hfill%
  \subcaptionbox{%
    $d = 2$%
  }[37mm]{%
    \includegraphics{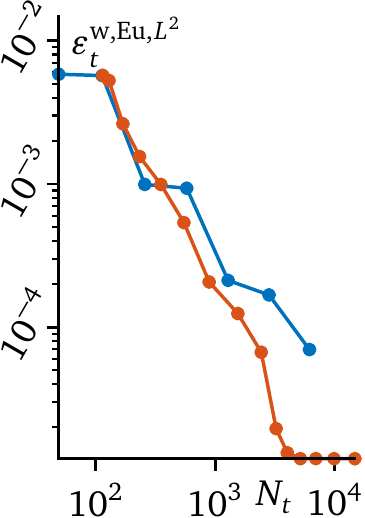}%
  }%
  \hfill%
  \subcaptionbox{%
    $d = 3$%
  }[37mm]{%
    \includegraphics{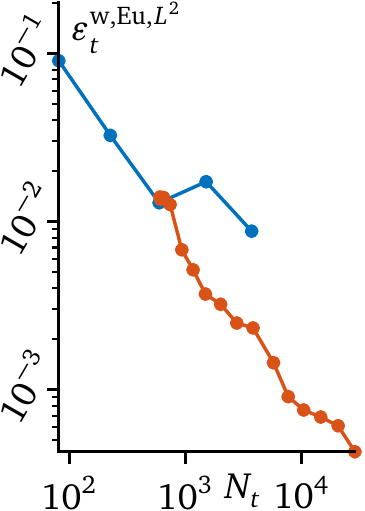}%
  }%
  \hfill%
  \subcaptionbox{%
    $d = 4$%
  }[37mm]{%
    \includegraphics{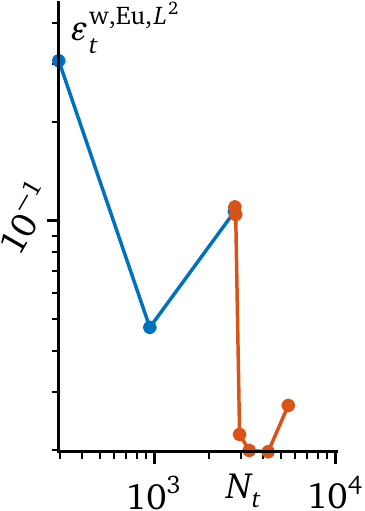}%
  }%
  \caption[Convergence of the weighted Euler equation error]{%
    Convergence of the $\Ltwo$ norm $\weightedeulererrorLtwo_t$
    of the weighted Euler equation error for $t = 0$ for
    regular sparse grids \emph{\textcolor{C0}{(blue)}} and
    spatially adaptive sparse grids \emph{\textcolor{C1}{(red)}.}
    The number $\ngp_t$ is the average number
    $\frac{1}{m_{\policy}} \sum_{j=1}^{m_{\policy}} \ngp_{t,j}$
    of grid points over all policy grids for $t = 0$,
    where $\ngp_{t,j}$ is the number of grid points
    of the $j$-th policy entry.%
  }%
  \label{fig:financeEulerError}%
\end{figure}

\Cref{fig:financeEulerError} shows the convergence of the
$\Ltwo$ norm $\weightedeulererrorLtwo_0$ of the
weighted Euler equation error for $t = 0$ for regular sparse grids
and spatially adaptive sparse grids
for the cases of $d = 1, \dotsc, 4$ stocks.
For this and the following plots,
the value function grid is left unchanged
(usually a slightly refined regular sparse grid),
while the average number $\ngp_t$ of policy grid points increases
with decreasing refinement threshold $\refinetol_t$.
This is because the value function grid does not seem to have a great influence
on the convergence of the Euler equation errors.
Compared to regular grids, the spatial adaptivity decreases the error by
two orders of magnitude in one dimension.
The gain is smaller for higher dimensionalities $d$,
but spatially adaptive grids still outperform regular grids.
For $d = 2$, we observe that the error saturates
at $\ngp_t \approx \num{4000}$ points just above $10^{-5}$.
This is most likely due to the parts
\ref{item:financeErrorExtrapolation} to
\ref{item:financeErrorRounding} of the error that are not influenced
by sparse grid interpolation.
In addition, convergence significantly decelerates starting with $d = 4$.
For $d = 4$, spatially adaptive sparse grids
are able to achieve a weighted Euler equation error of
$\weightedeulererrorLtwo_t \approx \num{2.0e-2}$ for $t = 0$
(with an average number $\ngp_0 = \num{4252}$ of policy grid points).
For $d = 5$, we are still able to achieve a small error of
$\weightedeulererrorLtwo_t \approx \num{1.9e-2}$ for $t = 0$
with spatially adaptive sparse grids with
an average number $\ngp_0 = \num{12572}$ of policy grid points.
While we cannot detect any convergence for this dimensionality yet,
this is still a major result as such high-dimensional models
could not be solved up to now with conventional methods.

\paragraph{Optimal policies in 2D and 5D}

\begin{figure}
  \subcaptionbox{%
    $\normcetvalueintp[1]_t$%
  }[48mm]{%
    \includegraphics{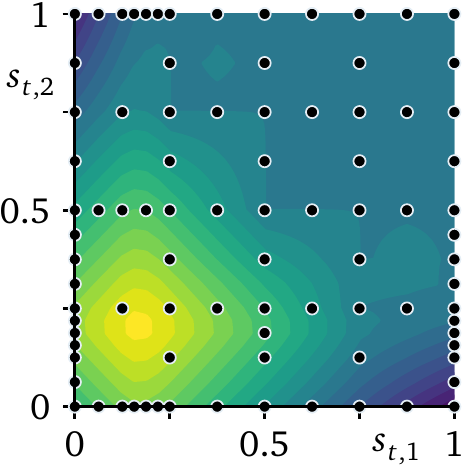}%
  }%
  \hfill%
  \subcaptionbox{%
    $\normbuy[\opt,\sparse,1]_{t,1}$%
  }[48mm]{%
    \includegraphics{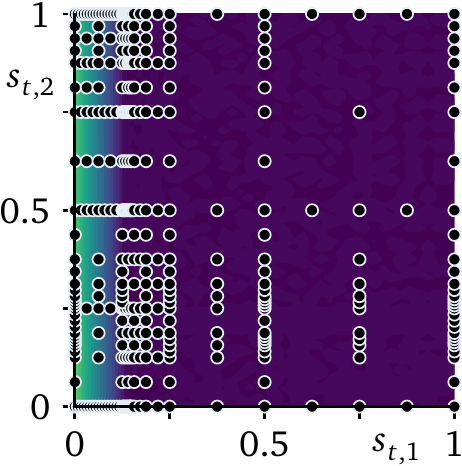}%
  }%
  \hfill%
  \subcaptionbox{%
    $\normsell[\opt,\sparse,1]_{t,1}$%
  }[48mm]{%
    \includegraphics{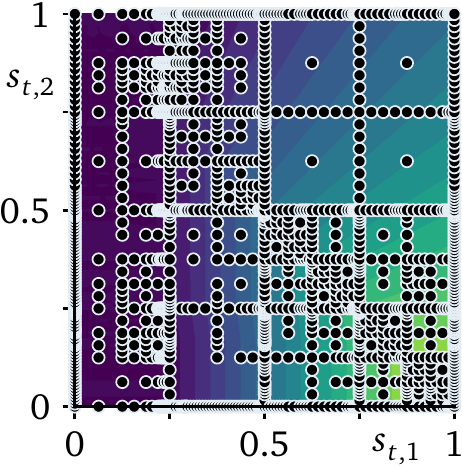}%
  }%
  \\[2mm]%
  \subcaptionbox{%
    $\normbond_t^{\opt,\sparse,1}$%
  }[48mm]{%
    \includegraphics{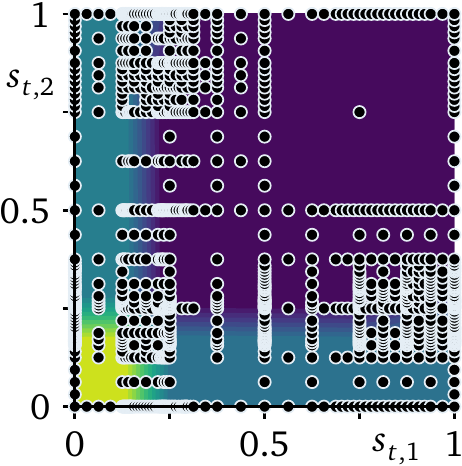}%
  }%
  \hfill%
  \subcaptionbox{%
    $\normbuy[\opt,\sparse,1]_{t,2}$%
  }[48mm]{%
    \includegraphics{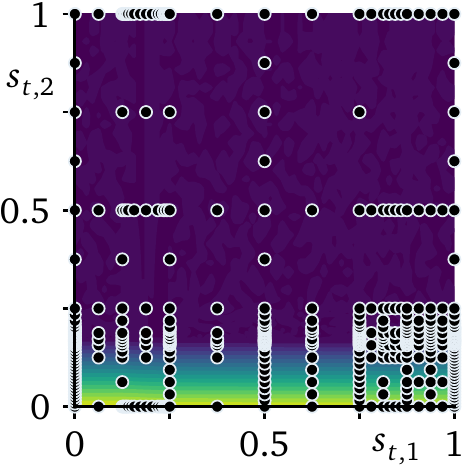}%
  }%
  \hfill%
  \subcaptionbox{%
    $\normsell[\opt,\sparse,1]_{t,2}$%
  }[48mm]{%
    \includegraphics{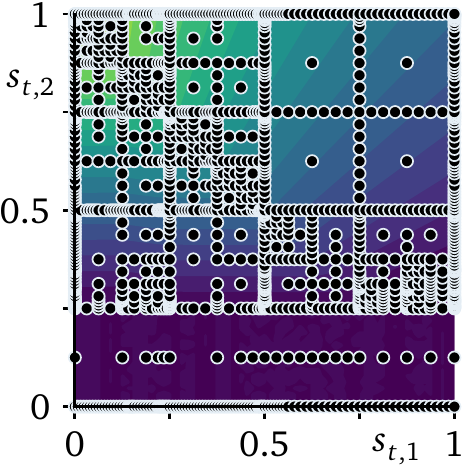}%
  }%
  \caption[Sparse grid solution for the two-dimensional TCP]{%
    Spatially adaptive sparse grid solution for the transaction costs problem
    \vspace{-0.15em}%
    with $d = 2$ stocks.
    \vspace{-0.15em}%
    Shown are the value function $\normcetvalueintp_t$ \emph{(top left)} and the
    optimal policy $\optnormpolicyintp_t$ for the initial time step $t = 0$,
    together with the corresponding grid points \emph{(dots).}
    The color coding is the same as in
    \cref{fig:financeSolution2DReference}.%
  }%
  \label{fig:financeSolution2DSparseGrid}%
\end{figure}

\begin{figure}
  \makebox[37mm][c]{%
    \hspace*{3.8mm}%
    \raisebox{-\height}{\includegraphics{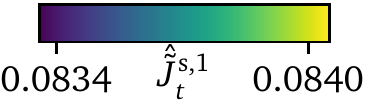}}%
  }%
  \hfill%
  \makebox[37mm][c]{%
    \hspace*{3.8mm}%
    \raisebox{-\height}{\includegraphics{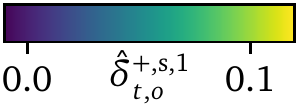}}%
  }%
  \hfill%
  \makebox[37mm][c]{%
    \hspace*{4.4mm}%
    \raisebox{-\height}{\includegraphics{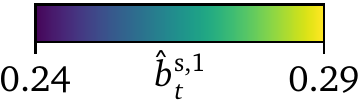}}%
  }%
  \hfill%
  \makebox[37mm][c]{%
    \hspace*{2.4mm}%
    \raisebox{-\height}{\includegraphics{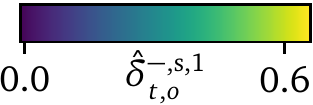}}%
  }%
  \\[1mm]%
  \makebox[37mm][c]{%
    \includegraphics{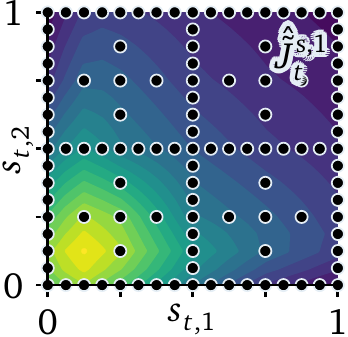}%
  }%
  \hfill%
  \makebox[37mm][c]{%
    \includegraphics{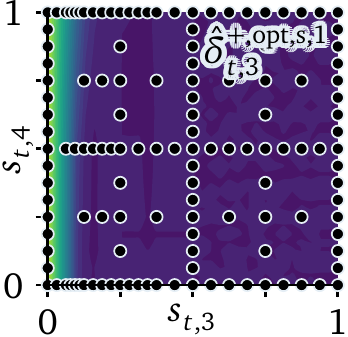}%
  }%
  \hfill%
  \makebox[37mm][c]{%
    \includegraphics{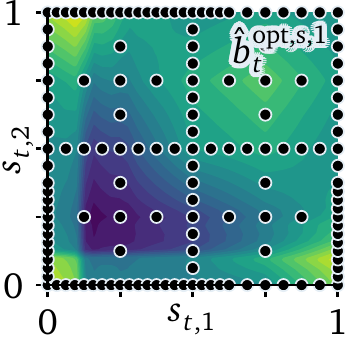}%
  }%
  \hfill%
  \makebox[37mm][c]{%
    \includegraphics{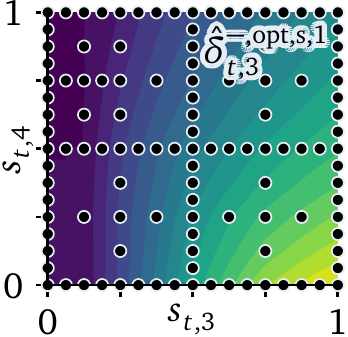}%
  }%
  \\[1mm]%
  \makebox[37mm][c]{%
    \includegraphics{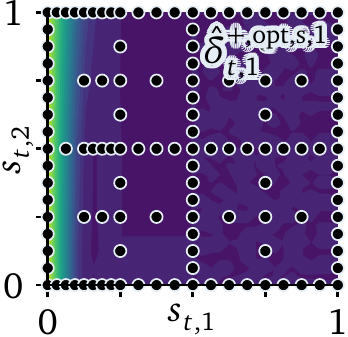}%
  }%
  \hfill%
  \makebox[37mm][c]{%
    \includegraphics{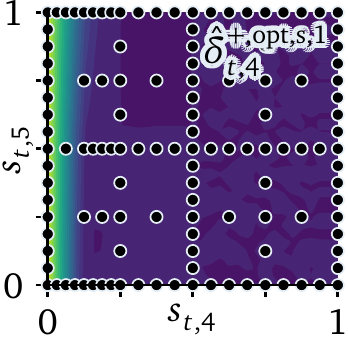}%
  }%
  \hfill%
  \makebox[37mm][c]{%
    \includegraphics{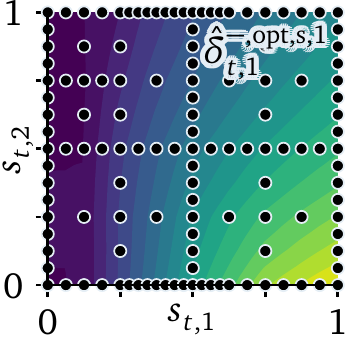}%
  }%
  \hfill%
  \makebox[37mm][c]{%
    \includegraphics{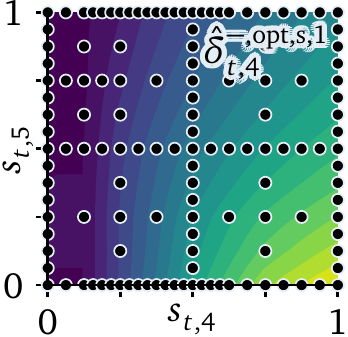}%
  }%
  \\[1mm]%
  \makebox[37mm][c]{%
    \includegraphics{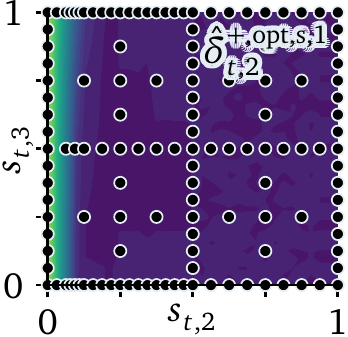}%
  }%
  \hfill%
  \makebox[37mm][c]{%
    \includegraphics{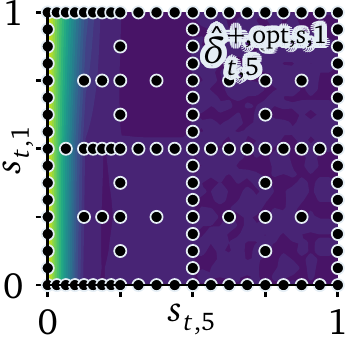}%
  }%
  \hfill%
  \makebox[37mm][c]{%
    \includegraphics{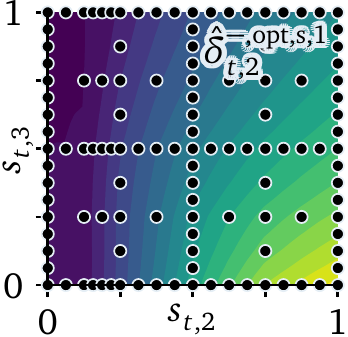}%
  }%
  \hfill%
  \makebox[37mm][c]{%
    \includegraphics{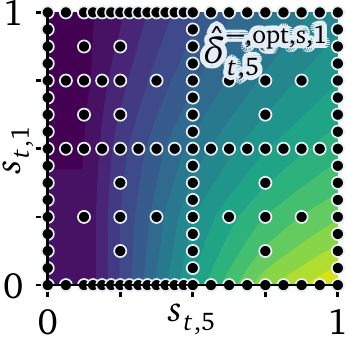}%
  }%
  \caption[Sparse grid solution for the five-dimensional TCP]{%
    Spatially adaptive sparse grid solution for the transaction costs problem
    \vspace{-0.15em}%
    with $d = 5$ stocks.
    \vspace{-0.15em}%
    Shown are slice plots of
    the value function $\normcetvalueintp_t$ \emph{(top left)} and
    the optimal policy $\optnormpolicyintp_t$
    for the initial time step $t = 0$,
    where for each function, a pair $(o_1, o_2)$
    of dimensions to be plotted was chosen,
    and the stock fractions $\stock_{t,o}$ of the other dimensions $o$
    are set to $0.1$.
    In addition, the corresponding grid points \emph{(dots)}
    are shown as the projection onto the
    $\stock_{t,o_1}$-$\stock_{t,o_2}$ plane.%
  }%
  \label{fig:financeSolution5DSparseGrid}%
\end{figure}

\Cref{fig:financeSolution2DSparseGrid,fig:financeSolution5DSparseGrid}
each display
the value function and the optimal policies corresponding to
sparse grid solutions for
$d = 2$ stocks with $\ngp_0 = \num{879}$ policy grid points or
$d = 5$ stocks with $\ngp_0 = \num{12572}$ policy grid points.
Obviously, most grid points are placed along the various kinks in the
policies.
Interestingly, experiments show that the surplus-based refinement
criterion does not place more grid points along the perfectly diagonal kink
caused by the cropping of the state space
(i.e., along $\sumfcn(\vstock_t) = 1$).
It is possible to circumvent this issue by either
transforming the domain (e.g., rotations as in \cite{Bohn18Optimally}) or
directly incorporating the distance to the diagonal into the
refinement criterion for the value function.
However, we refrain from doing so here as this does not seem to
drastically improve results.
Again, this might be due to the domination of the overall error by
the parts \ref{item:financeErrorExtrapolation} to
\ref{item:financeErrorRounding} that are not related to interpolation.

\paragraph{Pointwise error}

\begin{figure}
  \subcaptionbox{%
    $\ngp_t = 129$ (\num{5.3e-3})%
    \label{fig:financePointwiseError_1}%
  }[44mm]{%
    \clap{\includegraphics{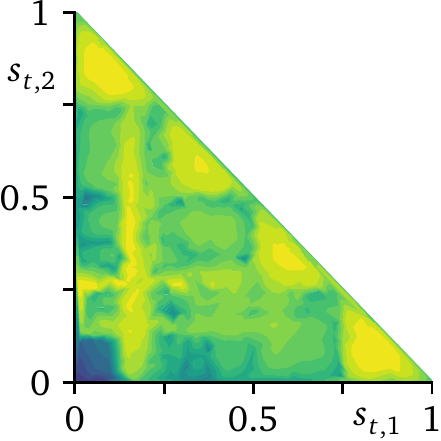}}%
  }%
  \subcaptionbox{%
    $\ngp_t = 889$ (\num{2.1e-4})%
    \label{fig:financePointwiseError_2}%
  }[44mm]{%
    \clap{\includegraphics{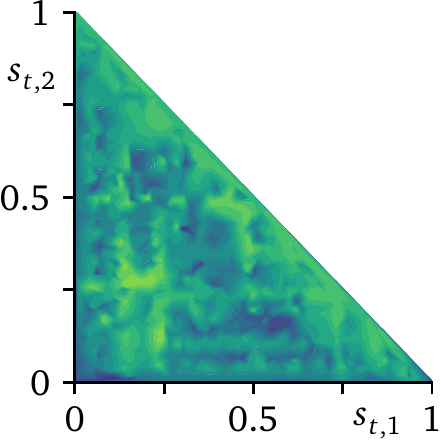}}%
  }%
  \subcaptionbox{%
    $\ngp_t = 5159$ (\num{1.2e-5})%
    \label{fig:financePointwiseError_3}%
  }[44mm]{%
    \clap{\includegraphics{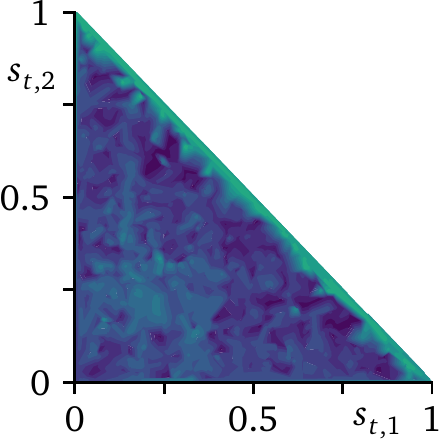}}%
  }%
  \hfill%
  \includegraphics{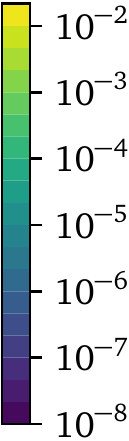}%
  \caption[Pointwise weighted Euler equation error for different grids]{%
    Pointwise weighted Euler equation error $\weightedeulererror_t(\vstock_t)$
    ($t = 0$) for the two-dimensional transaction costs problem and
    different spatially adaptive sparse grids.
    The $\Ltwo$ error $\weightedeulererrorLtwo_t$ is given in
    parentheses.%
  }%
  \label{fig:financePointwiseError}%
\end{figure}

Pointwise plots of the weighted Euler equation error
as in \cref{fig:financePointwiseError} for two stocks
reveal that there are two types of regions where the error is large:
The first type of region
is the neighborhood of the aforementioned diagonal boundary
$\sumfcn(\vstock_t) = 1$ of the uncropped region,
where the cropping distorts the error despite the weights.
The second type of region
are kinks of the optimal policy functions,
which is most visible for coarse grids
(e.g., \cref{fig:financePointwiseError_1}).
When increasing the number of grid points
(e.g., \cref{fig:financePointwiseError_1,fig:financePointwiseError_2}),
the error decreases quickly in the whole domain.

\paragraph{Monte Carlo simulation}

\begin{figure}
  \includegraphics{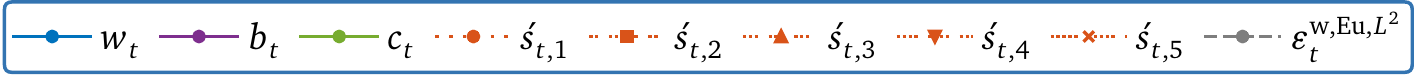}%
  \\[2mm]%
  \subcaptionbox{%
    $d = 3$ ($\ngp_0 = \num{28739}$)%
    \label{fig:financeSimulation_1}%
  }[48mm]{%
    \includegraphics{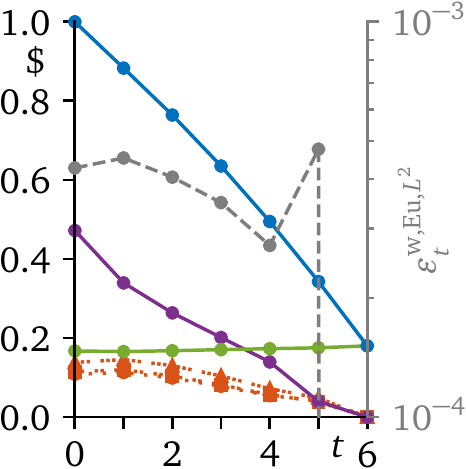}%
  }%
  \hfill%
  \subcaptionbox{%
    $d = 4$ ($\ngp_0 = \num{3343}$)%
    \label{fig:financeSimulation_2}%
  }[48mm]{%
    \includegraphics{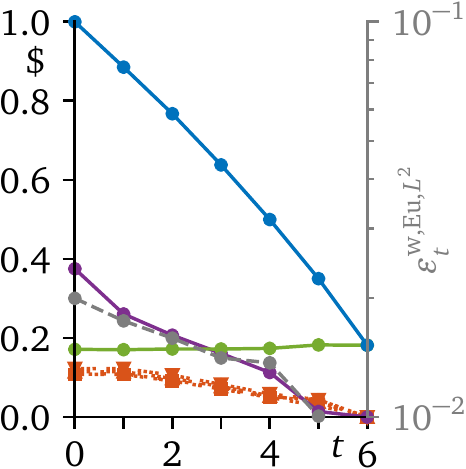}%
  }%
  \hfill%
  \subcaptionbox{%
    $d = 5$ ($\ngp_0 = \num{12572}$)%
    \label{fig:financeSimulation_3}%
  }[48mm]{%
    \includegraphics{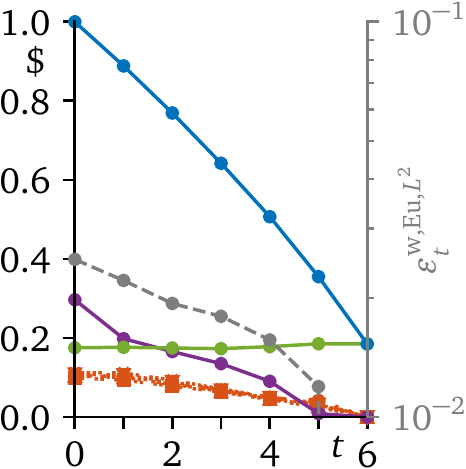}%
  }%
    \caption[Monte Carlo simulation of the TCP]{%
      Mean values of
      wealth $\wealth_t$
      \emph{\textcolor{C0}{(blue)},}
      unnormalized optimal bonds $\bond_t$
      \emph{\textcolor{C3}{(purple)},}
      unnormalized optimal consumption $\consume_t$
      \emph{\textcolor{C4}{(green)},} and
      unnormalized stock holdings
      $
        \acute{\vstock}_t
        \ceq (\vstock_t + \vnormbuy_t - \vnormsell_t) \wealth_t
      $
      after buying and selling
      \emph{\textcolor{C1}{(red)}}
      \vspace{0.1em}%
      in a Monte Carlo simulation of \num{10000} individuals,
      \vspace{-0.1em}%
      where we assume that $\wealth_0 = \$1$ for all individuals.
      In addition, the plots show the evolution of the
      $\Ltwo$ error $\weightedeulererrorLtwo_t$ over time $t$
      \emph{\textcolor{C8}{(gray, right axes)}.}%
    }%
    \label{fig:financeSimulation}%
\end{figure}

As explained in \cref{sec:828postProecessing},
we perform a multi-agent Monte Carlo simulation
and plot the resulting mean state and policy in \cref{fig:financeSimulation}
for $d = 3$, $4$, and $5$ stocks.
In addition, this figure contains the evolution of the
weighted Euler equation error $\weightedeulererrorLtwo_t$ over time.
We perform a two-part assessment of the simulated results.
First, consumption should ideally be constant over time
from a finance perspective.
We measure this by calculating the coefficients of variation
(ratio of the standard deviation of the $c_t$ values to their mean),
which is \SI{2.76}{\percent}, \SI{2.68}{\percent}, and \SI{2.58}{\percent}
for $d = 3$, $4$, and $5$, respectively.
This indicates that the variation of the consumption over time
is indeed small.
Second, we consider the so-called \term{Sharpe ratios}
\cite{Sharpe66Mutual}.
The ratios are stock fractions $\vstock$ that are determined
such that the excess stock return (compared to risk-free investment)
per unit of risk is maximized:%
\footnote{%
  The Sharpe ratios per se are derived for non-skewed
  stock return rate distributions.
  Our stock return rates are log-normally distributed and thus skewed,
  but the deviation should be small after six time steps.
  However, there are variants that
  take skewed distributions into account \cite{Mueller15Ansaetze}.%
}
{%
  \setlength{\abovedisplayskip}{9pt}%
  \setlength{\belowdisplayskip}{9pt}%
  \begin{equation}
    \vecargmax_{\vstock \in \clint{0, 1}^d}
    \frac{
      \tr{\*\mu_{\range{1}{d}}} \vstock - \bondreturn
    }{
      \sqrt{\tr{\vstock} \mat{\Sigma}_{\range{1}{d},\range{1}{d}} \vstock}
    },
  \end{equation}%
}%
where $\*\mu_{\range{1}{d}}$ and
$\mat{\Sigma}_{\range{1}{d},\range{1}{d}}$
are the first $d$ entries of $\*\mu$ and
the principal minor of order $d$ of $\mat{\Sigma}$ as given in
\cref{eq:financeStockReturnMeanCovariance}.
We compare these theoretical Sharpe ratios (left)
with the simulated stock fractions
$\acute{\stock}_{t,o}/\sumfcn(\acute{\vstock}_t)$ for $t = 0$ (right):
\begin{subequations}
  \setlength{\jot}{3pt}%
  \begin{align}
    d = 3\colon\,&
    (0.314, 0.302, 0.384\rlap{),}\hphantom{, 0.999, 0.999}\quad
    (0.300, 0.317, 0.383),\\
    d = 4\colon\,&
    (0.275, 0.185, 0.250, 0.289\rlap{),}\hphantom{, 0.999}\quad
    (0.239, 0.238, 0.253, 0.270),\\
    d = 5\colon\,&
    (0.275, 0.122, 0.176, 0.203, 0.223\rlap{),}\quad
    (0.199, 0.188, 0.197, 0.205, 0.212).
  \end{align}
\end{subequations}
The simulated stock fractions
match the predicted Sharpe ratios well for $d = 3$,
while the deviation for $d \ge 4$ is larger.
However, as the simulated stock fractions do not change much over time,
we may suspect that the skewness of the distribution of the stock return rates
limits the applicability of the Sharpe ratios to these cases.

\paragraph{Complexity and computation time}

\begin{figure}
  \includegraphics{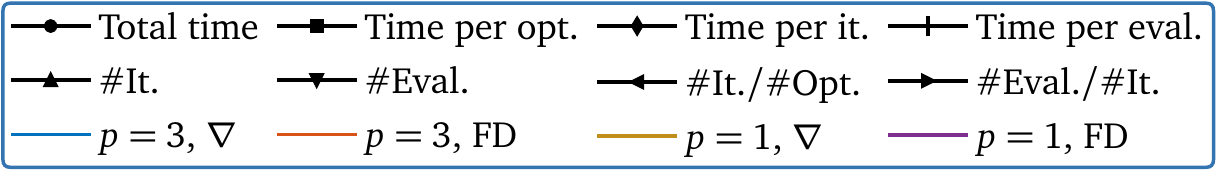}%
  \\[1mm]%
  \subcaptionbox{%
    $d = 1$%
  }[49.5mm]{%
    \includegraphics{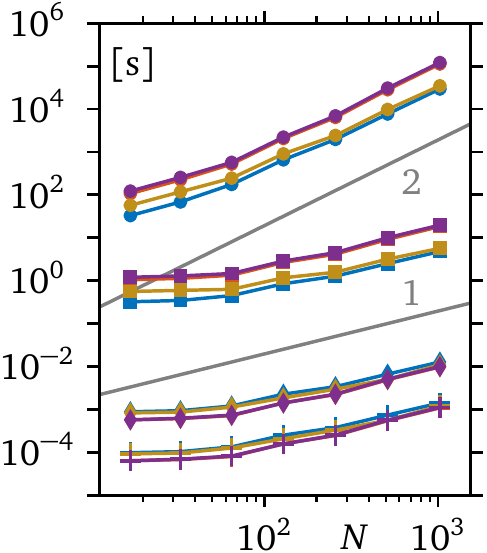}%
    \vspace*{0.3mm}\newline\hspace*{1.0mm}%
    \includegraphics{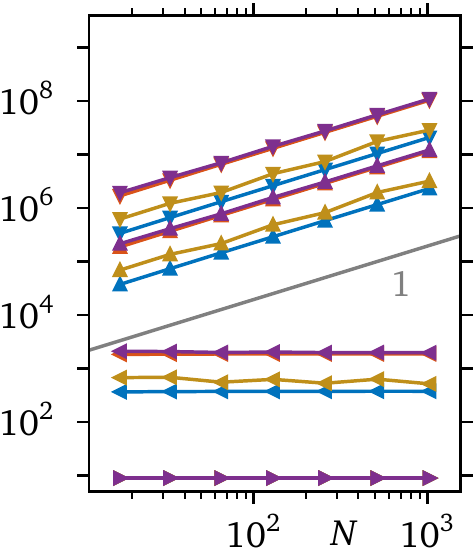}%
  }%
  \hfill%
  \subcaptionbox{%
    $d = 2$%
  }[49.5mm]{%
    \includegraphics{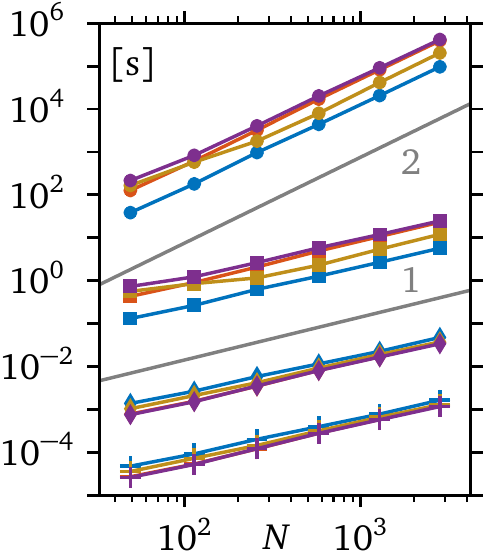}%
    \vspace*{0.3mm}\newline\hspace*{1.0mm}%
    \includegraphics{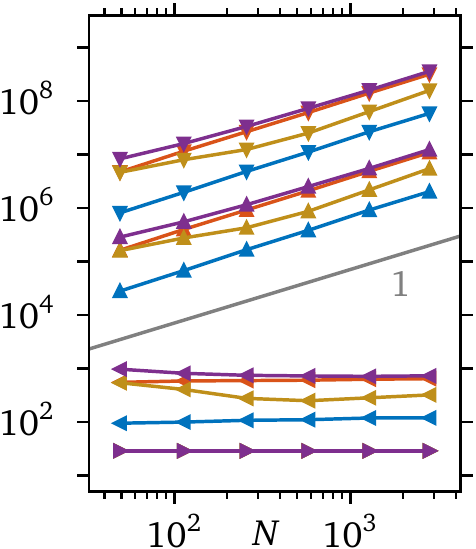}%
  }%
  \hfill%
  \subcaptionbox{%
    $d = 3$%
  }[49.5mm]{%
    \includegraphics{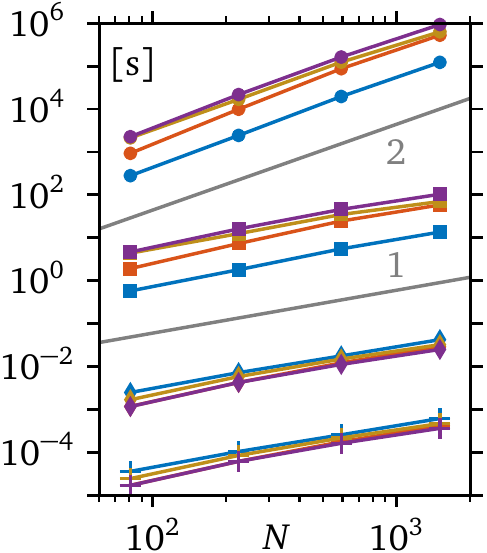}%
    \vspace*{0.3mm}\newline\hspace*{1.0mm}%
    \includegraphics{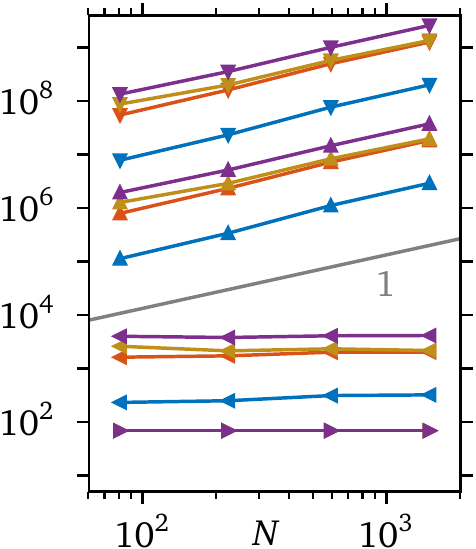}%
  }%
  \caption[Computation times and numbers of iterations for the TCP]{%
    Computation times \emph{(top)} and
    numbers of iterations and evaluations of the interpolant \emph{(bottom)}
    for the transaction costs problem on static regular sparse grids
    (i.e., without refinement).
    ``Total time'' is the serial time required to solve all
    emerging optimization problems.
    ``Time per opt.''\ is this time divided by the number
    $\mathrm{\#Opt.} = T\ngp$ of optimization problems.
    ``Time per it.''\ is the time divided by the number
    \#It.\ of optimization iterations,
    each of which is assumed to correspond to
    exactly one combined evaluation of objective function and gradient
    (the latter only if gradients are used).
    ``Time per eval.''\ is the time divided by
    the number \#Eval.\ of evaluations of
    the sparse grid interpolant and its gradient.
    The colors correspond to B-spline degrees $p = 3$ or $p = 1$ and
    to using gradients (``$\nabla$'') or finite differences (``FD'').%
  }%
  \label{fig:financeRuntime}%
\end{figure}

A complexity analysis reveals that the difficulty of solving
transaction cost problems quickly grows with the dimensionality $d$:
As shown in \cref{fig:structureSolveValueFunction},
the number of necessary arithmetic operations grows like
{%
  \setlength{\abovedisplayskip}{6pt}%
  \setlength{\belowdisplayskip}{6pt}%
  \begin{equation}
    \landauTheta{
      T
      \cdot \ngp_t
      \cdot \text{\#optimizer iterations}
      \cdot
      \underbrace{
        m_{\zeta}
        \cdot
        \overbrace{
          m_{\policy}
          \cdot
          \ngp_{t+1}
          \cdot m_{\state}
          \cdot p
        }^{\mathclap{\text{one evaluation of interpolant}}}
      }_{\mathclap{\text{one evaluation of objective gradient}}}\,
    },
  \end{equation}%
}%
where $m_{\state}, m_{\policy} \in \landauTheta{d}$ and
$m_{\zeta}, \ngp_t, \ngp_{t+1} \in \landauTheta{2^n n^{d-1}}$
if regular sparse grids of level $n$
without boundary points are used for state and stochastic grids
(due to $m_{\stochastic} = d$).
In addition,
the number of optimizer iterations is likely superlinear in $d$,
as this depends on the dimensionality $m_{\policy}$ of the search space
as well as on the number of multi-start points
(which also grows with $m_{\policy}$).
This means that the complexity is at least cubic in $d$,
quadratic in the average number $\ngp$ of employed state grid points, and
linear in the number $m_{\zeta}$ of quadrature points.
\Cref{fig:financeRuntime} confirms these observations with
experimental data.
For fixed $d$, the total time required by the optimization process
grows quadratically with the number $\ngp$ of grid points.
The time for one solution of the Bellman equation,
the time for one optimizer iteration, and
the time for one evaluation of the interpolant are all linear in $\ngp$,
as the number of optimizer iterations
is constant for fixed $d$.
If $d$ increases,
then the number of interpolant evaluations per optimizer iteration
(i.e., the number of quadrature points) increases as well.
Surprisingly, the number of optimizer iterations per grid point
and the time per evaluation are not monotonously increasing.
The latter observation might be due to vectorization effects.

\paragraph{Comparison to piecewise linear functions}

Hierarchical B-splines introduce two major benefits to
the solution of dynamic portfolio choice models.
The first benefit are the smooth objective functions:
When repeating the computations with piecewise linear functions (i.e., $p = 1$),
one obtains almost the same weighted Euler equation errors as in the cubic case
(except for the case of $d = 1$, where the error is one order of
magnitude greater than in the cubic case).
However, as we see in \cref{fig:financeRuntime},
the total computation time is several times larger
(e.g., more than five times for $d = 3$) for piecewise linear functions,
although evaluations are cheaper than for B-splines.
The main reason is that the number of required optimizer iterations
is for $p = 1$ almost seven times as high as in the cubic case,
since the optimizer has to deal with kinks in the objective function.
Experiments show that beginning with $d = 4$,
the total optimization time required to solve the transaction costs problem
is one whole order of magnitude shorter for cubic B-splines
than for piecewise linear functions.

\paragraph{Comparing exact gradients to finite differences}

The second benefit is the availability of exact gradients:
\Cref{fig:financeRuntime} also contains computation times of the
solution process
if we artificially do not use exact gradients of the objective functions,
but rather approximate them with finite differences.
For each evaluation of the objective gradient,
at least $m_{\policy}$ additional evaluations of the objective function
have to performed to compute the finite differences
($2m_{\policy}$ if central differences are used).
Consequently, while the resulting weighted Euler equation errors are
similar, the total optimization time roughly increases by a
factor of up to five if we do not use exact gradients.

\cleardoublepage

  \chapter{Conclusion}
\label{chap:90conclusion}

\noindent
Finally, we conclude the thesis by summarizing its results
and by giving an outlook on possible future work.
In particular, we highlight key contributions of the thesis to research,
give recommendations for future applications of the presented method, and
state possible downsides and limitations.

\vspace*{-0.5em}

\paragraph{Summary of the thesis}

The contribution of this thesis consisted of two major parts.
In the first part,
hierarchical B-splines on sparse grids were comprehensively presented
and embedded in a sparse grid framework with general
tensor product basis functions.
The advantage of this approach was that the framework could be reused
for different hierarchical bases
(as for the various spline bases derived in this thesis)
and that it clarified which properties only held for the
classical piecewise linear bases and not for other tensor product bases.
We saw that standard hierarchical B-splines suffer from
approximation issues near the boundary, and
we resolved these issues by incorporating not-a-knot boundary conditions
into the hierarchical B-spline basis.
In the further course of the thesis,
the focus was put on the algorithmic implications of the novel bases,
taking the hierarchization problem as an example.
We looked at requirements that had to be satisfied by grids and bases
to enable efficient hierarchization algorithms
such as breadth-first search and unidirectional principle,
for which we gave clear formulations and formal correctness proofs.
As a result, a whole ``zoo'' of hierarchical (B-)spline functions
has been derived in this thesis.
The main types were
standard hierarchical B-splines,
modified hierarchical B-splines,
hierarchical not-a-knot B-splines,
hierarchical fundamental splines, and
hierarchical weakly fundamental splines
(where the first two are not novel).
Modified, not-a-knot, and (weakly) fundamental splines could be combined
almost arbitrarily to tailor the ansatz functions to suit one's specific needs.

\pagebreak

The second, more practical part of the thesis was dedicated to transferring
the newly gained theoretical knowledge to
academic and real-life application test cases.
We verified that only with the new hierarchical not-a-knot conditions,
one is able to obtain the best possible order of convergence
$\landauO{\ms{n}^{p+1}}$ for interpolation
(B-spline degree $p$, fixed dimensionality $d$).
Using the Novak--Ritter criterion,
which was specifically designed for optimization,
we were able to achieve optimization gaps
that were for some test functions up to six orders of magnitude smaller
for cubic B-splines than for standard piecewise linear functions.
We transferred the Novak--Ritter criterion to uncertainty quantification
and obtained similarly strong results for the
propagation of fuzzy uncertainties with the fuzzy extension principle.
Furthermore, we successfully showed the suitability of hierarchical B-splines
for three real-world applications,
which are summarized in \cref{tbl:applicationSummary}.
\begin{table}
  \setnumberoftableheaderrows{1}%
  \begin{tabular}{%
    >{\kern\tabcolsep}=l<{\kern5mm}+l+l+l<{\kern\tabcolsep}%
  }
    \toprulec
    \headerrow
    Category&                   Topology opt.&      Biomechanics&      Finance\\
    \midrulec
    Interpolated quantities&    Elasticity tensors& Muscle forces&     Value functions\\
    SG dimensionality&          5&                  2&                 5\\
    \#Optimization variables&   \num{40000}&        2&                 11\\
    Time per evaluation&        \SI{30}{\second}&   \SI{30}{\minute}&  ---\\
    \#Eval. per opt. iteration& \num{8000}&         4&                 150\\
    Objective function type&    Non-linear&         Linear/non-linear& Non-linear\\
    Constraint function type&   Non-linear&         Non-linear/---&    Linear\\
    Optimization method&        SQP&                Augm. Lagrangian&  SQP\\
    \bottomrulec
  \end{tabular}
  \caption[%
    Summary of characteristics of the applications%
  ]{%
    Summary of characteristics of the applications presented in this thesis.
    The given values are rough example values that
    represent possible application test cases.%
  }%
  \label{tbl:applicationSummary}%
\end{table}%
First, by interpolating Cholesky factors of elasticity tensors,
we accomplished to efficiently solve
topology optimization problems in three spatial dimensions
with complex micro-cell structures.
Second, in the biomechanical application,
we dramatically reduced the computational time to solve
test scenarios by up to \SI{99}{\percent} by using
sparse grid surrogates with B-splines instead of the
exact continuum-mechanical model.
Third, we were able to solve dynamic portfolio choice problems
with five state variables and eleven policy variables
with unprecedented precision, as one could only speculate
how the solution looked like with state-of-the-art methods.
In all of these applications, the advantages of B-splines were made clear
by comparing the results to the classical piecewise linear basis.
The implementation of hierarchical B-splines of sparse grids is
publicly available as part of the sparse grid toolbox \sgpp
under a free and open-source license.%
\footnote{%
  \url{http://sgpp.sparsegrids.org/}%
}

\vspace*{\fill}
\pagebreak

\paragraph{Recommendations (advantages and disadvantages)}

Despite its broad applicability,
the presented method of B-splines on sparse grids is of course
not suited for all possible scenarios.
One must be able to sample the objective function at arbitrary
locations in some hyper-rectangle in order to use sparse grids;
a prescribed point cloud of scattered data does not suffice.
Moreover, the problem should not have more than ten dimensions,
since convergence notably slows down as the dimensionality grows,
although spatially adaptive approaches might still be feasible for
higher dimensionalities \cite{Pflueger10Spatially}.
In addition, the objective function should be ``as smooth as possible''
in order to benefit from higher-order B-splines.
This means ``continuous'' at the very least,
but twice continuous differentiability is more desirable.
The general rule is that the objective function should be
at least as smooth as the employed basis functions in order to obtain
optimal convergence results.
The concrete choice of basis (general type and degree) depends
on the application:
Not-a-knot B-splines are well-suited for objective functions with
dominating near-polynomial parts.
Fundamental splines may be used to accelerate the process of
hierarchization by enabling breadth-first search in quadratic time.
With weakly fundamental splines, this can be further reduced to
linear time using the unidirectional principle.
However, the additional grid points that have to be inserted
have to be taken into account as well.
The rule of thumb is that the more spatially adaptive a sparse grid is
(i.e., only few high-level grid points),
the more points have to be inserted.
In general, it does not hurt to try the different available
B-spline types and degrees,
since most function values can simply be reused
once the objective function has been sampled.

\vspace*{\fill}

\paragraph{Outlook and future work}

Finally, we briefly give suggestions for possible future work.
A major topic of interest is that of refinement criteria and adaptivity.
Besides the Novak--Ritter criterion,
there are other refinement criteria that are tailored to optimization
such as simultaneous optimistic optimization \cite{Wang14Bayesian}.
In addition, nested methods for hierarchical optimization could
use multiple interpolants with different resolutions
on different grids \cite{Delbos14Global}.
Criteria that directly incorporate constraints would improve
results in constrained optimization settings.
With respect to adaptivity,
there is also much work left to do.
This thesis focused on spatial adaptivity for its applications,
but there are interesting applications
that greatly benefit from dimensional adaptivity,
for example plasma physics \cite{Pflueger14EXAHD}.
Another key task would be the introduction of
$h$-$p$-adaptivity to B-splines on sparse grids,
which would greatly enhance the applicability of B-splines in
non-smooth scenarios.
As a simple special case,
one could investigate different B-spline degrees in different dimensions.
However, true $h$-$p$-adaptivity would allow to locally choose
both the spatial resolution $h$ and the B-spline degree $p$,
adapting them according to the local smoothness of the function.

\pagebreak

\begin{SCfigure}
  \begin{tikzpicture}
    \node[anchor=south west] at (0mm,0mm) {%
      \includegraphics[height=80mm]{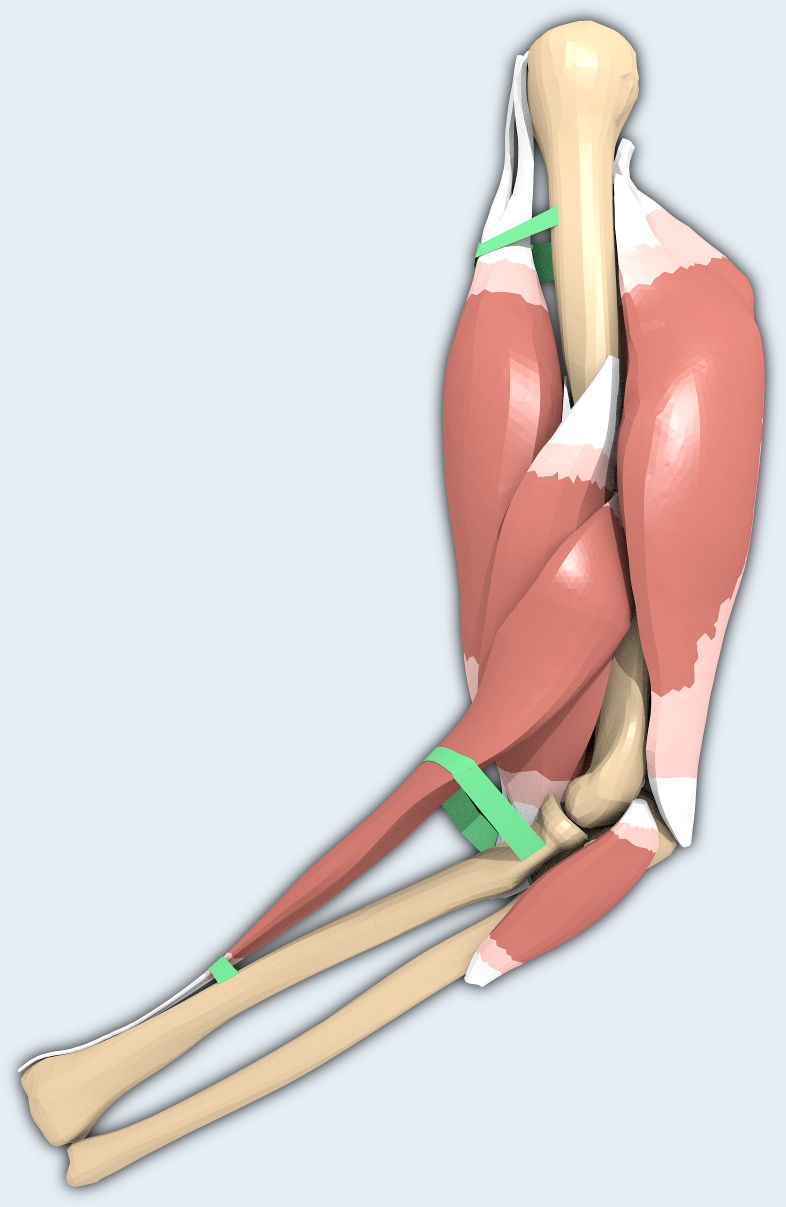}%
    };
    \node[anchor=west] (humerus)    at (0mm,75mm)  {\small{}Humerus};
    \node[anchor=west] (radius)     at (0mm,20mm)  {\small{}Radius};
    \node[anchor=east] (ulna)       at (40mm,10mm) {\small{}Ulna};
    \node[anchor=east] (triceps)    at (73mm,52mm) {\small{}Triceps};
    \node[anchor=west] (biceps)     at (0mm,52mm)  {\small{}Biceps};
    \node[anchor=west] (brachialis) at (0mm,44mm)  {\small{}Brachialis};
    \node[anchor=west] (brachioradialis) at (0mm,36mm) {%
      \small{}Brachioradialis%
    };
    \node[anchor=east] (anonceus)   at (73mm,20mm) {\small{}Anonceus};
    \draw (humerus)         -- (39mm,75mm);
    \draw (radius)          -- ($(0mm,10mm)!(radius)!(20mm,10mm)$);
    \draw (ulna)            -- (20mm,10mm);
    \draw (triceps)         -- (48mm,52mm);
    \draw (biceps)          -- (33mm,52mm);
    \draw (brachialis)      -- (36mm,44mm);
    \draw (brachioradialis) -- (36mm,36mm);
    \draw (anonceus)        -- (38mm,20mm);
  \end{tikzpicture}%
  \caption[Extended model of the human upper limb with five muscles]{%
    Extended model of the human upper limb with
    the three bones humerus, radius, and ulna
    \emph{\textcolor[RGB]{191,162,118}{(light brown)}} and
    the five muscles
    triceps, biceps, brachialis, brachioradialis, and anonceus
    \emph{\textcolor[RGB]{150,58,46}{(red)}.}
    Each muscle has associated tendon
    \emph{\textcolor[RGB]{128,128,128}{(white)}}
    and muscle-tendon complexes
    \emph{\textcolor[RGB]{255,160,160}{(light red)}.}
    Biceps and brachioradialis have each been fixed with one and two bands
    \emph{\textcolor[RGB]{80,200,100}{(green)}}
    to simulate the effect of the missing skin.
    Without bands, the muscles would unnaturally raise from the bones.%
  }%
  \label{fig:upperLimb5D}%
\end{SCfigure}

With regard to the application side, there are also quite a few
possibilities for future work.
The sparse grids in the biomechanical application
we considered in this thesis were only two-dimensional.
This is not in the range of dimensions in which sparse grids
demonstrate their full strength,
although the two-dimensional surrogates were already able to drastically
reduce the required computation time compared to full grids.
Currently, an extended model with five muscles and therefore
five-dimensional sparse grids is being considered
(see \cref{fig:upperLimb5D}).
Here, it will be mandatory to employ spatial adaptivity
to cope with the increased dimensionality.
In the application of topology optimization,
more complicated micro-cell models and more complicated
settings could be studied.
For example, the widths of the diagonal macro-bars could be constrained
\cite{Allaire16Towards}.
The dynamic portfolio choice models in the financial application
were quite limited.
For example,
there was no inheritance motive (i.e., bequest),
the model did not contain the individual's regular income, and
the model did not account for savings for large necessary investments
(e.g., cars or houses), which seems unnecessarily unrealistic.
Finally, one could consider many other real-world optimization problems
or other application fields of B-splines on sparse grids,
for instance, data mining or uncertainty quantification.
If objective gradients are available besides function values,
it might be feasible to directly incorporate the gradients into
the interpolation scheme \cite{Baar15Gradient}.

This extensive but by no means exhaustive list of possible future work
can be seen as an inspiration and starting point
for new and interesting applications of B-splines for sparse grids.

\cleardoublepage

  \appendix
  
  \chapter{Proofs}
\label{chap:a10proofs}

This chapter contains proofs that
were too long or too technical to include them in the main text.
For convenience, the corresponding propositions and theorems are repeated
before the proofs, using the same numbering as in their
original chapter.

\addtocontents{lop}{\protect\iffalse}
\isincluded{tex/document/20sparseGrids}{\section{Proofs for Chapter 2}
\label{sec:a11chapter2}

\disableornamentsfornextheadingtrue
\subsection{Proof of the Size of the Regular Sparse Grid with Coarse Boundaries}
\label{sec:a111proofGridSizeCoarseBoundary}

\propGridSizeCoarseBoundary*

\begin{proof}
  Note that the outer union in the definition of $\coarselevelset{n}{d}{b}$ in
  \eqref{eq:coarseBoundary2} is indeed disjoint.
  Therefore,
  \begin{equation}
    \setsize{\coarseregsgset{n}{d}{b}}
    = \sum_{\substack{\*l \in \nat^d\\\normone{\*l} \le n}}
    \setsize{\hiset{\*l}} +
    \hspace*{-5mm}\sum_{
      \substack{
        \*l \in \natz^d \setminus \nat^d\\
        (\normone{\vecmax(\*l, \*1)} \le n - b + 1) \lor
        (\*l = \*0)
      }
    }\hspace*{-10mm} \setsize{\hiset{\*l}}.
  \end{equation}
  The first sum is the number $\setsize{\interiorregsgset{n}{d}}$
  of interior grid points in $\regsgset{n}{d}$.
  The second sum can be split into summands
  with the same number $q$ of zero entries,
  which we count with
  $N_\*l \ceq \setsize{\{t \mid l_t = 0\}}
  = \normone{\vecmax(\*l, \*1)} - \normone{\*l}$,
  and the same level sum $m = \normone{\*l}$:
  \begin{equation}
    \setsize{\coarseregsgset{n}{d}{b}}
    = \setsize{\interiorregsgset{n}{d}} + 2^d +
    \sum_{q=1}^{d-1} \sum_{m=d-q}^{n-b-q+1}
    \sum_{
      \substack{
        \*l \in \natz^d\\
        N_\*l = q,\,\normone{\*l} = m}
    }\hspace*{-4mm} \setsize{\hiset{\*l}},
  \end{equation}
  where $2^d$ is the summand for $\*l = \*0$
  (number $\setsize{\hiset{\*0}}$ of corners of $\clint{\*0, \*1}$).
  The limits of the sum over $m$ are $d-q$,
  since there are $d-q$ entries $\ge 1$ in a level vector
  with $q$ zero entries, and $n-b-q+1$,
  since $m = \normone{\*l}
  = \normone{\vecmax(\*l, \*1)} - N_\*l
  \le n-b+1 - N_\*l
  = n-b-q+1$.
  
  \vspace*{2em}
  
  In general, the innermost summand $\setsize{\hiset{\*l}}$ equals
  $\setsize{\hiset{\*l}}
  = \prod_{\{t \mid l_t \ge 1\}} 2^{l_t-1} \cdot \prod_{\{t \mid l_t = 0\}} 2
  = 2^{\normone{\*l} - d + 2N_\*l}$.
  The number of innermost summands is given by
  \begin{equation}
    \setsize{\{\*l \in \natz^d \mid N_\*l = q,\; \normone{\*l} = m\}}
    = \binom{d}{q} \setsize{\{\*l \in \nat^{d-q} \mid \normone{\*l} = m\}}
    = \binom{d}{q} \binom{m-1}{d-q-1}.
  \end{equation}
  This can be seen by first putting $q$ zeros in $d$ places,
  for which there are $\binom{d}{q}$ possibilities, and then
  counting all positive vectors of length $d - q$ with level sum $m$,
  which can be done in $\binom{m-1}{d-q-1}$ ways.
  Thus,
  \begin{subequations}
    \begin{align}
      \setsize{\coarseregsgset{n}{d}{b}}
      &= \setsize{\interiorregsgset{n}{d}} + 2^d +
      \sum_{q=1}^{d-1} \binom{d}{q} \sum_{m=d-q}^{n-b-q+1}
      2^{m - d + 2q} \binom{m-1}{d-q-1}.\\[1em]
      \intertext{%
        After shifting the index $m \to (m + d - q)$ and slightly
        rearranging the terms, we obtain%
      }
      \cdots
      &= \setsize{\interiorregsgset{n}{d}} + 2^d +
      \sum_{q=1}^{d-1} 2^q \binom{d}{q} \sum_{m=0}^{n-d-b+1}
      2^m \binom{(d-q)-1+m}{(d-q)-1}.\\[1em]
      \intertext{%
        We can now use \thmref{lemma:numberOfGridPointsInterior} to
        conclude that%
      }
      \cdots
      &= \setsize{\interiorregsgset{n}{d}} +
      \sum_{q=1}^d 2^q \binom{d}{q}
      \setsize{\interiorregsgset{n-q-b+1}{d-q}}
    \end{align}
  \end{subequations}
  as desired.
\end{proof}

\vspace*{2em}

\breakpagebeforenextheadingtrue
\subsection{%
  Correctness Proof of the Construction of the Regular Sparse Grid
  with Coarse Boundaries%
}
\label{sec:a112proofInvariantCoarseBoundary}

\propInvariantCoarseBoundary*

\begin{proof}
  First, we show that every inserted level $\*l' \in \natz^t$ in the inner loop
  can be found on the right-hand side of \eqref{eq:coarseInvariant}.
  If $\*l' \ceq (\*l, 0)$
  is inserted for some $\*l \in \levelset^{(t-1)}$,
  then we have $\normone{\vecmax(\*l, \*1)} \le n-d+t-b$ or
  $\*l = \*0$ by \cref{line:algCoarseBoundary1} of
  \cref{alg:coarseBoundary}.
  In the first case, we have
  \begin{equation}
    \normone{\vecmax(\*l', \*1)}
    = \normone{\vecmax(\*l, \*1)} + 1
    \le n - d + t - b + 1,
  \end{equation}
  and in the second case $\*l' = \vec{0}$.
  In either case, $\*l'$ is contained in the \rhs of
  \eqref{eq:coarseInvariant}.
  
  If $\*l' \ceq (\*l, l_t)$ is inserted
  for some $\*l \in \levelset^{(t-1)}$ and
  $l_t \in \{1, \dotsc, l^\ast\}$, then there are,
  depending on whether $\*l \in \nat^{t-1}$, two cases:
  \begin{itemize}
    \item
    \mbox{If $\*l \in \nat^{t-1}$, then $\*l' \in \nat^t$ and}
    $\normone{\*l'} \le \normone{\*l} + l^\ast = n - d + t$
    due to \cref{line:algCoarseBoundary2},
    i.e., $\*l'$ is contained in the first set of the \rhs of
    \eqref{eq:coarseInvariant}.
    
    \item
    \mbox{If $\*l \notin \nat^{t-1}$, then $\*l' \notin \nat^t$ and}
    $\normone{\vecmax(\*l', \*1)}
    \le \normone{\vecmax(\*l, \*1)} + l^\ast
    = n - d + t - b + 1$
    due to \cref{line:algCoarseBoundary3},
    i.e., $\*l$ is contained in the second set of the \rhs of
    \eqref{eq:coarseInvariant}.
  \end{itemize}
  Thus, all levels that the algorithm inserts into $\levelset^{(t)}$
  can be found on the \rhs of \eqref{eq:coarseInvariant}.
  
  It remains to prove that all levels on the \rhs of
  \eqref{eq:coarseInvariant}
  are eventually inserted by the algorithm into $\levelset^{(t)}$.
  We prove this by induction over $t = 1, \dotsc, d$.
  For $t = 1$, the \rhs of \eqref{eq:coarseInvariant} equals
  $\{l \in \natz \mid l \le n-d+1\}$, which is just $\levelset^{(1)}$
  (see \cref{line:algCoarseBoundary6} of \cref{alg:coarseBoundary}).
  For the induction step $(t - 1) \to t$, we assume
  the validity of the induction hypothesis
  \begin{equation}
    \label{eq:proofCoarseInductionHypothesis}
    \begin{aligned}
      \levelset^{(t-1)}
      &= \{\*l \in \nat^{t-1} \mid
      \normone{\*l} \le n-d+t-1\} \dotcup {}\\
      &\hphantom{{}={}}
      \paren*{
        \{\*l \in \natz^{t-1} \setminus \nat^{t-1} \mid
        \normone{\vecmax(\*l, \*1)} \le n-d+t-b\} \cup
        \{\vec{0}\}
      }.
    \end{aligned}
  \end{equation}
  The \rhs of \eqref{eq:coarseInvariant} has three parts,
  so we check for elements $\*l' \in \natz^t$
  of each of the three sets that they are appended to $\levelset^{(t)}$
  eventually.
  
  First, let $\*l' = (\*l, l_t)$ be in the first set of the \rhs,
  i.e., $\*l' \in \nat^t$ (in particular $l_t \ge 1$) and
  $\normone{\*l'} \le n - d + t$.
  Note that $\*l$ will be encountered in the inner loop, as
  $\*l \in \nat^{t-1}$ and
  $\normone{\*l} = \normone{\*l'} - l_t \le n - d + t - 1$,
  which implies $\*l \in \levelset^{(t-1)}$ by the induction
  hypothesis \eqref{eq:proofCoarseInductionHypothesis}.
  Since $1 \le l_t \le l^\ast$
  (due to
  $l_t = \normone{\*l'} - \normone{\*l} \le n-d+t - \normone{\*l} = l^\ast$),
  the level $\*l'$ is inserted into $\levelset^{(t)}$ during the innermost loop
  in \cref{line:algCoarseBoundary4} of \cref{alg:coarseBoundary}.
  
  Second, let $\*l' = (\*l, l_t)$
  be in the second set of the \rhs, i.e., we have
  $\*l' \notin \nat^t$ and
  $\normone{\vecmax(\*l', \*1)} \le n-d+t-b+1$.
  Here, there are three cases:
  \begin{enumerate}
    \item
    $l_t \ge 1$:
    This implies $\*l \notin \nat^{t-1}$ and 
    $\normone{\vecmax(\*l, \*1)}
    = \normone{\vecmax(\*l', \*1)} - l_t
    \le n-d+t-b$.
    Consequently, $\*l \in \levelset^{(t-1)}$ by the induction hypothesis
    \eqref{eq:proofCoarseInductionHypothesis}.
    As $1 \le l_t \le l^\ast$
    (due to $l_t \le n-d+t-b+1 -
    \normone{\vecmax(\*l, \*1)} = l^\ast$),
    $\*l$ is added to $\levelset^{(t)}$ in \cref{line:algCoarseBoundary4}.
    
    \item
    $l_t = 0$ and $\*l \in \nat^{t-1}$:
    This implies $\normone{\*l} = \normone{\*l'}
    = \normone{\vecmax(\*l', \*1)} - 1
    \le n - d + t - b
    \le n - d + t - 1$ since $b \ge 1$.
    Again, by the induction hypothesis
    \eqref{eq:proofCoarseInductionHypothesis},
    $\*l$ is added to $\levelset^{(t)}$ in \cref{line:algCoarseBoundary5}
    due to
    $\normone{\vecmax(\*l, \*1)}
    = \normone{\*l} \le n - d + t - b$.
    
    \item
    $l_t = 0$ and $\*l \notin \nat^{t-1}$:
    This implies $\normone{\vecmax(\*l, \*1)}
    = \normone{\vecmax(\*l', \*1)} - 1
    \le n - d + t - b$.
    Again, by the induction hypothesis
    \eqref{eq:proofCoarseInductionHypothesis},
    $\*l$ is added to $\levelset^{(t)}$ in \cref{line:algCoarseBoundary5}.
  \end{enumerate}
  
  \noindent
  Third, let $\*l = (\vec{0}, 0) \in \natz^t$
  be in the third set of the \rhs.
  This level is appended in \cref{line:algCoarseBoundary5}
  to $\levelset^{(t)}$, since $\*l' = \vec{0} \in \natz^{t-1}$
  is in $\levelset^{(t-1)}$ by the
  induction hypothesis~\eqref{eq:proofCoarseInductionHypothesis}.
\end{proof}
}{}
\isincluded{tex/document/30BSplines}{\section{Proofs for Chapter 3}
\label{sec:a12chapter3}

\disableornamentsfornextheadingtrue
\subsection{Proof of the Linear Independence of Hierarchical B-Splines}
\label{sec:a121proofHierBSplineLinearlyIndependent}

\propHierBSplineLinearlyIndependent*

\begin{proof}
  The proof is rigorous for the common B-spline degrees of
  $p \in \{1, 3, 5, 7\}$.
  For higher degrees, the proof has to be viewed as a sketch.
  
  We follow the presentation in \cite{Valentin16Hierarchical} and
  prove the assertion by induction over $l \in \natz$.
  For $l = 0$, the B-splines $\bspl{0,i'}{p}$ with $i' \in \{0, 1\}$
  are linearly independent.
  For the induction step $(l-1) \to l$, let
  \begin{equation}
    \label{eq:proofPropHierBSplineLinearlyIndependent1}
    \sum_{l'=0}^l \sum_{i' \in \hiset{l'}} \surplus{l',i'} \bspl{l',i'}{p}
    \equiv 0
  \end{equation}
  be a linear combination of the zero function.
  We separate the summands of level $l$
  from the summands of coarser levels $l' < l$:
  \begin{equation}
    \sum_{i \in \hiset{l}} \surplus{l,i} \bspl{l,i}{p}
    =: g_1 \equiv g_2 \ceq
    -\sum_{l'=0}^{l-1} \sum_{i' \in \hiset{l'}} \surplus{l',i'} \bspl{l',i'}{p}.
  \end{equation}
  The right-hand side $g_2$ is smooth in every grid point
  $\gp{l,i}$ of level $l$ ($i \in \hiset{l}$),
  since these grid points are not knots of the hierarchical B-splines
  $\bspl{l',i'}{p}$ of level $l' < l$ ($i' \in \hiset{l'}$).
  This implies that the left-hand side $g_1$ must be smooth there as well:
  \begin{equation}
    \label{eq:proofPropHierBSplineLinearlyIndependent2}
    \underbrace{
      \sum_{i \in \hiset{l}} \surplus{l,i}
      \partialdiff_-^p \bspl{l,i}{p}(\gp{l,i'})
    }_{{} = \partialdiff_-^p g_1(\gp{l,i'})}
    =
    \underbrace{
      \sum_{i \in \hiset{l}} \surplus{l,i}
      \partialdiff_+^p \bspl{l,i}{p}(\gp{l,i'})
    }_{{} = \partialdiff_+^p g_1(\gp{l,i'})},\quad
    i' \in \hiset{l},
  \end{equation}
  where $\partialdiff_-^p$ and $\partialdiff_+^p$ denote the left and right
  derivative of order $p$, respectively.
  By repeated application of \eqref{eq:cardinalBSplineDerivative},
  one can show that
  \begin{equation}
    \partialdiff_-^p \cardbspl{p}(k + 1)
    = (-1)^k \binom{p}{k}
    = \partialdiff_+^p \cardbspl{p}(k),\quad
    k \in \integer,
  \end{equation}
  where $\binom{p}{k} = 0$ for $k < 0$ or $k > p$
  \cite{Hoellig13Approximation}.
  We can insert this relation into
  \eqref{eq:proofPropHierBSplineLinearlyIndependent2}
  and use \eqref{eq:uniformHierarchicalBSplineUV} to obtain
  \begin{alignat}{2}
    \sum_{i \in \hiset{l}} \surplus{l,i} (-1)^{k-1} \binom{p}{k-1}
    &= \sum_{i \in \hiset{l}} \surplus{l,i} (-1)^k \binom{p}{k},\quad
    &&i' \in \hiset{l},\quad
    k \ceq \frac{p+1}{2} + i' - i.
  \end{alignat}
  As $\binom{p}{k-1} + \binom{p}{k} = \binom{p+1}{k}$
  and $(-1)^k$ is constant for $i \in \hiset{l}$ when $i'$ is fixed,
  this is equivalent to
  \begin{equation}
    \label{eq:proofPropHierBSplineLinearlyIndependent3}
    \sum_{i \in \hiset{l}} \surplus{l,i}
    \binom{p+1}{\frac{p+1}{2} + i' - i} = 0,\quad
    i' \in \hiset{l}.
  \end{equation}
  
  This is a square system of linear equations
  whose system matrix
  $\mat{A}(p)$ is a banded symmetric Toeplitz matrix%
  \footnote{%
    The entries $A_{k,j}$ of a Toeplitz matrix $\mat{A}$
    solely depend on $k - j$, i.e.,
    $A_{k,j} = c_{k-j}$ for some vector $\*c$.%
  }
  of size
  $2^{l-1} \times 2^{l-1}$ with bandwidth $\ceil{\tfrac{p-1}{4}}$.
  The non-zero values of $\mat{A}(p)$ are tabulated for some degrees $p$
  in \cref{tbl:proofHierBSplineLinearlyIndependent}.
  For $p = 1, 3, 5, 7$,
  the corresponding matrices are diagonally dominant and therefore regular.
  For higher B-spline degrees $p$, the regularity of $\mat{A}(p)$ has
  to be shown differently.
  
  \begin{table}
    \setnumberoftableheaderrows{1}%
    \begin{tabular}{>{\kern\tabcolsep}=l<{\kern7mm}+c+c+c+c<{\kern\tabcolsep}}
      \toprulec
      \headerrow
      &$k = 0$&$k = 1$&$k = 2$&$k = 3$\\
      \midrulec
      $p = 1$&2&&&\\
      $p = 3$&6&1&&\\
      $p = 5$&20&6&&\\
      $p = 7$&70&28&1&\\
      $p = 9$&252&120&10&\\
      $p = 11$&924&495&66&1\\
      \bottomrulec
    \end{tabular}%
    \caption[%
      Non-zero matrix values in the proof of linear independence%
    ]{%
      Non-zero values $A_{j,j+k}(p)$ of the diagonals of $\mat{A}(p)$
      obtained in \eqref{eq:proofPropHierBSplineLinearlyIndependent3}.%
    }%
    \label{tbl:proofHierBSplineLinearlyIndependent}%
  \end{table}
  
  Due to the regularity of $\mat{A}(p)$, we infer from
  \eqref{eq:proofPropHierBSplineLinearlyIndependent3} that
  $\surplus{l,i} = 0$ for $i \in \hiset{l}$.
  According to \eqref{eq:proofPropHierBSplineLinearlyIndependent1},
  we obtain
  a linear combination of the zero function with the hierarchical
  B-splines of level $< l$, i.e.,
  \begin{equation}
    \sum_{l'=0}^{l-1} \sum_{i' \in \hiset{l'}} \surplus{l',i'} \bspl{l',i'}{p}
    = 0,
  \end{equation}
  which implies $\surplus{l',i'} = 0$ for all $l' < l$ and $i' \in \hiset{l'}$
  by the induction hypothesis.
  Thus, the hierarchical B-splines $\bspl{l',i'}{p}$
  ($l' \le l$, $i' \in \hiset{l'}$) are linearly independent.
\end{proof}
}{}
\isincluded{tex/document/40algorithms}{\section{Proofs for Chapter 4}
\label{sec:a13chapter4}

\disableornamentsfornextheadingtrue
\subsection{Combinatorial Proof of the Combination Technique}
\label{sec:a131proofCombiTechnique}

\addtocontents{lop}{\protect\fi}

\begin{definition}[binomial coefficient for integer parameters]
  \label{def:binomialCoefficient}
  The binomial coefficient $\binom{n}{k}$ is defined for
  $n \in \natz$ and $k \in \integer$ as
  \begin{equation}
    \binom{n}{k}
    \ceq
    \begin{cases}
      \frac{n (n - 1) \dotsm (n - (k-1))}{k!},&0 < k < n,\\
      1,&(k = 0) \lor (k = n),\\
      0,&(k < 0) \lor (k > n).
    \end{cases}
  \end{equation}
\end{definition}

\begin{lemma}[inclusion-exclusion counting lemma]
  \label{lemma:inclusionExclusionCountingLemma}
  For $a \in \natz$, $r \ge a$, and $s \in \integer$, we have
  \begin{equation}
    \sum_{q=0}^a (-1)^q \binom{a}{q} \binom{r-q}{s}
    = \binom{r-a}{s-a}.
  \end{equation}
\end{lemma}

\begin{proof}
  We apply the upper negation formula
  (see Eq. (5.14) of \cite{Graham94Concrete})
  to the second binomial of the \lhs:
  \begin{subequations}
    \begin{align}
      \sum_{q=0}^a (-1)^q \binom{a}{q} \binom{r-q}{s}
      &= (-1)^s \sum_{q=0}^a \binom{a}{0+q} \binom{(s-r-1)+q}{s} (-1)^q.\\
      \intertext{%
        This sum can be simplified using the identity
        in Eq.\ (5.24) of \cite{Graham94Concrete}
        (the sum has already been written in the same way as in
        \cite{Graham94Concrete}):%
      }
      \cdots
      &= (-1)^{s+a} \binom{s-r-1}{s-a}.\\
      \intertext{%
        Applying the upper negation formula again,%
      }
      \cdots
      &= \binom{r-a}{s-a},
    \end{align}
  \end{subequations}
  we obtain the desired quantity.
\end{proof}

\addtocontents{lop}{\protect\iffalse}

\propCombiTechniqueOne*

\begin{proof}
  Let $q = 0, \dotsc, d - 1$ and $\gp{\*l,\*i} \in \regsgset{n}{d}$, i.e.,
  $\normone{\*l} \le n$ and $\*i \in \hiset{\*l}$.
  Note that for $\*l' \in \natz^d$, we have
  $\fgset{\*l'} \ni \gp{\*l,\*i} \iff \*l' \ge \*l$.
  Hence,
  \begin{subequations}
    \begin{align}
      \setsize{
        \{\*l' \mid \normone{\*l'} = n - q,\; \fgset{\*l'} \ni \gp{\*l,\*i}\}
      }
      &= \setsize{
        \{\*l' \mid \normone{\*l'} = n - q,\; \*l' \ge \*l\}
      }\\
      &= \setsize{
        \{\*a \in \natz^d \mid \normone{\*a} = n - q - \normone{\*l}\}
      }\\
      \intertext{%
        by mapping $\*a \ceq \*l' - \*l$.
        The size of the last set is
        known as the number of \term{weak compositions}
        of $n - q - \normone{\*l}$ into $d$ parts
        and can be computed as%
      }
      \cdots
      &= \binom{n - q - \normone{\*l} + d - 1}{d - 1},
    \end{align}
  \end{subequations}
  see Theorem 2.2 of \cite{Bona15Introduction}.
  Now, we can use \cref{lemma:inclusionExclusionCountingLemma}
  with the values
  $a \ceq s \ceq d - 1$ and
  $r \ceq n - \normone{\*l} + d - 1$
  to conclude that the \lhs of the assertion
  \eqref{eq:combiTechniqueOne} equals
  \begin{equation}
    \sum_{q=0}^{d-1} (-1)^q \binom{d-1}{q} \cdot
    \binom{n - q - \normone{\*l} + d - 1}{d - 1}
    = \binom{n - \normone{\*l}}{0}
    = 1,
  \end{equation}
  proving the proposition.
\end{proof}

\addtocontents{lop}{\protect\fi}

\begin{shortlemma}[relation is equivalence relation]
  \label{lemma:combiTechniqueEquivalenceRelation}
  $\eq$ is an equivalence relation.
\end{shortlemma}

\begin{proof}
  We check reflexivity, symmetry, and transitivity of $\eq$:
  \begin{itemize}
    \item
    \emph{Reflexivity:}
    Using the same level $\*l' = \*l''$ implies
    $T_{\*l',\*l'} = \{t \mid l'_t < l_t\}$.
    For all $t \notin T_{\*l',\*l'}$, we have $l'_t \ge l_t$.
    Consequently, $\*l' \eq \*l'$.
    
    \item
    \emph{Symmetry:}
    We have
    $\*l' \eq \*l'' \iff \*l'' \eq \*l'$, since
    $T_{\*l',\*l''} = T_{\*l'',\*l'}$ and
    $\min\{l'_t, l''_t\} = \min\{l''_t, l'_t\}$.
    
    \item
    \emph{Transitivity:}
    Let $\*l' \eq \hat{\*l}$, $\hat{\*l} \eq \*l''$, and
    $t \notin T_{\*l',\*l''}$.
    From the definition of $T_{\*l',\*l''}$,
    it holds that either $l'_t \not= l''_t$ or $l'_t = l''_t \ge l_t$.
    As $l'_t = l''_t \ge l_t$ already implies $\min\{l'_t, l''_t\} \ge l_t$,
    we assume that $l'_t \not= l''_t$.
    Here, we have three cases:
    \begin{itemize}
      \item
      \mbox{\emph{Case 1:} ${l'_t \not= \hat{l}_t = l''_t}$.}
      $t \notin T_{\*l',\hat{\*l}}$ implies $l'_t \ge l_t$ and
      $l''_t = \hat{l}_t \ge l_t$.
      Therefore, $\min\{l'_t, l''_t\} \ge l_t$.
      
      \item
      \mbox{\emph{Case 2:} ${l'_t = \hat{l}_t \not= l''_t}$.}
      Analogously to the first case, we conclude $\min\{l'_t, l''_t\} \ge l_t$.
      
      \item
      \mbox{\emph{Case 3:} ${l'_t \not= \hat{l}_t \not= l''_t}$.}
      $t \notin T_{\*l',\hat{\*l}}$ implies $l'_t \ge l_t$ and
      $t \notin T_{\*l'',\hat{\*l}}$ implies $l''_t \ge l_t$.
      Hence, $\min\{l'_t, l''_t\} \ge l_t$.
    \end{itemize}
    Therefore, it holds that $\min\{l'_t, l''_t\} \ge l_t$
    for all $t \notin T_{\*l',\*l''}$, i.e., $\*l' \eq \*l''$.
  \end{itemize}
  This shows that $\eq$ is an equivalence relation.
\end{proof}

\addtocontents{lop}{\protect\iffalse}

\lemmaCombiTechniqueIdenticalValues*

\begin{proof}
  First, we note that $T_{\*l',\*l''} \not= \emptyset$.
  Otherwise, for $T_{\*l',\*l''} = \emptyset$,
  we have $\min\{l'_t, l''_t\} \ge l_t$
  for all $t = 1, \dotsc, d$, which implies $\*l' \ge \*l$, i.e.,
  $\fgset{\*l'} \ni \gp{\*l,\*i}$.
  This contradicts the fact that $\*l' \in L$, where $L$ is defined
  in \eqref{eq:combiTechniqueSpecialLevelSet}
  (which holds as our equivalence relation is only defined on $L$).
  Therefore, $T_{\*l',\*l''} \not= \emptyset$ must hold.
  Without loss of generality,
  we assume that $T_{\*l',\*l''} = \{1, \dotsc, m\}$
  for some $m \in \{1, \dotsc, d\}$.
  
  Let
  \begin{equation}
    S \ceq \gp{\*l,\*i} + \spn\{\stdbasis{1}, \dotsc, \stdbasis{m}\}
    = \{\gp{\*l,\*i} + \textstyle\sum_{t=1}^m c_t \stdbasis{t} \mid
    c_1, \dotsc, c_m \in \real\}
  \end{equation}
  be the $m$-dimensional affine subspace of $\real^d$
  through $\gp{\*l, \*i}$
  parallel to the dimensions $1, \dotsc, m$,
  where $\stdbasis{t}$ is the $t$-th standard basis vector.
  It holds that $S \cap \fgset{\*l'} = S \cap \fgset{\*l''}$ due to
  $l'_t = l''_t$ for $t \le m$.%
  \footnote{%
    In more detail:
    If we have an $\gp{\hat{\*l},\hat{\*i}} \in S \cap \fgset{\*l'}$,
    then $\fa{t \le m}{\hat{l}_t \le l'_t = l''_t}$ and
    $\fa{t > m}{\hat{l}_t = l_t \le l''_t}$, i.e.,
    $\hat{\*l} \le \*l''$ and therefore
    $\gp{\hat{\*l},\hat{\*i}} \in S \cap \fgset{\*l''}$.%
  }
  
  On this $m$-dimensional grid $S \cap \fgset{\*l'} = S \cap \fgset{\*l''}$,
  the full grid interpolants $\fgintp{\*l'}$ and $\fgintp{\*l''}$
  coincide, as both interpolate the function values given by
  the objective function $\objfun$:
  \begin{equation}
    \label{eq:proofCombiTechniqueIdenticalValues1}
    \restrictfcn{\fgintp{\*l'}}{S \cap \fgset{\*l'}}
    = \restrictfcn{\objfun}{S \cap \fgset{\*l'}}
    =  \restrictfcn{\fgintp{\*l''}}{S \cap \fgset{\*l'}}.
  \end{equation}
  However, this does not suffice to conclude
  $\fgintp{\*l'}(\gp{\*l,\*i}) = \fgintp{\*l''}(\gp{\*l,\*i})$,
  since $\gp{\*l,\*i} \notin \fgset{\*l'}$.
  
  To this end, we recall from \eqref{eq:interpFullGridMV} that
  \begin{equation}
    \fgintp{\*l'}
    = \sum_{\*i'=\*0}^{\*2^{\*l'}} \interpcoeff{\*l',\*i'}
    \basis{\*l',\*i'},\quad
    \interpcoeff{\*l',\*i'} \in \real.
  \end{equation}
  This implies that the $m$-variate restricted interpolant
  $\restrictfcn{\fgintp{\*l'}}{S \cap \clint{\*0, \*1}}$ can be written as
  \begin{subequations}
    \begin{align}
      (\restrictfcn{\fgintp{\*l'}}{S \cap \clint{\*0, \*1}})
      (\*x_{\range{1}{m}})
      &\mathrel{\righthphantom{=}{\ceq}}
      \sum_{\*i'_{\range{1}{m}}=\*0}^{\*2^{\*l'_{\range{1}{m}}}}
      \interpcoefftilde{\*l'_{\range{1}{m}},\*i'_{\range{1}{m}}}
      \basis{\*l'_{\range{1}{m}},\*i'_{\range{1}{m}}}(\*x_{\range{1}{m}}),\quad
      \*x_{\range{1}{m}} \in \clint{0, 1}^m,\\
      \interpcoefftilde{\*l'_{\range{1}{m}},\*i'_{\range{1}{m}}}
      &\ceq \sum_{\*i'_{\range{m+1}{d}}=\*0}^{\*2^{\*l'_{\range{m+1}{d}}}}
      \interpcoeff{\*l',\*i'}
      \basis{\*l'_{\range{m+1}{d}},\*i'_{\range{m+1}{d}}}%
      (\gp{\*l_{\range{m+1}{d}},\*i_{\range{m+1}{d}}})
    \end{align}
  \end{subequations}
  by factoring out tensor product factors corresponding to dimensions
  $m + 1, \dotsc, d$.
  The subscripts $\range{1}{m}$ and $\range{m+1}{d}$ denote
  the entries with respect to the dimensions $1, \dotsc, m$ and
  $m + 1, \dotsc, d$, respectively.
  As a result,
  both $\restrictfcn{\fgintp{\*l'}}{S \cap \clint{\*0, \*1}}$
  and, analogously,
  $\restrictfcn{\fgintp{\*l''}}{S \cap \clint{\*0, \*1}}$
  are interpolants of $\objfun$ in
  $\ns{\*l'_{\range{1}{m}}} = \ns{\*l''_{\range{1}{m}}}$.
  Due to \thmref{lemma:tensorProductLinearIndependence},
  it follows from \eqref{eq:proofCombiTechniqueIdenticalValues1}
  that they must be the same:
  \begin{equation}
    \restrictfcn{\fgintp{\*l'}}{S \cap \clint{\*0, \*1}}
    = \restrictfcn{\fgintp{\*l''}}{S \cap \clint{\*0, \*1}}.
  \end{equation}
  Consequently, $\fgintp{\*l'}(\gp{\*l,\*i}) = \fgintp{\*l''}(\gp{\*l,\*i})$
  as $\gp{\*l,\*i} \in S \cap \clint{\*0, \*1}$.
\end{proof}

\lemmaCombiTechniqueCharacterization*

\begin{proof}
  ``$\subset$'':
  Let $\*l' \in L_0$.
  We have to prove that
  $\fa{t \in T_{L_0}}{l'_t = l^\ast_t}$ and
  $\fa{t \notin T_{L_0}}{l'_t \ge l_t}$.
  The first statement is clear by the definition of $T_{L_0}$.
  Therefore, let $t \notin T_{L_0}$.
  By the definition of $T_{L_0}$,
  we have either
  $\ex{\hat{\*l} \in L_0}{l'_t \not= \hat{l}_t}$ or
  $\fa{\hat{\*l} \in L_0}{l'_t = \hat{l}_t \ge l_t}$.
  In the latter case, we obtain $l'_t \ge l_t$
  (e.g., by setting $\hat{\*l}$ to $\*l'$).
  In the former case, there is an $\hat{\*l} \in L_0$ such that
  $l'_t \not= \hat{l}_t$.
  Due to $\*l' \eq \hat{\*l}$
  (since $\*l'$ and $\hat{\*l}$ are both contained in the same
  equivalence class $L_0$) and $t \notin T_{\*l',\hat{\*l}}$,
  we have $\min\{l'_t, \hat{l}_t\} \ge l_t$.
  This implies $\fa{t \notin T_{L_0}}{l'_t \ge l_t}$, as desired.
  
  ``$\supset$'':
  Let $\*l' \in L$ such that
  $\fa{t \in T_{L_0}}{l'_t = l^\ast_t}$ and
  $\fa{t \notin T_{L_0}}{l'_t \ge l_t}$.
  Furthermore,
  let $\*l'' \in L_0$ be an arbitrary representative of $L_0$.
  We prove that $\*l' \eq \*l''$ (i.e., $\*l' \in L_0$).
  Note that $T_{L_0} \subset T_{\*l',\*l''}$,
  as $t \in T_{L_0}$ implies
  $l''_t = l^\ast_t < l_t$, which can be combined with $l'_t = l^\ast_t$
  to $l'_t = l''_t < l_t$, i.e., $t \in T_{\*l',\*l''}$.
  
  To prove the equivalence of $\*l'$ and $\*l''$,
  let $t \notin T_{\*l',\*l''}$,
  i.e., $t \notin T_{L_0}$.
  By assumption on $\*l'$, it holds $l'_t \ge l_t$.
  Hence, it remains to show that $l''_t \ge l_t$ as well.
  Again, by definition of $T_{L_0}$, we have either
  $\ex{\hat{\*l} \in L_0}{l''_t \not= \hat{l}_t}$ or
  $\fa{\hat{\*l} \in L_0}{l''_t = \hat{l}_t \ge l_t}$.
  In the second case, it holds $l''_t \ge l_t$.
  In the first case, there is an $\hat{\*l} \in L_0$ such that
  $l''_t \not= \hat{l}_t$.
  Due to $\*l'' \eq \hat{\*l}$
  (since $\*l''$ and $\hat{\*l}$ are both contained in the same
  equivalence class $L_0$) and $t \notin T_{\*l'',\hat{\*l}}$,
  we have $\min\{l''_t, \hat{l}_t\} \ge l_t$.
  In particular, $l''_t \ge l_t$.
  In total, we have $\min\{l'_t, l''_t\} \ge l_t$ for all
  $t \notin T_{\*l',\*l''}$, proving that $\*l'$ and $\*l''$ are equivalent,
  as asserted.
\end{proof}

\propCombiTechniqueZero*

\begin{proof}
  \Cref{lemma:combiTechniqueIdenticalValues}
  implies that the summands $\fgintp{\*l'}(\gp{\*l,\*i})$
  corresponding to levels $\*l'$ of the same equivalence class
  $L_0 \in \eqclasses{L}{\eq}$ are identical.
  Let $f_{L_0}$ denote the common function value.
  The sum in the \lhs of the assertion can be reordered to combine
  levels of the equivalence classes $L_0 \in \eqclasses{L}{\eq}$:
  \begin{subequations}
    \begin{align}
      &\sum_{q=0}^{d-1} (-1)^q \binom{d-1}{q} \cdot
      \sum_{\substack{\normone{\*l'} = n - q\\\fgset{\*l'} \notni \gp{\*l,\*i}}}
      \fgintp{\*l'}(\gp{\*l,\*i})\\
      \label{eq:proofCombiTechniqueIdenticalValues2}
      &= \sum_{L_0 \in \eqclasses{L}{\eq}} f_{L_0} \sum_{q=0}^{d-1}
      (-1)^q \binom{d-1}{q} \cdot
      \setsize{\{\*l' \in L_0 \mid \normone{\*l'} = n - q\}}.
    \end{align}
  \end{subequations}
  It now suffices to show that the inner sum vanishes
  for every equivalence class $L_0 \in \eqclasses{L}{\eq}$.
  
  To this end, we have to calculate
  $\setsize{\{\*l' \in L_0 \mid \normone{\*l'} = n - q\}}$
  for a fixed equivalence class $L_0$.
  Without loss of generality,
  let $T_{L_0} = \{1, \dotsc, m\}$ in the notation of
  \thmref{lemma:combiTechniqueCharacterization}
  with $1 \le m \le d$.
  Note that the case $m = 0$ is impossible:
  Otherwise, $T_{L_0} = \emptyset$ implies
  $\fa{\*l' \in L_0}{\*l' \ge \*l}$ by
  \cref{lemma:combiTechniqueCharacterization}, and
  as equivalence classes are non-empty,
  there is at least one $\*l' \in L_0$ with $\*l' \ge \*l$.
  However, this is equivalent to $\fgset{\*l'} \ni \gp{\*l,\*i}$,
  which contradicts $\*l' \in L$.
  Hence, we have $m > 0$.
  
  To enumerate all levels $\*l' \in L_0$ with $\normone{\*l'} = n - q$,
  we exploit the characterization of $L_0$ of
  \cref{lemma:combiTechniqueCharacterization}.
  For notational convenience, we define the vector
  \begin{equation}
    \hat{\*l}
    \ceq (l^\ast_1, \dotsc, l^\ast_m,\; l_{m+1}, \dotsc, l_d),
  \end{equation}
  where $\*l^\ast$ is given as in \cref{lemma:combiTechniqueCharacterization}.
  We show that $\*a \ceq (l'_t - l_t)_{t = m+1, \dotsc, d}$
  constitutes a bijection between
  \begin{equation}
    \{\*l' \in L_0 \mid \normone{\*l'} = n - q\}
    \quad\text{and}\quad
    \{\*a \in \natz^{d-m} \mid \normone{\*a} = n - q - \normone{\hat{\*l}}\}
    \colon
  \end{equation}
  \begin{itemize}
    \item
    Let $\*l' \in L_0$ with $\normone{\*l'} = n - q$.
    Then, $\fa{t=m+1,\dotsc,d}{l'_t - l_t \ge 0}$
    (by \cref{lemma:combiTechniqueCharacterization}), i.e.,
    $\*a \in \natz^{d-m}$, and
    \begin{equation}
      \normone{\*a}
      = \sum_{t=m+1}^d (l'_t - l_t)
      = \paren*{\normone{\*l'} - \sum_{t=1}^m l'_t} -
      \sum_{t=m+1}^d \hat{l}_t
      = n - q - \normone{\hat{\*l}}.
    \end{equation}
    
    \item
    Conversely, let $\*a \in \natz^{d-m}$ with
    $\normone{\*a} = n - q - \normone{\hat{\*l}}$.
    If we define $\*l'$ as
    \begin{equation}
      \*l'
      = (l^\ast_1, \dotsc, l^\ast_m,\;
      a_1 + l_{m+1}, \dotsc, a_{d-m} + l_d),
    \end{equation}
    then $\fa{t=1,\dotsc,m}{l'_t = l^\ast_t < l_t}$ and
    $\fa{t=m+1,\dotsc,d}{l'_t \ge l_t}$.
    By \cref{lemma:combiTechniqueCharacterization},
    we obtain $\*l' \in L_0$ and
    \begin{equation}
      \normone{\*l'}
      = \normone{\hat{\*l}} + \normone{\*a}
      = n - q.
    \end{equation}
  \end{itemize}
  This bijection implies that
  \begin{subequations}
    \begin{align}
      \setsize{\{\*l' \in L_0 \mid \normone{\*l'} = n - q\}}
      &= \setsize{
        \{\*a \in \natz^{d-m} \mid \normone{\*a} = n - q - \normone{\hat{\*l}}\}
      }.\\
      \intertext{%
        This is the number of weak decompositions of
        $n - q - \normone{\hat{\*l}}$ into $d - m$ parts:%
      }
      \cdots
      &= \binom{n - q - \normone{\hat{\*l}} + d - m - 1}{d - m - 1},
    \end{align}
  \end{subequations}
  see Theorem 2.2 of \cite{Bona15Introduction}.
  We insert this quantity into the inner sum of
  \eqref{eq:proofCombiTechniqueIdenticalValues2}:
  \begin{subequations}
    \label{eq:proofCombiTechniqueIdenticalValues3}
    \begin{align}
      &\sum_{q=0}^{d-1}
      (-1)^q \binom{d-1}{q} \cdot
      \setsize{\{\*l' \in L_0 \mid \normone{\*l'} = n - q\}}\\
      &= \sum_{q=0}^{d-1} (-1)^q \binom{d-1}{q} \cdot
      \binom{n - q - \normone{\hat{\*l}} + d - m - 1}{d - m - 1}.
    \end{align}
  \end{subequations}
  Again, we apply \thmref{lemma:inclusionExclusionCountingLemma}
  with the values $a \ceq d - 1$,
  $r \ceq n - \normone{\hat{\*l}} + d - m - 1$, and
  $s \ceq d - m - 1$ to infer that as claimed,
  \eqref{eq:proofCombiTechniqueIdenticalValues3} is equal to
  \begin{equation}
  \label{eq:proofCombiTechniqueIdenticalValues4}
    \binom{n - \normone{\hat{\*l}} - m}{-m}
    = 0
  \end{equation}
  by the convention for binomial coefficients in \cref{def:binomialCoefficient}
  as $-m < 0$.
  
  Note that for the calculation in
  \cref{%
    eq:proofCombiTechniqueIdenticalValues3,%
    eq:proofCombiTechniqueIdenticalValues4%
  }
  to be correct,
  we have to ensure that $n - \normone{\hat{\*l}} - m \ge 0$;
  otherwise, the binomial coefficients would not be well-defined.
  This is a direct consequence of the fact that $l^\ast_t < l_t$
  for all $t = 1, \dotsc, m$
  (see \cref{lemma:combiTechniqueCharacterization}) as
  \begin{equation}
    n - \normone{\hat{\*l}} - m
    = n - \sum_{t=1}^m l^\ast_t - \sum_{\mathclap{t=m+1}}^d l_t - m
    \ge n - \sum_{t=1}^m (l_t - 1) - \sum_{\mathclap{t=m+1}}^d l_t - m
    = n - \normone{\*l}
    \ge 0,
  \end{equation}
  where we have used $l^\ast_t \le l_t - 1$ for $t = 1, \dotsc, m$
  and $\normone{\*l} \le n$.
\end{proof}

\fillsubsectionornament
\subsection{Correctness Proof of the Method of Residual Interpolation}
\label{sec:a132proofResidualInterpolation}

\propInvariantResidualInterpolation*

\begin{proof}
  \newcommand*{\eqwithref}[1]{%
    \quad\;%
    \if\relax\detokenize{#1}\relax%
      \mathclap{=}%
    \else%
      \mathclap{\overset{\eqref{eq:propInvariantResidualInterpolation#1}}{=}}%
    \fi%
    \quad\;%
  }%
  We prove the assertion by induction over $j = 1, \dotsc, m$.
  We will need the following two equations
  that directly follow from the algorithm
  (\cref{%
    line:algResidualInterpolation3,%
    line:algResidualInterpolation2,%
    line:algResidualInterpolation4%
  }, respectively):%
  \begin{subequations}
    \begin{alignat}{3}
      \label{eq:propInvariantResidualInterpolation5}
      r_{\*l^{(j)}}^{(j-1)}(\gp{\*l,\*i})
      &= r^{(j-1)}(\gp{\*l,\*i}),\quad
      &&\*l \le \*l^{(j)},\;\;
      &&\*i \in \hiset{\*l},\\
      \label{eq:propInvariantResidualInterpolation6}
      r^{(j)}(\gp{\*l,\*i})
      &= r^{(j-1)}(\gp{\*l,\*i}) - r_{\*l^{(j)}}^{(j-1)}(\gp{\*l,\*i}),\quad
      &&\*l \in L,\;\;
      &&\*i \in \hiset{\*l}.
    \end{alignat}
  \end{subequations}
  
  \vspace*{-2em}
  \pagebreak
  
  \noindent
  \textbf{Induction base case:}
  For $j = 1$, there is nothing to show for
  \eqref{eq:propInvariantResidualInterpolation1}.
  \Cref{eq:propInvariantResidualInterpolation2}
  can be proven as follows:
  \begin{equation}
    r^{(1)}(\gp{\*l,\*i})
    \eqwithref{6}
    r^{(0)}(\gp{\*l,\*i}) - r_{\*l^{(1)}}^{(0)}(\gp{\*l,\*i})
    \eqwithref{5} 0,\qquad
    \*l \le \*l^{(1)},\; \*i \in \hiset{\*l}.
  \end{equation}
  \Cref{eq:propInvariantResidualInterpolation3}
  holds as $r_{\*l^{(1)}}^{(0)} = f^{\sparse,(1)}$
  (by \cref{line:algResidualInterpolation2} in
  \cref{alg:residualInterpolation}) and, therefore,
  \begin{align}
    r^{(1)}(\gp{\*l,\*i})
    \eqwithref{6}
    r^{(0)}(\gp{\*l,\*i}) - r_{\*l^{(1)}}^{(0)}(\gp{\*l,\*i})
    = \fcnval{\*l,\*i} - f^{\sparse,(1)}(\gp{\*l,\*i}),\quad
    \*l \in L,\; \*i \in \hiset{\*l}.
  \end{align}
  
  \noindent
  \textbf{Induction step case:}
  We show the three statements for the induction step $j \to (j + 1)$.
  \begin{itemize}
    \item
    \emph{Showing \eqref{eq:propInvariantResidualInterpolation1} for $j + 1$:}
    Let $j' = 1, \dotsc, j$, $\*l \le \*l^{(j')}$,
    and $\*i \in \hiset{\*l}$.
    Due to the ordering of the levels $\*l^{(1)}, \dotsc, \*l^{(m)}$,
    we can conclude from $j + 1 > j'$ that
    $\normone{\*l^{(j+1)}} \le \normone{\*l^{(j')}}$.
    This implies that there must be a $t' \in \{1, \dotsc, d\}$
    such that $l_{t'}^{(j+1)} \le l_{t'}^{(j')}$.
    Let $S$ be the line in $\real^d$ defined by
    \begin{equation}
      S
      \ceq \gp{\*l,\*i} + \spn\{\stdbasis{t'}\},
    \end{equation}
    where $\stdbasis{t'}$ is the $t'$-th standard basis vector.
    It holds that $S \cap \fgset{\*l^{(j+1)}} \subset \fgset{\*l^{(j')}}$.
    To show this, let $\gp{\*l',\*i'} \in S \cap \fgset{\*l^{(j+1)}}$
    be arbitrary (with $\*i' \in \hiset{\*l'}$).
    Then, $\fa{t \not= t'}{l'_t = l_t \le l_t^{(j')}}$
    (due to $\gp{\*l',\*i'} \in S$) and
    $l'_{t'} \le l_{t'}^{(j+1)} \le l_{t'}^{(j')}$
    (due to $\gp{\*l',\*i'} \in \fgset{\*l^{(j+1)}}$).
    This means that $\*l' \le \*l^{(j')}$, which implies that
    $\gp{\*l',\*i'} \in \fgset{\*l^{(j')}}$.
    As $\gp{\*l',\*i'}$ is arbitrary,
    this shows $S \cap \fgset{\*l^{(j+1)}} \subset \fgset{\*l^{(j')}}$.
    
    Thus, we infer
    \begin{equation}
      \label{eq:proofPropInvariantResidualInterpolation2}
      r_{\*l^{(j+1)}}^{(j)}(\gp{\*l',\*i'})
      \eqwithref{5}
      r^{(j)}(\gp{\*l',\*i'})
      \eqwithref{2}
      0,\quad
      \gp{\*l',\*i'} \in S \cap \fgset{\*l^{(j+1)}}
      \subset \fgset{\*l^{(j')}},\;
      i' \in \hiset{\*l'},
    \end{equation}
    with the induction hypothesis
    \eqref{eq:propInvariantResidualInterpolation2} for $j$.
    Unfortunately, this does not suffice to directly conclude that
    $r_{\*l^{(j+1)}}^{(j)}(\gp{\*l,\*i}) = 0$ as
    $\gp{\*l,\*i}$ is in general not contained in $\fgset{\*l^{(j+1)}}$.
    
    As in the proof of \cref{lemma:combiTechniqueIdenticalValues},
    we exploit the tensor product nature of the basis functions and
    restrict $r_{\*l^{(j+1)}}^{(j)}$ to $S \cap \clint{0, 1}$:
    \begin{subequations}
      \begin{align}
        (\restrictfcn{r_{\*l^{(j+1)}}^{(j)}}{S \cap \clint{0, 1}})(x_{t'})
        &\mathrel{\righthphantom{=}{\ceq}}
        \sum_{l'_{t'}=0}^{l_{t'}^{(j+1)}}
        \sum_{i'_{t'} \in \hiset{l'_{t'}}}
        \surplustilde[(j+1)]{l'_{t'},i'_{t'}}
        \basis{l'_{t'},i'_{t'}}(x_{t'}),\quad
        x_{t'} \in \clint{0, 1},\\
        \surplustilde[(j+1)]{l'_{t'},i'_{t'}}
        &\ceq \sum_{\*l'_{-t'}=\*0}^{\*l^{(j+1)}_{-t'}}
        \sum_{\*i'_{-t'} \in \hiset{\*l'_{-t'}}}
        \surplus[(j+1)]{\*l',\*i'}
        \basis{\*l'_{-t'},\*i'_{-t'}}(\gp{\*l_{-t'},\*i_{-t'}}),
      \end{align}
    \end{subequations}
    where the subscript $-t'$ indicates all entries but the $t'$-th.
    As a consequence, this shows that
    $\restrictfcn{r_{\*l^{(j+1)}}^{(j)}}{S \cap \clint{0, 1}} \in
    \ns{l_{t'}^{(j+1)}}$ is an interpolant of the zero function
    (by \eqref{eq:proofPropInvariantResidualInterpolation2}).
    Due to the linear independence of the univariate basis functions,
    we conclude
    \begin{equation}
      \restrictfcn{r_{\*l^{(j+1)}}^{(j)}}{S \cap \clint{0, 1}}
      \equiv 0.
    \end{equation}
    Consequently, we obtain
    $r_{\*l^{(j+1)}}^{(j)}(\gp{\*l,\*i}) = 0$
    as $\gp{\*l,\*i} \in S \cap \clint{\*0, \*1}$.
    
    \item
    \emph{Showing \eqref{eq:propInvariantResidualInterpolation2} for $j + 1$:}
    Let $j' = 1, \dotsc, j + 1$, $\*l \le \*l^{(j')}$,
    and $\*i \in \hiset{\*l}$.
    For the case $j' \le j$, we obtain
    \begin{equation}
      \label{eq:proofPropInvariantResidualInterpolation1}
      r^{(j+1)}(\gp{\*l,\*i})
      \eqwithref{6}
      r^{(j)}(\gp{\*l,\*i}) - r_{\*l^{(j+1)}}^{(j)}(\gp{\*l,\*i})
      = 0
    \end{equation}
    due to $r^{(j)}(\gp{\*l,\*i}) = 0$ by induction hypothesis
    (\cref{eq:propInvariantResidualInterpolation2}) and
    $r_{\*l^{(j+1)}}^{(j)}(\gp{\*l,\*i}) = 0$ as shown above
    (\cref{eq:propInvariantResidualInterpolation1} for $j + 1$).
    
    For the case $j' = j + 1$,
    \cref{eq:proofPropInvariantResidualInterpolation1}
    still holds as the difference between
    $r^{(j)}(\gp{\*l,\*i})$ and $r_{\*l^{(j+1)}}^{(j)}(\gp{\*l,\*i})$
    vanishes due to \eqref{eq:propInvariantResidualInterpolation5}
    for $j + 1$
    (here, we need $\*l \le \*l^{(j+1)}$).
    
    \item
    \emph{Showing \eqref{eq:propInvariantResidualInterpolation3} for $j + 1$:}
    Let $\*l \in L$ and $\*i \in \hiset{\*l}$.
    Then,
    \begin{subequations}
      \begin{align}
        r^{(j+1)}(\gp{\*l,\*i})
        &\eqwithref{6}
        r^{(j)}(\gp{\*l,\*i}) - r_{\*l^{(j+1)}}^{(j)}(\gp{\*l,\*i}).\\
        \intertext{%
          The first term can replaced with the induction hypothesis
          (\eqref{eq:propInvariantResidualInterpolation3} for $j$).
          For the second term, note that
          $r_{\*l^{(j+1)}}^{(j)}
          = \sum_{\*l' \in \levelset} \sum_{\*i' \in \hiset{\*l'}}
          \surplus[(j+1)]{\*l',\*i'} \basis{\*l',\*i'}
          = f^{\sparse,(j+1)} - f^{\sparse,(j)}$ by definition.
          Hence, we obtain as desired%
        }
        \cdots
        &\eqwithref{} (\fcnval{\*l,\*i} - f^{\sparse,(j)}(\gp{\*l,\*i})) -
        (f^{\sparse,(j+1)}(\gp{\*l,\*i}) - f^{\sparse,(j)}(\gp{\*l,\*i}))\\
        &\eqwithref{} \fcnval{\*l,\*i} - f^{\sparse,(j+1)}(\gp{\*l,\*i}).
      \end{align}
    \end{subequations}
  \end{itemize}
  This shows the validity of the statements in
  \eqref{eq:propInvariantResidualInterpolationStatements}
  for $j + 1$.
\end{proof}

\subsection{Correctness Proof of Hierarchization with Breadth-First Search}
\label{sec:a133proofBFS}

\propInvariantBFS*

\begin{proof}
  We start with two observations:
  
  \begin{itemize}
    \item
    First, due to the breadth-first search nature of \cref{alg:BFS} and the
    hierarchical relation \eqref{eq:directAncestor},
    all grid points with level sum $< q$ are \pop{}ped before
    the first point with level sum $\ge q$ is \pop{}ped.
    
    \item
    Second, after \pop{}ping all grid points with
    level sum $< q$, the output values of the grid points with
    level sum $\le q$ remain unchanged for the rest of the algorithm:
    If \cref{line:algBFS2} of the algorithm updates the output value of a point
    $(\*l, \*i)$ in the iteration of $(\*l', \*i') \in \liset$ with
    $\normone{\*l'} \ge q$, then \cref{line:algBFS1} implies
    $\*l \ge \*l'$ and thus, $\normone{\*l} \ge \normone{\*l'} \ge q$.
    However, $\normone{\*l} = q$ is not possible as
    this would imply that $\normone{\*l} = \normone{\*l'}
    \implies (\*l, \*i) = (\*l', \*i')$ by \cref{line:algBFS1},
    but $(\*l, \*i) = (\*l', \*i')$ is explicitly excluded in the
    \texttt{\algorithmicfor} loop of \cref{line:algBFS1}.
    Therefore, we must have $\normone{\*l} > q$.
    Hence, if a point $(\*l, \*i)$ with level sum $\ge q$ has been \pop{}ped,
    only surpluses of points with level sum $> q$ may be updated.
  \end{itemize}
  
  \noindent
  Now, we prove the asserted claim by induction over $q$.
  
  \noindent
  \textbf{Induction base case:}
  For $q = 0$, \cref{alg:BFS} sets $\linout{\*l,\*i}$ to
  $\fcnval{\*l,\*i}$ in \cref{line:algBFS3}.
  As the sum in \eqref{eq:propInvariantBFS} is empty,
  the claim is correct for $q = 0$.
  
  \noindent
  \textbf{Induction step case:}
  Let $\linout[(q)]{\*l',\*i'}$ and
  $\linout[(q+1)]{\*l',\*i'}$
  be the surpluses after \pop{}ping all
  grid points with level sum $< q$ and $< q + 1$, respectively.
  We show the induction step $q \to (q + 1)$, i.e.,
  we assume that the assertion is true for $q$
  and prove that after \pop{}ping all grid points with level sum $< q + 1$,
  it holds
  \begin{equation}
    \label{eq:proofPropInvariantBFS1}
    \linout[(q+1)]{\*l,\*i}
    = \fcnval{\*l,\*i} -
    \largesum{\normone{\*l'} < q+1} \linout[(q+1)]{\*l',\*i'}
    \fundbasis{\*l',\*i'}(\gp{\*l,\*i}),\quad
    (\*l, \*i) \in \liset,\;\;
    \normone{\*l} = q+1.
  \end{equation}
  Therefore, let $(\*l, \*i) \in \liset$ fulfill $\normone{\*l} = q+1$.
  The update in \cref{line:algBFS2} can safely be applied
  with all grid points $(\*l', \*i')$ with level sum $q$.
  The grid points $(\*l', \*i')$ that do not satisfy the relation in the set in
  \cref{line:algBFS1} do not contribute as
  $\fundbasis{\*l',\*i'}(\gp{\*l,\*i}) = 0$
  due to the necessary condition \eqref{eq:fundamentalPropertyImplicationMV}.
  By summing all updates from \cref{line:algBFS2}, we obtain
  \begin{subequations}
    \begin{align}
      \linout[(q+1)]{\*l,\*i}
      &= \linout[(q)]{\*l,\*i} -
      \largesum{\normone{\*l'} = q} \linout[(q)]{\*l',\*i'}
      \fundbasis{\*l',\*i'}(\gp{\*l,\*i}).\\
      \intertext{%
        After inserting the induction hypothesis
        for the first $\linout[(q)]{\*l,\*i}$, we have%
      }
      \cdots
      &= \paren*{
        \fcnval{\*l,\*i} -
        \largesum{\normone{\*l'} < q} \linout[(q)]{\*l',\*i'}
        \fundbasis{\*l',\*i'}(\gp{\*l,\*i})
      } -
      \largesum{\normone{\*l'} = q} \linout[(q)]{\*l',\*i'}
      \fundbasis{\*l',\*i'}(\gp{\*l,\*i})\\
      &= \hphantom{\Biggl(}\fcnval{\*l,\*i} -
      \largesum{\normone{\*l'} < q + 1} \linout[(q)]{\*l',\*i'}
      \fundbasis{\*l',\*i'}(\gp{\*l,\*i}).
    \end{align}
  \end{subequations}
  As noted above, we have
  $\fa{(\*l', \*i'),\, \normone{\*l'} < q + 1}{
    \linout[(q)]{\*l',\*i'} = \linout[(q+1)]{\*l',\*i'}
  }$
  (the values of points with level sum $< q + 1$
  do not change after \pop{}ping all points with level sum $< q$).
  This shows the induction claim \eqref{eq:proofPropInvariantBFS1}.
\end{proof}

\disableornamentsfornextheadingtrue
\subsection{%
  Proof for the Correctness of the Unidirectional Principle on
  Spatially Adaptive Sparse Grids%
}
\label{sec:a134proofCorrectnessUnidirectionalPrincipleSASG}

\lemmaChainExistenceSufficient*

\begin{proof}
  We prove the assertion by induction over $j = 0, \dotsc, d$.
  
  For $j = 0$, the operator $\upop{\emptyset}$ is by definition
  \eqref{eq:upopProduct} the identity operator $\idop$.
  The assumption $(\upop{\emptyset})_{\*k'',\*k'} \not= 0$
  implies that $\*k' = \*k''$, since the identity matrix is diagonal.
  Therefore, the chain $(\chain{0})$ from $\*k'$ to $\*k''$ is given by
  $\chain{0} = \*k'$, which is contained in $\liset$.
  
  For the induction step $j \to (j+1)$, we split \cref{eq:upopProduct},
  i.e.,
  \begin{equation}
    \upop{t_1,\dotsc,t_{j+1}}
    = \upop{t_{j+1}} \upop{t_j} \dotsm \upop{t_1}
    = \upop{t_{j+1}} \upop{t_1,\dotsc,t_j},
  \end{equation}
  and infer by assumption
  \begin{equation}
    \label{eq:proofLemmaChainExistenceSufficient1}
    0
    \not= (\upop{t_1,\dotsc,t_{j+1}})_{\*k'',\*k'}
    = \sum_{\*k \in \liset} (\upop{t_{j+1}})_{\*k'',\*k}
    (\upop{t_1,\dotsc,t_j})_{\*k,\*k'}.
  \end{equation}
  Consequently,
  there is at least one summation index $\*k$
  for which both factors do not vanish.
  The first factor $(\upop{t_{j+1}})_{\*k'',\*k}$ can by definition
  only be non-zero if $\*k \samepole{t_{j+1}} \*k''$.
  The second factor $(\upop{t_1,\dotsc,t_j})_{\*k,\*k'}$ being
  non-zero implies that by induction hypothesis,
  $\liset$ contains the chain $(\chain{0}, \dotsc, \chain{j})$
  from $\*k'$ to $\*k$ with respect to $(t_1, \dotsc, t_j)$.
  The combination of both statements leads to
  the chain $(\chain{0}, \dotsc, \chain{j}, \chain{j+1})$ from $\*k'$
  to $\*k''$ with respect to $(t_1, \dotsc, t_{j+1})$.
  All points of the chain are contained in $\liset$.
\end{proof}

\lemmaChainExistenceNecessary*

\begin{proof}
  Again, we prove the claim by induction over $j = 0, \dotsc, d$.
  
  For $j = 0$, the operator $\upop{\emptyset}$ is the identity operator.
  Therefore, the \lhs of \eqref{eq:lemmaChainExistenceNecessary} is one
  (due to $\chain{0} = \*k'$).
  The \rhs is by convention also one, as it is an empty product.
  
  \setlength{\abovedisplayskip}{9pt}%
  \setlength{\belowdisplayskip}{9pt}%
  For the induction step $j \to (j+1)$, we consider again
  \begin{equation}
    \label{eq:proofLemmaChainExistenceNecessary1}
    (\upop{t_1,\dotsc,t_{j+1}})_{\chain{j+1},\*k'}
    = \sum_{\*k \in \liset} (\upop{t_{j+1}})_{\chain{j+1},\*k}
    (\upop{t_1,\dotsc,t_j})_{\*k,\*k'}
  \end{equation}%
  similar to \eqref{eq:proofLemmaChainExistenceSufficient1}.
  Recall that
  \begin{subequations}
    \begin{alignat}{4}
      \chainuv{j}_t
      &= k''_t
      &&\;\;\text{and}\;\;
      &\chainuv{j+1}_t
      &= k''_t
      &&\quad\text{for } t \in \{t_1, \dotsc, t_j\},\\
      \label{eq:prooflemmaChainExistenceNecessary}
      \chainuv{j}_t
      &= k'_t
      &&\;\;\text{and}\;\;
      &\chainuv{j+1}_t
      &= k''_t
      &&\quad\text{for } t = t_{j+1},\\
      \chainuv{j}_t
      &= k'_t
      &&\;\;\text{and}\;\;
      &\chainuv{j+1}_t
      &= k'_t
      &&\quad\text{for } t \notin \{t_1, \dotsc, t_j, t_{j+1}\}.
    \end{alignat}
  \end{subequations}
  We now argue that all summands of
  \eqref{eq:proofLemmaChainExistenceNecessary1} vanish,
  except the summand with index $\chain{j}$.
  There are two cases for the summation index $\*k$,
  if we assume $\*k \not= \chain{j}$:
  \begin{itemize}
    \item
    If there is a $t \in \{t_1, \dotsc, t_j\}$ with $k_t \not= k''_t$,
    then we have $k_t \not= k''_t = \chainuv{j+1}_t$.
    Consequently, $\chain{j+1} \not\samepole{t_{j+1}} \*k$ and
    $(\upop{t_{j+1}})_{\chain{j+1},\*k} = 0$ due to \eqref{eq:upopEntries},
    i.e., the first factor of the $\*k$-th summand in
    \eqref{eq:proofLemmaChainExistenceNecessary1} vanishes.
    
    \item
    If there is a $t \notin \{t_1, \dotsc, t_j\}$ with $k_t \not= k'_t$,
    then the second factor $(\upop{t_1,\dotsc,t_j})_{\*k,\*k'}$
    of the $\*k$-th summand in
    \eqref{eq:proofLemmaChainExistenceNecessary1} vanishes.
    Indeed, if we assume the contrary,
    then \cref{lemma:chainExistenceSufficient} implies that there is
    a chain from $\*k'$ to $\*k$ with respect to $(t_1, \dotsc, t_j)$.
    However, by definition of the chain, this means that
    $\*k'$ and $\*k$ coincide in all other dimensions
    (which are not in $\{t_1, \dotsc, t_j\}$).
    This contradicts $k_t \not= k'_t$ and therefore
    $(\upop{t_1,\dotsc,t_j})_{\*k,\*k'}$ must vanish.
  \end{itemize}
  We infer that only the summand $\*k = \chain{j}$ remains in
  \eqref{eq:proofLemmaChainExistenceNecessary1}:
  \begin{equation}
    (\upop{t_1,\dotsc,t_{j+1}})_{\chain{j+1},\*k'}
    = (\upop{t_{j+1}})_{\chain{j+1},\chain{j}}
    (\upop{t_1,\dotsc,t_j})_{\chain{j},\*k'}.
  \end{equation}
  The first factor equals
  \begin{equation}
    (\upop{t_{j+1}})_{\chain{j+1},\chain{j}}
    = (\upopuv{t_{j+1}}{\eqclass{\chain{j+1}}{\samepole{t_{j+1}}}})%
    _{\chainuv{j+1}_{t_{j+1}},\chainuv{j}_{t_{j+1}}}
    = (\upopuv{t_{j+1}}{\eqclass{\chain{j+1}}{\samepole{t_{j+1}}}})%
    _{k''_{t_{j+1}},k'_{t_{j+1}}}
  \end{equation}
  by \cref{eq:upopEntries,eq:prooflemmaChainExistenceNecessary}
  (and due to $\chain{j} \samepole{t_{j+1}} \chain{j+1}$).
  The second factor equals
  \begin{equation}
    (\upop{t_1,\dotsc,t_j})_{\chain{j},\*k'}
    =
    (\upopuv{t_1}{\eqclass{\chain{1}}{\samepole{t_1}}})_{k''_{t_1},k'_{t_1}}
    \dotsm
    (\upopuv{t_j}{\eqclass{\chain{j}}{\samepole{t_j}}})_{k''_{t_j},k'_{t_j}}
  \end{equation}
  by induction hypothesis.
  Hence, the product of both factors is
  \begin{equation}
    (\upop{t_1,\dotsc,t_{j+1}})_{\chain{j+1},\*k'}
    =
    (\upopuv{t_1}{\eqclass{\chain{1}}{\samepole{t_1}}})_{k''_{t_1},k'_{t_1}}
    \dotsm
    (\upopuv{t_{j+1}}{\eqclass{\chain{j+1}}{\samepole{t_{j+1}}}})%
    _{k''_{t_{j+1}},k'_{t_{j+1}}}.
    \vspace*{-2.5em}
  \end{equation}
\end{proof}

\vspace*{1em}
\pagebreak

\propCorrectnessUPCharacterization*

\begin{proof}
  ``$\implies$'':
  Let the \up be correct for $\linop$ and $(t_1, \dotsc, t_d)$
  and $\*k', \*k'' \in \liset$ with $(\linop)_{\*k'',\*k'} \not= 0$.
  Then, we obtain
  \begin{equation}
    (\upop{t_1,\dotsc,t_d})_{\*k'',\*k'}
    = (\linop)_{\*k'',\*k'}
    \not= 0.
  \end{equation}
  By \cref{lemma:chainExistenceSufficient},
  this implies that $\liset$ contains the chain from $\*k'$ to $\*k''$
  with respect to $(t_1, \dotsc, t_d)$.
  
  ``$\impliedby$'':
  For the converse direction, we assume that there are chains
  from $\*k'$ to $\*k''$ with respect to $(t_1, \dotsc, t_d)$
  for all $\*k', \*k'' \in \liset$ with $(\linop)_{\*k'',\*k'} \not= 0$.
  Let $\*k', \*k'' \in \liset$ be arbitrary.
  There are two cases:
  \begin{itemize}
    \item
    ${(\linop)_{\*k'',\*k'} \not= 0}$:
    By assumption, $\liset$ contains the chain from $\*k'$ to $\*k''$
    with respect to $(t_1, \dotsc, t_d)$.
    We apply \cref{lemma:chainExistenceNecessary} with $j = d$
    to infer
    \begin{equation}
      (\upop{t_1,\dotsc,t_d})_{\*k'',\*k'}
      =
      (\upopuv{t_1}{\eqclass{\chain{1}}{\samepole{t_1}}})_{k''_{t_1},k'_{t_1}}
      \dotsm
      (\upopuv{t_d}{\eqclass{\chain{d}}{\samepole{t_d}}})_{k''_{t_d},k'_{t_d}}
      = (\linop)_{\*k'',\*k'}
    \end{equation}
    by the assumption \eqref{eq:tensorProductOperator} on
    the tensor product structure of $\linop$.
    
    \item
    ${(\linop)_{\*k'',\*k'} = 0}$:
    In this case, $(\upop{t_1,\dotsc,t_d})_{\*k'',\*k'}$ must vanish as well.
    Indeed, if we assume the contrary
    $(\upop{t_1,\dotsc,t_d})_{\*k'',\*k'} \not= 0$,
    then we can apply \cref{lemma:chainExistenceSufficient}
    to obtain that $\liset$ contains the chain from $\*k'$ to $\*k''$
    with respect to $(t_1, \dotsc, t_d)$.
    We conclude with \cref{lemma:chainExistenceNecessary} as in the first case
    that $(\upop{t_1,\dotsc,t_d})_{\*k'',\*k'} = (\linop)_{\*k'',\*k'} = 0$,
    which is a contradiction.
  \end{itemize}
  In any case, we obtain
  $(\upop{t_1,\dotsc,t_d})_{\*k'',\*k'} = (\linop)_{\*k'',\*k'}$,
  from which follows the correctness of the \up\punctfix{,}
  as $\*k'$ and $\*k''$ are arbitrary.
\end{proof}

\breakpagebeforenextheadingtrue
\subsection{Correctness Proof of Hermite Hierarchization}
\label{sec:a135proofHermiteHierarchization}

\propInvariantHermiteHierarchization*

\begin{proof}
  We prove the assertion by induction over $l = 0, \dotsc, n$.
  
  For the induction base case $l = 0$ and $i \in \{0, 1\}$, we have
  \begin{equation}
    \sum_{i'=0}^1
    \linout{0,i'} \deriv[q]{x}{\bspl[\wfs]{0,i'}{p}}(\gp{0,i})
    = \kronecker{q}{0} \cdot \fcnval{0,i} +
    \kronecker{q}{1} \cdot (\fcnval{0,1} - \fcnval{0,0})
    = \deriv[q]{x}{\fgintp{0}}(\gp{0,i})
  \end{equation}
  for $q = 0, \dotsc, \frac{p-1}{2}$ by
  \cref{%
    line:algHermiteHierarchization1,%
    line:algHermiteHierarchization3%
  }
  of \cref{alg:hermiteHierarchization}.
  
  For the induction step case $(l-1) \to l$,
  it suffices to show that
  \begin{equation}
    \label{eq:proofPropInvariantHermiteHierarchization1}
    \deriv[q]{x}{\fgintp{l-1}}(\gp{l,i})
    \overset{!}{=} \sum_{l'=0}^{l-1} \sum_{i' \in \hiset{l'}}
    \linout{l',i'} \deriv[q]{x}{\bspl[\wfs]{l',i'}{p}}(\gp{l,i}),\quad
    i = 0, \dotsc, 2^l,\;\;
    q = 0, \dotsc, \frac{p-1}{2}.
  \end{equation}
  Indeed, if \eqref{eq:proofPropInvariantHermiteHierarchization1} holds,
  then we obtain by
  \cref{line:algHermiteHierarchization6,line:algHermiteHierarchization8}
  of \cref{alg:hermiteHierarchization}
  \begin{subequations}
    \begin{align}
      \deriv[q]{x}{\fgintp{l}}(\gp{l,i})
      &= \deriv[q]{x}{\fgintp{l-1}}(\gp{l,i}) +
      \deriv[q]{x}{r^{(l)}_l}(\gp{l,i})\\
      &= \sum_{l'=0}^{l-1} \sum_{i' \in \hiset{l'}}
      \linout{l',i'} \deriv[q]{x}{\bspl[\wfs]{l',i'}{p}}(\gp{l,i}) +
      \sum_{i' \in \hiset{l}}
      \linout{l,i'} \deriv[q]{x}{\bspl[\wfs]{l,i'}{p}}(\gp{l,i})\\
      &= \sum_{l'=0}^l \sum_{i' \in \hiset{l'}}
      \linout{l',i'} \deriv[q]{x}{\bspl[\wfs]{l',i'}{p}}(\gp{l,i}),
    \end{align}
  \end{subequations}
  which is the desired relation
  \eqref{eq:propInvariantHermiteHierarchization}
  for level $l$.
  
  To prove \eqref{eq:proofPropInvariantHermiteHierarchization1},
  we separate two cases:
  \begin{itemize}
    \item
    \mbox{${i \notin \hiset{l}}$: In this case,}
    the ``true'' level of $\gp{l,i}$ is actually $\le l - 1$.
    Therefore, we can apply the induction hypothesis
    for \cref{eq:propInvariantHermiteHierarchization} to
    obtain \eqref{eq:proofPropInvariantHermiteHierarchization1}.
    
    \item
    \mbox{${i \in \hiset{l}}$: In this case,}
    we cannot directly apply the induction hypothesis,
    as it only holds for grid points of levels $\le l - 1$.
    However, we note that in
    \eqref{eq:proofPropInvariantHermiteHierarchization1},
    the term $\deriv[q]{x}{\fgintp{l-1}}(\gp{l,i})$
    is the $q$-th derivative of the Hermite interpolant of
    the data $\deriv[q']{x}{\fgintp{l-1}}(\gp{l,i\pm1})$
    ($q' = 0, \dotsc, \frac{p-1}{2}$),
    as determined in \cref{line:algHermiteHierarchization4}
    of \cref{alg:hermiteHierarchization}.
    The ``true'' level of the grid points $\gp{l,i\pm1}$ is
    actually $\le l - 1$ due to $i \in \hiset{l}$.
    Hence, we can apply the induction hypothesis
    for \cref{eq:propInvariantHermiteHierarchization}
    to conclude that the interpolated data of
    $\deriv[q]{x}{\fgintp{l-1}}(\gp{l,i})$ are given by
    \begin{equation}
      \deriv[q']{x}{\fgintp{l-1}}(\gp{l,i\pm1})
      = \deriv[q']{x}{
        \bracket*{
          \sum_{l'=0}^{l-1} \sum_{i' \in \hiset{l'}}
          \linout{l',i'} \bspl[\wfs]{l',i'}{p}
        }
      }(\gp{l,i\pm1}),\quad
      q' = 0, \dotsc, \frac{p-1}{2}.
    \end{equation}
    The linear combination in square brackets
    is a polynomial of degree $\le p$ on the interval
    $\clint{\gp{l,i-1}, \gp{l,i+1}}$
    by construction of the hierarchical basis functions
    $\bspl[\wfs]{l',i'}{p}$ ($l' = 0, \dotsc, l - 1$, $i' \in \hiset{l'}$).
    Due to the uniqueness of Hermite interpolation
    (\cref{lemma:hermiteInterpolation}),
    the interpolation polynomial of the data must
    coincide on $\clint{\gp{l,i-1}, \gp{l,i+1}}$
    with the term in square brackets.
    In particular, as $\gp{l,i} \in \clint{\gp{l,i-1}, \gp{l,i+1}}$,
    we obtain the claim \eqref{eq:proofPropInvariantHermiteHierarchization1}.
  \end{itemize}
  In both cases, we obtain the desired relation
  \eqref{eq:proofPropInvariantHermiteHierarchization1}.
\end{proof}

\corAlgHermiteHierarchizationCorrectness*

\begin{proof}
  By \cref{prop:invariantHermiteHierarchization} ($q = 0$), we have
  \begin{equation}
    \sum_{l'=0}^n \sum_{i' \in \hiset{l'}}
    \linout{l',i'} \bspl[\wfs]{l',i'}{p}(\gp{l,i})
    = \sum_{l'=0}^l \sum_{i' \in \hiset{l'}}
    \linout{l',i'} \bspl[\wfs]{l',i'}{p}(\gp{l,i})
    = \fgintp{l}(\gp{l,i}),\quad
    l \le n,\; i \in \hiset{l},
  \end{equation}
  as $\bspl[\wfs]{l',i'}{p}(\gp{l,i}) = 0$ if $l' > l$
  (weakly fundamental property \eqref{eq:weaklyFundamentalProperty}).
  \Cref{line:algHermiteHierarchization8} of
  \cref{alg:hermiteHierarchization} implies
  $\fgintp{l}(\gp{l,i})
  = \fgintp{l-1}(\gp{l,i}) + r^{(l)}_l(\gp{l,i})$, and by
  \cref{line:algHermiteHierarchization5,line:algHermiteHierarchization7},
  the second summand $r^{(l)}_l(\gp{l,i})$ equals
  $\fcnval{l,i} - \fgintp{l-1}(\gp{l,i})$,
  which cancels out the first summand, resulting in
  $\fgintp{l}(\gp{l,i}) = \fcnval{l,i}$.
  Combining these statements, we obtain
  \begin{equation}
    \sgintp(\gp{l,i})
    = \fcnval{l,i},\quad
    l \le n,\; i \in \hiset{l},
    \quad\text{where}\quad
    \sgintp
    \ceq \sum_{l'=0}^n \sum_{i' \in \hiset{l'}}
    \linout{l',i'} \bspl[\wfs]{l',i'}{p}.
  \end{equation}
  This means that $\sgintp$ is the correct hierarchical interpolant
  of the given function values
  (see \cref{eq:hierarchizationInterpolant}).
  Due to the uniqueness of hierarchical surpluses,
  the coefficients $\linout{l,i}$
  (which are the output of \cref{alg:hermiteHierarchization})
  must coincide with the surpluses $\surplus{l,i}$.
\end{proof}
}{}
\addtocontents{lop}{\protect\fi}

\cleardoublepage

  \chapter{Test Problems for Optimization}
\label{chap:a20testProblems}

In the following, we give formal definitions of the test problems
mentioned in \cref{sec:53testProblems}.
For each problem, we state the objective function
$\objfunscaled\colon \clint{\*a, \*b} \to \real$,
$\xscaled \mapsto \objfunscaled(\xscaled)$,
its domain $\clint{\*a, \*b}$
(where $\xscaled \in \clint{\*a, \*b}$),
the location $\xoptscaled \in \clint{\*a, \*b}$ of its global minimum, and
its minimal value $\objfunscaled(\xoptscaled)$
(and the constraint function
$\ineqconfunscaled\colon \clint{\*a, \*b} \to \real^{m_{\ineqconfunscaled}}$,
if any).
Plots of the test problems are given in
\cref{%
  fig:unconstrainedOptimizationProblem,%
  fig:constrainedOptimizationProblem%
}.

{
  \newcommand*{\centertestfunline}[1]{%
    \mathclap{\hphantom{\mathrm{(B.99a)}}#1}%
  }
  \let\xse\xscaledentry
  
  \section{Unconstrained Problems}
\label{sec:a21unconstrained}

\disableornamentsfornextheadingtrue
\vspace{-5mm}
\subsection{Bivariate Unconstrained Problems}
\label{sec:a211bivariateUnconstrained}

\paragraph{Branin02}

The function originates from \cite{Munteanu98Global}.
Compared to \cite{Munteanu98Global},
we changed the domain from $\clint{-5, 10} \times \clint{0, 15}$
to $\clint{-5, 15}^2$,
which seems more common in recent literature \cite{Gavana13Global}.
In addition, \cite{Munteanu98Global} uses the reciprocal function value,
while searching for the maximum instead of the minimum.
\vspace{-1.6em}

\begin{subequations}
  \begin{gather}
    \centertestfunline{
      \vspace*{-10mm}
      \testobjfunscaled{Bra02}(\xscaled)
      \ceq \paren*{
        -\,\frac{51\xse{1}^2}{40\pi^2} +
        \frac{5\xse{1}}{\pi} + \xse{2} - 6
      }^2 +
      \paren*{10 - \frac{5}{4\pi}} \cos(\xse{1}) \cos(\xse{2})
    }\\
    \centertestfunline{
      {} + \ln(\xse{1}^2 + \xse{2}^2 + 1) + 10,\hspace*{35mm}
    }\notag\\
    \centertestfunline{
      \xscaled \in \clint{-5, 15}^2,\quad
      \xoptscaled = (-3.196988424804, 12.52625788532),
    }\\
    \centertestfunline{
      \testobjfunscaled{Bra02}(\xoptscaled) = 5.558914403894
    }
  \end{gather}
\end{subequations}

\pagebreak

\paragraph{GoldsteinPrice}

This function originates from \cite{Goldstein71Descent},
where the function was stated without bounds for the optimization domain.
We took the domain $\clint{-2, 2}^2$ from \cite{Gavana13Global}.
In addition, we scaled the function values by the factor $10^{-4}$
for the sake of plotting.
\vspace{-1.6em}

\begin{subequations}
  \begin{gather}
    \centertestfunline{
      \testobjfunscaled{GoP}(\xscaled)
      \ceq 10^{-4} \cdot \paren*{
        1 + (\xse{1} + \xse{2} + 1)^2
        (19 - 14\xse{1} + 3\xse{1}^2 - 14\xse{2} +
        6\xse{1}\xse{2} + 3\xse{2}^2)
      }
    }\\
    \centertestfunline{
      \hspace*{25mm}
      {} \cdot
      \paren*{
        30 + (2\xse{1} - 3\xse{2})^2
        (18 - 32\xse{1} + 12\xse{1}^2 + 48\xse{2} -
        36\xse{1}\xse{2} + 27\xse{2}^2)
      },
    }\notag\\
    \centertestfunline{
      \xscaled \in \clint{-2, 2}^2,\quad
      \xoptscaled = (0, -1),
    }\\
    \centertestfunline{
      \testobjfunscaled{GoP}(\xoptscaled) = 3 \cdot 10^{-4}
    }
  \end{gather}
\end{subequations}

\paragraph{Schwefel06}

This function originates from \cite{Schwefel77Numerische}.
We changed the domain from $\clint{-3, 5} \times \clint{-1, 7}$ to
$\clint{-6, 4}^2$, such that the optimum point is not located at the
center of the optimization domain.
\vspace{-1.6em}

\begin{subequations}
  \begin{gather}
    \centertestfunline{
      \testobjfunscaled{Sch06}(\xscaled)
      \ceq \max(
        \abs{\xse{1} + 2\xse{2} - 7},
        \abs{2\xse{1} + \xse{2} - 5}
      ),
    }\\
    \centertestfunline{
      \xscaled \in \clint{-6, 4}^2,\quad
      \xoptscaled = (1, 3),\quad
      \testobjfunscaled{Sch06}(\xoptscaled) = 0
    }
  \end{gather}
\end{subequations}

\subsection{\texorpdfstring{$d$}{d}-Variate Unconstrained Problems}
\label{sec:a212dvariateUnconstrained}

\paragraph{Ackley}

The form of this function originates from \cite{Ackley87Connectionist},
where it was stated only for two variables.
We use the generalization to $d$ variables from \cite{Gavana13Global}.
The optimization domain $\clint{1.5, 6.5}^d$
was chosen such that it does not contain $\*0$,
where the gradient of the objective function becomes singular.
Otherwise, the function would not be continuously differentiable,
which would be a disadvantage for spline-based approaches
(see Schwefel06 and Schwefel22 for functions with discontinuous derivatives).
\vspace{-1.6em}

\begin{subequations}
  \begin{gather}
   \centertestfunline{
      \testobjfunscaled{Ack}(\xscaled)
      \ceq -20 \exp\paren*{-\frac{\norm[2]{\xscaled}}{5\sqrt{d}}} -
      \exp\paren*{\frac{1}{d} \sum_{t=1}^d \cos(2\pi \xse{t})} +
      20 + \econst,
    }\\
    \centertestfunline{
      \xscaled \in \clint{1.5, 6.5}^d,\quad
      \xoptscaled = 1.974451986484 \cdot \*1,\quad
      \testobjfunscaled{Ack}(\xoptscaled) = 6.559645375628
    }
  \end{gather}
\end{subequations}

\paragraph{Alpine02}

This function originates from \cite{Clerc99Swarm}.
We changed the domain from $\clint{0, 10}^d$ to $\clint{2, 10}^d$
to exclude the singularities of the derivative of the objective function
at $\xse{t} = 0$.
In addition, the author of \cite{Clerc99Swarm} searched for maximal points.
For minimization, we changed the sign of the objective function.
\vspace{-1.6em}

\begin{subequations}
  \begin{gather}
    \centertestfunline{
      \testobjfunscaled{Alp02}(\xscaled)
      \ceq -\prod_{t=1}^d \sqrt{\xse{t}} \sin(\xse{t}),\qquad
      \xscaled \in \clint{2, 10}^d,
    }\\
    \centertestfunline{
      \xoptscaled = 7.917052684666 \cdot \*1,\quad
      \testobjfunscaled{Alp02}(\xoptscaled) = -2.808131180070^d
    }
  \end{gather}
\end{subequations}

\paragraph{Schwefel22}

This function originates from \cite{Schwefel77Numerische}.
We changed the domain from $\clint{-10, 10}^d$ to
$\clint{-3, 7}^d$, such that the optimum point is not located at the
center of the optimization domain.
\vspace{-1.6em}

\begin{subequations}
  \begin{gather}
    \centertestfunline{
      \testobjfunscaled{Sch22}(\xscaled)
      \ceq \sum_{t=1}^d \abs{\xse{t}} +
      \prod_{t=1}^d \abs{\xse{t}},\qquad
      \xscaled \in \clint{-3, 7}^d,
    }\\
    \centertestfunline{
      \xoptscaled = \*0,\quad
      \testobjfunscaled{Sch22}(\xoptscaled) = 0
    }
  \end{gather}
\end{subequations}

  \section{Constrained Problems}
\label{sec:a22constrained}

\paragraph{G08}

This problem originates from \cite{Schoenauer93Constrained}.
We changed the domain from $\clint{0, 10}^2$
to $\clint{0.5, 2.5} \times \clint{3, 6}$
to increase the size of feasible region.
In addition, we use different frequencies for the sine terms
as in \cite{Gavana13Global}.
\vspace{-1.6em}

\begin{subequations}
  \begin{gather}
    \centertestfunline{
      \testobjfunscaled{G08}(\xscaled)
      \ceq -\frac{
        \sin^3(2\pi\xse{1}) \sin(2\pi\xse{2})
      }{
        \xse{1}^3 (\xse{1} + \xse{2})
      },\quad
      \testineqconfunscaled{G08}(\xscaled)
      \ceq \begin{pmatrix}
        \xse{1}^2 - \xse{2} + 1\\
        1 - \xse{1} + (\xse{2} - 4)^2
      \end{pmatrix},
    }\\
    \centertestfunline{
      \xscaled \in \clint{0.5, 2.5} \times \clint{3, 6},\quad
      \xoptscaled = (1.227971358337, 4.245373366474),
    }\\
    \centertestfunline{
      \testobjfunscaled{G08}(\xoptscaled) = -0.09582504141804
    }
  \end{gather}
\end{subequations}

\paragraph{G04Squared}

This problem is based on a problem from
\cite{Colville68Comparative} with the objective function
$\testobjfunscaled{G04}(\xscaled)
\ceq 5.3578547 \xse{3}^2 + 0.8356891 \xse{1} \xse{5} +
37.293239 \xse{1} - 40792.141$ and the same constraints
$\testineqconfunscaled{G04}(\xscaled) \ceq
\testineqconfunscaled{G04Sq}(\xscaled)$.
However, hierarchical cubic not-a-knot B-splines are able to exactly
represent the polynomial $\testobjfunscaled{G04}$ of coordinate degree two
on the whole domain $\clint{\*0, \*1}$,
if the level of the sparse grids is high enough,
see \thmref{cor:sparseGridRegularNAKPolynomials}.
Therefore, we modified the original G04 problem by squaring the
objective function.
To ensure that this does not change the location of the global minimum,
we added a constant before squaring such that the shifted function
is non-negative on $\clint{\*0, \*1}$.
\vspace{-1.6em}

\begin{subequations}
  \allowdisplaybreaks
  \begin{gather}
  \centertestfunline{
    \testobjfunscaled{G04Sq}(\xscaled)
    \ceq (5.3578547 \xse{3}^2 + 0.8356891 \xse{1} \xse{5} +
    37.293239 \xse{1} - 10120)^2,
  }\\
  \centertestfunline{
    \testineqconfunscaled{G04Sq}(\xscaled)
    \ceq 10^{-3} \scalebox{0.92}{$
      \begin{pmatrix}
        85334.407 + 5.6858 \xse{2} \xse{5} +
        0.6262 \xse{1} \xse{4} -
        2.2053 \xse{3} \xse{5} - 92000\\
        -85334.407 - 5.6858 \xse{2} \xse{5} -
        0.6262 \xse{1} \xse{4} +
        2.2053 \xse{3} \xse{5}\\
        80512.49 + 7.1317 \xse{2} \xse{5} +
        2.9955 \xse{1} \xse{2} +
        2.1813 \xse{3}^2 - 110000\\
        -80512.49 - 7.1317 \xse{2} \xse{5} -
        2.9955 \xse{1} \xse{2} -
        2.1813 \xse{3}^2 + 90000\\
        9300.961 + 4.7026 \xse{3} \xse{5} +
        1.2547 \xse{1} \xse{3} +
        1.9085 \xse{3} \xse{4} - 25000\\
        -9300.961 - 4.7026 \xse{3} \xse{5} -
        1.2547 \xse{1} \xse{3} -
        1.9085 \xse{3} \xse{4} + 20000
      \end{pmatrix},
    $}
  }\\
  \centertestfunline{
    \xscaled \in \clint{78, 102} \times \clint{33, 45} \times
    \clint{27, 45}^3,
  }\\
  \centertestfunline{
    \xoptscaled = (78, 33, 29.995256025682, 45, 36.775812905788),
  }\\
  \centertestfunline{
    \testobjfunscaled{G04Sq}(\xoptscaled) = 43.590737882363
  }
  \end{gather}
\end{subequations}

}

\cleardoublepage

  \longchapter{%
  Detailed Results for Topology Optimization%
}{%
  Detailed Results for Topology Optimization%
}{%
  Detailed Results (Topology Optimization)%
}
\label{chap:a30topoOptDetails}

\begin{figure}
  \subcaptionbox{%
    2D cross%
  }[72mm]{%
    \includegraphics{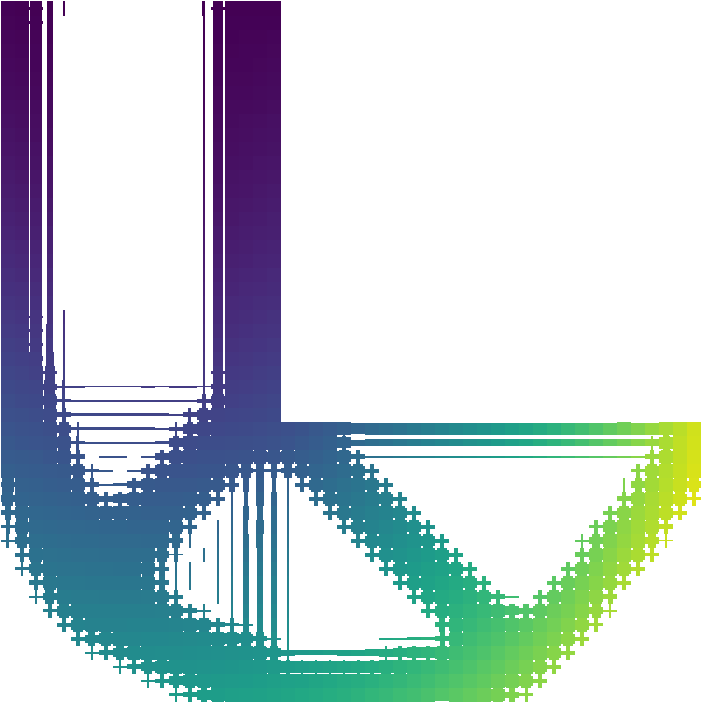}%
  }%
  \hfill%
  \subcaptionbox{%
    2D framed cross%
  }[72mm]{%
    \includegraphics{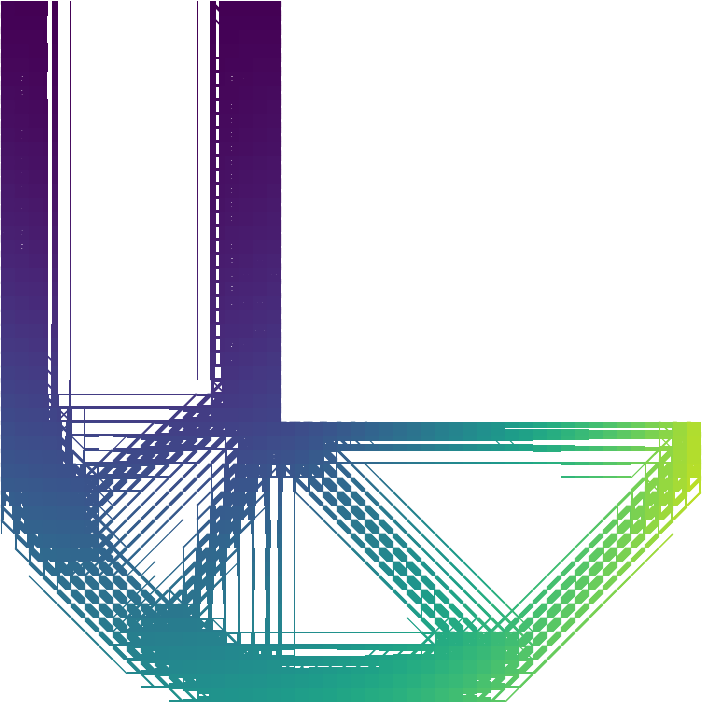}%
  }%
  \\[2mm]%
  \subcaptionbox{%
    2D sheared cross%
  }[72mm]{%
    \includegraphics{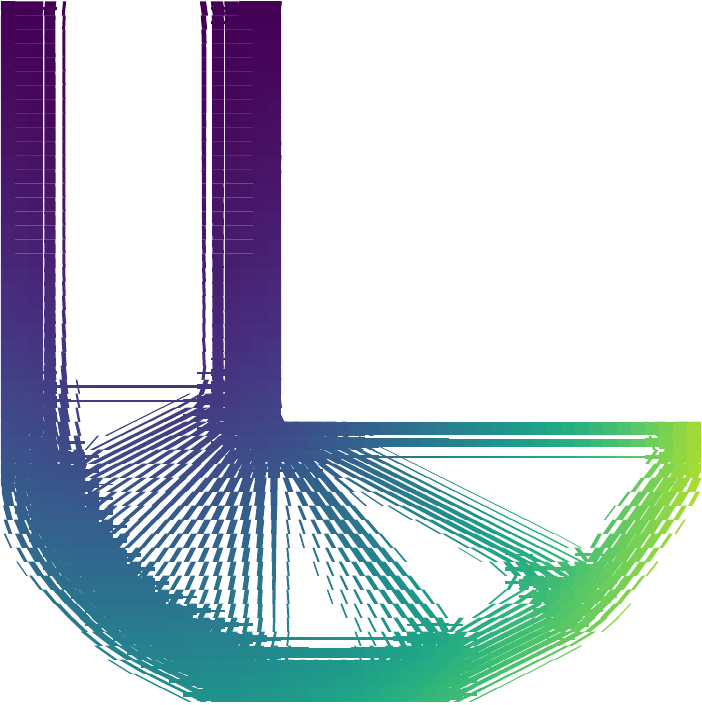}%
  }%
  \hfill%
  \subcaptionbox{%
    2D sheared framed cross%
  }[72mm]{%
    \includegraphics{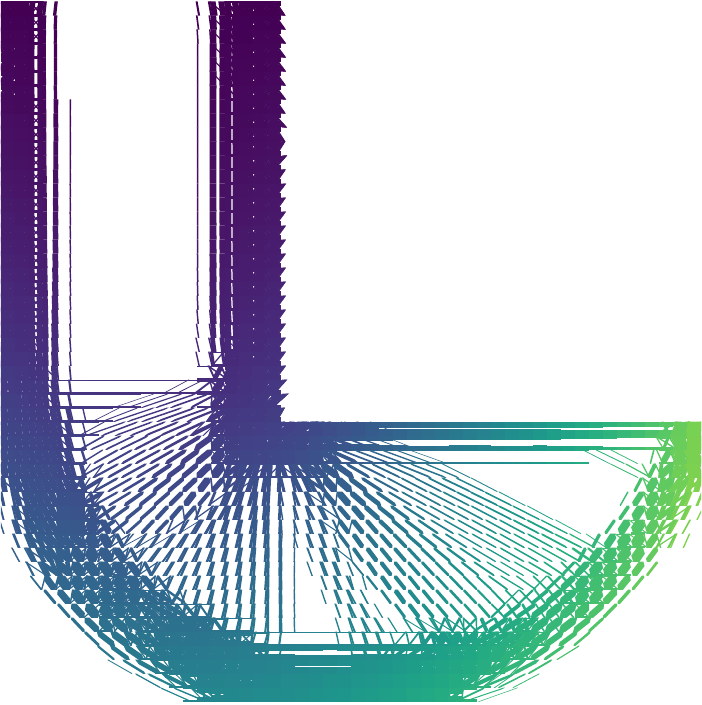}%
  }%
  \caption[Optimal structures in the 2D L-shape scenario]{%
    Topologically optimal structures in the 2D L-shape scenario
    for different micro-cell models using cubic B-splines
    (spatially adaptive grids with around \num{10000} points).
    The colors indicate the length of the displacement,
    where dark regions correspond to weak displacements and
    bright regions to strong displacements.
    The color map is the same as in
    \cref{fig:topoOptStructure2DCantilever}.
    Only bars with widths $\ge 0.1$ are shown.
    More details can be found in \cref{tbl:topoOptResultsDetailed}.%
  }%
  \label{fig:topoOptStructure2DLShape}%
\end{figure}

\begin{figure}
  \newcommand*{\myscale}{0.38}%
  \subcaptionbox{%
    3D cantilever, 3D cross%
  }[72mm]{%
    \includegraphics[scale=\myscale]{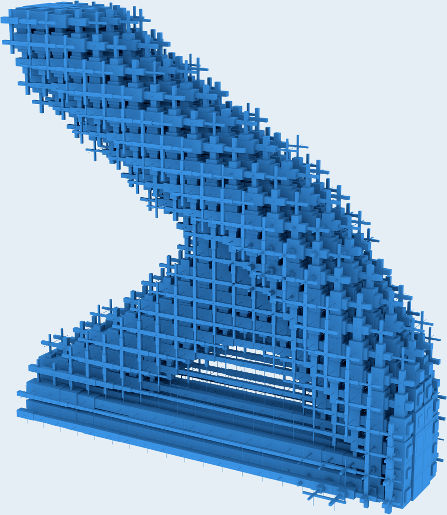}%
  }%
  \hfill%
  \subcaptionbox{%
    3D cantilever, 3D sheared cross%
  }[72mm]{%
    \includegraphics[scale=\myscale]{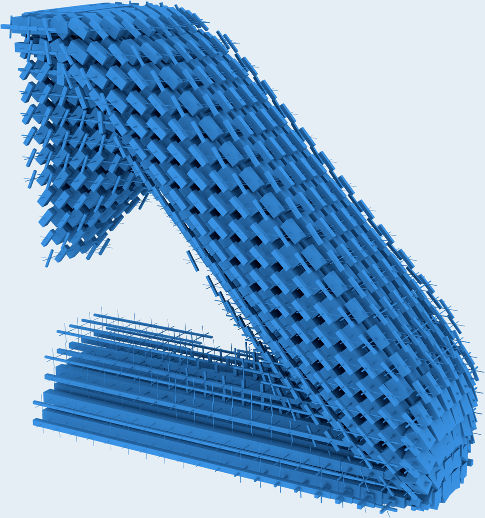}%
  }%
  \\[2mm]%
  \begin{tikzpicture}
    \node[anchor=south west] at (0mm,61mm) {%
      \includegraphics[scale=\myscale]{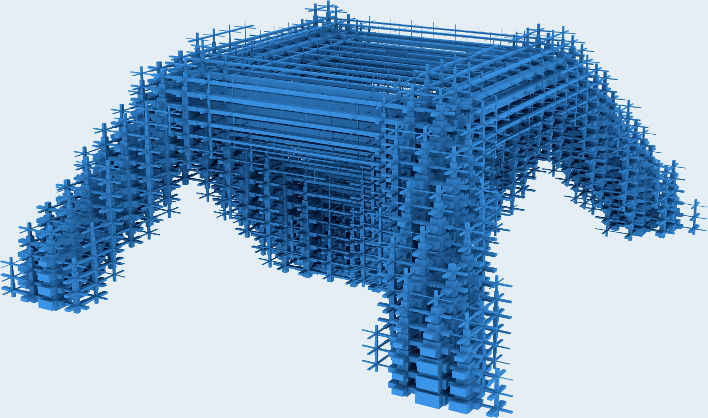}%
    };
    \node[anchor=south west] at (0mm,0mm) {%
      \includegraphics[scale=\myscale]{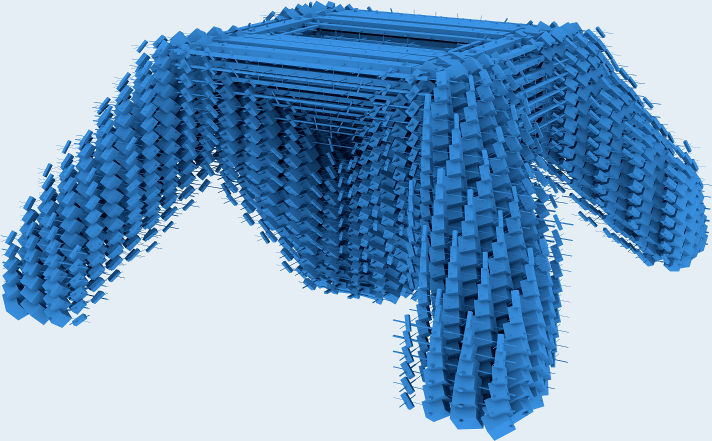}%
    };
    \node[anchor=south west,text width=50mm] at (-1mm,66mm) {%
      \subcaption{3D center-load, 3D cross}%
    };
    \node[anchor=south west,text width=60mm] at (-8mm,5mm) {%
      \subcaption{3D center-load, 3D sheared cross}%
    };
  \end{tikzpicture}
  \hspace*{6mm}
  \caption[Optimal structures in the 3D scenarios]{%
    Topologically optimal structures in the
    3D cantilever and center-load scenarios
    for different micro-cell models using cubic B-splines
    (spatially adaptive grids with around \num{10000} points).
    More details can be found in \cref{tbl:topoOptResultsDetailed}.%
  }%
  \label{fig:topoOptStructure3D}%
\end{figure}

\noindent
This appendix complements \cref{sec:64results}.
It contains visualizations of the topologically optimal structures
in the 2D L-shape and the 3D scenarios in
\cref{fig:topoOptStructure2DLShape} and \cref{fig:topoOptStructure3D},
respectively.
In addition, we report details of the corresponding optimization runs
in \cref{tbl:topoOptResultsDetailed} and
information about the employed spatially adaptive sparse grids in
\cref{tbl:topoOptResultsModels}.
The computation times were measured on a shared-memory computer
with 4x Intel Xeon E7-8880v3 (72 cores, 144 threads).

\begin{table}
  \sisetup{retain-zero-exponent=true}%
  \newcommand*{\cece}[1]{\multicolumn{1}{c}{#1}}%
  \newcommand*{\trcmp}{\compliance[\opt,\ast]}%
  \newcommand*{\apcmp}{\complianceintp[\opt,\ast]}%
  \newcommand*{\iic}{\multirow{-4}{*}{2D cantilever}}%
  \newcommand*{\iil}{\multirow{-4}{*}{2D L-shape}}%
  \newcommand*{\iiic}{\multirow{-2}{*}{3D cantilever}}%
  \newcommand*{\iiicl}{\multirow{-2}{*}{3D center-load}}%
  \setnumberoftableheaderrows{1}%
  \begin{tabular}{%
    >{\kern\tabcolsep}=l+l<{\kern5mm}+r+c+c+l+r<{\kern\tabcolsep}%
  }
    \toprulec
    \headerrow
    Scenario& Model& \#Iter.& $\trcmp$& $\apcmp$& \cece{O.-i. gap}& Time\\
    \midrulec
    &         2D-C&      547& 74.974&   74.974&   \num{3.67e-5}&    \hms{;5}\\
    &         2D-FC&     249& 70.816&   69.409&   \num{1.41e0}&     \hms{;14}\\
    &         2D-SC&    2196& 67.809&   67.804&   \num{5.21e-3}&    \hms{1;6}\\
    \iic&     2D-SFC&    749& 68.602&   65.201&   \num{3.40e0}&     \hms{;15}\\
    \midrulec
    &         2D-C&      289& 183.68&   183.68&   \num{9.06e-5}&    \hms{;3}\\
    &         2D-FC&     602& 177.51&   174.49&   \num{3.02e0}&     \hms{;17}\\
    &         2D-SC&    1609& 169.60&   169.60&   \num{7.33e-3}&    \hms{;33}\\
    \iil&     2D-SFC&    574& 174.55&   158.19&   \num{1.64e1}&     \hms{;8}\\
    \midrulec
    &         3D-C&       39& 247.60&   247.49&   \num{1.13e-1}&    \hms{;10}\\
    \iiic&    3D-SC&     608& 162.59&   159.33&   \num{3.25e0}&     \hms{3;17}\\
    \midrulec
    &         3D-C&       35& 169.27&   169.27&   \num{3.31e-3}&    \hms{;4}\\
    \iiicl&   3D-SC&    1026& 46.171&   45.571&   \num{6.00e-1}&    \hms{2;25}\\
    \bottomrulec
  \end{tabular}
  \caption[Details about topology optimization runs]{%
    Detailed information about the optimization runs corresponding to
    \cref{tbl:topoOptResultsModels} and
    \cref{%
      fig:topoOptStructure2DCantilever,%
      fig:topoOptStructure2DLShape,%
      fig:topoOptStructure3D%
    }, which employs cubic B-splines ($p = 3$) on the spatially
    adaptive grids listed in \cref{tbl:topoOptResultsModels}.
    From left to right, the columns contain
    the optimization scenario,
    the micro-cell model,
    the number of optimization iterations,
    the actual compliance value
    $\trcmp \ceq \compliance(\mcpoptappr{1}, \dotsc, \mcpoptappr{M})$,
    the approximated compliance value
    $\apcmp \ceq \complianceintp(\mcpoptappr{1}, \dotsc, \mcpoptappr{M})$
    as reported by the optimizer,
    the optimality-interpolation gaps $\abs{\trcmp - \apcmp}$, and
    the computation time of the online phase
    (without the time to generate the sparse grid data).%
  }%
  \label{tbl:topoOptResultsDetailed}%
\end{table}

\begin{table}
  \sisetup{retain-zero-exponent=true}%
  \newcommand*{\cece}[1]{\multicolumn{1}{c}{#1}}%
  \setnumberoftableheaderrows{1}%
  \begin{tabular}{%
    >{\kern\tabcolsep}=l<{\kern5mm}+c+r+l+l+c<{\kern\tabcolsep}%
  }
    \toprulec
    \headerrow
    Model&  $d$& $\ngp$&      \cece{Threshold}& \cece{Rel. err.}& Eval. time\\
    \midrulec
    2D-C&   2&   \num{10197}&    \num{2.15e-5}&    \num{2.24e-5}& \SI{6.96}{\second}\\
    2D-FC&  4&   \num{10502}&    \num{7.94e-1}&    \num{1.81e-2}& \SI{7.45}{\second}\\
    2D-SC&  3&   \num{10723}&    \num{4.64e-3}&    \num{1.49e-3}& \SI{8.72}{\second}\\
    2D-SFC& 5&   \num{10694}&    \num{5.01e0}&     \num{4.82e-2}& \SI{7.45}{\second}\\
    \midrulec
    3D-C&   3&    \num{9207}&    \num{7.94e-2}&    \num{3.18e-3}& \SI{33.5}{\second}\\
    3D-SC&  5&   \num{15389}&    \num{5.01e0}&     \num{4.95e-2}& \SI{40.0}{\second}\\
    \bottomrulec
  \end{tabular}
  \caption[%
    Details about spatially adaptive sparse grids for topology optimization%
  ]{%
    Detailed information about the spatially adaptive sparse grids used for
    \cref{tbl:topoOptResultsModels,tbl:topoOptResultsDetailed} and
    \cref{%
      fig:topoOptStructure2DCantilever,%
      fig:topoOptStructure2DLShape,%
      fig:topoOptStructure3D%
    }.
    The columns correspond to
    the micro-cell model,
    the number $d$ of micro-cell parameters,
    the number $\ngp$ of sparse grid points,
    the threshold $\refinetol$ used in the grid generation algorithm,
    the relative $\Ltwo$ spectral interpolation error
    $
      \normLtwoscaled{
        \vphantom{\big(}
        \norm[2]{\etensor({\cdot}) - \etensorcholintp({\cdot})}
      }
      /
      \normLtwoscaled{
        \vphantom{\big(}
        \norm[2]{\etensor({\cdot})}
      }
    $, and the
    time needed to evaluate the elasticity tensor $\etensor(\vgp{k})$
    at a single grid point $\vgp{k}$.%
  }%
  \label{tbl:topoOptModels}%
\end{table}

\cleardoublepage

  \printbibliography[%
  title=Bibliography,
  heading=myheading,
  postnote=mypostnote,
]

\iftoggle{partialCompileMode}{
  \cleardoublepage
}{
  \cleardoublepage
  \thispagestyle{empty}
  \null
  \cleardoubleevenpage
  
  \thispagestyle{empty}
  \newgeometry{
    bindingoffset=0mm,
    inner=0mm,
    outer=0mm,
    top=0mm,
    bottom=0mm,
  }
  
  \vspace*{\fill}
  
  \begin{center}
    \includegraphics[scale=2]{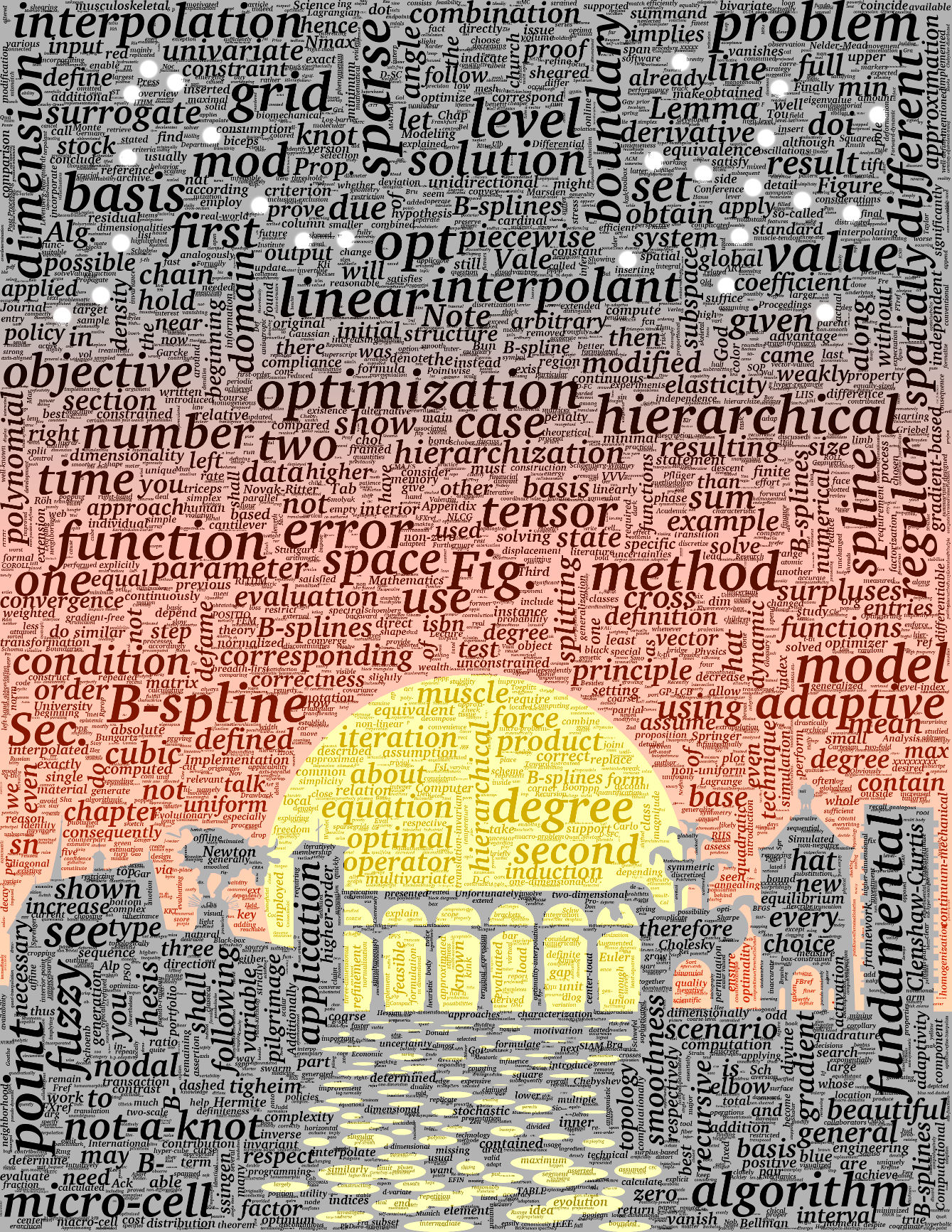}
  \end{center}
  
  \vspace*{\fill}
}

\end{document}